\documentclass[11pt,reqno]{amsart}
\usepackage[T1]{fontenc}
\usepackage{lmodern}
\usepackage{geometry}
\usepackage{amsmath,amssymb,mathtools,mathrsfs}
\usepackage{booktabs,longtable,array}
\usepackage{needspace}
\usepackage{enumitem}
\usepackage{xcolor}
\usepackage{tikz}
\usetikzlibrary{arrows.meta,positioning,fit,shapes.geometric}
\usepackage[hidelinks]{hyperref}
\usepackage{xr-hyper}

\numberwithin{equation}{section}
\setlist[itemize]{leftmargin=2em,itemsep=2pt,topsep=4pt}
\setlist[enumerate]{leftmargin=2.2em,itemsep=2pt,topsep=4pt}

\newtheorem{theorem}{Theorem}[section]
\newtheorem{proposition}[theorem]{Proposition}
\newtheorem{lemma}[theorem]{Lemma}
\newtheorem{corollary}[theorem]{Corollary}

\theoremstyle{definition}
\newtheorem{definition}[theorem]{Definition}
\newtheorem{remark}[theorem]{Remark}

\newcommand{\eps}{\varepsilon}
\newcommand{\epstar}{\varepsilon_*}
\newcommand{\T}{\mathbb T}
\newcommand{\R}{\mathbb R}
\newcommand{\Sph}{\mathbb S}
\newcommand{\N}{\mathbb N}
\newcommand{\Z}{\mathbb Z}
\newcommand{\1}{\mathbf 1}
\newcommand{\cD}{\mathcal D}
\newcommand{\cZ}{\mathcal Z}

\newcommand{\cL}{\mathcal L}
\newcommand{\Exc}{\operatorname{Exc}}

\newcommand{\Err}{\operatorname{Err}}
\newcommand{\Jac}{\operatorname{Jac}}
\newcommand{\dd}{\,\mathrm d}
\newcommand{\tensor}{\otimes}

\newcommand{\mainproofheading}[2]{%
  \par\medskip\noindent\textbf{Proof of #1.}\label{#2}\par\smallskip}
\newcommand{\routineproofinresource}[1]{%
  \par\smallskip\noindent\emph{Supplementary proof.}
  Online Resource~1, Supplementary Proof~\ref{#1}, gives the detailed
  chartwise calculation for the statement above.\par\smallskip}
\newcommand{\fallbackproofinresource}[1]{%
  \par\smallskip\noindent\emph{Complementary auxiliary route.}
  Online Resource~1, Supplementary Proof~\ref{#1}, gives the complete proof
  of this coarse or conditional estimate.\par\smallskip}
\newcommand{\prooflaterinmain}[1]{%
  \par\smallskip\noindent\emph{Proof.}
  The complete argument is given at Proof~\ref{#1}.\par\smallskip}

\newcommand{\mainstatementreference}[1]{%
  \par\medskip\noindent\emph{Main-text statement.}
  The statement of #1 appears in the main manuscript; the detailed proof
  follows here.\par\smallskip}
\newcounter{suppproof}[section]
\renewcommand{\thesuppproof}{\thesection.P\arabic{suppproof}}
\newcommand{\suppproofheading}[2]{%
  \refstepcounter{suppproof}\par\medskip\noindent
  \textbf{Supplementary Proof \thesuppproof\space (#1).}%
  \label{#2}\par\smallskip}
\makeatletter
\newcommand{\clearxrtocindents}{%
  \expandafter\let\csname r@tocindent-1\endcsname\relax
  \expandafter\let\csname r@tocindent0\endcsname\relax
  \expandafter\let\csname r@tocindent1\endcsname\relax
  \expandafter\let\csname r@tocindent2\endcsname\relax
  \expandafter\let\csname r@tocindent3\endcsname\relax}
\makeatother
\hypersetup{
  pdftitle={A Periodic Long-Time Boltzmann--Grad Limit in Every Fixed Dimension d at Least 4},
  pdfauthor={Zhixu Hua},
  pdfsubject={A periodic long-time Boltzmann--Grad limit for hard spheres in every fixed dimension d at least 4},
  pdfkeywords={Boltzmann--Grad limit, hard spheres, periodic dynamics, cumulants, kinetic theory}
}
\title[Periodic Boltzmann--Grad limit in every fixed $d\ge4$]
{A Periodic Long-Time Boltzmann--Grad Limit\\
in Every Fixed Dimension $d\ge4$}
\author{Zhixu Hua}
\thanks{ORCID: \href{https://orcid.org/0009-0004-0951-1170}{0009-0004-0951-1170}.}
\address{Changkong College, Nanjing University of Aeronautics and Astronautics,
Nanjing, Jiangsu, China}
\email{hzx001@nuaa.edu.cn}
\date{September 3, 2026}
\begin{document}
\begin{abstract}
We prove a periodic long-time Boltzmann--Grad limit for hard spheres in every
fixed spatial dimension $d\ge4$.  Under the assumptions of the main theorem,
the rescaled $s$-particle correlations converge in $L^1$ to the corresponding
tensor-product Boltzmann profile at rate $\varepsilon^{1/(400d)}$, uniformly
for $1\le s\le|\log\varepsilon|$ and $0\le t\le t_{\rm fin}$.  In particular,
when the activity and weighted solution bounds are fixed,
$t_{\rm fin}=O(\log|\log\varepsilon|)$.

The componentwise long-bond estimate used in dimensions two and three is
insufficient in higher dimension.  We replace it by a joint estimate for two
connected time sublayers, selecting the first two lower collisions as landing
roots.  Disjoint landing edges yield a direct tangential frame.  When the
edges overlap, the newly appearing particle line contains a collision-free
chord whose radial relative speed produces a rank-$(d-1)$ positive Schur
factor.  A positive Jacobi-network argument prevents focusing, while
bounded-degree finite-fibre coarea globalizes the estimate for arbitrary
nonnegative joint kernels.  A first-failure construction repairs periodic
double overlaps, and the resulting packet bound closes the exceptional
top-layer contribution.

A paid epoch restart for sealed complete joint kernels, using the fixed-word
multi-landing operator for every fixed integer $k\ge1$, yields the stated time
range.  Within the regular connected full two-sublayer packet class, every
packet with at least $2k-1$ physical lines admits a canonical $k$-birth flag
and the operator factor
$\varepsilon^{(d-1)(k-1)}\varepsilon_*^{-k(d-1)}$; $k=2$ is the first
gain-producing case.  A complementary full-chord estimate provides a coarse fallback.
\end{abstract}
\maketitle
\tableofcontents
\section{Introduction and main result}
\subsection{Problem and context}

The Boltzmann--Grad limit relates a dilute Newtonian hard-sphere gas to the
Boltzmann kinetic equation.  Its kinetic scaling goes back to Grad's
formulation \cite{Grad1958}; standard phase-space and dilute-gas backgrounds
are \cite{Spohn1991,CercignaniIllnerPulvirenti1994}.  Lanford's theorem gives
the limit for a short time \cite{Lanford1975}; a broad account of the
classical hierarchy and its geometric recollision problem is given in
\cite{GST2014}.  Deng, Hani and Ma
developed a long-time expansion on Euclidean space and then a periodic
version in dimensions two and three \cite{DHMlong,DHMtorus}.  Here we close
the corresponding periodic argument in every fixed dimension $d\ge4$.

Longer-time and fluctuation analyses of hard-sphere systems also develop
cluster, correlation-error, and tagged-particle structures that provide the
broader mathematical context \cite{BGS2016,PulvirentiSimonella2017,
BGSScluster2022,BGSS2023}.

The precise version-locked interfaces imported from
\cite{DHMlong,DHMtorus} are stated in Section~\ref{sec:setup} and documented
in Online Resource~1.

\subsection{Main theorem}

Put
\[
 \cD_s^\circ=\{(x_1,v_1,\dots,x_s,v_s):
 d_{\T^d}(x_i,x_j)>\eps\ \text{for }i\ne j\}.
\]
The grand-canonical ensemble and the rescaled correlations are defined in
Section~\ref{sec:setup}.

\begin{theorem}[Periodic long-time limit in fixed dimension $d\ge4$]
\label{thm:main}
Fix an integer $d\ge4$ and $\beta>0$.  There are constants
$c_0=c_0(d,\beta)>0$ and
$0<\eps_0=\eps_0(d,\beta)<e^{-e}$ with the following property.  For every
$0<\eps\le\eps_0$, let
$\alpha,A,t_{\rm fin}>0$ satisfy
\begin{equation}\label{eq:loglog-main}
 \max(1,\alpha)\max(1,A)
 \le c_0\sqrt{\log|\log\eps|},
 \qquad
 \max(1,\alpha)\max(1,A)\max(1,t_{\rm fin})
 \le c_0\log|\log\eps|.
\end{equation}
Let $n_0\ge0$ on $\T^d\times\R^d$, with $\int n_0=1$.  Suppose that the
hard-sphere Boltzmann equation
\begin{equation}\label{eq:boltzmann-main}
 (\partial_t+v\cdot\nabla_x)n
 =\alpha\int_{\R^d}\int_{\Sph^{d-1}}
 ((v-v_1)\cdot\omega)_+(n'n_1'-nn_1)\dd\omega\dd v_1
\end{equation}
has a nonnegative, mass-conserving mild solution $n$ on
$[0,t_{\rm fin}]$.  Precisely, if $\mathcal Q(n,n)$ denotes the collision
integral on the right of \eqref{eq:boltzmann-main} without the factor
$\alpha$, then, for every $0\le t\le t_{\rm fin}$, the following identity
holds for almost every $(x,v)$:
\begin{equation}\label{eq:boltzmann-mild}
 n(t,x,v)=n_0(x-tv,v)
 +\alpha\int_0^t
 \mathcal Q(n,n)(s,x-(t-s)v,v)\,\dd s.
\end{equation}
Its initial trace is
\[
 n(0,x,v)=n_0(x,v),
\]
satisfying
\begin{equation}\label{eq:weighted-assumptions}
 \sup_{0\le t\le t_{\rm fin}}
 \|e^{2\beta|v|^2}n(t)\|_\infty\le A,
 \qquad
 \|e^{2\beta|v|^2}\nabla_xn_0\|_\infty\le A.
\end{equation}
Start diameter-$\eps$ hard spheres from the grand-canonical ensemble of
activity $\alpha\eps^{-(d-1)}$ and define the $s$-point rescaled correlations
with factor $(\alpha^{-1}\eps^{d-1})^s$.  Then, uniformly
for $0\le t\le t_{\rm fin}$ and $1\le s\le|\log\eps|$,
\begin{equation}\label{eq:main-convergence}
 \left\|f_s(t)-n(t)^{\tensor s}\1_{\cD_s^\circ}\right\|_
 {L^1((\T^d\times\R^d)^s)}
 \le \eps^{1/(400d)}.
\end{equation}
The microscopic flow is understood outside its invariant singular null set.
\end{theorem}

After decreasing \(c_0\) once, this parameter set contains the previous
single-epoch range: its old hypothesis implies the first inequality in
\eqref{eq:loglog-main}, and the old \(O(\sqrt{\log|\log\eps|})\) product
bound implies the second because \(\log|\log\eps|\ge1\).  Thus the new
statement extends the time range without exchanging away an old admissible
triple.

Before constants growing with the number of time layers are absorbed, the
version-2 parameter is $1/(200d)$; the theorem states the resulting uniform
exponent.

\subsection{Scope of the main result}

The proof combines the version-locked long-time reduction of
Section~\ref{sec:setup} with four ingredients established here: the
higher-dimensional packet estimate, the periodic first-failure repair, the
exceptional-operator insertion, and the epoch restart.  The restart uses the
fixed-word result under its stated packet hypotheses with $k$ fixed
independently of $\eps$.  The sharp radial-conjugate route is complemented by
a coarse full-chord construction.  Local rank is globalized by bounded finite
fibres and branch-complete coarea for the same arbitrary nonnegative joint
kernel.

\subsection{Supplementary material}

\mbox{}\par
\smallskip\noindent\textbf{Online Resource~1.} Expanded routine coordinate verifications, source-version concordance, alternate packet charts, and supplementary parameter and compatibility tables.

Online Resource 1 is submitted with the manuscript and is part of its proof
record.  It expands routine coordinate verifications, records the exact
source-version concordance and equation-level compatibility checks, and
preserves independent coarse or conditional routes.  The positive-network
and finite-fibre modules are proved in
Sections~\ref{app:positive-network}--\ref{app:coarea}, and the parameter
arithmetic is proved in Appendix~\ref{app:parameters}.  The principal new load-bearing arguments are developed in the main
manuscript; Online Resource~1 collects routine coordinate verifications,
source concordance, and auxiliary coarse routes.

\subsection{Organization}

Section~\ref{sec:overview} gives a detailed proof overview.  The imported
reduction is isolated in Section~\ref{sec:setup}; the microscopic flow and initial correlations are established in
Section~\ref{sec:ensemble-flow}; the ordinary fixed-dimensional estimates
are proved in Section~\ref{sec:ordinary}; and
Section~\ref{sec:parallel-ov} proves the periodic repair.  The landing geometry and fixed-word transfer are developed in
Section~\ref{sec:two-landing}; Sections~\ref{app:positive-network}--
\ref{app:coarea} contain the positive-network and finite-fibre proof modules.
Their operator output is inserted in Section~\ref{sec:p38j}, followed by the
epoch restart in Section~\ref{sec:global}.  Section~\ref{sec:scope} records the final claim
boundaries.  Online Resource~1 contains expanded routine coordinate verifications,
source-version concordance, alternate packet charts, and supplementary
compatibility tables.
\section{Main Technical Theorem and Detailed Proof Strategy}
\label{sec:overview}

This section follows the mathematical dependency chain.  It explains the
geometric model, the periodic and global passages, and the exact point at
which each output
enters Theorem~\ref{thm:main}.

\subsection{The higher-dimensional obstruction}

The higher-dimensional issue is the power available from the long-bond
estimate.  Throughout the paper we write
\[
 L_\eps=|\log\eps|,
 \qquad L_\eps^\#=\lceil L_\eps\rceil\in\mathbb N,
 \qquad \cL_\eps=2+L_\eps.
\]
The integer envelope $L_\eps^\#$ is used only when a cutoff is a cardinality;
otherwise $L_\eps$ and $|\log\eps|$ retain their usual real-valued meaning.
We allow the absolute exponent on $\cL_\eps$ to change from line to line.
For a serial three-output long degree-$(3,3)$ component we prove the shell bound
\begin{equation}\label{eq:intro-sharp-shell}
 J_\rho\lesssim \cL_\eps^C
 \min\left\{\rho^2,
 \eps^{d-1}\mu^{-(d-1)}\rho^{-(d-3)}\right\}.
\end{equation}
At the critical output scale, the shell estimate falls short of the isolated
fixed-dimensional bound required by a componentwise continuation of the
low-dimensional proof.  This is the point at which the joint two-sublayer
estimate enters.

\begin{figure}[t]
\centering
\begin{tikzpicture}[
  box/.style={draw,rounded corners,align=center,minimum height=8mm,text width=2.45cm,font=\small},
  smallbox/.style={draw,rounded corners,align=center,minimum height=8mm,text width=2.35cm,font=\small},
  arr/.style={-{Latex[length=2mm]},thick},
  bad/.style={draw=gray!75,fill=gray!8},
  good/.style={draw=black,fill=blue!5}]
\node[box,bad] (c1) at (-5.0,1.0) {component 1\\separate estimate};
\node[box,bad] (c2) at (-5.0,-0.4) {component 2\\separate estimate};
\node[smallbox,bad] (prod) at (-1.85,0.3) {product of\\component bounds};
\draw[arr,gray!75] (c1.east) -- (prod.north west);
\draw[arr,gray!75] (c2.east) -- (prod.south west);
\node[align=center,text width=5.4cm,font=\small] at (-3.45,-1.65)
  {Low-dimensional route: the components are estimated separately.};

\node[box,good] (up) at (1.75,1.0) {connected upper\\time sublayer};
\node[box,good] (down) at (1.75,-0.4) {connected lower\\time sublayer};
\node[smallbox,good] (joint) at (5.0,0.3) {joint positive-kernel\\packet};
\draw[arr] (up.east) -- (joint.north west);
\draw[arr] (down.east) -- (joint.south west);
\node[align=center,text width=5.4cm,font=\small] at (3.4,-1.65)
  {Higher-dimensional route: two landing roots are recovered jointly.};
\end{tikzpicture}
\caption{Componentwise and joint proof architectures.  The right-hand route
keeps the complete nonnegative kernel intact and is the route used in the
main theorem.}
\label{fig:overview-componentwise-joint}
\end{figure}
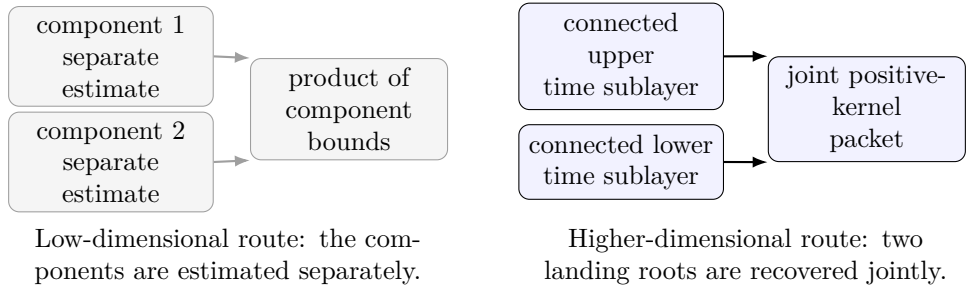

Figure~\ref{fig:overview-componentwise-joint} records the change in proof
architecture.  It is a dependency schematic whose boxes indicate the variables estimated
together.

\subsection{Main technical theorem}

The next statement is the canonical two-landing operator theorem.  Its packet
and output conventions are included here so that the estimate, its arbitrary
joint-kernel scope, its boundary-independent output domain, and its single
normalization are explicit before the detailed proof.  Section~\ref{sec:two-landing}
contains the proof; Sections~\ref{app:positive-network} and~\ref{app:coarea}
contain the two core closure modules used there.

Fix an integer packet cap $P_0$.  All constants in this section may depend on $P_0$;
the application in Section~\ref{sec:p38j} has
$|H|\le P_0=P_{\rm pkt}(d)$.
After replacing this cap by $2P_0$ once and relabelling it, the same symbol
$P_0$ bounds both the number of physical lines and the number of events,
since a packet of $P_0$ binary collision atoms meets at most $2P_0$ lines.
This simultaneous cap is the convention in every network and algebraic
multiplicity lemma below.  Fix one
half-open cell $\kappa$ for the time order, C-atom slot labels,
incoming signs, dyadic velocity sizes and torus lifts.  Write
$V_\eps=|\log\eps|^{C_0}$ for the velocity cutoff on this cell.  Every
function below is extended by zero outside the original support of the cell.
The cell is fixed during the calculation; its time and normal variables are
still integration variables.

\begin{definition}[Packet output operator]\label{def:packet-output}
For a packet $H$ and a cell $\kappa$, let
$\mathfrak B_{H,\kappa}$ be the product of the following live domains:
the one untranslated root position, all unpivoted transversal velocities,
all retained atom times and normals, and all external variables of the
complement.  Each position has a fixed fundamental-domain representative,
each time and normal remains in its original half-open cell, and each
velocity is integrated on $\R^d$ with the original zero extension.
If a consumer freezes a subcollection $y$ of the external variables, we
use the canonical product splitting
$\mathfrak B_{H,\kappa}=Y\times
\widehat{\mathfrak B}_{H,\kappa}$ and write $b$ for the remaining fibre
coordinate.  Under this splitting the output density remains in the branch measure, while
the conditional test function is the corresponding restriction of the
original joint kernel.

The recovery relation is obtained by the elastic, transport, tree-position
and collision-root substitutions constructed below.  The two
lower contacts are the first two lower atoms selected in
Lemma~\ref{lem:upper-tree-landings}.  The default sharp chart uses the
$2(d-1)$-velocity/two-time minor in \eqref{eq:maximal-minor-bound}, justified
for every packet cell by Theorem~\ref{thm:packet-fct}.  The maximal
$2d$-velocity disjoint chart of
Proposition~\ref{prop:disjoint-rooted-extension} and the full collision-free
new-line chord chart of
Proposition~\ref{prop:asymmetric-new-line-peeling} remain independent coarse
fallbacks.  The cell index $\kappa$ includes the corresponding finite
refinement.  For a nonnegative
Borel function $Q$,
Theorem~\ref{thm:appD-uniform-algebraic-multiplicity} proves that every
retained fibre has at most $M(d,P_0)$ regular solutions and supplies Borel
recovery branches $\Phi_{H,\kappa,j}$, $1\le j\le M(d,P_0)$.  An absent
branch is extended by zero.  Set
\begin{equation}\label{eq:Qsharp-definition}
 \int Q_{H,\kappa}^{\sharp}
 :=\sum_{j=1}^{M(d,P_0)}\int_{\mathfrak B_{H,\kappa}}
 Q(\Phi_{H,\kappa,j}(b))\,
 \mathbf 1_{\mathcal D_H}(\Phi_{H,\kappa,j}(b))
 W_{H,\kappa,j}(b)\dd b,
\end{equation}
where $W_{H,\kappa,j}$ consists only of the unspent collision weights and
the original cutoffs evaluated on branch $j$.  The domain
$\mathfrak B_{H,\kappa}$ is independent of the
values of the fixed boundary variables.  Dependence on those values occurs
only through the recovery branches.  The finite branch set is fixed by
$P_0$ and is counted as part of the output cell label.
\end{definition}

The symbol $J_{H,\kappa}$ in this section denotes the complete cell
operator, including integration in the source-boundary variables.  When a
consumer freezes such variables, its conditional meaning is the selected
Borel representative of
Definition~\ref{def:appD-selected-conditional-representative}: the original
conditional kernel is used off one fixed source-null set and is set to zero
on that set.  This measure-kernel convention is fixed before $Q$ is chosen, so the
exceptional set is independent of $Q$.
For a complete operator summed over the declared countable cell family, we
take the union of the cellwise null sets once before forming the sum; it is
still source-null and supplies one common selected representative.

\begin{theorem}[Two-landing-root operator]\label{thm:two-landing}
 Let $H=H_U\cup H_D$ be a connected full C-atom packet with $|H|\le P_0$
obtained by cutting it
as one free component.  Assume:
\begin{enumerate}[label=\textup{(\alph*)}]
\item $H_U$ and $H_D$ are connected in the atom graph;
\item every physical particle line of $H$ contains an atom in each
sublayer;
\item the two sublayers occupy disjoint globally ordered time intervals, so
$t_u<t_d$ for every $u\in H_U$ and $d\in H_D$.  Independently, on every
physical particle line all upper atoms precede all lower atoms; this
linewise interval property gives a unique internal bond crossing from
$H_U$ to $H_D$ on that line.  Neither assertion is used as an equivalent
reformulation of the other;
\item for $u\in H_U$ and $d\in H_D$,
$t_d-t_u\ge\epstar$ on the cell;
\item at least three physical particle lines meet $H$, and $H$ has no double
bond.
\end{enumerate}
Then there are two crossing bonds $e_1,e_2$, with distinct lower endpoints,
such that for every nonnegative joint test function $Q$ and every cell
$\kappa$,
\begin{equation}\label{eq:two-landing-cell}
 J_{H,\kappa}(Q)
 \le \cL_\eps^{C|H|}\eps^{-(d-1)}\eps^{2(d-1)}
       \epstar^{-2(d-1)}
 \int Q_{H,\kappa}^{\sharp}.
\end{equation}
Consequently, when $|H|\le P_{\rm pkt}(d)$,
\begin{align}
 J_H(Q)&\le \cL_\eps^C\eps^{d-1}\epstar^{-2(d-1)}
 \int Q_H^\sharp,\label{eq:two-landing-actual}\\
 \Exc(H)&\le\cL_\eps^C\eps^{2(d-1)}\epstar^{-2(d-1)}.
 \label{eq:two-landing-excess}
\end{align}
For a Source-B root set $R$, the normalization
\eqref{eq:root-normalization} gives
\begin{equation}\label{eq:two-landing-rooted}
 I_{H,R}(Q)\le
 \cL_\eps^C\eps^{(d-1)(|R|-1)}\eps^{2(d-1)}
       \epstar^{-2(d-1)}
 \int Q_H^\sharp.
\end{equation}
The output domain and its $M(d,P_0)$ branch slots in these inequalities are
independent of the concrete fixed boundary values.  Independently, the
full-chord charts prove every displayed estimate with the weaker loss
$\epstar^{-2d}$; this coarse route is retained as a fallback and a
coarse consistency check.
\end{theorem}

The every-$\kappa$ assertion above is an assertion for the complete
integrated operator.  Its fibrewise form holds for every fixed boundary
value only for the selected representative just specified; equivalently,
it holds for almost every boundary value for any original disintegration.
All fibrewise statements use this single selected boundary representative.

Section~\ref{app:coarea} gives the full regular-cell, submersion and
arbitrary-joint-kernel construction behind this theorem, including the
fibrewise preservation of the original time-order indicator.

The powers used by the theorem enter the main argument with the following
ownership.  The table records the factors in the displayed estimates and their roles in
the packet bound.
\begin{center}
\small
\begin{tabular}{@{}p{0.27\textwidth}p{0.25\textwidth}p{0.38\textwidth}@{}}
\toprule
Factor & Source & Role in the packet bound \\
\midrule
$\eps^{-(d-1)}$ & full connected component & unique baseline, counted once \\
$\eps^{2(d-1)}$ & two selected landing spheres & joint landing capacity \\
$\epstar^{-2(d-1)}$ & two retained long-separation blocks & sharp radial-conjugate loss \\
$\cL_\eps^{C|H|}$ & finite cells and recovery branches & bounded packet bookkeeping \\
$\eps^{(d-1)(k-1)}\epstar^{-k(d-1)}$ & fixed-$k$ birth flag & conditional multi-landing transfer \\
\bottomrule
\end{tabular}
\end{center}

\subsection{The two-sublayer mechanism and its consequences}

The main proof of Theorem~\ref{thm:main} is carried by Contributions I--III.
Contribution IV isolates the two-landing mechanism as a fixed-word module
proved under explicit packet hypotheses.  Section~\ref{sec:global} uses this
module through Corollary~\ref{cor:automatic-fixed-k-birth-production} within
that stated packet class.

\subsubsection{Two-landing packet geometry}

The analytic replacement acts jointly on two consecutive time sublayers.
Let $H=H_U\cup H_D$ be a connected collision packet such that both
sublayers are connected, every physical particle line meets both, no lower
atom is a parent of an upper atom, and the two time layers are separated by
\begin{equation}\label{eq:epstar-def}
 \epstar=\exp\big(-\sqrt{|\log\eps|}\big).
\end{equation}
After the preliminary double-bond case, at least three particle lines are
present.  The upper collisions contain a physical-particle spanning tree.
Every lower collision is consequently a co-tree contact.  More is true after
double bonds are excluded.  Take the first two lower collision atoms.  The
second must contain a particle line not present at the first; otherwise the
two consecutive atoms would be joined by both intervening particle bonds.
Choosing one of the two lines at the first atom then gives two crossing lines
with distinct last-upper endpoints.  Thus no unmarked lower contact occurs
before the selected second landing.

If the first two lower edges are disjoint, their incidence subspaces are
orthogonal, and chronological rooted-tree conormal elimination preserves the
resulting positive $2(d-1)$-dimensional tangential frame.  If they overlap,
write them as $\{a,b\}$ and $\{b,c\}$ in chronological order.  The new line
$c$ has no lower collision before the second landing.  From its last upper
contact $u_c$ to that landing $d_2$ it is therefore completely free.
After an exact zero-set row operation, retain the relative speed at $d_2$.
The radial-conjugate endpoint block has Jacobian
\[
 (t_{d_2}-t_{u_c})^{d-1}|g_2|
\]
and its literal physical Schur term is
$(t_{d_2}-t_{u_c})^{-1}A^*P_{g_2}A\ge0$.  The remaining lower constraint is
the one-speed frame at $d_1$ quotiented by its time direction.  The resulting
nonnegative rooted response is nonfocusing on its $(d-1)$-volume.
Thus, for every nonnegative joint test function and with no restriction on
the first lower time gap,
\begin{equation}\label{eq:intro-packet}
 J_H(Q)\le \cL_\eps^C\eps^{-(d-1)}\eps^{2(d-1)}
       \epstar^{-2(d-1)}
 \int Q^\sharp
\end{equation}
while allowing arbitrary upper co-tree contacts.  The factor
$\eps^{-(d-1)}$ is the single full-component normalization and is counted
once.  Formula \eqref{eq:intro-packet} acts on the joint positive kernel in a
fixed enlarged domain, with neither a boundary supremum nor a product of
marginal bounds.
Arbitrary upper incoming roots are included through their nonnegative
cylinder-curvature index terms, and all later lower roots are causally
neutral.  The full-chord $2d$-velocity construction is retained separately
and gives the weaker $\epstar^{-2d}$ loss.  The sharp chart closes
$\mathrm{FCT}(P_{\rm pkt}(d))$ with ordinary Cartesian recovery branches;
the radial row is used only in the local exterior-algebra proof.

\begin{figure}[t]
\centering
\begin{tikzpicture}[
  atom/.style={circle,draw,fill=white,minimum size=6mm,inner sep=0pt},
  line/.style={thick}, chord/.style={very thick,blue!65!black},
  every node/.style={font=\small}]
\node at (-3.7,2.0) {disjoint landings};
\draw[line] (-5.6,1.2) -- (-5.6,-1.25) node[below] {$a$};
\draw[line] (-4.2,1.2) -- (-4.2,-1.25) node[below] {$b$};
\draw[line] (-2.8,1.2) -- (-2.8,-1.25) node[below] {$c$};
\draw[line] (-1.4,1.2) -- (-1.4,-1.25) node[below] {$e$};
\node[atom] (dl1) at (-4.9,-0.25) {$d_1$};
\node[atom] (dl2) at (-2.1,-0.65) {$d_2$};
\draw[line] (-5.6,-0.25) -- (dl1) -- (-4.2,-0.25);
\draw[line] (-2.8,-0.65) -- (dl2) -- (-1.4,-0.65);
\draw[dashed] (-6.0,0.8) -- (-0.9,0.8) node[right] {$s$};

\node at (3.7,2.0) {overlapping landings};
\draw[line] (1.6,1.2) -- (1.6,-1.25) node[below] {$a$};
\draw[line] (3.4,1.2) -- (3.4,-1.25) node[below] {$b$};
\draw[chord] (5.2,0.8) -- (5.2,-1.25) node[below,text=black] {$c$};
\node[atom] (or1) at (2.5,-0.25) {$d_1$};
\node[atom] (or2) at (4.3,-0.65) {$d_2$};
\draw[line] (1.6,-0.25) -- (or1) -- (3.4,-0.25);
\draw[line] (3.4,-0.65) -- (or2) -- (5.2,-0.65);
\draw[dashed] (1.0,0.8) -- (5.7,0.8) node[right] {$s$};
\node[align=center,text=blue!65!black] at (5.5,0.02)
  {collision-free\\new-line chord};
\end{tikzpicture}
\caption{The exhaustive first-two-landing dichotomy.  The diagram records
particle-line incidence and chronological position only; the analytic frame
and its regular-cell hypotheses are those of
Theorem~\ref{thm:two-landing}.}
\label{fig:overview-landing-dichotomy}
\end{figure}
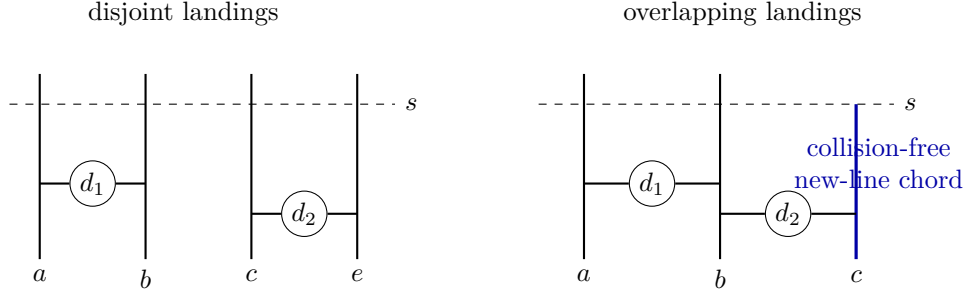

The first two lower landing edges are either disjoint or overlap in one
particle line; Figure~\ref{fig:overview-landing-dichotomy} displays exactly
this dichotomy.  In the second panel, the selected line $c$ first appears at
$d_2$, so its separator-to-$d_2$ segment is the collision-free chord used by
the radial-conjugate block.

\subsubsection{Periodic double-overlap repair}

The periodic obstruction also has a combinatorial part.  Source A already
observes that two torus overlap paths can form a double overlap segment and
lists the four
Source-B consumers that require modification.  We complete the arguments
left implicit there.  At the first failed no-double-overlap inference, the
two backward traces form a double support.  The legal cut gives a protected
degree-$(3,3)$ component or a singleton with two admissible contact roots.
At the higher endpoint the bottom pair gives the first admissible root.  At
the lower endpoint an overlap atom gives the lower root directly, while a
collision atom gives its earlier incoming quadratic companion.  The
one-incoming-root-per-lift lemma then separates the two lifts.  Permanent
original-incidence labels supply the certificate injections needed by the
cyclic, 2CONNUP and DOWN counts.

The joint estimate \eqref{eq:intro-packet} gives the modified alternative
$\mathrm{P3.8\mbox{-}J}$.  It supplies the second positive operator
interpoland required by the special version-2 top-layer consumer in
\cite{DHMlong}, with the entire complement inserted as $Q$ and with a larger
power margin.  This operator statement supplies exactly the form required by that consumer,
with the complete complement retained as one joint kernel.

\subsubsection{Long-time epoch restart}

The restart confines layer-dependent collision-word bookkeeping and
exceptional constants to individual blocks.

A block contains
$O(\sqrt{\log|\log\eps|})$ new layers, so its two layer-dependent constant
feedbacks retain the already verified subpower size.  At a block boundary,
the old nonempty cumulant is retained as one sealed joint kernel.  Contacts internal to
its influence block are resummed into the exact periodic hard-sphere group
before absolute values are taken; only paired C--O first contacts between
distinct blocks are exposed.  A common-radius, time-decreasing Gaussian closes
the remaining linear collision rate.  Every genuinely unmatched output root
is paid by a merge edge whose other leaves retain their source envelope, and
truncation defects retain their affected label set through the same sealed
interface.
A global cutoff tower of length $O(\log|\log\eps|)$ supplies every particle
order requested by the next block.  This separates the local exceptional
constant depth from the global cutoff depth and removes the former
$\exp(O(\mathfrak L^2))$ obstruction for the total number of layers.

\begin{figure}[t]
\centering
\begin{tikzpicture}[
  box/.style={draw,rounded corners,align=center,minimum height=9mm,text width=2.7cm},
  arr/.style={-{Latex[length=2mm]},thick},node distance=8mm]
\node[box] (epoch) {one-epoch operator\ and endpoint bounds};
\node[box,right=of epoch] (seal) {sealed complete\ joint kernel};
\node[box,right=of seal] (fresh) {fresh-reference\ conversion};
\node[box,right=of fresh] (next) {next epoch with\ paid defects};
\draw[arr] (epoch) -- (seal);
\draw[arr] (seal) -- (fresh);
\draw[arr] (fresh) -- (next);
\end{tikzpicture}
\caption{Epoch interface.  Internal contacts are first resummed inside the
sealed block; only the paired time-integrated interblock transfer is exposed
before the next epoch.}
\label{fig:overview-epoch-interface}
\end{figure}
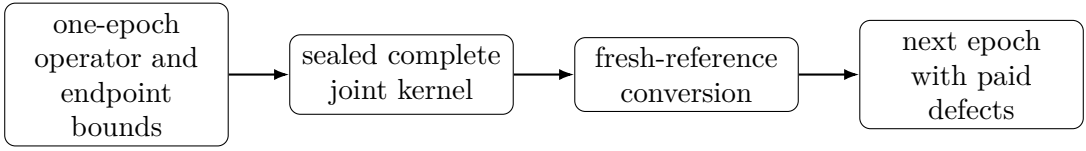

\subsubsection{Fixed-word multi-landing transfer}

Section~\ref{sec:transfer-principle} supplies a fixed-event-word theorem in
arbitrary dimension.  Within the regular
connected full two-sublayer packet class satisfying the hypotheses of
Theorem~\ref{thm:two-landing}, for every fixed $k$ any selected $k$-birth
flag, with one canonically chosen new line at each of its birth atoms, forms
simultaneous separator-to-landing ports.  A determinant-one
causal current clearing and a Moore--Penrose multiport Schur complement give
their joint $k(d-1)$-dimensional tangential frame.  The sphere normalization
then gives the exact capacity $\eps^{(d-1)(k-1)}$, and bounded-degree
finite-fibre coarea accepts an arbitrary joint kernel.  A packet with
$P\ge2k-1$ lines automatically has such a flag and gains
$\eps^{(d-1)(k-1)}\epstar^{-k(d-1)}$.  The case $k=2$ is a special case and
the first gain-producing one.  The long-time proof uses it on the smallest
exceptional packets and uses the largest admissible fixed \(k\) on larger
packets; the packet cap makes this a finite dimension-dependent choice.

\subsection{Proof architecture}

The proof is organized into six stages.  Each stage supplies the displayed
input for the stages that follow.
\begin{enumerate}[label=\textup{(\arabic*)}]
\item We construct the almost-everywhere periodic hard-sphere flow and expand
the initial grand-canonical correlations by an exact rooted-forest identity.
\item We establish the ordinary fixed-dimensional one-atom and serial
three-output estimates.  Their sharp shell bound also pinpoints why a direct
componentwise continuation of the low-dimensional argument is insufficient.
\item At the first periodic double-overlap failure, we assign a protected
certificate with permanent incidence labels.  This repairs the UP, cyclic,
2CONNUP and DOWN counts without charging the same configuration twice.
\item In a connected two-sublayer packet, a spanning tree in the upper
particle graph exposes two lower landing roots.  Disjoint roots are handled by
their joint incidence frame; overlapping roots are reduced by eliminating the
new particle line's collision-free chord.
\item Positive Jacobi-network coercivity prevents the rooted collision word
from focusing the landing frame.  Bounded-degree finite-fibre coarea then
turns the local submersion estimate into a bound for an arbitrary nonnegative
joint kernel, with every recovery branch included.
\item The joint packet is inserted directly into the exceptional top-layer
estimate.  Paid complete cumulants are then restarted across genuine blocks
by a paired stopped-contact partition identity; localized defects retain their
marked affected sets and cutoff remainders use one global tower.  The resulting
epoch ledger and final exponent bookkeeping yield Theorem~\ref{thm:main}.
\end{enumerate}

\begin{center}
\small
\begin{tabular}{@{}
>{\raggedright\arraybackslash}p{0.20\textwidth}
>{\raggedright\arraybackslash}p{0.31\textwidth}
>{\raggedright\arraybackslash}p{0.37\textwidth}
@{}}
\toprule
Stage & Verified output & Downstream use \\
\midrule
Imported architecture
& Section~\ref{sec:setup} and Online Resource 1, Section~\ref{app:source-interface}
& Version-locked structural reduction \\
Microscopic and ordinary inputs
& Sections~\ref{sec:ensemble-flow}--\ref{sec:ordinary}
& Flow, initial cumulants, and ordinary operator estimates \\
Periodic combinatorial repair
& Section~\ref{sec:parallel-ov}
& Periodic UP, cyclic, 2CONNUP, and DOWN inputs \\
Two-landing operator
& Theorems~\ref{thm:two-landing} and~\ref{thm:packet-fct},
with Sections~\ref{app:positive-network}--\ref{app:coarea}
& Arbitrary-kernel packet estimate \\
Top-layer insertion
& Proposition~\ref{prop:v2-consumer}
and Online Resource 1, Section~\ref{app:consumer-verification}
& Exceptional Source-B consumer \\
Epoch restart
& Proposition~\ref{prop:epoch-stacked-endpoint}
& Full-range decomposition and endpoint control \\
Final aggregation
& Section~\ref{sec:global}, Subsection~\ref{subsec:final-aggregation},
with arithmetic closure in Appendix~\ref{app:parameters}
& Theorem~\ref{thm:main} \\
\bottomrule
\end{tabular}
\end{center}

\subsection{The fixed-word transfer theorem and reusable mechanism}

Section~\ref{sec:transfer-principle} isolates the fixed-word chain of
landing-frame geometry, positive compliance and no focusing, and finite-fibre
arbitrary-kernel coarea developed in Theorem~\ref{thm:two-landing} and
Sections~\ref{app:positive-network}--\ref{app:coarea}.  Under its stated
packet-class hypotheses, this module is proved independently and yields
Corollary~\ref{cor:automatic-fixed-k-birth-production}, used in
Section~\ref{sec:global}.  Section~\ref{sec:scope} records its conditional,
fixed-word scope.

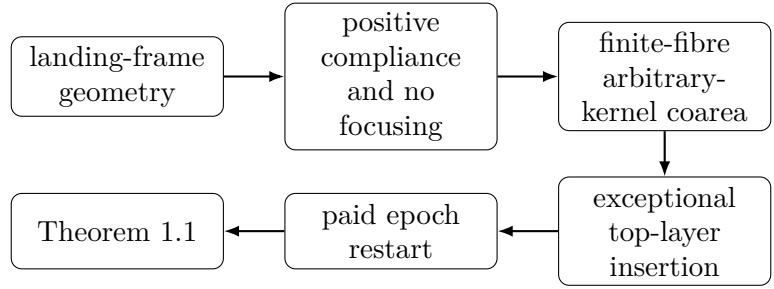
\begin{figure}[t]
\centering
\begin{tikzpicture}[
  box/.style={draw,rounded corners,align=center,minimum height=10mm,text width=2.55cm},
  arr/.style={-{Latex[length=2mm]},thick},
  node distance=6mm and 8mm]
\node[box] (frame) {landing-frame\\geometry};
\node[box,right=of frame] (positive) {positive compliance\\and no focusing};
\node[box,right=of positive] (coarea) {finite-fibre\\arbitrary-kernel coarea};
\node[box,below=of coarea] (consumer) {exceptional\\top-layer insertion};
\node[box,left=of consumer] (restart) {paid epoch\\restart};
\node[box,left=of restart] (main) {Theorem~\ref{thm:main}};
\draw[arr] (frame) -- (positive);
\draw[arr] (positive) -- (coarea);
\draw[arr] (coarea) -- (consumer);
\draw[arr] (consumer) -- (restart);
\draw[arr] (restart) -- (main);
\end{tikzpicture}
\caption{Load-bearing logical chain.  Periodic first-failure certificates
feed the exceptional insertion, while the microscopic, ordinary and imported
interfaces supply the surrounding reduction.}
\label{fig:overview-main-chain}
\end{figure}

Figure~\ref{fig:overview-main-chain} separates the local geometric output
from the two global passages needed for Theorem~\ref{thm:main}.
\section{Imported Reduction and Structural Inputs}
\label{sec:setup}

The microscopic objects and imported reduction are fixed here.  Together they
form the exact source interface for the geometric and hierarchy arguments that
follow.

Throughout the paper, \emph{Source B} denotes
arXiv:2408.07818v2 \cite{DHMlong}, and \emph{Source A} denotes
arXiv:2503.01800v1 \cite{DHMtorus}.  These are the exact imported proof
interfaces.  The mathematical interface remains Source-B v2.  Source-B v3 is used only as
a presentation benchmark, and Online Resource~1,
Section~\ref{app:source-interface}, records the equation-level concordance.

\subsection{Hard spheres and correlations}

For $N\ge0$ define the closed billiard table and its interior by
\[
 \overline{\cD}_N=\left\{z_N=(x_i,v_i)_{i=1}^N:
 d_{\T^d}(x_i,x_j)\ge\eps\text{ for }i\ne j\right\},
 \qquad
 \cD_N^\circ=\left\{d_{\T^d}(x_i,x_j)>\eps\text{ for }i\ne j\right\}.
\]
At a regular boundary point we identify the incoming and outgoing velocity
states by the reflection \eqref{eq:scattering-map}; equivalently, one may take
the postcollisional representative.  The boundary has phase measure zero, so
densities and correlation functions may be represented on $\cD_N^\circ$.
At a regular binary contact $x_i-x_j=\eps\omega$, the equal-mass elastic
reflection is
\begin{equation}\label{eq:scattering-map}
 v_i^+=v_i^- -((v_i^--v_j^-)\cdot\omega)\omega,
 \qquad
 v_j^+=v_j^- +((v_i^--v_j^-)\cdot\omega)\omega.
\end{equation}
It is an orthogonal reflection of $\R^{2d}$ and has absolute determinant one.

Set $a_\eps=\alpha\eps^{-(d-1)}$.  The grand partition function is
\begin{equation}\label{eq:partition}
 \cZ_\eps=\sum_{N\ge0}\frac{a_\eps^N}{N!}
 \int_{\overline{\cD}_N}\prod_{j=1}^N n_0(z_j)\dd z_N,
\end{equation}
and the initial density is
\begin{equation}\label{eq:grand-density}
 W_{0,N}(z_N)=\cZ_\eps^{-1}a_\eps^N
 \prod_{j=1}^Nn_0(z_j)\1_{\overline{\cD}_N}(z_N).
\end{equation}
Let $W_N(t)$ be the push-forward under the hard-sphere flow.  The rescaled
correlations are
\begin{equation}\label{eq:correlation-def}
 f_s(t,z_s)=(\alpha^{-1}\eps^{d-1})^s
 \sum_{n\ge0}\frac1{n!}\int W_{s+n}(t,z_{s+n})
 \dd z_{s+1}\cdots\dd z_{s+n}.
\end{equation}
The normalization cancels the activity of the $s$ displayed roots exactly.

In \eqref{eq:boltzmann-main},
\[
 v'=v-((v-v_1)\cdot\omega)\omega,
 \qquad
 v_1'=v_1+((v-v_1)\cdot\omega)\omega,
\]
and $n,n_1,n',n_1'$ are evaluated at $(x,v)$, $(x,v_1)$,
$(x,v')$, $(x,v_1')$, respectively.

\subsection{Collision molecules and associated operators}

We use the collision-history formalism of \cite[Definitions 2.1--2.6]{DHMtorus}
and the layered cumulant formalism of \cite[Sections 4--13, version 2]{DHMlong}.
The new argument needs only the structural facts listed below.

A collision atom has two top and two bottom slots.  Serial top/bottom slots
form a physical particle line.  A C-atom carries the reflection
\eqref{eq:scattering-map}; an overlap atom is velocity-transparent.  A
connected full component meeting $P$ physical particle lines has source
normalization $\eps^{-(d-1)P}$.  Each ordinary collision-tree contact pays one
position-sphere factor $\eps^{d-1}$.  Thus a spanning tree leaves exactly one
full-component baseline
\begin{equation}\label{eq:unique-baseline-setup}
 \eps^{-(d-1)P}\eps^{(d-1)(P-1)}=\eps^{-(d-1)}.
\end{equation}
For a connected full component $H$, its \emph{relative excess} is the
positive-operator norm after this unique baseline has been removed:
\begin{equation}\label{eq:relative-excess-definition}
 \Exc(H):=\inf\left\{A\ge0:
 J_H(Q)\le \eps^{-(d-1)}A\int Q_H^\sharp
 \text{ for every nonnegative Borel }Q\right\}.
\end{equation}
The infimum is taken after summing the fixed half-open cell decomposition.
Thus $\Exc(H)$ never includes a second full-component normalization, and an
estimate of the form $J_H\le\eps^{-(d-1)}G\int Q_H^\sharp$ says exactly
$\Exc(H)\le G$.

After times, slot labels, signs, dyadic velocity cells and torus lifts are
put into half-open cells, every molecule defines a positive associated
operator $J_H(Q)$.  Here $Q\ge0$ is the operator supplied by the unintegrated
complement.  The notation $Q^\sharp$ means that fixed boundary states are
substituted into $Q$ and all still-live output variables are integrated on a
fixed cutoff enlargement.  This enlargement may depend on the
global time, velocity and lift cutoffs, but not on the concrete boundary
values.  If $R$ is a set of root particle lines, the Source-B associated
integral satisfies the exact normalization identity
\begin{equation}\label{eq:root-normalization}
 I_{H,R}(Q)=\eps^{(d-1)|R|}J_H(Q).
\end{equation}

\subsection{The imported reduction}

The cluster expansion, time-layer decomposition and cutting identities are
separate from the quantitative diagonal.  The four logical levels below keep
that distinction explicit; in particular, every subpower factor is quantified.

\subsubsection{Replacement criterion}

\begin{proposition}[Equation-level replacement criterion]
\label{prop:source-exact-replacement}
For a fixed integer $d\ge4$, restrict to the high-recollision branch just described.  The new
two-landing construction may replace Source A, Proposition 3.8, in the proof
of \eqref{eq:source-trunc-output} only if all of the following assertions
hold.  Together with the unchanged high-particle route, these assertions
cover the full Source-B first-threshold selector.
\begin{enumerate}[label=\textup{(I\arabic*)}]
\item On every positive indicator, lift, slot, sign and dyadic cell, the
packet estimate is a fibrewise inequality for the complete $Q$ in
\eqref{eq:source-positive-Q}, with the original indicator
$\mathbf 1_{\mathcal D_M}$ and no supremum over fixed boundary values.
\item The packet is cut and protected by a cutting sequence admissible in
\eqref{eq:source-cutting-identity}; its complement has an ordered cleanup
for which the dependency property preceding
\eqref{eq:source-iterated-integral} remains true.
\item In the attached branch the packet supplies one and only one full
baseline, and the complement satisfies the count
\eqref{eq:source-attached-baseline}.  With this convention, the packet
version of \eqref{eq:source-attached-ledger} implies
\eqref{eq:source-top-second-interpoland}.
\item In the disjoint branch the top packet satisfies
\eqref{eq:source-disjoint-top}, including its $m$-time power and its
$r$-dependent gain, while the lower factor is exactly
\eqref{eq:source-disjoint-lower}.
\item The first interpoland \eqref{eq:source-top-first-interpoland} is
available by three source-level routes: the unchanged ordinary uncut
molecule estimate in the attached branch; the exact factorization
\eqref{eq:source-disjoint-factorization}, lower estimate
\eqref{eq:source-disjoint-lower}, and top rooted packet bound in the
disjoint nonempty branch; and the top rooted packet bound alone in the
disjoint empty branch.  All polylogarithmic factors and indicator counts are
absorbed in the displayed epsilon margins, and
Lemma~\ref{lem:source-two-regime-envelope} yields
\eqref{eq:source-top-target}.  Substitution into
\eqref{eq:source-top-error-reduction} then gives
\eqref{eq:source-trunc-output}.
\end{enumerate}
\end{proposition}

\begin{proof}
Items (I1)--(I2) are precisely the hypotheses needed to insert a new local
operator into the source cutting identity and its iterated positive
integration on the high-recollision branch.  In the attached case, (I3)
reproduces the only structural step in Source B that creates a strict positive
baseline power.  In the
disjoint case, (I4) reproduces the source factorization, including the empty
lower factor.  Under (I5), Lemma~\ref{lem:source-two-regime-envelope} converts
the two source estimates into \eqref{eq:source-top-target}; the error reduction
\eqref{eq:source-top-error-reduction} then closes exactly as in Source B.
The high-particle members of the supremum in
\eqref{eq:source-top-error-reduction} were already controlled by the
unchanged route preceding the proposition.  No other use of the exceptional
top combinatorial proposition occurs in that reduction.
\end{proof}

Proposition~\ref{prop:source-exact-replacement} states the replacement
criterion.  The global-coordinate, submersion, protected-cleanup and
fibrewise time-owner arguments below establish \textup{(I1)--(I5)}.

\begin{proposition}[Fixed structural Source-B reduction]
\label{prop:fixed-source-reduction}
Fix $d\ge2$, a positive integer $N_{\rm lay}$, a terminal time $T_*>0$, all
cutoffs and all starred constants.  Put $\tau_*=T_*/N_{\rm lay}$.  Suppose
that the initial cumulant, ordinary-component, cutting, weight, volume and
exceptional-top inputs satisfy the equation-level interfaces in
Online Resource~1, Section~\ref{app:source-interface}, with one layer-series parameter
$q_{\rm lay}<1$.  Then the version-2 expansion, at this fixed data, gives
the layer-endpoint identities \eqref{eq:source-fa}--
\eqref{eq:source-layer-error}, the terminal estimates
\eqref{eq:source-trunc-output} and \eqref{eq:source-iterated-error}, and the
exact terminal cumulant identity \eqref{eq:source-final-cumulant}.  This
statement concerns the selected endpoint $T_*$; it does not assert an
interpolation at interior times of that decomposition.
\end{proposition}

\begin{proof}
This is the structural content of Source B, Propositions~5.1--5.4 and its
Sections~8--13, in the source version specified in
Online Resource~1, Section~\ref{app:source-interface}.  The source correspondence table in
Online Resource~1, Section~\ref{appB:complete-call-matrix} records the variables,
normalizations and time owner at every call.  Positivity permits the full
unintegrated complement to remain one arbitrary kernel; no factorization of
that kernel is part of the proposition.  The exact initial-versus-restarted
error factors and the coefficient $2(d-1)$ are those in
\eqref{eq:source-iterated-error}.
\end{proof}

\begin{proposition}[Fixed-dimensional replacement principle]
\label{prop:reduction}
Use the normalizations, positive kernel \eqref{eq:source-positive-Q} and full
partial-order time domain \eqref{eq:source-time-domain}.  For the fixed
terminal data \((N_{\rm lay},T_*)\) of
Proposition~\ref{prop:fixed-source-reduction}, put
\[
 \tau_*:=T_*/N_{\rm lay},\qquad
 q_\beta:=C_\beta\max(1,\alpha)\max(1,A),\qquad
 \vartheta_*:=q_\beta\frac{\max(1,T_*)}{N_{\rm lay}},
\]
and assume
\begin{equation}\label{eq:fixed-owner-diagonal}
 0<\vartheta_*\le q_0<1,
 \qquad
 C_{\rm own}N_{\rm lay}\eps^{1/(72d)}\le1.
\end{equation}
For a fixed integer $d\ge4$, the hypotheses needed to apply
Proposition~\ref{prop:fixed-source-reduction} are exactly:
\begin{enumerate}[label=\textup{(R\arabic*)}]
\item the almost-everywhere periodic hard-sphere flow is invertible and
measure-preserving;
\item the initial grand-canonical cumulants obey
\eqref{eq:source-initial-cumulant} with logarithmic exponent
$C_{\rm init}^*(N_{\rm lay})$ from \eqref{eq:Cinit-definition}, and
$\|\Err_0\|_1\le\eps^{\Lambda_0/10}$;
\item all ordinary elementary kernels, the weight estimate and the
velocity-volume estimate hold fibrewise for arbitrary nonnegative $Q$ in
dimension $d$; cutoff, lift, dyadic and enlarged-domain losses are routed
through $C_{\rm ord}^*(N_{\rm lay})$ from \eqref{eq:Cord-definition}, and
the packet-free owner exhaustion of
Proposition~\ref{prop:packet-free-small-window-consumer} restores the full
$\vartheta_*^{\mathsf X}$ ordinary ledger at the fixed terminal window;
\item every bearing periodic failure of the Euclidean parallel-overlap
inference produces a distinct protected good component without changing the
cutting identity or its shadow map;
\item the exceptional high-recollision top layer has a rooted arbitrary-$Q$
packet of relative size at most
$\cL_\eps^{C_{\rm ord}^*}\eps^{2(d-1)}\epstar^{-2(d-1)}$, and all attached, disjoint,
cleanup and time-owner conditions (I1)--(I5) of
Proposition~\ref{prop:source-exact-replacement} hold.  Here a complete cell
operator is integrated in its source-boundary variables.  Whenever one of
the consumers is written fibrewise at a fixed boundary value, it uses the
single selected Borel conditional representative of
Definition~\ref{def:appD-selected-conditional-representative}; hence the
estimate holds at every value for that representative and agrees with the
original conditional operator almost everywhere.  For a complete countable
cell sum, the word single refers to the common zero extension obtained from
the union of the cellwise null sets, as fixed after
Definition~\ref{def:packet-output}.
\end{enumerate}
Under (R1)--(R5), all dimension-sensitive calls in the fixed structural
reduction are valid for the fixed integer $d\ge4$.
\end{proposition}

\begin{proof}
Parts 1--3 of \cite[Section~5]{DHMtorus} use (R1)--(R3).  Parts 4--5 use the
ordinary cutting algorithm, and their bearing periodic parallel-overlap
calls are precisely (R4).  In Part 6 the combinatorially unchanged
high-particle branch uses the ordinary UP estimate together with
Proposition~\ref{prop:packet-free-small-window-consumer}, both included in
(R3) with \(T=T_*\), whereas the
high-recollision branch uses (R5).  In the latter branch the fibrewise
inequality for the complete positive kernel permits the packet to be
integrated before its complement without a componentwise factorization.  The
attached and disjoint integer ledgers are exactly the conditions in
Proposition~\ref{prop:source-exact-replacement}.  No parameter diagonal is
used beyond the two explicit owner conditions
\eqref{eq:fixed-owner-diagonal}; their simultaneous proof for every
terminal rerun is deferred to Proposition~\ref{prop:sourceA-wrapper}.
\end{proof}

\begin{proposition}[Source-A growing-layer estimate]
\label{prop:sourceA-wrapper}
Fix an integer $d\ge4$ and $\beta>0$, and first fix the universal exponent
$\delta_0=1/20$.  Choose $\Gamma$, the fixed dominance exponents, $\kappa$
and $c_0$ in the order \eqref{eq:appE-choice-order}; in particular
$\delta_0$ is fixed before $\kappa$.  For
$0<\eps\le\eps_0(d,\beta)$ put $L_\eps=|\log\eps|$ and
\[
 N_{\rm lay}=\mathfrak L
 =\max\{1,\lfloor\kappa\sqrt{\log L_\eps}\rfloor\}.
\]
For every $(\alpha,A,t_{\rm fin})$ satisfying the one-epoch condition
\begin{equation}\label{eq:one-epoch-horizon}
 \max(1,\alpha)\max(1,A)\max(1,t_{\rm fin})
 \le c_0\sqrt{\log L_\eps},
\end{equation}
the Source-A Part-7 argument evaluates all constants before the terminal time
is chosen and gives
\begin{equation}\label{eq:reduction-wrapper-bounds}
 C_{\rm all}^*(\mathfrak L)\le L_\eps^{\delta_0},
 \qquad
 \vartheta_T:=
 \frac{C_\beta\max(1,\alpha)\max(1,A)\max(1,T)}{\mathfrak L}
 \le q_0<1,
 \qquad
 C_{\rm own}\mathfrak L\eps^{1/(72d)}\le1
\end{equation}
for every $T\in[0,t_{\rm fin}]$.  For $T>0$ the fixed structural reduction may
be run on $[0,T]$ with layer length $T/\mathfrak L$; the small-window
replacement \eqref{eq:source-weight-projection} is used in place of the source-only
formula \eqref{eq:source-weight-projection-original}, and
Proposition~\ref{prop:packet-free-small-window-consumer} supplies the
ordinary owner exhaustion needed to recover
$\vartheta_T^{\mathsf X}$.  All non-root
logarithmic losses in that run are bounded by
\begin{equation}\label{eq:structural-subpower-factor}
 \mathscr S_\eps:=L_\eps^{C\mathfrak L C_{\rm all}^*(\mathfrak L)}
 =\eps^{-o(1)},
\end{equation}
uniformly in $T$ and in the admissible triple.
\end{proposition}

\begin{proof}
The two genuine layer-dependent dominance conditions and every fixed-power
transition are established in Lemma~\ref{lem:appE-sourceA-Cj}.  The ordinary and
initial exponents are included by \eqref{eq:Call-definition} and
\eqref{eq:appE-Call-quadratic-growth}.  The parameter bound for
$\vartheta_T$ is monotone in $T$ and follows from the one-epoch hypothesis
after $c_0$ is chosen relative to $\kappa$.  The last, owner-diagonal
inequality in
\eqref{eq:reduction-wrapper-bounds} is
\eqref{eq:appE-owner-diagonal}; it follows from
\(\mathfrak L=O(\sqrt{\log|\log\eps|})\) and
\(\eta/\eta_*=\eps^{1/(72d)}\) and is independent of
\((\alpha,A,T)\).  Finally
\eqref{eq:appE-growing-subpower} proves
\eqref{eq:structural-subpower-factor}.  None of these thresholds depends on
a positive lower bound for $T/\mathfrak L$.  The endpoint $T=0$ is not
obtained by assigning a positive time power to a zero-dimensional
projection: every nonempty strict layer-time domain is empty, and the empty
molecule has $\mathsf X=0$.
\end{proof}

\begin{corollary}[Quantified diagonal reduction]
\label{cor:diagonal-reduction}
Assume that (R1)--(R5) of Proposition~\ref{prop:reduction} have been proved.
After the choices of Proposition~\ref{prop:sourceA-wrapper}, one
$\eps_0=\eps_0(d,\beta)>0$ works simultaneously for every admissible
$(\alpha,A,t_{\rm fin})$, every $T\in[0,t_{\rm fin}]$, every
$1\le s\le L_\eps$ and every nonempty $H\subset[s]$.  For the expansion
rerun with terminal time $T$,
\begin{align}
 \|E_H^T(T)\|_1
 &\le \mathscr S_\eps\eps^{1/(15d)}
 \left(\eps^{1/C_{14}^*}L_\eps^{C_{\rm all}^*}\right)^{|H|},
 \label{eq:reduction-cumulant}\\
 \|f_s^{{\rm err},T}(T)\|_1
 &\le\mathscr S_\eps\eps^{1/(20d)},
 \label{eq:reduction-truncation}\\
 \|f^{\mathcal A,T}(T)-n(T)\|_1
 &\le\mathscr S_\eps\eps^{7/(1000d)}.
 \label{eq:reduction-approx}
\end{align}
The associated iterated error is bounded by $\eps^{1/100}$.  At $T=0$ the
notation is fixed by
\begin{equation}\label{eq:T0-object-assignment}
 \widetilde f_s^0(0):=f_s(0),\qquad
 f_s^{{\rm err},0}(0):=0,\qquad
 \Err_{\mathfrak L}^0:=\Err_0,
 \qquad f^{\mathcal A,0}(0):=n_0,\qquad
 E_H^0(0):=E_H(0).
\end{equation}
Thus the initial forest remainder belongs exclusively to the iterated-error
term, and the endpoint construction requires no time-layer interpolation.
\end{corollary}

\begin{proof}
Apply Proposition~\ref{prop:fixed-source-reduction} to the dimension-$d$
interfaces and then Proposition~\ref{prop:sourceA-wrapper}.  The retained
endpoint exponents are $1/(15d)$ and $1/(20d)$, and the propagation
coefficient is $2(d-1)$.  The Source-A interpolation parameters are
$\theta=1/(200d)$ and $\theta_{\mathfrak L}=9/(1000d)$, whose averaged exponent is
$7/(1000d)$.  Per-root logarithmic losses are absorbed in the displayed
 $L_\eps^{C_{\rm all}^*}$, while all other occurrences are in
 $\mathscr S_\eps$.  It remains to check the iterated error here, without
 invoking the later assembly.  Put $a_0=-\log q_0>0$ and
 $Y_k=\Lambda_k^2A_k=\Lambda_k^{2+1/(10d)}$.  The hierarchy
 $\Lambda_k=\Lambda_{k+1}^{100d^2}$ gives, uniformly in
 $0\le\ell<\mathfrak L$,
 \[
  \sum_{k=\ell+1}^{\mathfrak L}Y_k\le2Y_{\ell+1}.
 \]
 Hence, for $1\le\ell<\mathfrak L$, the $\ell$th restarted term in
 \eqref{eq:source-iterated-error} is at most
 \[
  \exp\!\left(-\frac{a_0}{10}\Lambda_\ell
       +4(d-1)L_\eps\Lambda_{\ell+1}^{2+1/(10d)}\right)
  \le \exp\!\left(-\frac{a_0}{20}\Lambda_\ell\right)
  \le \eps^{101/100},
 \]
 after shrinking the same uniform $\eps_0$.  The last restarted term is
 $q_0^{\Lambda_{\mathfrak L}/10}\le\eps^{101/100}$, while the initial
 term is at most
 \[
  \exp\!\left[-L_\eps\left(\frac{\Lambda_0}{10}
            -4(d-1)\Lambda_1^{2+1/(10d)}\right)\right]
  \le\eps^{\Lambda_0/20}\le\eps^{101/100}.
 \]
 Since $\mathfrak L+1=O(\sqrt{\log L_\eps})$, their sum is at most
 $\eps^{1/100}$.  The uniform $T>0$ quantifier is already part of
Proposition~\ref{prop:sourceA-wrapper}: it applies the fixed structural
endpoint theorem separately on each interval $[0,T]$.  At $T=0$, use the
object assignment \eqref{eq:T0-object-assignment}.  The exact initial forest
expansion has the form \eqref{eq:source-final-cumulant} with
$\Err_{\mathfrak L}^0=\Err_0$; Corollary
\ref{cor:initial-cumulant-L1} gives \eqref{eq:reduction-cumulant},
\eqref{eq:source-initial-error} gives the iterated-error bound,
$f_s^{{\rm err},0}=0$ gives \eqref{eq:reduction-truncation}, and
$f^{\mathcal A,0}(0)=n_0=n(0)$ gives \eqref{eq:reduction-approx} with zero left
side.  Hence all four endpoint statements hold without double counting.
\end{proof}

The remaining sections prove (R1)--(R5), including every operator-level
condition in Proposition~\ref{prop:source-exact-replacement}.  Section~\ref{sec:global} and Appendix~\ref{app:parameters} establish the growing-layer estimate and perform the final
aggregation.
\section{Microscopic Flow and Initial Correlations}
\label{sec:ensemble-flow}

Every collision-history expansion below rests on two inputs: the
almost-everywhere microscopic group and the exact initial forest expansion.
They are established in this section.  The fixed-particle-number,
grand-canonical and reduced-Palm measures remain distinct, and the initial
cumulant proof includes the factorial resummation that cancels the grand
partition function.

\subsection{The fixed-$N$ billiard table and its singular set}

Write
\[
 Q_N=\left\{x_N\in(\T^d)^N:
 d_{\T^d}(x_i,x_j)\geq\eps\ \text{for }i\ne j\right\}.
\]
At a boundary point at which exactly one pair $i,j$ is in contact, the
condition $\eps<1/2$ selects a unique local lift
$x_i-x_j=m+\eps\omega$, $m\in\Z^d$, $\omega\in\Sph^{d-1}$.  The corresponding
regular incoming and outgoing states are identified by
\eqref{eq:scattering-map}.  Denote the resulting phase-space quotient by
$\mathsf X_N$.  Its Liouville measure is the measure induced by
$\dd x_N\dd v_N$; the boundary itself has zero Liouville measure.

\begin{definition}[Primary singular set]\label{def:primary-singular}
The primary singular set $\mathsf S_N^0$ consists of the states at which
at least one of the following occurs:
\begin{enumerate}[label=\textup{(S\arabic*)}]
\item two distinct boundary cylinders meet, so two or more contacts occur
at the same time;
\item a unique pair is in contact but
$(v_i-v_j)\cdot\omega=0$;
\item the same contact admits more than one torus lift.
\end{enumerate}
The last alternative is empty for $\eps<1/2$, but is retained to specify
the quotient without an implicit lift convention.
\end{definition}

On every bounded kinetic-energy set, $\mathsf S_N^0$ is a finite union of
smooth pieces of positive codimension.  The multiple-contact pieces have
codimension at least two in configuration space; the grazing condition is
one additional scalar equation on a regular boundary face.

\begin{lemma}[Finite regular histories]\label{lem:finite-history-cells}
Fix $N$, an energy bound $K$, a time $T$, and an integer $q$.  The set of
states of energy at most $K$ whose orbit has a specified regular history
of at most $q$ collisions on $[-T,T]$ is a countable union of relatively
open smooth cells.  On each cell the collision times, contact normals and
the map from initial to final state are smooth.  The subset whose history
ends at $\mathsf S_N^0$ has Liouville measure zero.
\end{lemma}

\begin{proof}
A history cell fixes the ordered colliding pairs and, whenever a free
flight crosses a fundamental-domain boundary, the finitely many torus
lifts allowed by $T\sqrt{2K}+1$.  Between contacts the motion is affine.
At a transversal contact the collision time is a simple zero of
\[
 |x_i-x_j+t(v_i-v_j)-m|^2-\eps^2,
\]
so the implicit-function theorem gives smooth dependence on the incoming
state.  Reflection is the smooth orthogonal map
\eqref{eq:scattering-map}.  Induction over the history proves the first
claim.  Pulling a positive-codimension primary singular piece back by the
smooth nonsingular history map preserves nullity.  There are only
countably many pair words and lift words.
\end{proof}

The preceding lemma does not by itself exclude infinitely many regular
collisions in finite time.  We isolate that input instead of hiding it in
the construction of the flow.

\begin{lemma}[No finite collision accumulation]
\label{lem:no-finite-accumulation}
Fix $N$ and a kinetic-energy bound.  A hard-sphere trajectory which avoids
multiple and grazing contacts has only finitely many contacts on every
compact time interval.
\end{lemma}

\begin{proof}
Suppose instead that regular collision times accumulate at a finite time
$t_*$.  Conservation of kinetic energy gives a common speed bound $V$.
On a terminal interval $I=(t_*-\delta,t_*)$, lift every continuous particle
path to $\R^d$ and choose $\delta>0$ so that
\begin{equation}\label{eq:terminal-common-lift-scale-d}
 2\eps+2V\delta<1.
\end{equation}
At a contact time $s\in I$ of particles $i,j$, there is a unique
$m_{ij}(s)\in\Z^d$ such that
\[
 |\widetilde x_i(s)-\widetilde x_j(s)-m_{ij}(s)|=\eps.
\]
If the same pair collides at $s,t\in I$, then
\[
 |m_{ij}(s)-m_{ij}(t)|\le2\eps+2V|s-t|<1.
\]
The left side is the norm of an integral vector, hence it is zero.  Thus
each pair which collides in $I$ has one fixed lattice label $m_{ij}$.

In $\R^{dN}$ retain only the finitely many pair walls
\[
 \mathcal B_{ij}:=\{X:|X_i-X_j-m_{ij}|\le\eps\}.
\]
They are globally Euclidean convex cylinders with smooth boundary.  The
lifted trajectory stays in their complement and has the usual specular
reflection at every collision.  After constant-speed reparametrization this
is a semi-dispersing billiard trajectory on the complete manifold $\R^{dN}$
with finitely many geodesically convex walls.  The no-accumulation theorem
\cite[Theorem~1, pp.~92--93]{BFK2002} gives only finitely many collisions in
$I$, a contradiction.  This common-lift proof uses exactly $\eps<1/2$ and
does not assert global convexity of a collision cylinder in the product
torus.

The quantitative Euclidean form of the same estimate is used later only
after the explicit short-window lifting in
Lemma~\ref{lem:long-gap-sublayers}; no global torus lift is assumed there.
\end{proof}

Alexander proves the corresponding full-measure all-time group and
Liouville theorem for finite hard spheres in bounded Euclidean piecewise
smooth tables
\cite[Theorems II.B.1--II.B.2, pp.~19--20]{Alexander1975}.  We do not
silently apply that Euclidean theorem to the torus.  Instead, the next proof
combines the manifold no-accumulation theorem above with the explicit
countable history cells of Lemma~\ref{lem:finite-history-cells}.  The
standard toral phase-space convention agrees with
\cite[Sections~1--2]{SimanyiSzasz1999}.

\begin{proposition}[Almost-everywhere periodic hard-sphere group]
\label{prop:ae-flow}
For fixed $N$ and $0<\eps<1/2$, there is a Borel invariant set
$\mathsf X_N^{\rm reg}\subset\mathsf X_N$ whose complement
$\Sigma_N$ has Liouville measure zero and maps
\[
 \Phi_N^t:\mathsf X_N^{\rm reg}\longrightarrow
 \mathsf X_N^{\rm reg},\qquad t\in\R,
\]
such that:
\begin{enumerate}[label=\textup{(F\arabic*)}]
\item $\Phi_N^0={\rm Id}$ and
$\Phi_N^{t+u}=\Phi_N^t\Phi_N^u$;
\item every trajectory has finitely many collisions on compact time
intervals and obeys free flight and \eqref{eq:scattering-map};
\item $\Phi_N^t$ is invertible, with inverse $\Phi_N^{-t}$, and velocity
reversal conjugates $\Phi_N^t$ to $\Phi_N^{-t}$;
\item $(\Phi_N^t)_\#(\dd x_N\dd v_N)=\dd x_N\dd v_N$.
\end{enumerate}
The exceptional set may be chosen as the forward and backward saturation
of the finite-history preimages of $\mathsf S_N^0$.  No additional
accumulation set is needed because Lemma~\ref{lem:no-finite-accumulation}
excludes that alternative before a singular endpoint.
\end{proposition}

\begin{proof}
First restrict to the kinetic-energy ball
\[
 \mathcal E_K:=\left\{\sum_i|v_i|^2\le K\right\}
\]
and an integer time window \([-T,T]\).  A maximal
trajectory that fails to cross the entire window either has infinitely many
regular collisions accumulating in finite time or reaches a first multiple
or grazing contact.  The first alternative is excluded by
Lemma~\ref{lem:no-finite-accumulation}.  Before the singular endpoint the
history is therefore finite.

Fix its pair word and lift word.  On the corresponding smooth history cell,
adjoin the possible hitting time $u\in[-T,T]$ as one extra variable and let
$\Psi(z,u)$ be the phase point obtained after the prescribed regular word
and the final free flight.  For fixed $u$ the map in $z$ is nonsingular by
Lemma~\ref{lem:finite-history-cells}.  In one fixed torus-lift chart write
\[
 g_{ij}(\Psi)=|x_i-x_j-m_{ij}|^2-\eps^2.
\]
At an intersection of two distinct collision walls, the differentials of
the two corresponding functions $g_{ij},g_{k\ell}$ are independent: for
disjoint pairs this is immediate, and for pairs sharing one label the two
unshared particle-position blocks separate the differentials.  Thus
\[
 (z,u)\longmapsto
 (g_{ij}(\Psi(z,u)),g_{k\ell}(\Psi(z,u)))
\]
has rank two on every regular multiple-contact stratum.  Its zero set in the
extended $(z,u)$ chart has codimension two and hence dimension at most one
less than the initial-state space.  Its projection to $z$ is locally a
Lipschitz image of that lower-dimensional set and therefore has Liouville
measure zero.

For grazing use instead the two functions
\[
 g_{ij}(\Psi),\qquad
 \dot g_{ij}(\Psi)
 =2(x_i-x_j-m_{ij})\cdot(v_i-v_j).
\]
On $g_{ij}=\dot g_{ij}=0$, the first differential has a nonzero position
component and the second has a nonzero relative-velocity component
$2(x_i-x_j-m_{ij})$; they are independent.  The same extended-chart
codimension-two and projection argument proves nullity.  Equivalently, this
is the usual zero-flux statement for the grazing part of the collision
section, now with the existential hitting time accounted for explicitly.

There are countably many pair/lift words, singular pair choices and smooth
history charts.  On each such chart split the terminal variable into the
forward and backward directions and restrict the extended zero set above to
the exact first-singular-hitting graph in that direction.  This graph is
Borel: the finite-word collision times and endpoint map are smooth, while
the event-order, simple-root, lift, no-earlier-hit and terminal-stratum
conditions are Borel.  Its projection to the initial state is injective,
because the first singular time in one fixed direction is unique.  The
Lusin--Souslin theorem \cite[Theorem~15.1]{Kechris1995} therefore makes the
exact projected image Borel.  Apply this separately on every labelled chart
and stratum, and include $\mathsf S_N^0\cap\mathcal E_K$ for a hit at time
zero.  The preceding codimension argument proves that the resulting exact
finite-window set $\Sigma_{K,T}^{\rm ex}\subset\mathcal E_K$ is null.  The
torus quotient causes no additional uncountable choice: compactness gives a
finite spatial chart cover at each event and the lift word is in the
countable set $(\Z^d)^q$.

Define, with no further Borel enlargement,
\[
 \Sigma_N:=\bigcup_{K,T\in\N}\Sigma_{K,T}^{\rm ex},
 \qquad
 \mathsf X_N^{\rm reg}:=\mathsf X_N\setminus\Sigma_N.
\]
This is a Borel conull set.  Every state has finite kinetic energy and every
finite time lies in an integer window.  If a maximal regular trajectory from
$\mathsf X_N^{\rm reg}$ had a finite endpoint, then
Lemma~\ref{lem:no-finite-accumulation} would leave only finitely many earlier
contacts; a nonsingular endpoint could be continued by free flight or one
regular reflection.  The endpoint would therefore be an exact singular hit,
contrary to the definition of $\Sigma_N$.  Thus every resulting trajectory
is nonsingular for all real times and has a finite history on compact
intervals.  If a shift $\Phi_N^s z$ hit a singular state at time $u$, uniqueness
and the elastic involution would make $z$ hit the same state at time $s+u$.
Hence $\mathsf X_N^{\rm reg}$ is invariant under every real trajectory shift;
no rational-shift completion is needed.

Free flight has determinant one.  At a regular collision the velocity
reflection is orthogonal, and its incoming and outgoing normal fluxes
satisfy
\[
 |(v_i^--v_j^-)\cdot\omega|
 =|(v_i^+-v_j^+)\cdot\omega|.
\]
Thus the collision cross-section measure is preserved.  Suspending it by
the free-flight time proves Liouville preservation on every finite
history cell.  The countable pair/lift/history partition and monotone
exhaustion in \(K,T\) give (F4).  The reflection is an
involution after velocity reversal, which proves (F1)--(F3).
\end{proof}

\subsection{The grand-canonical flow}

Let
\[
 \mathsf X_\#=\bigsqcup_{N=0}^{\infty}\{N\}\times\mathsf X_N,
 \qquad
 \mathsf X_\#^{\rm reg}=
 \bigsqcup_{N=0}^{\infty}\{N\}\times\mathsf X_N^{\rm reg}.
\]
The natural reference measure on this disjoint union has sector measure
$N!^{-1}\dd x_N\dd v_N$.  The countable union
\[
 \Sigma_\#=\bigsqcup_{N\ge0}\{N\}\times\Sigma_N
\]
is therefore null.  Define
\begin{equation}\label{eq:grand-flow}
 \Phi_\#^t(N,z_N)=(N,\Phi_N^t z_N)
 \quad\text{on }\mathsf X_\#^{\rm reg}.
\end{equation}

\begin{corollary}[Grand-canonical Liouville transport]
\label{cor:grand-flow}
The initial law \eqref{eq:grand-density} assigns probability zero to
$\Sigma_\#$, and \eqref{eq:grand-flow} is a measurable reversible group
on a set of full initial probability.  For every $N$,
\begin{equation}\label{eq:liouville-pushforward}
 W_N(t,z_N)=W_{0,N}(\Phi_N^{-t}z_N)
 \quad\text{for a.e. }z_N.
\end{equation}
Consequently the correlations \eqref{eq:correlation-def} are the
correlations of the unique sectorwise push-forward law.
\end{corollary}

\begin{proof}
Every $W_{0,N}$ is absolutely continuous with respect to sector
Liouville measure.  Hence
\[
 W_{0,N}(\Sigma_N)=0.
\]
Summing the nonnegative
sector masses proves the first assertion.  Measurability follows from
the countable finite-history cell decomposition.  Proposition
\ref{prop:ae-flow}(F4) gives \eqref{eq:liouville-pushforward}; Tonelli's
theorem then gives the correlation statement.
\end{proof}

The initial inhomogeneous grand-canonical law need not be invariant.
Invariance belongs to the Liouville reference measure, and
\eqref{eq:liouville-pushforward} is the required transport identity.

\subsection{Cavity insertion, void probability and the mean number}

Define the unscaled one-point intensity by
\begin{equation}\label{eq:unscaled-intensity}
 \rho_1(z)=\sum_{m\ge0}\frac1{m!}
 \int W_{m+1}(0,z,z_m)\dd z_m,
 \qquad f_1(0,z)=a_\eps^{-1}\rho_1(z).
\end{equation}
For $z=(x,v)$ define the cavity partition function
\begin{equation}\label{eq:cavity-partition}
 \cZ_\eps[z]=\sum_{m\ge0}\frac{a_\eps^m}{m!}
 \int_{\overline{\cD}_m}\prod_{j=1}^m n_0(z_j)
 \prod_{j=1}^m
 \1_{\{d_{\T^d}(x,x_j)>\eps\}}\dd z_m.
\end{equation}
This is a number depending on the proposed inserted position $x$; it is
not the original partition function.  Let $\mathbb P_\eps$ denote the
original grand-canonical law \eqref{eq:grand-density}, and put
\[
 \mathcal B_\eps(x)=B_{\T^d}(x,\eps)\times\R^d.
\]

\begin{lemma}[Exact insertion and Palm identities]
\label{lem:cavity-palm}
The one-point intensity $\rho_1$ satisfies
\begin{equation}\label{eq:exact-insertion}
 \rho_1(z)=a_\eps n_0(z)\frac{\cZ_\eps[z]}{\cZ_\eps}.
\end{equation}
Moreover,
\begin{equation}\label{eq:void-ratio}
 \frac{\cZ_\eps[z]}{\cZ_\eps}
 =\mathbb P_\eps\{N(\mathcal B_\eps(x))=0\}.
\end{equation}
Conditioned in the reduced-Palm sense on a particle at $z$, the remaining
$m$ particles instead have density
\begin{equation}\label{eq:reduced-palm}
 \frac{1}{\cZ_\eps[z]}\frac{a_\eps^m}{m!}
 \1_{\overline{\cD}_m}(z_m)
 \prod_{j=1}^m n_0(z_j)
 \prod_{j=1}^m\1_{\{d_{\T^d}(x,x_j)>\eps\}}.
\end{equation}
Thus \eqref{eq:void-ratio} is an original-law void probability, whereas
\eqref{eq:reduced-palm} is the normalized reduced-Palm law.
\end{lemma}

\begin{proof}
In the definition of the one-point intensity, choose the displayed
particle among the $m+1$ labels.  The factor $m+1$ cancels
$(m+1)!$, its activity and density give $a_\eps n_0(z)$, and all
remaining exclusions are precisely those in \eqref{eq:cavity-partition}.
This proves \eqref{eq:exact-insertion}.  Expanding the original-law event
that no centre lies in $\mathcal B_\eps(x)$ gives the same numerator and
proves \eqref{eq:void-ratio}.  Dividing the insertion numerator by its own
normalization gives \eqref{eq:reduced-palm}.
\end{proof}

\begin{remark}[Role of the mean-number estimate]
\label{rem:mean-number-role}
Proposition~\ref{prop:mean-number} is not an additional hypothesis in
Proposition~\ref{prop:reduction}.  It verifies that the activity convention
$a_\eps=\alpha\eps^{-(d-1)}$ indeed gives
$\eps^{d-1}\mathbb E N=\alpha+O(\alpha^2A\eps)$ and distinguishes the original
law from its reduced-Palm insertion.  The bearing initial input (R2) is the
forest estimate proved below.
\end{remark}

\begin{proposition}[Mean particle number]\label{prop:mean-number}
Under \eqref{eq:weighted-assumptions},
\begin{equation}\label{eq:mean-number}
 0\le\alpha-\eps^{d-1}\mathbb E_\eps N
 \le C_{\beta,d}A\alpha^2\eps.
\end{equation}
\end{proposition}

\begin{proof}
Equations \eqref{eq:exact-insertion}--\eqref{eq:void-ratio} first give
$0\le\rho_1(z)\le a_\eps n_0(z)$ and hence
\begin{align}
 1-\frac{\cZ_\eps[z]}{\cZ_\eps}
 &=\mathbb P_\eps\{N(\mathcal B_\eps(x))\ge1\}\notag\\
 &\le\mathbb E_\eps N(\mathcal B_\eps(x))
 =\int_{\mathcal B_\eps(x)}\rho_1(z')\dd z'
 \le a_\eps\int_{\mathcal B_\eps(x)}n_0(z')\dd z'.
 \label{eq:cavity-union-bound}
\end{align}
Since $\mathbb E_\eps N=\int\rho_1$, insertion yields
\begin{align}
 0\le a_\eps-\mathbb E_\eps N
 &\le a_\eps^2
 \int n_0(z)\int_{\mathcal B_\eps(x)}n_0(z')\dd z'\dd z.
 \label{eq:mean-before-scale}
\end{align}
Let $\rho_0(x)=\int_{\R^d}n_0(x,v)\dd v$.  The weighted assumption gives
$\|\rho_0\|_\infty\le C_{\beta,d}A$, and
$|B_{\T^d}(x,\eps)|=|B_d(1)|\eps^d$ for $\eps<1/2$.  The last integral
in \eqref{eq:mean-before-scale} is at most
$C_{\beta,d}A\eps^d$.  Multiply by $\eps^{d-1}$ and use
$a_\eps=\alpha\eps^{-(d-1)}$.
\end{proof}

\subsection{The truncated rooted-forest identity}

Fix the initial-layer hierarchy from
\eqref{eq:source-hierarchy}.  Write $N_{\rm lay}$ for the integer denoted by
$L$ in that source display, so that it is not confused with the logarithmic
scale $|\log\eps|$.  For fixed $d\ge4$ the recursion has the closed form
\begin{equation}\label{eq:initial-hierarchy-closed}
 A_0=(L_\eps^\#)^{(100d^2)^{N_{\rm lay}}},
 \qquad
 \Lambda_0=(L_\eps^\#)^{10d(100d^2)^{N_{\rm lay}}}.
\end{equation}
Choose an absolute $C_{\rm init,0}$, large enough for all the finite
forest encodings below, and define
\begin{equation}\label{eq:Cinit-definition}
 C_{\rm init}^*(N_{\rm lay})
 :=C_{\rm init,0}\bigl(1+(100d^2)^{N_{\rm lay}}\bigr).
\end{equation}
Then $\Lambda_0=A_0^{10d}$, $s\le A_0$, and, after one absolute enlargement
of $C_{\rm init,0}$,
\begin{equation}\label{eq:initial-hierarchy-size}
 A_0+\Lambda_0\le\cL_\eps^{C_{\rm init}^*(N_{\rm lay})}.
\end{equation}
Thus the exponent used in this subsection is not a constant with a hidden
dependence on the number of layers.  On the diagonal
$N_{\rm lay}=\mathfrak L\le\kappa\sqrt{\log|\log\eps|}$, it satisfies, for
every fixed $\eta>0$ and all sufficiently small $\eps$,
\begin{equation}\label{eq:Cinit-subpower}
 C_{\rm init}^*(\mathfrak L)
 \le C\exp\!\bigl(\log(100d^2)\,\kappa
              \sqrt{\log|\log\eps|}\bigr)
 \le |\log\eps|^\eta.
\end{equation}

For particle labels $i,j$ put
\[
 b_{ij}=\1_{\{d_{\T^d}(x_i,x_j)\le\eps\}}.
\]
The following is the initial, one-particle-per-cluster specialization of
the Source-B inclusion--exclusion lemma.  We state its construction to fix
the exact combinatorial coefficients.  The Mayer expansion and tree-graph
organization have their classical antecedents in
\cite{Penrose1963,BrydgesFederbush1978}; the particular stopped rooted-forest
map used here is proved directly, rather than imported from either identity.

\begin{lemma}[Truncated rooted-forest identity]
\label{lem:truncated-forest-identity}
Fix $s,n$ and an integer $\Lambda\ge2$.  There are finite families
$\mathfrak F_{s,\Lambda,s+n}$ and
$\mathfrak F^{\rm err}_{s,\Lambda,s+n}$ of tuples
$\mathcal F=(F,E,[s],C)$ and functions
$\widetilde\chi_{\mathcal F}$ such that
\begin{equation}\label{eq:finite-forest-identity}
 \prod_{1\le i<j\le s+n}(1-b_{ij})
 =\sum_{\mathcal F\in\mathfrak F_{s,\Lambda,s+n}
 \cup\mathfrak F^{\rm err}_{s,\Lambda,s+n}}
 \widetilde\chi_{\mathcal F}.
\end{equation}
Here $(F,E,[s])$ is a rooted forest, $[s]\subset F\subset[s+n]$,
and $C$ is a set of cross pairs for which $(F,E\cup C)$ is still a
forest.  In the main family every rooted tree has fewer than $\Lambda$
vertices.  In the error family the ordered search stops when the first
rooted tree reaches exactly $\Lambda$ vertices; the earlier trees have
fewer than $\Lambda$ vertices and the later roots are singletons.

For a main tuple,
\begin{equation}\label{eq:forest-chi-factor}
 \widetilde\chi_{\mathcal F}
 =\chi_{\mathcal F}(z_F)
 \prod_{i<j:\ i,j\notin F}(1-b_{ij}),
 \qquad
 |\chi_{\mathcal F}|
 \le\prod_{\{i,j\}\in E\cup C}b_{ij}.
\end{equation}
Its coefficient is invariant under increasing relabelling of the selected
nonroots and factorizes over a disjoint union with isolated root trees.  For
an error tuple, only the following domination is asserted and used:
\begin{equation}\label{eq:error-forest-domination}
 |\widetilde\chi_{\mathcal F}|
 \le \prod_{\{i,j\}\in E\cup C}b_{ij}
 \prod_{i<j:\ i,j\notin F}(1-b_{ij}).
\end{equation}
In particular, we do not assert that the exact error coefficient is
independent of the unselected variables.
\end{lemma}

\begin{proof}
This is the singleton-cluster specialization of Source B, version 2,
Lemma 4.27.  We record its fibre calculation because the distinction
between main and error forests is needed below.  Expand the left side of
\eqref{eq:finite-forest-identity} as
\begin{equation}\label{eq:mayer-graph-expansion}
 \sum_{G\subset K_{s+n}}(-1)^{|E(G)|}
 \prod_{\{i,j\}\in E(G)}b_{ij}.
\end{equation}
Apply the deterministic ordered breadth-first map
$\pi_{s,\Lambda}$ of that lemma to $G$.  Starting with the roots
$1,\ldots,s$, it assigns a reachable nonroot to the least root that reaches
it after the other roots are deleted; within a root class it orders first
by graph distance, then recursively by the first predecessor, and finally
by label.  The first predecessors form $E$.  The lexicographically first
acyclic cross pairs between root classes form $C$.  If a root class first
reaches $\Lambda$ vertices, only its first $\Lambda$ ordered vertices are
kept and all later root classes are singletons.  This gives exactly one
main or error tuple.

For a fixed output $\mathcal F$, partition all unordered pairs into the
pivotal set ${\rm PV}=E\cup C$, the forbidden set ${\rm FB}$ whose presence
would change the preceding ordered search, and the allowed set ${\rm AL}$.
The exact fibre characterization in Source B, version 2, Lemma 4.27 is
\[
 \pi_{s,\Lambda}(G)=\mathcal F
 \quad\Longleftrightarrow\quad
 {\rm PV}\subset E(G),\qquad {\rm FB}\cap E(G)=\varnothing.
\]
Summing \eqref{eq:mayer-graph-expansion} over that fibre therefore gives
\begin{equation}\label{eq:forest-fibre-product}
 \widetilde\chi_{\mathcal F}
 =\prod_{\{i,j\}\in{\rm PV}}(-b_{ij})
  \prod_{\{i,j\}\in{\rm AL}}(1-b_{ij}).
\end{equation}
Every pair with both endpoints outside $F$ is allowed.  In a main fibre,
every other allowed pair has both endpoints in $F$; extracting the
outside--outside product in \eqref{eq:forest-fibre-product} proves
\eqref{eq:forest-chi-factor}, including independence from the total number
of unselected vertices.  The definitions of ${\rm PV},{\rm FB},{\rm AL}$
depend only on the ordered rooted trees, which proves relabelling invariance
and factorization.  In an error fibre there may be allowed pairs joining an
unselected vertex to the truncated $\Lambda$-vertex tree.  Dropping those
factors $1-b_{ij}\le1$ and all other selected allowed factors gives
\eqref{eq:error-forest-domination}.  Finally the fibres partition all
finite graphs, proving \eqref{eq:finite-forest-identity}.
\end{proof}

\subsection{Exact resummation of the initial correlations}

From \eqref{eq:grand-density} and \eqref{eq:correlation-def}, for almost every
$z_s$,
\begin{multline}\label{eq:initial-correlation-exact}
 f_s(0,z_s)=\prod_{j=1}^s n_0(z_j)
 \sum_{n\ge0}\frac{a_\eps^n}{\cZ_\eps n!}
 \int\prod_{j=s+1}^{s+n}n_0(z_j)\\
 \times\prod_{1\le i<j\le s+n}(1-b_{ij})
 \dd z_{s+1}\cdots\dd z_{s+n}.
\end{multline}
Fix a canonical main-forest tuple with $k=|F|-s$ selected nonroots.  Its
increasing-label equivalence class has exactly $\binom nk$ representatives.
Using Lemma~\ref{lem:truncated-forest-identity},
\begin{equation}\label{eq:factorial-cancellation}
 \frac{a_\eps^n}{n!}\binom nk
 =\frac{a_\eps^k}{k!}\frac{a_\eps^{n-k}}{(n-k)!}.
\end{equation}
The coefficient $\chi_{\mathcal F}$ depends only on the $k$ selected
variables.  The remaining $n-k$ variables carry exactly their mutual
hard-core indicator in \eqref{eq:forest-chi-factor}; summing and
integrating them gives $\cZ_\eps$, which cancels the denominator in
\eqref{eq:initial-correlation-exact}.  We have therefore proved the exact
main-term identity
\begin{multline}\label{eq:initial-after-resummation}
 f_s^{\rm main}(0,z_s)=\prod_{j=1}^s n_0(z_j)
 \sum_{\mathcal F\in\mathfrak F_{s,\Lambda}}
 \frac{a_\eps^{|F|-s}}{(|F|-s)!}\\
 \times\int\prod_{p\in F\setminus[s]}n_0(z_p)
 \chi_{\mathcal F}(z_F)\prod_{p\in F\setminus[s]}\dd z_p,
\end{multline}
where the canonical family has node set $F=[|F|]$ after increasing
relabeling.  Thus no ratio of partition functions remains in a main initial
cumulant coefficient.

For an error fibre take absolute values before summing the unselected
particles and use \eqref{eq:error-forest-domination}.  The residual
outside--outside hard-core product again resums to $\cZ_\eps$, whereas all
selected--outside allowed factors have been discarded.  Hence the error
part of \eqref{eq:finite-forest-identity} defines $\Err_0$ and obeys
\begin{multline}\label{eq:initial-error-resummation}
 |\Err_0(z_s)|\le\prod_{j=1}^s n_0(z_j)
 \sum_{\mathcal F\in\mathfrak F^{\rm err}_{s,\Lambda}}
 \frac{a_\eps^{|F|-s}}{(|F|-s)!}\\
 \times\int\prod_{p\in F\setminus[s]}n_0(z_p)
 \prod_{\{i,j\}\in E\cup C}b_{ij}
 \prod_{p\in F\setminus[s]}\dd z_p.
\end{multline}

For a main forest, call a root active if its rooted tree is nontrivial or
if its singleton tree is incident to a cross pair in $C$.  Let $H$ be the
active set.  The inactive roots are isolated singleton components;
factorization in Lemma~\ref{lem:truncated-forest-identity} extracts
$\prod_{j\notin H}n_0(z_j)$.  Applying this extraction to
$f_s^{\rm main}$ defines $E_H(0,z_H)$ and gives
\begin{equation}\label{eq:initial-expansion}
 f_s(0,z_s)=\sum_{H\subset[s]}
 \prod_{j\notin H}n_0(z_j)E_H(0,z_H)+\Err_0(z_s),
\end{equation}
with $E_\varnothing=1$ and with $\Err_0$ bounded by the error-forest
majorant \eqref{eq:initial-error-resummation}.

\subsection{The corrected close-root selector}

The next lemma supplies a point that is easy to state incorrectly.  If
$S$ denotes only the singleton-tree roots, a root in $S$ is close to the
root of another tree, but that second root need not itself lie in $S$.
The selector must therefore include both endpoints.

\begin{lemma}[Close-root selector]
\label{lem:corrected-close-selector}
Let $\mathcal F=(F,E,H,C)$ be an active main forest and assume every
rooted tree has fewer than $\Lambda$ vertices.  There is a set
$Z=Z(\mathcal F)\subset H$ such that
\begin{enumerate}[label=\textup{(C\arabic*)}]
\item every $r\in Z$ has a distinct $r'\in Z$ satisfying
$d_{\T^d}(x_r,x_{r'})\le\Lambda\eps$ on the support of
$\chi_{\mathcal F}$;
\item with $k=|F|-|H|$, one has $k\ge|H\setminus Z|$.
\end{enumerate}
\end{lemma}

\begin{proof}
Let $S\subset H$ be the roots whose rooted trees are singletons.  Every
$r\in S$ is active, hence some cross pair joins $r$ to a vertex
$j(r)$ in a different rooted tree $F(\pi(r))$.  Choose the first such
pair in the fixed forest order and set
\begin{equation}\label{eq:corrected-Z-definition}
 Z=S\cup\{\pi(r):r\in S\}.
\end{equation}
The cross pair gives $d(x_r,x_{j(r)})\le\eps$.  The unique tree path from
$j(r)$ to $\pi(r)$ has at most $|F(\pi(r))|-1\le\Lambda-2$ edges, each
of length at most $\eps$ on the support bound
\eqref{eq:forest-chi-factor}.  Therefore
\[
 d_{\T^d}(x_r,x_{\pi(r)})\le(\Lambda-1)\eps.
\]
Every element of $S$ is paired with its $\pi(r)$, and every element of
$Z\setminus S$ was included as $\pi(r)$ for at least one $r\in S$.
This proves (C1).  Since $S\subset Z$, every root in $H\setminus Z$ has
a nontrivial rooted tree.  The disjoint nonroot sets of these trees give
\[
 k=\sum_{r\in H}(|F(r)|-1)\ge|H\setminus Z|,
\]
which is (C2).
\end{proof}

\begin{remark}[The Source-B close-root choice]
\label{rem:source-close-repair}
In the cited Source-B v2 argument, the proof takes $Z=S$ and establishes only closeness to a
root in $H$.  The minimal obstruction has two roots $r_1,r_2$: let
$F(r_1)=\{r_1\}$, let $F(r_2)$ contain the edge $r_2-j$, and choose the
cross pair $\{r_1,j\}$.  Then $S=\{r_1\}$, so the printed internal-$S$
condition asks $r_1$ to have a distinct partner in a one-point set and is
false, although $d(x_{r_1},x_{r_2})\le2\eps$.  Thus the issue is not merely
terminological.
Definition \eqref{eq:corrected-Z-definition} repairs the proof while
retaining both the printed close-root conclusion and the defect inequality
(C2).  All later generalized links may therefore have both endpoints in
$Z$, exactly as required by \eqref{eq:source-positive-Q}.  Exact source-file
locations are recorded in Online Resource~1, Section~\ref{appB:source-location-notes}.
\end{remark}

\subsection{Forest counting and the initial cumulant bound}

The spatial marginal
$\rho_0(x)=\int n_0(x,v)\dd v$ satisfies
\begin{equation}\label{eq:rho0-bound}
 \|\rho_0\|_\infty\le C_{\beta,d}A.
\end{equation}
Consequently a nonroot joined to a fixed parent costs
\begin{equation}\label{eq:one-mayer-link}
 a_\eps\int n_0(z_q)b_{pq}\dd z_q
 \le C_{\beta,d}\alpha A\eps.
\end{equation}
By \eqref{eq:loglog-main}, after reducing $\eps_0$ we may arrange
\begin{equation}\label{eq:mayer-small}
 C_{\rm init}^*(N_{\rm lay})C_{\beta,d}\alpha A\eps
 \le\tfrac12\eps^{1/2}.
\end{equation}
For the diagonal $N_{\rm lay}=\mathfrak L$, this reduction of $\eps_0$ is
uniform in the theorem parameters: by \eqref{eq:Cinit-subpower} and
\eqref{eq:loglog-main}, the logarithm of the coefficient multiplying $\eps$
is $o(|\log\eps|)$.  The same observation applies to the error-forest
geometric series below.

\begin{lemma}[Uniform enumeration of active forests]
\label{lem:active-forest-count}
Fix $H$, $k$ and $Z$.  The number of canonical active main forests with
root set $H$, $k$ nonroots and selector
\eqref{eq:corrected-Z-definition} equal to $Z$ is at most
\begin{equation}\label{eq:forest-count}
 \cL_\eps^{C_{\rm init}^*(N_{\rm lay})|H|}C^{|H|+k}k!.
\end{equation}
The same bound, summed over $Z$, holds for the error forests of a fixed
size.
\end{lemma}

\begin{proof}
Choose the positive tree sizes $n_r=|F(r)|$.  Their number is at most
$2^{|H|+k}$.  Allocate the $k$ nonroot labels among the trees in at most
\[
 \frac{k!}{\prod_{r\in H}(n_r-1)!}
\]
ways.  Cayley's formula and Stirling's inequality give
\[
 \prod_{r\in H}n_r^{n_r-2}
 \le C^{|H|+k}\prod_{r\in H}(n_r-1)!.
\]
Finally $E\cup C$ is a forest, so $|C|\le|H|-1$.  Each cross pair has at
most $|F|^2\le(\Lambda_0|H|)^2$ choices.  The hierarchy bound
\eqref{eq:initial-hierarchy-size} absorbs all cross-pair choices and the
choice of $Z$ into
$\cL_\eps^{C_{\rm init}^*(N_{\rm lay})|H|}$.  Multiplying the displayed bounds
proves \eqref{eq:forest-count}.  The map from an error output
to these counting data is
\[
 \mathcal F\longmapsto
 \bigl(r_*,(n_r)_{r\in H},
       (F(r)\setminus\{r\})_{r\in H},(E_r)_{r\in H},C\bigr),
\]
where $r_*$ is the first root stopped by the ordered search,
$n_{r_*}=\Lambda$, every earlier size lies in
$\{1,\ldots,\Lambda-1\}$, and every later root has size one.  The sets
$F(r)\setminus\{r\}$ are the allocated nonroot labels, $E_r$ are the rooted
labelled trees, and $C$ is the lexicographically selected acyclic cross-pair
set.  These data reconstruct the selected error tuple uniquely; interactions
with unselected vertices occur only in the discarded allowed factors of
\eqref{eq:error-forest-domination} and create no additional combinatorial
choice.  Thus the error family uses the identical allocation, Cayley, and
cross-pair bounds, on the stated restricted size vectors.  Summing over the
at most $|H|$ choices of $r_*$ is absorbed by
$\cL_\eps^{C_{\rm init}^*(N_{\rm lay})|H|}$.
\end{proof}

\begin{proposition}[Initial forest cumulants]
\label{prop:initial-cumulant}
For $s\le A_0$, the exact expansion \eqref{eq:initial-expansion} satisfies
\begin{equation}\label{eq:initial-EH}
 |E_H(0,z_H)|\le
 \cL_\eps^{C_{\rm init}^*(N_{\rm lay})|H|}
 \sum_{Z\subset H}\eps^{|H\setminus Z|/2}
 \prod_{p\in H}n_0(z_p)\,
 \1_{\rm close}^{H,Z}(z_H),
\end{equation}
where $\1_{\rm close}^{H,Z}$ asserts that every $r\in Z$ has a distinct
$r'\in Z$ within $\Lambda_0\eps$.  Moreover,
\begin{equation}\label{eq:initial-error}
 \|\Err_0\|_1\le\eps^{\Lambda_0/10}.
\end{equation}
\end{proposition}

\begin{proof}
Fix an active main forest and integrate its nonroot variables in decreasing
tree distance.  The cross-pair indicators may be discarded, while every
tree edge integrates one previously unintegrated child.  From
\eqref{eq:one-mayer-link},
\begin{equation}\label{eq:forest-integral}
 \frac{a_\eps^k}{k!}
 \int\prod_{q\in F\setminus H}n_0(z_q)
 |\chi_{\mathcal F}|\prod_q\dd z_q
 \le\frac{(C_{\beta,d}\alpha A\eps)^k}{k!}.
\end{equation}
Assign to the forest the set $Z(\mathcal F)$ from
Lemma~\ref{lem:corrected-close-selector}.  It supplies the close indicator
and $k\ge|H\setminus Z|$.  The $k!$ in
Lemma~\ref{lem:active-forest-count} cancels the denominator in
\eqref{eq:forest-integral}.  Since
$k\le(\Lambda_0-1)|H|$, the geometric sum and
\eqref{eq:mayer-small} give
\begin{align*}
 &\sum_{k\ge|H\setminus Z|}
 \cL_\eps^{C_{\rm init}^*(N_{\rm lay})|H|}C^{|H|+k}
 (C_{\beta,d}\alpha A\eps)^k\\
 &\hspace{35mm}\le
 \cL_\eps^{C_{\rm init}^*(N_{\rm lay})|H|}
 \eps^{|H\setminus Z|/2}.
\end{align*}
Multiplication by the unintegrated root factor
$\prod_{p\in H}n_0(z_p)$ and summation over $Z$ proves
\eqref{eq:initial-EH}.

For an error forest the stopped tree has $\Lambda_0$ vertices, hence
\begin{equation}\label{eq:error-k-lower}
 k=|F|-s\ge\Lambda_0-s\ge\Lambda_0-A_0.
\end{equation}
Integrate also the root variables, whose $n_0$ factors have integral one,
and repeat the enumeration.  Since the stopped construction has at most
$s\Lambda_0$ selected nodes,
\[
 \|\Err_0\|_1
 \le\cL_\eps^{C_{\rm init}^*(N_{\rm lay})s}
 \sum_{k=\Lambda_0-s}^{s(\Lambda_0-1)}
 (C_{\rm init}^*(N_{\rm lay})C_{\beta,d}\alpha A\eps)^k.
\]
The hierarchy $\Lambda_0=A_0^{10d}$, $s\le A_0$, and
\eqref{eq:mayer-small} imply, after reducing $\eps_0$,
that the right side is at most $\eps^{\Lambda_0/10}$.  This proves
\eqref{eq:initial-error}.
\end{proof}

\begin{corollary}[Integrated initial cumulant endpoint]
\label{cor:initial-cumulant-L1}
On the diagonal hierarchy, for every nonempty $H\subset[s]$,
$h:=|H|$, and $s\le L_\eps$, the initial cumulant satisfies
\begin{equation}\label{eq:initial-EH-L1-raw}
 \|E_H(0)\|_{L^1((\T^d\times\R^d)^H)}
 \le \eps^{h/3}
\end{equation}
after decreasing the common $\eps_0(d,\beta)$.  In particular, because
$C_{14}^*>4$,
\begin{equation}\label{eq:initial-EH-L1-consumer}
 \|E_H(0)\|_1
 \le \mathscr S_\eps\eps^{1/(15d)}
 \bigl(\eps^{1/C_{14}^*}L_\eps^{C_{\rm all}^*}\bigr)^h .
\end{equation}
\end{corollary}

\begin{proof}
Fix $Z\subset H$, $z:=|Z|$, in
\eqref{eq:initial-EH}.  If $z=0$, integration of the root densities gives
one.  If $z>0$, form on $Z$ the graph that joins $p\ne q$ when
$d_{\T^d}(x_p,x_q)\le\Lambda_0\eps$.  The close indicator says that this
graph has minimum degree at least one.  Hence every connected component has
at least two vertices.  If it has $c$ components, a spanning forest has
$z-c\ge z/2$ edges.  Choosing in every nonempty graph the
lexicographically first spanning forest gives the pointwise union bound over
at most
\begin{equation}\label{eq:initial-close-forest-count}
 (z+1)^{z-1}\le z^z\le h^h
\end{equation}
labelled forests; here $z\ge2$, and the first inequality follows by adding
one auxiliary root and applying Cayley's formula.

For one such forest, integrate its leaves successively.  Velocities first
give the spatial densities $\rho_0$.  By \eqref{eq:rho0-bound}, every tree
edge, with the parent fixed, costs at most
\begin{equation}\label{eq:initial-close-edge-cost}
 \sup_x\int_{d_{\T^d}(x,y)\le\Lambda_0\eps}
       \rho_0(y)\,\dd y
 \le C_{\beta,d}A(\Lambda_0\eps)^d .
\end{equation}
The one remaining root of each component integrates to one.  After reducing
$\eps_0$, the quantity in \eqref{eq:initial-close-edge-cost} is at most
one, so
\begin{equation}\label{eq:initial-close-L1}
 \int\prod_{p\in H}n_0(z_p)\,
       \1_{\rm close}^{H,Z}(z_H)\,\dd z_H
 \le h^h
 \bigl(C_{\beta,d}A\Lambda_0^d\eps^d\bigr)^{z/2}.
\end{equation}
Insert this in \eqref{eq:initial-EH} and sum over $Z$.  Since

\[
 \frac{h-z}{2}+2z=\frac h2+\frac{3z}{2}\ge\frac h2,
\]
we obtain
\begin{equation}\label{eq:initial-close-summed}
 \|E_H(0)\|_1
 \le
 \left[
  2h\,L_\eps^{C_{\rm init}^*(\mathfrak L)}
 \max\{1,(C_{\beta,d}A\Lambda_0^d)^{1/2}\}
  \eps^{1/2}
 \right]^h .
\end{equation}
Here $h\le s\le L_\eps$.  The diagonal growth bounds for
$C_{\rm init}^*(\mathfrak L)$, $A$, and $\Lambda_0$, proved in
Appendix~\ref{app:parameters}, make the bracket
$\eps^{1/2-o(1)}$, uniformly in the theorem parameters.  It is at most
$\eps^{1/3}$ after decreasing the same $\eps_0$, which proves
\eqref{eq:initial-EH-L1-raw}.  Finally
$1/3>1/C_{14}^*+1/(15d)$; the subpower factors in
$\mathscr S_\eps L_\eps^{C_{\rm all}^*h}$ only enlarge the right side of
\eqref{eq:initial-EH-L1-consumer}.  This proves the consumer form.
\end{proof}

Proposition~\ref{prop:ae-flow} and Corollary~\ref{cor:grand-flow} prove
(R1) of Proposition~\ref{prop:reduction}.  Proposition
\ref{prop:initial-cumulant}, including the repaired internal close-root
selector, proves (R2).

\section{Ordinary Fixed-Dimensional Operator Estimates}
\label{sec:ordinary}

The ordinary Source-B cutting algorithm uses the local analytic inputs proved
in this section.  Every estimate is an inequality for a
nonnegative, possibly nonfactorized test function.  We work after fixing
the C/O labels, slot orientation, signs of time differences, dyadic
velocity cells and torus lifts.  The number of such cells is absorbed in
\(\cL_\eps^C\).  Here and in every display of this section the exponent is
routed as follows:
\begin{equation}\label{eq:section04-Cord-routing}
 C_{\rm loc}\le C_{\rm ord}^*(N_{\rm lay}),
 \qquad
 \cL_\eps^{C_{\rm loc}}
 \le \cL_\eps^{C_{\rm ord}^*(N_{\rm lay})},
\end{equation}
with \(C_{\rm ord}^*\) defined in \eqref{eq:Cord-definition}.  This ledger
includes velocity cutoffs, torus-lift counts, dyadic shells, normal and
rotation charts, slot/sign/time-order refinements, and enlarged output
boxes.  Every occurrence of \(C^*\) in the weight and projection argument
means the fixed source coefficient \(C^*(N_{\rm lay})\), which is also
dominated by \(C_{\rm ord}^*\).  When ordinary calls are repeated across
layers, their remaining repetition is charged to the global factor
\(\mathscr S_\eps\) of \eqref{eq:structural-subpower-factor}.  Thus an
apparently absolute \(C\) below never hides a coefficient growing with the
layer number or the ordering complexity.

The route from local collision coordinates to the hierarchy has three steps:
the one-atom maps, the possible three-output orientations, and the combination
of their velocity volume with the time weights.  In the last step the
distinguished time stays outside the projected simplex.  That separation is
used when the exceptional packet enters in Section~\ref{sec:p38j}.

\subsection{The one-atom coordinate maps}

Fix a regular atom with bottom states $z_1,z_2$, top states $z_3,z_4$ and
time $t$.  If $z_1$ is fixed, the independent variables for a degree-three
C-atom are
\begin{equation}\label{eq:degree-three-chart}
 (t,\omega,v_2)\in I\times\Sph^{d-1}\times\R^d,
\end{equation}
and all other states are recovered by
\begin{align}
 x_2&=x_1+t(v_1-v_2)-\eps\omega,\notag\\
 v_3&=v_1-[(v_1-v_2)\cdot\omega]\omega,
 &x_3&=x_1+t(v_1-v_3),\notag\\
 v_4&=v_2+[(v_1-v_2)\cdot\omega]\omega,
 &x_4&=x_2+t(v_2-v_4).
 \label{eq:degree-three-recovery}
\end{align}
For an O-atom, $v_3=v_1$, $v_4=v_2$, $x_3=x_1$, and $x_4=x_2$ in
transported coordinates.  The contact-sphere formula is
\begin{equation}\label{eq:contact-sphere-basic}
 \eps^{-(d-1)}\int_{\R^d}
 \delta(|x_1-x_2+t(v_1-v_2)|-\eps)F(x_2)\dd x_2
 =\int_{\Sph^{d-1}}F(x_1+t(v_1-v_2)-\eps\omega)\dd\omega.
\end{equation}
The scattering part of \eqref{eq:degree-three-recovery} is an orthogonal
map on $(v_1,v_2)\in\R^{2d}$.

It follows directly from \eqref{eq:source-C-kernel} and
\eqref{eq:contact-sphere-basic} that the three one-atom operators are:
\begin{enumerate}[label=\textup{(O\arabic*)}]
\item with two admissible fixed ends, the contact time and lift are
discrete and the operator is at most substitution of norm one;
\item with one fixed end,
\begin{equation}\label{eq:degree-three-exact}
 \mathcal J_3(Q)=
 \int_I\int_{\Sph^{d-1}}\int_{\R^d}
 [(v_1-v_2)\cdot\omega]_-Q^\flat
 \dd v_2\dd\omega\dd t;
\end{equation}
\item with no fixed end, one also integrates $(x_1,v_1)$ and multiplies
by the single full-component factor $\eps^{-(d-1)}$.
\end{enumerate}
Here $Q^\flat$ means the original $Q$ composed with
\eqref{eq:degree-three-recovery}; no marginalization is performed.
After the Maxwellian and position cutoffs, the containing domains in
\eqref{eq:degree-three-chart} and in (O3) have volume $\cL_\eps^C$ and
are independent of the concrete fixed state.

The ordinary small-time, small-relative-velocity and small-position
alternatives follow in the same chart, but the two kinds of smallness have
different radii.  Fix $0<\upsilon<1/2$.  A restriction
$|t-t_*|\le\eps^\upsilon$ gives a $t$-section of length
$2\eps^\upsilon$.  A permitted nonserial restriction
$|v_i-v_j|\le\eps^\upsilon$ reduces, after resolving the four possible
outputs of the scattering map, either to
$|u\cdot\omega|\le\eps^\upsilon$ or to
$|u-(u\cdot\omega)\omega|\le\eps^\upsilon$; for fixed $\omega$ at least
one coordinate of $u=v_2-v_1$ is then confined to an interval of length
$C\eps^\upsilon$.

For a permitted transported position restriction the radius is instead
$\eps^{1-\upsilon}$.  Since all positions at one atom agree after transport
up to $O(\eps)$, it implies, for a fixed external $x_*$ or a permitted
nonserial pair,
\[
 |(x_1-x_*)+t(v_1-v_j)|\le 2\eps^{1-\upsilon}.
\]
On the dyadic cell $|t|\sim D$, integration in one coordinate of
$v_1-v_j$ costs $C\eps^{1-\upsilon}/D$, whereas the length of the same
time cell is $O(D)$.  Their product is $O(\eps^{1-\upsilon})$ on every
nonterminal cell.  The terminal cell
$|t|\le C\eps^{1-\upsilon}$ costs the same by its time length.  There are
only $O(|\log\eps|)$ cells, and
$\eps^{1-\upsilon}|\log\eps|^C\le
\eps^\upsilon|\log\eps|^{C'}$ because $\upsilon<1/2$.
Consequently, if $G$ denotes any one of these three good cells, the exact
arbitrary-kernel statement is
\begin{equation}\label{eq:one-atom-good-section}
 \mathcal J_3(\mathbf 1_GQ)
 =\int_{\Omega'}\mathbf 1_G(w)Q^\flat(w)\,\dd\mu_3(w),
 \qquad
 \mu_3(G)\le \cL_\eps^C\eps^\upsilon .
\end{equation}
Here $\Omega'$ is fixed by the global cutoffs and $\dd\mu_3$ includes the
collision flux and the chart density.  The first identity, rather than the
false inequality
$\int_GQ\le |G|\int_{\Omega'}Q$, is the positive-operator statement valid
for arbitrary $Q$.  The small scalar may be extracted only after the
complete source weight and Maxwellian estimates have supplied a pointwise
conditional majorant: if $Q^\flat\le M$ on the current conditional cell,
then
\[
 \mathcal J_3(\mathbf 1_GQ)
 \le \cL_\eps^C\eps^\upsilon M .
\]
For a degree-four atom, use the same pointwise recovery identity on its full
vector output and the same zero extension; its unique
$\eps^{-(d-1)}$ baseline is retained outside the excess.  The slot exclusions for O-atoms are exactly those in
\eqref{eq:source-elementary-operator}; a serial equality is never charged
as a small-variable gain.

\subsection{Two periodic contacts of one transported pair}

Let a relative state $(x,g)$ have contacts at $t_1,t_2$ with lifts
$m_1,m_2$:
\begin{equation}\label{eq:two-lifted-contacts}
 x+t_i g-m_i=\eps\omega_i,\qquad i=1,2.
\end{equation}
If both roots are incoming and distinct, Lemma~\ref{lem:one-incoming-root}
gives $m=m_1-m_2\ne0$.  Subtracting
\eqref{eq:two-lifted-contacts} and wedging with $g$ gives
\begin{equation}\label{eq:lattice-tube-wedge}
 |m\wedge g|\le 2\eps|g|.
\end{equation}
Since $|m|\ge1$ and $|g|\le\cL_\eps^C$, $g$ belongs to a tube about
$\R m$ with $d-1$ transverse radii at most $\cL_\eps^C\eps$.  For every
fixed admissible $m$ its $d$-dimensional volume is at most
\begin{equation}\label{eq:multiple-contact-tube}
 \cL_\eps^{C_d}\eps^{d-1}.
\end{equation}
The bounded time and velocity cells allow only $\cL_\eps^C$ choices of
$m$, so the same bound holds after summing lifts.

The same conclusion holds for two outgoing roots after reversing time and,
more generally, for any two contact times once their lifts are known to be
different; the wedge calculation itself does not use the contact sign.
At a C-atom reflection is orthogonal in the relative velocity and preserves
the tube volume.  Thus it applies equally to the incoming pair, the reflected
outgoing pair, and the temporally ordered outgoing--incoming pair carried by
a periodic return.  These are precisely the support conditions created by a
double bond or double overlap segment.

If both the incoming and reflected outgoing relative states have two
contacts, apply \eqref{eq:lattice-tube-wedge} before and after reflection.
Put $\delta=\cL_\eps^C\eps$ and, on a fixed first-lift tube, write
$g=ae+z$ with $|z|\le\delta$, $e\in\Sph^{d-1}$, and $|a|\le\cL_\eps^C$.
The part $|a|\le2\delta$ already gains an additional axial length
$O(\delta)$.  On $|a|>2\delta$, division by $|a|$ shows that the second
tube requires
\[
 \operatorname{dist}\bigl(R_\omega e,\mathbb R f\bigr)
 \le C\delta/|a|,
 \qquad R_\omega e=e-2(e\cdot\omega)\omega,
\]
for the unit direction $f$ of the second nonzero lattice vector.  Set
 $q=C\delta/|a|$ and
 $s=\operatorname{dist}(\mathbb Re,\mathbb Rf)$.  Away from the equatorial
 critical fibre, the elementary reflection map
 $F_e(\omega)=\mathbb R(R_\omega e)$ has at most four regular local inverse
 branches.  On the nondegenerate branches its singular values are bounded below by an
absolute constant.  At either branch approaching the equator
$e\cdot\omega=0$, one singular value remains comparable to one and the
other $d-2$ are comparable to $s$.  Hence, when $s>Cq$, a direct coordinate
box estimate gives preimage measure
$C_dq(q/s)^{d-2}\le C_dq$.  When $s\le Cq$, the same preimage is contained
in the union of the equatorial band $|e\cdot\omega|\le Cq$ and finitely many
$Cq$-caps about the nondegenerate inverse images; their measures are
$O_d(q)$ and $O_d(q^{d-1})$, respectively.  Thus the area formula on these
explicit branches gives the uniform coarse bound
\[
 \sigma_{\Sph^{d-1}}\{\omega:
   \operatorname{dist}(R_\omega e,\mathbb Rf)\le C\delta/|a|\}
 \le C_d\delta/|a|.
\]
Integrating $a$ from $2\delta$ to the velocity cutoff costs only
$C\delta|\log\delta|\le\cL_\eps^C\eps$.  Together with the small-$a$ part
and \eqref{eq:multiple-contact-tube}, this proves the double capacity
\begin{equation}\label{eq:double-periodic-capacity}
 \cL_\eps^{C_d}\eps^d.
\end{equation}
Tangencies are zero-flux boundaries of the half-open sign cells.

\subsection{The serial $j=3$ output: independent coordinates}

The only ordinary two-atom estimate not supplied by the preceding
$d$-dimensional density concerns a type-$\{33{\rm A}\}$ pair.  Cut the lower
atom as a two-fixed-end collision root and keep the upper atom in the
coordinates \eqref{eq:degree-three-chart}.  Let
\[
 u=v_2-v_1,\qquad
 y=(u\cdot\omega)\omega,\qquad w=u-y\in y^\perp.
\]
The selected output is $y$; the two fixed exterior states are denoted by
$z_1,z_7$.  On a fixed second-contact lift $m$, transport and contact give
\begin{equation}\label{eq:j3-second-contact}
 b_m(t_2)-(t_1-t_2)y=\eps\omega_2,\qquad
 b_m(t)=x_1-x_7+t(v_1-v_7)-m.
\end{equation}
Put $\Delta=t_1-t_2$ and assume the ordinary long cell
$|\Delta|\ge\mu$.

\begin{lemma}[Carleman chart in dimension $d$]
\label{lem:d4-carleman}
On either incoming sign chart and away from $y=0$,
\begin{equation}\label{eq:d4-carleman}
 |u\cdot\omega|\dd u\dd\omega
 =c_d|y|^{-(d-2)}\dd y\,
 \dd\mathcal H^{d-1}_{y^\perp}(w)
 =c_d|y|^{-(d-3)}\delta(y\cdot w)\dd y\dd w.
\end{equation}
The right side is locally integrable at $y=0$ after restriction to a
bounded $w$-box, and \eqref{eq:d4-carleman} extends there as a measure
identity.
\end{lemma}

\begin{proof}
Write $u=a\omega+w$, $w\perp\omega$, on each sign of $a$.  Then
$y=a\omega$, so polar coordinates in $y\in\R^d$ contribute
$|a|^{d-1}\dd a\dd\omega$, whereas the left side contributes
$|a|\dd a\dd\omega\dd w$.  Their ratio is $|y|^{-(d-2)}$.
Moreover
\[
 \int_{\R^d}\delta(y\cdot w)F(w)\dd w
 =|y|^{-1}\int_{y^\perp}F(w)\dd\mathcal H^{d-1}(w),
\]
which proves the delta form.  On a shell
$|y|\sim\rho$, the weighted $y$-mass is $O(\rho^2)$, which proves local
integrability.
\end{proof}

There are two complementary charts.  In the direct chart we use
$(t_1,y,w)$, retain the second-collision indicator and its selected root
inside $Q$, and then enlarge by dropping that indicator.  Lemma
\ref{lem:d4-carleman} gives, on $|y|\sim\rho$,
\begin{equation}\label{eq:j3-direct-shell}
 \mathcal J^{\rm dir}_{3,\rho}(Q)
 \le\cL_\eps^C\rho^2\int_{\Omega'_{\rm dir}}Q^{\rm dir}.
\end{equation}
The independent variables are the displayed $(t_1,y,w)$, with $w$
written in one of finitely many orthonormal frames of $y^\perp$.

In the transformed chart solve \eqref{eq:j3-second-contact} for
\begin{equation}\label{eq:j3-Phi}
 y=\Phi(t_2,\omega_2)
 =\frac{b_m(t_2)-\eps\omega_2}{\Delta}.
\end{equation}
With $a=v_1-v_7$ and $p=a+y$, direct differentiation gives
\begin{equation}\label{eq:j3-Phi-Jacobian}
 \left|\det D_{(t_2,\omega_2)}\Phi\right|
 =\frac{\eps^{d-1}|p\cdot\omega_2|}{|\Delta|^d}.
\end{equation}
Indeed the $d-1$ angular derivatives are $-\eps/\Delta$ on
$T_{\omega_2}\Sph^{d-1}$, while
$\partial_{t_2}y=p/\Delta$.  Hence the transformed measure is
\begin{equation}\label{eq:j3-transformed-density}
 c_d\eps^{d-1}
 \frac{|p\cdot\omega_2|}{|y|^{d-2}|\Delta|^d}
 \dd t_1\dd t_2\dd\omega_2
 \dd\mathcal H^{d-1}_{y^\perp}(w).
\end{equation}

For fixed $t_1$ and $\omega_2$, let
$g(t_2)=y(t_2,\omega_2)\cdot\omega_2$.  Then
\begin{equation}\label{eq:j3-g-derivative}
 g'(t_2)=\frac{p\cdot\omega_2}{\Delta}.
\end{equation}
On each sign-of-$\Delta$ and derivative-sign cell this is a monotone
fractional-linear coordinate.  A vector shell $|y|\sim\rho$ has
$g$-length at most $C\rho$, and therefore
\begin{equation}\label{eq:j3-time-cancellation}
 \int_{\{|y|\sim\rho\}}
 \frac{|p\cdot\omega_2|}{|\Delta|^d}\dd t_2
 \le C_d\mu^{-(d-1)}\rho.
\end{equation}
Together with $|y|^{-(d-2)}\sim\rho^{-(d-2)}$ this proves
\begin{equation}\label{eq:j3-transformed-shell}
 \mathcal J^{\rm tr}_{3,\rho}(Q)
 \le\cL_\eps^{C_d}\eps^{d-1}\mu^{-(d-1)}\rho^{-(d-3)}
 \int_{\Omega'_{\rm tr}}Q^{\rm tr}.
\end{equation}

The following cell decomposition shows why \eqref{eq:j3-transformed-shell}
is an arbitrary-$Q$ operator estimate rather than a $Q=1$ volume calculation.
Refine the
shell by half-open cells $|\Delta|\sim D$ and
$|p\cdot\omega_2|\sim\nu$, including a terminal cell
$\nu\le D\rho$.  On a positive $\nu$-cell,
\eqref{eq:j3-g-derivative} gives a $t_2$-section of length at most
\begin{equation}\label{eq:j3-t2-cell}
 C\min\{1,D\rho/\nu\}.
\end{equation}
The coefficient in \eqref{eq:j3-transformed-density} is constant within
a fixed factor on this cell.  Multiplying it by
\eqref{eq:j3-t2-cell} gives at most
 $C_d\eps^{d-1}\rho^{-(d-3)}D^{-(d-1)}$.  On each positive $\nu$-cell use $s=g/\rho$
on every monotone component.  Then $s$ ranges in a fixed interval, and
\eqref{eq:j3-g-derivative} transforms the density into a quantity bounded by
 $C_d\eps^{d-1}\rho^{-(d-3)}D^{-(d-1)}\,\dd s$.  In the terminal cell retain $t_2$ in its
original fixed interval and use $\nu\le D\rho$ directly.  Along with fixed
boxes for $t_1,\omega_2$ and the coordinates of $w$, these charts construct a
containing domain $\Omega'_{\rm tr}$ independent of $(z_1,z_7)$.  The
rank-drop set $p\cdot\omega_2=0$ has zero collision flux; no inverse Jacobian
is taken there.
The direct chart treats the $y=0$ endpoint.  Thus all substitutions are
valid for measurable $Q\ge0$ by Tonelli and the area formula.

Taking the better chart gives the announced exact shell estimate
\begin{equation}\label{eq:j3-shell}
 \mathcal J_{3,\rho}(Q)\le\cL_\eps^C
 \min\left\{\rho^2,
 \eps^{d-1}\mu^{-(d-1)}\rho^{-(d-3)}\right\}
 \int Q^\sharp.
\end{equation}
The crossover is $r=\eps/\mu$.  Dyadic summation, with the terminal ball
$|y|\le r$ assigned to the direct chart, gives
\begin{equation}\label{eq:j3-summed}
 \mathcal J_3(Q)\le\cL_\eps^C\eps^2\mu^{-2}
 \int Q^\sharp.
\end{equation}
At the ordinary threshold $\mu=\eps^{1/(8d)}$ this is
$\eps^{2-1/(4d)-o(1)}$, which is stronger than the ordinary credit
$\eps^{1/(8d)}$.

\subsection{The other $\{33\}$ orientations}

We give the two remaining type-A outputs in full because the
dimension-independent radial factor in the $j=4$ cone chart must be verified
rather than read formally from the low-dimensional source formula.  Resolve
the lift of the second contact first.  On a fixed lift, changing a
$d$-dimensional output
$v_j-v_1$ to $(t_2,\omega_2)$ has collision-coarea density
\begin{equation}\label{eq:safe-output-second-coarea}
 \eps^{d-1}\frac{|(v_j-v_7)\cdot\omega_2|}{|\Delta|^d}
 \dd t_2\dd\omega_2,
 \qquad |\Delta|\ge\mu.
\end{equation}
The second normal velocity is a numerator.  Bounded time and velocity cells
allow only $\cL_\eps^C$ lifts.

For $j=2$ the output is $u$.  Put $r=\eps/\mu$.  On $|u|\le4r$, do not use
\eqref{eq:safe-output-second-coarea}; drop the second-contact indicator and
use $|y|\le|u|$ to obtain
\begin{equation}\label{eq:j2-direct-ball}
 \int_{|u|\le4r}|y|\dd u\dd\omega
 \le C_d\int_{|u|\le4r}|u|\dd u
 \le C_dr^{d+1}\le C_d\eps^{d-1}\mu^{-d}.
\end{equation}
Here the last inequality uses the ordinary choice
$\mu=\eps^{1/(8d)}$; more generally it is valid whenever
$\eps^2\le\mu\le1$.  The estimate is used only on this range.
On $|u|>4r$, use \eqref{eq:safe-output-second-coarea}.  The first weight
$|y|$ and the second numerator are bounded by the velocity cutoff, while
 $|\Delta|^{-d}\le\mu^{-d}$.  All remaining coordinate boxes have
polylogarithmic volume, and hence
\begin{equation}\label{eq:j2-transformed}
 \mathcal J_{33A,2}^{\rm tr}(Q)
 \le\cL_\eps^{C_d}\eps^{d-1}\mu^{-d}
 \int_{\Omega'_2}Q^\sharp.
\end{equation}
For the arbitrary-kernel form of the direct branch, put
$u=r\widehat u$ with $r=\eps/\mu$.  The Jacobian is $r^d$, while
$|y|\le |u|\le4r$, so the pulled-back density has the factor $Cr^{d+1}$.
Compose $Q$ with this inverse scaling and the remaining recovery maps, and
extend it by zero to the fixed ball $|\widehat u|\le4$ and the fixed angular
boxes.  This proves the same bound for arbitrary $Q$, without extracting a
scalar support volume from its integral.

For $j=4$ the output is $w=u-y\in y^\perp$.  On $|w|\le4r$, the Carleman
formula gives
\begin{equation}\label{eq:j4-direct-ball}
 \int_{|y|\le\cL_\eps^C}|y|^{-(d-2)}\dd y
 \int_{\substack{w\in y^\perp\\|w|\le4r}}
       \dd\mathcal H^{d-1}(w)
 \le\cL_\eps^{C_d}r^{d-1}
 \le\cL_\eps^{C_d}\eps^{d-1}\mu^{-d}.
\end{equation}
This includes $w=0$ and $y=0$; for every $d\ge4$,
$|y|^{-(d-2)}\dd y$ is locally integrable.
For the arbitrary-kernel form, write $y=\sqrt{s}\,e$ and
$w=r\widehat w$ in the moving hyperplane $e^\perp$.  Then
\[
 |y|^{-(d-2)}\dd y=\tfrac12\dd s\,\dd e,
 \qquad
 \dd\mathcal H^{d-1}(w)=r^{d-1}
       \dd\mathcal H^{d-1}(\widehat w).
\]
Finite half-open moving-frame charts put $(s,e,\widehat w)$ in fixed boxes.
Composing $Q$ with the inverse maps and extending it by zero extracts the
factor $r^{d-1}$ from the same joint integral.  Thus
\eqref{eq:j4-direct-ball} is also an arbitrary-$Q$ operator bound.

On $|w|>4r$, the resolved second-contact identity can be written
\begin{equation}\label{eq:j4-wstar}
 w=w_*+\frac{\eps\omega}{\Delta},
 \qquad
 w_*:=\frac{A_m+t_2a-\eps\omega_2}{\Delta},
 \qquad a=v_1-v_7,
\end{equation}
where $A_m$ is fixed on the lift chart.  It follows that
$|w_*|\ge3r$.  Set, with its sign retained,
\[
 \lambda_{\rm sgn}:=\frac{\eps}{\Delta |w_*|},
 \qquad |\lambda_{\rm sgn}|\le\frac13 .
\]
Using the delta form of the
Carleman identity,
\[
 |y|^{-(d-3)}\delta(y\cdot w)\dd w\dd y,
\]
the relation $y\cdot w=0$ restricts $y$ to the cone
\[
 \Gamma_{e,\lambda_{\rm sgn}}
 =\{y:e\cdot y+\lambda_{\rm sgn}|y|=0\},
 \qquad e=w_*/|w_*|.
\]
The defining gradient is uniformly separated from zero because
$|\lambda_{\rm sgn}|\le1/3$.  We use Riemannian hypersurface measure on the
cone.  In any one of the fixed cone-axis charts, the Gram density and the
coarea factor are bounded above and below by absolute constants.  Combining
this chartwise coarea with \eqref{eq:safe-output-second-coarea} therefore
gives the following upper density (with a chart-independent constant):
\begin{equation}\label{eq:j4-cone-density}
 C_d\eps^{d-1}
 \frac{|(v_4-v_7)\cdot\omega_2|}
 { |\Delta|^d|w_*||y|^{d-3}}
 \dd t_1\dd t_2\dd\omega_2
 \dd\mathcal H^{d-1}_{\Gamma_{e,\lambda_{\rm sgn}}}(y).
\end{equation}
The factor $|y|^{-(d-3)}$ is integrable, uniformly in the cone parameters:
\begin{equation}\label{eq:j4-cone-radial}
 \sup_{e,|\lambda_{\rm sgn}|\le1/3}
 \int_{\Gamma_{e,\lambda_{\rm sgn}}\cap B(0,\cL_\eps^C)}
 |y|^{-(d-3)}\dd\mathcal H^{d-1}(y)
\le\cL_\eps^C.
\end{equation}
Indeed, polar coordinates on the $(d-1)$-dimensional cone leave the exact
radial density
\[
 r^{d-2}r^{-(d-3)}\dd r=r\dd r.
\]

It remains to control $|w_*|^{-1}$.  Since
$v_4-v_7=a+w$ and $|w-w_*|\le r$,
\begin{equation}\label{eq:j4-wstar-numerator}
 \frac{|(v_4-v_7)\cdot\omega_2|}{|w_*|}
 \le C+\frac{|a|}{|w_*|}.
\end{equation}
For fixed $(t_1,\omega_2)$ put
$z(t_2)=A_m+t_2a-\eps\omega_2=\Delta w_*$.  The transformed cell excludes
$|z|<3\eps$, and therefore
\begin{equation}\label{eq:j4-affine-log}
 \int\frac{|a|}{|w_*|}\dd t_2
 =\int\frac{|a||\Delta|}{|A_m+t_2a-\eps\omega_2|}\dd t_2
 \le\cL_\eps^C.
\end{equation}
Indeed, after splitting the fixed vector into parts parallel and
orthogonal to $a$, the nonzero-$a$ integral is bounded by
$\int (d_0^2+q^2)^{-1/2}\dd q$ away from a ball of radius $3\eps$;
the case $a=0$ is immediate.  This proves \eqref{eq:j4-affine-log}
uniformly, including tangency.

For the arbitrary-$Q$ form, refine into half-open cells
$|\Delta|\sim D$, $|w_*|\sim R$ and $|y|\sim\rho$.  The affine-line
section in $t_2$ has length at most
$\cL_\eps^C\min\{1,DR/|a|\}$, with the second entry omitted if $a=0$.
We record the pointwise normalization needed for an arbitrary kernel.  If
$|a|\le DR$, retain $t_2$ itself; then
$1+|a|/R\le1+D$.  If $|a|>DR$, each of the at most two affine-line
components is parametrized by
\[
 s=\frac{a\cdot z(t_2)}{DR|a|},
 \qquad
 \dd t_2=\frac{DR}{|a|}\dd s,
 \qquad
 \left(1+\frac{|a|}{R}\right)\dd t_2
 \le(1+D)\dd s .
\]
In both cases the normalized time coordinate lies in a fixed box.  On a
cone-axis chart write $y=r\vartheta$ and set $q=r^2/\rho^2$ on the radial
cell.  The weighted radial density becomes
\[
 r^{d-2}r^{-(d-3)}\dd r=r\dd r
   =\tfrac12\rho^2\dd q,
\]
and the angular variables range over fixed boxes with uniformly bounded
chart densities.  Since $D$ and $\rho$ are bounded by the velocity and time
cutoffs, the factors $1+D$ and $\rho^2$ are polylogarithmic.  Compose $Q$
with these inverse maps and extend it by zero to the fixed
time--radial--angular boxes.  Thus the pulled-back density is pointwise at
most $\cL_\eps^C\eps^{d-1}D^{-d}$, independently of the fixed boundary
states; no support volume is extracted through $Q$.  Equations
\eqref{eq:j4-cone-density}--\eqref{eq:j4-affine-log}, together with
$|\Delta|^{-d}\le\mu^{-d}$, give the transformed bound.

Combining the direct and transformed charts, both outputs satisfy
\begin{equation}\label{eq:j2j4-full-range}
 \mathcal J_{33A,j}(Q)\le
 \cL_\eps^C\eps^{d-1}\mu^{-d}\int Q^\sharp,
 \qquad j\in\{2,4\}.
\end{equation}
The zero-output and lower radial endpoints are in the direct cells; cone
axes are covered by the finite charts, and both first- and second-contact
grazing sets carry vanishing collision flux.

\begin{lemma}[Fixed-time type-A sphere operators]
\label{lem:fixed-time-33A-j2-j3}
Assume $0<\eps<\eps_0$, $\eps^2\le\mu\le1$, and fix both atom times
$(t_1,t_2)$ in an ordinary type-$\{33{\rm A}\}$ cell with
$|t_1-t_2|\ge\mu$.  On every fixed lift, the $j=2$, serial $j=3$, and $j=4$
conditional Radon measures have explicit selected Borel versions for which,
for every
measurable $Q\ge0$,
\begin{align}
 \mathcal J_{33A,2}^{\,t_1,t_2}(Q)
 &\le \cL_\eps^C\eps^{d-1}\mu^{-d}
       \int_{\Omega'_{2,{\rm ft}}}Q^\sharp,
 \label{eq:j2-fixed-time-sphere}\\
 \mathcal J_{33A,3}^{\,t_1,t_2}(Q)
 &\le \cL_\eps^C\eps\mu^{-2}
       \int_{\Omega'_{3,{\rm ft}}}Q^\sharp.
 \label{eq:j3-fixed-time-sphere}\\
 \mathcal J_{33A,4}^{\,t_1,t_2}(Q)
 &\le \cL_\eps^C\eps^{d-2}\mu^{-(d-1)}
       \int_{\Omega'_{4,{\rm ft}}}Q^\sharp.
 \label{eq:j4-fixed-time-sphere}
\end{align}
Here a finite constant-coefficient sum over lift, angular and moving-frame
charts is understood, as in \eqref{eq:source-elementary-operator}.  Each
$\Omega'_{j,{\rm ft}}$ is a fixed product of spherical coordinate patches
and cutoff boxes, independent of the two fixed exterior states and of the
fixed times.  In particular these are arbitrary-kernel operator estimates,
not estimates obtained by integrating either distinguished time or by
putting $Q=1$.
\end{lemma}

\routineproofinresource{supp:proof:lem:fixed-time-33A-j2-j3}

For type B, both atoms are C-atoms, the higher fixed end is top and the
lower fixed end is bottom.  There are exactly three slot orientations:
both fixed ends are serial with the common bond, exactly one is serial, or
neither is serial.  We first record a common raw measure, so that the three
calculations below are operator identities for an arbitrary joint kernel and
not merely volume computations with \(Q=1\).  For bottom states
\((z_a,z_b)\), top states \((z_c,z_d)\), time \(t\), and normal \(\omega\),
put
\begin{align}
 \mathscr C_t(z_a,z_b;z_c,z_d;\omega)
 &:=[(v_a-v_b)\cdot\omega]_+
 \delta^{(d)}(x_c-x_a+t(v_c-v_a))\notag\\
 &\quad\times\delta^{(d)}(x_d-x_b+t(v_d-v_b))
 \delta^{(d)}\!\left(v_c-v_a+[(v_a-v_b)\cdot\omega]\omega\right)\notag\\
 &\quad\times\delta^{(d)}\!\left(v_d-v_b-[(v_a-v_b)\cdot\omega]\omega\right)
 \delta^{(d)}(x_a-x_b+t(v_a-v_b)-\eps\omega).
 \label{eq:typeB-raw-C-kernel}
\end{align}
The fixed states are \(z_1,z_7\), and the raw variables are
\begin{equation}\label{eq:typeB-raw-variables}
 \mathbf r=(t_1,t_2,\omega,\eta,z_2,z_3,z_4,z_5,z_6)
 \in I^2\times(\Sph^{d-1})^2\times(\R^{2d})^5.
\end{equation}
Here and below \(\dd\mathbf r\) means Lebesgue measure in every \(z_i\)
and \(t_i\), and surface measure in \(\omega,\eta\).  After the two sphere
disintegrations, including the original single \(\eps^{-(d-1)}\)
full-component normalization, the three exact raw measures are
\begin{align}
 \mathcal J_{33B}^{(1)}(Q)
 &=\eps^{d-1}\!\int \mathscr C_{t_1}(z_1,z_2;z_3,z_4;\omega)
 \mathscr C_{t_2}(z_5,z_3;z_6,z_7;\eta)Q\,\dd\mathbf r,
 \label{eq:typeB-raw-case1}\\
 \mathcal J_{33B}^{(2)}(Q)
 &=\eps^{d-1}\!\int \mathscr C_{t_1}(z_1,z_2;z_3,z_4;\omega)
 \mathscr C_{t_2}(z_5,z_3;z_7,z_6;\eta)Q\,\dd\mathbf r,
 \label{eq:typeB-raw-case2}\\
 \mathcal J_{33B}^{(3)}(Q)
 &=\eps^{d-1}\!\int \mathscr C_{t_1}(z_1,z_2;z_3,z_4;\omega)
 \mathscr C_{t_2}(z_4,z_5;z_6,z_7;\eta)Q\,\dd\mathbf r.
 \label{eq:typeB-raw-case3}
\end{align}
Every transport and scattering delta is displayed in
\eqref{eq:typeB-raw-C-kernel}; in particular, the \(Q\) in these three
formulas is the original \(Q(z_1,\ldots,z_7,t_1,t_2)\), with no
factorization or supremum.  Put $\lambda_{\rm small}:=\eps^{1/(8d)}$.
The type-B nondegenerate cell is
\begin{equation}\label{eq:typeB-nondegeneracy}
 |t_1-t_2|\ge\lambda_{\rm small},\qquad
 |v_e-v_f|\ge\lambda_{\rm small}
\end{equation}
for every distinct velocity pair used as a Carleman or spherical pivot.
The three push-forwards below include their independent coordinates and fix
which rank drops are delegated to the complementary small-variable cells.

If both fixed ends are serial, the \(z_4,z_6\) blocks of
\eqref{eq:typeB-raw-case1} are unit-Jacobian transport/scattering
substitutions.  The two Carleman disintegrations, including their incoming
fluxes, leave respectively
\[
 w\in(v_7-v_3)^\perp,\qquad
 \theta\in(v_1-v_3)^\perp,
\]
and the factors \(|v_7-v_3|^{-(d-2)}\) and \(|v_1-v_3|^{-(d-2)}\).
After these substitutions the last square constraint and its pivot
variables are
\begin{equation}\label{eq:typeB-case1-pivot-map}
 G_1(x_3,v_3)=
 \binom{x_3-x_1+t_1(v_3-v_1)}
       {x_7-x_3+t_2(v_7-v_3)},
 \qquad p_1=(x_3,v_3)\in\R^{2d}.
\end{equation}
Block row subtraction gives the exact determinant
\begin{equation}\label{eq:typeB-case1-pivot-determinant}
 |\det D_{p_1}G_1|=|t_1-t_2|^d.
\end{equation}
Thus \(G_1=0\) has one solution on every nonzero time-difference cell.  On
each orthogonal-frame chart the independent variables are
\[
 (t_1,t_2,w,\theta)\in I^2\times
 (v_7-v_3)^\perp\times(v_1-v_3)^\perp,
\]
and the push-forward density multiplying \(Q\circ\Psi_1\) is, up to
universal frame-chart constants,
\begin{equation}\label{eq:typeB-first-density}
 \eps^{d-1}|t_1-t_2|^{-d}
 |v_7-v_3|^{-(d-2)}|v_1-v_3|^{-(d-2)}.
\end{equation}
Here \(v_3\) and every eliminated state are the explicit functions of
\((t_1,t_2,w,\theta)\) obtained from the two scattering rules.  Thus \(Q\)
is composed with the full reconstruction map \(\Psi_1\) and is not
averaged or supremized.  The only rank drops in this chart are
\(t_1=t_2\), \(v_7=v_3\), and \(v_1=v_3\), all excluded by
\eqref{eq:typeB-nondegeneracy}; their small half-open neighborhoods, rather
than the zero sets alone, are retained in the complementary good-support
cells.

If exactly one fixed end is serial, write \(\eta\) for the lower contact
normal, \(\alpha=(v_6-v_7)\cdot\eta\), and let
\(\theta\in(v_1-v_3)^\perp\).  In
\eqref{eq:typeB-raw-case2}, the \(z_4\) block is a unit-Jacobian
substitution and the upper Carleman map gives
\(|v_1-v_3|^{-(d-2)}\).  Resolve
\[
 v_6-v_7=\alpha\eta+\sigma,\qquad \sigma\in\eta^\perp.
\]
The tangential part of the remaining transport constraint uses the $d-1$
coordinates of \(\sigma\) as pivots, with exact determinant
\begin{equation}\label{eq:typeB-case2-tangential-determinant}
 |t_1-t_2|^{d-1}.
\end{equation}
Indeed, inverse elastic reconstruction at the lower atom gives
\(v_3=v_7+\sigma\), so the complete remaining vector constraint is
\begin{equation}\label{eq:typeB-case2-constraint-map}
 G_2(\sigma,\eta)
 :=\zeta+(t_1-t_2)\sigma-\eps\eta=0,
 \qquad \sigma\in\eta^\perp .
\end{equation}
Its \(\eta^\perp\)-projection is the $(d-1)$-dimensional pivot just
described; its last component is the scalar constraint
\(h_2(\eta)=\eta\cdot\zeta-\eps=0\), where
\begin{equation}\label{eq:typeB-zeta}
 \zeta=(x_7+t_1v_7)-(x_1+t_1v_1)
       =(t_1-t_2)(v_7-v_3)+\eps\eta.
\end{equation}
On \(h_2=0\),
\begin{equation}\label{eq:typeB-case2-sphere-jacobian}
 |\nabla_{\Sph^{d-1}}h_2|=\sqrt{|\zeta|^2-\eps^2}.
\end{equation}
The level is a $(d-2)$-sphere of Euclidean radius
\[
 r=(1-\eps^2|\zeta|^{-2})^{1/2}
\]
inside an affine translate of \(\zeta^\perp\), and its area element is
\(r^{d-2}\dd\eta^\perp\).  Division by
\eqref{eq:typeB-case2-sphere-jacobian} gives \(r^{d-3}/|\zeta|\).  Sphere coarea
therefore leaves the independent coordinates
\((t_1,t_2,\alpha,\eta^\perp,\theta)\), with
\(\eta^\perp\in\Sph^{d-2}\subset\zeta^\perp\), and density
\begin{equation}\label{eq:typeB-second-density}
 \eps^{d-1}\alpha_+|t_1-t_2|^{-(d-1)}|v_1-v_3|^{-(d-2)}
 \frac{(1-\eps^2|\zeta|^{-2})^{(d-3)/2}}{|\zeta|}
 \dd t_1\dd t_2\dd\alpha\dd\eta^\perp\dd\theta.
\end{equation}
On \eqref{eq:typeB-nondegeneracy},
\(|\zeta|\ge\lambda_{\rm small}^2-\eps\ge\lambda_{\rm small}^2/2\), provided
\(\eps\lambda_{\rm small}^{-2}\le1/2\).  The radial factor is at most one, and the
total inverse power is therefore at most
\(C_d\lambda_{\rm small}^{-(2d-1)}\).  Notice that
the density vanishes, rather than blows up, at the endpoint
\(|\zeta|=\eps\); the regular-level measures converge to zero there.  The
collapsed level is therefore defined by monotone regular-level exhaustion
and carries no atom.

If neither fixed end is serial, start from
\eqref{eq:typeB-raw-case3}.  Five position blocks are unit translations;
the two occurrences of \(x_4\) leave, with the fixed-end affine traces
\[
 y^*:=x_1+t_1v_1,\qquad z^*:=x_7+t_2v_7,
\]
\begin{equation}\label{eq:typeB-third-position}
 (t_2-t_1)v_4+y^*-z^*-\eps(\omega+\eta)=0.
\end{equation}
The velocity deltas first substitute \(v_3,v_6\).  Replacing the two source
fibre scalars by their negatives (a unit-Jacobian reparametrization), write
\(v_2=v_4-\alpha_1\omega\), \(v_5=v_4-\alpha_2\eta\), and impose
\begin{equation}\label{eq:typeB-case3-sphere-sections}
 (v_4-v_1)\cdot\omega=0,
 \qquad (v_4-v_7)\cdot\eta=0.
\end{equation}
Thus, before coarea, the complete remaining constraint map is
\begin{equation}\label{eq:typeB-case3-constraint-map}
 \mathcal G_3(v_4,\omega,\eta)
 :=
 \begin{pmatrix}
  (t_2-t_1)v_4+y^*-z^*-\eps(\omega+\eta)\\
  (v_4-v_1)\cdot\omega\\
  (v_4-v_7)\cdot\eta
 \end{pmatrix}=0 .
\end{equation}
The two spherical normal Jacobians are
\(|v_4-v_1|\) and \(|v_4-v_7|\).  Parametrize their equators by
\(\omega',\eta'\in\Sph^{d-2}\) and the two free normal scalars by
\(\alpha_1,\alpha_2\).  For fixed
\((t_1,t_2,\omega',\eta')\), equation
\eqref{eq:typeB-third-position} is the zero of
\begin{equation}\label{eq:typeB-third-map}
 \Xi(v_4)=v_4+
 \frac{y^*-z^*-\eps(R_{v_4-v_1}\omega'
                         +R_{v_4-v_7}\eta')}{t_2-t_1}.
\end{equation}
On a fixed rotation chart,
\(\|D\Xi-I\|\le C_d\eps\lambda_{\rm small}^{-2}\).  For sufficiently small
$\eps=\eps(d)$ one has
\(C_d\eps\lambda_{\rm small}^{-2}\le1/2\).  Then
\(1/2\le|\det D\Xi|\le2\).  Moreover
\[
 |\Xi(v)-\Xi(w)|\ge\tfrac12|v-w|
\]
on each convex velocity chart, so this pivot has at most one branch there.
In the adapted coordinates
\((v_4,\omega',\eta')\), the square part of
\eqref{eq:typeB-case3-constraint-map} therefore has determinant
\begin{equation}\label{eq:typeB-case3-full-pivot}
 |t_2-t_1|^d|\det D\Xi|\,
 |v_4-v_1|\,|v_4-v_7|.
\end{equation}
The two original collision fluxes are not part of the spherical coarea
Jacobian and must be retained.  Writing
\(\alpha_{j,+}:=\max\{\alpha_j,0\}\), the remaining push-forward density is
bounded by
\begin{equation}\label{eq:typeB-third-density}
 C_d\eps^{d-1}|t_1-t_2|^{-d}|v_4-v_1|^{-1}|v_4-v_7|^{-1}
 \alpha_{1,+}\alpha_{2,+}
 \dd t_1\dd t_2\dd\alpha_1\dd\alpha_2\dd\omega'\dd\eta'.
\end{equation}

In the three cases the inverse degrees are respectively \(3d-4\), \(2d-1\), and
\(d+2\).  Choose finite orthogonal-frame and rotation charts, and write every
moving $(d-1)$-plane in its \(\R^{d-1}\) frame.  Enumerate the finite atlas
$(U_k)_{k=1}^{K}$ and replace it by the subordinate half-open Borel
refinement
\[
 U_1^\circ=U_1,\qquad
 U_k^\circ=U_k\setminus\bigcup_{j<k}U_j\quad(2\le k\le K).
\]
Thus chart interiors never overlap in the branch sum; restricting a smooth
chart to $U_k^\circ$ changes neither its Gram bounds nor the area formula.
If \(V=\cL_\eps^C\), fixed
containing domains, chosen before the concrete values of \(z_1,z_7\), are
\begin{align}
 \Omega'_1&=I^2\times[-V,V]^{d-1}\times[-V,V]^{d-1},\notag\\
 \Omega'_2&=I^2\times[-V,V]\times\Sph^{d-2}\times[-V,V]^{d-1},\notag\\
 \Omega'_3&=I^2\times[-V,V]^2\times(\Sph^{d-2})^2.
 \label{eq:typeB-fixed-output-domains}
\end{align}
The branch kernel is \(Q\circ\Psi_i\) on the part whose full reconstruction
lies in the original half-open cell and is zero elsewhere.  Consequently
all free scalars and orthogonal coordinates stay in the original cutoff,
while the output domain itself is boundary-value independent.

The sets \(t_1=t_2\) and the zero relative speeds are not discarded: their
small half-open neighborhoods are precisely the ordinary small-time or
small-relative-speed alternatives of \eqref{eq:one-atom-good-section}.
Exact grazing carries zero flux.  Frame-chart boundaries are assigned once
by the half-open convention rather than assumed null in a conditional
fibre.  On the third containing domain the velocity cutoff gives
\(\alpha_{1,+}\alpha_{2,+}\le V^2\), which is absorbed in
\(\cL_\eps^C\); no flux is discarded.  The area formula and Tonelli's
theorem applied to the same \(Q\circ\Psi_i\) now give
\begin{equation}\label{eq:typeB-bound}
 \mathcal J_{33B}(Q)\le
 \cL_\eps^C\eps^{d-1}\lambda_{\rm small}^{-(3d-4)}\int Q^\sharp.
\end{equation}
These are the three coarea cases of
\cite[Proposition 9.2]{DHMlong}, reconstructed above for every fixed $d\ge4$.

For a $\{44\}$ component, choose the two free ends which become the fixed
ends of a type-A pair.  Work on one of the already selected half-open lift
cells and denote its lifted relative position again by $x_1-x_7$.  Put
$\lambda_{\rm close}:=\eps^{1-1/(8d)}$.  On
\begin{equation}\label{eq:44-close-cell}
 \max\{|x_1-x_7|,|v_1-v_7|\}\le\lambda_{\rm close},
\end{equation}
make the linear change
\[
 c_x=\tfrac12(x_1+x_7),\quad c_v=\tfrac12(v_1+v_7),\qquad
 \widehat r_x=\frac{x_1-x_7}{\lambda_{\rm close}},\quad
 \widehat r_v=\frac{v_1-v_7}{\lambda_{\rm close}}.
\]
Its Jacobian is $C_d\lambda_{\rm close}^{2d}$, while
$(\widehat r_x,\widehat r_v)$ ranges in a fixed product of unit balls and
the centre state remains in a fixed cutoff box.  For an arbitrary
nonnegative joint kernel, define $Q^\sharp$ by composing with this inverse
linear map and with the conditional type-A recovery, and extend it by zero
outside the fixed product box.  The type-A estimate is uniform in the two
end states, so Tonelli may be applied before the rescaled relative variables
are integrated.  Consequently the close-state volume is extracted from the
same joint integral, rather than from a supremum of $Q$.  The two full
baselines give
\begin{equation}\label{eq:44-bound}
 \mathcal J_{44}(Q)\le
 \cL_\eps^C\eps^{-2(d-1)}\lambda_{\rm close}^{2d}\int Q^\sharp.
\end{equation}
The two cutoffs are deliberately different: the inverse-pivot cutoff is
$\lambda_{\rm small}$, while the close-endpoint cutoff is
$\lambda_{\rm close}$.  Their exponents are recorded in
Table~\ref{tab:ordinary-exponent-check} below.

\subsection{The weight--time inequality with an arbitrary time owner}

The two auxiliary estimates are proved directly in the fixed dimension under consideration.  The reason is that the time variables used by an exceptional
packet are omitted from the later molecule integral.  What remains is not a
product simplex, but the existential projection of the original partial-order
domain.

An \emph{extended binary time tree} is a rooted atom tree lying in one layer,
oriented so that every child time is below its parent time.  The time-reversed
notion is called an upside-down extended binary tree.  For an atom $a$ in such
a tree $T$, let $\sigma_T(a)$ be the number of descendants of $a$, including
$a$, and put
\[
 C_T:=\prod_{a\in T}\sigma_T(a).
\]

\begin{lemma}[Tree simplex and all existential projections]
\label{lem:tree-time-owner}
Let $T$ be an extended binary time tree in an interval of length $\tau$.
For every $A\subset T$,
\begin{align}
 \int_{\mathcal D_T}1\dd t_T
 &=\tau^{|T|}C_T^{-1},
 \label{eq:tree-time-simplex}\\
 \int_{\widetilde{\mathcal D}_{T,A}}1\dd t_{T\setminus A}
 &\le\tau^{|T|-|A|}C_T^{-1}
 \frac{|T|!}{(|T|-|A|)!}.
 \label{eq:tree-projected-simplex}
\end{align}
Here $\widetilde{\mathcal D}_{T,A}$ consists precisely of the retained times
for which there exists a choice of the omitted $A$-times placing the complete
vector in $\mathcal D_T$.  Both statements hold for upside-down trees.
\end{lemma}

\begin{proof}
Translation and dilation reduce the proof to the interval $(0,1)$.  Let $r$
be the root and let $T_1,T_2$ be its nonempty child trees, allowing one of
them to be empty.  If $m_j=|T_j|$, integration first below $t_r$ gives
\[
 \int_0^1
 \frac{t_r^{m_1}}{C_{T_1}}
 \frac{t_r^{m_2}}{C_{T_2}}\dd t_r
 =\frac1{C_{T_1}C_{T_2}(m_1+m_2+1)}
 =C_T^{-1},
\]
which proves \eqref{eq:tree-time-simplex} inductively.

For the projection estimate write $A_j=A\cap T_j$ and $k_j=|A_j|$.  Suppose
first that $r\notin A$.  The induction hypothesis on the interval
$(0,t_r)$, followed by integration of $t_r$, yields
\begin{multline}
 \int_{\widetilde{\mathcal D}_{T,A}}1\dd t_{T\setminus A}
 \le C_{T_1}^{-1}C_{T_2}^{-1}
 \frac{m_1!m_2!}{(m_1-k_1)!(m_2-k_2)!}\\
 {}\times\frac1{m_1+m_2-k_1-k_2+1}.
 \label{eq:tree-projection-root-kept}
\end{multline}
The elementary injection obtained by interlacing two ordered lists gives
\begin{equation}\label{eq:factorial-interlacing}
 \frac{m_1!m_2!}{(m_1-k_1)!(m_2-k_2)!}
 \le
 \frac{(m_1+m_2)!}{(m_1+m_2-k_1-k_2)!}.
\end{equation}
Indeed, the left side counts ordered selections made separately from two
disjoint lists, and forgetting the separation injects them into ordered
selections from their union.  Combining
\eqref{eq:tree-projection-root-kept}--\eqref{eq:factorial-interlacing} and
$C_T=(m_1+m_2+1)C_{T_1}C_{T_2}$ proves the claim.

If $r\in A$, the root is not integrated.  Existence of an admissible root
time only enlarges each child interval to $(0,1)$, so the two induction
bounds give
\[
 C_{T_1}^{-1}C_{T_2}^{-1}
 \frac{m_1!m_2!}{(m_1-k_1)!(m_2-k_2)!}.
\]
Now $|A|=k_1+k_2+1$; the same interlacing inequality and the factor
$m_1+m_2+1$ in $C_T$ give exactly the asserted factorial.  Restoring
$\tau$ gives one factor $\tau$ for each retained variable.  Reversing all
inequalities proves the upside-down case.
\end{proof}

\begin{lemma}[Weight paid by the same tree simplex]
\label{lem:tree-weight-owner}
For an extended binary tree $T$, with $E_T^-$ the bottom ends of its leaf
atoms,
\begin{equation}\label{eq:tree-weight-proof}
 \prod_{a\in T}(1+|v_{a,1}-v_{a,2}|)
 \le C^{|T|}C_T^{5/6}
 \exp\left(C\sum_{e\in E_T^-}|v_e|^{4/3}\right).
\end{equation}
For an upside-down tree the same estimate holds with its top leaf ends.
\end{lemma}

\begin{proof}
Energy conservation at $a$ and all of its descendants gives
\begin{equation}\label{eq:descendant-energy}
 (1+|v_{a,1}-v_{a,2}|)^2
 \le 2+2\sum_{(b,e)\in\mathcal L(a)}|v_e|^2,
\end{equation}
where $\mathcal L(a)$ contains the at most two bottom ends at every leaf
below $a$.  Factor $\sigma_T(a)^{5/3}$ from the right-hand side.  Hölder's
inequality in the finite leaf sum and $1+Z^{3/2}\le Ce^Z$ give
\begin{multline}
 \prod_{a\in T}(1+|v_{a,1}-v_{a,2}|)^2\\
 \le C^{|T|}C_T^{5/3}
 \exp\left(
 \sum_{(b,e)}|v_e|^{4/3}
 \sum_{a\preceq b}\sigma_T(a)^{-10/9}\right).
 \label{eq:tree-weight-square}
\end{multline}
Along the ancestor chain of a fixed leaf, the positive integers
$\sigma_T(a)$ are strictly increasing.  Hence the inner sum is bounded by
$\sum_{n\ge1}n^{-10/9}$.  Taking square roots proves
\eqref{eq:tree-weight-proof}.  Energy and time reversal prove the second
statement.
\end{proof}

The application of these tree statements to the actual molecule is where the
forest and layer hypotheses enter.

\begin{lemma}[Alternating tree decomposition]
\label{lem:alternating-tree-decomposition}
Every one-layer atom tree can be partitioned into extended binary trees and
upside-down extended binary trees $T^{(j)}$ so that the exterior leaf ends
used in Lemma~\ref{lem:tree-weight-owner} are genuine ends of the original
one-layer tree.  Consequently, with
$C_{\rm lay}=\prod_j C_{T^{(j)}}$, the products of
\eqref{eq:tree-weight-proof} and the bounds
\eqref{eq:tree-projected-simplex} hold for the complete one-layer tree.
\end{lemma}

\begin{proof}
Choose any atom $a$.  Take all descendants of $a$; this is an extended
binary tree and all of its bottom leaves are original bottom ends.  Removing
it leaves tree components, each attached through one bottom end at its first
atom.  In each component take all ancestors of that first atom.  These form
upside-down extended binary trees whose top leaves are original top ends.
After removing them, repeat in the downward direction.  Each round removes
at least one atom, so the construction terminates.  A connecting bond is the
boundary end of exactly the two consecutive pieces which it joins and of no
other piece.  Hence its energy exponential is repeated at most twice,
uniformly in the number of atoms; replacing the numerical constant $C$ by
$2C$ absorbs this duplication before the later layer cancellation.  The time domains of the pieces are enlarged
when connecting inequalities are discarded, hence the product of the
projection bounds is an upper bound, which is the direction needed here.
\end{proof}

For a one-layer CH molecule, remove one collision from each independent
recollision cycle.  Source-B cluster bookkeeping permits at most
\(\Gamma\mathsf R\) such removals, and the remaining molecule is a forest.
Let \(D_{\rm cyc}\) be the set of removed collision atoms.  Their time
variables are accounted for before the forest estimate: if
\(a\in D_{\rm cyc}\setminus A\), discard its order constraints and
integrate \(t_a\) over its layer interval, giving the trivial factor
\(\tau\); if \(a\in D_{\rm cyc}\cap A\), then \(t_a\) is an omitted
existential witness and contributes no integration factor.  Consequently
the removed atoms contribute exactly
\[
 \tau^{|D_{\rm cyc}|-|D_{\rm cyc}\cap A|}
\]
to the enlarged projected-time bound.  Multiplying this by the forest
projection restores the exponent
\(|M_{\rm lay}|-|A\cap M_{\rm lay}|\) in
\eqref{eq:one-layer-time-full}; no deleted cycle time is silently counted
as a forest coordinate.
At a removed collision, energy conservation transfers every new boundary
velocity to the bottom particle velocities of that cluster.  The velocity
cutoff therefore costs at most $\cL_\eps^{C^*\mathsf R}$, while
Lemmas~\ref{lem:tree-time-owner}--\ref{lem:alternating-tree-decomposition}
give, for a scalar $C_{\rm lay}\ge1$,
\begin{align}
 \prod_{a\in M_{\rm lay}}(1+|v_{a,1}-v_{a,2}|)
 &\le C^{|M_{\rm lay}|}\cL_\eps^{C^*\mathsf R}
 C_{\rm lay}^{5/6}
 \exp\left(C\sum_{e\in E_{\rm lay}}|v_e|^{4/3}\right),
 \label{eq:one-layer-weight-full}\\
 \int_{\widetilde{\mathcal D}_{A,\rm lay}}1
 \dd t_{M_{\rm lay}\setminus A}
 &\le \tau^{|M_{\rm lay}|-|A\cap M_{\rm lay}|}C_{\rm lay}^{-1}
 \cL_\eps^{C^*|A\cap M_{\rm lay}|}.
 \label{eq:one-layer-time-full}
\end{align}
In the second line we used $|M|\le\cL_\eps^{C^*}$ only to absorb the
factorial in \eqref{eq:tree-projected-simplex}.  The exponent of $\tau$ is
the exact scaling dimension of the integrated variables.  In particular,
no lower bound on $\tau$ is used.

It remains to restore the parameter $\gamma$.  This must be done before
absorbing the $4/3$-power exponential.  Put
$\widehat v_e=\gamma^{1/2}v_e$.  Since $0<\gamma\le1$,
\begin{equation}\label{eq:gamma-weight-scaling}
 1+|v_{a,1}-v_{a,2}|
 \le \gamma^{-1/2}
 (1+|\widehat v_{a,1}-\widehat v_{a,2}|).
\end{equation}
Energy conservation is homogeneous, so the complete $\gamma=1$ tree
argument applies to $\widehat v$.  This produces
$\exp(C\sum|\widehat v_e|^{4/3})$, not
$\exp(C\sum|v_e|^{4/3})$; hence the elementary inequality
\begin{equation}\label{eq:scaled-four-thirds-absorption}
 e^{C|\widehat v|^{4/3}}\le C_0e^{|\widehat v|^2/20}
 =C_0e^{\gamma|v|^2/20}
\end{equation}
has a constant independent of $\gamma$.  Thus
\eqref{eq:gamma-weight-scaling}, rather than a false polynomial bound for
the unscaled exponential, is the sole source of
$\gamma^{-|M|/2}$.

Now multiply \eqref{eq:one-layer-weight-full} and
\eqref{eq:one-layer-time-full} over all layers and put
$C_M=\prod_{\rm lay}C_{\rm lay}$.  A velocity which is an end of a layer is
either a true end of $M$ or lies on a bond between two layers.  True top and
bottom energies are equal.  Every interlayer bond belongs to a particle line
which is incoming in some layer; there are $\mathsf U$ such incidences.
For a true end use \eqref{eq:scaled-four-thirds-absorption}.  For an
interlayer representative use instead
\[
 e^{C|\widehat v|^{4/3}}
 \le C^*e^{|\widehat v|^2/(20dL)},
\]
where $L$ is the fixed number of time layers and $C^*$ is allowed to
depend on it.  Iterated energy conservation bounds the representatives from
each layer by the true incoming energy; summing over at most $L$
layers leaves at most $\gamma(10d)^{-1}$ times that energy.  Enlarging the
harmless numerical constants gives
\begin{equation}\label{eq:all-layer-weight-full}
 \prod_{a\in M}(1+|v_{a,1}-v_{a,2}|)
 \le(C\gamma^{-1/2})^{|M|}\cL_\eps^{C^*\mathsf R}
 (C^*)^{\mathsf U}C_M^{5/6}
 \exp\left(\gamma\sum_{e\in E^-_{\rm end}}|v_e|^2\right).
\end{equation}
Set $B_M=C_M^{5/6}\ge1$.  Raising the product of
\eqref{eq:one-layer-time-full} to the power $5/6$ cancels $B_M$ and gives
the exact projection estimate \eqref{eq:source-weight-projection}.  Taking
$A=\varnothing$ gives \eqref{eq:source-weight-2}; equation
\eqref{eq:all-layer-weight-full} is \eqref{eq:source-weight-1}.  This proves
the full weight proposition for every atom set $A$, with exactly the
existential time projection required by the time owner.  No step used a
restriction on $d$ or a lower bound on the layer length.

\begin{lemma}[Distinguished-time projection]
\label{lem:distinguished-time-weight}
Let $R\subset M$ be a set of time variables which a local operator retains
as genuine output coordinates, and put $A=M\setminus R$.  Then
\begin{equation}\label{eq:distinguished-time-projection}
 B_M\left(\int_{\widetilde{\mathcal D}_{A}}1\,\dd t_R\right)^{5/6}
 \le \tau^{5|R|/6}\cL_\eps^{C^*|A|}.
\end{equation}
Writing \(V_R=\int_{\widetilde{\mathcal D}_{A}}1\,\dd t_R\), one has the
direct consequence
\begin{equation}\label{eq:distinguished-full-time}
 V_R
 \le B_M^{-6/5}\tau^{|R|}\cL_\eps^{C^*|A|},
 \qquad
 B_MV_R
 \le\tau^{|R|}\cL_\eps^{C^*|A|}.
\end{equation}
Here the harmless enlargement of \(C^*\) absorbs the exponent \(6/5\).
If \(R\) contains a distinguished subset \(R_0\) in bijection with \(q\)
new physical lines, then \(q\) of the \(|R|\) displayed powers of \(\tau\)
are owned by those lines.  No residual \(1/6\) power is delegated to a
different elementary integration.
\end{lemma}

\begin{proof}
Take $A=M\setminus R$ in \eqref{eq:source-weight-projection}.  This proves
\eqref{eq:distinguished-time-projection}.  Raise that inequality to the
power \(6/5\):
\[
 V_R\le
 B_M^{-6/5}\tau^{|R|}
 \cL_\eps^{(6/5)C^*|A|}.
\]
Multiplication by the unique weight factor \(B_M\) leaves
\(B_M^{-1/5}\le1\), which proves
\eqref{eq:distinguished-full-time}.  The last assertion uses only the
injection \(R_0\hookrightarrow R\) and the fact that these coordinates are
integrated over the displayed projection; it makes no claim that a time
eliminated by a contact-root coarea can later be restored without its
Jacobian.
\end{proof}

\subsection{Simultaneous velocity-volume control in fixed dimension}

It remains to prove the second auxiliary proposition uniformly in the fixed
boundary values.  Select the normal degree-three
components.  At a selected atom $a$, denote the fixed end by $e_a$, the free
end on the same side by $f_a$, and choose one of the remaining ends $p_a$
(the nonserial one when $a$ is an O-atom).  Trace $p_a$ backwards through the
shadow map.  If it reaches the fixed end of another selected C-atom, draw an
arrow to that atom.  Otherwise draw a free outgoing arrow terminating either
at an original free end or at a fixed end in a nonselected component.  Let
$\mathcal F=\{f_a\}$ and let $\mathcal G$ be the multiset of terminal ends.

For a C-atom the scattering rule is an orthogonal map on the pair of
velocities, and for an O-atom it is the identity.  At every selected atom,
the local energy identity implies
\[
 |\widetilde v_{f_a}|^2+
 \sum_{b\to a}|\widetilde v_{p_b}|^2
 \le \sum_{a\to c}|\widetilde v_{p_a}|^2,
\]
with a terminal arrow included on the right.  Summing over $a$ cancels every
internal arrow and gives
\begin{equation}\label{eq:velocity-arrow-contraction}
 \sum_{e\in\mathcal F}|\widetilde v_e|^2
 \le\sum_{g\in\mathcal G}|v_g|^2.
\end{equation}
Moreover the scattering formulas express the entire vector
$(\widetilde v_e)_{e\in\mathcal F}$ as a linear map
$P_\omega((v_g)_{g\in\mathcal G})$, whose coefficients depend only on the
fixed normals.  Equation \eqref{eq:velocity-arrow-contraction} says that
$\|P_\omega\|_{\ell^2\to\ell^2}\le1$.

There are at most $10|M|$ source and target vectors.  An original free end
$g$ lies in its prescribed $d$-dimensional box $|v_g|\le X_g$; a terminal
at a nonselected component lies in a $\cL_\eps^{C^*}$ box, and there are at
most $10K$ such terminals.  Hence the source domain is covered by at most
\begin{equation}\label{eq:source-unit-box-count}
 C^{|M|}\cL_\eps^{C^*K}
 \prod_{e\in E_{\rm free}(\mathcal M)}X_e^d
\end{equation}
unit coordinate boxes.  The image under a contraction of one unit box in
$\R^{d|\mathcal G|}$ is contained in a ball of radius
$C_d|\mathcal G|^{1/2}$ in $\R^{d|\mathcal F|}$.  The volume formula for this
ball and Stirling's lower bound for the gamma function cover it by
$C_d^{|\mathcal F|+|\mathcal G|}$ unit boxes.  Multiplying by
\eqref{eq:source-unit-box-count} produces a measurable containing set $Y$
with
\begin{equation}\label{eq:volume-d}
 |Y|\le\prod_{e\in E_{\rm free}(\mathcal M)}X_e^d
 C^{|M|}\cL_\eps^{C^*K}.
\end{equation}
The construction of $Y$ uses only the normals, dyadic sizes and cutting
sequence.  It is therefore uniform in the concrete fixed velocities and is
exactly \eqref{eq:source-volume} in the fixed dimension $d$.

The following table is the complete Source-B
slot/orientation correspondence used by the ordinary consumer.  ``Serial'' means
that the fixed incidence shares the connecting particle line; every O-atom
entry uses the nonserial slot required by the source definition.
\begin{center}
\scriptsize
\begin{tabular}{@{}p{0.13\textwidth}p{0.12\textwidth}p{0.25\textwidth}
 p{0.23\textwidth}p{0.17\textwidth}@{}}
\toprule
source class & selected slot and orientation & independent output chart
& rank-drop or endpoint routing & proved bound\\
\midrule
normal degree three
 & C/O admissible free slot
 & $(t,\omega,v_2)$ in \eqref{eq:degree-three-chart}
 & O serial equalities are excluded from the gain; small time, speed, or
   transported position goes to \eqref{eq:one-atom-good-section}
 & \eqref{eq:degree-three-exact}, \eqref{eq:one-atom-good-section}\\
\addlinespace
$\{33A\}$
 & $j=2$, full relative output $u$
 & direct $|u|\le4\eps/\mu$ or second-root
   $(t_2,\omega_2)$ chart
 & zero output is assigned to the direct ball
 & \eqref{eq:j2-transformed}, \eqref{eq:j2j4-full-range}\\
\addlinespace
$\{33A\}$
 & serial $j=3$, normal output $y$
 & Carleman $(t_1,y,w)$ and transformed
   $(t_1,t_2,\omega_2,w)$ charts
 & $y=0$ is direct; derivative zero has zero flux or lies in the terminal
   cell
 & \eqref{eq:j3-shell}--\eqref{eq:j3-summed}\\
\addlinespace
$\{33A\}$
 & $j=4$, tangential output $w$
 & direct $|w|\le4\eps/\mu$ or cone
   $(t_1,t_2,\omega_2,y)$ chart
 & cone axis and $y=0$ are covered by the finite direct/cone atlas
 & \eqref{eq:j4-cone-density}, \eqref{eq:j2j4-full-range}\\
\addlinespace
$\{33B\}$
 & both fixed ends serial
 & $(t_1,t_2,w,\theta)$
 & $t_1=t_2$ or either displayed relative speed small
 & \eqref{eq:typeB-first-density}, \eqref{eq:typeB-bound}\\
\addlinespace
$\{33B\}$
 & exactly one fixed end serial
 & $(t_1,t_2,\alpha,\eta^\perp,\theta)$
 & $|\zeta|$ small is routed to the time/speed alternatives
 & \eqref{eq:typeB-second-density}, \eqref{eq:typeB-bound}\\
\addlinespace
$\{33B\}$
 & neither fixed end serial
 & $(t_1,t_2,\alpha_1,\alpha_2,\omega',\eta')$
 & $D\Xi$ rank loss is excluded by
   $C_d\eps\lambda_{\rm small}^{-2}\le1/2$
 & \eqref{eq:typeB-third-density}, \eqref{eq:typeB-bound}\\
\addlinespace
$\{44\}$
 & two free ends selected as type-A fixed ends
 & $2d$-dimensional close relative state plus the type-A chart
 & complement is the ordinary nonclose branch
 & \eqref{eq:44-bound}\\
\bottomrule
\end{tabular}
\end{center}

The fixed-time powers which must dominate
$\eta_*:=\eps^{1/(9d)}$ are
\begin{table}[ht]
\centering
\begin{tabular}{@{}lll@{}}
\toprule
component & retained coefficient at the ordinary cutoffs & exponent of $\eps$\\
\midrule
$\{33{\rm A}\}$, $j=2$
 & $\eps^{d-1}\mu^{-d}$ & $d-1-1/8$\\
$\{33{\rm A}\}$, serial $j=3$
 & $\eps\mu^{-2}$ & $1-1/(4d)$\\
$\{33{\rm A}\}$, $j=4$
 & $\eps^{d-2}\mu^{-(d-1)}$ & $d-2-(d-1)/(8d)$\\
$\{33{\rm B}\}$
 & $\eps^{d-1}\lambda_{\rm small}^{-(3d-4)}$
 & $d-1-(3d-4)/(8d)$\\
$\{44\}$
 & $\eps^{-2(d-1)}\lambda_{\rm close}^{2d}$ & $7/4$\\
\bottomrule
\end{tabular}
\caption{Fixed-time ordinary exponents at the two distinct cutoffs.}
\label{tab:ordinary-exponent-check}
\end{table}
For every $d\ge4$ every exponent in the last column is strictly larger than
$1/(9d)$; logarithmic losses are therefore absorbed without using either
distinguished atom-time volume.

\begin{lemma}[Elementary retained-time estimates outside the type-A short gap]
\label{lem:ordinary-fixed-distinguished-times}
Let $\mathcal X$ be any one- or two-atom elementary component other than a
standalone degree-two root substitution, and let
$R_{\mathcal X}$ be any subset of its atom times which has been designated
as a new-line owner.  Except on a type-$\{33{\rm A}\}$ small-time cell, the variables in
$R_{\mathcal X}$ may be kept as genuine output coordinates until the global
time projection.  The local operator acts only on the complementary
variables and is bounded, uniformly in the fixed values of
$R_{\mathcal X}$, by the usual full-component baseline times
$\cL_\eps^C$.  Put $\eta=\eps^{1/(8d)}$ and
$\eta_*=\eps^{1/(9d)}$.  Moreover, if $\mathcal X$ is a two-atom support-good
$\{33{\rm B}\}$ or $\{44\}$ component, or a type-$\{33{\rm A}\}$ component
on $|t_1-t_2|\ge\eta$, its retained good factor $\eta_*$ is still available
with \emph{both} atom times fixed.  On a type-$\{33{\rm A}\}$ cell with
$|t_1-t_2|<\eta$ no pointwise fixed-time estimate is asserted: the original
two-time Radon block and its indicator are retained until the post-envelope
common-measure estimate in Lemma~\ref{lem:owner-marked-good-component}.
\end{lemma}

\begin{proof}
The standalone degree-two root is excluded from the statement; no time is
recovered from that norm-one substitution in this local lemma.  In the
one-fixed-line degree-three chart \eqref{eq:degree-three-chart}, fix the atom
time.  The angular and free-velocity coordinates remain in a cutoff box
independent of all boundary states;
the normal degree-three simultaneous-volume estimate is therefore uniform in
that time.  The same observation applies to a normalized one-atom
degree-four chart.  A one-atom component whose good alternative itself is a
small-time cell is not claimed to keep an independent good factor here; the
present statement retains only its ordinary baseline when that time is
fixed.

For $\{33{\rm A}\}$ on $|t_1-t_2|\ge\eta$,
Lemma~\ref{lem:fixed-time-33A-j2-j3} gives, with both times fixed and for an
arbitrary nonnegative kernel, the bounds
\begin{equation}\label{eq:j2j3-fixed-distinguished-times}
 \mathcal J_{33A,2}^{\,t_1,t_2}(Q)
 \le\cL_\eps^C\eps^{d-1}\mu^{-d}\int Q^\sharp,
 \qquad
 \mathcal J_{33A,3}^{\,t_1,t_2}(Q)
 \le\cL_\eps^C\eps\mu^{-2}\int Q^\sharp.
\end{equation}
These estimates use neither time volume.  The all-range latitude coarea
estimate \eqref{eq:j4-fixed-time-sphere}, including the direct small-$w$
cells, gives
\begin{equation}\label{eq:j4-fixed-distinguished-times}
 \mathcal J_{33A,4}^{\,t_1,t_2}(Q)
 \le\cL_\eps^C\eps^{d-2}\mu^{-(d-1)}\int Q^\sharp.
\end{equation}
At the ordinary value $\mu=\eta=\eps^{1/(8d)}$, the five exponents are exactly
those in Table~\ref{tab:ordinary-exponent-check}.  Hence every orientation
retains the final $\eta_*$ credit after logarithmic absorption.  If neither
time is distinguished one may of course use the
stronger affine-log form.  If only one time is distinguished, first apply
the selected fixed-both-times kernel estimate pointwise in the other time
and then use conditional Tonelli in that unmarked time.  The same displayed
density coefficient is retained; the actual half-open interval of the
unmarked time is merely appended to the fixed output box and is consumed
locally.  If both times are distinguished, neither time is present in the
local output box.  Thus the fixed-time argument never obtains its
$\eps$-gain from the length of a time interval that is later reused by the
global projection.

The complementary type-$\{33{\rm A}\}$ small-time cell is not covered by
these fixed-time inequalities.  It is kept on the original common two-time
measure, with its half-open indicator, for the post-envelope clause of
Lemma~\ref{lem:owner-marked-good-component}; no conditional representative
is bounded pointwise there.

For $\{33{\rm B}\}$, fix $t_1,t_2$ in the three exact densities
\eqref{eq:typeB-first-density}, \eqref{eq:typeB-second-density} and
\eqref{eq:typeB-third-density}.  All remaining inverse pivots are bounded by
$\lambda_{\rm small}^{-(3d-4)}$ on \eqref{eq:typeB-nondegeneracy}, and the half-open charts
in \eqref{eq:typeB-fixed-output-domains} integrate only the complementary
velocity, angular and scalar variables.  Thus
\begin{equation}\label{eq:typeB-fixed-distinguished-times}
 \mathcal J_{33B}^{\,t_1,t_2}(Q)
 \le\cL_\eps^C\eps^{d-1}\lambda_{\rm small}^{-(3d-4)}\int Q^\sharp.
\end{equation}
With $\lambda_{\rm small}=\eps^{1/(8d)}$ this is again smaller than $\eta_*$.  Finally,
the $\{44\}$ gain \eqref{eq:44-bound} is the $2d$-dimensional volume of the
close relative endpoint state and does not use either atom-time volume; it
 is unchanged when both times are fixed.  The complementary small-relative-
speed cells for these component types are kept half-open and are bounded by
the same fixed-time coarse estimates, without manufacturing a good credit.
Together with the explicit type-$\{33{\rm A}\}$ exception above, hence
every distinguished time remains an actual output of the local operator,
unless it is consumed later on that common two-time measure, and all asserted
uniform and good bounds follow on their stated regimes.
\end{proof}

\begin{proposition}[Ordinary estimate in every fixed dimension $d\ge4$]
\label{prop:ordinary-subset}
The degree-two, normal degree-three and degree-four, all ordinary good
$\{33{\rm A}\}$, all three $\{33{\rm B}\}$, $\{44\}$, one-pair multiple
contact, weight and simultaneous velocity-volume estimates required in
Source B, Propositions 9.1--9.4, hold in every fixed dimension $d\ge4$ as positive
operators on fixed enlarged domains.  The serial $j=3$ orientation is
bounded by \eqref{eq:j3-summed}; no stronger isolated estimate is used by
the ordinary consumer.
\end{proposition}

\begin{proof}
The one-atom and multiple-contact assertions are proved in the first two
subsections.  Equations \eqref{eq:j3-shell}--\eqref{eq:j3-summed} cover
the serial orientation, while \eqref{eq:j2j4-full-range},
\eqref{eq:typeB-bound} and \eqref{eq:44-bound} cover every remaining
two-atom orientation.  Equations \eqref{eq:source-weight-1}--
\eqref{eq:source-weight-2} and \eqref{eq:volume-d} prove the auxiliary
claims.  All decompositions are positive and all changes of variables
were performed before enlarging their domains, so the conclusion holds
for arbitrary measurable $Q\ge0$.  For the support-good alternatives this
means the indicator-preserving measure form
\eqref{eq:one-atom-good-section} and
\eqref{eq:source-elementary-operator}; their scalar smallness is used only
after the pointwise source envelope, as recorded in
Lemma~\ref{lem:owner-marked-good-component}.
\end{proof}

This proposition proves every local analytic assertion in (R3).  The global
small-window time ledger is not inferred from an arbitrary omitted set in
\eqref{eq:source-weight-projection}: Proposition
\ref{prop:packet-free-small-window-consumer} first exhausts the actual
new-line birth owners and then applies the present fixed-distinguished-time
estimates.  Together the two propositions prove (R3) of
Proposition~\ref{prop:reduction}.
\section{Periodic double-overlap repair}
\label{sec:parallel-ov}

Four steps of the Euclidean cutting argument use the absence of a double
overlap segment: the two conclusions of the UP proposition, the
cyclic-lower-component estimate, and the 2CONNUP estimate.  The final DOWN tail
requires the same replacement in the periodic setting.  Keeping the original
cutting order, we assign each periodic object canonically at the first failed
no-double-overlap inference and charge it at most once.

\subsection{Lift separation and permanent incidence labels}

The backward map $\operatorname{Ba}$ traces an active edge to its edge before
the current cut, and the shadow map $\operatorname{Sh}$ identifies a fixed
edge created by a cut with the free edge having the same transported state;
see Source B, Definition~8.6 and Propositions~8.7--8.8.  Source B, Proposition~8.10,
then supplies the positive-kernel cutting identity used below.  Iterating these
maps attaches to every current half-edge a permanent label
\begin{equation}\label{eq:original-incidence-label}
 \mathfrak i=(a,\mathfrak p,\varsigma),
\end{equation}
where $a$ is an atom of the molecule before cutting, $\mathfrak p$ is one of
the two physical particle lines through $a$, and $\varsigma$ is its top or
bottom side.  We call \eqref{eq:original-incidence-label} an
\emph{original incidence}.  A cut changes its current representative but not
its original incidence or its physical-line label.

\begin{lemma}[One incoming root for a fixed lift]
\label{lem:one-incoming-root}
Let $x,g\in\R^d$, $g\ne0$, and $m\in\Z^d$.  Away from tangency, the equation
\begin{equation}\label{eq:fixed-lift-contact}
 |x+tg-m|=\eps
\end{equation}
has at most one root satisfying the incoming condition
$(x+tg-m)\cdot g<0$.  Hence two distinct incoming roots of one freely
transported relative state have different torus lifts.
\end{lemma}

\begin{proof}
For $\phi_m(t)=|x+tg-m|^2-\eps^2$, strict convexity gives, at two
nongrazing roots $t_m^-<t_m^+$,
\[
 \phi_m'(t_m^-)<0<\phi_m'(t_m^+),
 \qquad \phi_m'(t)=2(x+tg-m)\cdot g.
\]
Thus only $t_m^-$ is incoming.  A double root is tangential and belongs to
the collision-flux null set already removed in
Proposition~\ref{prop:ae-flow}.
\end{proof}

\begin{definition}[Double support and protected certificate]
\label{def:protected-certificate}
A \emph{double support} is a pair of atoms $p,n$ joined by maximal overlap
segments on two distinct physical particle lines.  At the first cutting step
at which such a support invalidates a Source-B no-double-overlap inference,
we attach one immutable certificate identifier to the resulting elementary
object.  Its representative is one of
\begin{enumerate}[label=\textup{(\roman*)}]
\item a type-$\{33\mathrm A\}$ adjacent pair;
\item a good $\{3\}$ singleton;
\item a good $\{4\}$ singleton.
\end{enumerate}
The certificate stores the double support, its creation step, and the
original incidences of the representative.  The object is then protected:
its atoms and all their active representatives are removed before the
algorithm restarts.
\end{definition}

For a pair representative we retain six labelled bearing ports, one for each
of its at most six active external incidences; unused ports are filled with a
cemetery symbol.  A singleton has at most four such ports and is padded to
six.  The ordered port list is denoted by
\begin{equation}\label{eq:certificate-ports}
 \operatorname{Port}(C)=(\mathfrak i_1(C),\ldots,
 \mathfrak i_6(C)).
\end{equation}
Only original labels occur in this list, so it is unchanged by later cuts.

For uniform accounting, an ordinary $\{33\mathrm A\}$ output is also given
an immutable identifier and the same six-port record at its creation.  A
\emph{bearing record} means either such an ordinary pair record or a protected
periodic certificate.  Only the latter represents a failed periodic
inference, but all bounded-fibre maps below take values in the disjoint union
of the two kinds of bearing records.

\begin{lemma}[First-failure classification]
\label{lem:first-failure}
Run UP, DOWN, or 2CONNUP and stop at the first step where one of the
no-double-overlap inferences used in Source B fails.  If the step has already
been assigned to the Source-B exceptional ledger---a step-1 free cut or the
degree-four branch of step~2 in 2CONNUP---make no new combinatorial charge.
Otherwise the failed cut produces exactly one protected certificate in the
sense of Definition~\ref{def:protected-certificate}.  The conclusion is
invariant under TOP/BOTTOM reversal.
\end{lemma}

\begin{proof}
At the first failure, the two active paths whose seriality was identified in
the Euclidean proof have backward traces on two distinct original particle
lines.  Both paths have the cut atom $p$ as one endpoint and the atom $n$ at
which the competing serial paths meet as the other endpoint.  They therefore
form a double support.  This is precisely condition $(\triangle)$ in Source
A, Section~5, Part~4.

If the two cut-side terminals remain distinct, the adjacent legal pair cut
has one same-side external fixed incidence at each atom.  For UP these are
both top incidences and for DOWN both are bottom incidences.  The component is
therefore of type $\{33\mathrm A\}$, not type $\{33\mathrm B\}$.  If the two
terminals coincide, the cutting rule gives either the same pair or a
singleton whose two same-side edges are the two lines of the double support.
The first alternative gives (i).  We prove that the singleton is good.

The two original particle lines in the double support carry one common
freely transported relative state between the endpoint atoms.  Depending on
the C/O type and the side exposed by the cut, this is either the incoming
pair or the reflected outgoing pair at the representative singleton.  The
two endpoint atoms give two distinct contact times
\begin{equation}\label{eq:parallel-contacts}
 |x+t_pg-m_p|=\eps,
 \qquad |x+t_ng-m_n|=\eps,
 \qquad t_n<t_p.
\end{equation}
If $|g|\le\eps^{1/(8d)}$, the ordinary small-relative-speed alternative is
already a good singleton.  Suppose otherwise.  At a contact $t$ the other
root for the same lift, when it exists, is
\begin{equation}\label{eq:companion-root}
 t'=t-2\frac{(x+tg-m)\cdot g}{|g|^2}.
\end{equation}
The higher witness at $p$ is incoming.  If the lower endpoint $n$ is an
O-atom, velocity transparency makes $t_n$ incoming for the same selected
state.  If $n$ is a C-atom, that selected state is outgoing at $t_n$;
\eqref{eq:companion-root} supplies its incoming companion
 $s_n<t_n<t_p$ with the same lift as $n$.  Thus in both cases we have two
 distinct incoming roots, $t_p$ and $s_n$ (where $s_n=t_n$ in the O-case).
 In the C-case,
 $0<t_n-s_n\le2\eps/|g|\le2\eps^{1-1/(8d)}$, so the companion remains in the
 fixed enlargement of the original time cell and introduces no new
 unbounded lift or time sum.
 Lemma~\ref{lem:one-incoming-root} forces their lifts to differ.  The DOWN
case is the exact time reverse.  The sign-symmetric form of
\eqref{eq:multiple-contact-tube} now gives, on every lift and velocity cell,
a relative $\cL_\eps^C\eps^{d-1}$ section.  This is stronger than the
$\eps^{1/(8d)}$ excess required of a good singleton, and proves (ii) or (iii)
according to its active degree.  This is exactly the extra support in Source
A, Proposition~3.2(4), now with the two lifts and the retained side stated
explicitly.  Protection removes the source atom, so no later restart can
create a second identifier from the same failure.  The proof after reversing
top and bottom is identical.
\end{proof}

\begin{lemma}[Finite-state exhaustion of the four first-failure calls]
\label{lem:first-failure-finite-state}
The four calls in Source A, Section~5, Part~4, at which absence of a double
overlap is used are denoted by
\[
 \mathfrak c\in
 \{\mathsf U_1,\mathsf U_2,\mathsf{CY},\mathsf{2C}\}.
\]
They are, respectively, Source B, Proposition~10.2(1),
Proposition~10.2(2), Proposition~10.9, and Proposition~11.2.  At the first
failed inference record the finite state
\begin{equation}\label{eq:first-failure-state-tuple}
 \mathfrak s=(\mathfrak c,o,\varsigma,\delta,\chi,\theta).
\end{equation}
Their correspondence with the four Source-A calls is recorded in
Online Resource~1, Section~\ref{appB:source-location-notes}.
Here \(o\in\{\mathrm{old},\mathrm{new}\}\) is the priority-ledger flag,
\(\varsigma\in\{\mathrm{top},\mathrm{bottom}\}\) is the exposed side,
\(\delta\) is the active-degree signature of the elementary output,
\(\chi\in\{\mathrm{sep},\mathrm{coin}\}\) records whether the two cut-side
terminals are distinct or coincident, and
\(\theta\in\{\mathrm C,\mathrm O,*\}\) is the C/O type of the lower endpoint
(\(*\) when no root conversion is used).  Every realizable state is in the
following call-site list.  The list is an upper list: a row does not assert
that every syntactically listed new state occurs, only that no other state
can occur and that every listed state has the unique disposition in the
second table.

\begin{center}
\small
\begin{tabular}{@{}p{0.13\textwidth}p{0.25\textwidth}p{0.31\textwidth}
 p{0.23\textwidth}@{}}
\toprule
call & old states & new degree/terminal states & excluded state\\
\midrule
\(\mathsf U_1\)
& step--1 free-cut charge
& \(\mathcal R\)
& no other elementary degree signature\\
\(\mathsf U_2\)
& step--1 free-cut charge
& \(\mathcal R\)
& no other elementary degree signature\\
\(\mathsf{CY}\)
& step--1 free-cut charge
& \(\mathcal R\)
& no other elementary degree signature\\
\(\mathsf{2C}\)
& step--1 free cut or the degree--four branch of step~2
& \(\mathcal R\setminus\{(4,\mathrm{coin})\}\)
& a non-old degree--four output\\
\bottomrule
\end{tabular}
\end{center}
where
\begin{equation}\label{eq:first-failure-new-state-set}
 \mathcal R:=
 \{((3,3),\mathrm{sep}),((3,3),\mathrm{coin}),
        (3,\mathrm{coin}),(4,\mathrm{coin})\}.
\end{equation}
The entry \(\varsigma\) records the side.  For a singleton, \(\theta\)
records its lower type.  These entries complete the state, and the unique
ledger/output map is
\begin{center}
\small
\begin{tabular}{@{}p{0.13\textwidth}p{0.14\textwidth}p{0.14\textwidth}
 p{0.14\textwidth}p{0.35\textwidth}@{}}
\toprule
priority & \(\delta\) & \(\chi\) & \(\theta\) & unique disposition\\
\midrule
old & any old signature & any & any
& \(\mathcal L_{\rm old}\); no new certificate and no new charge\\
new & \((3,3)\) & sep or coin & \(*\)
& one same-\(\varsigma\)-side \(\{33\mathrm A\}\) record in
  \(\mathcal L_{33}\)\\
new & \(3\) & coin & O
& one protected good \(\{3\}\), using the two incoming roots\\
new & \(3\) & coin & C
& one protected good \(\{3\}\), after the incoming-companion substitution\\
new & \(4\) & coin & O
& one protected good \(\{4\}\), with its full baseline retained\\
new & \(4\) & coin & C
& one protected good \(\{4\}\), after the incoming companion and with its
  full baseline retained\\
\bottomrule
\end{tabular}
\end{center}
The last two rows are absent for the non-old \(\mathsf{2C}\) call.  TOP/BOTTOM
reversal changes only \(\varsigma\) and the time orientation; it does not
create another state.
\end{lemma}

\begin{proof}
At the four call sites, the source proofs use no-double-overlap in exactly
the following statements: the first degree drop in
\(\mathsf U_1\), the serial fixed-end test in \(\mathsf U_2\), the first
degree-two atom of a cyclic lower component in \(\mathsf{CY}\), and the
simultaneous break of two bearing overlap segments in \(\mathsf{2C}\).
Tracing the two active paths backwards at the first failure gives the same
two labelled support lines in all four cases, as in the first paragraph of
the proof of Lemma~\ref{lem:first-failure}.  Thus the call label changes no
local state: condition \((\triangle)\) holds and the only remaining data are
the last five entries of \eqref{eq:first-failure-state-tuple}.

Apply the Source-B priority rule before classifying the local output.  A
step--1 free cut already owns its degree-two charge.  In \(\mathsf{2C}\), the
degree-four branch of step~2 is likewise an old exception.  These and only
these states give the first row of the disposition table.  In particular,
a degree-four \(\mathsf{2C}\) state cannot enter a new periodic ledger.

Consider a non-old state.  The legal elementary-output list leaves only two
possibilities.  If two active atoms remain, each has degree three and each
has one external fixed incidence on the exposed side.  The signature is
\((3,3)\), and the output is \(\{33\mathrm A\}\), whether the two cut-side
terminals were kept distinct or identified by the local cut.  Opposite-side
fixed incidences, and hence \(\{33\mathrm B\}\), are impossible.  If one
active atom remains, the other atom is the norm-one root substitution and
the active singleton has degree three or four.  A singleton cannot have
\(\chi=\mathrm{sep}\), so these are precisely the last four elements of
\(\mathcal R\) after the C/O bit is supplied.  No \(\{44\}\) output belongs
to the elementary list at any of the four failed cuts.

For an O lower endpoint, transparency makes the two contact witnesses
incoming for one freely transported relative state.  For a C lower endpoint,
\eqref{eq:companion-root} replaces its outgoing witness by the incoming
companion.  Lemma~\ref{lem:one-incoming-root} then forces different lifts,
and \eqref{eq:multiple-contact-tube} gives the good singleton.  Its active
degree decides uniquely between \(\{3\}\) and \(\{4\}\); in the latter case
the good bit does not remove the unique full baseline.  Finally, the
priority classes are disjoint and protection removes the chosen active
representative.  Hence no state has two dispositions and no realizable
state is missing from the two tables.
\end{proof}

\begin{lemma}[Post-envelope measure of a periodic certificate]
\label{lem:periodic-certificate-measure}
Let $C$ be a protected periodic first-failure certificate in the sense of
Definition~\ref{def:protected-certificate} whose current analytic
representative is counted by $N_{\rm good}$.  Retain its
complete different-lift indicator through every arbitrary-kernel operation.
After the source weight and Maxwellian bounds have produced a pointwise
conditional majorant, the certificate's remaining relative-velocity
section has measure at most
\begin{equation}\label{eq:periodic-certificate-section}
 \cL_\eps^C\eps^{d-1}.
\end{equation}
This section estimate is uniform in every fixed time and boundary variable.
If the representative also carries a new-line owner, that time is not used
in \eqref{eq:periodic-certificate-section} and remains available for the
ordinary distinguished-time projection.  If the representative is a full
degree-four singleton, its normalization $\eps^{-(d-1)}$ is retained outside
\eqref{eq:periodic-certificate-section}; the tube is its good-support excess,
not a deletion of that baseline.
\end{lemma}

\begin{proof}
The two original lines in the stored double support give two contacts of one
freely transported relative state.  After the incoming companion conversion
in the C-case their lifts are different by
Lemma~\ref{lem:one-incoming-root}.  For every fixed pair of contact times and
all remaining boundary variables, \eqref{eq:lattice-tube-wedge} confines the
relative velocity to a tube with $d-1$ transverse radii
$O(\cL_\eps^C\eps)$ and a polylogarithmically bounded axial coordinate.
This proves \eqref{eq:periodic-certificate-section} for a pair as well as for
either singleton representative.  The calculation is a conditional measure
bound; it is not pulled through an arbitrary test kernel.  The owner time is
only a fixed parameter in this section, so Fubini leaves its original layer
interval untouched.  Finally the immutable identifier prevents the same
 tube from being assigned to two current representatives.
\end{proof}

\begin{theorem}[Exact full-pair velocity disintegration]
\label{thm:exact-full-pair-disintegration}
For the velocity assertions below, fix a finite legal Source-B cutting
branch, all atom times and normals, and all discrete cells.  On the
right-hand side of the exact cutting identity
\eqref{eq:source-cutting-identity}, first resolve the velocity-conservation
deltas at every C- and O-atom.  Let $\mathcal R$ be the physical root lines
and $\mathcal P\setminus\mathcal R$ the lines created by the actual birth
forest, write $n:=|\mathcal P|$, and put
\begin{equation}\label{eq:exact-global-input-fibre}
 \mathcal U
 :=\bigoplus_{p\in\mathcal R}\mathbb R^d_{v_p}
   \oplus
   \bigoplus_{q\in\mathcal P\setminus\mathcal R}\mathbb R^d_{h_q}
 \simeq\mathbb R^{dn}.
\end{equation}
Thus $u\in\mathcal U$ consists of one root velocity for each root line and
one independent birth velocity for each nonroot line.  Every bond/free-end velocity on the cut
molecule is an affine-linear function of $u$.  A fixed shadow is the same
function as its unique free shadow and is not another coordinate.  In
particular, distinct shadow classes at different trajectory segments are
not incorrectly counted as independent when a collision delta relates
them.  Suppose that the live local consumers are
listed in reverse cutting-dependency order and that, whenever the two
outputs of one C-junction have different eventual consumers, the two
outputs remain in one joint group until after the complete pair junction.
Then there is a finite Borel partition \((E_\lambda)_\lambda\) of
$\mathcal U$ and, on each piece, an affine bijection
\begin{equation}\label{eq:exact-global-velocity-map}
 \mathcal G_\lambda:u\longmapsto
 (U_\alpha)_{\alpha\in\mathscr A}\oplus U_{\rm rem}
 =:w
\end{equation}
with the following properties.
\begin{enumerate}[label=\textup{(\roman*)}]
\item The output blocks form a partition of the complete linear velocity
fibre: every input column and every output row occurs once.  A block
$U_\alpha$ contains only the \emph{pre-coarea linear velocity inputs} of
consumer $\alpha$.  Angular variables, centre variables, root times and
other nonlinear local coordinates are not rows of
$D\mathcal G_\lambda$.
\item At a C-atom on lines $p,q$ the event factor is the complete matrix
\begin{equation}\label{eq:exact-global-C-factor}
 \binom{\delta v_p^+}{\delta v_q^+}
 =\mathscr R_{pq}^{(2d)}
  \binom{\delta v_p^-}{\delta v_q^-},
 \qquad |\det\mathscr R_{pq}^{(2d)}|=1.
\end{equation}
At a birth of a child $q$, the coordinate $h_q$ already present in
\eqref{eq:exact-global-input-fibre} occupies the child input slot, and the
same square formula is used on $(\delta v_p^-,h_q)$.  An O-atom is the
identity.  No projected block $P_\omega$ or $I-P_\omega$ is inverted.
\item For the chronologically ordered event word, let
$\widehat{\mathscr R}_{j,\lambda}\in O(dn)$ be the identity on every line
except for the complete $2d$-dimensional pair at the $j$th C-event, where
it equals \eqref{eq:exact-global-C-factor}; it is the identity at an
O-event.  There are square block-triangular cutting/forest basis maps
$\mathfrak S_\lambda,\mathfrak T_\lambda\in GL(dn)$ such that
\begin{equation}\label{eq:exact-global-factorization}
 D\mathcal G_\lambda
 =\mathfrak T_\lambda
  \widehat{\mathscr R}_{r,\lambda}\cdots
  \widehat{\mathscr R}_{1,\lambda}
  \mathfrak S_\lambda,
\end{equation}
and
\begin{equation}\label{eq:exact-global-factor-bounds}
 C^{-n}\le|\det D\mathcal G_\lambda|\le C^n,
 \qquad
 \|D\mathcal G_\lambda\|+\|(D\mathcal G_\lambda)^{-1}\|
 \le C^n.
\end{equation}
\item For every nonnegative Borel $F$, including an arbitrary complete
complement kernel and the original time indicator, one has the exact
push-forward identity
\begin{multline}\label{eq:exact-global-pushforward}
 \int_{\mathcal U}F(u)\,\dd u\\
 =\sum_\lambda\int_{\mathcal G_\lambda(E_\lambda)}
 F(\mathcal G_\lambda^{-1}w)
 |\det D\mathcal G_\lambda|^{-1}
 \prod_{\alpha\in\mathscr A}\dd U_\alpha\,\dd U_{\rm rem}.
\end{multline}
The identity remains true with every fixed-shadow substitution and with
$\mathbf1_{\mathcal D_M}$ inside $F$; neither is enlarged or averaged.
\item Let $Z_\alpha$ be the component-owned nonvelocity variables of a local
consumer after every distinguished time has been held outside.  These
variables are disjoint between components in the Source-B iterated
integral and are not rows of the global velocity map.  A local nonlinear
chart may subsequently be applied to $(U_\alpha,Z_\alpha)$ while all other
blocks are fixed.  Its equidimensional full output
$W_\alpha=(w_\alpha,\zeta_\alpha)$ is a local coarea output and is
\emph{not} identified with a row group of
\eqref{eq:exact-global-velocity-map}.  Conditional area formula followed
by Tonelli therefore commutes with the global linear push-forward without
duplicating a velocity coordinate.
\end{enumerate}

Separately, before conditioning the marked atom time, the collision time at
a birth $b(q)$ can be retained through a standalone degree-two cleanup call.
More precisely, write
$z_{\rm loc}=(z_2,z_3,z_4)$ for the three nonfixed edge states of the
original one-fixed-line atom (with $z_1=z_p$ held as a parameter).  The
time below is the original indexed variable $t=t_{b(q)}$.  The original
normalized C-atom Radon measure has the one-fixed-line
disintegration
\begin{multline}\label{eq:exact-transparent-birth-relay}
 \eps^{-(d-1)}\int \Delta_{b(q)}(z,t)F(z,t)\,
       \dd z_{\rm loc}\,\dd t\\
 =\int_{\mathbb R_{t_{b(q)}}\times\Sph^{d-1}\times\mathbb R^d}
 [(v_p-h_q)\cdot\omega]_-\,
 F(\Psi_{b(q)}(t,\omega,h_q),t)
 \dd t\,\dd\omega\,\dd h_q .
\end{multline}
Here $\Psi_{b(q)}$ is the exact incoming one-fixed-line collision chart,
including the fixed shadow and lift.  Formula
\eqref{eq:exact-transparent-birth-relay} holds branchwise with
$\mathbf1_{\mathcal D_M}$ and the complete complement inside $F$.
Consequently the marked time and $h_q$ are genuine Lebesgue output
coordinates; the degree-two call is a relay, not their owner.  This is used
inside the unchanged global normalization of
\eqref{eq:source-cutting-identity}: it does not change
$|E_*|-2|M|$.  In particular, the unique additional $\eps^{-(d-1)}$ baseline
of a no-fixed-end full-$\{4\}$ component is separate from the contact
normalization displayed on the left and remains outside the relay.
\end{theorem}

\begin{proof}
We construct the map directly from the global velocity fibre.  Put every
coordinate in \eqref{eq:exact-global-input-fibre} into the fixed space
$\mathcal U\simeq\mathbb R^{dn}$ at the outset.  Before the creation of $q$,
the coordinate $h_q$ is carried unchanged by the identity until its birth.  At its birth, apply the embedded $dn$-dimensional matrix
$\widehat{\mathscr R}_{b(q),\lambda}$, whose only nonidentity block is
\eqref{eq:exact-global-C-factor} on the complete parent--child pair
$(v_p^-,h_q)$.  At every later C-event apply the analogous embedded matrix,
and at an O-event apply the identity.  Thus repeated collisions involving
the same physical line are represented by an \emph{ordered product} of
square matrices on the same $dn$-dimensional fibre, never by a direct sum of
event blocks.  Induction on the event word gives the middle factor of
\eqref{eq:exact-global-factorization} and a square orthogonal map at every
frontier.

Now perform the legal cuts.  A cut does not adjoin a velocity: it selects
affine evaluation functions of $u$ as free representatives, while every
fixed shadow is the same evaluation function used as a parameter.  In
reverse dependency order the elementary square basis changes are only:
(a) a permutation or renaming of $d$-dimensional coordinates; (b) a
complete pair-scattering block; (c) a centre/relative transformation of a
complete pair, with fixed nonzero determinant; and (d) on a reduced rooted
forest, the root together with the child-minus-parent differences.  After a
topological ordering, the matrix in (d) is block triangular with diagonal
blocks $\pm I_d$.  These operations construct
$\mathfrak S_\lambda,\mathfrak T_\lambda\in GL(dn)$.  Affine lift and cell
translations do not enter their derivatives, and the path formulas for the
inverse have coefficients bounded by $C^n$.  This proves
\eqref{eq:exact-global-factorization}--
\eqref{eq:exact-global-factor-bounds}.

Row exhaustion follows from the exact velocity-fibre recursion.  List the consumers as $\alpha_1,\ldots,\alpha_s$ in reverse
dependency order and let $V_{k-1}$ be the unconsumed velocity fibre just
before consumer $\alpha_k$ is exposed.  The elementary square basis changes
above give, at each step, the exact recursion
\begin{equation}\label{eq:exact-global-consumer-recursion}
 V_{k-1}\simeq U_{\alpha_k}\oplus V_k,
 \qquad
 \dim V_{k-1}=\dim U_{\alpha_k}+\dim V_k.
\end{equation}
If the boundary met at that step is a C-junction, inverse complete-pair
scattering first recovers \emph{both} incoming pre-states; the two output
slots remain one joint group until they are assigned to their possibly
different downstream consumers.  A fixed shadow is then only evaluation at
a parameter and removes no coordinate.  A marked degree-two relay likewise
has no output row of its own.  Iterating
\eqref{eq:exact-global-consumer-recursion} leaves $V_s=U_{\rm rem}$ and gives
\[
 \mathcal U\simeq
 \bigoplus_{k=1}^sU_{\alpha_k}\oplus U_{\rm rem},
 \qquad
 \sum_{k=1}^s\dim U_{\alpha_k}+\dim U_{\rm rem}=dn.
\]
Hence every row and every column occurs exactly once.  This is the required
ownership bijection and proves (i).

On each half-open Borel cell the affine change-of-variables theorem gives
\eqref{eq:exact-global-pushforward}; summing the disjoint cells gives the
full identity.  Since $F$ is arbitrary and nonnegative, fixed-shadow and
time-order indicators pass through pointwise.  Applying a local area
formula only after the corresponding linear block $U_\alpha$ has been
isolated, and adjoining only its separately owned nonvelocity block
$Z_\alpha$, proves (v); all partner blocks are parameters in that
equidimensional conditional formula.

Finally, integration of the contact delta in the free position gives
$\eps^{d-1}\dd\omega$ by \eqref{eq:contact-sphere-basic}.  The incoming flux
$[(v_p-h_q)\cdot\omega]_-$ is the explicit density already present in the
C-kernel \eqref{eq:source-C-kernel}; it is not obtained by differentiating
or integrating out the atom time.  The transport and scattering deltas then
resolve the remaining two edge states with unit determinant.  Thus the
$\eps^{-(d-1)}$ contact normalization cancels exactly, and the one-fixed-line
incoming chart is one-to-one on each sign/lift cell.  The area formula gives
\eqref{eq:exact-transparent-birth-relay}.  Lift, sign and time-cell
boundaries are assigned by the same half-open convention as the cutting
identity.  The grazing set and the zero-relative-velocity set carry zero
C-kernel Radon mass because of the explicit incoming-flux factor, so no
exceptional endpoint measure is added.  Keeping
$\mathbf1_{\mathcal D_M}F$ as the test function
 retains the original atom-time index pointwise and proves the last assertion
 without a separate time marginalization.  Since this is merely the
 one-fixed-line disintegration within the same cutting identity, it neither
 creates nor removes the extra full-component baseline.
\end{proof}

\begin{lemma}[Local free-coordinate injection for periodic certificates]
\label{lem:periodic-certificate-free-coordinate}
Let $\mathscr C$ be any family of protected periodic certificates produced
by one legal cutting sequence.  In reverse cutting order one can assign to
each $C\in\mathscr C$ its complete local non-distinguished Source-B block
$W_C=(w_C,\zeta_C)$.  Here $U_C$ denotes its pre-coarea linear velocity
input and $Z_C$ its separately owned nonvelocity input (nondistinguished
angular, spatial-centre, time and auxiliary variables).  The local chart is
the equidimensional map
\begin{equation}\label{eq:certificate-local-equidimensional-chart}
 (U_C,Z_C)\longmapsto W_C=(w_C,\zeta_C),
 \qquad \dim W_C=\dim U_C+\dim Z_C,
\end{equation}
where $w_C$ is a $d$-dimensional local free coordinate obtained from
$U_C$ and $\zeta_C$ comprises all remaining local output coordinates, such
that:
\begin{enumerate}[label=\textup{(\roman*)}]
\item the complete different-lift support indicator of $C$ is bounded by
the tube indicator in \eqref{eq:lattice-tube-wedge} for an affine image of $w_C$ with determinant
bounded above and below by absolute constants;
\item after all earlier components are fixed, $U_C$ is the complete linear
velocity input of the current elementary component: it contains every
genuinely integrated free velocity of that component and no fixed boundary
shadow.  The block $Z_C$ is its
component-owned nonvelocity input after every distinguished birth-owner time
has been held outside.  The pairs $(U_C,Z_C)$ are disjoint integration
blocks for distinct certificates.  The corresponding local coarea output is
$W_C$.  The complementary domain of $\zeta_C$ is a boundary-independent
$\cL_\eps^C$ box (or the norm-one
$\{2\}$ substitution in a temporary $\{33{\rm A}\}$ refinement).  Because
the current representative is classified as good, the whole block is
excluded from the Source-B normal degree-three family used in
Proposition~\ref{prop:ordinary-subset};
\item the map $C\mapsto U_C$ is injective at the level of global linear
velocity blocks by Theorem~\ref{thm:exact-full-pair-disintegration}; local
area formula then converts $(U_C,Z_C)$ to $W_C$ with all other blocks fixed.
Therefore, after a pointwise source majorant has been obtained, the
conditional complete-block measures multiply:
\begin{equation}\label{eq:periodic-certificate-product-measure}
 \int\prod_{C\in\mathscr C}\mathbf1_{\rm tube,C}
      \prod_{C\in\mathscr C}\dd w_C\dd\zeta_C
 \le \cL_\eps^{C|\mathscr C|}\eps^{(d-1)|\mathscr C|}.
\end{equation}
\end{enumerate}
All distinguished $S$-times belonging to $C$, not only its selected
bookkeeping time, and all degree-four normalizations are outside this
coordinate product.  The tube estimate is uniform in those fixed times, so
they remain genuine outputs for the global time projection.
\end{lemma}

\begin{proof}
For a stored first-failure record use the two labelled support lines, not an
arbitrary port.  For a $\{3\}$ singleton choose the nonserial free velocity in
the degree-three chart; subtracting its fixed shadow gives the stored relative
velocity.  For a $\{33{\rm A}\}$ pair first use the legal exact refinement
into its one-fixed-line $\{3\}$ member and its norm-one $\{2\}$ root
substitution, and choose the same free relative velocity on the $\{3\}$
member.  For a full $\{4\}$ singleton pass from one free velocity pair to
centre-of-mass and relative coordinates and choose the relative coordinate.
The latter change has constant nonzero determinant, and its separate
$\eps^{-(d-1)}$ source normalization is untouched.  In the O-case the relative
state is unchanged.  In the C-case first perform the incoming companion
substitution \eqref{eq:companion-root}; the scattering map is orthogonal, so
the stored relative velocity is an affine orthogonal image of the selected
local coordinate.  In each case extend this selected coordinate to a
complete square linear chart
\[
 U_C\longmapsto (w_C,r_C)
\]
by retaining the unused velocity coordinates (for a full pair, $r_C$
contains the centre coordinate).  The determinant and inverse of this chart
are bounded by absolute constants.  Taking
$\zeta_C=(r_C,Z_C)$, or the corresponding bounded-Jacobian Source-B local
coordinates, now gives the equidimensional chart
\eqref{eq:certificate-local-equidimensional-chart}.
Equation~\eqref{eq:lattice-tube-wedge} then gives the support domination in
(i) for every representative class and supplies the full-dimensional Jacobian used
below.

The positive cutting identity orders the component carrying a free shadow
before the later component which uses its fixed copy.  Apply
Theorem~\ref{thm:exact-full-pair-disintegration} before any local nonlinear
chart.  It assigns to the certificate one pre-coarea velocity block $U_C$;
the Source-B iterated integral assigns the disjoint nonvelocity input $Z_C$.
In reverse dependency order every state outside $U_C\oplus Z_C$ is held
fixed as a conditional parameter when the equidimensional local block $W_C$
is formed and integrated.

We prove the required pivot disjointness directly by induction over legal
cutting events.  Attach to a newly created certificate $C$ its protected
original atom set $A_C$ and its selected current free slot $\pi_C$.  At
creation, the priority rule assigns the cutting event to exactly one ledger,
and protection removes every atom in $A_C$ and every active representative
of $\pi_C$ before the next call.  Hence a later certificate $C'$ is created
from the active remainder and
\begin{equation}\label{eq:certificate-pivot-disjointness}
 A_C\cap A_{C'}=\varnothing,
 \qquad \pi_C\ne\pi_{C'}.
\end{equation}
A later legal cut partitions the active free slots.  When it breaks a bond,
the selected side receives the free copy and the other side receives only a
fixed shadow of that copy.  Equality of their transported values therefore
does not identify two integration coordinates: in reverse order the fixed
copy is a parameter and only the free copy is integrated.  This preserves
\eqref{eq:certificate-pivot-disjointness}, including when several protected
components share fixed boundary states in the same packet.  Protection of
the disjoint original atom sets $A_C$ also makes the local nonvelocity input
blocks $Z_C$ disjoint after the global velocity blocks have been assigned;
consequently their local output complements $\zeta_C$ are disjoint.  Shared
shadows are fixed parameters, not common integration variables.

The block $Z_C$ of a singleton or pair consists precisely of its continuous
nondistinguished angular, spatial-centre, time and auxiliary root variables
in the Source-B local operator.  Discrete cell and lift labels were already
fixed in the finite Borel branch and are not Lebesgue coordinates of $Z_C$.
After the pointwise envelope, its enlarged conditional domain
has measure at most $\cL_\eps^C$ by
\eqref{eq:source-elementary-operator}; a degree-two root substitution has
operator norm at most one.  The variables in $Z_C$ are therefore integrated together with
$w_C$, not deferred to the normal family.

Under a pair-to-singleton recut, the norm-one $\{2\}$ substitution retires
the old pair block and transfers the same single pivot $\pi_C$ to the good
$\{3\}$ singleton; no new identifier or free slot is created.  A shadow or
restart only renames the current representative of $\pi_C$, while the
protected original set $A_C$ is unchanged.  These are all transitions in
Lemma~\ref{lem:certificate-transition-table}, so the induction proves that
the selected pivots of distinct certificate identifiers remain distinct at
every stage.  If the certificates are listed in reverse dependency order,
the change from the corresponding original free slots to
$(U_C)_{C\in\mathscr C}\oplus U_{\rm rem}$ is first obtained by the square
global linear map of Theorem~\ref{thm:exact-full-pair-disintegration}.
Only then does the local conditional area formula form
$W_C=(w_C,\zeta_C)$ from its own equidimensional input $(U_C,Z_C)$; all
entries involving another block are fixed parameters.  The distinguished
$d$-dimensional local
determinants are uniformly nonzero.  Conditional Tonelli,
followed by Lemma~\ref{lem:periodic-certificate-measure} on each diagonal
block, now gives \eqref{eq:periodic-certificate-product-measure}.  Notice
that the argument uses disjoint global linear inputs before it uses
nonlinear local outputs; it does not place angular or auxiliary coordinates
in a velocity matrix.  This proves (ii)--(iii).
\end{proof}

\begin{center}
\small
\begin{tabular}{@{}p{0.14\textwidth}p{0.13\textwidth}p{0.17\textwidth}
 p{0.24\textwidth}p{0.22\textwidth}@{}}
\toprule
direction & lower exposed atom & terminals & incoming witnesses
& protected output\\
\midrule
UP
& C or O
& distinct
& no root conversion needed
& same-top-side \(\{33A\}\) pair\\
\addlinespace
UP
& O
& coincident
& \(t_n<t_p\) are incoming; equal lift is excluded by
  Lemma~\ref{lem:one-incoming-root}
& good \(\{3\}\) or good \(\{4\}\) singleton on two different lifts\\
\addlinespace
UP
& C
& coincident
& replace outgoing \(t_n\) by its incoming companion
  \(s_n<t_n<t_p\); equal lift is excluded
& good \(\{3\}\) or good \(\{4\}\) singleton; companion enlargement
  \(O(\eps^{1-1/(8d)})\)\\
\addlinespace
DOWN
& C or O
& distinct
& time reversal of the first row
& same-bottom-side \(\{33A\}\) pair\\
\addlinespace
DOWN
& O
& coincident
& time-reversed incoming witnesses
& good singleton on different lifts\\
\addlinespace
DOWN
& C
& coincident
& time-reversed companion construction
& good singleton on different lifts\\
\bottomrule
\end{tabular}
\end{center}

For the singleton signatures in
Lemma~\ref{lem:first-failure-finite-state}, the table exhausts the C/O and
exposed-side choices.  The two-atom signatures have already been sent
uniquely to the same-side \(\{33\mathrm A\}\) row of that lemma.  A purported
same-lift coincident-terminal singleton is not an additional output: after
both witnesses are put in incoming orientation it contradicts
Lemma~\ref{lem:one-incoming-root}.  The active degree alone decides whether
the singleton is recorded as good \(\{3\}\) or good \(\{4\}\).

\begin{figure}[t]
\centering
\begin{tikzpicture}[
 atom/.style={circle,draw,minimum size=5.5mm,inner sep=0pt},
 flow/.style={-{Latex[length=2mm]},thick},
 every node/.style={font=\scriptsize}]
 \node[align=center] at (0,1.48) {O-atom\\coincident terminal};
 \node[atom] (op) at (-1.4,0.55) {$p$};
 \node[atom] (on) at (-1.4,-0.55) {$n$};
 \node[atom] (oq) at (0.6,0) {$q$};
 \draw[flow] (op) to[bend left=16] node[above] {$m_p,t_p$} (oq);
 \draw[flow] (on) to[bend right=16] node[below] {$m_n,t_n$} (oq);
 \node[align=center] at (0,-1.25) {two incoming roots;\\same lift is impossible};

 \draw[dashed] (2,-1.55)--(2,1.65);
 \node[align=center] at (4.1,1.48) {C-atom\\coincident terminal};
 \node[atom] (cp) at (2.7,0.55) {$p$};
 \node[atom] (cn) at (2.7,-0.55) {$n$};
 \node[atom] (cq) at (4.5,0) {$q$};
 \node[atom] (cc) at (5.75,-0.75) {$q'$};
 \draw[flow] (cp) to[bend left=16] node[above] {$t_p$} (cq);
 \draw[flow] (cn) to[bend right=16] node[below] {$t_n$} (cq);
 \draw[flow,densely dotted] (cq) -- node[right] {companion $s_n$} (cc);
 \node[align=center] at (4.2,-1.35) {replace the outgoing witness;\\then apply the incoming-root test};
\end{tikzpicture}
\caption{Coincident-terminal first-failure alternatives.  The drawing fixes
the C/O and time-orientation convention; the two curves denote distinct
original particle-line traces, not two uses of one certificate.}
\label{fig:coincident-terminal-CO}
\end{figure}
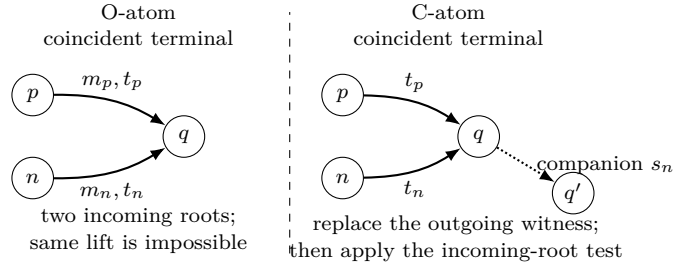

\begin{lemma}[Permanent ledger and recut invariance]
\label{lem:certificate-ledger}
There is a deterministic ledger with three disjoint classes:
\begin{equation}\label{eq:three-ledgers}
 \mathcal L_{\rm old},\qquad
 \mathcal L_{33},\qquad
 \mathcal L_{\rm par}.
\end{equation}
Here $\mathcal L_{\rm old}$ consists of the Source-B step-1 and degree-four
step-2 exceptions, $\mathcal L_{33}$ consists of ordinary adjacent-pair
outputs, and $\mathcal L_{\rm par}$ consists of the non-old first-failure
 certificates.  Every cutting event belongs to at most one class.  There are
two operations which must not be confused.  A \emph{permanent recut}
replaces a protected $\{33\mathrm A\}$ pair in the active source
decomposition by a good singleton and a norm-one $\{2\}$ component.  Its
identifier, owner and port list are transferred to the singleton.  If
$g_C\in\{0,1\}$ records whether that identifier was already counted by
$N_{\rm good}$ immediately before the recut, the atomic ledger update is
\begin{equation}\label{eq:simultaneous-ledger-transfer}
 \begin{split}
 N_{33A}&\longleftarrow N_{33A}-1,\qquad
 N_{\rm sing}\longleftarrow N_{\rm sing}+1,\\
 N_{\rm good}&\longleftarrow N_{\rm good}+(1-g_C).
 \end{split}
\end{equation}
Thus $g_C=1$ for a periodic good pair and for an ordinary pair already
declared support-good, whereas $g_C=0$ when the Source-B good partition
permanently turns a previously nongood ordinary pair into a good singleton.
By contrast, the exact $\{3\}+\{2\}$ refinement used only inside an owner or
conditional-measure calculation is temporary: it is coarsened before the
active decomposition is resumed and changes none of
$N_{33A},N_{\rm sing},N_{\rm good}$.  In particular, under every permanent
recut
$N_{\rm prot}:=N_{33A}+N_{\rm sing}$ and every previously proved assignment to
a certificate are invariant under recutting.  Here and below $N_{33A}$
counts every current $\{33\mathrm A\}$ pair representative, whether its
provenance ledger is ordinary or periodic, whereas $N_{\rm sing}$ counts
every protected identifier, ordinary or periodic, whose current
representative is a good singleton and no longer a pair.  The two summands
are therefore disjoint.  The ledger class
records provenance; the counter records the current analytic representative.
\end{lemma}

\begin{proof}
Order all atoms, slots and legal cuts once and for all.  At a cutting event,
first test membership in the old exceptional class.  If it is old, record it
there and do not make a new periodic charge.  Otherwise apply the ordinary
pair rule if it succeeds; only if that rule fails apply
Lemma~\ref{lem:first-failure}.  This priority rule makes
\eqref{eq:three-ledgers} disjoint.  The certificate identifier and the labels
in \eqref{eq:certificate-ports} refer to the molecule before any cut, whereas
a recut changes only the current representative.  Transferring the same
identifier proves every assertion in the lemma.
\end{proof}

The equivalent ordered priority rule is recorded in
Online Resource 1, Section~\ref{appA:certificate-priority}.

\begin{center}
\small
\begin{tabular}{@{}p{0.18\textwidth}p{0.14\textwidth}p{0.17\textwidth}
 p{0.24\textwidth}p{0.17\textwidth}@{}}
\toprule
transition & provenance & current representative & counter update
& full baseline\\
\midrule
old step-1 or degree-four step-2 exception
& \(\mathcal L_{\rm old}\)
& Source-B old object
& only the original Source-B \(N_4\)/old counter changes
& exactly as in Source B; no periodic charge\\
\addlinespace
ordinary adjacent-pair success
& \(\mathcal L_{33}\)
& protected \(\{33A\}\)
& \(N_{33A}{+}{=}1\); \(N_{\rm good}{+}{=}1\) only when its support is good
& none\\
\addlinespace
periodic first failure, pair output
& \(\mathcal L_{\rm par}\)
& protected good \(\{33A\}\)
& \(N_{33A}{+}{=}1,\ N_{\rm good}{+}{=}1\);
  \(N_{\rm sing},N_4\) unchanged
& none\\
\addlinespace
periodic first failure, degree-three singleton
& \(\mathcal L_{\rm par}\)
& protected good \(\{3\}\)
& \(N_{\rm sing}{+}{=}1,\ N_{\rm good}{+}{=}1\);
  \(N_{33A},N_4\) unchanged
& none\\
\addlinespace
periodic first failure, degree-four singleton
& \(\mathcal L_{\rm par}\)
& protected good \(\{4\}\)
& \(N_{\rm sing},N_{\rm good},N_4\) each increase by one
& exactly one \(\eps^{-(d-1)}\) full baseline; the periodic tube is the
  corresponding excess and does not erase it\\
\addlinespace
periodic good pair, or already-good ordinary pair, permanently recut
& identifier and owner unchanged
& good \(\{33A\}\to\) good \(\{3\}+\{2\}\)
& \(\Delta(N_{33A},N_{\rm sing})=(-1,+1)\);
  all other counters fixed
& none\\
\addlinespace
nongood ordinary pair made good by a permanent Source-B recut
& \(\mathcal L_{33}\) identifier and owner unchanged
& nongood \(\{33A\}\to\) good \(\{3\}+\{2\}\)
& \(\Delta N_{33A}=-1\),
  \(\Delta N_{\rm sing}=\Delta N_{\rm good}=+1\);
  \(N_{\rm prot},N_4\) fixed
& none\\
\addlinespace
temporary exact refinement for owner/measure computation
& identifier and owner unchanged
& \(\{33A\}\leftrightarrow\{3\}+\{2\}\) inside one local operator
& all five analytic counters fixed
& none\\
\addlinespace
protection, restart, or shadow transfer
& original ledger unchanged
& current representative renamed only
& every analytic counter unchanged
& none\\
\addlinespace
later attached/disjoint packet consumer
& not a new certificate event
& whole protected packet
& certificate counters are read, not incremented
& attached: one \(\eps^{-(d-1)}\) baseline; disjoint: rooted normalization only\\
\bottomrule
\end{tabular}
\end{center}

\begin{lemma}[State-transition completeness]
\label{lem:certificate-transition-table}
The preceding table exhausts every transition in the periodic UP, cyclic,
2CONNUP and DOWN consumers.  In particular, at every stage
\[
 N_{\rm prot}=N_{33A}+N_{\rm sing},
\]
the provenance classes remain disjoint, \(N_{\rm good}\) counts the current
analytic representative once, and no certificate transition duplicates a
full-component baseline.  A periodic good \(\{4\}\) creates exactly its one
source-required \(N_4\) charge and no second copy.
\end{lemma}

\begin{proof}
At a cut, the priority rule tests the old class, the ordinary pair rule and
the periodic first-failure rule in that order.  These are the first five
rows.  After an object exists, the only structural changes permitted by the
cutting algorithm are the simultaneous recut, protection/restart and shadow
renaming, which are the next two rows.  The attached/disjoint distinction is
made only after the protected special packet is inserted into the Source-B
consumer and therefore cannot change its certificate provenance.  The
displayed permanent counter changes are exactly
\eqref{eq:simultaneous-ledger-transfer}, including its pre-recut good bit;
the temporary refinement is an identity inside one positive local operator
and hence is not a state transition.  All other entries follow from the
definition of current representative, with the good-$\{4\}$ row carrying
the same $N_4$ charge as every Source-B full degree-four component.  Hence no
additional transition is possible.
\end{proof}

\subsection{The two local UP consumers}

\begin{proposition}[Periodic UP classification]
\label{prop:periodic-up-local}
The local UP conclusions in Source B, Proposition~10.2, parts (1)--(2), remain true for the
periodic molecule after every non-old failed no-double-\allowbreak overlap inference is
replaced by one protected certificate.  A component charged as a certificate
is not charged again as an ordinary good pair.
\end{proposition}

\begin{proof}
In part~(1), a degree-four atom becomes degree two either by the ordinary
Source-B mechanism or because the two paths that disappear at the same cut
are the two sides of a double support.  The former is unchanged and the
latter is Lemma~\ref{lem:first-failure}.  The invariant
$|E_*|-3|\mathcal M|$ is evaluated after removing the protected object, so
the original contradiction proof applies to the regular remainder.

In part~(2), the only new case is that the two free ends corresponding to the
fixed ends of a $\{33\mathrm A\}$ pair are serial through an overlap atom.
Their backward traces give a double support.  Recut the pair as in
Source~A, Part~4(2).  The singleton is good by
Lemma~\ref{lem:first-failure}, while the other atom is a normal $\{2\}$
component.  Lemma~\ref{lem:certificate-ledger} makes this a transfer rather
than a second charge.  All other branches are literally the Source-B fixed-
end partition.
\end{proof}

\subsection{Cyclic lower components: the first-certificate map}

Let $\mathcal M=\mathcal M_U\cup\mathcal M_D$ satisfy the hypotheses of
Source B, Proposition~10.9: $\mathcal M_U$ is a forest and no atom in
$\mathcal M_D$ is the parent of an atom in $\mathcal M_U$.  Denote by
$\operatorname{Cyc}(\mathcal M_D)$ the set of connected components of
$\mathcal M_D$ that contain a cycle.

\begin{lemma}[One original lower component per bearing port]
\label{lem:one-lower-component-per-port}
Fix an original incidence $\mathfrak i$ and its physical line
$\mathfrak p(\mathfrak i)$.  All atoms of $\mathcal M_D$ met by the lower
trace of this incidence belong to one original connected component of
$\mathcal M_D$.  This remains true after any sequence of cuts and restarts
when components and traces are interpreted by their original labels.
\end{lemma}

\begin{proof}
The no-$D$-to-$U$ order makes the $D$-portion of a physical line an interval
in its serial order.  Consecutive atoms on this interval are joined by the
intervening serial overlap path in $\mathcal M_D$, hence are in one connected
component.  A cut may delete a current representative of that path but does
not change its original atom set or the labels
\eqref{eq:original-incidence-label}; the assertion after restart is therefore
the same assertion about the pre-cut molecule.
\end{proof}

\begin{definition}[First-certificate map]
\label{def:first-certificate-map}
For $A\in\operatorname{Cyc}(\mathcal M_D)$ run UP in the common global
cutting order.  Let $n(A)$ be the first atom of $A$ that has degree two when
it is cut, unless an atom of $A$ has already entered a protected certificate.
In the first case let $p(A)$ be the atom whose cut first changes the degree of
$n(A)$ to two; in the second use the earlier protected certificate.  Let
$C(A)$ be the bearing record containing $p(A)$ or that earlier atom.  Finally
let $h(A)$ be the least port in \eqref{eq:certificate-ports} whose original
lower trace meets $A$.  Define
\begin{equation}\label{eq:first-certificate-map}
 \Phi_{\rm cyc}(A)=(C(A),h(A)).
\end{equation}
\end{definition}

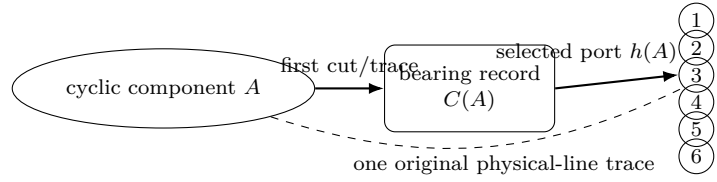
\begin{figure}[t]
\centering
\begin{tikzpicture}[
 every node/.style={font=\scriptsize},
 flow/.style={-{Latex[length=2mm]},thick}]
 \node[draw,ellipse,minimum width=2.35cm,minimum height=1.05cm] (A) at (0,0)
   {cyclic component $A$};
 \node[draw,rounded corners,minimum width=2.25cm,minimum height=1.15cm,
   align=center] (C) at (4.05,0) {bearing record\\$C(A)$};
 \foreach \y/\lab in {0.90/1,0.54/2,0.18/3,-0.18/4,-0.54/5,-0.90/6}{
   \node[circle,draw,minimum size=4.6mm,inner sep=0pt] (p\lab) at (7.05,\y) {$\lab$};}
 \draw[flow] (A) -- node[above] {first cut/trace} (C);
 \draw[flow] (C.east) -- (p3.west);
 \node at (5.60,0.48) {selected port $h(A)$};
 \draw[dashed] (A.south east) to[bend right=22]
   node[below,pos=0.56] {one original physical-line trace} (p3.south west);
\end{tikzpicture}
\caption{The cyclic-component injection.  The pair $(C(A),h(A))$ retains
the immutable original incidence.  One bearing port meets at most one
original lower component, while forgetting the port has fibre at most six.}
\label{fig:cyclic-bearing-port-injection}
\end{figure}

\begin{proposition}[Cyclic-component injection]
\label{prop:cyclic-certificate-injection}
The map \eqref{eq:first-certificate-map} is defined on every cyclic lower
component and is injective.  Consequently,
\begin{equation}\label{eq:lower-cyclic-count}
 N_{33A}+N_{\rm sing}\ge
 \frac16\,|\operatorname{Cyc}(\mathcal M_D)|.
\end{equation}
The statement and the map remain valid after protection, restart and the
recut transfer \eqref{eq:simultaneous-ledger-transfer}.
\end{proposition}

\begin{proof}
Choose a cycle in $A$ and the atom of that cycle cut last.  Its two neighbours
have already produced fixed ends, so some atom of $A$ has degree two when it
is cut.  Hence $n(A)$ exists unless $A$ has already met a certificate.

The proof of Source B, Proposition~10.9, now applies until the point where it
uses absence of double overlap.  If $p(A)$ is an overlap child of $n(A)$, or
is an overlap parent lying in $A$, the first-degree-two choice would force
both atoms to have degree three and the UP pair rule would already have
created a $\{33\mathrm A\}$ component.  Thus in the remaining ordinary case
$p(A)\in\mathcal M_U$, and the forest property of $\mathcal M_U$ again makes
$p(A)$ degree three at its cut; it belongs to the required pair.  If the
last inference fails, $p(A),n(A)$ form a double support and
Lemma~\ref{lem:first-failure} supplies a protected certificate.  In the
regular branch the immutable record attached to the ordinary pair plays the
same role.  Thus $C(A)$ is always
defined.  Its construction also shows that either it contains an atom of
$A$ or one of its bearing particle lines meets $A$, so $h(A)$ is defined.

Suppose $\Phi_{\rm cyc}(A)=\Phi_{\rm cyc}(A')=(C,h)$.  The original line of
port $h$ then meets both $A$ and $A'$.  Lemma~\ref{lem:one-lower-component-per-port}
gives $A=A'$.  Hence $\Phi_{\rm cyc}$ is injective.  Each bearing record has six
ports, so projection from $(C,h)$ to $C$ has fibres of size at most six and
\eqref{eq:lower-cyclic-count} follows.  Permanent original labels and
Lemma~\ref{lem:certificate-ledger} prove the last assertion.
\end{proof}

\subsection{The exact periodic 2CONNUP count}

Use exactly the Source-B sets.  Thus
$\mathcal M_{\rm 2conn}^1$ consists of lower degree-three atoms with one
uninterrupted upper adjacency, and $\mathcal M_{\rm 2conn}^2$ consists of
lower atoms adjacent to the upper molecule on two different particle lines.
Let $\mathcal L$ be the set of particle lines meeting both the upper molecule
and $\mathcal M_{\rm 2conn}$.  For $\mathfrak p\in\mathcal L$, let
$\sigma_{\mathfrak p}$ be the maximal overlap segment meeting both layers.
Define the auxiliary graph
\begin{equation}\label{eq:2conn-graph}
 G=(\mathcal L,E_G),\qquad
 \{\mathfrak p,\mathfrak q\}\in E_G
 \Longleftrightarrow
 \mathfrak p\cap\mathfrak q
 \text{ is an atom of }\mathcal M_{\rm 2conn}^2.
\end{equation}
The lower molecule is a forest, so $G$ is a forest: an auxiliary cycle lifts
to a lower atom cycle.  Two different lower intersections of the same two
lines would already give such a lower cycle.  Periodicity does not alter this
combinatorial fact.  Its new effect is instead the following: one upper cut
may break $\sigma_{\mathfrak p}$ and $\sigma_{\mathfrak q}$ simultaneously
even when $\{\mathfrak p,\mathfrak q\}\in E_G$.  The upper cut atom and their
lower intersection are then the endpoints of a double support.

For each $\mathfrak p$, let $b(\mathfrak p)$ be the first cutting event that
breaks $\sigma_{\mathfrak p}$.  Orient every edge of $G$ from the endpoint
with earlier breaking time to the endpoint with later breaking time.  If the
times agree, call the edge \emph{simultaneous}.  The priority rule of
Lemma~\ref{lem:certificate-ledger} gives the disjoint partition
\begin{equation}\label{eq:2conn-edge-partition}
 E_G=E_{\rm reg}\mathbin{\dot\cup}E_{\rm old}
       \mathbin{\dot\cup}E_{\rm par},
\end{equation}
where a simultaneous edge belongs to $E_{\rm old}$ when its break is a
Source-B step-1 or degree-four step-2 exception, and otherwise belongs to
$E_{\rm par}$.  Every edge in $E_{\rm par}$ creates one protected
$\{33\mathrm A\}$ or good $\{3\}$ certificate.  Indeed the non-old breaking
atom has degree three; hence the good-$\{4\}$ alternative in
Lemma~\ref{lem:first-failure} cannot occur here.

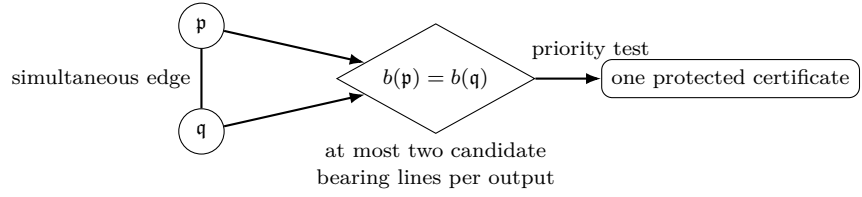
\begin{figure}[t]
\centering
\begin{tikzpicture}[
 every node/.style={font=\scriptsize},
 line/.style={circle,draw,minimum size=6mm,inner sep=0pt},
 flow/.style={-{Latex[length=2mm]},thick}]
 \node[line] (p) at (0,0.7) {$\mathfrak p$};
 \node[line] (q) at (0,-0.7) {$\mathfrak q$};
 \draw[thick] (p)--node[left] {simultaneous edge} (q);
 \node[draw,diamond,aspect=1.8] (b) at (3.1,0) {$b(\mathfrak p)=b(\mathfrak q)$};
 \node[draw,rounded corners,minimum width=2.3cm] (cert) at (7.0,0)
   {one protected certificate};
 \draw[flow] (p) -- (b);
 \draw[flow] (q) -- (b);
 \draw[flow] (b) -- (cert);
 \node[above] at (5.15,0.12) {priority test};
 \node[align=center] at (3.1,-1.15) {at most two candidate\\bearing lines per output};
\end{tikzpicture}
\caption{The periodic change in the 2CONNUP candidate map.  A simultaneous
auxiliary edge has one breaking event and creates at most one non-old
certificate; conversely that output owns at most its two endpoint
candidates.  This is the two-to-one fibre in
Lemma~\ref{lem:2conn-candidates}.}
\label{fig:2conn-simultaneous-fibre}
\end{figure}

\begin{lemma}[2CONN candidate accounting]
\label{lem:2conn-candidates}
Let $K=400d^2$.  Let $\mathcal L_0$ be the lines selected in step~1 of
2CONNUP, $\mathcal L_1$ all lines broken there, and let
$\mathcal P_1$ be the edges between $\mathcal L_1$ and its complement.  Order
$\mathcal L\setminus\mathcal L_1=(\mathfrak p_1,\ldots,\mathfrak p_s)$ by
breaking time and let $\mathcal P_2$ be the edges among these lines.  If
\begin{equation}\label{eq:2conn-cert-count-def}
 N_{\rm 2c}:=N_{33A}^{\rm 2c}+N_{{\rm sing},3}^{\rm 2c},
\end{equation}
then
\begin{align}
 |\mathcal L_1|&\le2|\mathcal L_0|,
 \label{eq:2conn-L01}\\
 |\mathcal P_1|&\ge K|\mathcal L_0|-2|\mathcal L_1|,
 \qquad
 |\mathcal P_1|+|\mathcal P_2|
 \ge |\mathcal M_{\rm 2conn}^2|-|\mathcal L_1|,
 \label{eq:2conn-P12}\\
 N_{\rm 2c}&\ge \frac12\max\left\{
 |\mathcal P_1|-|\mathcal L_1|,
 |\mathcal M_{\rm 2conn}^1|,
 K^{-1}|\mathcal P_2|
 \right\}-E_{\rm old}.
 \label{eq:2conn-candidate-lower}
\end{align}
Here $N_{{\rm sing},3}^{\rm 2c}$ counts only degree-three protected
singletons created by, or transferred from a protected pair created by, the
present 2CONN candidate call.  Thus an ordinary or periodic pair-to-singleton
transfer changes its analytic representative but not $N_{\rm 2c}$.
Here $E_{\rm old}$ is the number of old exceptional cutting events.  No
certificate counted in $N_{\rm 2c}$ is counted in $E_{\rm old}$.
\end{lemma}

\begin{proof}
The first two lines are the forest calculation of Source B,
Proposition~11.2.  Each selected line has at least $K$ incident edges, one
step-1 cut breaks at most two lines, and the subgraph induced by
$\mathcal L_1$ and its neighbours is a forest.  After step~1 every remaining
vertex has degree at most $K$, and every edge not lost inside
$\mathcal L_1$ lies in $\mathcal P_1\cup\mathcal P_2$, proving
\eqref{eq:2conn-P12}.

It remains to check the only periodic change in the candidate-to-output
map.  The three candidate families in the maximum are, exactly as in Source
B: lines in the external neighbourhood of $\mathcal L_1$, lines carrying an
atom of $\mathcal M_{\rm 2conn}^1$, and later endpoints of oriented edges in
$\mathcal P_2$.  At the breaking of a candidate line, either the Source-B
argument creates a $\{33\mathrm A\}$ pair, or an old exceptional event is
charged, or its breaking is simultaneous with an adjacent line.  In the last
case the two breaking paths and the lower intersection form a double support,
so Lemma~\ref{lem:first-failure} creates a protected $\{33\mathrm A\}$ or
good $\{3\}$ certificate.

One output contains at most two candidate bearing lines.  For a regular pair
this is the Source-B two-to-one bound.  For a periodic certificate it follows
because the triggering cut breaks at most the two lines in its unique
simultaneous auxiliary edge.  In the third family, a fixed later vertex has
at most $K$ earlier neighbours, which is the displayed $K^{-1}$.  Subtracting
the old events gives \eqref{eq:2conn-candidate-lower}.  The last claim is the
priority rule in Lemma~\ref{lem:certificate-ledger}.
\end{proof}

\begin{proposition}[Periodic 2CONNUP source inequality]
\label{prop:periodic-2connup}
For periodic 2CONNUP,
$N_{33B}=N_{44}=0$ and
\begin{equation}\label{eq:2conn-count}
 N_{\rm 2c}-50d^2N_4
 \ge c_d\bigl(|\mathcal M_{\rm 2conn}^1|
              +|\mathcal M_{\rm 2conn}^2|\bigr)
       -C_dN_{\rm full}(\mathcal M_U),
\end{equation}
where one may take $c_d=[9(400d^2+4)]^{-1}$ and a finite
$C_d\le 2+100d^2$.  The dual statement holds for 2CONNDN.
\end{proposition}

\begin{proof}
The fixed-end orientation argument in Source B is unchanged by protection,
so $N_{33B}=N_{44}=0$.  Put $K=400d^2$ and denote the maximum in
\eqref{eq:2conn-candidate-lower} by $Y$.  From
\eqref{eq:2conn-L01}--\eqref{eq:2conn-P12},
\[
 |\mathcal P_1|-|\mathcal L_1|
 \ge (K/2-3)|\mathcal L_1|,
 \qquad |\mathcal L_1|\le 3Y/K.
\]
Moreover $|\mathcal P_1|\le Y+|\mathcal L_1|\le2Y$ and
$|\mathcal P_2|\le KY$, whence
\begin{equation}\label{eq:2conn-Y-source}
 |\mathcal M_{\rm 2conn}^1|+|\mathcal M_{\rm 2conn}^2|
 \le (K+4)Y.
\end{equation}
The two unchanged Source-B component bounds are
\begin{equation}\label{eq:2conn-old-bounds}
 \max\{E_{\rm old},N_4\}
 \le N_{\rm full}(\mathcal M_U)+|\mathcal L_1|.
\end{equation}
Using $50d^2=K/8$, Lemma~\ref{lem:2conn-candidates} and
\eqref{eq:2conn-old-bounds} give
\begin{align*}
 N_{\rm 2c}-50d^2N_4
 &\ge \frac Y2-(1+K/8)
       \bigl(N_{\rm full}(\mathcal M_U)+|\mathcal L_1|\bigr)\\
 &\ge \left(\frac12-\frac3K-\frac38\right)Y
       -(1+K/8)N_{\rm full}(\mathcal M_U)\\
 &\ge \frac Y9-(1+K/8)N_{\rm full}(\mathcal M_U),
\end{align*}
because $K=400d^2\ge400$.  Replacing the penultimate coefficient by
$1/8$ would be false: its exact value is $1/8-3/K$.
Together with \eqref{eq:2conn-Y-source} this is
\eqref{eq:2conn-count}.  Reversing the layer order proves the dual result.
The coefficient $50d^2$ is exactly the Source-B coefficient; it
has not been weakened by the periodic excision.
\end{proof}

\subsection{DOWN tails and blockwise injection}

The last omitted assignment occurs in Source A, Part~5, and in the final
DOWN construction of Source B, Proposition~12.1.  We formulate it using the
original block labels appearing there.

\begin{definition}[Original recollisional blocks]
\label{def:recollisional-blocks}
In the subleading DOWN molecule write the C-atom sets as
$\mathcal M_j$, $j\in T$, and the connecting overlap atoms according to the
source tree $(T,E)$, rooted at $1$.  For every $j\ne1$, let
$\mathsf o_j$ be the unique incoming overlap atom on the edge from the
parent of $j$; this is the ``cute'' atom in the final-DOWN construction of
Source B, Proposition~12.1.  Assign that atom to its child and define
\begin{equation}\label{eq:down-original-block-partition}
 V_1:=\mathcal M_1,
 \qquad V_j:=\mathcal M_j\cup\{\mathsf o_j\}\quad(j\ne1),
 \qquad \mathcal B_j:=V_j,
 \qquad V=\mathop{\dot\bigcup}_{j\in T}V_j .
\end{equation}
This is a disjoint partition by original atom identifiers: the
 $\mathcal M_j$ are disjoint C-atom sets and each tree edge has exactly one
overlap atom, assigned to its unique child.  It is deliberately not defined
as a path closure.  Source B cuts exactly $\mathcal M_1$ in the root call
and $\mathcal M_j\cup\{\mathsf o_j\}$ in a nonroot call.  A block is
recollisional when its original C-atom interaction graph contains a cycle.
Blocks are ordered by those actual calls in the Source-B restart procedure.
\end{definition}

\begin{lemma}[Block owner and recut invariance]
\label{lem:down-block-injection}
For every recollisional original block $\mathcal B_j$, the localized DOWN
procedure produces a protected good object $C_j$ before moving to the next
block.  The map
\begin{equation}\label{eq:down-block-map}
 \Phi_{\rm down}:\mathcal B_j\longmapsto C_j
\end{equation}
is injective.  Its owner is the unique original block containing the first
bearing atom of $C_j$.  Protection, restart and
\eqref{eq:simultaneous-ledger-transfer} do not change this owner.
\end{lemma}

\begin{proof}
Cut $\mathcal B_j$ as free, including its unique incoming overlap atom when
$j\ne1$; for $j=1$ retain the distinguished output-end boundary of the
source construction.  A recollisional block then contains a cycle and
has exactly the one degree-three boundary atom specified by the rooted block
tree.  The Source-B DOWN proposition therefore produces a
$\{33\mathrm A\}$ component, but that statement alone does not say that the
pair is good.  Apply the finite positive partition in part~(2) of the same
Source-B UP/DOWN proposition to this pair and to the atom that produced its
first fixed end.  If every Euclidean no-double-overlap inference in that
argument is valid, the proposition selects an actual good output: either the
$\{33\mathrm A\}$ pair itself or its good $\{3\}$/$\{4\}$ fixed-end
producer.  Define $C_j$ to be that selected output and protect it.  If the
first such inference fails because the two supports use different periodic
lifts, Lemma~\ref{lem:first-failure} replaces precisely that step by a
periodic good certificate, which we instead call $C_j$ and protect.  Thus in
all positive cells $C_j$ is genuinely one of the good objects counted by
$N_{\rm good}$ and exists before the algorithm leaves $\mathcal B_j$.
The choices are made once per block, so their total positive partition has
at most $C^{|T|}\le C^{|\mathcal M|}$ cells and is already within the source
indicator allowance.

The boundary inputs and outputs just used are exhaustive:
\begin{center}
\small
\begin{tabular}{@{}p{0.18\textwidth}p{0.28\textwidth}p{0.42\textwidth}@{}}
\toprule
block/cell & boundary supplied to localized DOWN & protected output\\
\midrule
$j=1$
 & the distinguished source output end; no parent overlap atom
 & the first ordinary good output selected by DOWN, or the first periodic
 certificate replacing its failed no-double-overlap inference\\
$j\ne1$
 & the unique parent overlap atom $\mathsf o_j$, assigned to the child
 block in \eqref{eq:down-original-block-partition}
 & the same two alternatives, before the call returns to the parent block\\
ordinary successful cell
 & the literal free/fixed incidences created by Source-B cutting cases
 (0),(2)--(4)
 & a good $\{33{\rm A}\}$, $\{3\}$ or $\{4\}$ object selected by the
 finite DOWN partition\\
periodic first-failure cell
 & the same incidences together with the two stored original support lines
 & the unique certificate of Lemma~\ref{lem:first-failure}; a full
 $\{4\}$ baseline, when present, is retained\\
\bottomrule
\end{tabular}
\end{center}
There is no no-boundary or two-parent-boundary case: the rooted block tree
gives zero parent overlaps for $j=1$ and exactly one for $j\ne1$.

The localization is at the level of atom identifiers, which is the level
needed for injection.  Cut the atom set $V_j$ in
\eqref{eq:down-original-block-partition} as free.  Source B,
Definition~8.6 and the paragraph following it, state for each of cases
(2)--(4) that the result has one molecule with atom set $V_j$ and another
with the complementary atom set.  This remains true in case (2), where an
ov-segment can leave $V_j$ and return: the newly inserted shortcut bond has
both endpoints in $V_j$, whereas every intermediate O-atom remains in the
complement.  Hence no path-convexity assertion is required.  Inductively,
through the DOWN call,
\begin{enumerate}[label=\textup{(L\arabic*)}]
\item every current atom has an original identifier in $V_j$;
\item every elementary component or periodic first-failure certificate
produced in this call has all of its bearing atoms in $V_j$;
\item after a component is protected, every later cut acts on a subset of
the still-unprocessed original identifiers and never reintroduces an
identifier already removed.
\end{enumerate}
Indeed, the induction step cuts some $A\subset V_j$.  The same atom-set
statement gives the two new sets $A$ and $V_j\setminus A$.  A shortcut
bond or a matched pair changes incidences but creates no atom.  The periodic
companion record of Lemma~\ref{lem:first-failure} is attached to the first
active bearing atom and likewise introduces no atom.  Protection then
removes one of these disjoint atom subsets.  This proves (L1)--(L3).

In particular, the first bearing atom of every protected regular pair or
periodic certificate has a unique permanent original atom identifier in
$V_j$ and hence a unique owner.  The
sets $V_j$ are disjoint by \eqref{eq:down-original-block-partition}, so
objects with different owners cannot coincide.  Once $C_j$ is protected,
Source B's protection convention removes its entire current atom component
from all later cuts.  A later restart acts on a disjoint $V_k$ and therefore
cannot recreate an atom identifier of $C_j$.  Thus
\begin{equation}\label{eq:down-protected-atom-equalities}
 A(C_j)\subset V_j,\qquad
 V_j=A(C_j)\,\dot\cup\,(V_j\setminus A(C_j)),\qquad
 (V_j\setminus A(C_j))\cap\bigcup_{k\ne j}V_k=\varnothing .
\end{equation}
The first set is permanently removed and the second is the only set on
which the localized restart may continue; matched shadows change edge
states but none of the three atom sets in
\eqref{eq:down-protected-atom-equalities}.  Consequently
\eqref{eq:down-block-map} is injective, and
Lemma~\ref{lem:certificate-ledger} proves recut invariance.
\end{proof}

\begin{proposition}[Periodic final DOWN tail]
\label{prop:periodic-down-tail}
The leading and subleading final DOWN conclusions used in Source B,
Proposition~12.1, and Source A, Part~5, remain valid on $\T^d$.
In the leading case at least one protected good component is produced and no
new $\{4\}$ baseline occurs.  In the subleading case the number of distinct
protected good components is at least the number of recollisional original
blocks.
\end{proposition}

\begin{proof}
In the leading case fixing the output end turns the tree into a molecule with
one cycle.  Its components are nonfull and contain neither a degree-two atom
nor a C-bottom fixed end, exactly as in the source proof.  DOWN therefore
produces a first $\{33\mathrm A\}$ pair.  If its two bearing ends are regular,
the Source-B fixed-end argument makes it good.  If they are the two serial
sides of a double support, Lemma~\ref{lem:first-failure} produces a protected
good pair or singleton.  The full/nonfull flag is the structural
degree-four-incidence flag of the UP/DOWN decomposition and is fixed before
the ordinary/periodic good bit is attached.  Every ordinary or periodic
$\{4\}$ singleton is therefore full and carries exactly one $N_4$ baseline,
whether or not it is good.  Since every component in this leading cleanup is
nonfull, no $\{4\}$ singleton of either provenance can be created.  This
proves the asserted absence of a new baseline, not merely the absence of a
new bad $\{4\}$.

In the subleading case apply Lemma~\ref{lem:down-block-injection}.  It gives
one distinct protected good object for each recollisional block.  Since the
map is expressed in original block labels, later restarts and recuts cannot
merge two charges.  This is the exact distinctness assertion omitted in the
one-paragraph source tail.
\end{proof}

Propositions~\ref{prop:periodic-up-local},
\ref{prop:cyclic-certificate-injection},
\ref{prop:periodic-2connup}, and~\ref{prop:periodic-down-tail} prove the
periodic cutting input (R4).  The contribution relative to Source A is not
the discovery of double overlap segments---they are already identified in
its Section~5---but the different-lift admissibility argument, the permanent
certificate ledger, and the bounded-fibre assignments needed already in
dimension four and hence throughout the present fixed-$d\ge4$ argument.

\section{Two-Landing Packet Geometry and Proof of the Main Technical Theorem}
\label{sec:two-landing}

The exact two-landing operator has already been stated in
Section~\ref{sec:overview}.  This section proves it.  The proof first selects
the two landing roots, treats the disjoint and overlapping landing geometries,
establishes flux completion, and then invokes the positive-network and
finite-fibre modules of Sections~\ref{app:positive-network}
and~\ref{app:coarea}.  The complete nonnegative joint kernel remains intact
throughout; the full-component normalization and each retained owner time are
charged only once.  Independent full-chord and conditional charts are recorded
in Online Resource~1 and do not replace the sharp radial-conjugate route.

\subsection{Particle graph and two distinct landing atoms}

The atom graph has the atoms as vertices and the molecule bonds as edges.
The physical-particle graph $\Gamma_U$ has one vertex for every physical
particle line and one edge for every upper collision atom, joining the two
lines present at that collision.  Parallel collision edges are retained.

\begin{lemma}[Upper tree and landing bonds]
\label{lem:upper-tree-landings}
Under assumptions (a)--(e), $\Gamma_U$ is connected and spans all physical
particle lines.  The crossing bonds may be chosen as $e_1,e_2$ so that both
their upper endpoints $u_1,u_2$ and their lower endpoints $d_1,d_2$ are
distinct.  Moreover, $d_1<d_2$ may be taken to be the first two lower
collision atoms in the fixed strict order.  In particular, there is no
lower contact before $d_1$ or strictly between $d_1$ and $d_2$.
Their two upper collision edges form a forest and hence extend to
a collision spanning tree $S_U$ of $\Gamma_U$.  The extension may be chosen
by reverse chronological Kruskal: scan upper contacts from latest to earliest
and retain an edge exactly when it joins two current components.  With this
choice $u_1,u_2\in S_U$, and every unmarked upper contact has its two
endpoints joined by an $S_U$-path consisting entirely of later contacts.
Every lower collision is outside this tree.  The unordered particle-line
pairs of $d_1,d_2$ are different.
\end{lemma}

\begin{proof}
Two upper atoms joined by a molecule bond lie on a common physical particle
line, so the corresponding two edges of $\Gamma_U$ share a vertex.
Connectivity of the upper atom graph therefore places all its collision
edges in one connected edge component of $\Gamma_U$.  Assumption (b) says
that every physical particle line occurs at some upper atom, hence no vertex
is missing or isolated.  Thus $\Gamma_U$ is connected and spanning.

Let $d_1<d_2$ be the first two lower collision atoms.  Such a $d_2$ exists:
$d_1$ contains only two particle lines, whereas assumptions (b) and (e)
give a third line which also occurs in the lower sublayer.  Write $a,b$ for
the two lines at $d_1$.  If both lines at $d_2$ had already occurred in the
lower sublayer, they would have to be $a,b$.  Because no lower atom lies
strictly between $d_1$ and $d_2$, the two line bonds from $d_1$ to $d_2$
would then be parallel molecule bonds, contrary to the no-double-bond part
of (e).  Hence $d_2$ contains a new line; call it $c$ and choose its crossing
bond as $e_2$.  Its lower endpoint is $d_2$.

For a particle line $p$, denote its last upper atom by $u(p)$.  At most one
of $u(a),u(b)$ equals $u(c)$.  Indeed, if $u(a)=u(b)$, that C-atom already
uses its two slots on $a,b$ and cannot contain $c$; if $u(a)\ne u(b)$, the
claim is immediate.  Choose $p\in\{a,b\}$ with $u(p)\ne u(c)$ and take its
crossing bond as $e_1$.  Then
\[
 u_1=u(p)\ne u(c)=u_2,
 \qquad d_1\ne d_2.
\]
Thus the chosen crossing bonds have distinct upper and lower endpoints, and
the lower endpoints are precisely the first two lower atoms.

The two upper collision edges are not parallel.  Indeed, if they joined the
same two physical lines, the later of $u_1,u_2$ would be the last upper atom
on both lines, contradicting the fact that the two selected lines have the
distinct last upper endpoints $u_1,u_2$.  Hence these two edges form a
forest.

Scan all upper collision edges in decreasing time order and retain an edge
precisely when its endpoints are in different components of the already
retained forest.  The selected physical line at $u_i$ has no later upper
collision.  It is therefore an isolated vertex of the later-edge graph just
before $u_i$ is scanned, so the edge at $u_i$ must be retained.  At the end
the retained forest is connected because $\Gamma_U$ is connected, and hence
it is a spanning tree containing both selected edges.  If an edge is
rejected, its endpoints were already joined at that moment; the retained
joining path uses only edges scanned earlier, that is, contacts later in
physical time.  Call this reverse-Kruskal tree $S_U$.  Since every edge of
$S_U$ is an upper atom, every lower collision is a co-tree contact relative
to $S_U$.

Finally, the unordered particle pairs at $d_1,d_2$ are different: the latter
contains the new line $c$, while the former contains only $a,b$.  Thus their
two incidence covectors in the particle-line graph are nonparallel.
\end{proof}

The matching refinement also exposes an exact local determinant which is
useful for assessing possible replacements for
\eqref{eq:tree-exterior-coercivity}.  Its role is local; the global packet chart is constructed below.

\begin{lemma}[Two-sided free-segment bridge]
\label{lem:two-sided-bridge}
For either matched crossing line, put
$a_i=s-t_{u_i}>0$ and $b_i=t_{d_i}-s>0$.  Freeze the two companion states at
the endpoints of that line and retain its separator state $(X_i,V_i)$.  In
top/bottom contact representations compatible with the transversal cut, the
upper and lower vector residuals have the affine form
\[
 U_i=X_i-a_iV_i-Y_i^U,
 \qquad
 R_i=X_i+b_iV_i-Y_i^D.
\]
Consequently
\begin{equation}\label{eq:two-sided-bridge-determinant}
 \left|\det D_{(X_i,V_i)}(U_i,R_i)\right|
 = (a_i+b_i)^d=(t_{d_i}-t_{u_i})^d
 \ge \epstar^d.
\end{equation}
For the two matched lines the product of the two conditional determinants is
at least $\epstar^{2d}$.
\end{lemma}

\fallbackproofinresource{supp:proof:lem:two-sided-bridge}

\begin{remark}[Role of the bridge identity]
\label{rem:bridge-not-global}
In the full molecule, the companion states in
Lemma~\ref{lem:two-sided-bridge} remain coupled to the other collision
constraints.  Eliminating the remaining upper-tree
constraints and the rooted co-tree contacts differentiates those companions
with respect to $(X_i,V_i)$ and produces exactly the Schur terms in
\eqref{eq:tree-schur-row}.  Thus multiplying the two local determinants would
repeat the original fixed-boundary error.  A successful bridge-based proof
must exhibit one global ordering of all raw edge-state constraints in which
the companion blocks are genuinely triangular, or else control the resulting
Schur complement.  The argument below therefore controls the resulting Schur complement directly.
\end{remark}

There is, however, one genuinely noncollinear word on which the bridge
coordinates do form a complete global chart.  This gives an unconditional
test of the mechanism beyond the common-normal stratum considered later.

\begin{proposition}[Arbitrary-normal aligned three-line bridge]
\label{prop:aligned-three-line-bridge}
Suppose there are exactly three physical particle lines, labelled
$0,c,2$, and exactly four contacts.  In the upper sublayer the contacts are
$(0,c)$ at $t_{u_1}$ and then $(c,2)$ at $t_{u_2}$; in the lower sublayer
they are $(0,c)$ at $t_{d_1}$ and then $(c,2)$ at $t_{d_2}$.  There are no
other contacts.  The four collision normals and velocities may be arbitrary
on the regular incoming cell.

Put $h=\epstar/4$ and choose a fixed separator number $s$ satisfying
$t_{u_2}+h\le s\le t_{d_1}-h$; the fixed-grid construction in
\eqref{eq:separator-grid-choice} supplies such pieces.  Retain the separator
state of the centre line $c$ as a parameter and use the
separator positions and velocities $(x,u)$ and $(y,w)$ of lines $0$ and $2$
as the $4d$ scalar pivot variables.  Then the four vector contact residuals have
a Jacobian satisfying
\begin{equation}\label{eq:aligned-three-line-lower}
 \left|\det D_{(x,u,y,w)}(U_1,U_2,R_1,R_2)\right|
 \ge 2^{2d}h^{2d},
 \qquad h=\epstar/4.
\end{equation}
Consequently this complete packet, with an arbitrary nonnegative joint test
function, obeys the unconditional cell estimate
\begin{equation}\label{eq:aligned-three-line-operator}
 J_{H,\kappa}(Q)
 \le \cL_\eps^C\eps^{-3(d-1)}\eps^{4(d-1)}h^{-2d}
       \int Q_{H,\kappa}^{\sharp,\mathrm{br}}
 \le \cL_\eps^C\eps^{d-1}\epstar^{-2d}
       \int Q_{H,\kappa}^{\sharp,\mathrm{br}}.
\end{equation}
Here $Q_{H,\kappa}^{\sharp,\mathrm{br}}$ denotes the resulting zero-extended
bridge output.  Its domain is independent of the retained centre-line state.
\end{proposition}

\fallbackproofinresource{supp:proof:prop:aligned-three-line-bridge}

\begin{proposition}[All three-line tree words]
\label{prop:all-three-line-bridge}
Suppose there are exactly three physical particle lines, exactly two upper
contacts forming their collision tree, exactly two lower contacts containing
the two matched first landings, and no other contacts.  Then, for arbitrary
regular collision normals, one can choose two crossing-line separator states
so that the four vector contact residuals have a $4d$-dimensional pivot
minor bounded below by $2^{2d}h^{2d}$.  Hence the unconditional arbitrary-$Q$
bound \eqref{eq:aligned-three-line-operator} holds for every such word.
\end{proposition}

\fallbackproofinresource{supp:proof:prop:all-three-line-bridge}

\begin{proposition}[Three lines with arbitrary upstream roots]
\label{prop:three-line-upstream-rooted}
Assume that the packet has exactly three physical particle lines and no
double bond.  Choose the matched crossing lines and the reverse-Kruskal tree
as in Lemma~\ref{lem:upper-tree-landings}.  Allow arbitrary unmarked upper
contacts and arbitrary unmarked lower contacts strictly after $d_2$.  Since
$d_1,d_2$ are the first two lower atoms, there is no earlier lower contact.
Then, for arbitrary regular collision
normals,
\begin{equation}\label{eq:three-line-upstream-rooted-operator}
 J_{H,\kappa}(Q)
 \le \cL_\eps^C\eps^{d-1}\epstar^{-2d}
       \int Q_{H,\kappa}^{\sharp,3\mathrm r}.
\end{equation}
The output box is independent of the fixed boundary values.
The assertion about upstream roots is made in the raw separator/event-state
chart before any incoming-root substitution: chronological event states are
independent coordinates, the four selected residuals are then differentiated,
and only afterward are the upstream root blocks eliminated.  It is not a
claim of independence in the original initial-data coordinates.
\end{proposition}

\fallbackproofinresource{supp:proof:prop:three-line-upstream-rooted}

\begin{proposition}[All tree-only words, arbitrary normals]
\label{prop:all-tree-only-bridge}
Suppose the word consists only of an upper collision spanning tree and the
two selected lower first-landing contacts.  Then, for every number of
physical particle lines allowed by assumptions \textup{(a)--(e)}, the packet
obeys the unconditional arbitrary-$Q$ estimate
\begin{equation}\label{eq:all-tree-only-operator}
 J_{H,\kappa}(Q)
 \le \cL_\eps^C\eps^{d-1}\epstar^{-2d}
       \int Q_{H,\kappa}^{\sharp,\mathrm{tree}}.
\end{equation}
No relation among the collision normals is assumed.
\end{proposition}

\fallbackproofinresource{supp:proof:prop:all-tree-only-bridge}

\begin{proposition}[Disjoint landings with one-sided rooted contacts]
\label{prop:disjoint-rooted-extension}
Assume the two selected lower edges $f_1,f_2$ are disjoint.  Allow arbitrary
unmarked contacts in the upper sublayer and arbitrary unmarked lower contacts
strictly after $d_2$.  The first-two choice excludes an earlier lower
contact.  Then the
unconditional operator estimate
\begin{equation}\label{eq:disjoint-rooted-operator}
 J_{H,\kappa}(Q)
 \le \cL_\eps^C\eps^{d-1}\epstar^{-2d}
       \int Q_{H,\kappa}^{\sharp,\mathrm{disj}}
\end{equation}
holds, with a constant depending only on the packet cap.  The output domain
is fixed by the cell and is independent of the boundary values.
\end{proposition}

\fallbackproofinresource{supp:proof:prop:disjoint-rooted-extension}

The fact that $d_1\ne d_2$ has an analytic consequence: $e_1$ cannot be the
companion top edge at $d_2$, and $e_2$ cannot be the companion top edge at
$d_1$.  This makes the first block of the iterated coordinate map diagonal
before any other molecule variable is eliminated.

Fix the strict total order of lower atoms belonging to the half-open cell.
For a crossing bond $e$ write $\lambda(e)$ for its first lower atom.  If a
lower edge state is evaluated immediately above an atom $d$, define its
\emph{lower influence cone} to be the set of crossing-bond states on which
it depends after all strictly earlier lower elastic maps and transports have
been performed.

\begin{lemma}[Causal lower influence cones]
\label{lem:lower-influence-cone}
The state immediately above $d$ depends only on crossing bonds $e$ with
$\lambda(e)<d$, together with its own unchanged crossing state when the
physical line has $\lambda(e)\ge d$.  In particular, at
$d=\lambda(e)$ the companion top-edge state is independent of the velocity
and crossing position on $e$.  If $d_1<d_2$, then the companion state at
$d_1$ is independent of both selected states, while the companion state at
$d_2$ may depend on the state on $e_1$ but is independent of the state on
$e_2$.
\end{lemma}

\begin{proof}
Immediately below the transversal cut, every particle line carries only its
own crossing state.  Process lower atoms in the fixed strict order.  The
complete elastic and transport maps at the next atom use only the two states
arriving at that atom, and their outputs continue on those same two physical
lines.  Thus a crossing state can enter another line's influence cone only
at a collision involving its present descendant line.  Before
$\lambda(e)$, the line of $e$ has had no lower collision, so its information
cannot have left that line.  At $\lambda(e)$ the other top slot belongs to a
different physical line and consequently cannot depend on $e$.

If a selected state $e_j$ occurred as the companion at an earlier atom, that
earlier atom would already be $\lambda(e_j)$, contradicting the definition
of $d_j$.  Therefore no state whose first landing is later can influence an
earlier companion.  Applying this observation first at $d_1$ and then at
$d_2$ proves the final assertion.  The same induction applies to physical
positions because free transport propagates a state only along its own
particle line.
\end{proof}

\subsection{The two raw landing coordinates}

We first record two exact changes of representation.  At a C-atom the
elastic scattering map is an orthogonal involution on the velocity pair.
It preserves the contact sphere and the absolute flux.  Hence, on a fixed
incoming-sign cell, the atom distribution may be written with the contact
condition on either its bottom pair or its top pair; the other pair is then
recovered by the two velocity deltas.  We use the top-pair representation at
$d_1,d_2$.  More generally, in the global recovery below we use the bottom
pair at every upper atom and the top pair at every lower atom.  Thus the pair
on the transversal side of an atom is always the free pair, independent of
that atom's own time.

Second, in dimension $d$ the sphere identity is
\begin{equation}\label{eq:sphere-delta-identity}
 \delta(|R|-\eps)F(R)
 =\eps^{d-1}\int_{\Sph^{d-1}}\delta^{(d)}(R-\eps\Omega)
 F(\eps\Omega)\dd\Omega
\end{equation}
as distributions in $R\in\R^d$.  This identity is applied without fixing
$\Omega$; $\Omega$ remains a live variable in
\eqref{eq:Qsharp-definition}.

The separator is not the time of a root-substituted upper contact.  Put
$h=\epstar/4$ and fix the grid $h\mathbb Z$.  On every
time history choose the smallest grid point $s$ such that
\begin{equation}\label{eq:separator-grid-choice}
 \max_{u\in H_U}t_u+h\le s\le
 \min_{d\in H_D}t_d-h.
\end{equation}
The interval in \eqref{eq:separator-grid-choice} has length at least $2h$,
so such a point exists.  Refine the half-open time cells by the chosen grid
index.  These pieces are disjoint and retain their grid label, so summing
the packet estimate over them does not multiply the right-hand integral by
the number of grid points.  On each refined cell $s$ is a fixed external
number; in particular it has zero derivative under every contact-root or
landing-pivot substitution.  Assumptions (c)--(d) give
\begin{equation}\label{eq:separator-gap}
 t_u\le s-h,\qquad t_d-s\ge h
 \quad(u\in H_U,\ d\in H_D).
\end{equation}
Take the right trace of the hard-sphere state at time $s$.  Every particle
line is free between its last upper atom and this section.  For $i=1,2$,
write $(x_i,v_i)$ for the affine state on the crossing bond $e_i$ and make
the common-section shear
\begin{equation}\label{eq:packet-shear}
 (x_i,v_i,s)\longmapsto(X_i,v_i,s),
 \qquad X_i=x_i+s v_i.
\end{equation}
Its extended Jacobian is one although $s$ is a live time.  The complete
configuration and velocity vectors at the separator will be denoted by
$(X,V)$.

Let $Y_i(t_{d_i};X,V)$ be the companion trajectory arriving at $d_i$ from
the separator, and let $m_i$ be the fixed torus lift.  In the top-pair
representation at $d_i$, the contact residual is
\begin{equation}\label{eq:landing-residual}
 R_i=X_i+(t_{d_i}-s)v_i-Y_i(t_{d_i};X,V)-m_i.
\end{equation}
Put $\Delta_i=t_{d_i}-s$.  The variables $v_1,v_2$ are distinct.  Order the
landings so that $d_1<d_2$.  Conditional on the separator positions and on
the upper history, Lemma~\ref{lem:lower-influence-cone} gives the lower
triangular block
\[
 L_{21}^{(v)}:=D_{v_1}R_2.
\]
\begin{equation}\label{eq:raw-landing-Jacobian}
 D_{(v_1,v_2)}(R_1,R_2)
 =\begin{pmatrix}\Delta_1I_d&0\\
                   L_{21}^{(v)}&\Delta_2I_d\end{pmatrix},
 \qquad
 |\det|=|\Delta_1\Delta_2|^d\ge4^{-2d}\epstar^{2d}.
\end{equation}
The lower-left block includes every possible shared earlier lower history;
it is irrelevant to this conditional determinant.  It is not yet legitimate
to multiply this determinant by a later upper-tree position Jacobian.  The
exterior-coarea argument below controls the resulting Schur row.

If the separator positions and the upper history were external parameters,
one could apply \eqref{eq:sphere-delta-identity} at $d_1,d_2$ and integrate
first in $v_1$ and then in $v_2$.  The conditional substitution would be
\begin{equation}\label{eq:landing-velocity-solution}
 v_i=\Delta_i^{-1}
 \big(\eps\Omega_i-X_i+Y_i(t_{d_i};X,V)+m_i\big),
 \qquad i=1,2,
\end{equation}
where the data on the right for $i=2$ are evaluated after the $i=1$
substitution.  Thus \eqref{eq:landing-velocity-solution} does not assert
that the two companion states are mutually independent; it is the solution
of the triangular system \eqref{eq:raw-landing-Jacobian}.  Its conditional
inverse-Jacobian factor would be
\begin{equation}\label{eq:two-contact-payment}
 \eps^{2(d-1)}|\Delta_1\Delta_2|^{-d}
 \le4^{2d}\eps^{2(d-1)}\epstar^{-2d}.
\end{equation}
The two flux numerators are retained in the positive weight.  On the support
of $Q$ they are at most $2V_\eps$ each.  We do not perform this conditional
substitution before imposing the upper contacts: the joint coarea below
selects a maximal minor among all separator velocities.

\subsection{Exterior coarea and the collision-tree Schur form}

The correct joint object is not a prescribed $2d$-by-$2d$ velocity minor.
It is the $2d$-dimensional row Jacobian with respect to all separator
velocities.  For a linear map $A:E\to\R^k$ put
\begin{equation}\label{eq:normal-jacobian}
 \Jac_k(A):=\sqrt{\det(AA^*)}.
\end{equation}
Thus $\Jac_k(A)^2$ is the sum of the squares of all $k$-column minors in an
orthonormal basis.  A prescribed minor can contract after the upper
translations are eliminated even when the full row volume remains
transverse.

Let $\widehat{\mathcal Q}_{\mathcal P}$ be the hard-sphere configuration
space of the $P=|\mathcal P|$ particle lines, modulo common translation, with
the kinetic-energy metric.  At a regular collision let $\nu$ be the unit
normal pointing into the excluded collision cylinder (hence away from the
admissible exterior), let
$\mathcal K=-D\nu|_{\nu^\perp}\ge0$ be the shape operator of that cylinder,
and let $\mathcal R=I-2\nu\otimes\nu$ be the orthogonal reflection in
$\nu$.  Put \(E=T_X\widehat{\mathcal Q}_{\mathcal P}\).  For the incoming
configuration velocity \(V\in E\), with \(V\cdot\nu\ne0\), define the
event-time projection
\begin{equation}\label{eq:event-time-projection}
 \mathcal V:E\longrightarrow\nu^\perp,\qquad
 \mathcal V\xi
 :=\xi-\frac{\xi\cdot\nu}{V\cdot\nu}\,V .
\end{equation}
Thus \(\mathcal V\xi\) is the displacement at the shifted collision time:
\(\nu\cdot\mathcal V\xi=0\).  Its adjoint
\(\mathcal V^*:\nu^\perp\to E\) is taken in the kinetic metric.  The shape
operator \(\mathcal K:\nu^\perp\to\nu^\perp\) is extended by zero in the
 cylinder-axis and centre-of-mass directions.

For later network elimination this reflection is always kept on the full
reduced space.  If $a=\{i,j\}$, the normalized incidence
$B_a=2^{-1/2}(e_i-e_j)^*\otimes I_d$ gives
\[
 \mathcal R_a=I_E-2B_a^*(\omega_a\otimes\omega_a)B_a,
\]
whose two-line block is exactly
\eqref{eq:appC-pair-dimensional-collision}.  Likewise
$\mathcal H_a=B_a^*\widehat{\mathcal H}_aB_a$.  Thus no step below treats a
binary collision as two independent $d$-dimensional line gauges; the
$d$-dimensional objects are the free-line factor currents and the landing
frame.

\begin{lemma}[Conormal-frame identities]
\label{lem:conormal-frame-identities}
Let the columns of $(Z,W)$ be a $k$-frame of phase-space conormals, with
$Z$ its configuration part and $W$ its velocity part.  Under backward
adjoint pullback through a forward free flight of length $r\ge0$ and a
forward regular collision one has, after
the orthogonal collision identification,
\begin{align}
 (Z,W)&\longmapsto (Z,W+rZ),
 \label{eq:backward-free-conormal}\\
 (Z,W)&\longmapsto
 (\mathcal RZ+\mathcal H\mathcal R W,\mathcal RW),
 \qquad \mathcal H=2|V\cdot\nu|\,
 \mathcal V^*\mathcal K\mathcal V\ge0.
 \label{eq:backward-collision-conormal}
\end{align}
If the symmetric part
$G_{\rm s}:=\frac12(Z^*W+W^*Z)$ is positive semidefinite, these operations
keep it positive semidefinite and do not decrease $\det(W^*W)$.  No symmetry
of $Z^*W$ itself is required.
\end{lemma}

\begin{proof}
For tangent vectors the free map is
$(\delta X,\delta V)\mapsto(\delta X+r\delta V,\delta V)$.
At a collision the exact tangent map is
\[
 \delta X^+=\mathcal R\delta X^-,\qquad
 \delta V^+=\mathcal R\delta V^-+
 \mathcal R\mathcal H\delta X^-.
\]
This is also the full variable-speed saltation matrix used here, not only
its restriction to a unit-energy shell.  Indeed, differentiating the event
equation first gives
\(\delta t=-(\nu\cdot\delta X^-)/(\nu\cdot V^-)\), hence the boundary
displacement is precisely
\(\delta X^-+V^-\delta t=\mathcal V\delta X^-\).  Differentiating
\[
 V^+=V^--2(V^-\cdot\nu)\nu,\qquad
 D\nu[\eta]=-\mathcal K\eta\quad(\eta\in\nu^\perp)
\]
at that shifted event, and then transporting the output back to the fixed
reference time, gives
\[
 D\mathscr S=
 \begin{pmatrix}
  \mathcal R&0\\
  \mathcal R\mathcal H&\mathcal R
 \end{pmatrix},
 \qquad
 \mathcal H=2|V^-\cdot\nu|\mathcal V^*\mathcal K\mathcal V.
\]
To fix all four maps, put
$F_r=\left(\begin{smallmatrix}I&rI\\0&I\end{smallmatrix}\right)$ and
$M=D\mathscr S$.  With the same incoming/outgoing orthogonal
identification, the convention ledger is
\begin{equation}\label{eq:four-map-convention-ledger}
\begin{array}{c|c|c}
 &\text{free flight}&\text{collision}\\ \hline
D\Phi
 &\begin{pmatrix}I&rI\\0&I\end{pmatrix}
 &\begin{pmatrix}\mathcal R&0\\\mathcal R\mathcal H&\mathcal R\end{pmatrix}\\[3mm]
(D\Phi)^{-1}
 &\begin{pmatrix}I&-rI\\0&I\end{pmatrix}
 &\begin{pmatrix}\mathcal R&0\\-\mathcal H\mathcal R&\mathcal R\end{pmatrix}\\[3mm]
(D\Phi)^*
 &\begin{pmatrix}I&0\\rI&I\end{pmatrix}
 &\begin{pmatrix}\mathcal R&\mathcal H\mathcal R\\0&\mathcal R\end{pmatrix}\\[3mm]
(D\Phi)^{-*}
 &\begin{pmatrix}I&0\\-rI&I\end{pmatrix}
 &\begin{pmatrix}\mathcal R&-\mathcal R\mathcal H\\0&\mathcal R\end{pmatrix}.
\end{array}
\end{equation}
The conormal of a future residual, differentiated with respect to the past
state, is pulled back by $(D\Phi)^*$, the third row of
\eqref{eq:four-map-convention-ledger}; this gives
\eqref{eq:backward-free-conormal}--
\eqref{eq:backward-collision-conormal}.
This computation uses no normalization of \(|V^-|\).  Since
\(\mathcal K\ge0\), the displayed definition gives
\(\langle\mathcal H\xi,\xi\rangle
=2|V^-\cdot\nu|\langle\mathcal K\mathcal V\xi,
\mathcal V\xi\rangle\ge0\).  Centre-of-mass and cylinder-axis directions
lie in its kernel and cause no additional term.
Taking adjoints and pulling future conormals back gives
\eqref{eq:backward-free-conormal}--
\eqref{eq:backward-collision-conormal}; these are the matrix form of the
normal-vector identities in
\cite[Section 2, equations (2.2)--(2.4)]{ChernovSimanyi2007}.  Hence, up to
orthogonal conjugacy,
\[
 (G_{\rm s})_{\rm free}=G_{\rm s}+rZ^*Z,
 \qquad
 (G_{\rm s})_{\rm coll}=G_{\rm s}+W^*\mathcal R\mathcal H\mathcal R W.
\]
Both increments are positive semidefinite.  Moreover,
\[
 (W+rZ)^*(W+rZ)-W^*W
 =r(G+G^*)+r^2Z^*Z\ge0,
\]
while a collision replaces $W$ by $\mathcal RW$.  Loewner monotonicity of
the determinant proves the assertion.  A kernel of $\mathcal K$ gives
equality and causes no loss.  The hypothesis is only $G_{\rm s}\ge0$; the
antisymmetric part of $Z^*W$ is unchanged by both displayed increments and
is irrelevant to the velocity-volume conclusion.  A negative direction of
$G_{\rm s}$ is the genuine obstruction for a two-time flag.
\end{proof}

\begin{lemma}[Tree Schur lift of a conormal frame]
\label{lem:tree-schur-conormal-lift}
Assume the upper marked contacts form a collision tree and there is no
unmarked upper contact.  Let $(Z_0,W_0)$ be a $k$-frame of output conormals at
the separator, and let $\mathcal U$ be the upper tree residual.  Add the
unique linear combination of the conormals of $\mathcal U$ which annihilates
all chronological relative-cluster translation columns $h$.  Its remaining
velocity part is
\begin{equation}\label{eq:abstract-tree-schur-lift}
 W_{\rm Sch}=W_0-(D_V\mathcal U)^*(D_h\mathcal U)^{-*}Z_0.
\end{equation}
In chronological tree-cluster coordinates this lift is the composition of
the backward adjoint-pullback free-flight and collision maps
\eqref{eq:backward-free-conormal}--
\eqref{eq:backward-collision-conormal}.  Consequently, if
$\frac12(Z_0^*W_0+W_0^*Z_0)$ is positive semidefinite, then
\begin{equation}\label{eq:tree-schur-volume-monotonicity}
 \sqrt{\det(W_{\rm Sch}^*W_{\rm Sch})}
 \ge c(P_0)\sqrt{\det(W_0^*W_0)}.
\end{equation}
For the landing frame $Z_0=(D_h\mathcal R_D)^*$ and
$W_0=(D_V\mathcal R_D)^*$, one has $W_{\rm Sch}=\mathcal S^*$.
\end{lemma}

\begin{proof}
The algebraic formula is the block elimination identity.  Indeed, if
$\lambda$ is the coefficient frame of upper conormals, annihilation of the
configuration columns reads
\[
 Z_0+(D_h\mathcal U)^*\lambda=0,
 \qquad
 \lambda=-(D_h\mathcal U)^{-*}Z_0,
\]
and substitution in the velocity component gives
\eqref{eq:abstract-tree-schur-lift}.

It remains to identify this single block operation with chronological
conormal transport.  At every prefix of a collision tree, the next tree edge
joins two previously disjoint collision clusters.  Use their relative
cluster translation as the pivot column of that contact.  In these raw
coordinates the contact row has an identity configuration block.  Eliminating
it changes a conormal across the preceding free segment by
$(Z,W)\mapsto(Z,W+rZ)$; across the collision it applies the orthogonal
reflection and, when the normal is differentiated before being retained, the
shape-operator term in \eqref{eq:backward-collision-conormal}.  This is a
direct multiplication of the two local block matrices, not an appeal to an
abstract billiard-flow chart.  Explicitly, one edge step is
\begin{equation}\label{eq:tree-one-event-block-multiplication}
 \binom{Z'}{W'}
 =
 \underbrace{\begin{pmatrix}
  \mathcal R&\mathcal H\mathcal R\\0&\mathcal R
 \end{pmatrix}}_{\text{collision adjoint}}
 \underbrace{\begin{pmatrix}I&0\\rI&I\end{pmatrix}}_{\text{free adjoint}}
 \binom{Z}{W}
 =
 \binom{\mathcal RZ+\mathcal H\mathcal R(W+rZ)}
       {\mathcal R(W+rZ)} .
\end{equation}
In the raw cluster chart the same step is Gaussian elimination of the new
relative-translation column: if the new contact row is
$h_e+L_er$ and a later row is $Y_hh_e+Y_rr$, its Schur row is
\[
 Y_r-Y_hL_e=(Y_h,Y_r)
 \binom{-L_e}{I}.
\]
Differentiating the free flight and scattering relation gives exactly the
two factors in
\eqref{eq:tree-one-event-block-multiplication}; hence the Schur injection
$(-L_e,I)^T$ and the chronological adjoint pullback are the same matrix,
including the causal off-diagonal block.  Induction over the tree edges
proves the
claimed identification.  Passing between cluster translations and
orthonormal reduced particle coordinates has condition number bounded by
$C(P_0)$.

Lemma~\ref{lem:conormal-frame-identities} now preserves the positive cone and
does not decrease the velocity-frame volume at any local step.  The bounded
tree-coordinate changes give \eqref{eq:tree-schur-volume-monotonicity}.
Finally, inserting the landing values of $Z_0,W_0$ in
\eqref{eq:abstract-tree-schur-lift} gives the adjoint of
\eqref{eq:tree-schur-row}.
\end{proof}

\begin{lemma}[Rooted-tree conormal elimination]
\label{lem:rooted-tree-schur-lift}
Let an upper collision word contain a marked chronological spanning tree
$T$ and a set $A_U$ of unmarked regular incoming contacts.  Retain the times
and normals of $T$.  For $a\in A_U$, temporarily retain the $d$ local
variables $y_a=(t_a,\theta_a)$, where $\theta_a$ belongs to a regular normal
chart, and write $C_a=R_a-\eps\omega_a(\theta_a)$.  Jointly eliminate the
rows $C_{A_U}$ against $y_{A_U}$ and the marked-tree rows against their
chronological relative-cluster translations.  If $(Z_0,W_0)$ is a
$k$-frame of separator output conormals satisfying
$\frac12(Z_0^*W_0+W_0^*Z_0)\ge0$, then the velocity part $W_{\rm rt}$ of the
resulting Schur frame satisfies
\begin{equation}\label{eq:rooted-tree-volume-monotonicity}
 \sqrt{\det(W_{\rm rt}^*W_{\rm rt})}
 \ge c(P_0)\sqrt{\det(W_0^*W_0)}.
\end{equation}
Moreover the eliminated root block has exact absolute determinant
\begin{equation}\label{eq:upper-root-neutral-determinant}
 \eps^{(d-1)|A_U|}\prod_{a\in A_U}|g_a\cdot\omega_a|,
\end{equation}
up to the bounded densities of the normal charts.  After solving
$C_{A_U}=0$, $W_{\rm rt}$ is the adjoint of the total derivative obtained
from the physical incoming-root and scattering substitutions; in
particular it includes the derivatives of all substituted normals.
\end{lemma}

\begin{proof}
Use raw edge states and traverse the upper word backward from the separator.
At an unmarked contact write
$G_a=(D_{\theta_a}\omega_a)^*D_{\theta_a}\omega_a$ and
$E_a=D_{\theta_a}\omega_aG_a^{-1/2}$.  Pointwise normalize the tangent
columns by $\delta\widehat\theta_a=G_a^{1/2}\delta\theta_a$.  The local
constraint block in these normalized columns is
\[
 D_{(t_a,\theta_a)}C_a
   =\bigl[g_a\ \ -\eps E_a\bigr],
\]
where the $d-1$ columns of $E_a$ form an orthonormal basis of
$\omega_a^\perp$.  The column normalization has determinant
$J_a=(\det G_a)^{1/2}$, exactly the density in
$\dd\omega_a=J_a\dd\theta_a$; it is therefore retained and cancelled in the
same contact-sphere disintegration.  Resolving the first column along
$\omega_a$ gives
\[
 \left|\det D_{(t_a,\theta_a)}C_a\right|
 =\eps^{d-1}|g_a\cdot\omega_a|.
\]
Causality orders these blocks triangularly.  Hence their product is
\eqref{eq:upper-root-neutral-determinant}; exact grazing is excluded on the
regular cell and has zero collision weight.

The Schur complement of this local block is the derivative after imposing
$C_a=0$ by the implicit-function theorem.  Differentiating the identity
$R_a=\eps\omega_a$ shows that it is the derivative of the physical
incoming collision branch, rather than the fixed-time or fixed-normal
derivative.  On conormals its free-flight and collision factors are exactly
\eqref{eq:backward-free-conormal} and
\eqref{eq:backward-collision-conormal}.  Thus the normal derivative is
already contained in the positive shape-operator term $\mathcal H$; no
normal is frozen in this step.

Write $x$ for all older raw event-state variables,
$y_a=(t_a,\theta_a)$, and $\Phi_a(x,y_a)$ for the state immediately after
the collision.  On a regular sign/lift cell,
\begin{equation}\label{eq:root-IFT-event-matrix}
 D\Phi_a^{\rm phys}
 =
 D_x\Phi_a-D_{y_a}\Phi_a
 (D_{y_a}C_a)^{-1}D_xC_a
 =
 \begin{pmatrix}
  \mathcal R_a&0\\
  \mathcal R_a\mathcal H_a&\mathcal R_a
 \end{pmatrix}.
\end{equation}
The first equality is the literal IFT Schur formula.  The second follows by
differentiating the event time, giving
$\delta t=-(\nu_a\cdot\delta X)/(\nu_a\cdot V^-)$, and then
differentiating the reflected velocity at the shifted contact.  Taking the
adjoint gives
\[
 (D\Phi_a^{\rm phys})^*
 =
 \begin{pmatrix}
  \mathcal R_a&\mathcal H_a\mathcal R_a\\
  0&\mathcal R_a
 \end{pmatrix},
\]
which is exactly the collision factor in
\eqref{eq:tree-one-event-block-multiplication}.  Hence the time, normal and
all older-variable off-diagonal derivatives are accounted for by one
displayed matrix, not by a later identification.

At a marked tree contact, the next tree edge joins two components of the
marked-tree prefix.  Its relative-cluster translation gives an identity
configuration pivot in the raw state immediately adjacent to that contact.
Eliminating this pivot gives the same backward free-flight factor and the
fixed-normal orthogonal collision factor (or the same positive
$\mathcal H$ term if the normal chart is eliminated first).  Interspersed
unmarked contacts are changes between adjacent raw event states and do not
alter this identity pivot.  Induction over the complete chronological word
therefore identifies the one global Schur complement with a composition of
the two maps in Lemma~\ref{lem:conormal-frame-identities}.

That lemma preserves the positive symmetric part of $Z^*W$ and does not
decrease the velocity-frame
volume at every step.  Only the passage between raw cluster translations
and orthonormal reduced particle coordinates remains; its condition number
is bounded by $C(P_0)$.  This proves
\eqref{eq:rooted-tree-volume-monotonicity}.  Finally, the block determinant
identity says that eliminating all $C_a$ is identical to differentiating
after their incoming roots and normals have been substituted, which proves
the last assertion.  With the chronological cluster translations denoted by
$h$, the remaining velocity frame for the landing residual is therefore
exactly the adjoint of
$D_V\mathcal R_D-D_h\mathcal R_D(D_h\mathcal U)^{-1}D_V\mathcal U$.
\end{proof}

\begin{lemma}[Causal endpoint block on a collision-free line]
\label{lem:causal-free-line-endpoints}
Fix the raw separator coordinates before any contact-position row is
eliminated.  Let the physical line $c$ have last upper contact $u$ and first
lower contact $d$, so that it has no collision in $(u,d)$.  In the
separator-side outgoing representation at $u$ and the separator-side
incoming representation at $d$, there are companion trajectories
$Y_u,Y_d$, independent of $(X_c,V_c)$, for which the two vector contact
residuals are
\begin{equation}\label{eq:causal-free-line-residuals}
 \begin{split}
 C_u&=X_c+(u-s)V_c-Y_u-m_u-\eps\omega_u,\\
 C_d&=X_c+(d-s)V_c-Y_d-m_d-\eps\omega_d .
 \end{split}
\end{equation}
Here the lifts $m_u,m_d$ and the normal charts are fixed only as cell
labels; the normals themselves remain live variables.  Consequently
\begin{equation}\label{eq:causal-free-line-block}
 D_{(X_c,V_c)}(C_u,C_d)
 =\begin{pmatrix}I_d&(u-s)I_d\\ I_d&(d-s)I_d\end{pmatrix},
 \qquad
 \left|\det D_{(X_c,V_c)}(C_u,C_d)\right|=(d-u)^d.
\end{equation}
The assertion remains true when the companion at $d$ has passed through
arbitrary lower collisions strictly before $d$, provided none of those
collisions contains $c$.
\end{lemma}

\begin{proof}
Starting from the separator, reconstruct the upper word backward.  Between
$u$ and $s$ no collision contains $c$.  Hence every elastic reset and free
transport used to reconstruct the companion slot at $u$ acts only on
separator states other than $(X_c,V_c)$.  Induction over those contacts
shows that $Y_u$ is independent of $(X_c,V_c)$.  The
separator-side outgoing pair is essential here: the incoming velocity pair
at $u$ is obtained only after the collision involving $c$ has been inverted.
The contact position, however, is common to the incoming and outgoing
representatives.

Reconstruct the lower word forward from the separator.  Since $d$ is the
first lower contact of $c$, no earlier lower reset contains that line.
Exactly the same causal induction, now in increasing time, shows that the
companion trajectory $Y_d$ is independent of $(X_c,V_c)$.  This also proves
the final assertion and is the one-line version of
Lemma~\ref{lem:lower-influence-cone}.

On the fixed lift cell, the position on $c$ is the affine function
$X_c+(t-s)V_c$.  This gives \eqref{eq:causal-free-line-residuals}; subtracting
the first block row from the second gives
$\operatorname{diag}(I_d,(d-u)I_d)$ and proves
\eqref{eq:causal-free-line-block}.  No collision constraint has yet been
used in this differentiation, so the asserted independence is an identity
in genuinely free raw coordinates rather than an independence claimed
after a Schur substitution.
\end{proof}

\begin{lemma}[Positive free-chord reduction]
\label{lem:positive-free-chord-reduction}
Let one physical particle line have no collision in the open interval
$(u,d)$, $u<d$, and let its contact positions at the two endpoints be
constrained to companion ports $Y_u,Y_d$.  For any separator time
$s\in(u,d)$ its state $(X,V)$ is recovered from the two vector residuals
with
\begin{equation}\label{eq:free-chord-recovery}
 V=\frac{Y_d-Y_u}{d-u},
 \qquad
 X=\frac{(d-s)Y_u+(s-u)Y_d}{d-u},
\end{equation}
and the absolute Jacobian of those residuals with respect to $(X,V)$ is
$(d-u)^d$.

At the tangent level, eliminating $(X,V)$ is a positive Schur reduction.  In
the discrete Jacobi index form it replaces the free line by the port term
\begin{equation}\label{eq:free-chord-index-term}
 \frac1{d-u}|\xi_d-O_{u,d}\xi_u|^2\ge0,
\end{equation}
Here $O_{u,d}$ is the ordered orthogonal transport identifying the two
companion-line coordinate spaces.  In the complete reduced space the term is
the factor square
\[
 (d-u)^{-1}|A_\gamma\xi|^2
\]
from \eqref{eq:appC-network-form}, not a full-rank edge on $E$.  Hence
adjoining and then eliminating such a free chord preserves the
positive-index-form hypothesis in
Lemma~\ref{lem:rooted-tree-schur-lift}, including for a forest whose
components are grounded at the two chord endpoints.
\end{lemma}

\begin{proof}
In one common affine coordinate system the two endpoint equations are
$X+(u-s)V=Y_u$ and $X+(d-s)V=Y_d$.  Subtraction gives
\eqref{eq:free-chord-recovery} and the block determinant is $(d-u)^d$.
For a tangent endpoint field the recovered field on the open chord is the
affine interpolant
\[
 \xi(t)=\frac{d-t}{d-u}\xi_u+\frac{t-u}{d-u}\xi_d.
\]
It uniquely minimizes the free Jacobi energy
$\int_u^d|\dot\eta(t)|^2\,dt$ among fields with the same endpoints, and its
energy is $(d-u)^{-1}|\xi_d-\xi_u|^2$.  Collision identifications at the
ports are orthogonal, giving \eqref{eq:free-chord-index-term}.  Schur
complements of nonnegative quadratic forms are nonnegative on the grounded
range, which proves the last assertion.
\end{proof}

\begin{lemma}[Grounded positive-network one-port lift]
\label{lem:grounded-positive-network-lift}
Consider a finite time-expanded hard-sphere Jacobi network on at most
$P_0$ physical lines and with at most $P_0$ events.  Its local ingredients are:
\begin{enumerate}[label=\textup{(\roman*)}]
\item physical-line free segments of lengths $\ell_e>0$, whose
  $d$-dimensional drop maps $A_e$ are obtained after the adjacent pair
collision junctions are inserted;
\item fixed-normal elastic junctions with the full two-line block
\(\mathscr R_a^{(2d)}\) in
\eqref{eq:appC-pair-dimensional-collision};
\item incoming-root or collision-curvature shunts $\mathcal H_v\ge0$;
\item at most one already eliminated free chord, contributing the rank-$d$ factor
 $\ell^{-1}|A_\gamma\xi|^2$.
\end{enumerate}
Let a collision forest be grounded by fixing one translation in every
component.  Assume that the resulting rooted endpoint word is
trace-complete in the sense of
Definition~\ref{def:appC-trace-complete-word}, and suppose its reduced
incidence pivots are then eliminated
together with all regular incoming time--normal roots.  Let
$K_{\mathcal N}$ be the current-to-potential compliance of the resulting
grounded network.  There are dual covariant current/potential gauges
$Q_-,Q_+$, consisting only of full-pair orthogonal junctions and bounded
forest/incidence coordinate maps, such that
\begin{equation}\label{eq:grounded-network-compliance-factorization}
 K_{\mathcal N}=K_{\mathcal N}^*\ge0,
 \qquad
 s_{\min}(Q_\pm)\ge c(P_0),\quad
 \|Q_\pm\|\le C(P_0).
\end{equation}
Consequently, if the input frame is the one-speed frame $(Z,bZ)$, $b>0$,
then
\begin{equation}\label{eq:grounded-one-speed-volume}
 \widetilde K_{\mathcal N}:=
 Q_+K_{\mathcal N}Q_-
 =Q_-^*K_{\mathcal N}Q_-\ge0,
 \qquad
 W_{\rm Sch}=(bI+\widetilde K_{\mathcal N})Z,
 \qquad
 \Jac_k(W_{\rm Sch})\ge b^k\Jac_k(Z).
\end{equation}
This statement allows the endpoints of a positive chord to lie on either
grounded component; no assertion that a preferred raw bridge minor is a
positive update is used.
\end{lemma}

\begin{proof}
We give the finite-dimensional elimination argument because positivity of
the final network cannot be inferred merely from a diagram.  Split every
physical trajectory at all contact ports and put one tangent field $\xi_v$
at each resulting vertex.  After the orthogonal endpoint identifications,
the part of the discrete index form which has already been eliminated is
\begin{equation}\label{eq:grounded-network-index-form}
 \mathcal I(\xi)=
 \frac12\sum_e\ell_e^{-1}|A_e\xi|^2
 +\frac12\sum_v\langle\mathcal H_v\xi_v,\xi_v\rangle
 +\frac12\sum_\gamma
   \ell_\gamma^{-1}|A_\gamma\xi|^2 .
\end{equation}
Every summand is nonnegative.  The last sum is precisely the contribution
of the free chords in item~\textup{(iv)}, not a new collision assumption.

Choose one ground in every forest component and order the remaining port
fields in the same leaf/chronological order as the reduced-incidence
pivots.  If $I$ denotes a set of fields about to be eliminated and $B$ the
fields still on the boundary, the matrix of
\eqref{eq:grounded-network-index-form} has the block form
\[
 L=\begin{pmatrix}L_{II}&L_{IB}\\L_{BI}&L_{BB}\end{pmatrix}\ge0.
\]
The grounded forest part makes $L_{II}$ positive definite on the columns
actually pivoted.  Indeed, the chronological forest injection
\eqref{eq:appC-forest-coordinate-injection} has no grounded kernel, and every
chord factor is a nonnegative additional row; hence the extra-kernel
hypothesis of Lemma~\ref{lem:appC-grounded-invertibility} holds.  Gaussian
elimination replaces $L$ by its one-port Kron
stiffness
\begin{equation}\label{eq:grounded-kron-block}
 S_{\mathcal N}:=L/L_{II}
 =L_{BB}-L_{BI}L_{II}^{-1}L_{IB}\ge0,
\end{equation}
because, for every boundary field $z$,
\[
 \frac12\langle z,S_{\mathcal N}z\rangle
 =\min_y\mathcal I(y,z)\ge0.
\]
If a harmless zero mode is retained before grounding, the same formula with
the Moore--Penrose inverse holds on the grounded range; our chosen grounds
remove that mode and keep the displayed inverse ordinary.  In fact the
retained one-port block is positive definite.  The physical trace estimate
\eqref{eq:appC-boundary-trace-coercivity} uses the actual factor rows,
forest constraints, full pair junctions and grounds and proves that a
zero-energy minimizer must have zero retained boundary datum, even when the
chord is only rank $d$ and the curvature shunts have kernels.  Hence
$S_{\mathcal N}>0$.

We next identify the operator produced by this elimination.  Prescribe a
boundary current $j$ and solve the grounded Euler--Lagrange equations for
\eqref{eq:grounded-network-index-form}.  Equivalently, minimize the Thomson
energy
\begin{equation}\label{eq:grounded-network-thomson}
 \frac12\sum_e\ell_e|J_e|^2
 +\frac12\sum_v\langle\mathcal H_v^+q_v,q_v\rangle
\end{equation}
over covariantly transported edge currents with boundary divergence $j$;
here $q_v\in\operatorname{Ran}\mathcal H_v$ is the current leaving through
the shunt and its kernel carries no shunt current.  Strict positivity on the grounded range
gives a unique minimizing potential modulo no residual translation.  Its
boundary value is $K_{\mathcal N}j$, where $K_{\mathcal N}$ is the grounded
compliance, equivalently $K_{\mathcal N}=S_{\mathcal N}^{-1}$.  The
polarization identity for the minimum in
\eqref{eq:grounded-network-thomson} gives
\begin{equation}\label{eq:grounded-network-compliance-positive}
 \frac12\langle j,K_{\mathcal N}j\rangle
 =\min_{\operatorname{div}J=j}
   \frac12\left(\sum_e\ell_e|J_e|^2+\text{shunt energy}\right)\ge0.
\end{equation}
Thus $K_{\mathcal N}=K_{\mathcal N}^*\ge0$.  The stiffness and compliance
are not interchanged: a single edge of length $\ell$ has
$S_{\mathcal N}=\ell^{-1}I$ and $K_{\mathcal N}=\ell I$.  Formula
\eqref{eq:grounded-kron-block}, applied to the common full block matrix,
shows that eliminating any subset of internal ports gives the same
stiffness and hence the same compliance.  In particular, the order
``chord first, forest second'' agrees with the simultaneous Schur
complement; Lemma~\ref{lem:appC-chord-associativity} writes the exact
Schur-quotient identity.

Theorem~\ref{thm:appC-physical-schur-jacobi}, applied to the complete event
word before any quotient, identifies the two operators through the fully
defined action \eqref{eq:appC-physical-discrete-action} and the chronological
Green induction \eqref{eq:appC-word-Green-induction}.  Formula
\eqref{eq:appC-root-local-Schur-check} inserts every causal off-diagonal
derivative of an interspersed root into the later factor rows and the same
embedded shunt.  All bordered blocks are the explicit physical blocks displayed above.

The clean boundary calculation is made in current variables, so no
historical chord endpoint is expressed through an instantaneous phase
potential.  Prescribe the raw landing current $z_B$ and let
$j_B=Q_-z_B$ be its covariant representative.  The lower free interval
contributes $b|z_B|^2/2$.  All remaining physical segment currents,
the chord current and the nonnegative root-port currents have the Thomson
energy
\[
 \frac12\sum_e\ell_e|J_e|^2
 +\frac12\sum_v\langle H_v^+q_v,q_v\rangle
\]
under exact full-pair current conservation and boundary divergence $j_B$.
By \eqref{eq:grounded-network-compliance-positive}, its minimum is
\[
 \frac12\langle j_B,K_{\mathcal N}j_B\rangle.
\]
Therefore the complete raw dual boundary action is
\begin{equation}\label{eq:grounded-Thomson-one-speed-action}
 \mathcal T_b^{\rm raw}(z_B)=\frac12
 \left(b|z_B|^2+\langle Q_-z_B,K_{\mathcal N}Q_-z_B\rangle\right)
 =\frac12\langle z_B,(bI+\widetilde K_{\mathcal N})z_B\rangle .
\end{equation}
The discrete Green identity says that its derivative is exactly the
retained boundary potential.  This is the physical Schur landing conormal,
not merely a comparison quadratic form; Theorem
~\ref{thm:appC-physical-schur-jacobi} proves the equality before and after
the same Gaussian quotient.

Let $Q_-$ send the raw landing-current frame to the covariant current in
\eqref{eq:grounded-Thomson-one-speed-action}, and let $Q_+$ send the dual
covariant potential back to the raw retained Schur row.  If $P_{\mathcal N}$
is the square primal raw-to-covariant potential map, then, with the convention
of \eqref{eq:appC-exact-dual-boundary-gauges},
\[
 Q_-=P_{\mathcal N}^{-*},\qquad
 Q_+=P_{\mathcal N}^{-1}=Q_-^*,\qquad
 \langle Q_-z,\xi\rangle_{\rm cov}
   =\langle z,Q_+\xi\rangle_{\rm raw}.
\]
Their
only nonorthogonal factors are the normalized incidence and forest QR maps;
the bounds in
\eqref{eq:appC-forest-coordinate-injection}--
\eqref{eq:appC-forest-boundary-splitting} therefore give
\[
 s_{\min}(Q_\pm)\ge c(P_0),\qquad
 \|Q_\pm\|\le C(P_0).
\]
No free-segment or chord length occurs in $Q_\pm$: every such length is in
the positive Thomson functional.  Differentiating
\eqref{eq:grounded-Thomson-one-speed-action} now gives the identity in
\eqref{eq:grounded-one-speed-volume}.  Since
$\widetilde K_{\mathcal N}=Q_-^*K_{\mathcal N}Q_-\ge0$, the least singular
value of $bI+\widetilde K_{\mathcal N}$ is at least $b$, which proves the
displayed exterior estimate.  This also proves
\eqref{eq:grounded-network-compliance-factorization}.  The chord uses its
exact selector and unchanged physical weight; no one-particle projection
or reciprocal endpoint coordinate is used.
\end{proof}

\begin{proposition}[Asymmetric new-line peeling]
\label{prop:asymmetric-new-line-peeling}
Assume \textup{(a)--(e)} of Theorem~\ref{thm:two-landing}, and suppose the
first two lower edges overlap.  Write
\[
 f_{d_1}=\{a,b\},\qquad f_{d_2}=\{b,c\},
\]
after relabelling, so $c$ is the line which first appears at $d_2$.  Allow
arbitrary unmarked upper contacts and arbitrary lower contacts strictly
after $d_2$.  Then, without assuming $\mathrm{FCT}(P_0)$ or a bound on
$d_2-d_1$, every nonnegative joint kernel obeys
\begin{equation}\label{eq:asymmetric-new-line-operator}
 J_{H,\kappa}(Q)
 \le \cL_\eps^{C|H|}\eps^{d-1}\epstar^{-2d}
       \int Q_{H,\kappa}^{\sharp,\mathrm{ch}}.
\end{equation}
The output box is independent of the concrete fixed boundary values.
\end{proposition}

\fallbackproofinresource{supp:proof:prop:asymmetric-new-line-peeling}

\begin{corollary}[Positive lower-flag criterion]
\label{cor:positive-lower-flag}
Pull the $2d$ conormals of the two selected landing-position residuals
back to the separator through the exact lower collision history up to
$d_2$, solving every unmarked incoming contact on the way, and denote their
configuration and velocity parts by $(Z_D,W_D)$.  Suppose
\begin{equation}\label{eq:positive-lower-flag-hypothesis}
 \frac12(Z_D^*W_D+W_D^*Z_D)\ge0,
 \qquad
 \sqrt{\det(W_D^*W_D)}
 \ge c_L|\Delta_1\Delta_2|^d.
\end{equation}
Then arbitrary marked upper trees and arbitrary unmarked upper contacts
satisfy the unconditional packet estimate
\begin{equation}\label{eq:positive-lower-flag-operator}
 J_{H,\kappa}(Q)
 \le C(P_0)c_L^{-1}\cL_\eps^C
       \eps^{d-1}\epstar^{-2d}\int Q_{H,\kappa}^{\sharp,L}.
\end{equation}
Every unmarked lower contact used in forming $(Z_D,W_D)$ is neutral in the
measure through its own sphere--flux/root determinant.
\end{corollary}

\begin{proof}
Apply Lemma~\ref{lem:rooted-tree-schur-lift} to the frame
$(Z_D,W_D)$.  Equations
\eqref{eq:rooted-tree-volume-monotonicity} and
\eqref{eq:positive-lower-flag-hypothesis} give a $2d$-row separator
velocity Jacobian at least
$c(d,P_0)c_L|\Delta_1\Delta_2|^d$.  Cauchy--Binet chooses a $2d$-coordinate
minor.  The local determinant calculation
\eqref{eq:upper-root-neutral-determinant}, applied in the same chronological
raw chart on both sides of the separator, cancels every unmarked sphere and
flux factor.  The remaining sphere/normalization ledger is
\[
 \eps^{-(d-1)P}\eps^{(d-1)(P-1)}\eps^{2(d-1)}=\eps^{d-1}.
\]
Since $\Delta_i\ge h=\epstar/4$, joint coarea costs at most
$C(P_0)c_L^{-1}\epstar^{-2d}$.  Zero extension on the fixed half-open cell
defines $Q_{H,\kappa}^{\sharp,L}$ and preserves the arbitrary nonnegative
kernel throughout.
\end{proof}

For an oriented particle-line edge $f=(p,q)$ put
\begin{equation}\label{eq:incidence-covector}
 \beta_f x=x_p-x_q,\qquad
 b_f=\beta_f\tensor I_d:\R^{dP}\longrightarrow\R^d.
\end{equation}
Two different unordered edges give independent incidence rows.  Their
scalar Gram matrix is either
$\left(\begin{smallmatrix}2&\pm1\\\pm1&2\end{smallmatrix}\right)$ or
$2I_2$.

\begin{definition}[Short first-lower-gap branch]
\label{def:alternative-A}
After selecting the first two lower atoms \(d_1<d_2\), we say that
\emph{alternative \textup{(A)}} holds when
\begin{equation}\label{eq:alternative-A}
 0<t_{d_2}-t_{d_1}<\epstar .
\end{equation}
Its complement is the long first-lower-gap branch
\(t_{d_2}-t_{d_1}\ge\epstar\).  This terminology is used only in
Lemmas~\ref{lem:consecutive-lower-positive-flag}--\ref{cor:no-double-packet-direct};
the unconditional proof of Theorem~\ref{thm:two-landing} treats both
branches by the disjoint/overlap dichotomy and does not assume
alternative~\textup{(A)}.
\end{definition}

\begin{lemma}[Two consecutive lower landings form a positive flag]
\label{lem:consecutive-lower-positive-flag}
Assume alternative \textup{(A)} of Definition~\ref{def:alternative-A}, and
assume that $d_1<d_2$ are the only lower contacts up to time $d_2$.  Put
\[
 b=\Delta_1=d_1-s,\qquad q=d_2-d_1,
 \qquad B=\Delta_2=b+q.
\]
If $s$ is the smallest separator-grid point in
\eqref{eq:separator-grid-choice}, then
\begin{equation}\label{eq:lower-gap-ratio}
 b\ge2h,\qquad 0<q<4h,\qquad q/b<2.
\end{equation}
Let $(Z_D,W_D)$ be the $2d$ landing conormals pulled back from $d_1,d_2$
to $s$, including the fixed-normal reflection at $d_1$.  For arbitrary
distinct lower collision edges and arbitrary regular normals,
\begin{align}
 \frac12(Z_D^*W_D+W_D^*Z_D)&\ge c_d\,b I_{2d},
 \label{eq:lower-flag-symmetric-cone}\\
 \sqrt{\det(W_D^*W_D)}&\ge c_d(bB)^d.
 \label{eq:lower-flag-volume}
\end{align}
The constants are absolute and no relation between the two landing normals
is required.
\end{lemma}

\begin{proof}
Let $u_*:=\max_{u\in H_U}t_u$.  Since there is no lower contact before
$d_1$, the inter-sublayer gap is $d_1-u_*\ge\epstar=4h$.  The smallest grid
point satisfying $u_*+h\le s$ obeys $s<u_*+2h$, and hence $b>2h$ (the weak
inequality in \eqref{eq:lower-gap-ratio} is enough).  The two contacts are
consecutive lower atoms, so alternative \textup{(A)} gives
$q<\epstar=4h$.  This proves
\eqref{eq:lower-gap-ratio}.

Normalize an incidence row by $a_f=2^{-1/2}b_f$, so
$a_fa_f^*=I_d$.  Let $\mathcal R_1$ be the orthogonal configuration
reflection at $d_1$.  For output covectors $y_1,y_2\in\R^d$, set
\[
 z_1=a_{f_1}^*y_1,
 \qquad z_2=\mathcal R_1a_{f_2}^*y_2.
\]
Before a row operation, direct differentiation at fixed $d_1,d_2$ and
fixed retained normal $\omega_1$ gives the two conormal columns
\[
 (a_{f_1}^*y_1,\,b a_{f_1}^*y_1),
 \qquad
 (a_{f_2}^*y_2,\,(bI+q\mathcal R_1)a_{f_2}^*y_2).
\]
The difference
$\mathcal R_1a_{f_2}^*y_2-a_{f_2}^*y_2$ lies in
$\operatorname{Ran}a_{f_1}^*$, because
$\mathcal R_1-I=-2(a_{f_1}^*\omega_1)
(a_{f_1}^*\omega_1)^*$.  Add the corresponding multiple of the first $d$
landing rows to the second $d$.  This is a block-unitriangular row
operation of determinant one, and it changes the second pair above to
$(z_2,Bz_2)$.  Thus the row-equivalent exact frame is
\begin{equation}\label{eq:two-lower-frame-formula}
 Z_Dy=z_1+z_2,
 \qquad W_Dy=bz_1+Bz_2,
 \qquad y=(y_1,y_2).
\end{equation}

If the two lower edges are disjoint, then
$a_{f_1}\mathcal R_1a_{f_2}^*=0$.  If they share one endpoint, their
oriented scalar incidence product is $\sigma/2$, $\sigma\in\{-1,1\}$.
Writing
$\nu_1=a_{f_1}^*\omega_1$ and
$\mathcal R_1=I-2\nu_1\nu_1^*$ gives the exact cross block
\begin{equation}\label{eq:lower-incidence-reflection-cross}
 a_{f_1}\mathcal R_1a_{f_2}^*
 =\frac{\sigma}{2}(I_d-2\omega_1\omega_1^*).
\end{equation}
Thus this block is self-adjoint and all its eigenvalues $\lambda$ satisfy
$|\lambda|\le1/2$.

Diagonalize \eqref{eq:lower-incidence-reflection-cross}.  In each of its
$d$ scalar eigendirections, the quadratic form associated with the
symmetric part of $Z_D^*W_D$ is
\begin{equation}\label{eq:lower-cone-two-by-two}
 M_\lambda=
 \begin{pmatrix}
  b&\lambda(b+q/2)\\
  \lambda(b+q/2)&b+q
 \end{pmatrix}.
\end{equation}
With $r=q/b<2$,
\[
 \frac{\det M_\lambda}{b^2}
 \ge 1+r-\frac14(1+r/2)^2\ge\frac34,
 \qquad \operatorname{tr}M_\lambda\le4b.
\]
Hence $M_\lambda\ge(3/16)bI_2$, proving
\eqref{eq:lower-flag-symmetric-cone} in normalized incidence coordinates.

The corresponding block of $W_D^*W_D$ is
\[
 N_\lambda=
 \begin{pmatrix}
  b^2&\lambda bB\\
  \lambda bB&B^2
 \end{pmatrix},
 \qquad
 \det N_\lambda=b^2B^2(1-\lambda^2)
 \ge\frac34b^2B^2.
\]
Multiplying the $d$ blocks and taking a square root gives
$\sqrt{\det(W_D^*W_D)}\ge(3/4)^{d/2}(bB)^d$ in the normalized frame.
Returning from $a_f$ to $b_f$ changes only an absolute constant and proves
\eqref{eq:lower-flag-volume}.
\end{proof}

\begin{proposition}[All rooted words with a short first lower gap]
\label{prop:no-early-lower-root}
Under assumptions (a)--(e) and alternative \textup{(A)} of
Definition~\ref{def:alternative-A}, arbitrary
unmarked upper contacts and arbitrary lower contacts strictly after $d_2$
are allowed.  The first-two choice leaves no lower contact before $d_2$.
Then every bounded packet, for arbitrary regular collision
normals and arbitrary nonnegative $Q$, obeys
\begin{equation}\label{eq:no-early-lower-root-operator}
 J_{H,\kappa}(Q)
 \le \cL_\eps^C\eps^{d-1}\epstar^{-2d}
       \int Q_{H,\kappa}^{\sharp,0D}.
\end{equation}
The output domain is fixed by the half-open cell and not by its boundary
values.
\end{proposition}

\begin{proof}
Lemma~\ref{lem:consecutive-lower-positive-flag} gives the positive symmetric
cone and the $2d$-volume lower bound for the exact lower landing frame at
the separator.  Apply the rooted-tree elimination of
Lemma~\ref{lem:rooted-tree-schur-lift}.  The strengthened conormal identity
uses only the positive symmetric part, so the possibly nonzero
antisymmetric pairing between the two landing times causes no loss.  It
follows that
\[
 \Jac_{2d}(\mathcal S)\ge c(d,P_0)(\Delta_1\Delta_2)^d
 \ge c(d,P_0)h^{2d}.
\]
Every upper root contributes the neutral determinant
\eqref{eq:upper-root-neutral-determinant}, and later lower roots are absent
from the selected residuals by causality and are neutral by the same local
calculation.  Cauchy--Binet, joint coarea, and the sphere/normalization
identity
\[
 \eps^{-(d-1)P}\eps^{(d-1)(P-1)}\eps^{2(d-1)}=\eps^{d-1}
\]
therefore give \eqref{eq:no-early-lower-root-operator}.  All substitutions
remain inside one nonnegative integrand; zero extension on the retained
cell defines $Q_{H,\kappa}^{\sharp,0D}$.
\end{proof}

\begin{corollary}[Direct no-double packet bound]
\label{cor:no-double-packet-direct}
Every packet satisfying assumptions (a)--(e) and alternative \textup{(A)}
of Definition~\ref{def:alternative-A} obeys
\begin{equation}\label{eq:no-double-packet-direct}
 J_{H,\kappa}(Q)
 \le \cL_\eps^{C|H|}\eps^{d-1}\epstar^{-2d}
       \int Q_{H,\kappa}^{\sharp}
\end{equation}
for every nonnegative joint kernel $Q$.  This conclusion is independent of
$\mathrm{FCT}(P_0)$.
\end{corollary}

\begin{proof}
Choose the crossing lines by Lemma~\ref{lem:upper-tree-landings}.  Their
lower endpoints are the first two lower atoms, so the hypothesis of
Proposition~\ref{prop:no-early-lower-root} holds.  That proposition gives
\eqref{eq:no-double-packet-direct}; its recovery map is precisely the main
chart fixed in Definition~\ref{def:packet-output}.
\end{proof}

\begin{definition}[Flux-completed collision-tree condition $\mathrm{FCT}(P_0)$]
\label{def:collision-tree-coercivity}
Consider a regular hard-sphere collision word on $P\le P_0$ labelled
particle lines.  Before a separator time $s$, mark a chronological spanning
tree $T$ of collision contacts.  After $s$, mark two first-landing contacts
$d_1<d_2$ whose unordered particle pairs are different.  Other contacts may
be arbitrary co-tree contacts.  Fix the order, lifts, incoming signs and
half-open normal charts.  Retain the times and normals of the upper tree and
of $d_1,d_2$.  Every other contact is first replaced by its unique fixed-lift
incoming time root as in \eqref{eq:remaining-time-root}; its contact normal
and scattering map are then functions of the variables that remain.  All
derivatives below are the \emph{total derivatives after these root
substitutions}.  In particular they include the derivatives of every
substituted root and normal.  Put $\Delta_i=t_{d_i}-s>0$.

Let $V$ be the separator velocities and let $h\in\R^{d(P-1)}$ be the
chronological relative-cluster translation coordinates used in
Lemma~\ref{lem:rooted-tree-schur-lift}; the common translation is omitted.
Let $\mathcal U$ be the $d(P-1)$ marked-tree residuals and let
$\mathcal R_D=(R_1,R_2)$.  In the raw chronological chart
$D_h\mathcal U$ is block triangular with reduced tree-incidence diagonal,
and hence it and its inverse are bounded in terms of $P_0$.  Define
\begin{equation}\label{eq:tree-schur-row}
 \mathcal S:=D_V\mathcal R_D-
 D_h\mathcal R_D(D_h\mathcal U)^{-1}D_V\mathcal U.
\end{equation}
Let $\vartheta=(t_{d_1},t_{d_2})$.  Causality makes $\mathcal U$
independent of these two future times, so the reduced landing-time matrix is
\[
 G_{21}^{(t)}:=D_{t_{d_1}}R_2.
\]
\begin{equation}\label{eq:landing-time-columns}
 \mathcal T:=D_\vartheta\mathcal R_D
 =\begin{pmatrix}g_1&0\\G_{21}^{(t)}&g_2\end{pmatrix},
 \qquad \phi_i:=|g_i\cdot\omega_i|,
\end{equation}
where $g_i$ and $\omega_i$ are the incoming relative velocity and contact
normal at $d_i$.  The derivatives in \eqref{eq:landing-time-columns}, like
those in \eqref{eq:tree-schur-row}, are total derivatives after all earlier
root substitutions.

For every $2(d-1)$-element set $I$ of scalar separator-velocity coordinates,
write $\mathcal S_I$ for the corresponding $2d$-by-$2(d-1)$ submatrix and put
\begin{equation}\label{eq:mixed-normal-jacobian}
 \Jac_{2(d-1)\mid2}(\mathcal S,\mathcal T)
 :=\left(\sum_{|I|=2(d-1)}
       |\det[\mathcal S_I\ \mathcal T]|^2\right)^{1/2}.
\end{equation}
We say that $\mathrm{FCT}(P_0)$ holds if, for every such word on its regular
physical collision fibre,
\begin{equation}\label{eq:tree-exterior-coercivity}
 \Jac_{2(d-1)\mid2}(\mathcal S,\mathcal T)
 \ge c(d,P_0)\phi_1\phi_2|\Delta_1\Delta_2|^{d-1}
\end{equation}
for some $c(d,P_0)>0$ uniform in regular collision normals and all
within-layer time gaps.  All co-tree collisions are included through these
implicit root/scattering maps; none is frozen while differentiating
$\mathcal U$ or $\mathcal R_D$.
\end{definition}

\begin{lemma}[Exact time-wedge factorization]
\label{lem:time-wedge-factorization}
On a regular cell let $E_T=\operatorname{Ran}\mathcal T\subset\R^{2d}$, and
let $Q_T:\R^{2(d-1)}\to E_T^\perp$ be any orthogonal identification.  Then
\begin{equation}\label{eq:time-wedge-factorization}
 \Jac_{2(d-1)\mid2}(\mathcal S,\mathcal T)
 =\sqrt{\det(\mathcal T^*\mathcal T)}\,
   \Jac_{2(d-1)}(Q_T^*\mathcal S).
\end{equation}
Moreover the causal block form \eqref{eq:landing-time-columns} gives
\begin{equation}\label{eq:time-wedge-flux-lower}
 \sqrt{\det(\mathcal T^*\mathcal T)}\ge\phi_1\phi_2.
\end{equation}
In particular, $\mathrm{FCT}(P_0)$ follows from the flux-free
$2(d-1)$-dimensional spacetime-tangential estimate
\begin{equation}\label{eq:quotient-fct}
 \Jac_{2(d-1)}(Q_T^*\mathcal S)
 \ge c(d,P_0)|\Delta_1\Delta_2|^{d-1}.
\end{equation}
\end{lemma}

\begin{proof}
Choose an orthogonal matrix on the landing-row space whose first $2(d-1)$ rows
span $E_T^\perp$ and whose last two rows span $E_T$.  For every
$2(d-1)$-set $I$,
the transformed matrix $[\mathcal S_I\ \mathcal T]$ is block triangular, so
its determinant is the product of the determinant of
$Q_T^*\mathcal S_I$ and the two-dimensional determinant of $\mathcal T$ in
$E_T$.  Squaring, summing over $I$, and applying Cauchy--Binet gives
\eqref{eq:time-wedge-factorization}.

For \eqref{eq:time-wedge-flux-lower}, use the orthogonal output covectors
$(\omega_1,0)$ and $(0,\omega_2)$.  The corresponding two-by-two minor of
\eqref{eq:landing-time-columns} is lower triangular with diagonal entries
$g_1\cdot\omega_1$ and $g_2\cdot\omega_2$.  Its modulus is
$\phi_1\phi_2$, and the full two-row volume is at least the modulus of this
minor.  Combining the two identities proves the sufficient criterion.
\end{proof}

\begin{lemma}[Exact two-speed landing quotient]
\label{lem:two-speed-landing-quotient}
Assume that $d_1<d_2$ are the first two lower contacts and that their
unordered particle edges are different.  Put
\[
 b=t_{d_1}-s,\qquad B=t_{d_2}-s.
\]
After a block-unitriangular landing-row operation with determinant one and
uniformly bounded inverse, the two landing-time columns and the lower
conormal frame have the simultaneous form
\begin{equation}\label{eq:two-speed-time-diagonal}
 \widetilde{\mathcal T}
 =\begin{pmatrix}g_1&0\\0&g_2\end{pmatrix},\qquad
 \widetilde Z=(Z_1,Z_2),\qquad
 \widetilde W=(bZ_1,BZ_2).
\end{equation}
Here $Z_i:\R^d\to\widehat{\mathcal Q}_{\mathcal P}$ is an isometry in
normalized incidence coordinates, up to the harmless common convention for
the kinetic metric, and
\begin{equation}\label{eq:two-speed-cross-contraction}
 \|Z_1^*Z_2\|\le\frac12.
\end{equation}
If $K_i:\R^{d-1}\to g_i^\perp$ are orthogonal identifications, put
\[
 \overline Z=(Z_1K_1,Z_2K_2),\qquad
 \overline W=(bZ_1K_1,BZ_2K_2).
\]
Then
\begin{equation}\label{eq:two-speed-graph-form}
 \overline W=\overline Z
 \begin{pmatrix}bI_{d-1}&0\\0&BI_{d-1}\end{pmatrix},
 \qquad
 \frac12I_{2(d-1)}\le \overline Z^*\overline Z
 \le\frac32I_{2(d-1)},
\end{equation}
and consequently
\begin{equation}\label{eq:two-speed-initial-volume}
 \sqrt{\det(\overline W^*\overline W)}
 \ge 2^{-(d-1)}(bB)^{d-1}.
\end{equation}
Thus the flux-completed $2(d-1)$-dimensional landing quotient is uniformly
transverse at the separator for every first gap; the unresolved issue in a
long-gap overlapping word is solely whether the complete rooted upper
elimination preserves a fixed fraction of the two-speed volume
\eqref{eq:two-speed-initial-volume}.
\end{lemma}

\begin{proof}
Use the normalized incidence maps $a_{f_i}$ from
\eqref{eq:incidence-covector}, and let $\mathcal R_1$ be the orthogonal
configuration reflection at $d_1$.  Before a row operation the two
configuration/velocity conormal blocks are
\[
 (a_{f_1}^*,ba_{f_1}^*),\qquad
 (a_{f_2}^*,ba_{f_2}^*+(B-b)\mathcal R_1a_{f_2}^*).
\]
Since $\mathcal R_1-I$ has range in
$\operatorname{Ran}a_{f_1}^*$, there is a linear map $C:\R^d\to\R^d$
with $\|C\|\le2$ such that
\[
 \mathcal R_1a_{f_2}^*-a_{f_2}^*=a_{f_1}^*C.
\]
Add the corresponding $C^*$-multiple of the first $d$ landing rows to the
second $d$.  This gives
\[
 Z_1=a_{f_1}^*,\qquad Z_2=\mathcal R_1a_{f_2}^*,
 \qquad \widetilde W=(bZ_1,BZ_2).
\]

The same operation diagonalizes the time columns.  Indeed, differentiating
the second residual with respect to $t_{d_1}$ gives
$a_{f_2}(I-\mathcal R_1)V=-C^*a_{f_1}V=-C^*g_1$,
which is cancelled by the added first row; the $t_{d_2}$ column is unchanged.
This proves \eqref{eq:two-speed-time-diagonal}.  It also shows directly that
the $2(d-1)$-dimensional row quotient is obtained by restricting the first block
to $g_1^\perp$ and the second to $g_2^\perp$.

Distinct particle edges have normalized scalar incidence product of modulus
at most $1/2$.  Orthogonal reflection and restriction to the subspaces
$g_i^\perp$ do not increase its norm, proving
\eqref{eq:two-speed-cross-contraction}.  The Gram matrix of $\overline Z$ is
therefore
\[
 \begin{pmatrix}I_{d-1}&K_1^*Z_1^*Z_2K_2\\
 K_2^*Z_2^*Z_1K_1&I_{d-1}\end{pmatrix},
\]
whose spectrum lies in $[1/2,3/2]$.  Equation
\eqref{eq:two-speed-graph-form} follows, and multiplication of its $2(d-1)$
singular values gives \eqref{eq:two-speed-initial-volume}.
The row operation and its inverse are uniformly bounded, so the equivalent
orthogonal quotient in Lemma~\ref{lem:time-wedge-factorization} differs only
by an absolute constant.  No upper collision has been used in this
calculation, which proves the final scope assertion.
\end{proof}

\begin{lemma}[Fixed-normal transport of the two-speed quotient]
\label{lem:fixed-normal-two-speed-transport}
Let $(\overline Z,\overline W)$ be the $2(d-1)$-frame in
Lemma~\ref{lem:two-speed-landing-quotient}, so that
\[
 \overline W=\overline ZL,
 \qquad L=\operatorname{diag}(bI_{d-1},BI_{d-1})>0.
\]
Propagate this frame backward through an arbitrary prescribed upper
collision word while retaining every collision time and normal, so no
incoming-root normal is differentiated.  If $r$ is the sum of the intervening
free-flight lengths, then, after one orthogonal identification of the initial
and final reduced configuration spaces,
\begin{equation}\label{eq:fixed-normal-two-speed-transport}
 \overline Z_{\rm fn}=O\overline Z,
 \qquad
 \overline W_{\rm fn}=O\overline Z(L+rI_{2(d-1)})
\end{equation}
for an orthogonal map $O$.  In particular
\begin{equation}\label{eq:fixed-normal-two-speed-volume}
 \Jac_{2(d-1)}(\overline W_{\rm fn})
 =\Jac_{2(d-1)}(\overline Z)\det(L+rI_{2(d-1)})
 \ge\Jac_{2(d-1)}(\overline W)
 \ge2^{-(d-1)}(bB)^{d-1}.
\end{equation}
Thus neither a fixed collision reflection nor an arbitrary fixed-normal
co-tree word can focus the mixed landing quotient.  Any possible loss in
the rooted cell must be created by the curvature terms produced when one or
more incoming contact normals are differentiated.
\end{lemma}

\begin{proof}
For a retained normal the collision part of the adjoint-pullback conormal map
is $(Z,W)\mapsto(\mathcal RZ,\mathcal RW)$ with $\mathcal R$ orthogonal.
It preserves every identity of the form $W=ZL_0$.  A backward free flight of
length $\rho$ sends
\[
 (Z,ZL_0)\longmapsto(Z,Z(L_0+\rho I)).
\]
Induction through the word proves
\eqref{eq:fixed-normal-two-speed-transport}, because the scalar free-flight
increments commute with the right-hand matrix $L$.  Exterior volume is
unchanged by $O$, and
\[
 \det(L+rI_{2(d-1)})
 =(b+r)^{d-1}(B+r)^{d-1}\ge(bB)^{d-1}.
\]
Combine this with
\eqref{eq:two-speed-graph-form}--\eqref{eq:two-speed-initial-volume}.
\end{proof}

\begin{lemma}[Symplectic compliance factorization of the rooted lift]
\label{lem:rooted-compliance-factorization}
In orthonormal reduced-particle coordinates, write the complete backward
adjoint-pullback upper lift from the separator through the rooted upper word as a
block map
\[
 \binom{Z}{W}\longmapsto
 \begin{pmatrix}A_U&B_U\\C_U&D_U\end{pmatrix}
 \binom{Z}{W}.
\]
After the collision-by-collision orthogonal identifications, $D_U$ is
invertible and
\begin{equation}\label{eq:rooted-compliance-factorization}
 K_U:=D_U^{-1}C_U=K_U^*\ge0,
 \qquad
 W_{\rm rt}=D_U(W+K_UZ).
\end{equation}
Moreover $D_U$ does not decrease any exterior volume:
\begin{equation}\label{eq:rooted-D-exterior-expansion}
 \Jac_k(D_UF)\ge\Jac_k(F)
 \qquad(1\le k\le d(P-1)).
\end{equation}
The operator $K_U$ is the one-port grounded compliance of the time-expanded
orthogonal Jacobi network as seen from the separator port: a free segment of
length $\ell$ has resistance $\ell I$, and eliminating an incoming contact
adds its nonnegative cylinder curvature as a shunt.  If instead a set
$\partial$ of event ports is retained before the interior fields are
eliminated, denote the resulting multiport compliance by
$\mathcal K_U^\partial$.  If a covariant path $\gamma$ joins two ports
$u,v\in\partial$ and
\[
 j_u=-f,\qquad j_v=O_\gamma f,
\]
then
\begin{equation}\label{eq:rooted-compliance-path-bound}
 \langle j,\mathcal K_U^\partial j\rangle
 \le \Big(\sum_{e\in\gamma}\ell_e\Big)|f|^2.
\end{equation}
Here $j$ is extended by zero at all other ports and projected to the grounded
range.
\end{lemma}

\begin{proof}
The two local conormal factors, after absorbing their orthogonal reflections,
are the symplectic shears
\[
 \begin{pmatrix}I&0\\rI&I\end{pmatrix},\qquad
 \begin{pmatrix}I&\mathcal H\\0&I\end{pmatrix},
 \qquad r\ge0,\quad\mathcal H\ge0.
\]
Their product belongs to the positive symplectic semigroup.  Equivalently,
apply the complete lift first to a frame $(0,F)$.  Its initial symmetric
conormal Gram is zero, so
Lemma~\ref{lem:conormal-frame-identities} gives
\[
 \Jac_k(D_UF)\ge\Jac_k(F)
\]
for every $k$ and $F$.  In particular $s_{\min}(D_U)\ge1$; hence $D_U$ is
invertible and \eqref{eq:rooted-D-exterior-expansion} already holds.  We may
therefore define $K_U=D_U^{-1}C_U$ without a circular positivity argument.

Equivalently for the remaining identification,
split the word at all event times and assign to a port field $\xi$ the index
form
\begin{equation}\label{eq:rooted-discrete-index-form}
\mathcal I(\xi)=
 \sum_e\ell_e^{-1}|A_e\xi|^2
 +\sum_v\langle\mathcal H_v\xi_v,\xi_v\rangle.
\end{equation}
We now identify this form with the physical transfer blocks by an explicit
prefix induction.  Orient every segment, choose a retained boundary set
$\partial$, and put
$\Lambda=\operatorname{diag}_e(\ell_eI_d)$ and
$\mathbb A=(A_I\ A_\partial)$.  The complete grounded KKT system is
\begin{equation}\label{eq:rooted-compliance-KKT-matrix}
 \begin{pmatrix}
  -\Lambda&A_I&A_\partial\\
  A_I^*&H_{II}&H_{I\partial}\\
  A_\partial^*&H_{\partial I}&H_{\partial\partial}
 \end{pmatrix}
 \begin{pmatrix}J\\\xi_I\\\xi_\partial\end{pmatrix}
 =
 \begin{pmatrix}0\\0\\i_\partial\end{pmatrix},
\end{equation}
where $H=\operatorname{diag}_v\mathcal H_v$ after the orthogonal pair
identifications.  A free segment contributes the exact boundary relation
\[
 J_v=J_u,\qquad \xi_v=\xi_u+\ell J_u,
\]
whose transfer matrix is
$\left(\begin{smallmatrix}I&0\\\ell I&I\end{smallmatrix}\right)$.
A collision shunt contributes
\[
 \xi^-=\xi^+,\qquad J^-=J^++\mathcal H\xi^+,
\]
whose transfer matrix is
$\left(\begin{smallmatrix}I&\mathcal H\\0&I\end{smallmatrix}\right)$
in the same $(\text{current},\text{potential})$ ordering
$(J,\xi)=(Z,W)$ used by the adjoint frame.  These are precisely the two
physical shears displayed above.  First take
$\partial=\{\mathrm{sep}\}$.  If the transfer relation of the first $m$
events equals the one-port relation obtained by eliminating the first $m$
interior KKT fields, adjoining event $m+1$ composes both relations with the
same local two-by-two block.  The assertion holds for the empty prefix, so
induction proves it for the complete word, including every causal
off-diagonal block.  For a larger retained set $\partial$, cut the same
time-expanded network at those ports and use the same edge and shunt laws;
this defines the multiport KKT response, but it is not identified with the
single transfer quotient $D_U^{-1}C_U$.

The endpoint sign in that one-port quotient is worth recording explicitly.
Orient the upper trajectory edges in the backward direction, from the
separator towards the grounded upper roots.  Thus the transfer block sends
the separator data to the root data,
\begin{equation}\label{eq:rooted-transfer-endpoint-orientation}
 \binom{J_{\rm rt}}{\xi_{\rm rt}}
 =\begin{pmatrix}A_U&B_U\\C_U&D_U\end{pmatrix}
   \binom{J_{\rm sep}}{\xi_{\rm sep}},
 \qquad i_{\rm sep}=-J_{\rm sep}.
\end{equation}
The last identity is not an extra convention: $i_{\rm sep}$ is the KKT
current leaving the retained boundary into the exterior, whereas
$J_{\rm sep}$ is oriented from that boundary into the upper network.
Grounding gives $\xi_{\rm rt}=0$, and the second block row of
\eqref{eq:rooted-transfer-endpoint-orientation} therefore gives
\begin{equation}\label{eq:rooted-transfer-grounded-response}
 0=C_UJ_{\rm sep}+D_U\xi_{\rm sep},
 \qquad
 \xi_{\rm sep}=D_U^{-1}C_U\,i_{\rm sep}.
\end{equation}
This fixes both the current/potential ordering and the endpoint sign before
any appeal to symmetry or positivity.

Eliminating $J$ in
\eqref{eq:rooted-compliance-KKT-matrix} gives
\[
 L_U=\mathbb A^*\Lambda^{-1}\mathbb A+H.
\]
Grounding the common translation and the marked-tree cluster translations
makes the interior block strictly positive; eliminating $\xi_I$ gives the
boundary stiffness $S_U^\partial$, and its inverse is the multiport grounded
compliance
\[
 \mathcal K_U^\partial=(S_U^\partial)^{-1}.
\]
For the single retained separator port, the prefix identity says that this
same boundary relation, in transfer-block notation, is
\begin{equation}\label{eq:rooted-one-port-transfer-KKT-identity}
 K_U=\mathcal K_U^{\{\mathrm{sep}\}}=D_U^{-1}C_U.
\end{equation}
This is the promised matrix identification.  It does not identify
$D_U^{-1}C_U$ with the larger multiport matrix and does not follow merely
from the sign of the index form.

The free Euler--Lagrange current is
$\ell_e^{-1}A_e\xi$ in the $d$-dimensional factor space; at a rooted contact its jump is
$\mathcal H_v\xi_v$.  These are exactly the tangent equations dual to
\eqref{eq:backward-free-conormal}--
\eqref{eq:backward-collision-conormal}.  Grounding the common translation and
the marked-tree cluster translations therefore identifies the terminal
position response with the compliance of
\eqref{eq:rooted-discrete-index-form}.  In transfer-block notation that
response is $D_U^{-1}C_U$.  Symmetry and nonnegativity of the index form prove
$K_U=K_U^*\ge0$ and the first identity in
\eqref{eq:rooted-compliance-factorization}.  The second is the block identity
$C_UZ+D_UW=D_U(K_UZ+W)$.

The opening frame argument already proves
\eqref{eq:rooted-D-exterior-expansion} and the invertibility used in this
calculation.

Finally, retain the event-port set $\partial$ containing $u,v$ and first
delete the nonnegative shunts in
\eqref{eq:rooted-discrete-index-form}.  The Thomson principle expresses the
multiport compliance quadratic form
$\langle j,\mathcal K_U^\partial j\rangle$ as the minimum of
$\sum_e\ell_e|J_e|^2$ over currents with divergence $j$.  Send the transported
current $f$ only along $\gamma$; orthogonality keeps every edge norm equal to
$|f|$ and gives the right-hand side of
\eqref{eq:rooted-compliance-path-bound}.  Restoring the shunts can only
decrease the grounded inverse in Loewner order.  This proves the path bound.
\end{proof}

\begin{lemma}[Zero-set covariance and the relative new-line shear]
\label{lem:radial-chord-zero-set-shear}
Assume that the first two lower edges overlap and relabel them as
\[
 f_{d_1}=\{a,b\},\qquad f_{d_2}=\{b,c\},
\]
where $c$ is the line which first appears at $d_2$.  Let $u_c$ be its last
upper contact.  In the normalized incidence coordinates of
Lemma~\ref{lem:two-speed-landing-quotient}, let $C$ be the bounded map defined
by
\begin{equation}\label{eq:radial-row-C-definition}
 \mathcal R_1a_{f_2}^*-a_{f_2}^*=a_{f_1}^*C,
\end{equation}
and put
\begin{equation}\label{eq:radial-row-operation}
 \begin{aligned}
 \mathscr C_i&:=R_i-\eps\omega_i\quad(i=1,2),\\
 \widetilde R_2^{\rm zs}
 &:=\eps\omega_2+\mathscr C_2+C^*\mathscr C_1
   =R_2+C^*(R_1-\eps\omega_1),\\
 L_D&:=\begin{pmatrix}I_d&0\\ C^*&I_d\end{pmatrix}.
 \end{aligned}
\end{equation}
Thus the system
\(R_1=\eps\omega_1\),
\(\widetilde R_2^{\rm zs}=\eps\omega_2\) is equivalent to the two
original contact equations \(R_i=\eps\omega_i\), $i=1,2$.  On their common zero set, with the landing
normals retained, and as total derivatives after every earlier incoming-root
substitution,
\begin{equation}\label{eq:radial-row-total-derivatives}
 D_{t_{d_1}}\widetilde R_2^{\rm zs}=0,
 \qquad D_{t_{d_2}}\widetilde R_2^{\rm zs}=g_2,
 \qquad D_{(X_c,V_c)}\widetilde R_2^{\rm zs}
       =D_{(X_c,V_c)}R_2 .
\end{equation}

Let $w_2$ be the post-$d_1$ velocity of the companion line $b$ at the
second landing.  It is independent of $(X_c,V_c)$, $t_{d_1}$ and $t_{d_2}$
in this chart.  With
\[
 \tau_u=t_{u_c}-s,\qquad \tau_d=t_{d_2}-s,
 \qquad \ell=\tau_d-\tau_u=t_{d_2}-t_{u_c}>0,
\]
the change
\begin{equation}\label{eq:radial-relative-shear}
 g_2=V_c-w_2,
 \qquad \widehat X_c=X_c+\tau_uw_2
\end{equation}
is a fibrewise block-unitriangular change whose pivot-variable derivative
and inverse have identity diagonal blocks.  Its retained-variable
coefficients obey the existing $\cL_\eps^C$ cell cutoff and introduce no
Jacobian loss.  In these variables the endpoint rows are row-equivalent to
\begin{equation}\label{eq:radial-affine-endpoint-rows}
 U_{u_c}=\widehat X_c+\tau_ug_2-\Xi_u,
 \qquad
 \widetilde R_2^{\rm zs}=\widehat X_c+\tau_dg_2-\Xi_d,
\end{equation}
where $\Xi_u,\Xi_d$ are independent of $t_{d_2}$ and, on the common zero
fibre, of $t_{d_1}$.  In particular
\begin{equation}\label{eq:radial-first-time-zero-column}
 D_{t_{d_1}}(U_{u_c},\widetilde R_2^{\rm zs},|g_2|)=0.
\end{equation}
Both the row operation and the shear preserve every preceding Schur
complement.
\end{lemma}

\begin{proof}
Equation~\eqref{eq:radial-row-C-definition} is the exact row relation in the
proof of Lemma~\ref{lem:two-speed-landing-quotient}.  There
\[
 D_{t_{d_1}}R_2=-C^*g_1,
 \qquad D_{t_{d_1}}R_1=g_1,
\]
while $D_{t_{d_2}}R_1=0$ and $D_{t_{d_2}}R_2=g_2$.  The first landing does
not contain $c$, hence $D_{(X_c,V_c)}R_1=0$.  At fixed landing normals,
$D_p\mathscr C_i=D_pR_i$ for every pivot variable used above.  If $L_D$
depends on retained variables, the product-rule term vanishes on the actual
contact fibre, because it multiplies $(\mathscr C_1,\mathscr C_2)$ rather
than $(R_1,R_2)$:
\[
 D\bigl(L_D(\mathscr C_1,\mathscr C_2)^T\bigr)
 =(DL_D)(\mathscr C_1,\mathscr C_2)^T
   +L_DD(\mathscr C_1,\mathscr C_2)^T
 =L_DD(\mathscr C_1,\mathscr C_2)^T.
\]
This proves \eqref{eq:radial-row-total-derivatives}, including the total
root-substituted meaning of the derivatives.

The first two lower contacts are consecutive.  The companion line $b$ has no
lower collision after $d_1$ and before $d_2$, so its velocity there is the
fixed-normal elastic image at $d_1$ of separator velocities on lines $a,b$.
It therefore contains neither $(X_c,V_c)$ nor $t_{d_2}$.  Because $d_1$ is
the first lower contact and its normal is retained, free transport to $d_1$
changes positions but not the incoming or outgoing velocities, so $w_2$ is
also independent of $t_{d_1}$.  Earlier upper roots do not change this
claim: the forward lower reconstruction starts from the raw independent
separator states.  This is the velocity form of the causal independence in
Lemma~\ref{lem:causal-free-line-endpoints}.

Substitution of \eqref{eq:radial-relative-shear} in the two endpoint
residuals gives \eqref{eq:radial-affine-endpoint-rows}, with the fixed normal
terms absorbed into $\Xi_d$; the cancellation of
the $t_{d_1}$ column in \eqref{eq:radial-row-total-derivatives} gives the
last independence assertion.  The Jacobian of the change in
$(X_c,V_c)$ has identity diagonal blocks, so at fixed retained variables its
linear part and inverse have norm one.  Dependence of the translation on
the retained marked-tree time $t_{u_c}$ is already covered by the existing
$\cL_\eps^C$ time cutoff and changes no determinant or epsilon power.
Finally, for an
earlier pivot block,
\[
 \begin{pmatrix}A&B\\C_0&D\end{pmatrix}
 \longmapsto
 \begin{pmatrix}A&B-AR\\C_0&D-C_0R\end{pmatrix},
 \qquad
 (D-C_0R)-C_0A^{-1}(B-AR)=D-C_0A^{-1}B,
\]
which proves the Schur-covariance assertion.
\end{proof}

\begin{lemma}[Exact radial-conjugate chord block]
\label{lem:radial-conjugate-chord-block}
On a regular cell with $g_2\ne0$, retain $\rho_2=|g_2|$ and put
\[
 \widehat g_2=g_2/|g_2|,
 \qquad P_{g_2}=I_d-\widehat g_2\tensor\widehat g_2.
\]
The exact derivative of the endpoint map in
\eqref{eq:radial-affine-endpoint-rows} is
\begin{equation}\label{eq:radial-conjugate-matrix}
 \mathcal C_{\rm rad,d}
 :=D_{(\widehat X_c,g_2,t_{d_2})}
      (U_{u_c},\widetilde R_2^{\rm zs},\rho_2)
 =\begin{pmatrix}
    I_d&\tau_uI_d&0\\
    I_d&\tau_dI_d&g_2\\
    0&\widehat g_2^T&0
  \end{pmatrix},
\end{equation}
and
\begin{equation}\label{eq:radial-conjugate-determinants}
 |\det\mathcal C_{\rm rad,d}|=\ell^{d-1}|g_2|.
\end{equation}
With the energy row $g_2^T$ in place of $\widehat g_2^T$, the absolute
determinant is $\ell^{d-1}|g_2|^2$.

For every endpoint-drop variation $\eta\in\R^d$, fixing $\rho_2$ gives
\begin{equation}\label{eq:radial-tangent-solve}
 \delta g_2=\ell^{-1}P_{g_2}\eta,
 \qquad
 \delta t_{d_2}=|g_2|^{-2}g_2\cdot\eta,
 \qquad
 \eta\cdot\delta g_2=\ell^{-1}|P_{g_2}\eta|^2\ge0.
\end{equation}
If $T_\gamma$ is the exact total endpoint-trace map after
Lemma~\ref{lem:radial-chord-zero-set-shear}, $D=(-I_d,I_d)$ and
$A_\gamma=DT_\gamma$, then
\begin{equation}\label{eq:radial-columnwise-inverse}
 \Pi_g\mathcal C_{\rm rad,d}^{-1}
      \binom{T_\gamma}{0}
 =\ell^{-1}P_{g_2}A_\gamma.
\end{equation}
In the energy normalization the literal physical symmetric bordered block,
up to a simultaneous harmless sign convention, is
\begin{equation}\label{eq:radial-physical-KKT-block}
 \begin{pmatrix}
  -\ell I_d&g_2&A_\gamma\\
  g_2^T&0&0\\
  A_\gamma^*&0&H
 \end{pmatrix},
\end{equation}
and eliminating the radial chord current and multiplier gives exactly
\begin{equation}\label{eq:radial-positive-Schur-term}
 H_{\rm new}=H+\ell^{-1}A_\gamma^*P_{g_2}A_\gamma\succeq H.
\end{equation}
\end{lemma}

\begin{proof}
Subtract the first $d$ rows of \eqref{eq:radial-conjugate-matrix} from the
second $d$.  The only nontrivial determinant is
\[
 \det\begin{pmatrix}\ell I_d&g_2\\
               \widehat g_2^T&0\end{pmatrix}
 =-\ell^{d-1}|g_2|,
\]
which proves \eqref{eq:radial-conjugate-determinants}; the energy row
multiplies it by $|g_2|$.  The same row subtraction reduces the inverse
problem to
\[
 \ell\,\delta g_2+g_2\,\delta t_{d_2}=\eta,
 \qquad g_2\cdot\delta g_2=0.
\]
Orthogonal projection and contraction with $g_2$ give
\eqref{eq:radial-tangent-solve}.  Applying that solution columnwise with
$\eta=DT_\gamma$ proves \eqref{eq:radial-columnwise-inverse}.

The upper-left $(d+1)$-square block of
\eqref{eq:radial-physical-KKT-block} is the constrained energy form for this
same inverse problem, not a comparison matrix.  Its inverse has upper-left
block $-\ell^{-1}P_{g_2}$.  The block Schur formula therefore gives
\[
 H-[A_\gamma^*\ 0]
 \begin{pmatrix}-\ell I_d&g_2\\g_2^T&0\end{pmatrix}^{-1}
 \binom{A_\gamma}{0}
 =H+\ell^{-1}A_\gamma^*P_{g_2}A_\gamma,
\]
which proves \eqref{eq:radial-positive-Schur-term}.  The energy row and the
normalized speed row define the same tangent space and differ by the positive
factor $|g_2|$; hence no inverse speed enters a raw or covariant boundary
gauge.  The locus $g_2=0$ belongs to the existing zero-flux/rank-loss null
stratum and no unit vector is defined there.
\end{proof}

\begin{lemma}[Radial chord compatibility with the positive network]
\label{lem:radial-chord-network-compatibility}
In the overlapping packet word, replacing the full chord factor
$\ell^{-1}|A_\gamma\xi|^2$ by
\begin{equation}\label{eq:radial-network-factor}
 \ell^{-1}|E_{g_2}^*A_\gamma\xi|^2,
 \qquad E_{g_2}:\R^{d-1}\longrightarrow g_2^\perp
\end{equation}
on a finite orthogonal half-open atlas preserves the trace left inverse, the
grounded-kernel assertion and the physical Schur--Jacobi prefix identity.
The local map $E_{g_2}$ is not a global recovery variable.  The first
landing-time column remains $g_1$.
\end{lemma}

\begin{proof}
In the trace-completeness certificate of
Lemma~\ref{lem:appC-overlap-trace-certificate}, every frontier crossed by the
chord interval is reconstructed from the selectors of the physical lines
other than $c$.  The explicit triangular left inverse in
Lemma~\ref{lem:appC-uniform-trace-left-inverse} does not select or invert the
chord row.  It therefore remains valid if that row is deleted, and a fortiori
if the nonnegative rank-$(d-1)$ factor \eqref{eq:radial-network-factor} is
inserted.  Thus no new grounded kernel is possible.

For the prefix identity, all free-segment, complete-pair, incoming-root and
forest blocks are unchanged.  At the chord prefix,
Lemma~\ref{lem:radial-conjugate-chord-block} identifies the physical
off-diagonal inverse and the KKT elimination with the same bordered matrix
\eqref{eq:radial-physical-KKT-block}.  Radial-chord-first and forest-first
elimination are therefore Schur complements of one symmetric system, so
Schur associativity gives the same one-port response.  Equation
\eqref{eq:radial-first-time-zero-column} shows that none of these rows uses
$t_{d_1}$, proving the last assertion.  The finite atlas is used only to
factor the projector in \eqref{eq:radial-positive-Schur-term}; the physical
derivative and all global coarea coordinates remain Cartesian.
\end{proof}

\begin{lemma}[One-time wedge identity]
\label{lem:one-time-wedge}
Let $S:\R^N\to\R^d$, let $g\ne0$, and let
$Q_g:\R^{d-1}\to g^\perp$ be orthogonal.  Then
\begin{equation}\label{eq:one-time-wedge}
 \Jac_{d-1\mid1}(S,g)
 :=\left(\sum_{|I|=d-1}|\det[S_I\ g]|^2\right)^{1/2}
 =|g|\Jac_{d-1}(Q_g^*S).
\end{equation}
\end{lemma}

\begin{proof}
Rotate $g/|g|$ to the last coordinate in the row space.  Expansion along
the last column and Cauchy--Binet on the first $d-1$ rows give
\eqref{eq:one-time-wedge}.
\end{proof}

\begin{proposition}[Flux completion for disjoint first landings]
\label{prop:disjoint-radial-fct}
On every regular packet cell whose first two lower edges are disjoint,
\eqref{eq:tree-exterior-coercivity} holds.
\end{proposition}

\begin{proof}
For disjoint edges the cross incidence block in
Lemma~\ref{lem:two-speed-landing-quotient} is zero.  After quotienting the
two landing-time directions, its exact tangential frame therefore satisfies
\[
 \overline W=\overline Z
 \operatorname{diag}(\Delta_1I_{d-1},\Delta_2I_{d-1}),
 \qquad
 \tfrac12(\overline Z^*\overline W+
                 \overline W^*\overline Z)\ge0,
\]
and
\[
 \Jac_{2(d-1)}(\overline W)
 \ge c_d|\Delta_1\Delta_2|^{d-1}.
\]
Apply the total-derivative rooted conormal theorem
Lemma~\ref{lem:rooted-tree-schur-lift}, not merely the fixed-normal
transport lemma.  It retains every incoming-root curvature shunt and cannot
decrease this exterior volume.  Finally
Lemma~\ref{lem:time-wedge-factorization} supplies the two landing fluxes and
gives \eqref{eq:tree-exterior-coercivity}.
\end{proof}

\begin{proposition}[Radial flux completion for overlapping first landings]
\label{prop:overlap-radial-fct}
On every regular packet cell whose first two lower edges overlap,
\eqref{eq:tree-exterior-coercivity} holds.
\end{proposition}

\begin{proof}
Use the notation of Lemma~\ref{lem:radial-chord-zero-set-shear}.  The
reverse-Kruskal tree contains $u_c$; remove it and write
$F=T\setminus\{u_c\}$.  Let $A$ be the contacts outside
$T\cup\{d_1,d_2\}$.  Before any quotient, take a square augmented derivative
with row order equivalent to
\begin{equation}\label{eq:radial-full-row-order}
 (C_A,U_F,U_{u_c},\widetilde R_2^{\rm zs},\rho_2,R_1)
\end{equation}
and column order equivalent to
\begin{equation}\label{eq:radial-full-column-order}
 (y_A,h_F,\widehat X_c,g_2,t_{d_2},V_I,t_{d_1}),
 \qquad |I|=d-1.
\end{equation}
The root blocks have $d|A|$ rows and columns; $F$ has $P-2$ edges and
contributes $d(P-2)$; the radial chord contributes $2d+1$; and the first
landing contributes $d=(d-1)+1$.  Thus
\eqref{eq:radial-full-row-order} and
\eqref{eq:radial-full-column-order} have the same dimension.

First eliminate every $y_a$, using the exact chronological determinant
\eqref{eq:unrooted-causal-determinant}.  To track the $t_{d_1}$ column, split
$A=A_U\mathbin{\dot\cup}A_+$ into upper contacts and lower contacts strictly
after $d_2$; there is no other lower contact before $d_2$ because $d_1,d_2$
are the first two lower contacts.  In chronological root coordinates,
causality and \eqref{eq:unrooted-causal-determinant} give invertible diagonal
blocks $B_U,B_+$ and
\[
 D_{(y_U,y_+)}(C_U,C_+)
 =\begin{pmatrix}B_U&0\\ B_{+U}&B_+\end{pmatrix},
 \qquad
 D_{t_{d_1}}(C_U,C_+)=\binom{0}{\beta_+}.
\]
For the retained rows
\[
 R_{\rm ret}:=(U_F,U_{u_c},\widetilde R_2^{\rm zs},\rho_2,R_1),
\]
future lower-root variables cannot affect earlier rows, so
$D_{(y_U,y_+)}R_{\rm ret}=(R_U,0)$.  Hence their root-elimination Schur
correction to the $t_{d_1}$ column is exactly
\[
 (R_U,0)
 \begin{pmatrix}B_U^{-1}&0\\ *&B_+^{-1}\end{pmatrix}
 \binom{0}{\beta_+}=0.
\]
Thus future lower roots may depend on $t_{d_1}$, but they do not alter the
retained $t_{d_1}$ column.  Keep the forest variables live and
eliminate the radial block before the grounded forest quotient.  Its
determinant is \eqref{eq:radial-conjugate-determinants}, and its literal
Schur update is the nonnegative factor
\eqref{eq:radial-positive-Schur-term}.  By
Lemma~\ref{lem:radial-chord-network-compatibility}, the remaining physical
forest is trace-complete without using the full chord row.  Eliminate its
internal fields and apply Proposition~\ref{prop:appC-positive-network-closure}
to the first-landing current frame restricted to $g_1^\perp$.
Equation~\eqref{eq:radial-first-time-zero-column} makes the radial retained
rows independent of $t_{d_1}$, the upper rows are independent by causality,
and $D_{t_{d_1}}R_1=g_1$.  Together with the vanishing Schur correction just
proved, the remaining time column is therefore exactly $g_1$.
Lemma~\ref{lem:one-time-wedge} then gives
\begin{equation}\label{eq:radial-first-one-time-bound}
 \Jac_{d-1\mid1}(S_1,g_1)
 \ge c(d,P_0)\phi_1\Delta_1^{d-1}.
\end{equation}
These are successive exact block eliminations of the same square matrix, not
a product of unrelated lower bounds.  Therefore its nonneutral part obeys
\begin{equation}\label{eq:radial-augmented-full-determinant}
 |\det J_{\rm aug}|
 \ge c(d,P_0)
 \left(\prod_{a\in A}\eps^{d-1}|g_a\cdot\omega_a|\right)
 \phi_1|g_2|\Delta_1^{d-1}\ell^{d-1}.
\end{equation}

Rotate only the local $g_2$-velocity block so that its first direction is
$\widehat g_2$ and expand along $D\rho_2$.  The remaining determinant uses
exactly $d-1$ directions in $g_2^\perp$, the $d-1$ still-unpivoted
separator-velocity directions selected in
\eqref{eq:radial-first-one-time-bound}, and $t_{d_2},t_{d_1}$.  The first
set belongs to the radial block and is removed before the second
Cauchy--Binet selection, so the two velocity-column sets are disjoint.
Undoing the bounded shear and rotation turns this determinant into a
component of the same total root-substituted matrix
$[\mathcal S_I\ \mathcal T]$ in
Definition~\ref{def:collision-tree-coercivity}; the global recovery system remains in the original Cartesian variables.  Dividing the exact neutral root and forest
pivot factors in \eqref{eq:radial-augmented-full-determinant}, and using
\[
 |g_2|\ge\phi_2,
 \qquad \ell=t_{d_2}-t_{u_c}\ge t_{d_2}-s=\Delta_2,
\]
gives
\[
 \Jac_{2(d-1)\mid2}(\mathcal S,\mathcal T)
 \ge c(d,P_0)\phi_1\phi_2
       |\Delta_1\Delta_2|^{d-1}.
\]
This is \eqref{eq:tree-exterior-coercivity}.  Exact grazing, singular normal
charts and $g_2=0$ lie in the already declared zero-flux/rank-loss null
strata.
\end{proof}

\begin{theorem}[Flux completion for every first-two-landing packet]
\label{thm:packet-fct}
Every regular packet cell whose first two lower landing edges are distinct
and which satisfies assumptions \textup{(a)--(e)} of
Theorem~\ref{thm:two-landing} satisfies
$\mathrm{FCT}(P_0)$.  Equivalently,
\begin{equation}\label{eq:packet-fct-unconditional}
 \Jac_{2(d-1)\mid2}(\mathcal S,\mathcal T)
 \ge c(d,P_0)\phi_1\phi_2
       |\Delta_1\Delta_2|^{d-1}.
\end{equation}
\end{theorem}

\begin{proof}
The first two lower unordered edges are distinct.  They are either disjoint,
when Proposition~\ref{prop:disjoint-radial-fct} applies, or share exactly one
line, when Proposition~\ref{prop:overlap-radial-fct} applies.  These cases
exhaust the packet class.
\end{proof}

\begin{proposition}[Conditional covariant-port bridge criterion]
\label{prop:covariant-rooted-bridge}
Choose the first two lower landings and their matched crossing lines by
Lemma~\ref{lem:upper-tree-landings}.  Allow arbitrary regular unmarked upper
contacts and arbitrary lower contacts strictly after the second landing.
Assume that after all other root and tree-position eliminations the two
selected bridge rows have the exact port-diagonal form
\eqref{eq:covariant-rooted-bridge-block}, with no additional diagonal
self-response, and that $P,Q$ are the compliance compressions described
below.  Then the total root-substituted separator-velocity Schur row
$\mathcal S$ in \eqref{eq:tree-schur-row} satisfies
\begin{equation}\label{eq:covariant-rooted-bridge-volume}
 \Jac_{2d}(\mathcal S)\ge c(d,P_0)h^{2d},
 \qquad h=\epstar/4.
\end{equation}
The criterion itself places no restriction on the first lower gap, the
overlap of the two lower edges, or the retained collision normals.
\end{proposition}

\fallbackproofinresource{supp:proof:prop:covariant-rooted-bridge}

\begin{corollary}[Packet consequence of the covariant-port criterion]
\label{cor:conditional-long-gap-rooted-packet}
Assume that every long-first-gap overlapping cell on at least four particle
lines with upper co-tree contacts satisfies the port-diagonal hypothesis of
Proposition~\ref{prop:covariant-rooted-bridge}.  Then every packet satisfying
assumptions \textup{(a)--(e)} of Theorem~\ref{thm:two-landing} obeys
\begin{equation}\label{eq:conditional-long-gap-rooted-packet}
 J_{H,\kappa}(Q)
 \le \cL_\eps^{C|H|}\eps^{d-1}\epstar^{-2d}
       \int Q_{H,\kappa}^{\sharp,\mathrm{port}}
\end{equation}
for every nonnegative joint kernel $Q$.  The output box is independent of the
concrete fixed boundary values.
\end{corollary}

\fallbackproofinresource{supp:proof:cor:conditional-long-gap-rooted-packet}

\begin{remark}[Why the port-diagonal hypothesis remains open]
\label{rem:covariant-port-self-response}
Lemma~\ref{lem:rooted-compliance-factorization} proves the covariant passive
compliance and its sharp balanced-current path bound.  It does not delete the
compliance self-blocks created when earlier marked-tree position rows are
eliminated.  Algebraically, a Schur complement of a passive multiport may
have nonzero diagonal compliance blocks, and its individual off-diagonal
blocks need not be controlled by the two bare flight spans.  Nor does a
single preferred bridge minor inherit positivity from the full operator:
positive semidefiniteness concerns the complete quadratic form, not each
rectangular subblock or a polar-aligned correction.
Any attempt to prove the conditional covariant-port shortcut must therefore
control the full Pl\"ucker vector of the response or use a partial/multiple
root-normal exchange atlas; one preferred $2$-by-$2$ bridge minor is
insufficient.  This condition is not used by the main overlap proof, which
uses the full positive-network exterior frame of Proposition
\ref{prop:asymmetric-new-line-peeling} and the finite-branch coarea theorem.
\end{remark}

\begin{lemma}[Unrooted neutral completion]
\label{lem:unrooted-neutral-completion}
Let $A$ be the set of contacts outside the upper tree and
$\{d_1,d_2\}$.  Before making their incoming-root substitutions, retain for
every $a\in A$ its time and $d-1$ local normal coordinates, denoted together
by $y_a\in\R^d$, and replace its scalar sphere delta by the vector residual
$C_a=R_a-\eps\omega_a$.  After eliminating only the upper-tree positions,
put $C_A=(C_a)_{a\in A}$ and $y_A=(y_a)_{a\in A}$.  On a fixed strict
chronological cell,
\begin{equation}\label{eq:unrooted-causal-determinant}
 \left|\det D_{y_A}C_A\right|
 =\eps^{(d-1)|A|}\prod_{a\in A}|g_a\cdot\omega_a|,
\end{equation}
up to the bounded density of the chosen normal charts.

For every $2(d-1)$-set $I$ of separator-velocity coordinates, the square
Jacobian of $(C_A,\mathcal R_D)$ with respect to
$(y_A,V_I,t_{d_1},t_{d_2})$ satisfies
\begin{multline}\label{eq:unrooted-schur-equivalence}
 \left|\det D_{(y_A,V_I,t_{d_1},t_{d_2})}
       (C_A,\mathcal R_D)\right|\\
 =\eps^{(d-1)|A|}\prod_{a\in A}|g_a\cdot\omega_a|
   \left|\det[\mathcal S_I\ \mathcal T]\right|.
\end{multline}
Here $(\mathcal S,\mathcal T)$ are the total root-substituted derivatives in
Definition~\ref{def:collision-tree-coercivity}.  Thus
$\mathrm{FCT}(P_0)$ is equivalently a lower bound for a completely explicit
unrooted joint Jacobian; all additional factors in
\eqref{eq:unrooted-schur-equivalence} are neutral in the molecule measure.
\end{lemma}

\fallbackproofinresource{supp:proof:lem:unrooted-neutral-completion}

\begin{lemma}[$(d-1)$-normal root exchange]
\label{lem:three-normal-root-exchange}
Use the notation of Lemma~\ref{lem:unrooted-neutral-completion}, assume
$P\ge4$, and choose one contact $a\in A$.  Write
$y_b=(t_b,\theta_b)$, where $\theta_b\in\R^{d-1}$ is a local normal coordinate.
Before solving the contact at $a$, eliminate the rows
$(C_{A\setminus\{a\}},\mathcal U)$ against the columns
$(y_{A\setminus\{a\}},h)$.  Denote the resulting $3d$-row map by
\[
 \mathcal F_a=(C_a^{\rm red},R_1^{\rm red},R_2^{\rm red})
\]
and keep $\theta_a$ fixed.  Put
\begin{equation}\label{eq:root-exchange-jacobian}
 \Jac_{3d-3\mid d-1}^{(a)}
 :=\left(\sum_{|I|=3d-3}
 \left|\det D_{(t_a,V_I,t_{d_1},t_{d_2})}\mathcal F_a\right|^2
 \right)^{1/2},
\end{equation}
where $I$ ranges over the scalar separator-velocity coordinates.
For every such $I$ the corresponding full square minor satisfies
\begin{multline}\label{eq:root-exchange-block-identity}
 \left|\det D_{(y_{A\setminus\{a\}},h,t_a,V_I,t_{d_1},t_{d_2})}
 (C_{A\setminus\{a\}},\mathcal U,C_a,\mathcal R_D)\right|\\
 =c_T\eps^{(d-1)(|A|-1)}
   \prod_{b\in A\setminus\{a\}}|g_b\cdot\omega_b|
 \left|\det D_{(t_a,V_I,t_{d_1},t_{d_2})}\mathcal F_a\right|,
\end{multline}
where $c_T$ is the reduced tree-incidence determinant and is bounded above
and below in terms of $P_0$.  In particular the right-hand side contains no
factor $\eps^{d-1}$ from the selected normal $\theta_a$.
\end{lemma}

\fallbackproofinresource{supp:proof:lem:three-normal-root-exchange}

\begin{lemma}[Partial normal root-exchange atlas]
\label{lem:partial-normal-root-exchange}
In the setting of Lemma~\ref{lem:three-normal-root-exchange}, choose
$J\subset\{1,\ldots,d-1\}$ with $|J|=k$, $0\le k\le d-1$, in a regular
half-open normal chart with uniformly conditioned metric
$\theta_a=(\theta_a^1,\ldots,\theta_a^{d-1})$.
Keep $\theta_a^J$ as output coordinates and pivot the complementary normal
coordinates $\theta_a^{J^c}$.  After the causal row operations which isolate
those $d-1-k$ tangent components of $C_a$, and after eliminating
$(C_{A\setminus\{a\}},\mathcal U)$ as before, denote the remaining
$(2d+1+k)$-row map by $\mathcal F_{a,J}$.  Define
\begin{equation}\label{eq:partial-root-exchange-jacobian}
 \Jac_{2(d-1)+k\mid d-1}^{(a,J)}
 :=\left(\sum_{|I|=2(d-1)+k}
 \left|\det D_{(t_a,V_I,t_{d_1},t_{d_2})}
               \mathcal F_{a,J}\right|^2\right)^{1/2}.
\end{equation}
For every $I$ in this sum, the corresponding full square minor has modulus
\begin{multline}\label{eq:partial-root-exchange-block-identity}
 c_T\eps^{(d-1)(|A|-1)}
 \prod_{b\in A\setminus\{a\}}|g_b\cdot\omega_b|\,
 \eps^{d-1-k}\times{}\\
 \left|\det D_{(t_a,V_I,t_{d_1},t_{d_2})}
                  \mathcal F_{a,J}\right|,
\end{multline}
up to the uniformly bounded density of the selected normal chart.  For
$k=0$ the last determinant factors further as
\begin{equation}\label{eq:partial-root-standard-endpoint}
 \left|\det D_{(t_a,V_I,t_{d_1},t_{d_2})}
                  \mathcal F_{a,\varnothing}\right|
 =|g_a\cdot\omega_a|\,
   \left|\det[\mathcal S_I\ \mathcal T]\right|.
\end{equation}
For $k=d-1$, \eqref{eq:partial-root-exchange-jacobian} and
\eqref{eq:partial-root-exchange-block-identity} agree, up to the same
uniformly bounded chart-density factor, with
\eqref{eq:root-exchange-jacobian} and
\eqref{eq:root-exchange-block-identity}.
\end{lemma}

\fallbackproofinresource{supp:proof:lem:partial-normal-root-exchange}

\begin{remark}[The gain of a partial exchange]
\label{rem:partial-root-exchange-gain}
The selected collision sphere and flux contribute
$\eps^{d-1}|g_a\cdot\omega_a|$, whereas
\eqref{eq:partial-root-exchange-block-identity} spends only
$\eps^{d-1-k}$.  Thus a $k$-normal chart with $k\ge1$ retains the net factor
$\eps^k$, provided its last Jacobian contains the selected root flux.
For every fixed $C$,
\[
 \eps^k\epstar^{-C}=o(1).
\]
Consequently a partial exchange may pay any fixed additional
$\epstar^{-C}$ loss and still dominate the required packet estimate.  The
sharp $\epstar^{-2d}$ lower bound is needed only in the standard $k=0$
chart.  Iterating the same causal row operation gives product charts for
several roots, with $k$ replaced by the total number of retained normal
directions; only charts using at most the available separator-velocity
dimension are kept.
\end{remark}

\begin{remark}[What the root exchange would buy]
\label{rem:root-exchange-gain}
If, on a bounded family of legal charts, one could choose $a$ so that
\begin{equation}\label{eq:root-exchange-open-bound}
 \Jac_{3d-3\mid d-1}^{(a)}
 \ge c(P_0)|g_a\cdot\omega_a|\phi_1\phi_2
 \ell_a^{d-1}|\Delta_1\Delta_2|^{d-1}
\end{equation}
for an upper scale $\ell_a\ge c\epstar$, then joint coarea would retain the
selected collision-sphere factor $\eps^{d-1}$.  Indeed, all contacts
$b\ne a$ are neutral by
\eqref{eq:root-exchange-block-identity}, whereas the collision numerator at
$a$ supplies $\eps^{d-1}|g_a\cdot\omega_a|$ and its $d-1$ normal coordinates
remain in the output.  The complete ledger would be
\[
 \eps^{-(d-1)P}\eps^{(d-1)(P-1)}\eps^{2(d-1)}\eps^{d-1}=\eps^{2(d-1)}.
\]
Because $\epstar=\exp(-\sqrt{|\log\eps|})$, this surplus absorbs any fixed
additional power of $\epstar^{-1}$ and would be stronger than the packet
bound required below.

Identity \eqref{eq:root-exchange-block-identity} is unconditional, but
\eqref{eq:root-exchange-open-bound} is not proved here.  It can fail for a
blind choice of $a$: when the selected co-tree edge and the two overlapping
lower edges lie on the same three particle lines, the associated incidence
directions are dependent and a repeated-pair exchange can be rank deficient.
Thus the presence of an upper root alone does not prove the stronger
$(d-1)$-normal alternative \eqref{eq:root-exchange-open-bound}.  It does not
exclude the finite atlas in
Lemma~\ref{lem:partial-normal-root-exchange}: on the rank-deficient repeated-
pair chart the standard mixed volume may remain coercive, or a partial
normal exchange may carry the missing Pl\"ucker component.  A valid proof
must establish that joint direct-or-partial-exchange inequality, or else
choose an external root quantitatively and peel the external marked-tree
branches.  This is the root/tree pivot-exchange version of the same bearing
condition recorded in Remark~\ref{rem:ct-obstruction}.
\end{remark}

\begin{proposition}[Tree-only common-normal stratum]
\label{prop:common-normal-fct}
Suppose the collision word consists only of the marked upper tree and the
two marked lower contacts, and every retained collision normal is parallel
to one fixed unit vector $n$.  Then \eqref{eq:tree-exterior-coercivity}
holds.  No collinearity of the velocities is required.
\end{proposition}

\fallbackproofinresource{supp:proof:prop:common-normal-fct}

\begin{remark}[Why a preferred pure bridge block is insufficient]
\label{rem:ct-obstruction}
With $X$ fixed, Lemma~\ref{lem:lower-influence-cone} and Cauchy--Binet do
give the unconditional lower bound
\begin{equation}\label{eq:lower-row-volume}
 \Jac_{2d}(D_V\mathcal R_D)\ge|\Delta_1\Delta_2|^d.
\end{equation}
Also, reinstating the time and $d-1$ normal coordinates of an unmarked
contact gives the exact local determinant
\begin{equation}\label{eq:root-four-jacobian}
 \eps^{d-1}|g\cdot\omega|,
\end{equation}
which is cancelled by its sphere factor and collision flux.  Finally,
eliminating the upper tree cluster-translation rows adds the multiplier
\[
 -(D_h\mathcal U)^{-*}(D_h\mathcal R_D)^*
\]
to the landing conormals, and its remaining velocity part is precisely
$\mathcal S^*$.

These facts alone do not prove $\mathrm{FCT}(P_0)$.  For a prescribed pair
of nonconsecutive lower landings, the one-frame identities do not force the
symmetric part of the combined $Z^*W$ to stay positive after an intervening
lower root.  The antisymmetric part is harmless, but the missing positive
cone prevents a direct use of Lemma~\ref{lem:conormal-frame-identities}.
Thus these one-frame facts do not by themselves prove the
$2(d-1)$-velocity/two-time estimate for an arbitrary prescribed pair of
nonconsecutive lower landings.  Proposition~\ref{prop:common-normal-fct} and
the tree-only calculations above remain useful strata, but positivity of a
complete response does not control a preferred rectangular bridge block.

For the first-two-landing packet class, the missing exterior component is
instead supplied by the radial-conjugate atlas.  Disjoint edges are covered
by Proposition~\ref{prop:disjoint-radial-fct}; on an overlapping cell,
Proposition~\ref{prop:overlap-radial-fct} uses the collision-free new line,
retains its speed, and combines the resulting rank-$(d-1)$ positive chord
factor with the remaining one-time wedge.  Theorem~\ref{thm:packet-fct}
therefore closes $\mathrm{FCT}(P_0)$ for the application and gives the sharp
$\epstar^{-2(d-1)}$ loss.  The earlier full-chord construction in
Proposition~\ref{prop:asymmetric-new-line-peeling} remains an independent
unconditional $\epstar^{-2d}$ fallback.
\end{remark}

\begin{lemma}[Joint packet submersion]
\label{lem:joint-upper-landing-minor}
Fix a regular first-two-landing packet sign, lift, time-order and velocity
cell.  After the complete elastic, transport and unmarked-contact root
substitutions specified in Definition~\ref{def:collision-tree-coercivity}, let
$h\in\R^{d(P-1)}$ be upper-tree cluster translations and let
$V\in\R^{dP}$ be all separator velocities.  For the map, with the two
landing times kept live,
\begin{equation}\label{eq:joint-constraint-map}
 (h,V,\vartheta)\longmapsto(\mathcal U,\mathcal R_D),
\end{equation}
the tree rows are a submersion in $h$.  On their zero fibre the remaining
$2d$-row mixed map has
\begin{equation}\label{eq:joint-minor-lower-bound}
 \Jac_{2(d-1)\mid2}(\mathcal S,\mathcal T)
 \ge c(d,P_0)\phi_1\phi_2|\Delta_1\Delta_2|^{d-1}.
\end{equation}
Consequently the full joint map is a submersion on every nongrazing cell.
It has a finite half-open partition on each piece of which $2(d-1)$ scalar
velocity coordinates together with the two landing times form a pivot minor
of size at least
\begin{equation}\label{eq:maximal-minor-bound}
 \binom{dP}{2(d-1)}^{-1/2}c(d,P_0)
 \phi_1\phi_2|\Delta_1\Delta_2|^{d-1}.
\end{equation}
\end{lemma}

\begin{proof}
Orient the tree from a root particle and use the chronological relative
cluster translations constructed in Lemma~\ref{lem:rooted-tree-schur-lift}.
In raw event coordinates each marked edge has an identity pivot when it
first joins two marked-tree clusters, while intervening root contacts are
invertible physical collision branches.  Consequently
$D_h\mathcal U$ is block triangular with reduced incidence diagonal.
The passage to orthonormal cluster coordinates and its inverse are bounded
because $P\le P_0$, so $D_h\mathcal U$ and its inverse are bounded in terms
of $P_0$.
Integrating the tree deltas first is exact and
leaves the Schur row displayed in \eqref{eq:joint-minor-lower-bound}.
Theorem~\ref{thm:packet-fct} and
Definition~\ref{def:collision-tree-coercivity} give that bound.

For the coordinate statement use Cauchy--Binet:
\[
 \Jac_{2(d-1)\mid2}(\mathcal S,\mathcal T)^2
 =\sum_{I\subset\{1,\ldots,dP\},\ |I|=2(d-1)}
 |\det[\mathcal S_I\ \mathcal T]|^2.
\]
Choose the lexicographically first maximal minor and use half-open
inequalities to break ties.  The pieces are measurable, disjoint and finite.
The implicit-function theorem supplies the $2(d-1)$ velocity and two time pivots.
The two inverse factors $\phi_i^{-1}$ are cancelled by the two marked
collision fluxes.  Exact grazing has zero collision measure, and
\eqref{eq:root-four-jacobian} shows that no inverse grazing factor is hidden
at any unmarked contact.
\end{proof}

\subsection{Global recovery and the submersion statement}

The raw landing pivot and the upper cluster translations are treated as one
map.  Lemma~\ref{lem:joint-upper-landing-minor} supplies its local Schur
complement bound, and
Theorem~\ref{thm:appD-uniform-algebraic-multiplicity} supplies global
finite-fibre recovery.  The next lemma records the resulting complete measure
chart.

\begin{lemma}[Global packet chart and iterated submersion]
\label{lem:global-packet-chart}
On every nongrazing cell $\kappa$, the complete molecule constraint measure
admits the following elimination order:
\begin{enumerate}[label=\textup{(\roman*)}]
\item the nontransversal edge velocities are recovered from the one
transversal velocity on each physical particle line by complete
$2d$-dimensional elastic maps;
\item the nontransversal edge positions are recovered from the transversal
positions by the transport deltas;
\item every contact not in $S_U$ and not equal to $d_1,d_2$ is evaluated by
its incoming contact-time root on the fixed lift.
\item jointly, the $|\mathcal P|-1$ upper tree vector contacts and the two
landing vector contacts recover the relative cluster translations and
 the $2(d-1)$ scalar velocity coordinates selected by
 \eqref{eq:maximal-minor-bound}, together with $t_{d_1},t_{d_2}$; the
 remaining separator velocities and one common translation stay free.
\end{enumerate}
The maps in (i)--(ii) have absolute Jacobian one.  In (iii), the contact flux
is the absolute derivative of the scalar distance and cancels the
one-dimensional coarea denominator; there is at most one incoming root per
lift.  By \eqref{eq:joint-minor-lower-bound}, the joint inverse Jacobian in
(iv) is at most
\[
C(d,P_0)(\phi_1\phi_2)^{-1}|\Delta_1\Delta_2|^{-(d-1)}.
\]
Its sphere factors are
$\eps^{(d-1)(P-1)}\eps^{2(d-1)}$, and the two marked collision weights contain
$\phi_1\phi_2$.  Consequently the fluxes cancel exactly, the complete
constraint map is a submersion on each nongrazing cell, and its only
long-gap power loss is $C(d,P_0)\epstar^{-2(d-1)}$.

The maximal mixed minor supplies only local transversality.  Global recovery
uses Theorem~\ref{thm:appD-uniform-algebraic-multiplicity}: the complete
square pivot system has at most $M(d,P_0)$ regular solutions per retained
fibre, all placed in fixed branch slots.  Hence this optional chart also has
no implicit multiplicity-one assumption.
\end{lemma}

Let
\begin{equation}\label{eq:packet-time-partition}
 \mathcal T_{\rm ret}=S_U,\qquad
 \mathcal T_{\rm piv}=\{d_1,d_2\},\qquad
 \mathcal T_{\rm root}=H\setminus
 (\mathcal T_{\rm ret}\cup\mathcal T_{\rm piv}).
\end{equation}
The following ledger specifies every delta distribution used in the proof
and what remains live in $Q_{H,\kappa}^{\sharp}$.
\begin{center}
\scriptsize
\begin{tabular}{p{0.24\textwidth}p{0.25\textwidth}p{0.21\textwidth}p{0.20\textwidth}}
\toprule
constraint family & pivot variables & inverse Jacobian / factor & variables retained \\
\midrule
elastic velocity deltas
& nontransversal velocity pairs
& $1$
& one crossing velocity per physical line \\
transport deltas
& nontransversal edge positions
& $1$
& crossing positions and live times \\
remaining scalar contacts
& $t_a$, $a\in\mathcal T_{\rm root}$
& coarea denominator cancelled by flux
& selected incoming roots \\
joint tree/landing map
& relative cluster translations, one maximal $2(d-1)$-coordinate velocity minor,
and $t_{d_1},t_{d_2}$
& \(\begin{gathered}
 C(d,P_0)\\[-2pt]
 {}\cdot\eps^{(d-1)|S_U|+2(d-1)}\\[-2pt]
 {}\cdot|\Delta_1\Delta_2|^{-(d-1)}
 \\[-1pt]\text{after marked-flux}\\[-2pt]\text{cancellation}
 \end{gathered}\)
& other velocities, tree times, all retained normals, one root position \\
\bottomrule
\end{tabular}
\end{center}
Since $|S_U|=|\mathcal P|-1$, the last row has the factor stated later in
the normalization calculation.  The retained time domain is not replaced by
$\prod_{a\in\mathcal T_{\rm ret}}I_a$: after the root and pivot
substitutions it is
\begin{equation}\label{eq:packet-retained-time-domain}
 \left\{t_{\mathcal T_{\rm ret}}:
  \bigl(t_{\mathcal T_{\rm ret}},
 t_{\mathcal T_{\rm piv}}^{\rm piv},
  t_{\mathcal T_{\rm root}}^{\rm in}\bigr)\in\mathcal D_H\right\}.
\end{equation}
Every $t_a^{\rm in}$ is the selected incoming root on its fixed lift, while
$t_{d_i}^{\rm piv}$ is recovered by the mixed implicit map.  Thus the
retained domain is a subset of, and carries explicit witnesses for, the
existential projection of the original Source-B domain.

\begin{proof}
The crossing bonds supplied by assumption (c), one on each physical particle
line, form a transversal cut: removing them separates the upper and lower
directed histories and meets every particle line exactly once.  Take their
velocities as free variables.  Starting from this cut and moving upward,
the two bottom velocities at an upper C-atom determine the complete top pair
by the inverse elastic reflection.  Moving downward, the top pair at a lower
C-atom determines the complete bottom pair.  The molecule is acyclic in its
directed time order, so every nontransversal velocity is reached exactly
once.  Each pair map is orthogonal and has determinant of absolute value one.
At $u_i$ use the top pair as the recovered pair and at $d_i$ use the bottom
pair, so all crossing velocities remain transversal variables.  The collision normal is
either a retained sphere variable or the normalized contact residual on the
support of a scalar root measure; in either description the complete
reflection is an orthogonal pair map.  This proves (i).

With velocities and times treated as parameters, each transport delta has
the form $\delta^{(d)}(x'-x+t(v'-v))$.  Traverse the same two directed
histories from the transversal cut.  Each delta is a translation in the new
position variable and therefore has determinant one.  This proves (ii).

For every contact outside $S_U\cup\{d_1,d_2\}$, take the incoming pair on the
causal side on which that pair has already been recovered.  Before imposing
the contact, it is independent of that atom's own time and has the fixed-lift
equation
\begin{equation}\label{eq:remaining-time-root}
 |a+tg-m|=\eps,\qquad (a+tg-m)\cdot g<0.
\end{equation}
The derivative of $|a+tg-m|$ at a nongrazing root is $g\cdot\omega$.
The atom kernel contains $[g\cdot\omega]_-$, so one-dimensional coarea gives
an operator of norm at most one.  Lemma~\ref{lem:one-incoming-root} gives at
most one incoming root on the fixed lift.  Process these roots in their
causal order on each side of the separator.  Their scattering maps may enter
the coefficients of the later joint map, but
Lemma~\ref{lem:joint-upper-landing-minor} was formulated after precisely
these root substitutions and includes all such co-tree collisions.  This
proves (iii).

What remains is the rectangular vector constraint map
\eqref{eq:joint-constraint-map}.  Apply the sphere identity
\eqref{eq:sphere-delta-identity} at every $a\in S_U$ and at $d_1,d_2$.
Lemma~\ref{lem:joint-upper-landing-minor} permits joint coarea in the
$d(P-1)$ cluster translations and, on each maximal-minor cell, $2(d-1)$ selected
scalar velocity coordinates and the two landing times.  It leaves the other
separator velocities, one common translation, the upper-tree times and all
displayed normals, and has, after cancellation by the two marked fluxes,
factor
\[
 C(d,P_0)\eps^{(d-1)(P-1)+2(d-1)}|\Delta_1\Delta_2|^{-(d-1)}
 \le C(d,P_0)\eps^{(d-1)(P-1)+2(d-1)}\epstar^{-2(d-1)}.
\]
The common torus translation is integrated over one fixed fundamental
domain.  This proves (iv) and the asserted submersion.

At every substitution the full time-order factor is retained.  It becomes
$\mathbf 1_{\mathcal D_H}$ in \eqref{eq:Qsharp-definition}, and the retained
times therefore lie in \eqref{eq:packet-retained-time-domain}.  The coarea
measure of a submersion is unique; hence smooth approximate identities
converge to this measure independently of the displayed elimination order.
Exact grazing has zero flux, and simultaneous-collision sets and ordinary
geometric chart boundaries are null.  Dyadic, lift, sign, and maximal-minor
ties are instead retained and assigned to exactly one Borel half-open cell.
This proves the lemma.
\end{proof}

\subsection{Arbitrary positive kernels and normalization}

\begin{proof}[Proof of Theorem~\ref{thm:two-landing}]
The unordered first two lower edges are different by
Lemma~\ref{lem:upper-tree-landings}.  Theorem~\ref{thm:packet-fct} and
Lemma~\ref{lem:global-packet-chart} give on both the disjoint and overlapping
cells
\[
 J_{H,\kappa}(Q)
 \le \cL_\eps^{C|H|}\eps^{d-1}\epstar^{-2(d-1)}
       \int Q_{H,\kappa}^{\sharp},
\]
which is exactly \eqref{eq:two-landing-cell} because
$\eps^{-(d-1)}\eps^{2(d-1)}=\eps^{d-1}$.  In both cases all operations are performed on one
nonnegative integrand.  Elastic and transport maps have unit Jacobian; every
unmarked root is cancelled by its own sphere--flux factor; and the marked
tree plus the two landing spheres contribute
\[
 \eps^{-(d-1)P}\eps^{(d-1)(P-1)}\eps^{2(d-1)}
       \epstar^{-2(d-1)}
 =\eps^{-(d-1)}\eps^{2(d-1)}\epstar^{-2(d-1)}.
\]
The two chart constructions retain the full complement inside $Q$ and
neither factorize a marginal nor take a supremum over boundary data.  Their
local maximal minors are completed by the uniform algebraic multiplicity
theorem: every regular solution occurs in exactly one of the fixed
$M(d,P_0)$ branch slots in \eqref{eq:Qsharp-definition}.  Since $d,P_0$ are fixed
before $\eps$, this contributes only to $C(d,P_0)$ and changes no power in the
displayed ledger.

The number of slot, sign, lift and
dyadic cells is $\cL_\eps^{C|H|}$.  For
$|H|\le P_{\rm pkt}(d)$ it is absorbed into
$\cL_\eps^C$.  The separator-grid refinement is different: its pieces are
disjoint and the grid label is retained in the output, so their
$Q_{H,\kappa}^{\sharp}$ integrals add with no cardinality factor.  This proves
\eqref{eq:two-landing-actual} and
\eqref{eq:two-landing-excess}.  Finally,
\eqref{eq:root-normalization} proves
\eqref{eq:two-landing-rooted}.  For an independent fallback, disjoint cells
may instead use Proposition~\ref{prop:disjoint-rooted-extension}, while
overlapping cells may use
Proposition~\ref{prop:asymmetric-new-line-peeling}.  Their full $2d$-velocity
payments reproduce the weaker $\epstar^{-2d}$ estimate without using the
radial-conjugate chart.
\end{proof}

\subsection{Core positive-network and finite-fibre interfaces}

The following statements are retained in the main manuscript because they
carry the exact interface from positive compliance and local physical Schur
rows to global finite-fibre recovery.  Their complete proofs, including all
definitions and construction details, remain in Sections~\ref{app:positive-network}--\ref{app:coarea}.
\begin{proposition}[One-port positivity and compliance]
\label{prop:appC-positive-network-closure}
Consider a grounded rooted chronological forest with one retained boundary
port and at most $P_0$ physical lines and junctions.  Allow at most one
restored collision-free chord and any nonnegative collision-curvature ports,
 represented by the factor edge and full pair junctions above, assume that
 the endpoint word is trace-complete in the sense of
 Definition~\ref{def:appC-trace-complete-word}, and
 eliminate all internal position fields.  In the hard-sphere application the
retained port is the full reduced space $E$; a landing is a frame
$Z:\mathbb R^k\to E$.  Let $S_{\mathcal N}$ be
the one-port Dirichlet-to-Neumann stiffness and let
$K_{\mathcal N}:=S_{\mathcal N}^{-1}$ be its grounded
current-to-potential compliance in covariant boundary coordinates.  Let
$Q_-$ be the raw-to-covariant landing-current gauge and let $Q_+$ be the
dual covariant-to-raw potential gauge.  Then
\begin{equation}\label{eq:appC-network-closure}
 K_{\mathcal N}=K_{\mathcal N}^*\ge0,
 \qquad
 s_{\min}(Q_\pm)\ge c(P_0),\qquad
 \|Q_\pm\|\le C(P_0)
\end{equation}
in covariant orthonormal port coordinates.  If $b>0$ and
$F:\mathbb R^k\to E$ is a landing-current frame, the discrete Green identity
gives the exact physical Schur landing row
\begin{equation}\label{eq:appC-network-CD-identification}
 \widetilde K_{\mathcal N}:=
 Q_+K_{\mathcal N}Q_-=Q_-^*K_{\mathcal N}Q_-\ge0,
 \qquad
 W_{\rm Sch}=(bI+\widetilde K_{\mathcal N})F,
 \qquad
 \operatorname{Vol}_k(W_{\rm Sch})
 \ge b^k\operatorname{Vol}_k(F).
\end{equation}
\end{proposition}
\prooflaterinmain{supp:proof:prop:appC-positive-network-closure}
\begin{theorem}[Physical Schur--Jacobi identity for an arbitrary fixed word]
\label{thm:appC-physical-schur-jacobi}
Fix a regular hard-sphere event word $W$ with at most $P_0$ physical lines
and at most $P_0$ events.  Choose unmarked incoming roots $A$, a
chronological marked forest $F$, a collision-free chord $c$ with its two
endpoint constraints, and one retained landing port.  Assume directly that
the chord line has no collision between its two endpoint constraints, that
the rooted forest/chord endpoint word is trace-complete in the sense of
Definition~\ref{def:appC-trace-complete-word}, and that all declared root,
forest, chord and landing pivots are regular.  Roots and forest events may be
interspersed arbitrarily.  These are intrinsic hypotheses on $W$ and make no
reference to a later proposition.  In the full-chord chart let
$p_<=(y_A,h_F,X_c,V_c)$ and
$\mathcal C_<=(C_A,U_F,U_c,R_2)$.  In the radial-conjugate chart of
Section~\ref{appC:radial-conjugate-restoration}, replace these by
\[
 p_<^{\rm rad}=(y_A,h_F,\widehat X_c,g_2,t_{d_2}),
 \qquad
 \mathcal C_<^{\rm rad}
 =(C_A,U_F,U_c,\widetilde R_2^{\rm zs},\rho_2).
\]
Use $(p_<,\mathcal C_<)$ for either square pivot block, let $p_B$ be all
retained port columns, and let $\mathcal C_B=R_1$.  With every
causal derivative retained, set
\begin{equation}\label{eq:appC-full-physical-J}
 J_W=D_{(p_<,p_B)}(\mathcal C_<,\mathcal C_B)
 =\begin{pmatrix}A_W&B_W\\ C_W&D_W\end{pmatrix},
 \quad
 \begin{cases}
 A_W=D_{p_<}\mathcal C_<,\\
 B_W=D_{p_B}\mathcal C_<,\\
 C_W=D_{p_<}\mathcal C_B,\\
 D_W=D_{p_B}\mathcal C_B.
 \end{cases}
\end{equation}
On every regular cell $A_W$ is invertible and the physical retained block is
exactly
\begin{equation}\label{eq:appC-physical-Schur-exact}
 S_{\rm phys}=D_W-C_WA_W^{-1}B_W.
\end{equation}
Equivalently, for a landing multiplier $\mu$ solve the physical bordered
adjoint equations
\begin{equation}\label{eq:appC-physical-bordered-adjoint}
 A_W^*\lambda+C_W^*\mu=0,
 \qquad
 w_{\rm raw}=B_W^*\lambda+D_W^*\mu .
\end{equation}
Then, entry by entry,
\begin{equation}\label{eq:appC-bordered-adjoint-Schur-row}
 \lambda=-A_W^{-*}C_W^*\mu,
 \qquad
 w_{\rm raw}=S_{\rm phys}^*\mu .
\end{equation}

At every fixed-normal junction retain the full pair block
\eqref{eq:appC-pair-dimensional-collision} and eliminate that orthogonal
junction only after both incident line fields have been inserted.  Orient
every physical free segment and every restored chord.  After these
orthogonal junction eliminations their drop maps assemble
to the factor incidence
\begin{equation}\label{eq:appC-factor-incidence}
 \mathbb A_W:\bigoplus_vE_v\longrightarrow
       \bigoplus_eG_e,
 \qquad
 G_e=\begin{cases}
  \mathbb R^{d-1},&e=\gamma_{\rm rad}\text{ in the radial chord chart},\\
  \mathbb R^d,&\text{otherwise}.
 \end{cases}
\end{equation}
Thus $\mathbb A_W$ contains the pair mixing
\eqref{eq:appC-pair-dimensional-collision}; it is not, in general, a stack
of independent single-line differences $\xi_v-O_e\xi_u$.  At every
incoming-root port use the
normal pointing into the excluded collision cylinder (not into the
admissible exterior configuration domain), for which
$\mathcal K=-D\nu|_{\nu^\perp}\ge0$ and put
\begin{equation}\label{eq:appC-root-shunt}
 H_v=2|V_v\cdot\nu_v|\,\mathcal V_v^*\mathcal K_v\mathcal V_v\ge0.
\end{equation}
Then the Hessian of the physical second-variation equations, in the same
covariant gauge and with the same forest/chord variables, is
\begin{equation}\label{eq:appC-word-Jacobi-Hessian}
 L_W=\mathbb A_W^*\operatorname{diag}_{e}
       (\ell_e^{-1}I_{G_e})\mathbb A_W
     +\operatorname{diag}_{v}(H_v).
\end{equation}
More precisely, put
$\Lambda=\operatorname{diag}_e(\ell_eI_{G_e})$, split the endpoint fields as
$\xi=(\xi_I,\xi_B)$, and write
\[
 \mathbb A_W=(A_I\ A_B),
 \qquad
 H_W=\begin{pmatrix}H_{II}&H_{IB}\\H_{BI}&H_{BB}\end{pmatrix}.
\]
Before the edge currents are eliminated the common bordered network system
is the square symmetric system
\begin{equation}\label{eq:appC-word-KKT-system}
 \mathscr K_W^{I,B}
 \begin{pmatrix}j\\\xi_I\\\xi_B\end{pmatrix}
 =\begin{pmatrix}0\\0\\i_B\end{pmatrix},
 \qquad
 \mathscr K_W^{I,B}=
 \begin{pmatrix}
  -\Lambda&A_I&A_B\\
  A_I^*&H_{II}&H_{IB}\\
  A_B^*&H_{BI}&H_{BB}
 \end{pmatrix}.
\end{equation}
The first row is the collection of edge laws
$\ell_ej_e=(\mathbb A_W\xi)_e$; the next two rows are internal and boundary
Kirchhoff balance.  Eliminating $j$ gives exactly
$L_W\xi=(0,i_B)$.  The prefix relation induction in the proof constructs
$\mathbb A_W$ and $H_W$ directly from the physical event word and proves
that the boundary response of \eqref{eq:appC-word-KKT-system} is the bordered
physical relation \eqref{eq:appC-physical-bordered-adjoint}.
In particular, every
causal off-diagonal derivative created by an interspersed root is carried by
the later factor rows; no unnamed bordered matrix is used.
After grounding one translation in every free component, let $I$ be the
internal fields and $B$ the retained port.  The same Gaussian quotient as
in \eqref{eq:appC-physical-Schur-exact} gives
\begin{equation}\label{eq:appC-same-word-quotient}
 S_{\mathcal N}=L_{BB}-L_{BI}L_{II}^{-1}L_{IB},
 \qquad K_{\mathcal N}=S_{\mathcal N}^{-1}\ge0.
\end{equation}
Let $\mathsf P_W:E_{\rm raw}\to E_{\rm cov}$ be the square primal potential
coordinate map obtained from the complete pair junctions and the forest
splitting
\eqref{eq:appC-forest-coordinate-injection}--
\eqref{eq:appC-forest-boundary-splitting}.  Define, rather than merely name,
the two boundary gauges by
\eqref{eq:appC-exact-dual-boundary-gauges}.  They are square and satisfy
\eqref{eq:appC-positive-word-coordinate-bounds}.  For the one-speed landing
input, the raw frame $Z_1$ produces the covariant current $Q_-Z_1$.  The
internal network contributes the raw potential
$Q_+K_{\mathcal N}Q_-Z_1$, while the lower free interval contributes $bZ_1$.
The exact physical identity is
\begin{equation}\label{eq:appC-physical-network-final-identity}
 S_{\rm phys}^*
 =(bI+\widetilde K_{\mathcal N})Z_1,
 \qquad
 \widetilde K_{\mathcal N}
 :=Q_+K_{\mathcal N}Q_-=Q_-^*K_{\mathcal N}Q_-\ge0.
\end{equation}
Here the star is the Hilbert adjoint, so both sides of
\eqref{eq:appC-physical-network-final-identity} are maps from the $d$
landing rows to the retained port-column space.
Thus the positive network is the Schur quotient of the same physical word,
not a comparison network.
\end{theorem}
\prooflaterinmain{supp:proof:thm:appC-physical-schur-jacobi}
\begin{theorem}[Physical multiport Schur--Jacobi identity for birth segments]
\label{thm:appC-birth-multiport-schur-jacobi}
Use the hypotheses and notation of
Lemma~\ref{lem:appC-causal-birth-current-clearing}.  In the remainder word
retain every physical free segment other than the interiors of the
$\gamma_i$, every complete pair junction, every marked-forest constraint,
every complementary incoming-root shunt, and the declared common-
translation quotient.  Its Jacobi Hessian is
\begin{equation}\label{eq:appC-birth-remainder-Hessian}
 L_{\rm rem}=\mathbb A_{\rm rem}^*
       \Lambda_{\rm rem}^{-1}\mathbb A_{\rm rem}+H_{\rm rem},
 \qquad H_{\rm rem}\ge0.
\end{equation}
Here the prefix construction gives
$H_{\rm rem}=\mathbb C_{\rm rem}^*\mathbb C_{\rm rem}$ with
$\mathbb C_{\rm rem}$ the direct stack of the transported physical
incoming-root shunt factors.  Put $K_{\rm rem}:=L_{\rm rem}^{+}$.

Let $\iota_i:\mathbb R^d\to\mathcal Y$ be the $i$th block inclusion and
$F_i:=F\iota_i$.  With the extended factor
$\mathbb B_{\rm rem}$ from
\eqref{eq:appC-birth-extended-factor}, there are augmented physical
edge--shunt flows $\mathcal J_i$ such that
\begin{equation}\label{eq:appC-birth-current-divergence}
 F_i=\mathbb B_{\rm rem}^*\mathcal J_i,
 \qquad 1\le i\le k.
\end{equation}
Consequently the cleared physical injection satisfies
\begin{equation}\label{eq:appC-birth-current-range}
 F(\mathcal Y)\subseteq\operatorname{Ran}L_{\rm rem},
 \qquad K_{\rm rem}=K_{\rm rem}^*\ge0
       \quad\hbox{on this range}.
\end{equation}
The physical cut relation and the remainder Kirchhoff relation are the
two rows of the same symmetric bordered system
\begin{equation}\label{eq:appC-birth-common-radial-KKT}
\mathscr K_{\Gamma}^{\rm rad}=
\begin{pmatrix}
  D_\Delta&\mathcal G&F^*\\
  \mathcal G^*&0&0\\
  F&0&-L_{\rm rem}
\end{pmatrix}.
\end{equation}
In particular, after quotienting the $k$ landing-time directions, choose
orthogonal half-open frames
\[
 E_i:\mathbb R^{d-1}\longrightarrow g_i^\perp,
 \qquad E=\operatorname{diag}_{i}E_i.
\]
On the regular set $g_i\ne0$, the exact tangential response is
\begin{equation}\label{eq:appC-birth-tangential-response}
 M_{\Gamma,\perp}
 :=E^*D_\Delta E+(FE)^*K_{\rm rem}(FE)
 \succeq \operatorname{diag}_{i}(\Delta_iI_{d-1}),
\end{equation}
and hence
\begin{equation}\label{eq:appC-birth-tangential-determinant}
 \det M_{\Gamma,\perp}
 \ge\prod_{i=1}^k\Delta_i^{d-1}.
\end{equation}

Let $\mathcal S_k$ and $\mathcal T_k$ be the root- and forest-Schur
landing derivatives in the original Cartesian velocity and landing-time
coordinates.  Then
\begin{equation}\label{eq:appC-birth-mixed-frame}
 \Jac_{k(d-1)\mid k}(\mathcal S_k,\mathcal T_k)
 \ge c(d,P_0,k)
       \prod_{i=1}^k |g_i|\Delta_i^{d-1}
 \ge c(d,P_0,k)
       \prod_{i=1}^k |g_i\cdot\omega_i|\Delta_i^{d-1}.
\end{equation}
No comparison of the $\Delta_i$ is used.
\end{theorem}
\prooflaterinmain{supp:proof:thm:appC-birth-multiport-schur-jacobi}
\begin{theorem}[Uniform algebraic multiplicity and branch atlas]
\label{thm:appD-uniform-algebraic-multiplicity}
On every regular packet cell with \(|H|\le P_0\), and for every fixed output
value \(b\), the pivot fibre
\begin{equation}\label{eq:appD-pivot-fibre}
 \mathfrak F_{H,\kappa}(b)
 :=\{p:\mathscr F_{H,\kappa}(p;b)=0,\ (p,b)
       \text{ belongs to the declared cell}\}
\end{equation}
has at most
\begin{equation}\label{eq:appD-multiplicity-bound}
 M(d,P_0):=(4P_0+4)^{(2d+1)P_0+2d}
\end{equation}
regular points.  Moreover the regular zero set has Borel recovery branches
\begin{equation}\label{eq:appD-recovery-branches}
 p=\Phi_{H,\kappa,j}(b),\qquad 1\le j\le M(d,P_0),
\end{equation}
where an absent branch is assigned \(\dagger\).  The number and the domain
of the branch slots are independent of all concrete boundary values.
\end{theorem}
\prooflaterinmain{supp:proof:thm:appD-uniform-algebraic-multiplicity}
\begin{corollary}[Uniform multiplicity for fixed-$k$ birth recovery]
\label{cor:appD-fixed-k-birth-multiplicity}
Fix $k$ and use the simultaneous birth-landing recovery system of
Theorem~\ref{thm:birth-flag-frame-production}.  If the word has $P\le P_0$
physical lines and $|H|\le P_0$ events, then its number of pivot variables
and equations is at most
\begin{equation}\label{eq:appD-fixed-k-birth-variable-count}
 N_k(d,P_0):=d(P_0-1+k)+(d+1)P_0,
\end{equation}
every equation has degree at most $D_0(P_0)=4P_0+4$, and every regular
physical output fibre has at most
\begin{equation}\label{eq:appD-fixed-k-birth-multiplicity}
 M_k(d,P_0):=(4P_0+4)^{N_k(d,P_0)}
\end{equation}
points.  These points admit $M_k(d,P_0)$ boundary-independent Borel branch
slots
\begin{equation}\label{eq:appD-fixed-k-birth-branches}
 p=\Psi^{(k)}_{H,\kappa,j}(\xi,b,\zeta),
 \qquad 1\le j\le M_k(d,P_0),
\end{equation}
padded by $\dagger$.
\end{corollary}
\prooflaterinmain{supp:proof:cor:appD-fixed-k-birth-multiplicity}
\begin{definition}[Selected Borel conditional representative]
\label{def:appD-selected-conditional-representative}
For a fixed pair $(H,\kappa)$, let $\lambda_Y$ be the product reference
measure of the complement and fixed-boundary variables $y$.  The kernel
$K_{H,\kappa}$ below disintegrates the source-normalized packet measure
defining $J_{H,\kappa}$ (that is, the measure
\eqref{eq:appD-packet-collision-measure} after the source normalization
specified there).  Write the
fixed branch product as
$Y\times\widehat{\mathfrak B}_{H,\kappa}^{\rm br}$, so that $b$ below
denotes only the branch coordinates not contained in $y$.  The packet
collision measure, its exact branch push-forward
\eqref{eq:appD-exact-pushforward-measure}, and the dominating output measure
\eqref{eq:appD-dominating-output-measure} are Radon measures on standard
Borel product spaces.  Fix Borel disintegrations
\[
 y\longmapsto K_{H,\kappa}(y,\dd z),
 \qquad
 y\longmapsto K^{\rm ex}_{H,\kappa}(y,\dd b,\dd j),
 \qquad
 y\longmapsto K^\sharp_{H,\kappa}(y,\dd b,\dd j)
\]
with respect to $\lambda_Y$, once and for all and independently of the
test function.  Every recovery branch fixes the $Y$-coordinate (the
complement and fixed-boundary variables are not pivots), so the exact
push-forward identity is an identity over the same base $Y$.  Because all
three fibre spaces are standard Borel, choose
countable determining algebras in them.  On each regular compact exhaustion,
apply the integrated exact push-forward identity and the measure domination
to those algebras; take the countable union of the resulting null sets and
then pass to the monotone limit.  The monotone-class theorem gives one Borel
$\lambda_Y$-null set,
still denoted by $N_\partial$, outside which, simultaneously for every
Borel set $A$ in the collision fibre and every Borel set $B$ in the branch
fibre,
\begin{align}
 K_{H,\kappa}(y,A)
 &=\int
   \1_A\!\left(\Phi_{H,\kappa,j}(y,b)\right)
   K^{\rm ex}_{H,\kappa}(y,\dd b,\dd j),
 \label{eq:appD-conditional-exact-kernel}\\
 K^{\rm ex}_{H,\kappa}(y,B)
 &\le
 C(P_0)\cL_\eps^{C|H|}\eps^{d-1}\epstar^{-\gamma_\kappa}
 K^\sharp_{H,\kappa}(y,B).
 \label{eq:appD-conditional-dominating-kernel}
\end{align}
Enlarge this set once more by the null set of
Lemma~\ref{lem:appD-conditional-null-strata}, and put
\begin{align}
 \overline K_{H,\kappa}(y,\dd z)
 &:=\1_{N_\partial^c}(y)K_{H,\kappa}(y,\dd z),\label{eq:appD-selected-K}\\
 \overline K^\sharp_{H,\kappa}(y,\dd b,\dd j)
 &:=\1_{N_\partial^c}(y)K^\sharp_{H,\kappa}(y,\dd b,\dd j).
 \label{eq:appD-selected-Ksharp}
\end{align}
For every nonnegative Borel $Q$, define
\begin{align*}
 \overline J_{H,\kappa}^{y}(Q)
 &:=\int Q(y,z)\,\overline K_{H,\kappa}(y,\dd z),\\
 \overline J_{H,\kappa}^{\sharp,y}(Q)
 &:=\int Q\!\left(y,\Phi_{H,\kappa,j}(y,b)\right)\,
        \overline K^\sharp_{H,\kappa}(y,\dd b,\dd j).
\end{align*}
For an absent branch the last integrand is defined to be zero.  Notice that
$K^\sharp$ already contains the density
$\1_{\mathcal D_H}W^{\rm rem}$ from
\eqref{eq:appD-dominating-output-measure}; using $Q^\sharp$ as its test
function would count that density twice.
We call these the selected conditional operators.  The unadorned
$J_{H,\kappa}$ always denotes the complete operator after integration in
$y$.
\end{definition}
\begin{proposition}[Arbitrary-kernel packet coarea]
\label{prop:appD-arbitrary-kernel-coarea}
On every regular cell, for every fixed source-boundary value $y$ and every
nonnegative joint \(Q\), the selected representatives of
Definition~\ref{def:appD-selected-conditional-representative} satisfy
\begin{equation}\label{eq:appD-arbitrary-kernel-coarea}
 \overline J_{H,\kappa}^{\,y}(Q)
 \le
 C(P_0)\cL_\eps^{C|H|}
 \eps^{-(d-1)}\eps^{2(d-1)}\epstar^{-\gamma_\kappa}
 \overline J_{H,\kappa}^{\sharp,y}(Q).
\end{equation}
Equivalently, the original conditional disintegrations satisfy this bound
for $\lambda_Y$-almost every $y$.  After integration in $y$, the same
inequality holds for the complete operator $J_{H,\kappa}$ by
\eqref{eq:appD-selected-representative-integral} and monotone exhaustion.
The branch domain on the right is independent of the concrete value of
$y$.
\end{proposition}
\prooflaterinmain{supp:proof:prop:appD-arbitrary-kernel-coarea}
\begin{lemma}[Indicator-preserving fibre inequality]
\label{lem:appD-fibre-time-owner}
Disintegrate the two selected measure kernels of
Definition~\ref{def:appD-selected-conditional-representative}, together with
the auxiliary exact-branch kernel
$\1_{N_\partial^c}K^{\rm ex}_{H,\kappa}$, once more in the $S$-variables.
Using countable determining algebras before choosing $Q$, fix one common
Borel null set $N_S$, independent of $Q$, outside which the conditional
exact push-forward identity and conditional measure domination hold, and
zero-extend all three kernels on $N_S$.  Also include in $N_S$ the Tonelli
null set of those $S$-values for which the remaining fibre of
$N_\partial$ has nonzero reference measure.  Consequently, off $N_S$, the
conditional integral of the selected output kernel equals the bare
$Q^\sharp$ branch integral displayed below.
Then, for every fixed value of the \(S\)-variables and every nonnegative
$Q$,
\begin{multline}\label{eq:appD-fibre-time-owner}
 \overline J_{H,\kappa}^{\,t_S}\!\left(
   Q\,\1_{\mathcal D_H}\right)(t_S)\\
 \le
 C(P_0)\cL_\eps^{C|H|}
 \eps^{-(d-1)}\eps^{2(d-1)}\epstar^{-\gamma_\kappa}
 \1_{N_S^c}(t_S)
 \1_{\widetilde{\mathcal D}_{H,S}}(t_S)
 \sum_{j=1}^{M(d,P_0)}
 \int Q_{H,\kappa}^{\sharp}(b,j;t_S)\dd b .
\end{multline}
Here the left side is the selected collision representative, while the
factor $\1_{N_S^c}$ makes the displayed branch-output side its simultaneous
selected representative.  The
projection on the right is exactly
\eqref{eq:appD-time-projection}; no enlarged time domain replaces it.
\end{lemma}
\prooflaterinmain{supp:proof:lem:appD-fibre-time-owner}
\subsection{Fixed-word multi-landing transfer and birth-flag production}
\label{sec:transfer-principle}

The fixed-$k$ abstraction has two layers: the source-specific construction of
a uniformly transverse landing frame, and the dimension-independent chain of
root cancellation, finite algebraic multiplicity and arbitrary-kernel coarea.
Theorem~\ref{thm:fixed-word-multilanding-transfer} gives the conditional
transfer once a landing frame is available.
Theorem~\ref{thm:birth-flag-frame-production} constructs the birth-flag
frame in the present packet class, and
Corollary~\ref{cor:automatic-fixed-k-birth-production} is the form used in
Section~\ref{sec:global}.
These results do not constitute an unconditional extension to the full model.

We fix the coarea notation independently of the Source-B molecule language.
A finite event word is denoted by $H$, and $|H|$ is its number of
events, not its number of physical lines.  On a half-open regular cell
$\chi$, let $p\in\mathfrak P_{H,\chi}\subset\mathbb R^N$ be the declared
pivot variables, $b\in\mathfrak B_{H,\chi}$ the retained packet variables,
$\zeta\in\mathfrak Z_{H,\chi}$ the retained complement variables, and let
$\xi\in\mathfrak X_{H,\chi}$ be the concrete source-boundary parameter.  The
boundary space is equipped with its source measure $\dd\rho(\xi)$.  All
three pivot/retained boxes, with their Lebesgue--surface product measures,
are fixed independently of $\xi$.

The packet is defined by a square constraint map
\[
 \mathcal C_{H,\chi}(p;b,\zeta,\xi)=0\in\mathbb R^N,
\]
not by assuming that the retained variables are the value of a globally
invertible map of $p$.  With the delta interpreted as the regular coarea
measure, put
\begin{equation}\label{eq:abstract-J-definition}
 J_{H,\chi}^{\xi}(Q):=
 \int Q(p,b,\zeta)
 \1_{\mathcal D_{H,\chi}^{\xi}}(p,b,\zeta)
 W_{H,\chi}(p,b,\zeta;\xi)
 \delta^{(N)}(\mathcal C_{H,\chi})
 \dd p\dd b\dd\nu(\zeta).
\end{equation}
Here $W_{H,\chi}$ contains the declared normalization, sphere measures,
collision fluxes, Boltzmann weights and cutoffs.  The kernel $Q\ge0$ is an
arbitrary Borel function of all packet and complement variables and need not
factor through $(b,\zeta)$.

The physical recovery relation is
\begin{equation}\label{eq:abstract-recovery-relation}
 \mathscr R_{H,\chi}
 :=\{(\xi,b,\zeta,p):(p,b,\zeta)\in\mathcal D_{H,\chi}^{\xi},
       \ \mathcal C_{H,\chi}(p;b,\zeta,\xi)=0\}.
\end{equation}
If $p=\Psi_{H,\chi,j}(\xi,b,\zeta)$ is a recovery branch, let
$W_{H,\chi,j}^{\rm rem}(\xi,b,\zeta)$ be the uncancelled positive density after
the tree, root and landing changes of variables, and define
\begin{equation}\label{eq:abstract-Qsharp-definition}
 Q_{H,\chi,j}^{\sharp,\xi}(b,\zeta)
 :=Q(\Psi_{H,\chi,j}(\xi,b,\zeta),b,\zeta)
       \1_{\mathscr R_{H,\chi}}
       (\xi,b,\zeta,\Psi_{H,\chi,j}(\xi,b,\zeta))
       W_{H,\chi,j}^{\rm rem}(\xi,b,\zeta),
\end{equation}
extended by zero when the branch is absent.  Thus
$\int Q_{H,\chi,j}^{\sharp,\xi}$ below always means integration over the fixed
domain $\mathfrak B_{H,\chi}\times\mathfrak Z_{H,\chi}$ against
$\dd b\dd\nu(\zeta)$, and no supremum over boundary data is implicit.

\begin{theorem}[Fixed-word multi-landing transfer]
\label{thm:fixed-word-multilanding-transfer}
Fix integers $d\ge2$, $P_{\rm line}\ge2$, $W_{\rm ev}\ge1$, and $k\ge1$.
Consider a regular cell $\chi$ of a connected hard-sphere packet whose word
$H$ contains $|H|\le W_{\rm ev}$ events and meets
$P\le P_{\rm line}$ physical lines.  Assume that the packet is full and its
normalization inside $W_{H,\chi}$ is exactly $\eps^{-(d-1)P}$, and that:
\begin{enumerate}[label=\textup{(T\arabic*)}]
\item $P-1$ marked contacts form an upper spanning tree, $k$ further
contacts are selected as landings, and chronological relative-cluster
coordinates give a square tree block $B_T$ satisfying
\[
 |\det B_T|\ge c_T(d,P_{\rm line})>0;
\]
\item every unmarked contact has a transversal incoming time--normal root,
whose exact Jacobian is paired with its own sphere measure and collision
flux;
\item after the marked tree and all earlier root blocks are eliminated, the
joint $kd$-row landing derivative has
\begin{equation}\label{eq:abstract-k-landing-rank}
 \Jac_{kd}(\mathcal S_{\rm land})\ge\lambda>0;
\end{equation}
\item after extrinsic normal variables are introduced,
the constraint $\mathcal C_{H,\chi}=0$ is represented by $N$ polynomial
equations in the same number of pivot coordinates, where
\[
 N\le N_0(d,P_{\rm line},W_{\rm ev},k),
 \qquad
 \deg\mathcal C_{H,\chi}\le D_0(d,P_{\rm line},W_{\rm ev},k).
\]
Their coefficients are Borel
(indeed semialgebraic) functions of $(\xi,b,\zeta)$.  Every physical pivot in
\eqref{eq:abstract-recovery-relation} is a zero of this square system and has
nonzero full pivot Jacobian.  Polynomial zeros failing the original sphere,
sign, time-order, lift or half-open cell conditions are discarded by those
Borel conditions; no physical branch is omitted;
\item the relation \eqref{eq:abstract-recovery-relation} is a Borel subset of
the displayed standard Borel product, the retained variables range in the
fixed boundary-independent domains specified above, and the cell family has
cardinality at most $\cL_\eps^{C|H|}$.  Maximal-minor and chart ties are
assigned by a disjoint half-open
Borel partition rather than discarded.  The complement of the declared
regular cells is null for the packet measure or is covered by a separately
stated estimate for every $\xi\notin N_\partial$, where
$N_\partial$ is Borel and $\rho(N_\partial)=0$.
\end{enumerate}
Then for every $\xi\notin N_\partial$ and every nonnegative Borel kernel
$Q$ coupling the packet to an arbitrary complement, there is one
boundary-independent enlarged domain and
at most
\begin{equation}\label{eq:abstract-multilanding-multiplicity}
 M_{\rm br}(d,P_{\rm line},W_{\rm ev},k):=
 D_0(d,P_{\rm line},W_{\rm ev},k)^{
       N_0(d,P_{\rm line},W_{\rm ev},k)}
\end{equation}
Borel recovery branches such that
\begin{equation}\label{eq:abstract-multilanding-operator}
 J_{H,\chi}^{\xi}(Q)
 \le C(d,P_{\rm line},W_{\rm ev},k)\,\cL_\eps^{C|H|}
 \eps^{(d-1)(k-1)}\lambda^{-1}
 \sum_{j=1}^{M_{\rm br}(d,P_{\rm line},W_{\rm ev},k)}
 \int Q_{H,\chi,j}^{\sharp,\xi}(b,\zeta)\dd b\dd\nu(\zeta).
\end{equation}
The same statement holds for a finite atlas of successive landing blocks
provided the product of their square Jacobians is at least $\lambda$ and
all blocks belong to one common triangular pivot system.
\end{theorem}

\begin{proof}
Use the exact sphere identity in dimension $d$ at the $P-1$ tree contacts
and the $k$ selected landings.  Their surface factors contribute
$\eps^{(d-1)(P-1+k)}$.  The packet normalization contributes
$\eps^{-(d-1)P}$, leaving exactly
\begin{equation}\label{eq:abstract-multilanding-normalization}
 \eps^{-(d-1)P}\eps^{(d-1)(P-1+k)}
 =\eps^{(d-1)(k-1)}.
\end{equation}
At every unmarked contact, the incoming time--normal determinant cancels
pointwise with that contact's sphere and flux factors by (T2).  These
 cancellations occur inside the same nonnegative integrand and do not require
 factorization of $Q$.
The tree block in (T1) is eliminated next; its inverse determinant is at
most $c_T(d,P_{\rm line})^{-1}$ and is absorbed into the displayed constant.

By Cauchy--Binet, (T3) supplies a finite atlas of $kd$-column square minors;
on each half-open maximal-minor cell the inverse coarea cost is at most
$C\lambda^{-1}$.  Local invertibility alone would not control the number of
preimages.  Instead retain all event equations.  By (T4), every physical
real solution is a regular isolated complex zero of the same bounded-degree
square system.  The system may have nonphysical or positive-dimensional
complex components, but neither causes a problem: only isolated zeros enter
the bound, and the original Borel physical conditions discard spurious real
zeros.
The affine B\'ezout degree inequality
\cite[Theorem~1]{Heintz1983}, read together with its published corrigendum
\cite{Heintz1985}, bounds the zero-dimensional components of a square system
of $N\le N_0$ equations in $\mathbb C^N$, each of degree at most $D_0$, by
$D_0^N\le D_0^{N_0}$.  It therefore gives
\eqref{eq:abstract-multilanding-multiplicity}.  The Lusin--Novikov theorem
\cite[Theorem~18.10]{Kechris1995} enumerates the finite Borel fibres.  By
(T5) they live in the common fixed output domain and every tie has a unique
 Borel owner.  Add the branch label to the output domain, substitute each
 recovery map into both the remaining positive density and the original
 joint kernel $Q(p,b,\zeta)$, and extend the resulting function by zero
 outside its original cell.  The finite-fibre coarea formula applied
 directly to the delta measure in \eqref{eq:abstract-J-definition}, Tonelli,
and
\eqref{eq:abstract-multilanding-normalization} prove
\eqref{eq:abstract-multilanding-operator}.  For successive square blocks,
the derivative of the common recovery map is block lower triangular in
chronological order, and its determinant is the product of the diagonal
block determinants; the preceding coarea and zero-extension argument then
applies to that single joint map.
 The cell-cardinality bound in (T5) gives the displayed
 $\cL_\eps^{C|H|}$ after the disjoint cells are summed; for one fixed cell
 that factor may be omitted.
\end{proof}

For later use, if $A:\mathbb R^N\to\mathbb R^{kd}$ and
$T:\mathbb R^k\to\mathbb R^{kd}$, set
\begin{equation}\label{eq:abstract-general-mixed-jacobian}
 \Jac_{k(d-1)\mid k}(A,T)
 :=\left(\sum_{\substack{I\subset\{1,\ldots,N\}\\|I|=k(d-1)}}
       |\det[A_I\ T]|^2\right)^{1/2}.
\end{equation}
This definition is invariant under orthogonal changes of the $A$-columns
and under determinant-one changes of all $kd$ output rows.

\begin{corollary}[Flux-completed fixed-word transfer]
\label{cor:fixed-word-flux-multilanding-transfer}
In Theorem~\ref{thm:fixed-word-multilanding-transfer}, replace
\textup{(T3)} by the following flux-completed alternative.  The $k$ selected
landing times are pivot columns, $k(d-1)$ further pivot columns are ordinary
 Cartesian velocity coordinates, and, with
$\phi_i=|g_i\cdot\omega_i|$, there is a positive Borel function
$\lambda_{\rm fr}$ on the regular recovery relation such that the joint landing
rows satisfy pointwise
\begin{equation}\label{eq:abstract-flux-completed-frame}
 \Jac_{k(d-1)\mid k}(\mathcal S_{\rm land},\mathcal T_{\rm land})
 \ge \lambda_{\rm fr}\prod_{i=1}^k\phi_i .
\end{equation}
Assume also that the physical selected-landing density contains exactly the
corresponding product $\prod_i\phi_i$.  Write
$M_{\rm br}:=M_{\rm br}(d,P_{\rm line},W_{\rm ev},k)$.  Then the conclusion is
\begin{equation}\label{eq:abstract-flux-multilanding-operator}
\begin{aligned}
 J_{H,\chi}^{\xi}(Q)
 &\le C(d,P_{\rm line},W_{\rm ev},k)\cL_\eps^{C|H|}
       \eps^{(d-1)(k-1)}\\
 &\quad\times\sum_{j=1}^{M_{\rm br}}\int
 \frac{Q_{H,\chi,j}^{\sharp,\xi}(b,\zeta)}
 {\lambda_{\rm fr}(\xi,b,\zeta,
    \Psi_{H,\chi,j}(\xi,b,\zeta))}
 \dd b\dd\nu(\zeta).
\end{aligned}
\end{equation}
An absent branch is interpreted as zero.  In particular, if
$\lambda_{\rm fr}\ge\underline\lambda_{\rm fr}>0$ uniformly, the right-hand
side is bounded by the same expression with
$\underline\lambda_{\rm fr}^{-1}$ pulled in front of the sum.
All other conclusions of the theorem, including arbitrary joint kernels and
boundary-independent branch slots, remain unchanged.
\end{corollary}

\begin{proof}
Use the same common square recovery system as in the theorem, now with the
$k$ landing-time columns and $k(d-1)$ Cartesian velocity columns in the
landing block.  Cauchy--Binet and
\eqref{eq:abstract-flux-completed-frame} cost
$\lambda_{\rm fr}^{-1}\prod_i\phi_i^{-1}$ pointwise on each recovery
branch.  Each inverse flux is cancelled
pointwise by the flux of that same selected collision.  No flux belonging to
another contact is used.  The sphere and normalization ledger is unchanged:
\[
 \eps^{-(d-1)P}\eps^{(d-1)(P-1+k)}
 =\eps^{(d-1)(k-1)}.
\]
The bounded-degree system has the same number of pivot variables because
$k(d-1)+k=kd$.  The finite-fibre, Borel enumeration, zero extension and
arbitrary-kernel parts of the proof of
Theorem~\ref{thm:fixed-word-multilanding-transfer} therefore apply verbatim
to the same nonnegative integrand.
\end{proof}

\begin{corollary}[The present packet as the first noncomponentwise case]
\label{cor:d4-transfer-instance}
For $\rho$-almost every source-boundary parameter, the packet already proved
in Theorem~\ref{thm:two-landing} satisfies the hypotheses of
Theorem~\ref{thm:fixed-word-multilanding-transfer} with
\[
 (d,k,P_{\rm line},W_{\rm ev})=(d,2,2P_0,P_0),
 \qquad \lambda\ge c(d,P_0)\epstar^{2d}.
\]
 Consequently its abstract transfer estimate is valid on those slices and,
 after Tonelli, for the integrated source operator:
\[
 J_H(Q)\le C(P_0)\cL_\eps^{C|H|}
 \eps^{d-1}\epstar^{-2d}\int Q_H^\sharp,
\]
and the gain $\eps^{d-1}$ is the structural identity
$(d-1)(k-1)=d-1$.
\end{corollary}

\begin{proof}
Here $|H|\le P_0$ is an atom/event cap, so the word contains at most $P_0$
events and meets at most $2P_0$ physical lines.  The two landing-frame
 constructions give the stated $\lambda$.  Lemma
\ref{lem:appD-bounded-degree-word}, Theorem
\ref{thm:appD-uniform-algebraic-multiplicity} and Proposition
 \ref{prop:appD-arbitrary-kernel-coarea} verify the square recovery system,
 Borel branches and arbitrary-kernel output measure.  Substitution in
\eqref{eq:abstract-multilanding-operator} gives the display.  Logically, this
corollary is a structural reclassification of the already proved packet and
is not used to prove Theorem~\ref{thm:two-landing}; hence the references to
 Section~\ref{app:coarea} create no circular dependence.  The almost-everywhere boundary
 qualifier is exactly Lemma~\ref{lem:appD-conditional-null-strata}; Tonelli
 is used before the source-boundary variables are integrated, and no
 pointwise assertion is made on its exceptional set.
\end{proof}

\subsubsection{Birth-flag production for arbitrary fixed multiplicity}

\begin{definition}[Canonical lower birth flag]
\label{def:canonical-lower-birth-flag}
Assume that every physical line $p$ of a two-sublayer packet occurs in the
lower sublayer.  In the strict lower atom order put
\[
 \lambda_D(p):=\min\{d\in H_D:p\text{ is incident to }d\},
 \qquad B(d):=\{p:\lambda_D(p)=d\}.
\]
The nonempty $B(d)$ are the lower birth atoms.  The canonical $k$-birth
flag consists of the chronological first $k$ such atoms
$d_1<\cdots<d_k$ and, at each $d_i$, the line $c_i\in B(d_i)$ with the
smallest permanent physical-line label.  Put
\[
 \Delta_i=t_{d_i}-s,
 \quad g_i=\hbox{incoming relative velocity at }d_i,
 \quad \phi_i=|g_i\cdot\omega_i|.
\]
\end{definition}

\begin{lemma}[Birth count, causal segments and sharpness]
\label{lem:lower-birth-count}
If the packet has $P$ physical lines, the nonempty sets $B(d)$ partition
those lines and
\begin{equation}\label{eq:lower-birth-count}
 P=\sum_{d:B(d)\ne\varnothing}|B(d)|
   \le2b_D(H),
 \qquad b_D(H)\ge\left\lceil\frac P2\right\rceil.
\end{equation}
Consequently $P\ge2k-1$ guarantees the canonical $k$-birth flag.  Its
selected lines are distinct, and $c_i$ has no lower collision in
$(s,t_{d_i})$.  The selection is deterministic and cellwise Borel and
introduces no combinatorial summation.

The bound is sharp for first-appearance counting: take pairwise disjoint
lower atoms which each introduce two new lines and then add later bridge
atoms to connect their components.
\end{lemma}

\begin{proof}
Every line has a unique first lower atom, so the $B(d)$ are disjoint and
exhaust all $P$ lines.  A binary atom has two physical-line slots and hence
$|B(d)|\le2$, proving \eqref{eq:lower-birth-count}.  Different birth sets
are disjoint, so their selected representatives are distinct.  By the
definition of $\lambda_D(c_i)$, the line $c_i$ has no lower atom before
$d_i$; since the separator lies after every upper atom, its segment from
$s$ to $d_i$ is a genuine collision-free physical segment.  Strict atom
order and permanent line labels make the rule Borel on every half-open
cell.  The displayed construction has exactly one birth atom per two new
lines before its later bridges, proving sharpness.
\end{proof}

\begin{lemma}[Causal reverse-Kruskal upper tree]
\label{lem:birth-reverse-Kruskal-tree}
Under assumptions \textup{(a)--(b)} of
Theorem~\ref{thm:two-landing}, scan all upper collision edges from latest to
earliest and retain an edge precisely when it joins two distinct components
of the already retained forest.  The result is a spanning tree $T_K$ on the
$P$ physical lines, with exactly $P-1$ edges.  Moreover, for every rejected
upper contact $a$, the endpoints of $a$ are joined in $T_K$ by a path all of
whose contacts are strictly later than $a$.  The choice is deterministic on
each strict-order cell and hence is cellwise Borel.
\end{lemma}

\begin{proof}
The upper atom graph is connected by assumption \textup{(a)}.  Two adjacent
upper atoms share a physical line, and assumption \textup{(b)} rules out an
isolated physical-line vertex.  Thus the upper physical-line multigraph is
connected and spanning.  Reverse Kruskal on a connected finite multigraph
ends with a spanning tree, so it has $P-1$ edges.  When $a$ is rejected, its
endpoints are already connected by retained edges scanned earlier, which are
exactly contacts later in physical time.  The strict atom order fixes every
decision and proves the Borel assertion.
\end{proof}

For selected landing residuals write, in the top-pair representation,
\begin{equation}\label{eq:k-birth-landing-residuals}
 R_i=X_{c_i}+\Delta_iV_{c_i}
       -Y_i(t_{d_i};X,V)-m_i,
 \qquad \mathcal R_k=(R_1,\ldots,R_k).
\end{equation}
The actual selected contact equations are
\begin{equation}\label{eq:k-birth-contact-equations}
 C_i:=R_i-\eps\omega_i=0,
 \qquad C=(C_1,\ldots,C_k),
 \qquad \Omega=(\omega_1,\ldots,\omega_k).
\end{equation}
The selected normals are retained variables in the recovery chart.  Hence
the velocity and landing-time derivatives of $C$ equal those of
$\mathcal R_k$; all zero-set row operations below are nevertheless applied
to $C$, not to the nonzero vector $R=\eps\Omega$.
Let $h$ be the upper-tree relative-cluster translations and let $\mathcal U$
be the marked upper-tree residuals.  After all complementary incoming roots
are substituted, put
\begin{equation}\label{eq:k-birth-Schur-derivatives}
 \mathcal S_k=D_V\mathcal R_k-
 D_h\mathcal R_k(D_h\mathcal U)^{-1}D_V\mathcal U,
 \qquad
 \mathcal T_k=D_{(t_{d_1},\ldots,t_{d_k})}\mathcal R_k.
\end{equation}
The mixed Jacobian below is the one in
\eqref{eq:abstract-general-mixed-jacobian}, with $A=\mathcal S_k$ and
$T=\mathcal T_k$.

\begin{lemma}[Exact fixed-$k$ time-wedge factorization]
\label{lem:k-time-wedge-factorization}
Let $A:\mathbb R^N\to\mathbb R^{kd}$ and let
$T:\mathbb R^k\to\mathbb R^{kd}$ have rank $k$.  If
$Q_T:\mathbb R^{k(d-1)}\to(\operatorname{Ran}T)^\perp$ is orthogonal, then
\begin{equation}\label{eq:k-time-wedge-factorization}
 \Jac_{k(d-1)\mid k}(A,T)
 =\sqrt{\det(T^*T)}\,
       \Jac_{k(d-1)}(Q_T^*A).
\end{equation}
In particular, for
$\mathcal G\tau=(g_i\tau_i)_{i=1}^k$ and
$E=\operatorname{diag}_iE_i$ with
$E_i:\mathbb R^{d-1}\to g_i^\perp$ orthogonal,
\begin{equation}\label{eq:k-block-time-wedge}
 \operatorname{Ran}E=(\operatorname{Ran}\mathcal G)^\perp,
 \qquad
 \sqrt{\det(\mathcal G^*\mathcal G)}=\prod_{i=1}^k|g_i|.
\end{equation}
\end{lemma}

\begin{proof}
Choose an orthogonal output basis whose first $k(d-1)$ vectors span
$(\operatorname{Ran}T)^\perp$ and whose last $k$ vectors span
$\operatorname{Ran}T$.  Every square matrix $[A_I\ T]$ is then block
triangular.  Its last diagonal block has Gram determinant
$\det(T^*T)$; squaring and summing over $I$ gives
\eqref{eq:k-time-wedge-factorization} by Cauchy--Binet.  The two identities
in \eqref{eq:k-block-time-wedge} follow block by block.
\end{proof}

\begin{theorem}[Birth-flag multi-landing frame production]
\label{thm:birth-flag-frame-production}
Fix $d\ge2$, a packet cap $P_0$, and $k\ge1$.  Let
$H=H_U\cup H_D$ be a regular connected full hard-sphere packet satisfying
the two-sublayer assumptions \textup{(a)--(d)} of
Theorem~\ref{thm:two-landing}, with at most $P_0$ physical lines and events.
Suppose that a $k$-birth flag $(d_i,c_i)_{i=1}^k$ has been selected, at most
one line at each birth atom.  Use the reverse-chronological Kruskal upper
collision tree $T_K$ of Lemma~\ref{lem:birth-reverse-Kruskal-tree}, substitute every
complementary contact by its incoming time--normal root, and form
\eqref{eq:k-birth-Schur-derivatives}.  Then
\begin{equation}\label{eq:k-birth-frame-production}
 \Jac_{k(d-1)\mid k}(\mathcal S_k,\mathcal T_k)
 \ge c(d,P_0,k)\prod_{i=1}^k
       \phi_i\Delta_i^{d-1}.
\end{equation}
The estimate is uniform in all regular normals and all within-layer time
gaps.  The selected contacts need not be consecutive or disjoint, later
companions may depend on earlier selected states, and no factorization of a
joint complement kernel is assumed.
\end{theorem}

\begin{proof}
Lemma~\ref{lem:lower-influence-cone} gives
$D_{(X_{c_i},V_{c_i})}Y_i=0$ immediately before $d_i$ and shows that later
selected residuals may depend only on earlier selected states.  Thus the
landing derivative is causal, but we do not discard its off-diagonal
blocks.  Apply the determinant-one current clearing of
Lemma~\ref{lem:appC-causal-birth-current-clearing}.  It simultaneously
turns the landing-time columns into
$\mathcal G\tau=(g_i\tau_i)_i$ and makes the actual current on each selected
birth segment equal to its own cleared multiplier.

All upper-tree, incoming-root, companion-history and unselected lower
fields now form the single remainder network of
Theorem~\ref{thm:appC-birth-multiport-schur-jacobi}.  Its Moore--Penrose
compliance is used only on the physical current range; no trace-completeness
or artificial grounding of a disconnected component is assumed.  The
theorem gives
\[
 M_{\Gamma,\perp}\succeq
 \operatorname{diag}_{i}(\Delta_iI_{d-1}),
 \qquad
 \det M_{\Gamma,\perp}\ge\prod_i\Delta_i^{d-1}.
\]
 The same common square physical/KKT system supplies the $k$ time columns.
Lemma~\ref{lem:k-time-wedge-factorization} gives their block-diagonal volume
$\prod_i|g_i|$.  The determinant-one landing
row operation, exact dual forest gauges, and Cauchy--Binet in the original
Cartesian velocity columns therefore give
\[
 \Jac_{k(d-1)\mid k}(\mathcal S_k,\mathcal T_k)
 \ge c(d,P_0,k)\prod_i|g_i|\Delta_i^{d-1}.
\]
Since $|g_i|\ge\phi_i$, this is
\eqref{eq:k-birth-frame-production}.  The argument uses the diagonal
$D_\Delta$ and never compares two time gaps.
\end{proof}

\begin{corollary}[Automatic arbitrary fixed $k$ production]
\label{cor:automatic-fixed-k-birth-production}
Under the assumptions of Theorem~\ref{thm:birth-flag-frame-production}, if
the packet has $P\ge2k-1$ physical lines, its canonical flag satisfies, for
every nonnegative joint kernel $Q$ and for every
$\xi\notin N_\partial$, where $\rho(N_\partial)=0$, the following bounds.
On recovery branch $j$ write
$\Delta_{i,j}^{\sharp}(b,\zeta)
:=\Delta_i(\Psi^{(k)}_{H,\kappa,j}(\xi,b,\zeta),b,\zeta)$ on a present
branch, set $\Delta_{i,j}^{\sharp}=1$ on an absent branch, and define
\begin{equation}\label{eq:fixed-k-Qsharp-definition}
 Q_{H,\kappa,j}^{\sharp,k,\xi}(b,\zeta)
 :=Q(\Psi^{(k)}_{H,\kappa,j}(\xi,b,\zeta),b,\zeta)
   \1_{\mathscr R_{H,\kappa}}
   (\xi,b,\zeta,\Psi^{(k)}_{H,\kappa,j}(\xi,b,\zeta))
   W_{H,\kappa,j}^{\rm rem}(\xi,b,\zeta),
\end{equation}
extended by zero when the branch is absent.  The integration below is over
the same fixed $(b,\zeta)$ domain as in
\eqref{eq:abstract-Qsharp-definition}.
\begin{align}
 J_{H,\kappa}^{\xi}(Q)
 &\le C(d,P_0,k)\cL_\eps^{C|H|}
 \eps^{(d-1)(k-1)}
 \sum_{j=1}^{M_k(d,P_0)}\int
 \left(\prod_{i=1}^k
       (\Delta_{i,j}^{\sharp})^{-(d-1)}\right)
 Q_{H,\kappa,j}^{\sharp,k,\xi}\dd b\dd\nu(\zeta),
 \label{eq:automatic-fixed-k-birth-operator}\\
 J_{H,\kappa}^{\xi}(Q)
 &\le C(d,P_0,k)\cL_\eps^{C|H|}
 \eps^{(d-1)(k-1)}\epstar^{-k(d-1)}
 \sum_{j=1}^{M_k(d,P_0)}\int
 Q_{H,\kappa,j}^{\sharp,k,\xi}\dd b\dd\nu(\zeta).
 \label{eq:automatic-fixed-k-birth-gap}
\end{align}
Here $k$ and $P_0$ are fixed before $\eps$, and the output contains all
regular algebraic recovery branches in boundary-independent slots.
Integrating \eqref{eq:automatic-fixed-k-birth-gap} against
$\dd\rho(\xi)$ and using Tonelli gives the full operator estimate
\begin{equation}\label{eq:automatic-fixed-k-birth-integrated}
 J_{H,\kappa}(Q)
 \le C(d,P_0,k)\cL_\eps^{C|H|}
 \eps^{(d-1)(k-1)}\epstar^{-k(d-1)}
 \sum_{j=1}^{M_k(d,P_0)}\int
 Q_{H,\kappa,j}^{\sharp,k,\xi}
 \dd b\dd\nu(\zeta)\dd\rho(\xi).
\end{equation}
where $J_{H,\kappa}(Q):=\int J_{H,\kappa}^{\xi}(Q)\dd\rho(\xi)$.
With the shorthand
\begin{equation}\label{eq:fixed-k-aggregate-Qsharp}
 \int Q_{H,\kappa}^{\sharp,k}
 :=\sum_{j=1}^{M_k(d,P_0)}\int
 Q_{H,\kappa,j}^{\sharp,k,\xi}
 \dd b\dd\nu(\zeta)\dd\rho(\xi),
\end{equation}
the full operator estimate has the compact form
\begin{equation}\label{eq:automatic-fixed-k-birth-compact}
 J_{H,\kappa}(Q)
 \le C(d,P_0,k)\cL_\eps^{C|H|}
 \eps^{(d-1)(k-1)}\epstar^{-k(d-1)}
 \int Q_{H,\kappa}^{\sharp,k}.
\end{equation}
For every fixed $k\ge2$ the displayed uniform prefactor is a genuine gain:
if $L=|\log\eps|$, then
\[
 \eps^{(d-1)(k-1)}\epstar^{-k(d-1)}
 =\exp\{- (d-1)((k-1)L-k\sqrt L)\}=o(1).
\]
\end{corollary}

\begin{proof}
Lemma~\ref{lem:lower-birth-count} supplies the canonical flag and
Theorem~\ref{thm:birth-flag-frame-production} verifies the flux-completed
frame hypothesis of
Corollary~\ref{cor:fixed-word-flux-multilanding-transfer} with
$\lambda_{\rm fr}=c(d,P_0,k)\prod_i\Delta_i^{d-1}$.  The $P-1$ upper-tree
spheres,
$k$ landing spheres, and packet normalization give exactly
\[
 \eps^{(d-1)(P-1)}\eps^{k(d-1)}
 \eps^{-(d-1)P}=\eps^{(d-1)(k-1)}.
\]
Every unmarked contact cancels its own incoming-root determinant, and every
selected flux $\phi_i$ cancels only its own factor in
\eqref{eq:k-birth-frame-production}.  All changes of variables act on the
same nonnegative integrand.

The reverse-Kruskal upper spanning tree and its uniform translation block
verify \textup{(T1)}, while the original incoming root charts verify
\textup{(T2)}.  The word, lift, sign, normal-frame and maximal-minor cells
are the existing disjoint half-open semialgebraic atlas, so their retained
domains and tie ownership verify \textup{(T5)}.  The canonical birth flag is
one deterministic rule on each word cell and introduces no extra cell sum.

 Corollary~\ref{cor:appD-fixed-k-birth-multiplicity} gives the explicit
degree, variable-count and boundary-independent Borel branch bounds for this
same recovery system.  The radial frames are local proof frames and add no
polynomial variable.  The arbitrary joint kernel is substituted branchwise
and is never marginalized.

Finally the separator construction gives
$\Delta_i\ge h=\epstar/4$.  Absorbing $4^{k(d-1)}$ into the fixed constant
proves \eqref{eq:automatic-fixed-k-birth-gap}.
\end{proof}

For $k=1$ the normalization power is neutral and the bound is
$\epstar^{-(d-1)}$ times the one-port output; the Schur identity is the
one-port formula in
\eqref{eq:appC-physical-network-final-identity}.  For $k=2$,
\eqref{eq:automatic-fixed-k-birth-gap} is
$\eps^{d-1}\epstar^{-2(d-1)}$, exactly the sharp scale in
\eqref{eq:two-landing-actual}.  The zero-set formulation
\eqref{eq:radial-row-operation} verifies that the old radial two-landing
chart and the new two-port clearing use the same physical contact fibre.
Thus the two-port result is the first special case with an operator gain.
The next case is worth recording explicitly.

\begin{corollary}[Three-birth packet gain]
\label{cor:three-birth-packet-gain}
If $P\ge5$, then, for every $\xi\notin N_\partial$, the canonical
three-birth flag satisfies
\begin{equation}\label{eq:three-birth-packet-gain}
 J_{H,\kappa}^{\xi}(Q)
 \le C(d,P_0)\cL_\eps^{C|H|}
 \eps^{2(d-1)}\epstar^{-3(d-1)}
 \sum_{j=1}^{M_3(d,P_0)}\int
 Q_{H,\kappa,j}^{\sharp,3,\xi}\dd b\dd\nu(\zeta).
\end{equation}
Here $M_3(d,P_0):=M_k(d,P_0)|_{k=3}$.
Its integrated full-operator version is the case $k=3$ of
\eqref{eq:automatic-fixed-k-birth-integrated}.
\end{corollary}

\begin{proof}
Take $k=3$ in
Corollary~\ref{cor:automatic-fixed-k-birth-production}.
\end{proof}

Theorem~\ref{thm:birth-flag-frame-production} is used in the paid-cumulant
restart through the adaptive bounded choice \eqref{eq:adaptive-fixed-k}; the
two-landing estimate remains sufficient in the smallest packet and supplies
an independent coarse consistency check.  No value of \(k\) grows with \(\eps\).

\subsubsection{Extensions reduced to new geometry}

Theorem~\ref{thm:fixed-word-multilanding-transfer} separates proved
transfer from open geometric input.
\begin{enumerate}[label=\textup{(\roman*)}]
\item Any construction of $k$ landing residuals with a uniform frame
immediately yields the capacity predicted by the tree/landing normalization.
Theorem~\ref{thm:birth-flag-frame-production} solves the birth-flag class for
every fixed $k$, and
Corollary~\ref{cor:automatic-fixed-k-birth-production} produces such a flag
whenever $P\ge2k-1$.  Exhaustive production beyond birth flags or for other
collision geometries remains open.
\item The compliance proof uses orthogonal free transport, nonnegative shape
operators, finite-dimensional Schur reduction, and the rooted coercivity
hypothesis encoded in the physical network theorem.  Consequently it is a
candidate module for flat Bravais tori and fixed-event semi-dispersing
scatterer words, but only after their reflection law, port maps, sign
convention and grounded coercivity have been verified.  A complete kinetic
limit would additionally require the appropriate collision kernel and source
expansion.
\item If the collision and transport laws are algebraic of uniformly bounded
degree, the finite-fibre part survives with modified $D_0,N_0$.  For smooth
nonalgebraic scatterers, the local coarea and compliance modules remain, but
a different global multiplicity theorem is required.
\item The method is compatible with other positive cluster expansions: the
packet estimate accepts an arbitrary nonnegative complement kernel.  Thus a
future application may replace the Source-B consumer without reopening the
landing geometry, provided its normalization and time-owner ledgers match
the hypotheses of the transfer theorem.
\end{enumerate}

The reusable novelty is consequently the three-module chain
\[
 \begin{gathered}
 \text{landing-frame geometry}
 \quad\Longrightarrow\quad
 \text{positive compliance/no focusing}
 \\
 \Longrightarrow\quad
 \text{finite-fibre arbitrary-kernel coarea}.
 \end{gathered}
\]
The first module is explicitly dimension- and model-specific.  Positive
compliance is portable only under the geometric hypotheses just listed, and
finite-fibre arbitrary-kernel coarea is portable under the bounded-degree,
physical-branch completeness and semialgebraic/Borel hypotheses (T4)--(T5).
The theorem
makes these conditional interfaces, rather than an unconditional extension
claim, the reusable content.
\section{Positive Jacobi Networks and Physical Schur Identities}
\label{app:positive-network}

This section proves the finite-dimensional positive-network assertion used
in Lemma~\ref{lem:grounded-positive-network-lift}.  The argument separates
three facts that a diagrammatic presentation can obscure:
positivity of the grounded Dirichlet form, positivity of its Thomson
current-to-potential compliance, and the one-speed exterior lower bound for
the physical Schur landing row.  Kron reduction proves the first, the dual
current principle proves the second without tracing a historical chord
endpoint through an instantaneous phase pair, and the discrete Green
identity plus square covariant gauges proves the third.  Together they
justify the free-chord elimination in the overlapping-landing chart.
Our terminology follows the classical network elimination of Kron
\cite{Kron1939}; a modern graph-theoretic treatment, including Schur
complements and effective boundary response, is \cite{DorflerBullo2013}.

Throughout this section all vector spaces are finite-dimensional real
Hilbert spaces.  Adjoints and orthogonality refer to the kinetic-energy
inner product.  In the hard-sphere application a chronological port carries
the complete reduced configuration space
\(E=T_X\widehat{\mathcal Q}_{\mathcal P}\), of dimension \(d(P-1)\).
The $d$-dimensional object is instead the landing-current frame
\(Z_1:\mathbb R^d\to E\).  This distinction is essential: a binary
collision mixes the two incident particle lines and is not an independent
orthogonal gauge on either $d$-dimensional line port.

\subsection{Exterior volume and the positive conormal cone}

If $F:\R^k\to E$ is a column frame, write
\begin{equation}\label{eq:appC-frame-volume}
 \operatorname{Vol}_k(F)
 :=\sqrt{\det(F^*F)}=\|\bigwedge^kF\|.
\end{equation}
This quantity is zero exactly when the columns are dependent.  We shall use
the following elementary consequences of the min--max principle.

\begin{lemma}[Exterior monotonicity]
\label{lem:appC-exterior-monotonicity}
Let $A,B:\R^k\to E$.
\begin{enumerate}[label=\textup{(\roman*)}]
\item If $A^*A\ge B^*B\ge0$, then
$\operatorname{Vol}_k(A)\ge\operatorname{Vol}_k(B)$.
\item If $S=S^*\ge bI_E$ with $b>0$, then
\[
 \operatorname{Vol}_k(SF)\ge b^k\operatorname{Vol}_k(F).
\]
\item If $P:E\to E$ is invertible and
$\|P^{-1}\|\le c^{-1}$, then
$\operatorname{Vol}_k(PF)\ge c^k\operatorname{Vol}_k(F)$.
\end{enumerate}
\end{lemma}

\begin{proof}
Part~\textup{(i)} follows by multiplying the ordered eigenvalues of the two
positive $k$-by-$k$ Gram matrices.  For part~\textup{(ii)}, the least
singular value of $S$ is at least $b$, so the least singular value of its
restriction to every $k$-plane is at least $b$.  Part~\textup{(iii)} is the
same argument with the least singular value of $P$.
\end{proof}

A conormal frame is a pair $(Z,W)$ of maps $\R^k\to E$.  Define its
symmetric mixed Gram matrix by
\begin{equation}\label{eq:appC-positive-cone}
 \mathcal G(Z,W):=\frac12(Z^*W+W^*Z).
\end{equation}
We call $(Z,W)$ \emph{positive} if $\mathcal G(Z,W)\ge0$.  This
does not require $Z^*W$ itself to be symmetric.

For $r\ge0$, an orthogonal $O:E\to E$, and $H=H^*\ge0$, introduce the two
backward adjoint-pullback elementary maps
\begin{align}
 \mathsf F_r(Z,W)&=(Z,W+rZ),
 \label{eq:appC-free-generator}\\
 \mathsf C_{O,H}(Z,W)&=(OZ+HO W,OW).
 \label{eq:appC-collision-generator}
\end{align}
The first is backward free transport.  The second is a fixed-normal
orthogonal collision followed by its nonnegative cylinder-curvature shear.
The order of $H$ and $O$ in
\eqref{eq:appC-collision-generator} is a convention: replacing $H$ by an
orthogonal conjugate gives the other common convention.

\begin{lemma}[Local cone and exterior invariance]
\label{lem:appC-local-positive-generators}
Both maps \eqref{eq:appC-free-generator}--
\eqref{eq:appC-collision-generator} preserve the cone
$\mathcal G\ge0$.  On this cone neither map decreases
$\operatorname{Vol}_k(W)$.
\end{lemma}

\begin{proof}
For free transport,
\begin{align*}
 \mathcal G(Z,W+rZ)&=\mathcal G(Z,W)+rZ^*Z,\\
 (W+rZ)^*(W+rZ)-W^*W
 &=2r\mathcal G(Z,W)+r^2Z^*Z\ge0.
\end{align*}
For a collision, orthogonality gives
\[
 \mathcal G(OZ+HO W,OW)
 =\mathcal G(Z,W)+W^*O^*HO W\ge0,
\]
whereas the new velocity frame is $OW$ and has the same Gram matrix as
$W$.  Lemma~\ref{lem:appC-exterior-monotonicity} completes the proof.
\end{proof}

Thus every word in the generators is exterior-noncontracting on positive
frames.  What remains is to show that rooted constraint elimination and a
positive free chord produce precisely such a word, rather than an arbitrary
Schur complement.

\subsection{The block Riccati invariant}

Write a linear conormal relation on $E\oplus E$ in block form
\begin{equation}\label{eq:appC-block-relation}
 \binom{Z_{\rm out}}{W_{\rm out}}
 =\mathsf T\binom{Z}{W},
 \qquad
 \mathsf T=\begin{pmatrix}A&B\\ C&D\end{pmatrix}.
\end{equation}
The velocity row is $CZ+DW$.  Whenever $D$ is invertible, set
\begin{equation}\label{eq:appC-Riccati-block}
 K(\mathsf T):=D^{-1}C.
\end{equation}
Then
\begin{equation}\label{eq:appC-velocity-factorization}
 W_{\rm out}=D\bigl(W+K(\mathsf T)Z\bigr).
\end{equation}

The generators have matrices
\begin{equation}\label{eq:appC-generator-matrices}
 \mathsf F_r=\begin{pmatrix}I&0\\ rI&I\end{pmatrix},
 \qquad
 \mathsf C_{O,H}=\begin{pmatrix}O&HO\\0&O\end{pmatrix}.
\end{equation}
We use right multiplication because it corresponds to appending the next
local operation on the input side of the adjoint-pullback word.

\begin{lemma}[Riccati updates]
\label{lem:appC-Riccati-updates}
Suppose $D$ is invertible and $K=D^{-1}C=K^*\ge0$.
\begin{enumerate}[label=\textup{(\roman*)}]
\item For $r\ge0$,
\[
 K(\mathsf T\mathsf F_r)=K+rI.
\]
\item For $O$ orthogonal and $H=H^*\ge0$,
\begin{equation}\label{eq:appC-collision-Riccati}
 K(\mathsf T\mathsf C_{O,H})
 =O^*(I+KH)^{-1}KO
 =O^*K^{1/2}(I+K^{1/2}HK^{1/2})^{-1}K^{1/2}O.
\end{equation}
In particular, the new $D$ is invertible and the new Riccati block is
self-adjoint nonnegative.
\end{enumerate}
\end{lemma}

\begin{proof}
Right multiplication by $\mathsf F_r$ replaces $(C,D)$ by
$(C+rD,D)$, proving part~\textup{(i)}.  Right multiplication by the
collision generator replaces the velocity row by
\[
 (CO,(CH+D)O)=\bigl(DKO,D(I+KH)O\bigr).
\]
The nonzero spectra of $KH$ and $K^{1/2}HK^{1/2}$ agree away from zero and
are nonnegative.  Hence $I+KH$ is invertible.  Direct substitution gives
the first expression in \eqref{eq:appC-collision-Riccati}.  The identity
\[
 (I+KH)^{-1}K
 =K^{1/2}(I+K^{1/2}HK^{1/2})^{-1}K^{1/2}
\]
first holds when $K$ is invertible and then follows in general by replacing
$K$ by $K+\delta I$ and taking $\delta\downarrow0$.  This proves symmetry
and nonnegativity.
\end{proof}

\begin{proposition}[Positive-word factorization]
\label{prop:appC-positive-word}
Let $\mathsf T$ be any finite word in the generators
\eqref{eq:appC-generator-matrices}.  Then its block $D$ is invertible and
\begin{equation}\label{eq:appC-positive-word-properties}
 K(\mathsf T)=D^{-1}C=K(\mathsf T)^*\ge0,
 \qquad
 \operatorname{Vol}_k(DF)\ge\operatorname{Vol}_k(F)
\end{equation}
for every frame $F$.
\end{proposition}

\begin{proof}
Start with the identity word, for which $D=I$ and $K=0$.  The first
assertion follows by induction from
Lemma~\ref{lem:appC-Riccati-updates}.

For the exterior assertion, apply the word to the frame $(Z,W)=(0,F)$.
It is positive because its mixed Gram is zero.  Its final velocity frame is
$DF$.  Lemma~\ref{lem:appC-local-positive-generators}, iterated through the
word, shows that its exterior volume cannot decrease.
\end{proof}

The two conclusions in
\eqref{eq:appC-positive-word-properties} are logically distinct.  The
Riccati calculation proves that the position contribution cannot cancel a
one-speed input, while the conormal calculation proves that the prefactor
$D$ cannot destroy exterior volume.

\begin{corollary}[One-speed lower bound]
\label{cor:appC-one-speed-word}
Under the hypotheses of Proposition~\ref{prop:appC-positive-word}, for
$b>0$ one has
\begin{equation}\label{eq:appC-one-speed-word}
 CZ+DbZ=D(bI+K)Z,
 \qquad
 \operatorname{Vol}_k\bigl(D(bI+K)Z\bigr)
 \ge b^k\operatorname{Vol}_k(Z).
\end{equation}
\end{corollary}

\begin{proof}
The factorization is \eqref{eq:appC-velocity-factorization}.  Since
$K=K^*\ge0$, the least singular value of $bI+K$ is at least $b$.
Apply first Proposition~\ref{prop:appC-positive-word} and then
Lemma~\ref{lem:appC-exterior-monotonicity}.
\end{proof}

\subsection{Covariant graph Dirichlet forms}

We now separate two boundary operators which have reciprocal physical
dimensions.  The Dirichlet-to-Neumann operator is a \emph{stiffness}; the
operator occurring in the velocity row of the backward pullback relation is
its Neumann-to-Dirichlet \emph{compliance}.  Keeping these names distinct is
essential already for a single free edge.

Let $\mathcal V$ be a finite set of ports.  Each port carries a copy $E_v$
of a fixed Euclidean space $E_0$.  An oriented free segment $e=(u,v)$ has
length $\ell_e>0$ and an orthogonal identification
$O_e:E_u\to E_v$.  A port may also carry a shunt
$H_v=H_v^*\ge0$.  We also allow a physical-line factor edge
$\gamma$ with a finite-dimensional current space $G_\gamma$ and a linear
drop map
\[
 A_\gamma:\bigoplus_{v\in\mathcal V}E_v\longrightarrow G_\gamma.
\]
In the application $G_\gamma=\mathbb R^d$ and $A_\gamma$ is built from the
exact particle selector $S_i$ on the zero-centre-of-mass representative.
It is distinct from the normalized pair difference $B_e$ in
\eqref{eq:appC-normalized-pair-incidence}.  For a field
$\xi=(\xi_v)_{v\in\mathcal V}$ define
\begin{equation}\label{eq:appC-network-form}
 \mathcal I(\xi)
 =\frac12\sum_{e=(u,v)}\ell_e^{-1}
    |\xi_v-O_e\xi_u|^2
  +\frac12\sum_\gamma\ell_\gamma^{-1}|A_\gamma\xi|^2
  +\frac12\sum_{v\in\mathcal V}\langle H_v\xi_v,\xi_v\rangle.
\end{equation}
Parallel segments and loops after orthogonal identification are allowed.
The second sum is required for a chord carried by one physical line.  It is
rank $d$ in the full reduced space and must not be replaced by a fictitious
full-rank $E_0$ edge.

Choose a boundary set $B\subset\mathcal V$ and an internal set
$I=\mathcal V\setminus B$.  We ground one port in every connected component
which otherwise has the constant covariant zero mode.  Equivalently, delete
the corresponding field coordinates.  Let $L$ be the symmetric matrix of
$\mathcal I$ on the remaining coordinates.

\begin{lemma}[Grounded invertibility]
\label{lem:appC-grounded-invertibility}
Suppose every internal connected component meets either a boundary port, a
ground, or a shunt which is positive on its residual covariant constant
mode, and suppose the stacked full-edge and physical factor-incidence rows
have no additional kernel after those grounds are imposed.  Then $L_{II}$
is positive definite.
\end{lemma}

\begin{proof}
If $y^*L_{II}y=0$, every nonnegative summand in
\eqref{eq:appC-network-form} vanishes.  Hence $y_v=O_e y_u$ along every
full segment, $A_\gamma y=0$ at every factor edge,
$H_v^{1/2}y_v=0$ at every shunt, and the field is zero at
every boundary or ground adjacent to the internal component.  Transporting
these identities through the connected component gives $y=0$ under the
stated hypothesis.
\end{proof}

\begin{lemma}[Kron positivity and associativity]
\label{lem:appC-Kron-positivity}
Under the hypotheses of Lemma~\ref{lem:appC-grounded-invertibility}, the
 effective boundary stiffness
\begin{equation}\label{eq:appC-Kron-response}
 S_B:=L_{BB}-L_{BI}L_{II}^{-1}L_{IB}
\end{equation}
is self-adjoint nonnegative.  It is characterized by
\begin{equation}\label{eq:appC-Dirichlet-principle}
 \frac12\langle S_Bz,z\rangle
 =\min_y\mathcal I(y,z).
\end{equation}
If the internal variables are divided into any finite number of groups,
successive Schur complementation in any admissible order gives the same
$S_B$.
\end{lemma}

\begin{proof}
Complete the square:
\begin{align*}
 \mathcal I(y,z)
 &=\frac12\langle L_{II}(y+L_{II}^{-1}L_{IB}z),
             y+L_{II}^{-1}L_{IB}z\rangle\\
 &\quad+\frac12\langle
   (L_{BB}-L_{BI}L_{II}^{-1}L_{IB})z,z\rangle.
\end{align*}
This proves \eqref{eq:appC-Dirichlet-principle} and positivity.  The
associativity assertion follows either by repeating the completed-square
identity or by observing that every elimination order computes the same
minimum over all internal fields.
\end{proof}

Orthogonal edge identifications cause no sign problem.  On a spanning
forest they can be removed by choosing a root frame and transporting it
along the unique forest paths.  A co-tree chord then retains its orthogonal
holonomy, but its term in \eqref{eq:appC-network-form} remains a square.

\begin{lemma}[Gauge invariance]
\label{lem:appC-gauge-invariance}
Let $U_v:E_v\to E_0$ be arbitrary orthogonal port gauges and replace
$\xi_v$ by $U_v^*\eta_v$.  Then the network form retains the form
\eqref{eq:appC-network-form}, with
\[
 O_e\longmapsto U_vO_eU_u^*,
 \qquad H_v\longmapsto U_vH_vU_v^*.
\]
The Kron stiffness is changed only by orthogonal boundary conjugation; its
grounded compliance transforms by the same conjugation.
Consequently its positivity and every exterior-volume inequality are gauge
independent.
\end{lemma}

\begin{proof}
Substitution gives the displayed transformations term by term.
Schur complementation commutes with block-diagonal orthogonal conjugation,
and exterior volume is orthogonally invariant.
\end{proof}

\subsection{Boundary stiffness and compliance}

For a quadratic form $\frac12\langle L\xi,\xi\rangle$, define the boundary
current by $j_B=(L\xi)_B$.  At the minimizing internal field,
$(L\xi)_I=0$, so the Dirichlet-to-Neumann law is
\begin{equation}\label{eq:appC-boundary-current}
 j_B=S_Bz.
\end{equation}
This equation maps potential to current.  It does \emph{not} identify the
backward conormal shear with $S_B$.

\begin{lemma}[One-port stiffness--compliance duality]
\label{lem:appC-stiffness-compliance-duality}
Let $S_B$ be \eqref{eq:appC-Kron-response}.  On the grounded boundary range
define
\begin{equation}\label{eq:appC-compliance-definition}
 K_B:=S_B^+,
\end{equation}
where $S_B^+$ is the Moore--Penrose inverse.  Then $K_B=K_B^*\ge0$,
$z=K_Bj$ is the unique grounded boundary potential generated by
$j\in\operatorname{Ran}S_B$, and
\begin{equation}\label{eq:appC-Thomson-duality-one-port}
 \frac12\langle j,K_Bj\rangle
 =\min_{S_Bz=j}\frac12\langle S_Bz,z\rangle.
\end{equation}
If the retained one-port block is positive definite, $K_B=S_B^{-1}$.
Translation zero modes may equivalently be quotiented before inversion.
\end{lemma}

\begin{proof}
The spectral theorem gives all assertions on
$\operatorname{Ran}S_B=(\ker S_B)^\perp$.  In an eigenbasis, a positive
stiffness eigenvalue $s$ contributes potential $j/s$ and energy
$j^2/(2s)$; zero modes are removed by grounding or annihilated by the
range condition.  This is \eqref{eq:appC-Thomson-duality-one-port}.
\end{proof}

\begin{lemma}[Three exact one-port benchmark identities]
\label{lem:appC-three-benchmark-models}
In a common orthonormal port frame the following identities hold.
\begin{enumerate}[label=\textup{(\roman*)}]
\item A single edge of length $\ell$ to ground has
$S=\ell^{-1}I$ and $K=\ell I$.
\item Two series edges of lengths $\ell_1,\ell_2$ have
$S=(\ell_1+\ell_2)^{-1}I$ and $K=(\ell_1+\ell_2)I$.
\item An edge of length $\ell$ in parallel with a shunt $H\ge0$ has
\begin{equation}\label{eq:appC-edge-shunt-benchmark}
 S=\ell^{-1}I+H,
 \qquad K=(\ell^{-1}I+H)^{-1}
 =(I+\ell H)^{-1}\ell I.
\end{equation}
The last formula is exactly the compliance Riccati update
$K\mapsto(I+KH)^{-1}K$ after the free update $0\mapsto\ell I$.
\end{enumerate}
\end{lemma}

\begin{proof}
The first identity is immediate from \eqref{eq:appC-network-form}.  For two
series edges, minimizing
$|y|^2/(2\ell_1)+|z-y|^2/(2\ell_2)$ gives
$y=\ell_1(\ell_1+\ell_2)^{-1}z$ and the stated effective energy.  A shunt
adds $H$ to the Dirichlet stiffness, proving (iii); inversion gives the
displayed compliance and agrees with the Riccati formula by direct
multiplication.  These models fix the length, sign, and series/shunt
conventions used below.
\end{proof}

\begin{definition}[Trace-complete rooted word]
\label{def:appC-trace-complete-word}
For a particle set $C$ put
\[
 E_C:=\left\{(x_i)_{i\in C}\in(\mathbb R^d)^C:
                  \sum_{i\in C}x_i=0\right\},
 \qquad S_i x=x_i.
\]
A grounded chronological hard-sphere word is called \emph{trace-complete}
if, after the complete pair junctions have been substituted, its endpoint
fields admit a raw rooted ordering with the following properties.
\begin{enumerate}[label=\textup{(TC\arabic*)}]
\item Every initial forest component $C$ has a declared ground particle
  $g_C$, and the corresponding block of $\mathbb Gq$ is $S_{g_C}q$.
\item At a chronological forest merge $e$, the new relative-cluster
  translation $h_e$ occurs in the physical row in the form
  \begin{equation}\label{eq:appC-trace-merge-row}
   u_e=T_eh_e+V_eq_{<e},
   \qquad
   s_{\min}(T_e)\ge c(P_0),\quad
   \|T_e\|\le C(P_0),
  \end{equation}
  where $q_{<e}$ contains only coordinates already preceding $h_e$ in the
  rooted ordering.  In particular, an unresolved retained boundary field is
  not placed in $q_{<e}$.
\item The frontier requirement has the following two separately verifiable
  parts.
  \begin{enumerate}[label=\textup{(TC3\alph*)},leftmargin=*]
  \item Whenever a new reduced frontier $x_\alpha\in E_{C_\alpha}$ is
    reached from an already reached frontier, the ordinary physical-line
    factor rows contain
    \begin{equation}\label{eq:appC-trace-selector-row}
     d_\alpha=\Sigma_{I_\alpha}
       (x_\alpha-\mathscr R_\alpha x_{\operatorname{par}(\alpha)}),
     \qquad
     \Sigma_{I_\alpha}x=(S_ix)_{i\in I_\alpha}.
    \end{equation}
    If the slab is crossed by the restored chord, then
    $I_\alpha=C_\alpha\setminus\{c_\alpha\}$ with $c_\alpha$ its physical
    chord line; otherwise any one selector may be omitted.  The restored
    chord row itself is not used as a full-rank frontier row.
  \item The transport $\mathscr R_\alpha$ in \textup{(TC3a)} is a product
    only of complete-pair orthogonal junctions on the reached frontier.
    Hence its reduced-frontier extension is orthogonal and the inverse
    selector reconstruction has norm bounded by $C(P_0)$, with no inverse
    projected collision block.
  \end{enumerate}
\item The retained trace is a fixed trace of a reached frontier and satisfies
  $\|\operatorname{tr}_B\|\le C(P_0)$.
\end{enumerate}
Every selector, forest-incidence and path matrix in this definition is the
corresponding scalar graph matrix tensored with $I_d$.  Hence its rank is the
scalar rank multiplied by $d$, and its nonzero singular values are repeated
$d$ times.  The definition is a verifiable coverage certificate, not a
synonym for graph connectivity.
\end{definition}

\begin{lemma}[Uniform trace left inverse]
\label{lem:appC-uniform-trace-left-inverse}
Let $W$ be trace-complete, with at most $P_0$ physical lines and events and
at most one restored chord.  There is a linear map $\mathcal L_W$, depending
on the fixed word but satisfying $\|\mathcal L_W\|\le C(P_0)$, such that
\begin{equation}\label{eq:appC-explicit-trace-left-inverse}
 \operatorname{tr}_B
 =\mathcal L_W
   \begin{pmatrix}\mathbb A_W\\ \mathbb U_F\\ \mathbb G\end{pmatrix}.
\end{equation}
Consequently, for $z=\operatorname{tr}_Bq$,
\begin{equation}\label{eq:appC-certified-boundary-trace}
 \|z\|\le C(P_0)\bigl(
   \|\mathbb A_Wq\|+\|\mathbb U_Fq\|+\|\mathbb Gq\|\bigr).
\end{equation}
\end{lemma}

\begin{proof}
We construct the left inverse.  First fix a forest component $C$, a ground
$g\in C$, and an oriented rooted spanning tree $F_C$.  The raw difference
map
\[
 \mathcal D_{F_C,g}x
 :=\bigl(S_gx,(\sqrt2B_ex)_{e\in F_C}\bigr)
\]
has the explicit path inverse
\begin{equation}\label{eq:appC-rooted-path-inverse}
 x_i=x_g+\sum_{e\in\operatorname{path}(g,i)}
             \sigma_{i,e}\sqrt2B_ex,
 \qquad \sigma_{i,e}\in\{-1,1\}.
\end{equation}
Thus $\|\mathcal D_{F_C,g}^{-1}\|\le C(|C|)$.  Formula
\eqref{eq:appC-trace-merge-row} is the chronological physical version of
this path reconstruction: its new diagonal is $T_e$, while every causal
root contribution lies in an earlier column.  Forward substitution therefore
recovers every rooted-cluster coordinate from $(\mathbb Gq,\mathbb U_Fq)$
with norm at most $C(P_0)$.

Next let $c\in C$.  The selector stack
\[
 \Sigma_{\widehat c}:E_C\longrightarrow
          (\mathbb R^d)^{C\setminus\{c\}},
 \qquad \Sigma_{\widehat c}x=(x_i)_{i\ne c},
\]
is square.  Its inverse is the explicit map
\begin{equation}\label{eq:appC-missing-selector-inverse}
 (\Lambda_{\widehat c}y)_i=y_i\ (i\ne c),
 \qquad
 (\Lambda_{\widehat c}y)_c=-\sum_{i\ne c}y_i,
 \qquad
 \|\Lambda_{\widehat c}\|\le\sqrt{|C|}.
\end{equation}
Because $\mathscr R_\alpha$ is orthogonal, the inverse of
$\Sigma_{\widehat c}\mathscr R_\alpha$ is
$\mathscr R_\alpha^*\Lambda_{\widehat c}$ with the same bound.  Hence
\eqref{eq:appC-trace-selector-row} gives the exact recurrence
\begin{equation}\label{eq:appC-frontier-trace-recurrence}
 x_\alpha=\mathscr R_\alpha x_{\operatorname{par}(\alpha)}
       +\Lambda_{\widehat c}d_\alpha .
\end{equation}
In a chord slab the missing line is precisely the chord particle, so the
other physical-line rows already recover the entire zero-centre-of-mass
frontier.  The rank-$d$ chord row is an additional row and is never
inverted on $E_C$.

Select from $(\mathbb A_W,\mathbb U_F,\mathbb G)$ the ground rows, the merge
rows in rooted order, and at every new frontier the rows in
\eqref{eq:appC-trace-selector-row}.  Against the corresponding raw rooted
coordinate order the resulting square matrix is block lower triangular:
\begin{equation}\label{eq:appC-trace-triangular-matrix}
 \mathcal T_W=
 \begin{pmatrix}
  D_{\rm grd}&0&0&\cdots\\
  *&T_1&0&\cdots\\
  *&*&\Sigma_{I_1}\mathscr R_1&\cdots\\
  *&*&*&\ddots
 \end{pmatrix}.
\end{equation}
Equations \eqref{eq:appC-rooted-path-inverse},
\eqref{eq:appC-missing-selector-inverse}, and
\eqref{eq:appC-trace-merge-row} give
$\|\mathcal T_W^{-1}\|\le C(P_0)$ by forward substitution through at most
$P_0$ blocks.  If $\mathcal P_W$ selects these rows from the complete
residual stack, set
\begin{equation}\label{eq:appC-trace-left-inverse-formula}
 \mathcal L_W:=\operatorname{tr}_B\mathcal T_W^{-1}\mathcal P_W.
\end{equation}
This proves \eqref{eq:appC-explicit-trace-left-inverse} and
\eqref{eq:appC-certified-boundary-trace}.
\end{proof}

\mainstatementreference{Proposition~\ref{prop:appC-positive-network-closure}}

\mainproofheading{Proposition~\ref{prop:appC-positive-network-closure}}{supp:proof:prop:appC-positive-network-closure}
\begin{proof}
Split every free segment at every collision port.  The one-port stiffness
is positive definite, and we record why this is stronger than internal
invertibility.  Let \(\operatorname{tr}_Bq=z\) be the retained boundary
datum and let \(\mathbb Gq=0\) denote the declared ground coordinates.
 Lemma~\ref{lem:appC-uniform-trace-left-inverse} gives the trace estimate
 \begin{equation}\label{eq:appC-boundary-trace-coercivity}
  \|z\|\le C(P_0)\bigl(
        \|\mathbb A_Wq\|+\|\mathbb U_Fq\|+\|\mathbb Gq\|\bigr).
 \end{equation}
 In particular, this estimate is supplied by the square triangular matrix
 \eqref{eq:appC-trace-triangular-matrix}; it is not inferred from the internal
 estimate \eqref{eq:appC-actual-grounded-kernel}, and it never inverts the
 rank-$d$ chord row on the full reduced port.

Now fix \(z\ne0\) and let \(q_z\) be the internal minimizer in the
Dirichlet principle.  If its energy were zero, every factor drop
\(A_e q_z\) and every \(H_v^{1/2}q_{z,v}\) would vanish.  The physical forest
rows and grounds vanish because they define the admissible affine space.
Equation \eqref{eq:appC-boundary-trace-coercivity} would give \(z=0\), a
contradiction.  Consequently
\begin{equation}\label{eq:appC-boundary-stiffness-strict}
 \langle S_{\mathcal N}z,z\rangle
 =2\min_{\operatorname{tr}_Bq=z}\mathcal I(q)>0
 \qquad(z\ne0).
\end{equation}
This proof includes kernel shunts: their zero directions simply contribute
no energy, but they cannot create an additional zero mode because the
factor/forest/ground trace estimate already kills it.  Hence
\(S_{\mathcal N}>0\) and \(K_{\mathcal N}=S_{\mathcal N}^{-1}\) exists.

On the resulting graph prescribe a covariantly transported boundary current
$j$.  Along an edge
$e$ let $J_e$ be its current.  A curvature port $H_v\ge0$ is represented by
an auxiliary current $q_v\in\operatorname{Ran}H_v$ with dual energy
$\langle H_v^+q_v,q_v\rangle$, where $H_v^+$ is the Moore--Penrose inverse;
the component in $\ker H_v$ cannot leave through that port.  Minimize
\begin{equation}\label{eq:appC-Thomson-functional}
 \frac12\sum_e\ell_e|J_e|^2+
 \frac12\sum_v\langle H_v^+q_v,q_v\rangle
\end{equation}
under covariant current conservation at every internal vertex and boundary
divergence $j$.  Grounding removes the common-potential ambiguity.  The
minimum is a nonnegative quadratic form in $j$; polarization therefore
defines a self-adjoint nonnegative operator $K_{\mathcal N}$ by
\[
 \min\eqref{eq:appC-Thomson-functional}
 =\frac12\langle j,K_{\mathcal N}j\rangle.
\]
This is the Thomson dual of the Dirichlet principle
\eqref{eq:appC-Dirichlet-principle}, equivalently
Lemma~\ref{lem:appC-stiffness-compliance-duality} after internal Dirichlet
elimination.  A singular shunt causes no singular boundary compliance:
$H_v^+$ is used only on $\operatorname{Ran}H_v$, while its kernel simply
carries no shunt current.

Varying the minimizing current and summing edgewise
$\langle\Delta\xi,J\rangle$ gives the discrete Green identity: the covariant
boundary potential equals $K_{\mathcal N}j$.  For a raw one-speed landing
current $z$ put $j=Q_-z$.  The collision-free lower interval contributes the
independent raw dual energy $b|z|^2/2$.  Hence the complete raw current action
and its derivative are
\begin{equation}\label{eq:appC-one-speed-Thomson-action}
 \mathcal T_b^{\rm raw}(z)
 =\frac12\left(b|z|^2+\langle Q_-z,K_{\mathcal N}Q_-z\rangle\right)
 =\frac12\langle z,(bI+\widetilde K_{\mathcal N})z\rangle,
 \qquad
 D_z\mathcal T_b^{\rm raw}=(bI+\widetilde K_{\mathcal N})z.
\end{equation}
The chord causes no phase-trace ambiguity here: its current $J_\gamma$
remains an independent variable of \eqref{eq:appC-Thomson-functional} with
the positive term $\ell_\gamma|J_\gamma|^2/2$ until the global current
minimization.  No historical endpoint potential is expressed through an
instantaneous phase pair.

Let $Q_-$ send a raw landing-current frame to the covariant current $j$, and
let $Q_+$ send the dual covariant potential back to the raw retained Schur
row.  We construct them without using the conclusion of the physical
Schur theorem below.  For every normalized pair incidence $B_e$, the
injection $B_e^*:G_e\to E$ is isometric; complete it by Gram--Schmidt to an
orthogonal square basis matrix $[B_e^*,R_e]$ of $E$.  At a collision use the
complete orthogonal pair matrix, never either one-particle projection.  At
a forest merge use the square matrix
\[
 [R_F,R_B]\begin{pmatrix}G_F&0\\0&I_{E_B}\end{pmatrix},
\]
where $R_F^{\rm raw}=R_FG_F$ and $E=\operatorname{Ran}R_F\oplus E_B$.
The entries before QR are normalized incidences and triangular forest rows
with diagonal signs.  For at most $P_0$ lines there are only finitely many
ordered forest types, so this square matrix and its inverse are bounded by
constants depending only on $P_0$.  Composing at most $P_0$ such primal
potential factors defines $\mathsf P_W$.  The current and potential gauges are
defined by the exact dual-pairing identity
\begin{equation}\label{eq:appC-exact-dual-boundary-gauges}
 Q_-:=\mathsf P_W^{-*},\qquad Q_+:=\mathsf P_W^{-1}=Q_-^*,
 \qquad
 \langle Q_-z,\xi\rangle_{\rm cov}
   =\langle z,Q_+\xi\rangle_{\rm raw}.
\end{equation}
All other factors are orthogonal.  Therefore
\begin{equation}\label{eq:appC-positive-word-coordinate-bounds}
 s_{\min}(Q_\pm)\ge c_Q(P_0)>0,
 \qquad \|Q_\pm\|\le C_Q(P_0).
\end{equation}
No edge length, $C_{\rm ch}^{-1}$, or one-particle collision projection
occurs in either gauge.

For later raw-coordinate bookkeeping, the separate square junction fibre
matrix is
\begin{equation}\label{eq:appC-junction-shear-cost}
 J_a=\begin{pmatrix}I_{2d}&0\\ \mathscr R_a^{(2d)}&I_{2d}\end{pmatrix},
 \qquad
 s_{\min}(J_a)=\sqrt{\frac{3-\sqrt5}{2}}=:s_0.
\end{equation}
It is used only when the final physical frame is compared with the raw
cluster determinant; it is not a rectangular lift and is not part of the
chord current elimination.

Take the raw current in \eqref{eq:appC-one-speed-Thomson-action} to be
$z=F\zeta$, so its covariant representative is $Q_-F\zeta$.  The discrete
Green identity identifies the derivative of the raw action with the physical
retained Schur row.  Thus, as an identity of maps on the same landing
coefficient space,
\[
 W_{\rm Sch}=(bI+\widetilde K_{\mathcal N})F.
\]
Since $\widetilde K_{\mathcal N}\ge0$, the square factor has
$s_{\min}(bI+\widetilde K_{\mathcal N})\ge b$.  The exterior-power
singular-value inequality yields
\[
 \operatorname{Vol}_k(W_{\rm Sch})
 \ge b^k\operatorname{Vol}_k(F).
\]
This is \eqref{eq:appC-network-CD-identification}; no inverse chronological
block or endpoint-potential trace is used.  It proves the second assertion
of \eqref{eq:appC-network-closure}.
\end{proof}

If raw cluster translations rather than kinetic-orthonormal coordinates
are used, the incidence coordinate changes at the input and output have
condition numbers bounded in terms of the number $P$ of physical lines.
Consequently \eqref{eq:appC-network-closure} becomes
\begin{equation}\label{eq:appC-coordinate-cost}
 \operatorname{Vol}_k(W_{\rm Sch})
 \ge c(P)^k b^k\operatorname{Vol}_k(F).
\end{equation}
Together with the junction and forest factors in
\eqref{eq:appC-positive-word-coordinate-bounds}, this gives the constant
denoted $c(P_0)$ in Section~6.

\subsection{Exact free-chord restoration}

Let a physical line have no collision in $(u,d)$, $u<d$.  Its separator
state $(X,V)$ and its two endpoint values satisfy
\begin{equation}\label{eq:appC-chord-endpoints}
 X+(u-s)V=Y_u,
 \qquad X+(d-s)V=Y_d.
\end{equation}
The coefficient matrix in $(X,V)$ is
\[
 \begin{pmatrix}I&(u-s)I\\I&(d-s)I\end{pmatrix},
\]
whose determinant has absolute value $(d-u)^{\dim E_0}$.  Eliminating
$(X,V)$ therefore has neither multiplicity nor a hidden singular branch.

At the tangent level the unique interpolant is
\begin{equation}\label{eq:appC-affine-chord-field}
 \xi(t)=\frac{d-t}{d-u}\xi_u+
        \frac{t-u}{d-u}\xi_d.
\end{equation}
Its free index energy is
\begin{equation}\label{eq:appC-affine-chord-energy}
 \int_u^d|\dot\xi(t)|^2\dd t
 =\frac1{d-u}|\xi_d-\xi_u|^2.
\end{equation}
Let $S_u,S_d$ be the exact, unnormalised physical-particle selectors from the complete reduced fields
at the two endpoint ports to this physical line.  After the pair-collision
junctions have been transported orthogonally, the endpoint drop is
\[
 A_\gamma\xi:=S_d\xi_d-O_{u,d}S_u\xi_u\in\mathbb R^d.
\]
Thus the chord contributes, with no rescaling of its physical length,
$(d-u)^{-1}|A_\gamma\xi|^2$, precisely a rank-$d$ factor term in
\eqref{eq:appC-network-form}.

\begin{lemma}[Chord-first equals chord-last]
\label{lem:appC-chord-associativity}
Let a grounded forest network contain the endpoints of the collision-free
line \eqref{eq:appC-chord-endpoints}.  The following two procedures give the
same reduced boundary relation:
\begin{enumerate}[label=\textup{(\roman*)}]
\item eliminate $(X,V)$ first, insert the term
\eqref{eq:appC-affine-chord-energy}, and then eliminate the forest fields;
\item retain the affine chord fields, eliminate the forest fields first,
and eliminate the chord fields last.
\end{enumerate}
The common relation belongs to the positive network class of
Proposition~\ref{prop:appC-positive-network-closure}.
\end{lemma}

\begin{proof}
Order the variables of the \emph{same} complete quadratic form as chord
interior $c$, forest interior $f$, and retained boundary $b$:
\begin{equation}\label{eq:appC-common-chord-forest-matrix}
 \mathbb L=
 \begin{pmatrix}
  L_{cc}&L_{cf}&L_{cb}\\
  L_{fc}&L_{ff}&L_{fb}\\
  L_{bc}&L_{bf}&L_{bb}
 \end{pmatrix}.
\end{equation}
The chord block is positive definite with its endpoint values fixed, and
the relevant forest block is positive definite after grounding.  Block
Gaussian elimination gives the quotient identity
\begin{equation}\label{eq:appC-Schur-quotient-identity}
 (\mathbb L/L_{cc})/(\mathbb L/L_{cc})_{ff}
 =\mathbb L/L_{\{c,f\}}
 =(\mathbb L/L_{ff})/(\mathbb L/L_{ff})_{cc}.
\end{equation}
Thus chord-first and chord-last eliminate
\eqref{eq:appC-common-chord-forest-matrix}, not merely two formally similar
energies.  Formula \eqref{eq:appC-affine-chord-energy} identifies the first
quotient with the advertised nonnegative covariant square, so
Proposition~\ref{prop:appC-positive-network-closure} applies.
\end{proof}

\subsection{Exact radial-conjugate free-chord restoration}
\label{appC:radial-conjugate-restoration}

The full chord above remains the coarse packet pivot.  On the overlapping
first-two-landing cell, the sharper chart retains the speed of the incoming
relative velocity at the second landing.  After the zero-set row operation
\eqref{eq:radial-row-operation} and the shear
\eqref{eq:radial-relative-shear}, its exact endpoint derivative is
\begin{equation}\label{eq:appC-radial-chord-matrix}
 \mathcal C_{\rm rad,d}=
 \begin{pmatrix}
  I_d&\tau_uI_d&0\\
  I_d&\tau_dI_d&g_2\\
  0&\widehat g_2^T&0
 \end{pmatrix},
 \qquad
 |\det\mathcal C_{\rm rad,d}|=\ell^{d-1}|g_2|.
\end{equation}
The last row is only the normalized determinant gauge.  The physical
symmetric constraint uses the energy row $g_2^T$ and the bordered block
\begin{equation}\label{eq:appC-radial-chord-KKT-block}
 \begin{pmatrix}
  -\ell I_d&g_2&A_\gamma\\
  g_2^T&0&0\\
  A_\gamma^*&0&H
 \end{pmatrix},
 \qquad A_\gamma=(-I_d,I_d)T_\gamma.
\end{equation}
Its upper-left inverse has current block
$-\ell^{-1}P_{g_2}$, and its exact Schur quotient is
\begin{equation}\label{eq:appC-radial-chord-Schur-term}
 H+\ell^{-1}A_\gamma^*P_{g_2}A_\gamma.
\end{equation}
On a finite orthogonal atlas choose
$E_{g_2}:\R^{d-1}\to g_2^\perp$.  Then
\eqref{eq:appC-radial-chord-Schur-term} is the factor contribution
\begin{equation}\label{eq:appC-radial-factor-row}
 \ell^{-1}(E_{g_2}^*A_\gamma)^*(E_{g_2}^*A_\gamma).
\end{equation}
Thus the factor-current space of this one local chord is $\R^{d-1}$; all
ordinary physical free segments and the preserved full-chord fallback keep
their $\R^d$ current spaces.  The local orthogonal frame does not enter a
recovery map.  Equations
\eqref{eq:appC-radial-chord-matrix}--
\eqref{eq:appC-radial-factor-row} are exactly the physical/KKT dictionary,
not a positive square added after elimination.

The trace left inverse is unchanged by this rank reduction.  Indeed, the
square selector matrix in
\eqref{eq:appC-trace-triangular-matrix} never selects the chord row, and
\eqref{eq:appC-missing-selector-inverse} reconstructs the missing physical
line from all the other line selectors.  Consequently the non-chord
selector/forest/ground subsystem already has trivial grounded kernel.  The
radial factor is an additional nonnegative row.  Schur-complement
associativity applied to the one bordered system
\eqref{eq:appC-radial-chord-KKT-block} proves radial-chord-first equals
forest-first.  This concrete atlas bypasses the arbitrary-row realization
problem; it does not contradict the no-go warning that positivity of a full
multiport response need not control a preferred pure physical bridge block.

This lemma is stronger than saying that a preferred raw minor receives a
positive matrix update.  Such a statement is generally coordinate
dependent and need not hold.  What is invariant is positivity of the Kron
stiffness, its grounded compliance, and the Riccati block in covariant
coordinates.

\subsection{Hard-sphere/positive-network block dictionary}
\label{appC:hard-sphere-dictionary}

We now identify every abstract block used above with a block of the physical
constraint derivative.  This subsection is redundant for the finite
dimensional algebra but removes any reliance on a verbal assertion that the
Schur complement is the Jacobi network.

\paragraph{Conventions.}
For an oriented particle pair \(e=(i,j)\), set
\[
 \beta_e(x_1,\ldots,x_P)=x_i-x_j,\qquad
 b_e=\beta_e\otimes I_d .
\]
Reversing \(e\) changes both the residual and its conormal by a minus sign
and therefore changes no Gram matrix or absolute determinant.  The
following table fixes all other conventions.

\paragraph{The pair-mixing block.}
Let $E$ be the kinetic-orthogonal complement of common translation in
$(\mathbb R^d)^P$, and define the normalized incidence
\begin{equation}\label{eq:appC-normalized-pair-incidence}
 B_e:=2^{-1/2}(e_i-e_j)^*\otimes I_d:E\longrightarrow\mathbb R^d.
\end{equation}
Thus $B_eB_e^*=I_d$ and $B_e^*B_e$ is the orthogonal projection onto the
relative-coordinate subspace of the pair.  If $P_\omega=\omega\otimes
\omega$, the fixed-normal equal-mass reflection on the two line variables is
\begin{equation}\label{eq:appC-pair-dimensional-collision}
 \mathscr R_e^{(2d)}=
 \begin{pmatrix}I_d-P_\omega&P_\omega\\
                 P_\omega&I_d-P_\omega\end{pmatrix}
 =U_e^*\begin{pmatrix}I_d&0\\0&I_d-2P_\omega\end{pmatrix}U_e,
 \qquad
 U_e=2^{-1/2}\begin{pmatrix}I_d&I_d\\I_d&-I_d\end{pmatrix}.
\end{equation}
Its extension to $E$ is
\begin{equation}\label{eq:appC-full-reduced-reflection}
 \mathscr R_e=I_E-2B_e^*P_\omega B_e.
\end{equation}
Equations \eqref{eq:appC-pair-dimensional-collision}--
\eqref{eq:appC-full-reduced-reflection} are orthogonal involutions.  They
show explicitly why a collision is not two independent $d$-dimensional
line gauges.

With the cylinder-inward normal, the curvature shear has range in the same
relative subspace.  In normalized relative coordinates it has the form
\begin{equation}\label{eq:appC-embedded-root-shunt}
 H_e=B_e^*\widehat H_eB_e,
 \qquad
 \widehat H_e=2|B_eV\cdot\omega|\,
       \widehat{\mathcal V}_e^*\widehat{\mathcal K}_e
       \widehat{\mathcal V}_e\ge0,
\end{equation}
with the harmless factors of $\sqrt2$ fixed by
\eqref{eq:appC-normalized-pair-incidence}.  Equivalently, without relative
coordinates, $H_e=2|V\cdot\nu_e|\mathcal V_e^*\mathcal K_e\mathcal V_e$.
The exact tangent saltation block on $E\oplus E$ is
\begin{equation}\label{eq:appC-full-saltation-block}
 D\mathscr S_e=
 \begin{pmatrix}\mathscr R_e&0\\
                 \mathscr R_eH_e&\mathscr R_e\end{pmatrix},
 \qquad
 (D\mathscr S_e)^*=
 \begin{pmatrix}\mathscr R_e&H_e\mathscr R_e\\
                 0&\mathscr R_e\end{pmatrix}.
\end{equation}

\begin{center}
\small
\begin{tabular}{@{}p{0.19\textwidth}p{0.31\textwidth}p{0.39\textwidth}@{}}
\toprule
object & convention & invariant consequence\\
\midrule
forward collision
 & \(V^+=V^--2(V^-\cdot\nu)\nu\)
 & configuration reflection
   \(\mathcal R=I-2\nu\otimes\nu\)\\
shape operator
 & \(\nu\) points into the excluded cylinder and
   \(\mathcal K=-D\nu|_{\nu^\perp}\ge0\)
 & the saltation shear is
   \(\mathcal H=2|V\cdot\nu|\mathcal V^*\mathcal K\mathcal V\ge0\)\\
backward conormal
 & adjoint pullback $(D\Phi)^*$ of the forward tangent map
 & free flight is \((Z,W)\mapsto(Z,W+rZ)\)\\
conormal/network variables
 & $Z=$ boundary current, $W=$ boundary potential after homogeneous
   transport
 & a free edge of length $r$ has compliance $rI$, not stiffness
   $r^{-1}I$\\
edge current
 & positive from tail to head
 & reversing an edge sends \(J_e\mapsto-J_e\) and leaves
   \(\ell_e|J_e|^2\) unchanged\\
Schur complement
 & eliminate internal columns:
   \(L/L_{II}=L_{BB}-L_{BI}L_{II}^{-1}L_{IB}\)
 & equals the Dirichlet-to-Neumann stiffness $S_B$\\
one-port inverse
 & $K_B=S_B^{-1}$ on the grounded port
 & equals the Neumann-to-Dirichlet compliance in the conormal velocity row\\
port gauge
 & \(\eta_v=U_v\xi_v\), \(U_v\) orthogonal
 & both $S_B$ and $K_B$ change by boundary orthogonal conjugacy only\\
ground
 & delete one common-translation coordinate in each free component
 & makes the internal Dirichlet block positive definite\\
\bottomrule
\end{tabular}
\end{center}

\paragraph{Rows, columns, and elimination order.}
Let \(A\) be the unmarked incoming contacts, \(T\) the retained upper tree,
and \(d_1,d_2\) the first two lower contacts.  In the overlap case remove
the last upper contact \(u_c\) of the new line and write
\(F=T\setminus\{u_c\}\).  After the unit-Jacobian elastic and transport
reconstructions, use the row order
\begin{equation}\label{eq:appC-physical-row-order}
 (C_A,\ U_F,\ (U_{u_c},R_{d_2}),\ R_{d_1})
\end{equation}
and the pivot-column order
\begin{equation}\label{eq:appC-physical-column-order}
 (y_A,\ h_F,\ (X_c,V_c),\ v_*).
\end{equation}
Here \(y_A=(t_a,\theta_a)_{a\in A}\), \(h_F\) are chronological
relative-cluster translations, and \(v_*\) is the $d$-coordinate minor
chosen from the remaining separator velocities.  Successive block
elimination gives
\begin{equation}\label{eq:appC-physical-block-factorization}
 \left|\det D_p\mathcal C\right|
 =
 \left|\det A_{\rm root}\right|
 \left|\det B_F\right|
 \left|\det C_{\rm ch}\right|
 \left|\det S_1\right|,
\end{equation}
where
\begin{align}
 A_{\rm root}
  &=\operatorname{Tri}_{a\in A}
       [\,g_a\ \ -\eps E_{\omega_a}\,],
 &|\det A_{\rm root}|
  &=\eps^{(d-1)|A|}\prod_{a\in A}|g_a\cdot\omega_a|,
 \label{eq:appC-root-block-dictionary}\\
 B_F
  &=\operatorname{Tri}_{e\in F}(\sigma_e I_d),
 &|\det B_F|&=1
 \quad\hbox{in raw cluster coordinates},\label{eq:appC-tree-block-dictionary}\\
 C_{\rm ch}
  &=\begin{pmatrix}
      I_d&(t_{u_c}-s)I_d\\
      I_d&(t_{d_2}-s)I_d
    \end{pmatrix},
 &|\det C_{\rm ch}|&=(t_{d_2}-t_{u_c})^d.
 \label{eq:appC-chord-block-dictionary}
\end{align}

The forest coordinates in this display are not informal single-line
ports.  If $\mathbf B_F:E\to(\mathbb R^d)^F$ is the stack of the normalized
incidences \eqref{eq:appC-normalized-pair-incidence}, chronological cluster
translation supplies a raw injection
\begin{equation}\label{eq:appC-forest-coordinate-injection}
 R_F^{\rm raw}:(\mathbb R^d)^F\longrightarrow E,
 \qquad
 \mathbf B_FR_F^{\rm raw}=T_F,
\end{equation}
where $T_F$ is block lower triangular with diagonal
$\sigma_eI_d$, $\sigma_e\in\{-1,1\}$.  After one translation is grounded in
each forest component, $T_F$ is square and invertible.  A QR factorization
$R_F^{\rm raw}=R_FG_F$ gives an isometric injection $R_F$ and a matrix
$G_F$ whose condition number, together with that of $T_F$, is bounded only
by $P_0$.  Thus
\begin{equation}\label{eq:appC-forest-boundary-splitting}
 E=\operatorname{Ran}R_F\mathbin{\oplus}E_B,
 \qquad \xi=R_Fh+R_Bz,
\end{equation}
for an orthogonal boundary injection $R_B$.  Equations
\eqref{eq:appC-forest-coordinate-injection}--
\eqref{eq:appC-forest-boundary-splitting}, rather than an independent
$d$-dimensional gauge on each particle line, are the precise meaning of
the block $B_F$ in \eqref{eq:appC-physical-block-factorization}.
The final block \(S_1\) is a selected square minor of the $d$-row physical
Schur derivative.  Its full row frame, before selecting columns, is
\begin{equation}\label{eq:appC-final-block-dictionary}
 S_1^*
 =(bI+\widetilde K_{\mathcal N})Z_1,
 \qquad b=t_{d_1}-s.
\end{equation}
Here $Z_1$ is the retained raw boundary-current frame (the position
conormal), $Q_-Z_1$ is its covariant current, the internal network contributes
the raw potential $Q_+K_{\mathcal N}Q_-Z_1$, and the independent lower free
interval contributes the raw potential $bZ_1$.  Thus
\eqref{eq:appC-network-CD-identification} is exactly the identity
needed in the last factor of
\eqref{eq:appC-physical-block-factorization}; it is not an analogy between
two unrelated matrices.

\mainstatementreference{Theorem~\ref{thm:appC-physical-schur-jacobi}}

\mainproofheading{Theorem~\ref{thm:appC-physical-schur-jacobi}}{supp:proof:thm:appC-physical-schur-jacobi}
\begin{proof}
First perform no elimination.  The local matrix dictionary, in covariant
coordinates, is
\begin{center}
\small
\begin{tabular}{@{}p{0.18\textwidth}p{0.36\textwidth}p{0.36\textwidth}@{}}
\toprule
event & physical tangent/conormal equations & network operation\\
\midrule
free segment $e=(u,v)$
 & $A_e\xi=\ell_ej_e$
 & row $-\ell_ej_e+(\mathbb A_W\xi)_e=0$\\
fixed-normal collision on $a=\{i,j\}$
 & $(\xi_i^+,\xi_j^+)=\mathscr R_a^{(2d)}
      (\xi_i^-,\xi_j^-)$, with the adjoint current law
 & one orthogonal $2d$-dimensional junction, eliminated only after both
   incident lines are present\\
incoming root
 & $\xi^-=\mathscr R_a\xi^+$ and
   $j^-=\mathscr R_aj^++H_a\mathscr R_a\xi^+$
 & complete pair junction followed by the nonnegative shunt $H_a$\\
marked forest merge
 & the new relative translation enters as $\sigma_e h_e$,
   $\sigma_e\in\{-1,1\}$
 & unit triangular pivot followed by vertex identification\\
free chord
 & $A_\gamma\xi=\ell_\gamma j_\gamma$ in its
   $d$-dimensional line-current space
 & one rank-$d$ factor row; Schur elimination in either order\\
radial-conjugate chord
 & the literal bordered system
   \eqref{eq:appC-radial-chord-KKT-block}
 & factor row $E_{g_2}^*A_\gamma$ in $\mathbb R^{d-1}$;
   exact Schur term \eqref{eq:appC-radial-chord-Schur-term}\\
\bottomrule
\end{tabular}
\end{center}
Every displayed physical equation is a full vector equation; hence the
table retains its off-diagonal derivatives with respect to all earlier
ports.

On one free segment, the tangent equations are
$A_e\xi=\ell_e j_e$ and current conservation is obtained by varying the
endpoint field.  Their quadratic primitive is
$\frac1{2\ell_e}|A_e\xi|^2$.  At a collision the simultaneous orthogonal
change \eqref{eq:appC-pair-dimensional-collision} is first made on the two
incident fields; substituting it into the adjacent drop rows is exactly the
construction of $\mathbb A_W$ in \eqref{eq:appC-factor-incidence}.  At an
incoming collision root, transport the outgoing current and field to the
incoming fixed-normal gauge, writing
$\widetilde j^+:=\mathscr R_a j^+$ and
$\xi_v:=\xi^-=\mathscr R_a\xi^+$.  The saltation calculation with the
cylinder-inward normal then gives the current jump
$j^- -\widetilde j^+=H_v\xi_v$ with $H_v$ in
\eqref{eq:appC-root-shunt}; its primitive is
$\frac12\langle H_v\xi_v,\xi_v\rangle$.  A fixed-normal elastic collision
changes the complete pair gauge by the orthogonal matrix
\eqref{eq:appC-pair-dimensional-collision}.  Summing these local
primitives over the word gives exactly
$\frac12\langle L_W\xi,\xi\rangle$ with $L_W$ in
\eqref{eq:appC-word-Jacobi-Hessian}.  This proves the matrix identity before
elimination, including signs and transposes.

The sign in the root row is fixed directly by the root equation.
With this cylinder-inward unit normal $\nu$, write the specular reflection as
$V^+=V^--2(V^-\!\cdot\nu)\nu$ and use
$D\nu|_{\nu^\perp}=-\mathcal K$.  Differentiating this identity, transporting
both sides to the same fixed-normal gauge, and pairing the tangential
position variation with the conormal velocity variation gives
\[
 \langle\xi_v,j^- -\widetilde j^+\rangle
 =2|V^-\!\cdot\nu|\,
   \langle\mathcal V\xi,\mathcal K\mathcal V\xi\rangle.
\]
Polarization yields precisely $j^- -\widetilde j^+=H_v\xi_v$ with
$H_v=2|V^-\cdot\nu|\mathcal V^*\mathcal K\mathcal V\ge0$.  Reversing the
normal reverses both the signed crossing velocity and $D\nu$; the displayed
 positive-shunt convention is specifically the cylinder-inward one just
 fixed.

We now give a completely defined variational identification.  Let
$\mathscr X_W$ be the Hilbert space of all endpoint fields
$q_{i,m}\in\mathbb R^d$ on the physical line segments of the fixed word,
modulo common translation.  At a retained forest contact its homogeneous
contact row is imposed as a hard linear constraint; equivalently, its
chronological relative-cluster translation has already been solved by
$R_F^{\rm raw}$.  At a fixed-normal collision $a=\{i,j\}$ impose
the exact junction relation
\begin{equation}\label{eq:appC-collision-junction-relation}
 \binom{q_{i,a}^+}{q_{j,a}^+}
 =\mathscr R_a^{(2d)}
   \binom{q_{i,a}^-}{q_{j,a}^-}.
\end{equation}
All other line coordinates pass unchanged.  For a physical free segment
$e$ from one endpoint to the next, let $A_e:\mathscr X_W\to\mathbb R^d$
be its exact endpoint drop after the junction relations
\eqref{eq:appC-collision-junction-relation} are substituted.  Hence an
outgoing drop row $A^+q_a^+$ becomes
$A^+\mathscr R_a^{(2d)}q_a^-$, including both incident lines.  Stacking these
rows, the chord row, and no others defines $\mathbb A_W$ in
\eqref{eq:appC-factor-incidence}.

Define the discrete action on the grounded affine space by
\begin{equation}\label{eq:appC-physical-discrete-action}
 \mathcal I_W(q)=\frac12\sum_e\ell_e^{-1}|A_e q|^2
 +\frac12\sum_{a\in A}\langle H_aq_a,q_a\rangle,
\end{equation}
where every $H_a$ is the embedded operator
\eqref{eq:appC-embedded-root-shunt}.  Its first variation is exactly the
second row of \eqref{eq:appC-word-KKT-system}, while
$j_e=\ell_e^{-1}A_e q$ is its first row.  Thus the KKT system is the
Euler--Lagrange system of a fully specified action.

We now connect this action to the physical bordered adjoint by an equality of
linear boundary relations.  For a chronological prefix $W_j$, let
$\mathcal R_j^{\rm phys}$ be the raw current/potential relation obtained by
Schur-eliminating the physical pivot rows in that prefix, and let
$\mathcal R_j^{\rm KKT}$ be the boundary relation of the corresponding
prefix of \eqref{eq:appC-word-KKT-system}.  Let $\mathsf P_j$ be the product of the
complete-pair and square forest primal coordinate maps through that prefix.
We prove the invariant
\begin{equation}\label{eq:appC-prefix-relation-invariant}
 (z_{\rm raw},w_{\rm raw})\in\mathcal R_j^{\rm phys}
 \quad\Longleftrightarrow\quad
 (\mathsf P_j^{-*}z_{\rm raw},\mathsf P_jw_{\rm raw})
       \in\mathcal R_j^{\rm KKT}.
\end{equation}
Our multiplication convention is fixed by the recursion
\begin{equation}\label{eq:appC-prefix-primal-gauge-recursion}
 \mathsf P_0=I,\qquad
 \mathsf P_{j+1}=\begin{cases}
  \mathscr R_a\mathsf P_j,&\text{complete pair junction at a fixed collision or root},\\
  P_{F,e}\mathsf P_j,&\text{forest merge},\\
  \mathsf P_j,&\text{free segment or chord},
 \end{cases}
\end{equation}
where
$P_{F,e}=[R_e^{\rm qr},R_{B,e}]\operatorname{diag}(G_e,I)$ is the square
completion of the raw forest injection at that merge.  Thus no length or
chord endpoint inverse enters $\mathsf P_j$, and at the final prefix
$Q_-=\mathsf P_W^{-*}$ and $Q_+=\mathsf P_W^{-1}$ exactly as in
\eqref{eq:appC-exact-dual-boundary-gauges}.
Both sides are relations, rather than matrices of artificially equal size.
The following five computations prove the induction and specify every change
of variables.

\emph{Free segment.}
For a segment of length $\ell$, the backward physical adjoint block and the
KKT edge law are, respectively,
\begin{equation}\label{eq:appC-prefix-free-matrix}
 \binom{j^-}{\xi^-}
 =\begin{pmatrix}I&0\\ \ell I&I\end{pmatrix}
   \binom{j^+}{\xi^+},
 \qquad
 -\ell j+\xi_v-\xi_u=0.
\end{equation}
Thus both say that the current is constant and the potential drop is
$\ell j$.

\emph{Fixed-normal collision and incoming root.}
For a root write $y_a=(t_a,\theta_a)$ and set
\[
 G_a=(D_{\theta_a}\omega_a)^*D_{\theta_a}\omega_a.
\]
Normalize the tangent differential by
$E_{\omega_a}=D_{\theta_a}\omega_aG_a^{-1/2}$.  Thus
$[g_a\ -\eps E_{\omega_a}]$ is the root block after the pointwise column
change $\delta\widehat\theta_a=G_a^{1/2}\delta\theta_a$; the chart density
$(\det G_a)^{1/2}$ is retained in the surface measure, and the Schur
complement is invariant under this invertible pivot-column change.  On a regular root
$g_a\cdot\omega_a\ne0$, and the inverse used here is explicit: for
$r\in\mathbb R^d$,
\begin{equation}\label{eq:appC-root-pivot-explicit-inverse}
 [g_a\ -\eps E_{\omega_a}]^{-1}r
 =\binom{
    (\omega_a\cdot r)/(g_a\cdot\omega_a)}{
    \begin{aligned}
    \eps^{-1}\bigl(&E_{\omega_a}^*g_a\,
      (\omega_a\cdot r)/(g_a\cdot\omega_a)\\[-2pt]
      &{}-E_{\omega_a}^*r\bigr)
    \end{aligned}}.
\end{equation}
Multiplication by $[g_a\ -\eps E_{\omega_a}]$ returns separately the
$\omega_a$ and $\omega_a^\perp$ parts of $r$.  Substitution of this inverse
into the constrained scattering derivative gives the literal IFT Schur matrix
\begin{equation}\label{eq:appC-root-IFT-saltation-identity}
 D\Phi_a^{\rm root}
 =D_x\Phi_a-D_{y_a}\Phi_a(D_{y_a}C_a)^{-1}D_xC_a
 =\begin{pmatrix}
    \mathscr R_a&0\\ \mathscr R_aH_a&\mathscr R_a
   \end{pmatrix}.
\end{equation}
The second equality follows by substituting
$V^+=V^--2(V^-\cdot\nu)\nu$ and
$D\nu|_{\nu^\perp}=-\mathcal K$ into the first matrix: its four blocks are
$\mathscr R_a$, $0$, $\mathscr R_a
 (2|V^-\cdot\nu|\mathcal V^*\mathcal K\mathcal V)$, and
$\mathscr R_a$.  Taking the adjoint gives, without a sign convention left
implicit,
\begin{equation}\label{eq:appC-root-adjoint-KKT-matrix}
 \binom{j^-}{\xi^-}
 =\begin{pmatrix}
    \mathscr R_a&H_a\mathscr R_a\\0&\mathscr R_a
   \end{pmatrix}
   \binom{j^+}{\xi^+},
 \quad
 \xi^-=\mathscr R_a\xi^+,
 \quad
 j^-=\mathscr R_aj^++H_a\mathscr R_a\xi^+.
\end{equation}
This is exactly the complete-pair KKT junction followed by the shunt $H_a$.
At a fixed-normal collision one sets $H_a=0$.

\emph{Forest merge.}
In the raw chronological coordinates write one physical row as
\[
 U_e=T_eh_e+V_er.
\]
Eliminating it is the explicit injection
\begin{equation}\label{eq:appC-forest-Schur-injection}
 h_e=-T_e^{-1}V_er,
 \qquad
 \binom{h_e}{r}=R_er,
 \qquad
 R_e=\binom{-T_e^{-1}V_e}{I}.
\end{equation}
For every later row $Y=(Y_h\ Y_r)$ the physical Schur row and its adjoint
are
\begin{equation}\label{eq:appC-forest-row-congruence}
 Y_r-Y_hT_e^{-1}V_e=YR_e,
 \qquad (YR_e)^*=R_e^*Y^*.
\end{equation}
The KKT quadratic form restricted to the same hard constraint is
$L\mapsto R_e^*LR_e$.  Hence the two eliminations use the identical
off-diagonal matrix $V_e$.  QR factorization of $R_e$ only changes $\mathsf P_j$ by
the square forest factor in
\eqref{eq:appC-positive-word-coordinate-bounds}.

\emph{Restored chord.}
Let
\begin{equation}\label{eq:appC-general-chord-trace-blocks}
 \begin{aligned}
 \mathcal T_\gamma&:\mathscr X_W\to\mathbb R^{2d},
 &\mathcal D&=(-I_d\ I_d),
 &A_\gamma&=\mathcal D\mathcal T_\gamma,\\
 \Pi_V&=(0\ I_d),
 &C_\gamma&=\begin{pmatrix}I&aI\\I&cI\end{pmatrix},
 &\ell_\gamma&=c-a>0.
 \end{aligned}
\end{equation}
Direct inversion gives the dimensionally exact identity
\begin{equation}\label{eq:appC-general-chord-velocity-identity}
 \Pi_VC_\gamma^{-1}\mathcal T_\gamma
 =\ell_\gamma^{-1}\mathcal D\mathcal T_\gamma
 =\ell_\gamma^{-1}A_\gamma.
\end{equation}
After all preceding pivots, let $Z_j$ be the complete current landing frame.
Causal endpoint reconstruction gives the two physical off-diagonal blocks
\begin{equation}\label{eq:appC-general-chord-offdiagonal-blocks}
 J_{\gamma,I}=-\mathcal T_\gamma,
 \qquad
 J_{B,\gamma}=Z_j^*\mathcal T_\gamma^*
      \mathcal D^*\Pi_V.
\end{equation}
Their physical Schur contraction is therefore
\begin{align}
 -J_{B,\gamma}C_\gamma^{-1}J_{\gamma,I}
 &=Z_j^*\mathcal T_\gamma^*\mathcal D^*\Pi_V
       C_\gamma^{-1}\mathcal T_\gamma\notag\\
 &=Z_j^*\bigl(\ell_\gamma^{-1}A_\gamma^*A_\gamma\bigr).
 \label{eq:appC-general-chord-Schur-KKT-identity}
\end{align}
On the KKT side, the chord block is exactly
\begin{equation}\label{eq:appC-general-chord-KKT-block}
 \begin{pmatrix}-\ell_\gamma I_d&A_\gamma\\
                 A_\gamma^*&0\end{pmatrix}.
\end{equation}
Eliminating its $d$-dimensional current gives the same positive operator
$\ell_\gamma^{-1}A_\gamma^*A_\gamma$.  Thus the chord endpoint determinant
$|\det C_\gamma|=\ell_\gamma^d$ is a physical Gaussian pivot, whereas
$\ell_\gamma^{-1}A_\gamma^*A_\gamma$ is the KKT Hessian contribution; neither
is absorbed into the boundary gauge.

In the radial-conjugate chart, replace this one prefix calculation, and only
this calculation, by the exact derivative
\eqref{eq:appC-radial-chord-matrix}.  With $T_\gamma$ understood as the
total endpoint trace after the shear, the columnwise physical inverse is
\[
 \Pi_g\mathcal C_{\rm rad,d}^{-1}
       \binom{T_\gamma}{0}
 =\ell^{-1}P_{g_2}A_\gamma.
\]
The adjoint off-diagonal contraction and the KKT chord step are therefore
the same bordered system \eqref{eq:appC-radial-chord-KKT-block}, and both
produce exactly
$\ell^{-1}A_\gamma^*P_{g_2}A_\gamma$.  No old full-chord term is included in
the remaining block $H$, so the contribution is counted once.  This proves
the prefix invariant for the radial chart; all later prefix steps are
identical to the full-chord proof.

\emph{Boundary step.}
For a landing multiplier $\mu$, put
\begin{equation}\label{eq:appC-prefix-boundary-current}
 z_{\rm raw}=Z_1\mu,
 \qquad i_B=Q_-Z_1\mu.
\end{equation}
The internal KKT response is $K_{\mathcal N}i_B$, and the independent lower
free interval contributes the raw potential $bZ_1\mu$.  By
\eqref{eq:appC-exact-dual-boundary-gauges}, the internal network contributes
the raw potential $Q_+K_{\mathcal N}Q_-Z_1\mu$.  The total raw output is
\begin{equation}\label{eq:appC-prefix-boundary-output}
 w_{\rm raw}=(bI+\widetilde K_{\mathcal N})Z_1\mu.
\end{equation}
The physical bordered relation
\eqref{eq:appC-bordered-adjoint-Schur-row} gives
$w_{\rm raw}=S_{\rm phys}^*\mu$.  Since the five displayed local identities
prove \eqref{eq:appC-prefix-relation-invariant} at every prefix, comparison
for arbitrary $\mu$ proves
\eqref{eq:appC-physical-network-final-identity}.

The grounded kernel condition is also explicit.  Let \(\mathbb U_F\) be the
stack of the \emph{physical} retained-forest tangent rows after all earlier
root substitutions.  These are precisely the hard constraints already
imposed in the definition of \(\mathscr X_W\); they are not extra factor
edges silently added to \(\mathbb A_W\).  At a chronological forest merge
their exact coordinate form is
\begin{equation}\label{eq:appC-physical-forest-row-coordinate}
 \mathbb U_Fq=T_Fh+T_{\rm old}q_{\rm old}+T_Bz .
\end{equation}
The matrices \(T_{\rm old},T_B\) are obtained by differentiating the same
physical contact residual; thus they contain every causal contribution from
earlier roots.  The new relative-cluster column has the invertible diagonal
\(T_F\) from \eqref{eq:appC-forest-coordinate-injection}.

We now prove the kernel assertion by induction through the actual word.
On an active component \(C\), stack all its exact particle selectors:
\[
 \mathbf S_C=(S_i)_{i\in C}:E_C\longrightarrow
              (\mathbb R^d)^C .
\]
On the zero-centre-of-mass space \(E_C\), this stack is injective and its
left-inverse norm is bounded by a function of \(|C|\le P_0\).  Hence the
vanishing of all physical-line factor rows in one free slab propagates the
entire reduced component field, not merely one particle coordinate.  A
fixed-normal event then propagates it by the invertible pair reflection
\eqref{eq:appC-pair-dimensional-collision}.  At a forest merge, if the old
component fields, the retained boundary field \(z\), and the physical row
\(\mathbb U_Fq\) vanish, equation
\eqref{eq:appC-physical-forest-row-coordinate} and invertibility of \(T_F\)
force the new cluster translation \(h\) to vanish.  The two chord endpoint
components have already been reached by this induction without using a
chord row.  Hence either the rank-$d$ full-chord row or the rank-$(d-1)$
radial row is only an additional nonnegative constraint and is not needed to
kill the grounded kernel.  Starting from the
declared ground in every free component and iterating these four steps
therefore reaches every endpoint field of \(W\).

The same induction, with triangle inequalities, gives the quantitative
estimate
\begin{equation}\label{eq:appC-actual-grounded-kernel}
 \|q\|\le C(P_0)\bigl(\|\mathbb A_Wq\|
                 +\|\mathbb U_Fq\|+\|z\|\bigr),
 \qquad
 \ker\mathbb A_W\cap\ker H_W\cap\ker\mathbb U_F\cap\{z=0\}
 =\{0\}.
\end{equation}
Here the constant accumulates the selector left inverses, the forest
\(T_F^{-1}\) bounds, and at most \(P_0\) invertible pair junctions.  The
nonnegative shunt term is harmless and has been retained in the displayed
kernel because it is part of \(L_W\).
Consequently, on the already constrained tangent space \(\mathscr X_W\),
the actual factor incidence---not graph connectivity alone---makes
\(L_{II}\) positive definite.

At an incoming root, if $y_a=(t_a,\theta_a)$ and $z$ denotes every earlier
column, substitution of the physical root changes every later row by the
same exact Schur derivative
\begin{equation}\label{eq:appC-root-local-Schur-check}
 D_zR_{\rm fut}-D_{y_a}R_{\rm fut}
       (D_{y_a}C_a)^{-1}D_zC_a.
\end{equation}
Differentiating $C_a=0$ identifies its adjoint with
\eqref{eq:appC-full-saltation-block}; its symmetric current jump is exactly
$H_a$.  Consequently all causal off-diagonal derivatives created by the
root are inserted into the later factor rows $A_e$ and the same shunt
$H_a$, rather than being discarded.

Induct now over the actual chronological word.  A free interval adds its
square $|A_e q|^2/(2\ell_e)$ and the free Green boundary term; a fixed-normal
collision substitutes \eqref{eq:appC-collision-junction-relation} and its
adjoint current law; an incoming root uses
\eqref{eq:appC-root-local-Schur-check}; a forest merge is the injective
coordinate substitution \eqref{eq:appC-forest-coordinate-injection}; and a
chord is the exact minimization
\eqref{eq:appC-affine-chord-energy}.  Each step preserves the identity
\begin{equation}\label{eq:appC-word-Green-induction}
 \delta\mathcal I_W(q)[\eta]
 =\sum_{v\in I}\langle (L_Wq)_v,\eta_v\rangle
   +\langle i_B,\eta_B\rangle
\end{equation}
for every test field $\eta$ satisfying the same junction relations.  This
is a chronological induction, so roots, merges and the chord may be
interspersed arbitrarily.

After the internal forest fields are minimized, uniqueness of the KKT
solution follows from the trace estimate
\eqref{eq:appC-certified-boundary-trace} together with the grounded kernel
estimate \eqref{eq:appC-actual-grounded-kernel}.  Equality with the physical
Schur derivative is not inferred from uniqueness: it is the prefix-by-prefix
matrix identity \eqref{eq:appC-prefix-relation-invariant}, whose root,
forest, and chord off-diagonal blocks are respectively
\eqref{eq:appC-root-IFT-saltation-identity},
\eqref{eq:appC-forest-row-congruence}, and
\eqref{eq:appC-general-chord-Schur-KKT-identity}.

The physical Gaussian operations are the unit block
matrices
\begin{equation}\label{eq:appC-explicit-Gaussian-matrices}
 E_L=\begin{pmatrix}I&0\\-C_WA_W^{-1}&I\end{pmatrix},
 \qquad
 E_R=\begin{pmatrix}I&-A_W^{-1}B_W\\0&I\end{pmatrix},
 \qquad
 E_LJ_WE_R=\begin{pmatrix}A_W&0\\0&S_{\rm phys}\end{pmatrix}.
\end{equation}
Chronological root, forest and chord ordering makes $A_W$ block lower
triangular with the invertible diagonal blocks
\eqref{eq:appC-root-block-dictionary}--
\eqref{eq:appC-chord-block-dictionary}.  Grounding makes $L_{II}>0$ because
each internal component is joined to a ground or to the retained port.
Finally the second row of \eqref{eq:appC-word-KKT-system}, tested against an
arbitrary field $\eta$, gives the exact discrete Green identity
\[
 \sum_e\langle j_e,(\mathbb A_W\eta)_e\rangle
 +\sum_v\langle H_v\xi_v,\eta_v\rangle
 =\langle i_B,\eta_B\rangle.
\]
Solving the internal equations shows that boundary potential and current
obey $i_B=S_{\mathcal N}\xi_B$ and
$\xi_B=K_{\mathcal N}i_B$.  The complete pair junctions and the forest QR
splitting identify the raw physical current with $i_B=Q_-Z_1$; by duality
the corresponding raw Schur row is obtained by applying $Q_+$.  The lower
collision-free interval has length $b$ and contributes the independent
raw potential $bZ_1$.  Hence the Green boundary derivative is
$(bI+\widetilde K_{\mathcal N})Z_1$, which is
\eqref{eq:appC-physical-network-final-identity}.  This conclusion uses the
same Gaussian quotient on both sides and never identifies an instantaneous
historical chord endpoint with a phase-frontier potential.
\end{proof}

\subsection{Simultaneous birth-segment ports}
\label{appC:simultaneous-birth-ports}

The one-port identity above has a finite-port extension when the ports are
created by first lower appearances.  The point which makes the extension
work is not positivity alone.  Raw later landing multipliers can send
adjoint current through earlier birth segments.  We first remove this
causal mixing by one determinant-one row operation and only then take the
multiport Schur quotient.

\begin{lemma}[Causal clearing of birth-segment currents]
\label{lem:appC-causal-birth-current-clearing}
Let $d_1<\cdots<d_k$ be lower contacts and, for each $i$, let $c_i$ be a
distinct physical line whose first lower contact is $d_i$.  At most one
line is selected at each $d_i$.  Let $\gamma_i$ be the collision-free
segment of $c_i$ from the separator $s$ to $d_i$, put
\[
 \Gamma:=\{\gamma_1,\ldots,\gamma_k\},
 \qquad
 \Delta_i=t_{d_i}-s,
 \qquad
 \mathcal Y=\bigoplus_{i=1}^k\mathbb R^d,
 \qquad
 D_\Delta=\operatorname{diag}_{i}(\Delta_iI_d),
\]
and write $R=(R_1,\ldots,R_k)$ for the selected landing displacements after
all complementary incoming-root substitutions.  Put
\begin{equation}\label{eq:appC-birth-true-contact-vector}
 \Omega=(\omega_1,\ldots,\omega_k),
 \qquad C:=R-\eps\Omega.
\end{equation}
Thus the physical selected-contact fibre is $C=0$, equivalently
$R=\eps\Omega$.  Every derivative in this lemma is the resulting total
derivative, and the selected normals are retained unless stated otherwise.

For a raw landing multiplier $\mu=(\mu_i)\in\mathcal Y$, propagate the
bordered adjoint current backward through the complete lower word, using
the full $2d$-dimensional pair junction at every collision.  If $y_j$ is
the total current carried by $\gamma_j$, then
\begin{equation}\label{eq:appC-birth-current-transfer}
 y=\mathsf H\mu,
 \qquad
 \mathsf H=
 \begin{pmatrix}
 I&H_{12}&\cdots&H_{1k}\\
 0&I&\cdots&H_{2k}\\
 \vdots&&\ddots&\vdots\\
 0&\cdots&0&I
 \end{pmatrix}.
\end{equation}
Thus $\mathsf H$ is block upper unitriangular.  Set
\begin{equation}\label{eq:appC-birth-clearing-row-map}
 \mathsf U:=\mathsf H^{-*},
 \qquad \widetilde C:=\mathsf U C,
 \qquad
 \widetilde R^{\rm zs}:=\eps\Omega+\widetilde C
 =\eps\Omega+\mathsf U(R-\eps\Omega).
\end{equation}
Then $\mathsf U$ is block lower unitriangular and $\det\mathsf U=1$, while
\[
 \widetilde C=0\Longleftrightarrow C=0,
 \qquad
 \widetilde R^{\rm zs}=\eps\Omega
     \Longleftrightarrow R=\eps\Omega.
\]
The multiplier $y$ of the cleared constraint $\widetilde C$ induces the raw
multiplier $\mu=\mathsf U^*y=\mathsf H^{-1}y$, and the vector of its currents
on the selected segments is exactly $y$.  In particular, a multiplier
supported in the $i$th cleared block has current $y_i$ on $\gamma_i$ and
zero current on every $\gamma_j$, $j\ne i$.

Cut all the $\gamma_i$ and keep the remainder endpoint potentials $\xi$
independent.  Let
\[
 \mathcal G:\mathbb R^k\longrightarrow\mathcal Y,
 \qquad \mathcal G\tau=(g_i\tau_i)_{i=1}^k,
\]
where $g_i$ is the incoming relative velocity at $d_i$.  There is a
physical covariant current injection
$F:\mathcal Y\to\mathscr X_{\rm rem}$ such that the exact cut virtual-work
relation is
\begin{equation}\label{eq:appC-cleared-cut-derivative}
 D_{(V_c,\vartheta,\xi)}\widetilde C_{\rm cut}
 =D_{(V_c,\vartheta,\xi)}\widetilde R^{\rm zs}_{\rm cut}
   =\bigl(D_\Delta\ \ \mathcal G\ \ F^*\bigr),
 \qquad
 V_c=(V_{c_1},\ldots,V_{c_k}),
 \quad \vartheta=(t_{d_1},\ldots,t_{d_k}).
\end{equation}
Here $V_{c_i}$ are the original Cartesian separator velocities; no speed
or radial coordinate is introduced into the global recovery map.
If $\mathsf U$ depends on retained parameters, then on the actual common
contact fibre $C=0$
\begin{equation}\label{eq:appC-clearing-product-rule}
 D_z\widetilde C=(D_z\mathsf U)C+\mathsf U D_zC
                   =\mathsf U D_zC.
\end{equation}
Consequently the clearing preserves the physical zero set, every full
square determinant, and the coarea density.  In particular,
\(\delta^{(kd)}(\widetilde C)=\delta^{(kd)}(C)\) because
\(|\det\mathsf U|=1\).  The auxiliary cleared rows are used only in this
local Jacobian proof; the algebraic recovery system of Section~\ref{app:coarea} continues
to use the original polynomial equations $C_i=0$.
\end{lemma}

\begin{proof}
A multiplier at $d_i$ can propagate only through contacts earlier than
$d_i$.  It may therefore enter $\gamma_j$ only when $j\le i$.  Its current
on its own first-lower segment is $\mu_i$ with coefficient $I_d$.  Since
only one incident line is selected at a birth contact, no second selected
segment occurs in this diagonal block.  This proves
\eqref{eq:appC-birth-current-transfer}; the statements about $\mathsf U$
and the cleared segment currents follow by inversion.

It remains to identify all three blocks in
\eqref{eq:appC-cleared-cut-derivative}.  Pair the complete linearized word
with the bordered adjoint solution generated by the raw multiplier
$\mu=\mathsf H^{-1}y$.  On an internal free segment the endpoint work
cancels by current conservation.  At a fixed-normal collision it cancels
because the tangent and current laws are adjoints of the same complete
orthogonal pair matrix
\eqref{eq:appC-pair-dimensional-collision}.  At an incoming root, the
root multiplier cancels the variation of its contact equation; this is
the adjoint form of
\eqref{eq:appC-root-IFT-saltation-identity}.  At a marked forest row the
same cancellation is
\eqref{eq:appC-forest-row-congruence}.  We do not set the cut remainder
potential to zero: its uncancelled boundary work is, by definition,
$\langle Fy,\xi\rangle$.

The selected line $c_i$ has no lower collision before $d_i$.  Hence a
variation $\delta V_{c_i}$ changes the drop across $\gamma_i$ by
$\Delta_i\delta V_{c_i}$, while the companion endpoint before $d_i$ is
independent of $V_{c_i}$ by
Lemma~\ref{lem:lower-influence-cone}.  Moving the landing time changes the
same endpoint residual by $g_i\delta t_{d_i}$.  This is a simultaneous
moving-junction identity: for a raw multiplier $\mu$, the work of moving
$d_j$ is paired with the actual segment current
$(\mathsf H\mu)_j$, including contributions from every later raw landing,
not merely with the $j$th raw row.  Because the cleared current on
$\gamma_i$ is $y_i$, the remaining boundary work is therefore
\[
 \sum_i\langle y_i,\Delta_i\delta V_{c_i}
                         +g_i\delta t_{d_i}\rangle
 +\langle Fy,\xi\rangle.
\]
Dualizing this identity gives
\eqref{eq:appC-cleared-cut-derivative}, since the selected normals are fixed
in the displayed pivot derivatives and hence $D C=D R$ in those columns.
The last assertion is the product rule
\eqref{eq:appC-clearing-product-rule}; $\det\mathsf U=1$ makes the delta and
determinant statements exact.
\end{proof}

\begin{lemma}[Finite-prefix cut/KKT identification]
\label{lem:appC-birth-finite-prefix-invariant}
Order the rooted and forest-reduced physical word chronologically, cut the
selected birth segments when their lower endpoints are reached, and use the
true cleared contact equations $\widetilde C=\mathsf U C$ from
Lemma~\ref{lem:appC-causal-birth-current-clearing}.  After every finite
prefix $W_j$, let $\mathscr X_j$ be the covariant potential space after all
complete-pair and forest substitutions already encountered.  If
$\mathscr A_j$ is the set of incoming roots already encountered, let
$T_{a,j}:\mathscr X_j\to E_a$ be the exact covariant trace at root $a$ and
define the stacked physical shunt factor
\[
 \mathbb C_jq:=(H_a^{1/2}T_{a,j}q)_{a\in\mathscr A_j},
 \qquad H_j:=\mathbb C_j^*\mathbb C_j\ge0.
\]
There are
operators
\[
 \mathbb A_j:\mathscr X_j\longrightarrow\mathscr E_j,
 \quad F_j:\mathcal Y_j\longrightarrow\mathscr X_j,
 \quad D_j=D_j^*,
 \quad \mathcal G_j:\mathbb R^{I_j}\longrightarrow\mathcal Y_j
\]
and an exact square raw-to-covariant potential gauge $\mathsf P_j$ for which the
physical bordered relation at that prefix is equivalent to
\begin{align}
 \rho_j&=D_jy+\mathcal G_j\tau+F_j^*\xi,
       \label{eq:appC-birth-prefix-output}\\
 \sigma_j&=\mathcal G_j^*y,
       \label{eq:appC-birth-prefix-time-row}\\
 i_j&=F_jy-\mathbb A_j^*\mathbf j-\mathbb C_j^*r,
       \label{eq:appC-birth-prefix-kirchhoff}\\
 \Lambda_j\mathbf j&=\mathbb A_j\xi,
 \qquad r=\mathbb C_j\xi .
       \label{eq:appC-birth-prefix-edge-shunt-laws}
\end{align}
Here $I_j$ is the set of selected landings already encountered,
$\mathcal Y_j=\bigoplus_{i\in I_j}\mathbb R^d$, $\mathbf j$ contains every
retained physical edge current, and the $a$th component of $r$ is the
genuine factor current $H_a^{1/2}T_{a,j}\xi$ in the $a$th incoming-root
shunt.  In particular no root-shunt current is set to
zero.  The raw and covariant dual variables satisfy, at every prefix,
\begin{equation}\label{eq:appC-birth-prefix-exact-gauges}
 z_j^{\rm cov}=\mathsf P_j^{-*}z_j^{\rm raw},
 \qquad w_j^{\rm cov}=\mathsf P_jw_j^{\rm raw}.
\end{equation}

At the terminal prefix,
\begin{equation}\label{eq:appC-birth-prefix-terminal-blocks}
 D_j=D_\Delta,
 \qquad \mathcal G_j=\mathcal G,
 \qquad F_j=F,
 \qquad \mathbb C_j=\mathbb C_{\rm rem},
 \qquad
 L_j:=\mathbb A_j^*\Lambda_j^{-1}\mathbb A_j
       +\mathbb C_j^*\mathbb C_j=L_{\rm rem}.
\end{equation}
Put
\begin{equation}\label{eq:appC-birth-extended-factor}
 \mathbb B_{\rm rem}q
 :=\binom{\Lambda_{\rm rem}^{-1/2}\mathbb A_{\rm rem}q}
              {\mathbb C_{\rm rem}q},
 \qquad L_{\rm rem}=\mathbb B_{\rm rem}^*\mathbb B_{\rm rem}.
\end{equation}
For each block inclusion $\iota_i:\mathbb R^d\to\mathcal Y$, the physical
prefix induction produces an augmented current map $\mathcal J_i$ such that
\begin{equation}\label{eq:appC-birth-extended-current-divergence}
 F\iota_i=\mathbb B_{\rm rem}^*\mathcal J_i.
\end{equation}
Consequently
\begin{equation}\label{eq:appC-birth-prefix-range}
 F(\mathcal Y)\subseteq\operatorname{Ran}\mathbb B_{\rm rem}^*
 =\operatorname{Ran}L_{\rm rem}.
\end{equation}

Let $\widehat Z_\Gamma:\mathcal Y\to E_{\rm cov}$ be the terminal covariant
Riesz current frame generated by the cleared unit cut currents, and let
$\mathsf P_\Gamma$ be the terminal gauge in
\eqref{eq:appC-birth-prefix-exact-gauges}.  Define the raw Cartesian frame by
\begin{equation}\label{eq:appC-birth-raw-current-frame}
 Z_\Gamma:=\mathsf P_\Gamma^*\widehat Z_\Gamma,
 \qquad \mathsf P_\Gamma^{-*}Z_\Gamma=\widehat Z_\Gamma.
\end{equation}
Then, with $K_{\rm rem}=L_{\rm rem}^{+}$, the terminal physical derivative
obeys the literal identity
\begin{equation}\label{eq:appC-birth-prefix-terminal-identity}
 \boxed{\ (\mathsf U\mathcal S_k)Z_\Gamma
   =D_\Delta+F^*K_{\rm rem}F.\ }
\end{equation}
Moreover
\begin{equation}\label{eq:appC-birth-prefix-frame-bound}
 \Jac_{k(d-1)}(Z_\Gamma E)\le C(d,P_0,k).
\end{equation}
Every matrix in these identities is obtained from the actual bordered
physical derivative; none is a comparison-network response.
\end{lemma}

\begin{proof}
We prove the assertion by induction over the finite word.  At the empty
prefix $\mathsf P_0=I$ and all displayed spaces and maps are zero.  Suppose the
relations hold immediately before the next operation.

\emph{Free transport.}
For a free segment $e$ of length $\ell_e$, the physical conormal block and
the KKT edge law are exactly
\[
 \binom{\mathbf j^-}{\xi^-}
 =\begin{pmatrix}I&0\\ \ell_eI&I\end{pmatrix}
  \binom{\mathbf j^+}{\xi^+},
 \qquad -\ell_e\mathbf j_e+(\mathbb A_j\xi)_e=0,
\]
as in \eqref{eq:appC-prefix-free-matrix}.  Thus the endpoint works cancel,
the current is unchanged, and the new factor row is precisely the physical
drop row.  For a designated selected segment its physical lengths are
accumulated in the corresponding direct block; every other row is adjoined
to $\mathbb A_j$.

\emph{Complete-pair collision.}
Insert both incident $d$-fields before eliminating the junction.  Potentials
and dual currents are transported by the same orthogonal involution
$\mathscr R_a^{(2d)}$ from
\eqref{eq:appC-pair-dimensional-collision}.  Hence all incident endpoint
 works cancel simultaneously, including cross-line terms, and
 $\mathsf P_{j+1}=\mathscr R_a\mathsf P_j$.  No one-line projection is inverted or omitted.
 Under this same orthogonal identification, every pre-existing row of
 $\mathbb A_j$ and $\mathbb C_j$ and every boundary-work block of $F_j$ is
 pulled back before the next physical row is appended.

\emph{Incoming root.}
The exact IFT computation
\eqref{eq:appC-root-IFT-saltation-identity} and its adjoint
\eqref{eq:appC-root-adjoint-KKT-matrix} give
\[
 \xi^- =\mathscr R_a\xi^+,
 \qquad
 \mathbf j^- =\mathscr R_a\mathbf j^+
       +H_a\mathscr R_a\xi^+.
\]
The orthogonal pair part updates $\mathsf P_j$ exactly as above, while the second
term appends the genuine factor row $H_a^{1/2}T_{a,j+1}$ to
$\mathbb C_{j+1}$.  Thus all derivatives of
later rows created by the root Schur substitution occur in the updated
$\mathbb A_j,\mathbb C_j,F_j$; no causal off-diagonal block is discarded.

\emph{Forest merge.}
For the physical row $U_e=T_eh_e+V_er$, elimination is the injection
$R_e=(-T_e^{-1}V_e,I)^T$.  Equation
\eqref{eq:appC-forest-row-congruence} gives simultaneously
\[
 Y\longmapsto YR_e,
 \quad
 \mathbb A_j\longmapsto\mathbb A_jR_e,
 \quad
 \mathbb C_j\longmapsto\mathbb C_jR_e,
 \quad
 F_j^*\longmapsto F_j^*R_e.
\]
The dual current map is $R_e^*$ and the square QR completion is the factor
$P_{F,e}$ in \eqref{eq:appC-prefix-primal-gauge-recursion}.  Hence the raw
and covariant pairings in
\eqref{eq:appC-birth-prefix-exact-gauges} remain exactly equal.

\emph{Selected birth cut and causal clearing.}
When $d_i$ is reached, remove the interior of $\gamma_i$ but retain its two
endpoint fields and its current as the $i$th port.  The free drop is
$\Delta_i y_i$, motion of the lower endpoint contributes
$g_i\,\delta t_{d_i}$, and all remaining uncancelled endpoint work is, by
definition, $\langle F_i y_i,\xi\rangle$.  Here the direct term records the
already accumulated identity
$\sum_{e\subset\gamma_i}\ell_e=\Delta_i$; no second copy of a free-segment
contribution is added at the cut.  A later raw landing may carry
current through an earlier selected segment; the exact triangular law is
$y=\mathsf H\mu$.  Multiplying the true equations $C=0$ by
$\mathsf U=\mathsf H^{-*}$ changes the multiplier to
$\mu=\mathsf H^{-1}y$ and therefore makes the selected cut currents exactly
$y$.  Equation \eqref{eq:appC-clearing-product-rule} shows that this
operation has no parameter-derivative remainder on the contact fibre.
This proves \eqref{eq:appC-birth-prefix-output}--
\eqref{eq:appC-birth-prefix-edge-shunt-laws} after the cut, including the
complete off-diagonal contribution in $F_i$.

These five updates exhaust every operation in the rooted/forest-reduced
word, so induction proves the expanded prefix relation and
\eqref{eq:appC-birth-prefix-terminal-blocks}.  At the terminal prefix the
remainder has no external current other than the selected cut-port
injections, so $i_j=0$.  Apply the relation to the cleared multiplier
$y=\iota_i a$.  Before using any inverse,
the physical bordered adjoint gives edge currents $\mathbf j_i(a)$, a
remainder potential $\xi_i(a)$, and the actual shunt current
$r_i(a)=\mathbb C_{\rm rem}\xi_i(a)$.  The terminal Kirchhoff row is
\[
 F_i a=\mathbb A_{\rm rem}^*\mathbf j_i(a)
       +\mathbb C_{\rm rem}^*r_i(a).
\]
With
\[
 \mathcal J_i(a):=
 \binom{\Lambda_{\rm rem}^{1/2}\mathbf j_i(a)}{r_i(a)},
\]
this is exactly
\eqref{eq:appC-birth-extended-current-divergence}.  In finite dimension
$\operatorname{Ran}\mathbb B^*=\operatorname{Ran}(\mathbb B^*\mathbb B)$;
equivalently both spaces are $(\ker\mathbb B)^\perp$.  This proves
\eqref{eq:appC-birth-prefix-range} without deleting a root-shunt current.

The last three equations in the terminal prefix system give
\[
 Fy-L_{\rm rem}\xi=0,
 \qquad \xi=K_{\rm rem}Fy+\xi_0,
 \qquad \xi_0\in\ker L_{\rm rem}.
\]
The range identity makes this solvable for every $y$, and
$F^*\xi_0=0$.  Substitution in the first row gives
$D_\Delta+F^*K_{\rm rem}F$.  Here $K_{\rm rem}=L_{\rm rem}^+$ is taken in
the final covariant Hilbert quotient; it is never transported by a naive
Moore--Penrose congruence through a nonorthogonal forest gauge.  If
$W_\Gamma=(\mathsf U\mathcal S_k)^*$, exact raw/covariant duality gives,
for every $y,z\in\mathcal Y$,
\[
 \langle Z_\Gamma y,W_\Gamma z\rangle_{\rm raw}
 =\langle y,D_\Delta z\rangle
   +\langle Fy,K_{\rm rem}Fz\rangle.
\]
Since the right side is symmetric, dualizing proves
\eqref{eq:appC-birth-prefix-terminal-identity} with the stated orientation.
Finally, $\widehat Z_\Gamma$ is, by definition, the orthogonal direct sum of
the cleared unit cut-current blocks.  The gauge $\mathsf P_\Gamma$ contains only at
most $P_0$ complete orthogonal pair maps and bounded square forest factors.
The landing-row clearing acts on $\mathsf U\mathcal S_k$, not on this
raw/covariant boundary gauge, and the root curvature remains in the shunt
factor.  The coordinate bounds
\eqref{eq:appC-exact-dual-boundary-gauges}--
\eqref{eq:appC-positive-word-coordinate-bounds} therefore give
\eqref{eq:appC-birth-prefix-frame-bound}.
\end{proof}

\mainstatementreference{Theorem~\ref{thm:appC-birth-multiport-schur-jacobi}}

\mainproofheading{Theorem~\ref{thm:appC-birth-multiport-schur-jacobi}}{supp:proof:thm:appC-birth-multiport-schur-jacobi}
\begin{proof}
The finite-prefix physical/KKT induction in
Lemma~\ref{lem:appC-birth-finite-prefix-invariant} retains the genuine
root-shunt current as well as every physical edge current.  Its terminal
identity \eqref{eq:appC-birth-extended-current-divergence} is precisely
\eqref{eq:appC-birth-current-divergence}.  Since
$L_{\rm rem}=\mathbb B_{\rm rem}^*\mathbb B_{\rm rem}$ and
$\operatorname{Ran}\mathbb B_{\rm rem}^*
=\operatorname{Ran}(\mathbb B_{\rm rem}^*\mathbb B_{\rm rem})$ in finite
dimension, this proves \eqref{eq:appC-birth-current-range}.  Equivalently,
if $q\in\ker L_{\rm rem}$ then $\mathbb B_{\rm rem}q=0$ and
$\langle Fy,q\rangle=0$ for all $y$.  Thus the Moore--Penrose quotient below
does not discard a compatibility equation.

Before eliminating the remainder potential, the first row of
\eqref{eq:appC-birth-common-radial-KKT} is exactly
\eqref{eq:appC-cleared-cut-derivative}: the selected free segments give
$D_\Delta$, the landing-time columns give $\mathcal G$, and the remainder
boundary work gives $F^*$.  The last row is the Euler--Lagrange/Kirchhoff
equation $Fy-L_{\rm rem}\xi=0$ of
\eqref{eq:appC-birth-remainder-Hessian}, with the potential orientation
chosen consistently with the cut work.  The middle row is the symmetric
energy-row landing-time constraint.  At every chronological prefix,
the free, complete-pair, root and forest blocks are respectively
\eqref{eq:appC-prefix-free-matrix},
\eqref{eq:appC-pair-dimensional-collision},
\eqref{eq:appC-root-IFT-saltation-identity}, and
\eqref{eq:appC-forest-row-congruence}.  Hence
\eqref{eq:appC-birth-common-radial-KKT} is a quotient of the same physical
bordered derivative, entry by entry, including all causal off-diagonal
terms.  It is not a comparison network.

Since $FE$ lies in the range of $L_{\rm rem}$, eliminate the remainder
field with $K_{\rm rem}=L_{\rm rem}^{+}$ and restrict the current to
$\ker\mathcal G^*=\operatorname{Ran}E$.  This gives exactly
\eqref{eq:appC-birth-tangential-response}.  Positivity and
\eqref{eq:appC-birth-tangential-determinant} follow by congruence and
eigenvalue monotonicity.

It remains to relate this coefficient-space determinant to the physical
mixed Jacobian.  The clearing matrix $\mathsf U$ acts on all $kd$ landing
rows, has determinant one, and is applied before the root and forest Schur
quotients.  Thus for every $k(d-1)$-set $I$ of Cartesian velocity columns,
\[
 \det[(\mathsf U\mathcal S_k)_I\ \mathsf U\mathcal T_k]
 =\det[(\mathcal S_k)_I\ \mathcal T_k].
\]
In the cleared rows $\mathsf U\mathcal T_k=\mathcal G$.  The orthogonal
complement of this time-column range is $\operatorname{Ran}E$.  Let
$Z_\Gamma$ be the raw Cartesian Riesz current frame defined, with its exact
raw/covariant orientation, in
\eqref{eq:appC-birth-raw-current-frame}.  The finite-prefix theorem gives
\eqref{eq:appC-birth-prefix-terminal-identity}; restricting it to
$\operatorname{Ran}E$ gives the literal identities
\begin{equation}\label{eq:appC-birth-physical-column-identity}
 (\mathsf U\mathcal S_k)Z_\Gamma
 =D_\Delta+F^*K_{\rm rem}F=:M_\Gamma,
 \qquad
 (E^*\mathsf U\mathcal S_k)(Z_\Gamma E)=M_{\Gamma,\perp}.
\end{equation}
The first equality contains the same free, complete-pair, incoming-root,
forest, and selected-cut operations as $\mathcal S_k$, including the true
root shunts and every causal off-diagonal block.  It is the Green pairing of
the physical landing conormal with its exact dual Riesz frame, not an
identification of unrelated gauges.  The bound
\eqref{eq:appC-birth-prefix-frame-bound} is
\begin{equation}\label{eq:appC-birth-column-volume-bound}
 \Jac_{k(d-1)}(Z_\Gamma E)\le C(d,P_0,k).
\end{equation}
Cauchy--Binet applied to the second equality in
\eqref{eq:appC-birth-physical-column-identity}, in the original Cartesian
velocity coordinates, therefore gives
\begin{align*}
 \det M_{\Gamma,\perp}
 &\le \Jac_{k(d-1)}(E^*\mathsf U\mathcal S_k)
       \Jac_{k(d-1)}(Z_\Gamma E)\\
 &\le C(d,P_0,k)
       \Jac_{k(d-1)}(E^*\mathsf U\mathcal S_k),
\end{align*}
and hence
\[
 \Jac_{k(d-1)}(E^*\mathsf U\mathcal S_k)
 \ge c(d,P_0,k)\det M_{\Gamma,\perp}.
\]
 The exact time-wedge identity
\eqref{eq:k-time-wedge-factorization} gives
the additional factor
$\sqrt{\det(\mathcal G^*\mathcal G)}=\prod_i|g_i|$.
Together with \eqref{eq:appC-birth-tangential-determinant}, this proves the
first inequality in \eqref{eq:appC-birth-mixed-frame}; the second uses
$|g_i|\ge|g_i\cdot\omega_i|$.

The common system also contains the local radial block
\[
 \begin{pmatrix}\Delta_iI_d&g_i\\g_i^*&0\end{pmatrix}.
\]
Its energy-row determinant has modulus
$\Delta_i^{d-1}|g_i|^2$.  Replacing the last row by the normalized-speed
row $\widehat g_i^*$ divides this by $|g_i|$ and gives
$\Delta_i^{d-1}|g_i|$.  Eliminating the energy-row block gives the current
part $\Delta_i^{-1}P_{g_i}\succeq0$, exactly as in
Lemma~\ref{lem:radial-conjugate-chord-block}.  Since all these blocks occur
in the one matrix \eqref{eq:appC-birth-common-radial-KKT}, this is a joint
radial augmentation and not a product of unrelated local estimates.  No
$\widehat g_i$ is defined on $g_i=0$; those points, grazing points, chart
boundaries and rank-loss loci belong to the existing zero-flux/null strata.
\end{proof}

\begin{remark}[One-, two-, and three-port regressions]
\label{rem:appC-birth-port-regressions}
For $k=1$, $\mathsf H=I$ and the Schur response is
$\Delta_1I+F_1^*L_{\rm rem}^{+}F_1$; after the one-port boundary gauges this
is the identity $(bI+\widetilde K_{\mathcal N})Z_1$ in
\eqref{eq:appC-physical-network-final-identity}.  For $k=2$, all overlap
terms occur inside the single congruence
$(FE)^*L_{\rm rem}^{+}(FE)\ge0$; no oblique projection and no product of two
separately reduced chord determinants is used.  For $k=3$, writing
$\mathsf H_{12},\mathsf H_{13},\mathsf H_{23}$ for the three strict upper
blocks, the three cleared basis multipliers pull back to raw multipliers
\[
 (z,0,0),\qquad (-\mathsf H_{12}z,z,0),\qquad
 ((\mathsf H_{12}\mathsf H_{23}-\mathsf H_{13})z,
   -\mathsf H_{23}z,z).
\]
Multiplication by $\mathsf H$ gives selected-segment currents
$(z,0,0)$, $(0,z,0)$, and $(0,0,z)$, respectively.  Moreover, writing
$D=\operatorname{diag}(\Delta_1I_{d-1},\Delta_2I_{d-1},
\Delta_3I_{d-1})$ gives
\[
 M_{\Gamma,\perp}=D+(FE)^*L_{\rm rem}^{+}(FE)\succeq D,
 \qquad
 \det M_{\Gamma,\perp}\ge
 (\Delta_1\Delta_2\Delta_3)^{d-1}.
\]
This is the first new gain-producing case beyond the two-port theorem.  The
same argument contains every fixed $k$ without an induction on principal
minors.
\end{remark}

\begin{lemma}[Block dictionary and elimination commutation]
\label{lem:appC-block-dictionary}
With the row and column orders
\eqref{eq:appC-physical-row-order}--\eqref{eq:appC-physical-column-order},
the physical hard-sphere derivative and the covariant network reduction fit
the commutative diagram
\begin{equation}\label{eq:appC-elimination-diagram}
\begin{array}{ccc}
D_p(C_A,U_F,U_c,R_2,R_1)
 &\xrightarrow{\quad /\,A_{\rm root},B_F\quad}&
D_{(X_c,V_c,v_*)}(U_c,R_2,R_1)\\[2mm]
\Big\downarrow\hbox{\scriptsize adjoint pullback}
&&\Big\downarrow\hbox{\scriptsize eliminate \(C_{\rm ch}\)}\\[2mm]
\hbox{chronological conormal word}
&\xrightarrow{\quad\hbox{\scriptsize Thomson compliance reduction}\quad}&
(bI+\widetilde K_{\mathcal N})Z_1 .
\end{array}
\end{equation}
Both paths produce the same final Schur row and the determinant
factorization \eqref{eq:appC-physical-block-factorization}.
\end{lemma}

\begin{proof}
At an unmarked contact, differentiation in
\((t_a,\theta_a)\) gives the block in
\eqref{eq:appC-root-block-dictionary}.  Chronological causality makes these
blocks triangular.  At a retained tree contact, the new relative-cluster
translation occurs for the first time with coefficient \(\sigma_e I_d\);
this gives \eqref{eq:appC-tree-block-dictionary}.  The new physical line
has no collision between \(u_c\) and \(d_2\), so its two endpoint rows in
\((X_c,V_c)\) give \eqref{eq:appC-chord-block-dictionary}.

Apply Theorem~\ref{thm:appC-physical-schur-jacobi} to this word.  Its
matrix \eqref{eq:appC-full-physical-J} contains every causal off-diagonal
block, while \eqref{eq:appC-explicit-Gaussian-matrices} performs the exact
physical elimination.  The theorem identifies that quotient with the
Jacobi Hessian quotient and gives
\eqref{eq:appC-final-block-dictionary}.  The diagonal of $A_W$ is precisely
$A_{\rm root},B_F,C_{\rm ch}$, so its determinant is their product.
This proves the diagram and \eqref{eq:appC-physical-block-factorization}
without suppressing any causal derivative.
\end{proof}

\paragraph{A fully expanded four-particle test word.}
Put
\[
 \tau_0=t_{13}<\tau_1=t_{12}<\tau_2=t_{34}<\tau_3=t_{23}<s
 <\tau_4=t_{14}<\tau_5=t_{24}
\]
and take
\[
 r=13,\qquad u_1=12,\qquad u_2=34,\qquad u_c=23,
 \qquad d_1=14,\qquad d_2=24.
\]
Thus \(F=\{12,34\}\), the free chord is the particle--\(2\) segment from
\(23\) to \(24\), and \(13\) is the single unmarked root.  We now give a
finite matrix algorithm for this word; in particular, there is no anonymous
off-diagonal block.

Let
\[
 E_{4,d}=\{(q_1,q_2,q_3,q_4)\in(\mathbb R^d)^4:
             q_1+q_2+q_3+q_4=0\},
 \qquad S_iq=q_i,
 \qquad B_{ij}=2^{-1/2}(S_i-S_j).
\]
Write \(\Pi_{E_{4,d}}\) for the kinetic-orthogonal projection onto
\(E_{4,d}\), a space of dimension $3d$.  The subscript $4$ counts physical
particles, not the ambient spatial dimension.
For \(a=ij\) let
\[
 \mathscr R_{ij}=I_{E_{4,d}}-2B_{ij}^*P_{\omega_{ij}}B_{ij}.
\]
Enumerate the event fields \(q_{a^-},q_{a^+}\in E_{4,d}\) in the displayed
chronological order and impose, literally,
\begin{equation}\label{eq:appC-four-particle-junction-list}
 q_{a^+}=\mathscr R_aq_{a^-},
 \qquad a\in\{13,12,34,23,14\}.
\end{equation}
For every maximal physical-line segment \(e=(i;a^+,b^-)\), other than the
particle--\(2\) segment \(\gamma=(2;23^+,24^-)\), set
\[
 A_e q:=S_iq_{b^-}-S_iq_{a^+}\in\mathbb R^d,
 \qquad \ell_e=t_b-t_a.
\]
After the two contact equations defining the eliminated particle--\(2\)
state have been used, its exact homogeneous drop row is
\begin{equation}\label{eq:appC-four-particle-chord-row}
 A_\gamma q=S_4q_{24^-}-S_3q_{23^-},
 \qquad \ell_\gamma=\tau_5-\tau_3.
\end{equation}
Indeed the two affine endpoint equations are
\(X_2+(\tau_3-s)V_2=X_3(\tau_3)+\eps\omega_{23}\) and
\(X_2+(\tau_5-s)V_2=X_4(\tau_5)+\eps\omega_{24}\);
their homogeneous difference is exactly
\eqref{eq:appC-four-particle-chord-row}.  \(S_3,S_4\) are the
unnormalised particle selectors.  Therefore the coefficient of
\(|A_\gamma q|^2\) is the physical weight \(\ell_\gamma^{-1}\), with no
hidden factor of two.

At the root \(13\), the shunt on \(E_{4,d}\) is the completely specified block
\begin{equation}\label{eq:appC-four-particle-root-shunt}
 H_{13}=B_{13}^*\widehat H_{13}B_{13},
 \qquad
 \widehat H_{13}
 =2|B_{13}V(\tau_0^- )\cdot\omega_{13}|\,
   \widehat{\mathcal V}_{13}^*\widehat{\mathcal K}_{13}
   \widehat{\mathcal V}_{13}\ge0.
\end{equation}
Let \(Q_{13}q=q_{13^-}\), stack the rows \(A_e\) and \(A_\gamma\) in
\(\mathbb A_4\), and put
\begin{equation}\label{eq:appC-four-particle-Hessian}
 \Lambda_4=\operatorname{diag}
   \bigl((\ell_eI_d)_{e\ne\gamma},\ell_\gamma I_d\bigr),
 \qquad
 L_4=\mathbb A_4^*\Lambda_4^{-1}\mathbb A_4
       +Q_{13}^*H_{13}Q_{13}.
\end{equation}
Equations \eqref{eq:appC-four-particle-junction-list}--
\eqref{eq:appC-four-particle-Hessian} are an entry-by-entry construction of
the factor incidence and its Hessian: each row contains exactly two selector
blocks, after the displayed reflection substitutions are made.

The forest coordinate matrix can be written without a
graphical convention.  With particle \(2\) grounded in the component
\(\{1,2\}\) and particle \(4\) grounded in \(\{3,4\}\), let
\begin{equation}\label{eq:appC-four-particle-forest-matrix}
 R_F^{\rm raw}(h_{12},h_{34})
 =\Pi_{E_{4,d}}(\sqrt2\sigma_{12}h_{12},0,
             \sqrt2\sigma_{34}h_{34},0).
\end{equation}
Inserted in the local pre-collision gauges at \(12\) and \(34\), this gives
\[
 \binom{B_{12}}{B_{34}}R_F^{\rm raw}
 =T_F:=\operatorname{diag}(\sigma_{12}I_d,\sigma_{34}I_d).
\]
This identity is the diagonal \(h\)-part of the actual physical forest row
\(\mathbb U_F=(D U_{12},D U_{34})\), after the root \(13\) has been
substituted; its causal root dependence occurs in the other columns, not in
the newly introduced cluster-translation diagonal.

To form the actual internal and boundary injections, first order all endpoint
coordinates lexicographically by
\((13^-,13^+,12^-,12^+,34^-,34^+,23^-,23^+,14^-,14^+,24^-)\),
substitute the five junction equations
\eqref{eq:appC-four-particle-junction-list}, and delete the two grounded
particle-coordinate blocks.  Let $z\in E_{4,d}$ be the complete reduced
separator field immediately before the lower word, in its kinetic
orthonormal zero-centre-of-mass coordinates, and let $y$ contain all
remaining endpoint coordinates other than $h=(h_{12},h_{34})$.
To make this a fixed smooth chart, use the finitely many complement minors:
subdivide the regular cell by the lexicographically first nonzero complement
minor and use the corresponding endpoint rows for $y$ on each subcell.
There are only finitely many such measurable smooth subcells, and all
formulas below are identical on them.  Thus the retained boundary port is
the full $3d$-dimensional space $E_{4,d}$, not a $d$-dimensional landing
projection.  The chronological derivative recursion below computes unique
matrices \(T_I,T_B\) such that the \emph{already imposed physical forest
constraint} is
\begin{equation}\label{eq:appC-four-particle-physical-forest-row}
 \mathbb U_Fq=T_Fh+T_Iy+T_Bz=0.
\end{equation}
Use the displayed inverse of \(T_F\) to substitute
\[
 h=-T_F^{-1}(T_Iy+T_Bz).
\]
Ordinary left-to-right substitution in the stated endpoint ordering now
gives a unique full-column-rank matrix
\[
 q=R_Iy+R_Bz.
\]
This is the promised executable construction of $R_I,R_B$; in particular,
the forest constraint is part of their definition rather than an additional
kernel assumption.  Here $R_B:E_{4,d}\to\mathscr X_{4,d}$ has full column rank.
The four stiffness blocks are
\begin{equation}\label{eq:appC-four-particle-Kron-blocks}
 L_{II}=R_I^*L_4R_I,\quad L_{IB}=R_I^*L_4R_B,
 \quad L_{BI}=R_B^*L_4R_I,\quad L_{BB}=R_B^*L_4R_B.
\end{equation}
If \(L_{II}y=0\), every segment drop, the chord drop, and the root quadratic
form vanish, while \eqref{eq:appC-four-particle-physical-forest-row} holds
because it was used to construct \(R_I\).  Starting at the two declared
grounds, the physical forest rows \(12,34\), then the junction \(23\), and
finally the segment rows to \(14,24\), force \(R_Iy=0\); full column rank
gives \(y=0\).  Hence \(L_{II}>0\), and the
boundary stiffness and compliance are the explicit matrices
\begin{equation}\label{eq:appC-four-particle-boundary-Schur}
 S_B^{(4)}=L_{BB}-L_{BI}L_{II}^{-1}L_{IB}>0,
 \qquad K_4=(S_B^{(4)})^{-1}
\end{equation}
on the grounded full boundary space $E_{4,d}$.
Indeed, if a nonzero boundary datum \(z\) had zero minimized energy, all
factor drops would vanish.  With
\eqref{eq:appC-four-particle-physical-forest-row} and the two grounds already
imposed, the same explicit propagation \(12,34,23,14,24\) used above would
force \(R_Bz=0\), contradicting the full column rank of \(R_B\).  Thus this
strict inequality is not inferred from \(L_{II}>0\) alone.

We finally display the physical derivative and its elimination algorithm.
Use the row vector
\[
 \mathbf C_4=(C_{13},U_{12},U_{34},U_{23},R_{24},R_{14})
\]
and the column vector
\[
 \mathbf p_4=(y_{13},h_{12},h_{34},X_2,V_2,v_*),
 \qquad y_{13}=(t_{13},\theta_{13}).
\]
Combining the middle two chord rows, the complete \(5\)-by-\(5\) block
matrix is
\begingroup\scriptsize
\begin{equation}\label{eq:appC-four-particle-matrix}
 J_4=D_{\mathbf p_4}\mathbf C_4=
 \begin{pmatrix}
 A_{13}&D_{h_{12}}C_{13}&D_{h_{34}}C_{13}
    &D_{(X_2,V_2)}C_{13}&D_{v_*}C_{13}\\
 D_{y_{13}}U_{12}&D_{h_{12}}U_{12}&D_{h_{34}}U_{12}
    &D_{(X_2,V_2)}U_{12}&D_{v_*}U_{12}\\
 D_{y_{13}}U_{34}&D_{h_{12}}U_{34}&D_{h_{34}}U_{34}
    &D_{(X_2,V_2)}U_{34}&D_{v_*}U_{34}\\
 D_{y_{13}}(U_{23},R_{24})&D_{h_{12}}(U_{23},R_{24})
    &D_{h_{34}}(U_{23},R_{24})&C_{\rm ch}
    &D_{v_*}(U_{23},R_{24})\\
 D_{y_{13}}R_{14}&D_{h_{12}}R_{14}&D_{h_{34}}R_{14}
    &D_{(X_2,V_2)}R_{14}&D_{v_*}R_{14}
 \end{pmatrix},
\end{equation}
\endgroup
where
\[
 A_{13}=[\,g_{13}\ \ -\eps E_{\omega_{13}}\,],
 \qquad
 C_{\rm ch}=\begin{pmatrix}I_d&(\tau_3-s)I_d\\
                             I_d&(\tau_5-s)I_d\end{pmatrix}.
\]
Every derivative in \eqref{eq:appC-four-particle-matrix} is computed by the
following literal chronological recursion.  Between event times \(t<t'\),
\[
 D_\beta X_i(t')=D_\beta X_i(t)+(t'-t)D_\beta V_i(t),
 \qquad D_\beta V_i(t')=D_\beta V_i(t).
\]
At a collision \(a=ij\), positions are continuous and, writing
\(g=V_i^--V_j^-\),
\begin{align*}
 D_\beta V_i^+
 &=D_\beta V_i^--[D_\beta g\cdot\omega_a]\,\omega_a
   -[g\cdot D_\beta\omega_a]\,\omega_a
   -[g\cdot\omega_a]D_\beta\omega_a,\\
 D_\beta V_j^+
 &=D_\beta V_j^-+[D_\beta g\cdot\omega_a]\,\omega_a
   +[g\cdot D_\beta\omega_a]\,\omega_a
   +[g\cdot\omega_a]D_\beta\omega_a,
\end{align*}
and \(D_\beta V_k^+=D_\beta V_k^-\) for \(k\notin\{i,j\}\).  Finally
\(D_\beta C_{ij}=D_\beta X_i^--D_\beta X_j^-
-\eps D_\beta\omega_{ij}\).  These three lines determine every named block
of \eqref{eq:appC-four-particle-matrix}, including all interspersed-root
off-diagonal terms.

For a numerical or symbolic check, set \(J^{(0)}=J_4\).  For \(m=1,\ldots,4\)
let \(P_m=J^{(m-1)}_{mm}\) and, for \(i,j>m\), define
\begin{equation}\label{eq:appC-four-particle-Schur-algorithm}
 J^{(m)}_{ij}=J^{(m-1)}_{ij}
 -J^{(m-1)}_{im}P_m^{-1}J^{(m-1)}_{mj}.
\end{equation}
Chronological cluster coordinates give, exactly,
\[
 P_1=A_{13},\qquad P_2=\sigma_{12}I_d,\qquad
 P_3=\sigma_{34}I_d,\qquad P_4=C_{\rm ch},\qquad
 S_{14}=J^{(4)}_{55}.
\]
Thus all causal blocks in the old schematic triangular display are retained
by the recurrence \eqref{eq:appC-four-particle-Schur-algorithm}.

We next give a determinant/exterior consistency calculation which is independent of
Theorem~\ref{thm:appC-physical-schur-jacobi} and of
Proposition~\ref{prop:appC-positive-network-closure}.  Before choosing the
$d$ columns \(v_*\), use the exact physical landing conormal
\[
 \widehat Z_{14}=(S_1-S_4)^*:\mathbb R^d\longrightarrow E_{4,d},
 \qquad
 \widehat Z_{14}^*\widehat Z_{14}=2I_d,
 \qquad
 \operatorname{Vol}_d(\widehat Z_{14})=2^{d/2}.
\]
Pull the \(d_1=14\) conormal backward through its collision-free lower
interval to the separator.  There it is exactly the one-speed frame
\((\widehat Z_{14},b\widehat Z_{14})\), \(b=\tau_4-s\).  Do not solve the
particle--\(2\) endpoint equations for \(V_2\) in the frame comparison.  Put
\[
 a=\tau_3-s,\qquad c=\tau_5-s,\qquad
 \ell_{\rm ch}=c-a=\tau_5-\tau_3.
\]
The fourth diagonal pivot and its inverse are, entry by entry,
\begin{equation}\label{eq:appC-four-particle-chord-inverse}
 C_{\rm ch}=\begin{pmatrix}I_d&aI_d\\ I_d&cI_d\end{pmatrix},
 \qquad
 C_{\rm ch}^{-1}=\frac1{\ell_{\rm ch}}
 \begin{pmatrix}cI_d&-aI_d\\-I_d&I_d\end{pmatrix},
 \qquad |\det C_{\rm ch}|=\ell_{\rm ch}^d.
\end{equation}
Let
\[
 \Pi_V=(\,0\ \ I_d\,),\qquad
 \mathcal D=(\,-I_d\ \ I_d\,),
\]
and let \(\mathcal T_{35}\) be the exact two-endpoint trace from the current
four-particle reduced frontier to the endpoint pair, including the full
pair-collision transports.  Thus
\(\mathcal A_{2,35}:=\mathcal D\mathcal T_{35}\) is precisely the transported
selector difference in \eqref{eq:appC-four-particle-chord-row}.  The inverse
above gives the elementary identity
\begin{equation}\label{eq:appC-four-particle-velocity-response}
 \Pi_VC_{\rm ch}^{-1}\mathcal T_{35}
 =\ell_{\rm ch}^{-1}\mathcal D\mathcal T_{35}
 =\ell_{\rm ch}^{-1}\mathcal A_{2,35}.
\end{equation}

This also identifies the fourth Schur contraction, including its sign and
its dimensions.  After the first three causal substitutions let
\(\mathcal Z_{14}^{(3)}:\mathbb R^d\to E_{4,d}\) be the exact current landing
frame of the retained (14)-row, including all earlier complete-pair
transports.  The endpoint residual is
\(C_{\rm ch}(X_2,V_2)-\mathcal T_{35}W\), whereas the endpoint velocity
response reaches the physical landing row only after the pullback
\((\mathcal Z_{14}^{(3)})^*\).  Consequently
\begin{equation}\label{eq:appC-four-particle-fourth-contraction}
 J^{(3)}_{45}=-\mathcal T_{35},\qquad
 J^{(3)}_{54}=(\mathcal Z_{14}^{(3)})^*
    \mathcal T_{35}^*\mathcal D^*\Pi_V.
\end{equation}
Here $J^{(3)}_{45}$ is $2d$-by-$3d$,
$J^{(3)}_{54}$ is $d$-by-$2d$, and their contraction is the required
$d$-by-$3d$ landing row.  Define separately the $3d$-by-$3d$ KKT
Hessian
\begin{equation}\label{eq:appC-four-particle-KKT-chord-Hessian}
 H_{2,35}:=\mathcal T_{35}^*\mathcal D^*\Pi_V
 C_{\rm ch}^{-1}\mathcal T_{35}
 =\ell_{\rm ch}^{-1}\mathcal A_{2,35}^*\mathcal A_{2,35}\ge0.
\end{equation}
The term subtracted in
\eqref{eq:appC-four-particle-Schur-algorithm} is therefore
\begin{align}
 -J^{(3)}_{54}C_{\rm ch}^{-1}J^{(3)}_{45}
 &=(\mathcal Z_{14}^{(3)})^*\mathcal T_{35}^*\mathcal D^*\Pi_V
        C_{\rm ch}^{-1}\mathcal T_{35}\notag\\
 &=(\mathcal Z_{14}^{(3)})^*H_{2,35}.
 \label{eq:appC-four-particle-Schur-is-shunt}
\end{align}
Thus direct differentiation gives the positive KKT current response
\eqref{eq:appC-four-particle-KKT-chord-Hessian}, while the physical Jacobian
receives its landing-frame pullback.  The two objects have different
dimensions and are not identified.  The occurrence of
\(\ell_{\rm ch}^{-1}\) is confined to the positive current action.  The
independent Gaussian determinant contribution remains
\(|\det C_{\rm ch}|=\ell_{\rm ch}^d\); neither factor is silently absorbed
into the other.
We next obtain the exterior consistency identity directly in current space, without
using the general network theorem.  Retain every edge current, including
the chord current, in the finite four-particle Thomson functional obtained
from \eqref{eq:appC-four-particle-Hessian}.  After internal current balance
and one translation ground are imposed, minimization defines a symmetric
positive semidefinite compliance $K_4\ge0$ on the retained boundary
current.  The exact discrete Green identity gives
\begin{equation}\label{eq:appC-four-particle-direct-current-identity}
 W_{\rm full}:=(S_{14}^{\rm full})^*
 =(bI_{E_{4,d}}+\widetilde K_4)\widehat Z_{14},
 \qquad
 \widetilde K_4:=Q_{4,+}K_4Q_{4,-}
 =Q_{4,-}^*K_4Q_{4,-}\ge0.
\end{equation}
Here $K_4,Q_{4,-},Q_{4,+}:E_{4,d}\to E_{4,d}$ and
$Q_{4,+}=Q_{4,-}^*$ by the same exact dual-pairing convention as
\eqref{eq:appC-exact-dual-boundary-gauges}; thus every factor in the displayed
identity acts in the same $3d$-dimensional full boundary space.  The
square gauges are obtained by writing out the five complete-pair orthogonal
junctions, an orthogonal completion of the normalized $14$-incidence, and
the two-edge forest QR matrix.  Consequently
\[
 s_{\min}(Q_{4,-})\ge c_{4,-}>0,
 \qquad s_{\min}(Q_{4,+})\ge c_{4,+}>0,
 \qquad c_{\rm reg}:=\min\{1,c_{4,-}c_{4,+}\}>0,
\]
and these constants contain no chord length.  Since $\widetilde K_4\ge0$,
the exterior singular-value inequality and
\(\operatorname{Vol}_d(\widehat Z_{14})=2^{d/2}\) give
\begin{equation}\label{eq:appC-four-particle-direct-exterior-check}
 \operatorname{Vol}_d(W_{\rm full})
 \ge2^{d/2}b^d
 \ge2^{d/2}c_{\rm reg}^db^d.
\end{equation}
This is a direct calculation in the explicit four-particle KKT system; it
does not invoke Theorem~\ref{thm:appC-physical-schur-jacobi} or
Proposition~\ref{prop:appC-positive-network-closure}.

The four explicit Gaussian steps
\eqref{eq:appC-four-particle-Schur-algorithm}, applied to the same KKT
system, leave precisely the full separator-velocity row
\((S_{14}^{\rm full})^*=W_{\rm full}\).  This is the four-particle instance
of the discrete Green identity, derived here from the displayed finite
matrices.  The only changes between the physical
orthonormal frame and the raw cluster/selector coordinates are precisely the
junction/forest gauges just displayed.  The chord endpoint pivot
is not among these coordinate changes: its determinant
\(|\det C_{\rm ch}|=\ell_{\rm ch}^d\) is retained as the fourth diagonal
Gaussian pivot, whereas its effect on the landing conormal is already the
positive current action \eqref{eq:appC-four-particle-KKT-chord-Hessian}.
Hence neither
\(C_{\rm ch}^{-1}\) nor \(\ell_{\rm ch}^{-1}\) enters the frame comparison.
The displayed least-singular-value bounds give
\(\Jac_d(S_{14}^{\rm full})\ge2^{d/2}c_{\rm reg}^db^d\).
Choose \(v_*\) to be the lexicographically first maximal $d$-column minor
of this \(d\)-by-\(3d\) row.  Cauchy--Binet then gives the independent
square-minor bound
\begin{equation}\label{eq:appC-four-particle-minor-check}
 |\det S_{14}|
 \ge \frac{2^{d/2}c_{\rm reg}^d}{\sqrt{\binom{3d}{d}}}\,(\tau_4-s)^d.
\end{equation}
Together with the already displayed diagonal pivots,
\begin{equation}\label{eq:appC-four-particle-full-determinant-check}
 |\det J_4|
 \ge
 \eps^{d-1}|g_{13}\cdot\omega_{13}|\,
 (\tau_5-\tau_3)^d
 \frac{2^{d/2}c_{\rm reg}^d}{\sqrt{\binom{3d}{d}}}\,
 (\tau_4-s)^d .
\end{equation}
Thus this particular word checks the determinant and exterior conclusions
without feeding the general Schur--Jacobi theorem back into its own example.

The full compliance identity just proved is
\begin{equation}\label{eq:appC-four-particle-physical-block}
 (S_{14}^{\rm full})^*
 =(bI_{E_{4,d}}+\widetilde K_4)\widehat Z_{14},
 \qquad b=\tau_4-s.
\end{equation}
The determinant lower bound was obtained from this explicit current action,
not from the general theorem.  Consequently
\begin{align*}
 |\det A_{13}|&=\eps^{d-1}|g_{13}\cdot\omega_{13}|,\\
 |\det T_F|&=1,\\
 |\det C_{\rm ch}|&=(\tau_5-\tau_3)^d,\\
 |\det S_{14}|&\ge
 \dfrac{2^{d/2}c_{\rm reg}^d}{\sqrt{\binom{3d}{d}}}(\tau_4-s)^d.
\end{align*}
Changing a pair orientation signs the corresponding incidence row and does
not change any absolute determinant, Hessian, or exterior estimate.

\subsection{Application to the overlapping landing cell}

We finish by spelling out how the abstract objects correspond to the proof
of Proposition~\ref{prop:asymmetric-new-line-peeling}.  Let
$d_1=\{a,b\}$ and $d_2=\{b,c\}$ be the first two lower atoms, and let $u_c$
be the last upper atom on $c$.  Remove $u_c$ from the reverse-Kruskal upper
tree.  The remaining graph is a forest with two components.  Ground the
first at the separator coordinate of $c$ and the second at the chosen
reference particle.  The reduced incidence matrix is square and invertible
on all other relative translations.

The $(u_c,d_2)$ contact rows determine $(X_c,V_c)$ with determinant
$(t_{d_2}-t_{u_c})^d$.  By
Lemma~\ref{lem:causal-free-line-endpoints}, their two companion endpoint
fields do not use $(X_c,V_c)$.  Eliminating this state therefore inserts
exactly the chord term \eqref{eq:appC-affine-chord-energy}; no derivative of
the first landing has been discarded.

\begin{lemma}[Coverage certificate for the overlapping word]
\label{lem:appC-overlap-trace-certificate}
The regular overlapping word just constructed is trace-complete in the
sense of Definition~\ref{def:appC-trace-complete-word}.
\end{lemma}

\begin{proof}
We verify the four clauses of the definition, numbering the two separate
matrix requirements in its frontier clause as \textup{(TC3a)} and
\textup{(TC3b)}.
\begin{enumerate}[label=\textup{(\arabic*)}]
\item Removing $u_c$ from the reverse-Kruskal spanning tree produces exactly
two rooted components.  The separator coordinate of $c$ and the chosen
reference particle are retained raw particle coordinates, so their selector
rows are the two blocks of $\mathbb G$.  No quotient representative is being
used as a ground row.
\item Order the remaining Kruskal edges from the two roots towards their
leaves, before passing to the solved $(h,z)$ split.  In raw endpoint-particle
coordinates the contact at the first edge reaching a new cluster contains
that cluster translation for the first time, with derivative
$\sigma_eI_d$.  Its other endpoint coordinate has already been reached along
the rooted tree.  A complete-pair root substitution is an invertible change
of the two already present endpoint blocks and does not change this new
diagonal.  Thus the selected raw forest matrix is block lower triangular with
$T_e=\sigma_eI_d$ and has the form
\eqref{eq:appC-trace-merge-row}.  Only after this reconstruction is complete
do we change variables to $q=R_Iy+R_Bz$; in those later coordinates the same
row can contain a genuine $T_Bz$.  The trace proof therefore does not claim
$T_B=0$ and does not solve an equation having the unknown $z$ on its own
right-hand side.
\item For \textup{(TC3a)}, split the time-expanded word at every collision.
The only physical
state removed by the endpoint pivot $(u_c,d_2)$ is the state of line $c$ on
its collision-free interval.  Every other particle-line segment and its
factor row remains present.  Hence on every slab crossed by that interval
the selected rows are
\[
 (S_i(x_\alpha-\mathscr R_\alpha
        x_{\operatorname{par}(\alpha)}))_{i\ne c}
 =\Sigma_{\widehat c}
   (x_\alpha-\mathscr R_\alpha
        x_{\operatorname{par}(\alpha)}).
\]
Equation \eqref{eq:appC-missing-selector-inverse} is their inverse.  Outside
the chord interval all particle rows are present and one arbitrary selector
may be omitted.  This proves selector coverage without using
$A_\gamma$ as a full-rank row.
\item For \textup{(TC3b)}, at every collision, the two incident line
coordinates are changed by
the complete matrix $\mathscr R_a^{(2d)}$ and all other coordinates are
unchanged.  Its extension $\mathscr R_a$ to the reduced frontier is
orthogonal, so the inverse in item~\textup{(3)} is precisely
$\mathscr R_a^*\Lambda_{\widehat c}$ and has no grazing or chord-length
factor.
\item The first-landing row $d_1$ is retained only after the lower
chronological reconstruction has reached its separator frontier.  Its
one-port boundary field is the complete reached reduced frontier $E$, not the
$d$-dimensional landing projection.  Thus $\operatorname{tr}_B$ is the
identity in kinetic orthonormal coordinates (and the fixed raw-to-orthonormal
map in raw coordinates), so
$\|\operatorname{tr}_B\|\le C(P_0)$.  The normalized incidence $Z_1$ is only
the $d$-dimensional current frame inserted into that full port.
\end{enumerate}
These five checks are respectively \textup{(TC1)}, \textup{(TC2)},
\textup{(TC3a)}, \textup{(TC3b)} and \textup{(TC4)}.  They provide the row
selector $\mathcal P_W$ and triangular
matrix $\mathcal T_W$ used in
\eqref{eq:appC-trace-left-inverse-formula}; no additional topological
assumption is being inserted.
\end{proof}

The remaining first-landing conormal at the separator is
\begin{equation}\label{eq:appC-landing-one-speed}
 (Z_1,W_1)=(Z_1,bZ_1),
 \qquad b=t_{d_1}-s>0.
\end{equation}
The map $Z_1$ is the normalized incidence covector of the pair $\{a,b\}$,
tensored with $I_d$, and hence
\begin{equation}\label{eq:appC-landing-incidence-volume}
 \operatorname{Vol}_d(Z_1)\ge c(d,P_0)>0.
\end{equation}
All earlier upper free intervals are generators
\eqref{eq:appC-free-generator}; all fixed-normal resets are generators
\eqref{eq:appC-collision-generator}; all incoming-root eliminations supply
nonnegative curvature shunts; and the restored chord is a positive network
square.  Proposition~\ref{prop:appC-positive-network-closure} and
Corollary~\ref{cor:appC-one-speed-word} therefore give
\begin{equation}\label{eq:appC-final-landing-volume}
 \operatorname{Vol}_d(W_{1,\rm Sch})
 \ge c(d,P_0)b^d.
\end{equation}
Since the separator is chosen with
$b\ge\varepsilon_*/4$, the right-hand side is at least
$c(d,P_0)\varepsilon_*^d$.

There are two graph possibilities, and neither changes the proof.  If the
chord endpoints lie in different forest components, the chord reconnects
them.  If they lie in the same component, it creates a positive co-tree
term while the other grounded component is eliminated independently.  In
the first case the grounded Dirichlet matrix is positive because the chord
connects the two grounds through the boundary; in the second it is block
positive because each component already has its own ground.  Thus no
unstated topological alternative remains.

Exact grazing, coincident collision times, and rank loss of a normal chart
are excluded from the regular cells.  On those cells all matrices above
depend smoothly on the original variables.  Their exceptional complements
have collision-flux measure zero, and the nonnegative estimates extend to
the full integral by monotone exhaustion.  Consequently
\eqref{eq:appC-final-landing-volume} is an analytic coarea bound, not a
generic-transversality assertion and not a conclusion drawn from numerical
sampling.

\section{Finite-Fibre Coarea and Arbitrary Joint Kernels}
\label{app:coarea}

The measure-theoretic part of Theorem~\ref{thm:two-landing} is developed here.
We construct the cellwise measure, state the precise submersion result for
products of collision constraints, and verify that every substitution acts
on the same nonnegative kernel.  The construction first identifies a smooth
submersion on each regular cell and then uses the coarea measure of its zero
set; it does not define a product of delta distributions by formal iteration.
We use the coarea
theorem in the form of \cite[Section~3.2]{Federer1969}.  The only
real-algebraic inputs are the elementary affine B\'ezout bound and
semialgebraic selection; see also \cite{BasuPollackRoy2006} for a systematic
treatment.

\subsection{Regular half-open cells}

Fix a packet \(H\) with \(P\) physical particle lines and a strict total
order of its atoms.  A cell index \(\kappa\) contains the following finite
data:
\begin{enumerate}[label=\textup{(\roman*)}]
\item the strict atom order and all C-atom incoming/outgoing slot labels;
\item one sign choice for each scalar collision flux;
\item dyadic cells for every relative speed and for every velocity retained
as a free variable;
\item one normal chart from a fixed finite atlas of \(\Sph^{d-1}\) at each atom;
\item the torus lift attached to every transported contact residual;
\item a rooted upper spanning tree and, when needed, the lexicographically
first maximal coarea minor;
\item half-open tie-breaking inequalities for all of the preceding choices.
\end{enumerate}
The half-open convention makes ordinary chart/dyadic cells disjoint up to
their null geometric boundaries.  Maximal-minor ties are different: they
may have positive measure and are not deleted.  If the squared minors are
$m_1,\ldots,m_N$ in lexicographic order, use the Borel pieces
\[
 E_i=\{m_i>m_j\ (j<i),\ m_i\ge m_j\ (j>i)\}.
\]
They are disjoint and cover every tie point exactly once.

On the velocity cutoff \( |v|\le V_\eps=\cL_\eps^{C_0}\) and a time interval
of length at most \(\cL_\eps^{C_0}\), every displacement has size at most
\(\cL_\eps^{C}\).  Hence the number of admissible torus lifts for one
contact is at most \(\cL_\eps^C\).  The packet cap \( |H|\le P_0\) in the
application therefore gives at most \(\cL_\eps^{C(P_0)}\) lift words.
Normal charts, sign choices and rooted-tree choices cost only \(C(P_0)\).
Without substituting the packet cap, the uniform cell count is
\begin{equation}\label{eq:appD-cell-count}
 \#\{\kappa\}\le C^{|H|}\cL_\eps^{C|H|}.
\end{equation}
This is why Theorem~\ref{thm:two-landing} states
\(\cL_\eps^{C|H|}\) at cell level and writes \(\cL_\eps^C\) only after
\( |H|\le P_{\rm pkt}(d)\) has been used.

Let \(\mathfrak X_{H,\kappa}\) be the product of the raw position, velocity,
time and normal-chart variables on the cell, together with the fixed-end
and complement variables.  We remove the following sets:
\begin{enumerate}[label=\textup{(\alph*)}]
\item exact grazing, where a collision flux vanishes;
\item equality of two atom times;
\item a simultaneous contact of two distinct unordered particle pairs;
\item rank loss of a pivot block that the cell declares nonzero.
\end{enumerate}
Normal, velocity, lift and maximal-minor boundaries are not removed from a
conditional fibre.  They are assigned to exactly one Borel half-open cell;
only geometric normal-chart boundaries, when used for a surface atlas, may
also be chosen surface-null.
The remaining set $\mathfrak X_{H,\kappa}^{\rm reg}$, intersected with one
$E_i$, is Borel and is exhausted by compact sets on which every declared
inverse Jacobian is bounded.  It need not be relatively open along a minor
tie.  At each point the selected minor is nonzero on an ambient open
neighborhood; local coarea applies there and its Radon measure is then
restricted back to $E_i$.  The common-Radon-measure and conditional
fixed-boundary zero-measure verification for (a)--(c), and the exhaustive
routing of (d), are proved below in
Lemma~\ref{lem:appD-conditional-null-strata}; they are not inferred only
from an unconditional almost-everywhere flow statement.

\subsection{One collision as a vector submersion}

The first elementary identity is distributional but unambiguous.

\begin{lemma}[Sphere disintegration in fixed dimension]
\label{lem:appD-sphere-disintegration}
For every nonnegative Borel function \(F\) on \(\R^d\),
\begin{equation}\label{eq:appD-sphere-disintegration}
 \int_{\R^d}\delta(|R|-\eps)F(R)\dd R
 =
 \eps^{d-1}\int_{\Sph^{d-1}}F(\eps\omega)\dd\omega.
\end{equation}
Equivalently,
\[
 \delta(|R|-\eps)F(R)
 =\eps^{d-1}\int_{\Sph^{d-1}}\delta^{(d)}(R-\eps\omega)
 F(\eps\omega)\dd\omega
\]
as Radon measures in \(R\).
\end{lemma}

\begin{proof}
Use polar coordinates \(R=r\omega\).  Since
\(\dd R=r^{d-1}\dd r\,\dd\omega\), integration against
\(\delta(r-\eps)\) gives \eqref{eq:appD-sphere-disintegration}.  The
equivalent vector-delta identity follows by testing against a compactly
supported continuous function and applying Tonelli.
\end{proof}

Let $\theta\mapsto\omega(\theta)$ be a regular normal chart, put
\[
 G_\omega=(D_\theta\omega)^*D_\theta\omega,
 \qquad J_\omega=(\det G_\omega)^{1/2},
 \qquad E_\omega=D_\theta\omega\,G_\omega^{-1/2}.
\]
Thus $E_\omega:\R^{d-1}\to\omega^\perp$ is an orthonormal moving frame, while
$D_\theta\omega$---not $E_\omega$---is the actual chart differential, and
$\dd\omega=J_\omega\dd\theta$.  A freely transported relative contact residual has
the form
\begin{equation}\label{eq:appD-one-root-map}
 \Psi(t,\theta)=x+tg-m-\eps\omega(\theta).
\end{equation}

\begin{lemma}[Exact incoming-root block]
\label{lem:appD-exact-root-block}
At a regular root of \eqref{eq:appD-one-root-map},
\begin{equation}\label{eq:appD-exact-root-determinant}
 \left|\det D_{(t,\theta)}\Psi\right|
 =\eps^{d-1}|g\cdot\omega|J_\omega(\theta).
\end{equation}
On an incoming-sign cell the root is unique for a fixed lift.  Consequently,
for every nonnegative \(F\),
\begin{equation}\label{eq:appD-root-cancellation}
 \int \eps^{d-1}|g\cdot\omega|\,
 \delta^{(d)}(\Psi(t,\theta))F(t,\theta)
 J_\omega(\theta)\dd t\,\dd\theta
 =
 F(t_{\rm in},\theta_{\rm in})
\end{equation}
when the declared incoming root belongs to the cell, and the integral is
zero otherwise.
\end{lemma}

\begin{proof}
The derivative columns are $g$ and $-\eps D_\theta\omega$.  Factoring
$D_\theta\omega=E_\omega G_\omega^{1/2}$ and resolving $g$ into its normal
and tangential parts gives
\[
 |\det[g\ \ -\eps D_\theta\omega]|
 =\eps^{d-1}|g\cdot\omega|J_\omega.
\]
Uniqueness is Lemma~\ref{lem:one-incoming-root}.  The ordinary change of
variables theorem on the normal chart then proves
\eqref{eq:appD-root-cancellation}; the identical $J_\omega$ in the surface
measure and in the determinant cancels exactly.  To globalize, choose a
finite atlas $U_1,\ldots,U_N$ and its half-open Borel disjointification
$V_1=U_1$, $V_j=U_j\setminus\bigcup_{i<j}U_i$.  Chart boundaries have
surface measure zero, so every regular normal is counted exactly once and
no partition-of-unity multiplicity is introduced.
\end{proof}

Thus an unmarked collision is neutral: its sphere factor, collision flux
and time--normal inverse Jacobian cancel exactly.  This is not a bound by an
inverse grazing factor.  Exact grazing is simply absent from the collision
measure.

\subsection{Coarea, finite fibres, and branch labels}

We record the form of coarea used throughout the packet proof.

\begin{lemma}[Square block substitution]
\label{lem:appD-square-substitution}
Let \(U\subset\R^n\times\R^m\) be open, write variables as \((p,z)\), and
let \(G:U\to\R^n\) be \(C^1\).  Suppose that on a half-open cell every zero
of \(G(\,\cdot\,,z)\) is simple, there is at most one such zero
\(p=\pi(z)\), and
\[
 |\det D_pG|\ge\lambda>0.
\]
Then for every nonnegative Borel function \(Q\),
\begin{equation}\label{eq:appD-square-substitution}
 \int_U Q(p,z)\delta^{(n)}(G(p,z))\dd p\,\dd z
 \le
 \lambda^{-1}\int_{\pi_zU} Q(\pi(z),z)\dd z .
\end{equation}
The same formula with the exact determinant in the integrand is an
equality.
\end{lemma}

\begin{proof}
First take \(Q\) bounded and compactly supported inside a regular compact
exhaustion.  The implicit-function theorem and the change of variables
\(p\mapsto G(p,z)\) give the formula for each fixed \(z\); Tonelli then
integrates in \(z\).  Approximate an arbitrary nonnegative Borel function
from below by simple functions and exhaust the regular cell.  Monotone
convergence proves the stated form.  No product of independently defined
delta distributions occurs.
\end{proof}

For a rectangular map \(G:U\to\R^k\), \(k\le n\), its normal Jacobian is
\[
 \Jac_k(D_pG)=\sqrt{\det(D_pG(D_pG)^*)}.
\]
Cauchy--Binet gives
\begin{equation}\label{eq:appD-Cauchy-Binet}
 \Jac_k(D_pG)^2
 =\sum_{I\subset\{1,\ldots,n\},\,|I|=k}
  |\det D_{p_I}G|^2.
\end{equation}

\begin{lemma}[Finite-fibre square substitution]
\label{lem:appD-finite-fibre-substitution}
Let \(U\subset\R^n\times\R^m\) be open, write the variables as
\((p,z)\), and let \(G:U\to\R^n\) be \(C^1\).  Suppose that every zero of
\(G(\,\cdot\,,z)\) is simple, that

\[
 \#\{p:(p,z)\in U,\ G(p,z)=0\}\le M,
 \qquad |\det D_pG|\ge\lambda>0
\]

on the zero set.  Then the zeros admit Borel enumerations
\(p=\pi_j(z)\), \(1\le j\le M\), with an absent branch assigned the value
\(\dagger\), and for every nonnegative Borel \(Q\),
\begin{equation}\label{eq:appD-finite-fibre-substitution}
 \int_U Q(p,z)\delta^{(n)}(G(p,z))\dd p\dd z
 \le \lambda^{-1}\sum_{j=1}^{M}
 \int Q(\pi_j(z),z)\1_{\{\pi_j(z)\ne\dagger\}}\dd z.
\end{equation}
With \(|\det D_pG|\) in the integrand, the corresponding formula is an
equality.
\end{lemma}

\begin{proof}
On a compact regular exhaustion the implicit-function theorem writes the
zero set locally as finitely many graphs.  The area formula, followed by
Tonelli, gives for almost every \(z\) the sum over all roots in that fibre,
with weight \(|\det D_pG|^{-1}\).  A finite Borel fibre can be enumerated
by repeatedly taking the lexicographically least remaining point; the graph
of each selection is Borel.  This gives at most \(M\) branch slots.  The
determinant lower bound proves the inequality.  Monotone exhaustion and
monotone convergence remove compact support.  This lemma does
not infer \(M=1\) from a nonzero derivative.
\end{proof}

\begin{lemma}[Deterministic maximal-minor selection]
\label{lem:appD-maximal-minor-atlas}
Suppose \(\Jac_k(D_pG)\ge\lambda\) on the zero set.  Partition it by
declaring the lexicographically first minor attaining the maximum in
\eqref{eq:appD-Cauchy-Binet}.  The pieces are Borel and disjoint, and on
each piece the selected minor satisfies
\begin{equation}\label{eq:appD-maximal-minor-lower}
 |\det D_{p_I}G|
 \ge \binom nk^{-1/2}\lambda.
\end{equation}
Each selected minor is a local submersion.  This conclusion is local: it
does not assert that the selected coordinates have one zero in a complete
retained fibre.  Lemma~\ref{lem:appD-finite-fibre-substitution} may be
applied only after a separate fibre bound has been proved.
\end{lemma}

\begin{proof}
The minors are continuous.  The sets $E_i$ displayed at the start of this
section are therefore Borel, disjoint and exhaustive, including
positive-measure ties.  The largest square in a sum of
\(\binom nk\) nonnegative squares is at least the average, proving
\eqref{eq:appD-maximal-minor-lower}.  The implicit-function theorem supplies
local graph charts.  No global injectivity follows: for example
\(p\mapsto\sin p\) on \((0,4\pi)\) has nonzero derivative at each zero but
several zeros in one fibre.
\end{proof}

The packet proof uses maximal-minor selection in the coarse disjoint landing
chart, in the remaining coarse $d$-row overlap chart, and in the default FCT
chart.  In all three uses the missing global information is supplied by the
uniform algebraic fibre theorem below, not by the implicit-function theorem.

\subsection{Raw packet coordinates}

Choose a separator time \(s\) between the upper and lower sublayers.  On
every physical line retain its right trace \((X_p,V_p)\) at \(s\).  Quotient
the common translation by retaining one reference position.  The raw
separator space is
\begin{equation}\label{eq:appD-separator-space}
 \T^d\times\R^{d(P-1)}\times\R^{dP}.
\end{equation}
The transformation from any edge trace to
\eqref{eq:appD-separator-space} is an affine shear and has Jacobian one.

Given the separator states and the complete list of atom times and normals,
reconstruct the upper word backward and the lower word forward:
\begin{enumerate}[label=\textup{(\roman*)}]
\item between contacts use affine free transport;
\item at a C-atom use the complete elastic pair involution;
\item at an O-atom retain the transparent velocity pair;
\item use the fixed lift only in the position residual, never as a
continuous variable.
\end{enumerate}
Every elastic map is orthogonal on the $2d$-dimensional velocity pair and
has absolute determinant one.  Every transport-position delta is a
translation in the newly reconstructed position and also has determinant
one.  Hence all nontransversal edge states are uniquely reconstructed from
the separator states with no density loss.

Let \(\mathcal U_T\) be the vector of contact residuals at the
\(P-1\) reverse-Kruskal upper-tree atoms.  Use chronological
relative-cluster translations \(h\in\R^{d(P-1)}\) as their pivot variables.
Before any collision substitution, the derivative is block triangular:
\begin{equation}\label{eq:appD-tree-incidence-block}
 D_h\mathcal U_T=
 \begin{pmatrix}
  B_1\otimes I_d&0&\cdots&0\\
  *&B_2\otimes I_d&\cdots&0\\
  \vdots&\vdots&\ddots&\vdots\\
  *&*&\cdots&B_{P-1}\otimes I_d
 \end{pmatrix},
\end{equation}
where every scalar \(B_j\) is \(\pm1\) after the new cluster translation is
normalized.  Therefore
\begin{equation}\label{eq:appD-tree-incidence-determinant}
 |\det D_h\mathcal U_T|=1
\end{equation}
in raw cluster coordinates.  Passing to kinetic-orthonormal relative
particle coordinates changes this by a factor between \(c(P)\) and \(C(P)\).

Every unmarked upper contact is solved by the incoming-root block of
Lemma~\ref{lem:appD-exact-root-block}.  The reverse-Kruskal property says
that its two endpoints are already connected through later retained tree
contacts.  Consequently its causal root substitution does not change the
rank of \eqref{eq:appD-tree-incidence-block}; it only inserts the positive
curvature shear described in
Lemma~\ref{lem:rooted-tree-schur-lift}.

On the lower side, the first two atoms \(d_1<d_2\) have no lower atom before
\(d_1\) or between them.  If \(v_i\) is the separator velocity on the line
first used at \(d_i\), the residuals satisfy
\begin{equation}\label{eq:appD-lower-causal-form}
 \begin{split}
 R_1&=(t_{d_1}-s)Z_1v_1+F_1,\\
 R_2&=(t_{d_2}-s)Z_2v_2+F_2(v_1),
 \end{split}
\end{equation}
on disjoint landing edges.  The functions \(F_i\) may depend on every
retained complement variable and on the full upper history.  The essential
point is that \(F_1\) does not depend on \(v_2\).  This follows by forward
induction over the lower word: no event before \(d_1\) exists, and the only
event before \(d_2\) is \(d_1\).  Orthogonal slot changes merely replace
\(Z_i\) by normalized incidence isometries.  Hence
\begin{equation}\label{eq:appD-disjoint-lower-determinant}
 \left|\det D_{(v_1,v_2)}(R_1,R_2)\right|
 =|t_{d_1}-s|^d|t_{d_2}-s|^d.
\end{equation}

When the landing edges overlap, write them as
\(\{a,b\}\) and \(\{b,c\}\), with \(c\) new at \(d_2\).  The raw companion
independence in Lemma~\ref{lem:causal-free-line-endpoints} gives instead
\begin{equation}\label{eq:appD-overlap-first-block}
 \left|\det
 D_{(X_c,V_c)}(U_{u_c},R_{d_2})\right|
 =|t_{d_2}-t_{u_c}|^d.
\end{equation}
After this substitution, the remaining first-landing row is a one-speed
frame.  Section~\ref{app:positive-network} gives
\begin{equation}\label{eq:appD-overlap-second-block}
 \Jac_d(W_{1,\rm Sch})\ge c(d,P_0)|t_{d_1}-s|^d.
\end{equation}
Equations
\eqref{eq:appD-overlap-first-block}--\eqref{eq:appD-overlap-second-block}
give the two successive local submersions used on an overlapping cell.  The
first is globally affine in \((X_c,V_c)\); the second need not be globally
injective after the incoming-root variables have been eliminated.  We do
not assume that it is.

\subsection{Uniform global fibre multiplicity}

We now prove the global statement which is not supplied by a normal
Jacobian lower bound.  It is useful to keep the normal vectors extrinsic.
For every unmarked contact introduce
\((t_a,\omega_a)\in\R\times\R^d\)
and append the equation \(|\omega_a|^2-1=0\).  The $d$ contact equations
together with this sphere equation replace the $d$ local
time--normal-chart equations.  On a fixed normal chart the two systems
have the same regular roots.  More precisely,
\begin{equation}\label{eq:appD-extrinsic-sphere-measure}
 \int_{\Sph^{d-1}}F(\omega)\dd\omega
 =2\int_{\R^d}F(\omega)
       \delta(|\omega|^2-1)\dd\omega,
\end{equation}
and, at a contact with relative velocity \(g_a\),
\begin{equation}\label{eq:appD-extrinsic-root-determinant}
 \left|\det D_{(t_a,\omega_a)}
   \big(R_a-\eps\omega_a,|\omega_a|^2-1\big)\right|
 =2\eps^{d-1}|g_a\cdot\omega_a|.
\end{equation}
Indeed, expand the last row \(2\omega_a^*\) against the $d$ columns
\(-\eps I_d\).  Thus the factor $2$ in the surface measure cancels the
factor $2$ in the extrinsic root determinant, and polynomialization
changes neither the collision measure nor the root cancellation.

Let \(p_{H,\kappa}\) contain all continuous pivot variables in one of the
three packet charts: the chronological upper-tree translations; the
selected landing velocity coordinates (and, in the overlap chart, the
free-chord variables); and the pairs \((t_a,\omega_a)\) at every unmarked
contact.  Let \(b\) contain all remaining continuous variables and the
concrete fixed-boundary values.  Write
\begin{equation}\label{eq:appD-full-polynomial-system}
 \mathscr F_{H,\kappa}(p;b)=0
\end{equation}
for all selected vector contact residuals, all unmarked vector contact
residuals, and the extrinsic sphere equations.  It is a square system.

\begin{lemma}[Bounded degree of a fixed hard-sphere word]
\label{lem:appD-bounded-degree-word}
For \(|H|\le P_0\), every coordinate of
\(\mathscr F_{H,\kappa}\) is a real polynomial in \((p,b)\) of degree at
most
\begin{equation}\label{eq:appD-degree-bound}
 D_0(P_0):=4P_0+4.
\end{equation}
The number \(N_{H,\kappa}\) of pivot coordinates and equations satisfies
\(N_{H,\kappa}\le(2d+1)P_0+2d\).  These bounds are independent of
\(\eps\), the torus lifts, the dyadic cells and the fixed-boundary values.
\end{lemma}

\begin{proof}
At the separator the state coordinates have degree one.  Free transport
replaces \(x\) by \(x+(t'-t)v\), increasing total degree by at most one.
At a C-atom the complete elastic involution is
\[
 v_i^+=v_i^--[(v_i^--v_j^-)\cdot\omega_a]\omega_a,
 \qquad
 v_j^+=v_j^-+[(v_i^--v_j^-)\cdot\omega_a]\omega_a,
\]
so one collision increases total degree by at most two; an O-atom does not
increase it.  Induction through at most \(P_0\) atoms therefore bounds every
velocity by degree \(2P_0+1\) and every transported contact residual by
degree \(2P_0+2\).  The sphere equation has degree two, so
\eqref{eq:appD-degree-bound} is a safe common bound.  Torus lifts enter only
as constant translations.

Because the packet collision graph is connected, it contains a spanning
tree on its $P$ physical-line vertices.  Hence it has at least $P-1$
binary-event edges, and therefore
\begin{equation}\label{eq:appD-line-event-count}
 P\le |H|+1\le P_0+1.
\end{equation}
There are \(d(P-1)+2d=d(P+1)\) tree/landing scalar pivots in the disjoint
chart.  The overlap chart replaces $2d$ of them by \((X_c,V_c)\) and has
the same total.  Each unmarked time--normal block uses $d+1$ extrinsic
coordinates.  The deliberately loose estimate
\(d(P+1)+(d+1)|H|\le(2d+1)P_0+2d\) covers every chart, including the default FCT
choice.  The radial row is eliminated locally and does not enter this
Cartesian polynomial variable count.
\end{proof}

\mainstatementreference{Theorem~\ref{thm:appD-uniform-algebraic-multiplicity}}

\mainproofheading{Theorem~\ref{thm:appD-uniform-algebraic-multiplicity}}{supp:proof:thm:appD-uniform-algebraic-multiplicity}
\begin{proof}
Complexify the square polynomial system
\eqref{eq:appD-full-polynomial-system}.  At a regular real solution the
square pivot derivative is invertible: its block determinant is the product
of the declared upper-tree block, the selected landing minor, the incoming
root blocks and the extrinsic sphere-chart factors.  Hence that solution is
an isolated complex zero.  The affine B\'ezout inequality bounds the number
of isolated complex zeros, even when the full complex zero set has other
positive-dimensional components, by the product of the equation degrees.
This is the affine B\'ezout degree inequality of
\cite[Theorem~1]{Heintz1983}, with its published correction
\cite{Heintz1985}; it bounds the zero-dimensional components without
assuming that the entire zero set is zero dimensional.
Lemma~\ref{lem:appD-bounded-degree-word}
therefore gives
\[
 \#\mathfrak F_{H,\kappa}(b)
 \le D_0(P_0)^{N_{H,\kappa}}
 \le D_0(P_0)^{(2d+1)P_0+2d}
 =M(d,P_0),
\]
which is \eqref{eq:appD-multiplicity-bound}.

The cell conditions are semialgebraic explicitly: strict time order and
flux signs are polynomial inequalities; a lift is a fixed integral
translation; speed cells are inequalities between quadratic polynomials;
the normal atlas chooses the first maximal signed coordinate of
$\omega\in\R^d$ by half-open polynomial inequalities together with
$|\omega|^2=1$; and maximal-minor tie breakers compare squared polynomial
minors in a fixed lexicographic order.  Thus the
regular graph
\(\{(p,b):\mathscr F_{H,\kappa}(p;b)=0,\det D_p\mathscr F\ne0\}\)
is semialgebraic and therefore Borel in a product of standard Borel spaces.
The Lusin--Novikov theorem
\cite[Theorem~18.10]{Kechris1995} enumerates every countable Borel fibre by
Borel graphs.  To obtain exactly the fixed slots used here, order each
finite fibre lexicographically and let slot $j$ be its $j$th point.  The
graph of this selection is Borel: it is obtained recursively from the
Lusin--Novikov graphs by deleting previously selected points and taking the
least remaining point.  The degree bound makes all slots with
$j>M(d,P_0)$ empty.  Padding the first $M(d,P_0)$ slots with $\dagger$ yields
\eqref{eq:appD-recovery-branches}.  Neither the degree bound nor the number of slots
depends on polynomial coefficients, so neither depends on the boundary
values, lifts or \(\eps\).
\end{proof}

\mainstatementreference{Corollary~\ref{cor:appD-fixed-k-birth-multiplicity}}

\mainproofheading{Corollary~\ref{cor:appD-fixed-k-birth-multiplicity}}{supp:proof:cor:appD-fixed-k-birth-multiplicity}
\begin{proof}
Let $A\subset H$ be the unmarked contacts.  The square system consists of
the original $d(P-1)$ reverse-Kruskal upper-tree residuals, the $kd$ scalar
physical selected-contact equations
$C_i=R_i-\eps\omega_i=0$, and, for every $a\in A$, its $d$ contact residuals
and the extrinsic equation $|\omega_a|^2-1=0$.  The selected $\omega_i$ are
retained landing-normal variables, so $D C_i=D R_i$ in the velocity/time
pivot columns.  In particular the square system contains none of the
clearing matrix $\mathsf U$, the auxiliary frames $E_i$, or the
Moore--Penrose matrix $L_{\rm rem}^{+}$; those objects prove rank but are not
part of the polynomial recovery map.

The reverse-Kruskal upper tree contributes $d(P-1)$ scalar translation
pivots, the $k$ landing blocks contribute
$k(d-1)+k=kd$ velocity/time pivots, and each unmarked contact contributes
exactly $d+1$ extrinsic time--normal pivots.  Thus
$N=d(P-1)+kd+(d+1)|A|$, which proves
\eqref{eq:appD-fixed-k-birth-variable-count}.  The evolution is the same
fixed hard-sphere word as in Lemma~\ref{lem:appD-bounded-degree-word}; merely
selecting more of its contact equations does not increase their polynomial
degree.  Thus the degree is bounded by $D_0(P_0)$.  By Cauchy--Binet and
\eqref{eq:k-birth-frame-production}, at least one original Cartesian
$k(d-1)$-minor is no smaller than
\[
 \binom{dP}{k(d-1)}^{-1/2}
 \Jac_{k(d-1)\mid k}(\mathcal S_k,\mathcal T_k).
\]
The index $\kappa$ chooses the first maximal such minor, with equality ties
assigned by a finite disjoint half-open rule.  The complete pivot determinant
is the product of the tree block, this mixed landing minor, the unmarked root
blocks, and the extrinsic sphere blocks.  It is nonzero at every regular
physical solution, which is therefore an isolated complex zero of the
resulting square polynomial system.  The
affine B\'ezout argument in the proof of
Theorem~\ref{thm:appD-uniform-algebraic-multiplicity} gives
\eqref{eq:appD-fixed-k-birth-multiplicity}, and the same semialgebraic
maximal-minor partition and Lusin--Novikov enumeration give the asserted
Borel branch slots.  Their domains, including the original time-order and
lift conditions, are independent of concrete boundary values.  The same
exhaustion by lower bounds on the chosen minor and incoming fluxes covers
every transversal point; grazing roots have zero collision flux and the
remaining chart boundaries are lower-dimensional semialgebraic sets.  Hence
the singular/null treatment is independent of the joint kernel $Q$.
\end{proof}

The theorem concerns the regular graph only.  This does not leave an
uncontrolled singular mass in the packet integral.  On each physical packet
cell the conormal lower bounds prove full rank of the declared contact map;
the exhaustion used below restricts the smallest declared minor and every
incoming flux to be at least $1/n$.  The union over $n$ covers every
transversal point.  Grazing roots carry zero collision flux, while geometric
chart boundaries form lower-dimensional semialgebraic sets.  Maximal-minor
ties remain assigned to the Borel sets $E_i$ and are not removed.  Monotone
convergence therefore removes the exhaustion without assigning positive
collision measure to a singular algebraic branch.

\begin{remark}[Why this is not an implicit-function argument]
Theorem~\ref{thm:appD-uniform-algebraic-multiplicity} does not claim that a
maximal minor is globally one-to-one.  It permits several roots in the same
retained fibre, enumerates all of them, and pays a constant depending only
on the fixed packet cap.  This is the precise global input missing from a
local coarea calculation.
\end{remark}

\subsection{The zero-set measure}

For the original two-landing packet charts, let
\(\mathcal C_{H,\kappa}\) collect:
\begin{enumerate}[label=\textup{(\roman*)}]
\item the \(d(P-1)\) upper-tree contact residuals;
\item the two $d$-dimensional physical selected-contact residuals
      $C_i=R_i-\eps\omega_i$;
\item the $d$-dimensional residual of every unmarked contact, paired with
its incoming time--normal variables.
\end{enumerate}
The raw coordinate lemmas above put its derivative into a block triangular
form after an admissible permutation of rows and columns.  The unmarked
blocks have determinants \(\eps^{d-1}|g_a\cdot\omega_a|\); the selected
upper-tree and landing blocks have the lower bounds stated above.  Thus
\begin{equation}\label{eq:appD-total-rank}
\operatorname{rank}D\mathcal C_{H,\kappa}=d|H|
\end{equation}
on every regular cell.  In particular, zero is a regular value.

All Jacobians in the exact collision measure below are intrinsic.  More
precisely, each time--normal factor carries the product metric
\(\dd t^2+g_{\mathbb S^{d-1}}\), every sphere tangent determinant is computed in
an orthonormal frame \(E_\omega\), and the Hausdorff measure on the zero set
is induced by that product Riemannian metric.  In a coordinate chart
\(\theta\),
\[
 \left|\det[\,g\ \ -\eps D_\theta\omega\,]\right|
 =\eps^{d-1}|g\cdot\omega|J_\omega,
 \qquad \dd\omega=J_\omega\dd\theta.
\]
Thus the chart Gram density is carried by both surface measure and the
coordinate determinant and cancels exactly.  In
\eqref{eq:appD-packet-collision-measure},
\(\Jac_{d|H|}(D\mathcal C_{H,\kappa})\) therefore means the Riemannian
normal Jacobian, not a determinant relative to bare chart Lebesgue measure.
For extrinsic polynomialization one instead uses the ambient Euclidean
determinant together with the factor \(2\) in
\eqref{eq:appD-extrinsic-sphere-measure}, exactly as in
\eqref{eq:appD-extrinsic-root-determinant}.  The intrinsic coarea density
and the extrinsic polynomial fibre count consequently describe the same
collision measure.

\begin{definition}[Packet collision measure]
\label{def:appD-packet-collision-measure}
On a regular cell, the packet collision measure is the coarea measure
\begin{equation}\label{eq:appD-packet-collision-measure}
 \dd\nu_{H,\kappa}
 :=
 \frac{\dd\mathcal H^{\dim\mathfrak X-d|H|}
       \restriction\{\mathcal C_{H,\kappa}=0\}}
      {\Jac_{d|H|}(D\mathcal C_{H,\kappa})}
 \prod_{a\in H}\eps^{d-1}|g_a\cdot\omega_a|,
\end{equation}
with the source normalization multiplied afterward.  In the
arbitrary-kernel convention below, every Maxwellian, cutoff, remainder, and
generalized-link factor is absorbed into the single nonnegative \(Q\).
\end{definition}

\begin{lemma}[Conditional null strata and rank routing]
\label{lem:appD-conditional-null-strata}
Let \(y\) denote all complement and fixed-boundary variables, equipped with
the product reference measure used by the source integral.  There is one
Borel null set \(N_{\partial}\) in the \(y\)-space such that, for every
\(y\notin N_{\partial}\), the following assertions hold in the conditional
packet collision measure.
\begin{enumerate}[label=\textup{(\roman*)}]
\item exact grazing has zero measure;
\item equality of two distinct atom times, including a simultaneous contact,
has zero measure;
\item every rank loss is either grazing, belongs to a complementary
small-variable cell declared earlier, or is absent by an explicit packet
minor lower bound.  Maximal-minor ties are not rank losses and remain in the
Borel pieces \(E_i\).
\end{enumerate}
\end{lemma}

\begin{proof}
First keep \(y\) as a variable in
\(\mathfrak X_{H,\kappa}\).  The compact exhaustions of
\eqref{eq:appD-packet-collision-measure} are Radon measures on this common
product space, so their increasing limit admits the ordinary Tonelli
disintegration in \(y\).  We prove nullity before taking that
disintegration.

For an unmarked root the exact block determinant is
\(\eps^{d-1}|g_a\cdot\omega_a|J_{\omega_a}\).  On
\(|g_a\cdot\omega_a|\ge1/n\), the collision numerator cancels this
determinant and leaves a density uniformly bounded as \(n\to\infty\).
After the standard finite analytic stratification of each semialgebraic
recovery branch, \((g_a\cdot\omega_a)^2\) is analytic on every stratum.
Either it vanishes identically,
in which case that branch is absent from every regular exhaustion, or its
zero set has Lebesgue measure zero in the retained coordinates.  Hence the
increasing regular measure has no grazing atom.  The same dichotomy applies
to every collision flux retained in the output density and proves (i) in
the common measure.

For (ii), let \(e\ne f\) be two distinct unordered particle pairs and write
their lifted contact functions as
\[
 q_e(X)=|b_eX-m_e|-\eps,
 \qquad q_f(X)=|b_fX-m_f|-\eps.
\]
At a common contact their configuration conormals are
\(b_e^*\omega_e\) and \(b_f^*\omega_f\).  They are linearly independent:
if the edges are disjoint this is immediate, while if they share one
particle, an endpoint belonging only to \(e\) kills the coefficient of the
other conormal and hence forces the coefficient of \(b_e^*\omega_e\) to
vanish; the endpoint belonging only to \(f\) then forces the second
coefficient to vanish.  Consequently \((q_e,q_f)\) is a codimension-two
submersion on the full raw configuration space.  A simultaneous contact is
therefore codimension one inside either ordinary collision cross-section
and has zero measure for the common Radon collision measure by the surface
coarea formula.
If the unordered pair is the same, two different lifts cannot occur at one
time when \(2\eps<1\), since their lattice difference would have norm at most
\(2\eps\); the same lift is the same event and is excluded by the strict word
order.  This proves that the boundary \(t_a=t_b\) contributes no mass in
the common measure.

Finally, the upper-tree relative-translation block is an invertible reduced
incidence matrix.  The disjoint landing block, overlap chord block and
remaining overlap Schur row have the explicit positive lower bounds in
Section~\ref{sec:two-landing}; the default FCT cell retains only a selected
minor of a row Jacobian already bounded below.  The elementary cells with a
small time difference or small relative speed are assigned to their stated
good-support estimates and are never deleted.  These alternatives exhaust
all declared pivots and prove (iii) on the regular support.

Finally apply Tonelli to the indicators of the null strata in the common
Radon measure.  Their conditional masses vanish outside one Borel
source-boundary null set \(N_{\partial}\), which may be chosen
simultaneously because the atom, lift, chart, and minor families are
countable (indeed finite on a packet cell).  This is the required
conditional conclusion.  The almost-every qualifier is essential:
deliberately fixed external trajectories can force two independent
deterministic contact times to coincide, but those boundary values form
precisely the null set just removed and are never subjected to a boundary
supremum.
\end{proof}

The quantifier in Lemma~\ref{lem:appD-conditional-null-strata} is part of
the interface: removal of grazing and simultaneous-contact strata is valid
for almost every fixed source-boundary value, and therefore for the
integrated source operator.  It is not a geometric statement for an
arbitrarily prescribed exceptional boundary trajectory.  In order not to
silently turn this almost-everywhere fact into a boundary supremum, we now
fix the conditional representative used everywhere below.

For each regular packet cell, record its chart cost by
\[
 \gamma_\kappa:=
 \begin{cases}
  2(d-1),&\text{on the default radial FCT atlas},\\
  2d,&\text{on the main/coarse $2d$ atlas}.
 \end{cases}
\]

\mainstatementreference{Definition~\ref{def:appD-selected-conditional-representative}}

The definition is simultaneous for all $Q$: it modifies a Borel measure
kernel, not one exceptional set chosen after seeing the test function.
Moreover,
\begin{equation}\label{eq:appD-selected-representative-integral}
 \int \overline J_{H,\kappa}^{y}(Q)\,\lambda_Y(\dd y)
 =J_{H,\kappa}(Q),
\end{equation}
and, by \eqref{eq:appD-output-kernel},
\[
 \int \overline J_{H,\kappa}^{\sharp,y}(Q)\,\lambda_Y(\dd y)
 =
 \sum_j\int Q^\sharp_{H,\kappa}(b,j)\,\dd b .
\]
Thus the analogous equality holds for the dominating output operator,
because $\lambda_Y(N_\partial)=0$.  Thus zero extension changes no source
integral.  It only supplies a total Borel representative at boundary values
where the conditional operator was never canonically defined.

Formula \eqref{eq:appD-packet-collision-measure} is independent of the order
in which admissible square blocks are used.  Indeed, two orders are two
local parametrizations of the same regular level set, and the coarea
formula is coordinate invariant.  The Schur complement appearing after an
early block substitution is exactly the determinant quotient for the same
full derivative.

On the singular complement define the measure by monotone exhaustion:
restrict first to cells where every nonzero declared minor and every flux
used as a pivot exceeds \(\delta\), and every declared positive time gap
exceeds its lower threshold by \(\delta\); then let
\(\delta\downarrow0\).  Geometric normal-atlas boundaries are surface-null.
Dyadic, lift, sign, and maximal-minor ties are instead assigned to one Borel
piece and are retained throughout the exhaustion.  The collision-flux-null
sets contribute zero.  This procedure proves existence of the limiting
measure and rules out order-dependent delta products.

\subsection{Construction of the branch-complete enlarged kernel}

Let \(Q\ge0\) be a Borel function of all packet and complement variables.
It may couple different atoms, may contain the exact time-order indicator,
and need not factor across the packet boundary.  Apply the substitutions in
the following fixed order:
\begin{enumerate}[label=\textup{(\arabic*)}]
\item reconstruct all nontransversal velocities by elastic involutions;
\item reconstruct all nontransversal positions by transport translations;
\item solve every unmarked incoming root in its time--normal cell;
\item solve the upper-tree relative translations jointly with the selected
landing block, or use the two successive overlap blocks;
\item retain every variable not used as a pivot and retain the finite
algebraic-branch label from
Theorem~\ref{thm:appD-uniform-algebraic-multiplicity}.
\end{enumerate}
The resulting maps are the complete family
\(\Phi_{H,\kappa,j}\), \(1\le j\le M(d,P_0)\), from
\eqref{eq:appD-recovery-branches}.  In particular, no single-valued global
recovery is being inferred from a local Jacobian.

We next define the output at the level of Radon measures.  Let
\[
 U=H\setminus(T\cup D),\qquad T=S_U,\qquad D=\{d_1,d_2\},
\]
where \(U\) is the set of unmarked incoming-root atoms, \(T\) is the marked
upper tree, and \(D\) is the selected landing pair.  Put
\(\phi_a=|g_a\cdot\omega_a|\).  We fix the arbitrary-kernel convention
once: every Maxwellian, cutoff, generalized-link indicator, remainder, and
complement factor is part of \(Q\).  Thus \(W^0\equiv1\) below; the symbol
is retained only to make the collision-weight ledger readable.  In
particular, no factor already present in \(Q\) is duplicated in the output
measure.  If \(p\) denotes the complete square pivot list and \(b\) all
retained variables, define on the branch space
\begin{align}
 \dd\mu^{\rm ex}_{H,\kappa}(b,j)
 &:=
 \1_{\kappa}(\Phi_{H,\kappa,j}(b))
 \1_{\mathcal D_H}(\Phi_{H,\kappa,j}(b))\notag\\
 &\quad\times
 \frac{\eps^{-(d-1)P}\prod_{a\in H}\eps^{d-1}\phi_a}
 {|\det D_p\mathcal C_{H,\kappa}(\Phi_{H,\kappa,j}(b))|}
 W^0(\Phi_{H,\kappa,j}(b))\,\dd b .
 \label{eq:appD-exact-pushforward-measure}
\end{align}
In \eqref{eq:appD-exact-pushforward-measure}, $p$ and
$D_p\mathcal C_{H,\kappa}$ use the intrinsic $d$-variable block
$(t_a,\theta_a)$ at an unmarked root.  The $(d+1)$-variable extrinsic block
$(t_a,\omega_a)$ together with $|\omega_a|^2-1$ is used only to prove the
global algebraic multiplicity.  Equations
\eqref{eq:appD-extrinsic-sphere-measure}--
\eqref{eq:appD-extrinsic-root-determinant} show that passage between the two
descriptions inserts and cancels the same factor $2$ and does not change the
density above.
The density is zero when the branch is absent or its recovered point is
outside the original half-open cell.  The equality form of the finite-fibre
area formula gives, for every nonnegative Borel \(Q\),
\begin{equation}\label{eq:appD-exact-pushforward-identity}
 J_{H,\kappa}(Q)
 =\sum_{j=1}^{M(d,P_0)}\int
 Q(\Phi_{H,\kappa,j}(b))\,\dd\mu^{\rm ex}_{H,\kappa}(b,j)
\end{equation}
on every regular compact exhaustion.  Thus the branch sum is the exact
push-forward of the original collision measure, not a convention for
choosing one preimage.

Here is the complete atomwise determinant ledger.  ``Retained'' means that
the factor belongs to the density on the output side of
\eqref{eq:appD-exact-pushforward-measure}.
\begin{center}
\small
\begin{tabular}{@{}p{0.18\textwidth}p{0.22\textwidth}p{0.22\textwidth}
 p{0.25\textwidth}@{}}
\toprule
operation/atom & pivot determinant & collision numerator & output effect\\
\midrule
elastic pair reconstruction
 & \(1\) & none & no density change\\
transport-position reconstruction
 & \(1\) & none & no density change\\
\(a\in U\), incoming root
 & \(\eps^{d-1}\phi_a\) & \(\eps^{d-1}\phi_a\)
 & exact cancellation; \(t_a,\omega_a\) are eliminated\\
\(a\in T\), tree contact
 & reduced-incidence block, total \(J_T\)
 & \(\eps^{d-1}\phi_a\)
 & \(\eps^{d-1}\) remains; \(\phi_a\) enters \(W^{\rm rem}\)\\
\(a\in D\), main $2d$-power chart
 & joint block \(J_D\) & \(\eps^{d-1}\phi_a\)
 & the two \(\eps^{d-1}\)'s and both landing fluxes remain;
   \(J_D^{-1}\le C\epstar^{-2d}\)\\
\(a\in D\), default radial FCT chart
 & \(J_D\ge c\phi_{d_1}\phi_{d_2}|\Delta_1\Delta_2|^{d-1}\)
 & \(\eps^{2(d-1)}\phi_{d_1}\phi_{d_2}\)
 & the landing fluxes cancel; the inverse gap cost is
   \(C\epstar^{-2(d-1)}\)\\
source line normalization
 & none & \(\eps^{-(d-1)P}\)
 & with \(\eps^{(d-1)(P-1)}\eps^{2(d-1)}\), leaves \(\eps^{d-1}\)\\
\bottomrule
\end{tabular}
\end{center}
Indeed, the root rows contribute
\(\prod_{a\in U}\eps^{d-1}\phi_a\) to both numerator and determinant.
Since \(|T|=P-1\) and \(|D|=2\), chronological block elimination on
either packet chart gives the pointwise identity
\begin{equation}\label{eq:appD-exact-density-factorization}
 \frac{\eps^{-(d-1)P}\prod_{a\in H}\eps^{d-1}\phi_a}
 {|\det D_p\mathcal C_{H,\kappa}|}
 =\eps^{d-1}
 \frac{\prod_{a\in T\cup D}\phi_a}{|J_TJ_D|}.
\end{equation}
There is no factor from an unmarked root on the right-hand side.  On a
main chart set \(W^{\rm rem}_{H,\kappa,j}\) equal to
\(W^0\circ\Phi_{H,\kappa,j}=1\), times the recovered-cell indicator and
the fluxes \(\prod_{a\in T\cup D}\phi_a\).  On the default radial FCT chart
omit \(\phi_{d_1}\phi_{d_2}\), since those two factors have been cancelled
against \(J_D\).  In both cases omit every root flux cancelled above; the
full $2d$ chart remains the coarse fallback.
Define
\begin{equation}\label{eq:appD-dominating-output-measure}
 \dd\mu^\sharp_{H,\kappa}(b,j)
 :=
 \1_{\mathcal D_H}(\Phi_{H,\kappa,j}(b))
 W^{\rm rem}_{H,\kappa,j}(b)\,\dd b .
\end{equation}
The tree and landing lower bounds give the chart-specific measure domination
\begin{equation}\label{eq:appD-exact-measure-domination}
 \mu^{\rm ex}_{H,\kappa}
 \le C(P_0)\cL_\eps^{C|H|}\eps^{d-1}\epstar^{-\gamma_\kappa}
       \mu^\sharp_{H,\kappa}.
\end{equation}
Thus the default radial FCT atlas has the sharp loss
\(\epstar^{-2(d-1)}\), while the independent main/coarse atlas has the
weaker loss \(\epstar^{-2d}\).  No coarse-cell determinant is used to claim
the sharp FCT loss.
This is a statement on the fixed product branch domain; concrete boundary
values occur only inside the zero-extended density.  The displayed ledger
also shows that no collision flux or determinant used here is available
for a second time or velocity estimate.

The output domain \(\mathfrak B_{H,\kappa}\) is the product of:
\begin{enumerate}[label=\textup{(\alph*)}]
\item one untranslated root position;
\item all unpivoted separator velocities;
\item all retained atom times and normals;
\item the fixed-end and complement variables;
\item the finite cell and algebraic-branch labels, which are summed rather
than integrated.
\end{enumerate}
Each continuous factor is enlarged to its fixed fundamental domain, dyadic
box, time interval or normal chart.  This product is chosen before the
concrete fixed-boundary values are inserted.

\begin{definition}[Zero-extended output kernel]
\label{def:appD-output-kernel}
Put
\(\mathfrak B_{H,\kappa}^{\rm br}
=\mathfrak B_{H,\kappa}\times\{1,\ldots,M(d,P_0)\}\), with counting
measure on the second factor.  For
\((b,j)\in\mathfrak B_{H,\kappa}^{\rm br}\) set
\begin{equation}\label{eq:appD-output-kernel}
 Q_{H,\kappa}^{\sharp}(b,j)
 :=
 Q(\Phi_{H,\kappa,j}(b))\,
 \1_{\mathcal D_H}(\Phi_{H,\kappa,j}(b))\,
 W_{H,\kappa,j}^{\rm rem}(b)
\end{equation}
when the branch exists and every recovered variable lies in its original
half-open cell, and set it to zero otherwise.  Here
\(W_{H,\kappa,j}^{\rm rem}\) is precisely the chart-dependent density
defined after \eqref{eq:appD-exact-density-factorization}: it contains
\(W^0\circ\Phi_{H,\kappa,j}=1\), the recovered-cell indicator, and exactly
the collision fluxes not cancelled by a row determinant.  Every
Maxwellian, source cutoff, generalized link, and complement factor remains
inside \(Q\circ\Phi_{H,\kappa,j}\).  In particular,
\[
 Q^\sharp_{H,\kappa}(b,j)\,\dd b
 =Q(\Phi_{H,\kappa,j}(b))\,\dd\mu^\sharp_{H,\kappa}(b,j).
\]
\end{definition}

The definition keeps the same time indicator, rather than replacing it by
one.  Boundary data enter only through the recovery branches and not through
the box \(\mathfrak B_{H,\kappa}^{\rm br}\).  The branch set has the fixed
cardinality \eqref{eq:appD-multiplicity-bound}.  Hence
\(\int Q_{H,\kappa}^{\sharp}\) is a legitimate uniform enlargement in the
sense required by the Source-B weight and time consumers.

\mainstatementreference{Proposition~\ref{prop:appD-arbitrary-kernel-coarea}}

\mainproofheading{Proposition~\ref{prop:appD-arbitrary-kernel-coarea}}{supp:proof:prop:appD-arbitrary-kernel-coarea}
\begin{proof}
For $y\notin N_\partial$, apply
\eqref{eq:appD-conditional-exact-kernel} and
\eqref{eq:appD-conditional-dominating-kernel}, first to bounded simple
functions and then by monotone convergence to the original joint kernel.
This gives
\[
 \overline J_{H,\kappa}^{y}(Q)
 \le C(P_0)\cL_\eps^{C|H|}\eps^{d-1}\epstar^{-\gamma_\kappa}
 \overline J_{H,\kappa}^{\sharp,y}(Q).
\]
Since \(\eps^{d-1}=\eps^{-(d-1)}\eps^{2(d-1)}\), this is
\eqref{eq:appD-arbitrary-kernel-coarea}.  The branch set contains every
regular preimage by
Theorem~\ref{thm:appD-uniform-algebraic-multiplicity}; the cell, lift, sign,
and dyadic factors are those in \eqref{eq:appD-cell-count}.  Finally,
Tonelli and the fibrewise monotone exhaustion of
Lemma~\ref{lem:appD-conditional-null-strata} extend the inequality to
almost every boundary fibre without changing the joint kernel or the exact
time indicator.  If $y\in N_\partial$, both selected conditional operators
in \eqref{eq:appD-arbitrary-kernel-coarea} are zero by
\eqref{eq:appD-selected-K}--\eqref{eq:appD-selected-Ksharp}.  This proves
the asserted every-$y$ statement for the selected representatives.  Finally
integrate in $y$ and use
\eqref{eq:appD-selected-representative-integral} to recover the complete
cell operator.
\end{proof}

\subsection{Fibrewise time projection}

Let $M_t$ be the indexed set of all packet-time coordinates, and let
$R_{\rm birth}$ be the marked upper-tree birth-contact times: when a
marked edge first joins two chronological clusters, its time is kept and
its relative cluster translation, not that time, is the pivot.  Let $S$
consist of the variables fixed by the Source-B consumer together with
$R_{\rm birth}$, and put $A:=M_t\setminus S$.  Thus $A$ contains the two
landing times as well as the remaining packet times used in incoming
time--normal root blocks.  On the default radial FCT atlas the landing times
are recovered pivot variables; on the main/coarse $2d$ atlases they remain
retained branch coordinates integrated in the $A$-fibre.  Define the exact existential
projection
\begin{equation}\label{eq:appD-time-projection}
 \widetilde{\mathcal D}_{H,S}
 :=
 \{t_S:\ \exists\,t_A\text{ such that }(t_A,t_S)\in\mathcal D_H\}.
\end{equation}

\mainstatementreference{Lemma~\ref{lem:appD-fibre-time-owner}}

\mainproofheading{Lemma~\ref{lem:appD-fibre-time-owner}}{supp:proof:lem:appD-fibre-time-owner}
\begin{proof}
Outside the fixed exceptional set, apply the conditional exact
push-forward identity and conditional measure domination fixed in the
statement, with \(t_S\) held fixed and with
\(\1_{\mathcal D_H}\) included in \(Q\).  If the
branch integrand is nonzero, the complete packet-time values on that branch,
recovered or retained according to the chosen atlas, witness membership of \(t_S\) in
\(\widetilde{\mathcal D}_{H,S}\).  If no witness exists, both sides vanish.
All coarea pivots, the algebraic multiplicity bound and the fixed branch
slots are fibrewise uniform.  On the exceptional set both selected
conditional sides are zero.  Since this set is source-null, the subsequent
integration is unchanged.  This proves the inequality for every $t_S$ as a
statement about the selected representative, without asserting a geometric
bound on a deliberately prescribed exceptional trajectory.
\end{proof}

This is the precise statement needed before applying the Source-B estimate
for the scale-correct projection.  The variables in $R_{\rm birth}$ are not
recovered branch variables: they remain literal coordinates on the right
of \eqref{eq:appD-fibre-time-owner}.  Hence
Lemma~\ref{lem:distinguished-time-weight} supplies their honest time powers.
The eliminated incoming-root times receive only their root coarea Jacobians.
On the default radial FCT atlas the two landing times are recovered pivots
and receive only their declared landing-pivot Jacobian; on the main/coarse
$2d$ atlases they remain ordinary retained time coordinates in $\dd b$, while
the selected velocity/chord pivot block pays the chart's $J_D$.  No Tonelli
or finite-multiplicity argument is used to manufacture additional powers of
$\tau$.

\subsection{Multiplicity and boundary checks}

We close with the possible sources of hidden multiplicity.
\begin{enumerate}[label=\textup{(\roman*)}]
\item A fixed lift and incoming sign have at most one contact root by
Lemma~\ref{lem:one-incoming-root}.
\item Normal charts are finite and are made disjoint by the half-open Borel
sets $V_j=U_j\setminus\bigcup_{i<j}U_i$; their boundaries are null.
\item Elastic pair maps are involutions, so choosing the incoming/outgoing
slot label removes the two-branch ambiguity.
\item The upper-tree relative translations are locally recovered because
the reduced incidence matrix is invertible.
\item The disjoint landing Schur row and the remaining overlap landing row
may have more than one solution after root substitution.  Every such
solution is a regular point of the square polynomial system and is included
in exactly one branch slot by
Theorem~\ref{thm:appD-uniform-algebraic-multiplicity}.
\item The overlapping chord system is affine in \((X_c,V_c)\) and has one
solution for fixed values of all other pivot variables; the full packet
fibre is nevertheless enumerated rather than declared single-valued.
\item The radial augmentation proves the FCT row bound locally but is not a
global recovery variable.  The resulting rectangular Cartesian FCT chart
uses a deterministic maximal minor only for
local transversality.  Its possible global branches obey the same uniform
algebraic bound and are all retained.
\end{enumerate}
Therefore the constant in
\eqref{eq:appD-arbitrary-kernel-coarea} contains a controlled factor at
most \(M(d,P_0)\), never an unproved multiplicity-one assertion.  Because
\(P_0=P_{\rm pkt}(d)\) is fixed before \(\eps\), this is a genuine
\(C(d,P_0)\) constant
and has no epsilon or time cost.  This completes the global
measure-theoretic construction of the packet operator.

\section{Insertion into the long-time expansion}
\label{sec:p38j}

On every fixed packet cell, Theorem~\ref{thm:two-landing} is unconditional.
For disjoint first lower edges the sharp conclusion follows from Proposition
\ref{prop:disjoint-radial-fct}; for overlapping edges it follows from
Proposition~\ref{prop:overlap-radial-fct}.  Theorem~\ref{thm:packet-fct}
closes the cellwise rank argument and
Lemma~\ref{lem:global-packet-chart} supplies the complete measure chart.
Propositions~\ref{prop:disjoint-rooted-extension} and
\ref{prop:asymmetric-new-line-peeling} remain independent coarse fallbacks
with the weaker $\epstar^{-2d}$ loss.  The sharp theorem is inserted first
into the Source-A cutting class and then into the two Source-B consumers
recorded in Online Resource 1, Section~\ref{app:source-interface}.  Three steps connect the
local packet geometry to the long-time hierarchy: a strengthened special-case
induction, the outer Proposition 3.8 cutting algorithm, and the fibrewise
positive-operator transfer.  The resulting operator statement remains joint;
it does not become a componentwise factorization.

\subsection{The strengthened special-case class}

Recall the real-valued envelopes
\begin{equation}\label{eq:G-G0}
 \widehat G(q)=(2q+1)^{10}\binom q2(32q^{3/2})^{q^2},
 \qquad
 \widehat G_0(q)=(2q+1)^5\binom q2(32q^{3/2})^{q^2},
\end{equation}
and define the integer cutoffs
\[
 G(q)=\lceil\widehat G(q)\rceil,
 \qquad G_0(q)=\lceil\widehat G_0(q)\rceil.
\]
For integer $|M|$, the condition $|M|\ge G(q)$ is equivalent to the source
condition $|M|\ge\widehat G(q)$; the rounding of $G_0$ only enlarges the
integer sublayer threshold.
The special-case input class $\mathfrak C(q)$ consists of physically
realizable, one-layer, full C-molecules $M$ satisfying:
\begin{enumerate}[label=\textup{(C\arabic*)}]
\item $M$ has at most $q$ bottom free ends and at least $G(q)$ atoms;
\item $M$ has no double bond;
\item $M$ has no strongly degenerate primitive pair in the sense of Source
A, Definition 4.3;
\item a linear extension of the atom time order and all preliminary positive
indicator cells have been fixed.
\end{enumerate}
In the only calls made by the outer Source-A algorithm, $M$ is a connected
component of $X^+(A)$ selected below, or its TOP/BOTTOM reverse.  Write
$q=\#_{\rm conn}^+(A)$ for the total number of cross bonds.  Every connected
component of $X^+(A)$ contains an $X_0^+(A)$ atom and hence has at least two
bonds to the connected transversal set $A$.  If $c$ is the number of these
components, contracting $A$ and each component shows that the cross bonds
alone contribute $q-c$ independent cycles.  Since $c\le q/2$ and the full
molecule is a tree plus $\gamma\le2\Gamma$ bonds,
\begin{equation}\label{eq:application-size-bound}
 \frac q2\le q-c\le\rho(\mathcal M)=\gamma\le2\Gamma,
 \qquad q\le4\Gamma.
\end{equation}
Set
\begin{equation}\label{eq:C4-packet-cap}
 P_{\rm pkt}(d):=4\max_{3\le q\le\lfloor4\Gamma(d)\rfloor}G_0(q).
\end{equation}
Thus the packet selected below has dimension-only size.  The ambient set
$M$ need not have bounded size; all its time-order and indicator partitions
remain inside the source allowance $C^{|M|}\cL_\eps^C$.

\begin{lemma}[Hereditary free induced class]
\label{lem:hereditary-special-class}
Let $M\in\mathfrak C(q)$ and let $B$ be a connected component of a
time-convex sublayer.  Cut $B$ as free.  If $B$ meets $q'<q$ physical
particle lines and $|B|\ge G(q')$, then the resulting full molecule belongs
to $\mathfrak C(q')$.
\end{lemma}

\begin{proof}
Cutting $B$ as free turns every boundary bond at $B$ into a free end, so it
creates no fixed end in $B$ and preserves the number of incidences at every
atom.  A physical particle line meets a time-convex sublayer in an interval;
hence it contributes at most one top and one bottom boundary end.  The
number of bottom free ends is therefore at most the number $q'$ of lines.

A double bond in $B$ would be a double bond in $M$.  Suppose a strongly
degenerate primitive pair were created in $B$.  Its two atoms, slot pattern
and phase-space restriction are unchanged, because replacing an exterior
bond by a free end preserves the slot and degree.  Time convexity also
prevents a directed descendant path between atoms of $B$ from leaving $B$
and re-entering it.  Thus the descendant and primitive relations between
atoms of $B$ are the restrictions of the corresponding relations in $M$.
The pair would already be a strongly degenerate primitive pair in $M$, a
contradiction.  Physical realizability is inherited by restricting the
original collision history to $B$.  The size assumption finishes the proof.
\end{proof}

\begin{lemma}[One-sided boundary and protected restart]
\label{lem:protected-band-restart}
Let the ambient molecule be connected.  (For a disconnected source term the
statement is applied separately to each connected component.)  Let $I$ be a
consecutive interval in the fixed atom order, and suppose that
membership in $I$ is an interval on every physical particle line.  Let $B$
be a connected component of the induced molecule on $I$, cut $B$ as free,
and suppose a connected full submolecule $K\subset B$ has already been cut
and protected.  Then the remaining atoms admit the following ordered
cleanup.
\begin{enumerate}[label=\textup{(\roman*)}]
\item Every component of $B\setminus K$ produced by the recursive cutting
has a fixed end and hence is nonfull.
\item Let $B_-$ be the union of the components in the strictly earlier atom
interval having a bond to $B$.  Every component of $B_-$ has a fixed end on
its top boundary and no bottom fixed end; \textbf{DOWN} applies to it and
creates no full component.
\item After $B\cup B_-$ has been processed, every remaining connected
component has a fixed end and no top fixed end; \textbf{UP} applies and
creates no full component.
\item In the resulting cutting order, the live domain of a component depends
only on variables belonging to components earlier in the cutting order.
Consequently reverse analytic integration has exactly the dependency
property required before \eqref{eq:source-iterated-integral}.
\end{enumerate}
The same statement holds with time reversed and \textbf{UP}/\textbf{DOWN}
interchanged.  In particular, $K$ is the unique full component in this
restart.
\end{lemma}

\begin{proof}
The only cutting cases used below are cases (0), (2), (3), and (4) of the
cited Source-B v2 cutting definition:
\begin{center}
\small
\begin{tabular}{@{}p{0.25\textwidth}p{0.31\textwidth}p{0.34\textwidth}@{}}
\toprule
case & present boundary configuration & effect used here\\
\midrule
case (0)
 & no O-atoms on the crossing segment
 & one free/fixed matched pair\\
case (2)
 & two selected terminals with intervening exterior O-atoms
 & shortcut inside the selected atom set; O-atoms stay outside\\
case (3)
 & one selected terminal and one exterior C-terminal
 & selected free end and exterior fixed end\\
case (4)
 & one selected terminal and one pre-existing terminal end
 & selected free representative and fixed mate\\
\bottomrule
\end{tabular}
\end{center}
The accompanying Source-B fixed/free mate proposition states that a fixed
copy depends on the earlier free representative, and not conversely.  The
precise source labels and file locations are collected in
Online Resource~1, Section~\ref{appB:source-location-notes}; they fix the C/O and free/fixed
conventions used in every step of the argument.

The Source cutting operation replaces every bond crossing the chosen set by
a paired free end on the set being cut and a fixed end on its complement.
Thus cutting $B$ as free makes the induced molecule on $B$ full, but every
component on the other side of a crossing bond nonfull.  The recursive
cutting inside $B$ has the same property: once $K$ is protected, every
component of $B\setminus K$ is separated from $K$ by a cut bond and receives
its fixed mate.  This proves (i).

Consider a particle line meeting $B_-$.  If it subsequently meets $B$, time
convexity of $I$ makes the last bond from $B_-$ to $B$ a top-boundary bond of
the earlier component.  No cut made inside $B$ can create a bottom fixed end
in $B_-$, since that would require the same particle line to leave the
earlier interval and later re-enter it.  By definition each component of
$B_-$ has at least one bond to $B$, so it has a fixed end and is nonfull.
These are precisely the input conditions of \textbf{DOWN}; Source A's
one-sided algorithm then produces no degree-four full component.  This gives
(ii).

Let $C$ be a component left after $B\cup B_-$ is processed.  Choose a path
in the original connected molecule from $C$ to the processed set and take
the first bond of the path which crosses into it.  That bond has been cut,
so its endpoint in $C$ is fixed.  It cannot be a top fixed end: if the path
crossed from a later atom into an unprocessed earlier atom, then either that
earlier atom lay in the same earlier component attached to $B$ and hence in
$B_-$, or the path would leave and re-enter a time interval on one physical
line.  Both are impossible.  Hence $C$ is a nonfull no-top-fixed-end input
for \textbf{UP}.  This proves (iii), including uniqueness of the full
component.

Finally order the components by the moment at which they are cut.  A new
fixed end is, by definition, the shadow of the paired free end in a component
already cut.  Therefore its value may depend only on variables of that
earlier component.  No later component is used to define an earlier free
end.  Reversing the order for integration gives exactly the triangular
dependency used in the proof of the cutting identity
\eqref{eq:source-cutting-identity} and in
\eqref{eq:source-iterated-integral}.  This proves (iv).
\end{proof}

The boundary orientation and integration order in the preceding lemma are
summarized by the following exhaustive table.  Here ``earlier'' and
``later'' refer to the fixed atom order, while analytic integration runs in
the reverse of the last column.
\begin{center}
\small
\begin{tabular}{@{}p{0.18\textwidth}p{0.22\textwidth}p{0.20\textwidth}
 p{0.28\textwidth}@{}}
\toprule
region & fixed boundary supplied by & admissible cleanup & dependency order\\
\midrule
protected \(K\subset B\)
& cut from \(B\setminus K\)
& none; keep as the unique full packet
& packet variables precede every shadow they create\\
\(B\setminus K\)
& mate of a bond cut from \(K\)
& recursive nonfull cleanup
& each new component uses only an already cut mate\\
earlier band \(B_-\)
& top boundary bond into \(B\)
& \textbf{DOWN}; no bottom fixed end
& after \(B\), before the later complement\\
later complement
& bottom boundary bond into the processed set
& \textbf{UP}; no top fixed end
& last in cutting order, first in analytic integration\\
\bottomrule
\end{tabular}
\end{center}

\begin{lemma}[Packet-size and partition ledger]
\label{lem:packet-size-partition-ledger}
In every recursive call with $q$ physical particle lines, the returned
packet has fewer than $4G_0(q)$ atoms.  In the selected Source-A
Main-case-1 component, $q_B\le q_{\rm tot}\le4\Gamma(d)$, and hence
$|K|\le P_{\rm pkt}(d)$ as in
\eqref{eq:C4-packet-cap}.  Across the entire recursion, all added cells are
time-order, sign, lift, dyadic, or ordinary good/bad support cells, and their
number is at most $C^{|M|}\cL_\eps^C$.  No recursive call duplicates a
physical cell.
Here $q_{\rm tot}$ is the total number of $X^+(A)$--$A$ cross bonds and
$q_B$ is the number belonging to the selected component.
\end{lemma}

\begin{proof}
At the terminal call $K$ is the union of two complete sublayers, each of
size smaller than $2G_0(q)$.  A deficient-sublayer call passes to a component
with $q'<q$ and returns its packet unchanged, so the cap can only decrease
under recursion.  For the outer bound, every component of $X^+(A)$ contains
an $X_0^+(A)$ atom and therefore has at least two cross bonds to $A$.  If
$c$ is the number of components and
$q_{\rm tot}=\#_{\rm conn}^+(A)$, then $c\le q_{\rm tot}/2$.
After contracting $A$ and those components, the cross bonds contain
$q_{\rm tot}-c\ge q_{\rm tot}/2$ independent cycles.  Since the full molecule is a tree plus
 $\gamma\le2\Gamma(d)$ edges,
\[
 q_{\rm tot}/2\le q_{\rm tot}-c\le\rho(\mathcal M)=\gamma\le2\Gamma(d),
\]
which proves $q_B\le q_{\rm tot}\le4\Gamma(d)$ and the packet cap.

Fixing a strict order gives at most $C^{|M|}$ order pieces on the bounded
special call; incoming signs and slot choices have the same form.  The lift
and dyadic partitions have at most $\cL_\eps^{C|K|}$ pieces and $|K|=O_d(1)$.
let $N(q,n)$ be the largest number of descendants of a
half-open cell in a call with $q$ live lines and $n$ atoms.  If the call is
terminal, $N(q,n)\le C^n\cL_\eps^{C G_0(q)}$.  Otherwise its recursive
children have $q_i<q$, disjoint atom sets of sizes $n_i$ with
$\sum_i n_i\le n$, and the outer sign/order split costs at most $C^n$.
Thus
\begin{equation}\label{eq:packet-cell-recurrence}
 N(q,n)\le C^n\cL_\eps^{C G_0(q)}
 \max_{\substack{q_i<q\\\sum_i n_i\le n}}
 \prod_i N(q_i,n_i).
\end{equation}
Induction on $q$ gives
$N(q,n)\le C_q^n\cL_\eps^{C_q}$.  Here $q<4\Gamma(d)$ at the only special
call, so $C_q$ is a dimension constant and may be renamed $C$.  Every
child is a disjoint induced component, and every cleanup cell is defined by
half-open indicators.  Hence the product in
\eqref{eq:packet-cell-recurrence} neither duplicates atoms nor counts a
physical history twice.  This proves the asserted
$C^{|M|}\cL_\eps^C$ bound.
\end{proof}

\begin{lemma}[Short-window component lifting]
\label{lem:short-window-component-lifting}
Let $J\subset\R$ be an interval and let $x_1,\ldots,x_{q'}:J\to\T^d$
be hard-ball trajectories with $|v_i|\le V$.  Suppose their collision graph
on $J$ is connected and
\begin{equation}\label{eq:short-window-injectivity-condition}
 (q'-1)(\eps+2|J|V)<\operatorname{inj}(\T^d)=\frac12.
\end{equation}
Then the trajectories admit lifts
$\widetilde x_i:J\to\R^d$ such that every contact, including every contact
off a chosen spanning tree, is an ordinary Euclidean hard-ball contact.
Throughout $J$ all lifted pair differences are the unique short
representatives of their torus differences.  Different collision-connected
components may be lifted independently.
\end{lemma}

\begin{proof}
Choose a spanning tree of the collision graph and a root line.  Lift the
root trajectory on $J$.  Traverse the tree away from the root.  When the
edge from an already lifted line $i$ to a new line $j$ records a collision
at time $t_{ij}$, the covering-space lift of $x_j$ is unique after choosing
the integral translate for which
$|\widetilde x_i(t_{ij})-\widetilde x_j(t_{ij})|=\eps$.  This constructs
global lifts on $J$ for all vertices.

For a tree edge and any $t\in J$, free transport between events and velocity
reflection give
\[
 |\widetilde x_i(t)-\widetilde x_j(t)|
 \le\eps+2V|t-t_{ij}|
 \le\eps+2V|J|.
\]
Summing along the unique tree path gives, for every pair $a,b$ and every
$t\in J$,
\begin{equation}\label{eq:lifted-tree-diameter}
 |\widetilde x_a(t)-\widetilde x_b(t)|
 \le(q'-1)(\eps+2V|J|)<\frac12.
\end{equation}
If an off-tree pair collides, its torus difference has a unique
representative of norm $\eps<1/2$.  The difference of the already chosen
lifts is another representative and, by
\eqref{eq:lifted-tree-diameter}, lies in the injectivity ball.  Uniqueness
forces the two representatives to coincide.  The same argument at
noncollision times identifies torus and Euclidean separation, so no
spurious overlap is introduced.  Independent components have no contact
constraint relating their integral translates.
\end{proof}

\begin{lemma}[Long-gap sublayers]
\label{lem:long-gap-sublayers}
Let $M\in\mathfrak C(q)$.  Apart from a terminal remainder containing fewer
than $G_0(q)$ atoms, its fixed time order can be divided into consecutive
sublayers $M_j$ satisfying
\begin{equation}\label{eq:sublayer-size-gap}
 G_0(q)\le|M_j|<2G_0(q),
 \qquad
 |t_a-t_b|\ge\epstar
 \quad(a\in M_j, b\in M_{j+1}).
\end{equation}
Every physical particle line intersects $M_j$ in an interval.  There are at
least two complete sublayers.
\end{lemma}

\begin{proof}
Starting at one end of the fixed time order, put the next $G_0(q)$ atoms into
the current sublayer and continue adding atoms while the next consecutive
time gap is smaller than $\epstar$.  If this extension added another
$G_0(q)$ atoms, those added collisions would occur in an interval $J$ of
 length at most $G_0(q)\epstar$.  Split their collision graph into nontrivial
connected components; isolated lines contribute no added collision.  If one
such component has $q'$ lines, then
\[
 (q'-1)\bigl(\eps+2|J|V_\eps\bigr)
 \le q\bigl(\eps+2G_0(q)\epstar V_\eps\bigr)=o(1).
\]
For small $\eps$ this is below the injectivity radius of the torus.
Lemma~\ref{lem:short-window-component-lifting} therefore makes the entire
component, including all off-tree contacts, a genuine Euclidean hard-ball
history on $J$.  Different components may use different integral
translations because no collision constraint joins them.

The Burago--Ferleger--Kononenko Euclidean equal-mass/equal-radius bound
\cite[Theorem~1.3, p.~699]{BFK1998} gives at most
$(32(q')^{3/2})^{(q')^2}$ collisions in that component.  Since $q'\ge2$,
this is no more than
$\binom{q'}2(32(q')^{3/2})^{(q')^2}$.  There
are at most $q$ components and $q'\le q$, so the total number of added
collisions is at most
\[
 q\binom q2(32q^{3/2})^{q^2}<G_0(q),
\]
contrary to the assumed extension.  Hence a complete sublayer has fewer
than $2G_0(q)$ atoms.  A boundary is created only when the extension stops,
namely at a gap at least $\epstar$, which proves the second assertion in
\eqref{eq:sublayer-size-gap}.  No short-gap assertion is made
inside the initial block of $G_0(q)$ atoms.

Atom membership along one physical particle line is monotone in the fixed
time order, so the interval property follows.  Finally
 $\widehat G(q)=(2q+1)^5\widehat G_0(q)$, while
 $G_0(q)\le2\widehat G_0(q)$ for $q\ge2$.  Hence
 $G(q)>4G_0(q)$.  Consequently, after
discarding one terminal remainder, at least two complete sublayers remain.
The case $q=1$ is empty because a collision atom contains two distinct
physical particle lines.
\end{proof}

\begin{lemma}[Dimension-free one-sided cleanup]
\label{lem:dimension-free-up-down}
Let a connected C-molecule have no top fixed end and apply the finite
\textbf{UP} deletion rule: first remove degree-two atoms; then choose a
lowest degree-three atom (or a lowest atom if only degree-four atoms remain),
and in its descendant set remove the highest remaining atom, paired with an
adjacent degree-three parent or child exactly when such a pair exists.
Then the elementary outputs contain no $\{33{\rm B}\}$ or $\{44\}$
component, at most one degree-four component, and a degree-four component
occurs if and only if the input is full.  The time-reversed assertion holds
for \textbf{DOWN} when there is no bottom fixed end.  This statement holds in
every dimension.
\end{lemma}

\begin{proof}
The rule uses only atom degree, the parent--child order and whether an end is
fixed.  A two-atom output is created only when the chosen atom and its
adjacent parent or child both have degree three; the chronological choice
places their two exterior ends on the type-A sides, so a $\{33{\rm B}\}$
pair is never produced.  Degree-four atoms are never paired, excluding
$\{44\}$.  Once a singleton degree-four atom is removed, every remaining
connected component receives from that cut a fixed end on the processed
side; the same invariant holds after every later removal.  Thus there is at
most one degree-four singleton, and it can occur only before any fixed end is
present, namely when the original component is full.  Conversely, if the
input is full, repeated degree-two and degree-three removals either end in
the unique degree-four singleton or create one at the first stage at which
only degree-four atoms remain.  This proves the assertion by induction on
the number of atoms.  Reversing the parent--child order proves the DOWN
case.  No step refers to a velocity space or to $d$.
\end{proof}

\begin{lemma}[Packet induction for the special class]
\label{lem:packet-special-induction}
For every connected $M\in\mathfrak C(q)$ arising in
\eqref{eq:application-size-bound}, there is a cutting sequence and a
connected protected full component $K\subset M$ such that:
\begin{enumerate}[label=\textup{(\roman*)}]
\item $K=K_U\cup K_D$ satisfies all hypotheses of
Theorem~\ref{thm:two-landing};
\item $|K|<4G_0(q)\le P_{\rm pkt}(d)$;
\item every component of $M\setminus K$ is nonfull and is cut by the
ordinary UP or DOWN algorithm into $\{2\}$-, $\{3\}$- and
$\{33{\rm A}\}$-components, including their certified good recuts;
\item $K$ is the unique full component produced by the sequence;
\item the component order satisfies the dependency condition preceding
\eqref{eq:source-iterated-integral}.
\end{enumerate}
The extra positive decomposition has at most $C^{|M|}\cL_\eps^C$ cells.
\end{lemma}

\begin{proof}
We induct on $q$.  The classes $\mathfrak C(1)$ and $\mathfrak C(2)$ are
empty.  Indeed, with two physical particle lines, two consecutive collision
atoms are joined along both lines and form a double bond; without a double
bond there is at most one atom, whereas $G(2)>1$.

Let $q\ge3$ and form the complete sublayers of
Lemma~\ref{lem:long-gap-sublayers}.  First suppose some complete sublayer is
disconnected or misses a physical particle line.  It has at most $q$
connected components: the particle-line interval property prevents one
line from meeting two components of the same sublayer, so choosing one line
from each component is injective.  Choose a largest component $B$.  It meets $q'<q$
lines, and
\begin{equation}\label{eq:G-induction-inequality}
 |B|\ge q^{-1}G_0(q)>G(q-1)\ge G(q').
\end{equation}
division of the two real envelopes in
\eqref{eq:G-G0} gives
\begin{align*}
 \frac{q^{-1}\widehat G_0(q)}{\widehat G(q-1)}
 &=\frac1{q-2}\frac{(2q+1)^5}{(2q-1)^{10}}
 32^{2q-1}\frac{q^{3q^2/2}}{(q-1)^{3(q-1)^2/2}}\\
 &>4.
\end{align*}
for $q\ge3$; the last inequality follows already after replacing the final
ratio by $q^{3(2q-1)/2}$ and
$(2q+1)^5/(2q-1)^{10}$ by $(2q-1)^{-5}$.
Since $G_0(q)\ge\widehat G_0(q)$ and
$G(q-1)\le2\widehat G(q-1)$, this proves
$q^{-1}G_0(q)>G(q-1)$ as used in
\eqref{eq:G-induction-inequality}.
Lemma~\ref{lem:hereditary-special-class} applies to $B$.  Cut $B$ as free
and use the induction hypothesis.  Protect the returned packet $K$ and all
ordinary components in $B\setminus K$.  The hypotheses of
Lemma~\ref{lem:protected-band-restart} are now satisfied.

 Let $B_-$ be the union of those connected components of the union of all
 strictly earlier sublayers that have a molecule bond to $B$.  This is the
 set denoted $\mathcal M''$ in the first case of Source A's special-case
 induction.  Time convexity on every particle line has two consequences.
 First, every boundary bond of $B_-$ created by cutting $B$ is on its top
 side; hence $B_-$ has no bottom fixed end.  Second, a component of the
 earlier-sublayer union not contained in $B_-$ has no bond to $B_-$: such a
 bond would join it to the same connected component used in the definition
 of $B_-$.  Thus every component of $B_-$ receives a fixed end from $B$ and
 is nonfull.  Apply DOWN to all of them.

 Now consider a component $C$ of the remaining molecule.  Connectivity of
 the original molecule gives a first bond on a path from $C$ to
 $B\cup B_-$.  It becomes a fixed end of $C$, so $C$ is nonfull.  Such a
 boundary bond cannot point from a later processed atom back into an
 unprocessed earlier component: all earlier-sublayer components attached to
 $B$ were included in $B_-$, and components connected to them inside the
 earlier union were included at the same time.  Hence the remaining molecule
 has no top fixed end.  Apply UP componentwise.
 Lemma~\ref{lem:dimension-free-up-down} gives no $\{33{\rm B}\}$ or
 $\{44\}$ component and gives a degree-four component if and only if its
 input component is full.  Thus this
 cleanup contains no degree-four component and proves (iii)--(v); the last
assertion is also the dependency conclusion of
Lemma~\ref{lem:protected-band-restart}.

It remains to consider the case in which every complete sublayer is
connected and meets every physical particle line.  Choose two consecutive
complete sublayers at the initial end, opposite the terminal remainder, and call their
union $K=K_U\cup K_D$.  The interval property gives exactly one crossing
bond on each line, and the long boundary gap gives the separation in
Theorem~\ref{thm:two-landing}.  Each sublayer is connected and meets every
line.  Since the no-double hypothesis excludes the two-line case, at least
three lines occur.  Hence Theorem~\ref{thm:two-landing} applies.  Also
$|K|<4G_0(q)$.

Cut $K$ as free and protect it.  Every component of the one-sided remainder
receives a fixed end from this cut and has no fixed end on the opposite
side, so UP or DOWN cuts it into the components in (iii), without a
degree-four component.  The cut makes $K$ full and makes every complement
component nonfull, proving uniqueness.  The time-reversed and forward
forms of Lemma~\ref{lem:protected-band-restart} give the required shadow
dependency order.  All decompositions made here are
time-order, sign, lift or ordinary good/bad partitions.  By
\eqref{eq:application-size-bound} their number is
 $C^{|M|}\cL_\eps^C$ by
Lemma~\ref{lem:packet-size-partition-ledger}.  This closes the induction.
\end{proof}

\begin{lemma}[Dimension-free transversal cutting counts]
\label{lem:dimension-free-transversal-counts}
Let $A$ be a connected transversal atom set in a connected full C-molecule with
no double bond, and use the graph-theoretic TRANSUP and MAINTRUP deletion
rules with the sets $X_0^+(A)$ and $Y_1$.  Their output counts satisfy
\begin{align}
 \textup{TRANSUP:}\quad&
 N_{33{\rm B}}=N_{44}=0,\qquad N_4=1,\qquad
 N_{33{\rm A}}\ge \tfrac12|X_0^+(A)|-\rho(A)-1,
 \label{eq:dimension-free-transup}\\
 \textup{MAINTRUP--1:}\quad&
 N_{33{\rm B}}=N_{44}=0,\qquad N_4=1,\qquad
 N_{\rm good}\ge\tfrac1{10}|Y_1|-\rho(A)-1,
 \label{eq:dimension-free-maintrup-one}\\
 \textup{MAINTRUP--2:}\quad&
 N_{44}=0,\qquad N_{33{\rm B}}=N_{33{\rm B},{\rm good}},
 \notag\\[-2mm]
 &N_{33{\rm A}}+N_{33{\rm B}}
 \ge\tfrac1{10}\bigl(q_c-10^5E\bigr),\qquad N_4\le E.
 \label{eq:dimension-free-maintrup-two}
\end{align}
where $E=|Y_1|+\rho(A)+|X_0^+(A)|$ and $q_c$ is the number of
$A$--$X^+(A)$ bonds.  The same statements hold after reversing top and
bottom.  These are graph counts and hold in every dimension.
\end{lemma}

\begin{proof}
For TRANSUP, the causal fixed-end invariant in
Lemma~\ref{lem:dimension-free-up-down} excludes type B and $\{44\}$ outputs
and leaves exactly the initial full singleton.  If $q_2$ atoms of $A$ end in
degree-two components, invariance under a free cut of
$|E_*(A)|-3|A|$ gives $q_2\le\rho(A)$.  For each
$x\in X_0^+(A)$, mark the later-cut of its two neighbors in $A$.  One
neighbor receives at most two marks.  Except for the unique degree-four
atom and the $q_2$ degree-two atoms, a marked neighbor is paired by the rule
into a $\{33{\rm A}\}$ component.  This proves
\eqref{eq:dimension-free-transup}.

In MAINTRUP--1 the same fixed-end invariant gives the first three counts.
For these component counts the exact incidence identity is
$N_2+N_{33{\rm A}}=\rho(A)$.  Each remaining weakly-degenerate singleton, or the
later member of a weakly-degenerate pair, is a certified good component;
one output is charged by at most ten members of $Y_1$.  Removing the unique
degree-four output proves \eqref{eq:dimension-free-maintrup-one}.

For MAINTRUP--2, if $E=0$, then $X_0^+(A)=\varnothing$ and the closure
recursion gives $X^+(A)=\varnothing$, so $q_c=0$ and the claim is trivial.
Assume henceforth that $E\ge1$.  First remove $X_0^+(A)$ and the induced degree-two cascade,
and then join the components of $Y_1$ and of a cycle-breaking set by paths
of length at most four until their mutual distance exceeds four.  If $Z_0$
is the cascade, the incidence inequality gives $|Z_0|\le3E$, while the
path-joining construction gives $|S|\le200E$.  Thus at most
$2(|X_0^+(A)|+|Z_0|+|S|)\le408E$ crossing bonds are deleted, at most $E$
full singletons are created, and the remaining forest has at most
$1+4(|X_0^+(A)|+|Z_0|+|S|)\le817E$ components when $E\ge1$.
It is proper and has no uncertified weak degeneracy.  Greedily pair the
surviving $A$--$X^+(A)$ connections along this forest.  One pair uses at
most ten connections and at most one connection per forest component is
left over.  The combined loss is at most $1225E<10^5E$.  Thus at
least $(q_c-10^5E)/10$ type-A or type-B pairs are produced.  A type-B pair
lies wholly in the proper, non-weakly-degenerate part of $A$ and is therefore
on its certified good cell.  This proves
\eqref{eq:dimension-free-maintrup-two}.  Every argument used only incidence,
order and fixed-end data.
\end{proof}

\begin{lemma}[All-dimensional Main-case-2 count]
\label{lem:main-case-2-d}
Fix $d\ge4$ and assume that the preliminary strongly-degenerate
primitive-pair branch has been removed.  Thus every $\{33{\rm A}\}$ output
is either already good or is positively recut into one degree-two atom and
one good degree-three atom, without increasing $N_4$.  Use Source A's
notation for the connected transversal
sets $A_j$, the sets $X(A_j)$, $X_0^+(A)$ and $Y_1$, the cycle rank
$\rho$, and the crossing-bond count $q_c=\#_{\rm conn}^+(A)$.  Put
\begin{equation}\label{eq:main-case-2-cutoffs}
 D=K_1:=(60d)^{60d},\qquad
 K_{j+1}:=\bigl(60d\,G((60dK_j)^{60d})\bigr)^{60d},
 \qquad
 \Gamma>(60d)^{60d}K_D .
\end{equation}
The dimension-free cutting constructions below produce alternative~\textup{(2)} of
Theorem~\ref{thm:p38j}; more precisely,
\begin{equation}\label{eq:main-case-2-target}
 \frac{N_{\rm good}}{10d}-dN_4>100d^2.
\end{equation}
No ordinary-component estimate, packet estimate, or dimension-dependent
collision integral is used in this count.
\end{lemma}

\begin{proof}
We rerun the two exhaustive graph branches.  Suppose first that $A_j$ exists
and $X_0(A_j)\ne\varnothing$ for every $j\le D$.  Put
$A_0=\varnothing$, so that $X(A_0)=\varnothing$, and
\[
 B_j:=A_j\setminus(A_{j-1}\cup X(A_{j-1})).
\]
The rank-preserving construction of $A_j$ orders
$B_j=(a_1^j,\ldots,a_{m_j}^j)$ so that $a_i^j$ has exactly one bond to
$A_{j-1}\cup X(A_{j-1})\cup\{a_1^j,\ldots,a_{i-1}^j\}$.  Choose
$x_j\in X_0(A_j)$.  At least one of its two neighbors in $A_j$ lies in
$B_j$, since otherwise $x_j$ would already belong to $X(A_{j-1})$.

Starting with $j=1$, cut the atoms $a_i^j$ in increasing order.  Whenever
$a_i^j$ is adjacent to $x_j$ and $x_j$ has degree three at that step, cut
$\{a_i^j,x_j\}$ instead of $a_i^j$; after $B_j$ is exhausted, cut the remaining
atoms of $X(A_j)$ from the lowest ones, which then have degree two.  Proceed
only afterward to $B_{j+1}$.  Thus the cuts for distinct $j$ use disjoint new
blocks, and each $x_j$ belongs to a distinct $\{33{\rm A}\}$ component:
if $x_j\in X_0^+(A_j)$ its two exterior ends are both on top, and if
$x_j\in X_0^-(A_j)$ they are both on bottom.  The preliminary recut
assumption turns each such component into one good output without increasing
$N_4$.  The fixed-end invariant leaves only the initial full
degree-four component.  Hence
\[
 N_{\rm good}\ge D,
 \qquad N_4=1,
 \qquad
 \frac{N_{\rm good}}{10d}-dN_4
 \ge \frac{(60d)^{60d}}{10d}-d>100d^2.
\]

Otherwise the construction reaches the source's terminal-rank branch.
Since $\rho(A_1)=0<K_1$ while the full molecule has
$\rho(\mathcal M)=\gamma>\Gamma>K_D$, there is a least $j<D$ for which
$\rho(A_j)<K_j$ and $\rho(A_{j+1})\ge K_{j+1}$.
For the least relevant index $j$, set $A=A_j$.  The cycle-rank Euler identity
and the definition of the closure $X(A)$ give
\begin{equation}\label{eq:main-case-2-rank-step}
 \rho(A)<K_j,\qquad
 K_{j+1}\le \rho(A\cup X(A))\le \rho(A)+10|X(A)|,
\end{equation}
and, after exchanging top and bottom if necessary,
$|X^+(A)|\ge K_{j+1}/40$.  The negation of Main case 1 and the definition
of $K_{j+1}$ then give
\begin{equation}\label{eq:main-case-2-qc}
 q_c>(60d)^{60d}K_j.
\end{equation}
Let
\[
 T:=(30d)^{30d}K_j.
\]
The three alternatives in
Lemma~\ref{lem:dimension-free-transversal-counts} make the following
exhaustive choice.

If $|X_0^+(A)|\ge T$, \eqref{eq:dimension-free-transup} gives
\[
 N_{\rm good}\ge \frac{|X_0^+(A)|}{2}-\rho(A)-1
 >\frac T2-K_j-1,\qquad N_4=1.
\]
If $|X_0^+(A)|<T$ but $|Y_1|\ge T$,
\eqref{eq:dimension-free-maintrup-one} gives
\[
 N_{\rm good}\ge \frac{|Y_1|}{10}-\rho(A)-1
 >\frac T{10}-K_j-1,\qquad N_4=1.
\]
Since $K_j\ge K_1=(60d)^{60d}$, either lower bound implies
\eqref{eq:main-case-2-target} for $d\ge4$.

It remains that $|X_0^+(A)|<T$ and $|Y_1|<T$.  Put
\[
 E:=|Y_1|+\rho(A)+|X_0^+(A)|<2T+K_j\le3T.
\]
Equation~\eqref{eq:dimension-free-maintrup-two} gives, with every
$\{33{\rm B}\}$ component good; the preliminary recut condition makes each
counted $\{33{\rm A}\}$ contribute one good output without increasing $N_4$.
Therefore
\begin{equation}\label{eq:main-case-2-option-2}
 N_{\rm good}\ge\frac{q_c-10^5E}{10},
 \qquad N_4\le E.
\end{equation}
If $E=0$, \eqref{eq:main-case-2-qc} and $K_j\ge K_1$ immediately imply
\eqref{eq:main-case-2-target}.  If $E\ge1$, then
\begin{equation}\label{eq:main-case-2-ratio}
 q_c>R_dE,\qquad
 R_d:=\frac{(60d)^{60d}}{3(30d)^{30d}}
      =\frac{(120d)^{30d}}3.
\end{equation}
Consequently
\[
 \frac{N_{\rm good}}{10d}-dN_4
 >\left(\frac{R_d-10^5}{100d}-d\right)E>100d^2,
\]
where the last elementary inequality already holds at $d=4$ and its
left-hand side is increasing for integer $d\ge4$.  This proves
\eqref{eq:main-case-2-target} in every branch.  The argument used only the
source's structural rank and cutting counts, so it is independent of the
analytic upgrade of the elementary components.
\end{proof}

\subsection{The full Proposition 3.8 replacement}

\begin{theorem}[Closed-endpoint joint alternative
$\mathrm{P3.8\mbox{-}J}$]
\label{thm:p38j}
Fix an integer $\Gamma>(60d)^{60d}K_D$, with $D$ and $K_D$ defined in
\eqref{eq:main-case-2-cutoffs}.  Let $\mathcal M$ be a full one-layer
C-molecule whose bonds form a tree of
at most $|\log\eps|^{C^*}$ atoms plus $\gamma$ additional bonds, with
$\Gamma<\gamma\le2\Gamma$, and suppose that the remaining graph-theoretic
input hypotheses of Source A, Proposition~3.8 (fullness, one-layer order and
the tree-plus-excess-bond class) hold.  No analytic conclusion or ambient
dimension restriction is imported from that proposition.  After at most
$C^{|\mathcal M|}\cL_\eps^{C^*}$ positive indicator pieces, one of Source
A's alternatives (1), (2), or a locally certified instance of (3) holds, or
there is a protected packet $K$ such that
\begin{enumerate}[label=\textup{(\roman*)}]
\item $K$ is the unique full component and $|K|\le P_{\rm pkt}(d)$;
\item every component of $\mathcal M\setminus K$ is nonfull and has an
ordinary UP/DOWN cleanup with no degree-four component;
\item the packet and cleanup are one admissible cutting sequence for
\eqref{eq:source-cutting-identity}, ordered with the dependency property
preceding \eqref{eq:source-iterated-integral};
\item for the complete positive kernel $Q$ of
\eqref{eq:source-positive-Q}, with the original time indicator retained,
\begin{equation}\label{eq:p38j-packet}
 J_K(Q)\le\cL_\eps^C\eps^{-(d-1)}\eps^{2(d-1)}\epstar^{-2(d-1)}
 \int Q_K^\sharp.
\end{equation}
\end{enumerate}
This fourth outcome is an operator replacement for the Source-B-v2 top
consumer; it is not a componentwise factorization into at most ten excesses.
\end{theorem}

\begin{proof}
We first justify the closed upper endpoint without extending the printed
Source-A proposition.  The parameter order
\eqref{eq:appE-choice-order} fixes an integer Source-B threshold
\(\Gamma\) and, at the same stage, the Source-A threshold
\[
 \widehat\Gamma:=\Gamma+\frac12 .
\]
The recollision number \(\gamma\) is an integer.  Hence
\begin{equation}\label{eq:strict-endpoint-embedding}
 \Gamma<\gamma\le2\Gamma
 \quad\Longrightarrow\quad
 \widehat\Gamma<\gamma<2\widehat\Gamma .
\end{equation}
We use the parameter \(\widehat\Gamma\) to run the same finite graph
partition, including when \(\gamma=2\Gamma\), but do not import its analytic
conclusion from the source's $d\in\{2,3\}$ ambient theorem.  The branches of
that partition are rerun below: the one-sided cleanup is
Lemma~\ref{lem:dimension-free-up-down}, Main case~2 is
Lemma~\ref{lem:main-case-2-d}, and every elementary analytic output is
supplied by Proposition~\ref{prop:ordinary-subset} or by the packet theorem.
The constant \(\Gamma\) is chosen large enough that \(\widehat\Gamma\)
satisfies the fixed graph threshold, and every graph/order constant is
enlarged once to the maximum of the constants at \(\Gamma\) and
\(\widehat\Gamma\).  Independently, the incidence count
\(q/2\le\gamma\le2\Gamma\) still gives \(q\le4\Gamma\), so the packet cap
\eqref{eq:C4-packet-cap} is unchanged.  Thus no equality case is inferred
from an inspection of a strict source proof.

A strongly degenerate primitive pair gives alternative (1), with its
$\{44\}$ close-state gain supplied by \eqref{eq:44-bound}.  The two
double-bond graph branches are detected by the same degree/order test, which
is dimension-free; their $d$-dimensional analytic estimates are respectively
the double capacity \eqref{eq:double-periodic-capacity} and
the ordinary two-atom estimates in Proposition~\ref{prop:ordinary-subset}.
Their excess factors are at most
$\eps^{d-1/2}$ or the product
$\eps^{d-1}\eps^{1/(8d)}$, both stronger than
$\eps^{d-1+1/(15d)}$.  These are locally certified instances of alternative
(3).

After those branches, the molecule has no double bond and no strongly
degenerate primitive pair.  We do not invoke Source A's separate
long-triangle branch: its serial $d$-dimensional output is precisely the
componentwise estimate replaced by the joint packet.  Long triangles remain
inside the special class, whose hereditary and restart lemmas do not assume
their absence.

In Source A's Main case 1, let $B$ be a largest connected component of
$X^+(A)$, with ties resolved by the fixed atom order, and let $q_B$ be the
number of its bottom bonds to $A$.  If $c$ is the number of components of
$X^+(A)$, then $c\le q/2$, where $q=\#_{\rm conn}^+(A)$.  If $c>1$, every
other component contains an $X_0^+(A)$ atom and hence contributes at least
two of the $q$ cross bonds, so $q_B\le q-2$.  Moreover
\[
 |B|\ge \frac{|X^+(A)|}{c}\ge\frac{2G(q)}q
 >G(q-1)\ge G(q_B).
\]
Here the strict inequality follows from
\eqref{eq:G-induction-inequality}; if $c=1$, instead $B=X^+(A)$,
$q_B=q$, and the defining Main-case-1 bound gives $|B|\ge G(q_B)$.
The function $G$ is increasing, and the child-closure in the definition of
$X^+(A)$ makes the bottom end of every physical line of $B$ one of these
$B$--$A$ bonds; thus $q_B$ also bounds the number of lines of $B$.
After cutting $B$ as free, it is a connected member of
$\mathfrak C(q_B)$: its bottom free ends are precisely among those $q_B$
bonds, and physical realizability, the no-double condition and the absence
of strongly degenerate primitive pairs are inherited.  Use
Lemma~\ref{lem:packet-special-induction} on $B$ in place of Source A,
Proposition 4.4.  The returned packet $K$ is full, its complement inside
$B$ is nonfull, and $|K|<4G_0(q_B)\le P_{\rm pkt}(d)$.

We spell out the subsequent order.  The connected transversal set $A$
receives a fixed end from the cut of $B$ and has no bottom fixed end, so it is
processed by \textbf{DOWN}.  Let $U$ be a component of the unprocessed part
of $A^+$, including a component of $X^+(A)\setminus B$.  If $U$ has a bond
to $A$, the cut of $A$ supplies its fixed end.  Otherwise a path in its
original $A^+$ component to $A$ must cross $B$, so the earlier cut of $B$
already supplied a fixed end to $U$.  The child-closure defining $X^+(A)$
makes this a bottom fixed end, while transversality leaves no top fixed end;
thus \textbf{UP} applies in either case.  Finally each component of $A^-$
receives a fixed end from the processed set $A\cup A^+$ and has no bottom
fixed end, so \textbf{DOWN} applies.
Lemma~\ref{lem:dimension-free-up-down} gives no full degree-four output from
any of these nonfull inputs.
The paired free/fixed shadows point from $B$ to $A$ and to the separated
$A^+$ components, from $A$ to the remaining $A^+$ components, and from
$A$ or $A^+$ to $A^-$.  Hence the complete order also satisfies item (iii).  Thus no second
full component and no additional degree-four baseline is created.

Source A's Main case 2 never calls the special-case proposition.
Lemma~\ref{lem:main-case-2-d} reruns every TRANSUP and MAINTRUP branch with
the dimension-dependent target and gives alternative~(2).
All source indicator partitions are retained, and the new cell count in
Lemma~\ref{lem:packet-size-partition-ledger} is within the permitted factor.
Finally \eqref{eq:p38j-packet} is
Theorem~\ref{thm:two-landing} with $|K|\le P_{\rm pkt}(d)$.  This proves every
branch.
\end{proof}

Since $\epstar=\exp(-\sqrt{|\log\eps|})$, for every $\eta>0$,
\begin{equation}\label{eq:p38j-subpower}
 \cL_\eps^C\epstar^{-2(d-1)}\le\eps^{-\eta}
\end{equation}
once $\eps$ is sufficiently small.  Taking $\eta=1/(100d)$ gives
\begin{equation}\label{eq:packet-excess-d}
 \sigma_{K,d}^{\rm sharp}:=\cL_\eps^C\eps^{2(d-1)}\epstar^{-2(d-1)}
 \le \eps^{2(d-1)-1/(100d)}.
\end{equation}
The normalized packet estimate itself also contains the unique factor
$\eps^{-(d-1)}$; consequently it is neither
claimed nor needed to imply Source A's componentwise third alternative
\eqref{eq:source-A-third-alternative}.  Its role is instead to replace that
branch directly at the two Source-B operator estimates.  We now prove
this downstream statement, including the baseline and time ledgers.

\subsection{The fibrewise time owner}

We first make the new-line count behind the time owner literal.  Let
$\mathcal P_{\rm top}$ be the physical particle lines in the exceptional
top cluster.  Thus $|\mathcal P_{\rm top}|=r+m$.  If $r>0$, let
$\mathcal R_{\rm top}$ be its $r$ distinguished root lines; if $r=0$,
defer to the matching construction the choice of one line as the common-
translation representative and let
$\mathcal R_{\rm top}$ consist of that line.  In both cases
\begin{equation}\label{eq:top-new-line-set}
 \mathcal N_{\rm top}:=
 \mathcal P_{\rm top}\setminus\mathcal R_{\rm top},
 \qquad
 |\mathcal N_{\rm top}|=m-\iota .
\end{equation}
Write $\mathcal P(K)$ for the physical lines meeting the protected packet
$K$ and define
\begin{equation}\label{eq:PK-definition}
 P_K:=|\mathcal P(K)|.
\end{equation}
The upper physical-particle collision graph of $K$ has the reverse-
chronological Kruskal spanning tree from
Lemma~\ref{lem:upper-tree-landings}.  Denote it by $T_K$, retain the time of
each of its $P_K-1$ contact edges, and call this set of literal time
coordinates $R_{\rm birth}$.  For $p\in\mathcal N_{\rm top}$ let $b(p)$
denote its unique creation atom in the actual upper source birth forest.

The lower molecule has its source new-line set
$\mathcal N_{\rm low}$, with
$|\mathcal N_{\rm low}|=\mathsf X$.  Before making any local refinement,
retain the actual birth-forest creation atom of every line in
$\mathcal N_{\rm low}$ and of every
$p\in\mathcal N_{\rm top}$ whose creation atom $b(p)$ lies outside $K$.
Call this atom set $\mathfrak F_{\rm birth}$.  Distinct new lines have
distinct creation atoms, even when two such atoms later belong to the same
two-atom elementary component.

We separately select the atoms which will receive a one-atom joint
owner/good estimate, and one bookkeeping representative from each periodic
certificate which carries a birth time.  Group $\mathfrak F_{\rm birth}$ by
immutable identifier.  Only ordinary support-good one-atom $\{3\}$ and
full-$\{4\}$ identifiers are eligible for the joint estimate.  In the fixed
reverse forest order choose their unique atom.  For each periodic
certificate choose at most one of its birth atoms only to identify the
certificate's once-only tube; all of its times remain in the regular
projection.  Let $\mathfrak F_{3,\rm own}$ be the selected ordinary
degree-three family, $\mathfrak F_{4,\rm full,own}$ the selected ordinary
full-$\{4\}$ family, and $\mathfrak F_{\rm per,own}$ the periodic
representatives,
and put
\[
 \mathfrak F_{\rm own}:=
 \mathfrak F_{3,\rm own}\,\dot\cup\,
 \mathfrak F_{4,\rm full,own}\,\dot\cup\,
 \mathfrak F_{\rm per,own}\subset\mathfrak F_{\rm birth}.
\]
All other birth atoms remain literal time owners but are paid by the regular
existential projection; they do not receive a second good or tube factor.
In particular every birth atom lying in any ordinary two-atom
$\{33{\rm A}\}$, $\{33{\rm B}\}$ or $\{44\}$ component is left in this
regular class; no unsupported cut of such a component is made.  Outside a
type-$\{33{\rm A}\}$ short-gap cell, Lemma
~\ref{lem:ordinary-fixed-distinguished-times} keeps all of those times as
genuine outputs while obtaining every complementary local estimate.  On the
short-gap cell the original two-time block is instead retained intact until
the post-envelope common-measure step; neither time is first conditioned out.
For a selected periodic good-$\{4\}$ its $\eps^{-(d-1)}$ baseline remains in
$N_4$, and for a selected ordinary full-$\{4\}$ the same unique baseline is
retained.  Define
\begin{equation}\label{eq:S0-definition}
 S_0:=\{t_a:a\in\mathfrak F_{\rm birth}\}.
\end{equation}
Thus $S_0$ is a set of literal, pairwise distinct creation times, not a set
of components.  The selected subfamily $\mathfrak F_{\rm own}$ contains at
most one atom from each good identifier; $S_0$ itself can contain two times
from one $\{33{\rm A}\}$ component.  Neither degree-four class is recut.  The arbitrary source
kernel is not integrated or majorized at this matching stage.  After the
complete source envelope has been formed, Lemma
~\ref{lem:owner-marked-good-component} treats the ordinary support-good
owners, including ordinary good-$\{4\}$ owners relative to their retained
baselines, on their original conditional measures.  Periodic
certificate owners are kept separate and retain their ordinary projected
time coordinate.  In particular, no component is reclassified and no
$N_{\rm good}$ credit is deleted.

\begin{lemma}[Retention of the original birth charts]
\label{lem:birth-chart-retention}
Fix the actual source birth forest and any set $B$ of its creation atoms.
The positive cutting and shadow substitutions may be ordered so that, for
every $b(p)\in B$, the original atom time $t_{b(p)}$ remains a genuine
Lebesgue output coordinate until the distinguished-time projection, except
that the atom times of a type-$\{33{\rm A}\}$ short-gap block may be consumed
together only after the complete pointwise source majorant has been formed.
This
remains true if the current representative of $b(p)$ would otherwise be a
standalone degree-two root substitution.  Distinct new lines use distinct
blocks in one square global birth-velocity coordinate system; after a
C-atom the \emph{entire colliding pair}, rather than either particle alone,
is transported by the full orthogonal scattering matrix.  Consequently all
elementary estimates outside the declared short-gap exception may be applied
with the times in $B$ held fixed.  On the exception the same literal time
coordinates remain on one common Radon block, and no estimate consumes a
time which is later counted again by the projection.
\end{lemma}

\begin{proof}
Work before any elementary root substitution.  Root each source-birth tree
at an already present line (or at its declared common-translation
representative).  At the creation edge $b(p)$ the parent state is already
present and the new physical line has its independent $d$-dimensional birth
variable $v_p=h_p$ (already included in the global input fibre).
The exact collision chart is therefore the one-fixed-line chart
\eqref{eq:degree-three-chart}, with coordinates
$(t_{b(p)},\omega_{b(p)},v_p)$; for the initial full component the same
statement holds after its declared root line is used, while its unique full
baseline is left outside.  The equality used here is the branchwise Radon
identity \eqref{eq:exact-transparent-birth-relay}, with the complete
complement and $\mathbf1_{\mathcal D_M}$ inside its arbitrary test
function.  The contact-sphere Jacobian together with the explicit flux
factor in \eqref{eq:source-C-kernel} gives the displayed chart density, and
the original time coordinate is retained pointwise; no marginal time
argument or new determinant is introduced.

Use the fixed $dn$-dimensional input fibre
\eqref{eq:exact-global-input-fibre}: all root velocities and all independent
birth innovations $h_p$ are coordinates from the outset.  A future $h_p$ is
carried by the identity before its creation event; it is never represented
by a zero particle slot.  Root inputs may subsequently be conditioned on as
parameters, but this conditioning does not change the dimension of the
global square velocity map.

Suppose the birth of $p$ occurs by collision with the already active parent
$q$ and normal $\omega$.  At that event the embedded complete-pair matrix
\eqref{eq:exact-global-C-factor} acts on $(v_q^-,h_p)$.  Conditional on the
old inputs, its $h_p$-column is the $d$-dimensional isometric injection
\begin{equation}\label{eq:birth-mixed-isometry}
 J_{q,p,\omega}h_p
 :=\mathscr R_{qp}^{(2d)}\binom{0}{h_p}
 =\binom{P_\omega h_p}{(I-P_\omega)h_p},
 \qquad J_{q,p,\omega}^*J_{q,p,\omega}=I_d.
\end{equation}
The first component is a parent-line \emph{outgoing} direction and is not a
fixed root input.  Old variations enter the same pair as
$\mathscr R_{qp}^{(2d)}(\delta v_q,0)$, orthogonally to the new column.
Equation~\eqref{eq:birth-mixed-isometry} is only a $d$-to-$2d$ column
injection; no square determinant is attributed to it.  The determinant-one
statement belongs instead to the complete $2d$-dimensional event factor and,
globally, to the ordered product of its embedded $dn$-dimensional copies in
\eqref{eq:exact-global-factorization}.  At an O-atom the full vector is
unchanged, and every later C-atom again acts on the complete active pair.
This also covers repeated collisions on the same line and never uses the
rank-$(d-1)$ single-particle block $I-P_\omega$ as a pivot.

The transported range of the $h_p$ column is a legitimate mixed
$d$-dimensional innovation subspace, but it is \emph{not}, in general, a
$d$-dimensional particle-velocity output of the current Source-B
component: one collision may put its normal part on the parent line and its
tangential part on the child line.  We retain the input coefficient $h_p$
but never consume it separately from the complete downstream
free-coordinate system.  A later consumer using either transported partner
is resolved first, and inverse complete-pair scattering recovers both
partner pre-states before the earlier birth chart is reached.  Thus
$(t_{b(p)},\omega_{b(p)},h_p)$ remains a chart on every fixed-root affine
fibre with unit conditional velocity Jacobian, without identifying $h_p$
with a particle slot.

If the local cleanup would call the degree-two collision-root operator on a
marked birth atom, apply
\eqref{eq:exact-transparent-birth-relay} instead of the substituted
degree-two operator.  It leaves $t_{b(p)}$ and $h_p$ as literal Lebesgue
coordinates and propagates the complete outgoing pair by
\eqref{eq:exact-global-C-factor}, whose $h_p$ column is
\eqref{eq:birth-mixed-isometry}; the source normalization, fixed/free
shadow equality, full-$\{4\}$ baseline and
$\mathbf1_{\mathcal D_M}$ are identical on both sides of that equality.
The resulting positive conditional kernel has the same $\cL_\eps^C$
cutoff/flux bound as the original birth chart.  It neither creates a
degree-two free output nor spends $h_p$.  The coefficient is assigned by
the ownership bijection of
Theorem~\ref{thm:exact-full-pair-disintegration}, together with every
innovation coupled to it by later full-pair scatterings.  This is an exact
Radon disintegration of the same positive source measure, not a claim that
the substituted Source-B degree-two operator itself has a Lebesgue time
output.  The local $h_p$-column has Gram determinant one by
\eqref{eq:birth-mixed-isometry}, while the complete square velocity
Jacobian is the one in
\eqref{eq:exact-global-factorization}--
\eqref{eq:exact-global-factor-bounds}.

When two marked birth atoms occur in the same two-atom component, retain both
times and both distinct child-velocity blocks in that component.  Outside a
type-$\{33{\rm A}\}$ short-gap cell, Lemma
~\ref{lem:ordinary-fixed-distinguished-times} proves the required ordinary
operator bounds uniformly with both times fixed.  On that short-gap cell
the original block and its indicator are retained for the later common-time
integration, without claiming a pointwise fixed-time bound.  Lemma
~\ref{lem:periodic-certificate-measure} is uniform in all fixed times for a
periodic certificate.  Unmarked degree-two roots may still be substituted
with norm at most one.  Induction over the finite reverse forest proves that
every marked time survives as asserted.
\end{proof}

\begin{lemma}[Global innovation partition and one-use law]
\label{lem:global-innovation-partition}
Fix a cutting branch, all atom times and normals, and all discrete cells.
Let \(\mathcal H_{\rm all}\) be the complete velocity fibre on the
right-hand side of the exact cutting identity
\eqref{eq:source-cutting-identity}, including the root inputs and the birth
innovations \(h_p\) of Lemma~\ref{lem:birth-chart-retention}.  There is a
finite positive disintegration, compatible with the reverse cutting order,
whose velocity pivots give a direct sum
\begin{equation}\label{eq:global-innovation-direct-sum}
 \mathcal H_{\rm all}
 =\mathcal H_{\rm ord}\oplus
   \mathcal H_{\rm per}\oplus
   \mathcal H_{\rm nor}\oplus
   \mathcal H_{\rm rem}.
\end{equation}
Here \(\mathcal H_{\rm ord}\) and \(\mathcal H_{\rm per}\) contain the velocity
parts of, respectively, all ordinary good blocks and all protected periodic
blocks; \(\dim\mathcal H_{\rm nor}=d|\mathcal F|\), where \(\mathcal F\) is
the normal degree-three family of Proposition~\ref{prop:ordinary-subset}; and
\(\mathcal H_{\rm rem}\) contains root/envelope and unused variables.  A
marked standalone degree-two birth is a transparent relay and has no
summand of its own.  Its innovation is assigned, through a complete
full-pair coordinate group, to exactly one of the four displayed summands.
Each summand is
consumed once and only once by the correspondingly named operation.

The direct sum is obtained by pulling back disjoint \emph{output coordinate
groups} through the complete square velocity transformation.  At every
C-atom the whole $2d$-dimensional colliding pair is one orthogonal block;
neither of its $d$-dimensional projections is used as an exchange pivot.
Thus no small principal angle or inverse exchange determinant occurs.  All
additional forest/centre-of-mass coordinate factors are at most
\(C^{|M|}\cL_\eps^C\).  Conditional on the other three summands, the selected
normal \emph{output coordinates} $U_{\rm nor}\simeq
(\mathbb R^d)^{|\mathcal F|}$ still range in a boundary-independent set
\(Y_{\rm rem}\subset U_{\rm nor}\) satisfying
\begin{equation}\label{eq:normal-volume-on-innovation-complement}
 |Y_{\rm rem}|
 \le \prod_{e\in E_{\rm free}(\mathcal M)}X_e^d
       C^{|M|}\cL_\eps^{C^*K}.
\end{equation}
Thus a mixed marked-birth direction is spent by its genuine downstream
owner, and no ordinary or periodic good coordinate is counted again in the
simultaneous normal volume.
\end{lemma}

\begin{proof}
Apply Theorem~\ref{thm:exact-full-pair-disintegration} to the present legal
cutting branch.  Its consumer set $\mathscr A$ is partitioned into ordinary
good, periodic good, normal degree-three and residual consumers.  For a
periodic certificate the assigned linear block is $U_C$ of
Lemma~\ref{lem:periodic-certificate-free-coordinate}; the nonlinear block
$W_C=(w_C,\zeta_C)$ is formed only afterward from the equidimensional input
$(U_C,Z_C)$ by its local conditional area formula.  Thus no angular, centre
or auxiliary variable is a row of the global velocity matrix.

Writing $u$ for the root-plus-birth input and
$w=(U_{\rm ord},U_{\rm per},U_{\rm nor},U_{\rm rem})$, the exact theorem
gives, on every finite Borel branch,
\begin{equation}\label{eq:global-innovation-square-matrix}
 D_uw=\mathfrak T
 \widehat{\mathscr U}_J\cdots
 \widehat{\mathscr U}_1\mathfrak S,
\end{equation}
where every $\widehat{\mathscr U}_j$ is a $dn$-dimensional embedded event
matrix: it is the identity except on the complete $2d$-dimensional pair of
the $j$th C-event, and is the identity at an O-event.  The product is
chronologically ordered; event factors are not put in direct sum, since two
events may use the same particle line.  Every nontrivial physical line meets
an atom, so $n\le2|M|$ (any isolated identity coordinate is left in
$U_{\rm rem}$).  After enlarging the absolute constant, this gives
\begin{equation}\label{eq:global-innovation-square-bounds}
 C^{-|M|}\le |\det D_uw|\le C^{|M|},
 \qquad
 \|(D_uw)^{\pm1}\|\le C^{|M|}\cL_\eps^C.
\end{equation}
Equation~\eqref{eq:exact-global-pushforward}, with the complete positive
kernel and $\mathbf1_{\mathcal D_M}$ as $F$, is the asserted finite positive
disintegration.  The ownership induction in that theorem says that every
row and column occurs once.  Pulling the four disjoint output groups back by
$D_uw$ proves the direct sum
\eqref{eq:global-innovation-direct-sum}.  A marked degree-two birth has no
separate row: by \eqref{eq:exact-transparent-birth-relay} its $h_p$ remains
live and the same ownership bijection puts it into exactly one downstream
group.  All subsequent volumes are taken in these output coordinates; the
single global density $|\det D_uw|^{-1}$ in
\eqref{eq:exact-global-pushforward} is paid once and is not reintroduced in
an individual consumer.

It remains only to verify that the normal consumer keeps the already proved
volume bound after the other three blocks are conditioned.  The containing
set $Y$ in \eqref{eq:source-volume}, constructed explicitly in
\eqref{eq:volume-d}, depends only on the fixed normals, the dyadic sizes and
the cutting sequence, and is uniform in every concrete fixed velocity.  For
each value of the ordinary, periodic and residual blocks, reverse
triangularity in \eqref{eq:exact-global-consumer-recursion} makes every
boundary shadow external to the normal family a function of those fixed
blocks alone.  Shadows internal to the normal family are exactly the
internal arrows already treated jointly in the construction of $Y$; no
purported fixed boundary is allowed to vary with $U_{\rm nor}$.  The
admissible section of the selected normal free velocities is therefore a
measurable subset $Y_{\rm rem}\subset Y$; conditioning cannot enlarge it.
Hence
\[
 |Y_{\rm rem}|\le |Y|
 \le \prod_{e\in E_{\rm free}(\mathcal M)}X_e^d
       C^{|M|}\cL_\eps^{C^*K},
\]
which is \eqref{eq:normal-volume-on-innovation-complement} without a new
projection inverse or a second volume argument.  Finally apply the local
ordinary and periodic area formulas blockwise as permitted by part (v) of
Theorem~\ref{thm:exact-full-pair-disintegration}; for a periodic block the
equidimensional input is $(U_C,Z_C)$, not $U_C$ alone.  Conditional Tonelli
then consumes every input block exactly once.
\end{proof}

\begin{lemma}[New-line/time-owner matching]
\label{lem:new-line-time-owner-matching}
The marked tree and the ordinary complement may be selected so that there
are maps with the following properties.
\begin{enumerate}[label=\textup{(\roman*)}]
\item There is an anchor $p_K\in\mathcal P(K)$ at which $T_K$ may be
rooted, such that either $p_K\in\mathcal R_{\rm top}$ or
$b(p_K)\notin K$.  There is a bijection
\begin{equation}\label{eq:packet-merge-bijection}
 \psi_K:\mathcal P(K)\setminus\{p_K\}
 \longrightarrow R_{\rm birth}
\end{equation}
which sends a line to the retained time of its parent edge in the rooted
chronological merge tree $T_K$.
\item Put
\begin{equation}\label{eq:NK-definition}
 \mathcal N_K
 :=\{p\in\mathcal N_{\rm top}\cap\mathcal P(K):
       b(p)\in K\}.
\end{equation}
Then $\psi_K|_{\mathcal N_K}$ is injective.  There is also an injection
\begin{equation}\label{eq:outside-owner-injection}
 \psi_0:
 \mathcal N_{\rm low}\,\dot\cup\,
       (\mathcal N_{\rm top}\setminus\mathcal N_K)
 \longrightarrow S_0 .
\end{equation}
\item The atoms carrying $\operatorname{im}\psi_0$ lie outside $K$, whereas
the atoms carrying $\operatorname{im}(\psi_K|_{\mathcal N_K})$ lie in $K$.
Thus the two images are disjoint, and
\begin{equation}\label{eq:complete-owner-cardinality}
 |S_0|+|R_{\rm birth}|
 \ge \mathsf X+m-\iota .
\end{equation}
In the attached branch $r>0$ this is
$|S_0|+|R_{\rm birth}|\ge\mathsf X+m$.  In the disjoint branch the lower
and top images separate and give respectively
\begin{equation}\label{eq:disjoint-owner-cardinality}
 |S_{0,\rm low}|\ge\mathsf X,
 \qquad
 |S_{0,\rm top}|+|R_{\rm birth}|\ge m-\iota .
\end{equation}
Thus the two disjoint subcases require and receive, respectively,
\[
 \begin{array}{c|c|c}
  &\text{lower owner count}&\text{top owner count}\\ \hline
  r>0\ (\iota=0)&\mathsf X&m\\
  r=0\ (\iota=1)&\mathsf X&m-1.
 \end{array}
\]
\end{enumerate}
\end{lemma}

\begin{proof}
First fix the source birth forest in the lower molecule.  A genuinely new
particle line has a unique creation collision in this forest.  For every
$p\in\mathcal N_{\rm low}$ put
\begin{equation}\label{eq:lower-literal-birth-owner}
 \psi_0(p):=t_{b(p)}
\end{equation}
and enter the atom $b(p)$ in $\mathfrak F_{\rm birth}$.  Distinct new lines
have distinct birth-forest edges, so these literal times are distinct.  This
assignment is made before looking at the elementary-component type and
therefore cannot lose a birth edge when a two-atom component is refined.

Now run the fixed reverse forest order on the immutable identifiers which
meet these birth atoms.  For an ordinary support-good one-atom identifier,
put its unique birth atom in the appropriate ordinary summand of
$\mathfrak F_{\rm own}$.  For a periodic certificate choose the first birth
atom only as its bookkeeping representative in
$\mathfrak F_{\rm per,own}$; every periodic birth time remains projected.
For a protected periodic good-$\{4\}$ or the
unique ordinary full-$\{4\}$ output of Source B,
Proposition~10.2(1), use the singleton without recutting and retain its
unique $\eps^{-(d-1)}$ baseline in $N_4$.  A nongood ordinary full singleton is
not selected and its time stays in the regular projection.  Every ordinary
two-atom $\{33{\rm A}\}$, $\{33{\rm B}\}$ or $\{44\}$ identifier is likewise
ineligible for a one-atom joint estimate and all of its birth times stay in
that projection, except that a type-$\{33{\rm A}\}$ short-gap cell is
reserved for the common two-time estimate after the pointwise envelope.
Lemma~\ref{lem:birth-chart-retention} preserves the literal coordinates, and
Lemma~\ref{lem:ordinary-fixed-distinguished-times} applies on all remaining
two-atom cells.  Thus neither an exact degree-two substitution nor a
fixed-time coarea estimate consumes a time which will later be projected.
Thus every new line has a time owner, while every one-atom good identifier
receives at most one joint good/owner estimate.  This constructs
the $\mathcal N_{\rm low}$ part of $\psi_0$.

For the top layer use the actual acyclic source birth forest.  If $K$ meets
a distinguished root line, take such a line as $p_K$.  Otherwise the finite
set $\mathcal P(K)$ contains a line whose creation atom is outside $K$:
indeed, following birth parents strictly backwards from any line of
$\mathcal P(K)$ must leave the finite set of atoms of $K$.  Take the first
line before that exit as $p_K$.  In the rootless case this is also the
deferred common-translation representative.  Thus $p_K\notin\mathcal N_K$.

Root $T_K$ at $p_K$.  Every nonanchor packet line has one parent edge and
different lines have different parent edges.  Reading the
reverse-chronological Kruskal construction in its merge order, each such
edge is the unique edge which attaches the child side to the already rooted
side.  Assigning the retained time of that edge to the child gives the
chronological merge bijection \eqref{eq:packet-merge-bijection}.  Since
$\mathcal N_K\subset\mathcal P(K)\setminus\{p_K\}$, its restriction is the
first injection in (ii).

Now take $p\in\mathcal N_{\rm top}\setminus\mathcal N_K$.  By definition
its actual creation atom $b(p)$ lies outside the protected packet.  Enter
$b(p)$ in $\mathfrak F_{\rm birth}$ and define
$\psi_0(p)=t_{b(p)}$.  If its good identifier has no selected joint owner,
apply the same fixed-order selection as in the lower construction; if that
identifier was already selected, this new literal time is left in
$S_{\rm reg}$.  Theorem~\ref{thm:p38j}(ii), whose proof uses
Lemma~\ref{lem:packet-special-induction} on the selected connected component,
says that every component
outside $K$ is nonfull, so any selected local representative
here is an ordinary degree-three singleton or a periodic bookkeeping
representative, never an ordinary full-$\{4\}$ singleton.  Creation
atoms of distinct new lines are distinct,
and lower-layer atoms are disjoint from top-layer atoms even in the attached
case (attachment identifies only a boundary line), so the union map
$\psi_0$ is injective.  The $\psi_K$-owners are atoms of $K$, whereas every
$\psi_0$-owner lies outside $K$; hence their images are disjoint.  Finally use
$|\mathcal N_{\rm low}|=\mathsf X$ and
\eqref{eq:top-new-line-set}.  More explicitly, the domain of $\psi_0$ has
cardinality
\[
 \mathsf X+(m-\iota)-|\mathcal N_K|,
\]
whereas $\psi_K(\mathcal N_K)\subset R_{\rm birth}$ has cardinality
$|\mathcal N_K|$.  Disjointness of the two images therefore gives
\[
 |S_0|+|R_{\rm birth}|
 \ge \mathsf X+(m-\iota)-|\mathcal N_K|+|\mathcal N_K|
 =\mathsf X+m-\iota.
\]
Thus the count does not identify the required top credit with the possibly
smaller number $P_K-1$: the top new lines outside $\mathcal N_K$ are
assigned instead to $S_0$.  In the disjoint branch the two source
molecules and their forests use disjoint variables, giving
\eqref{eq:disjoint-owner-cardinality}.
\end{proof}

\begin{lemma}[Post-envelope owner/good intersection]
\label{lem:owner-marked-good-component}
Put
\[
 \eta=\eps^{1/(8d)},\qquad \eta_*=\eps^{1/(9d)}.
\]
Fix \(t_{\rm fin}>0\), \(\mathfrak L\in\mathbb N\) and
\(\alpha,A\ge0\), and write
\begin{equation}\label{eq:owner-good-bartau}
 \tau:=\frac{t_{\rm fin}}{\mathfrak L},
 \qquad
 \bar\tau:=\frac{\max(1,t_{\rm fin})}{\mathfrak L},
 \qquad
 q_\beta:=C_\beta\max(1,\alpha)\max(1,A),
 \qquad \vartheta:=q_\beta\bar\tau .
\end{equation}
Assume
\begin{equation}\label{eq:owner-good-diagonal-hypothesis}
 C_{\rm own}\mathfrak L\frac{\eta}{\eta_*}\le1,
\end{equation}
where the fixed constant $C_{\rm own}$ dominates every conditional chart,
support-partition and envelope constant used below.
Let $\mathcal C$ be an owner-marked $\{3\}$- or ordinary
full-$\{4\}$ component whose good status is one of the ordinary Source-B
support alternatives, rather than a periodic double-support certificate.
Keep its original indicator and its owner time through every arbitrary-
kernel operation.  First apply the collision-weight, dyadic and Maxwellian
envelope bounds while leaving the normal degree-three velocity family
unintegrated.  At this stage Proposition~\ref{prop:ordinary-subset} is used
only to identify its boundary-independent containing set, not to extract
that set's volume.  After these operations have produced the complete
pointwise conditional source majorant, the joint conditional contribution of
$\mathcal C$, including its unique owner time and the rate/activity/envelope
factor assigned to that owner, is bounded by
\begin{equation}\label{eq:owner-good-joint-bound}
 \cL_\eps^C\vartheta\eta_*
\end{equation}
times the same remaining source majorant.  No scalar good gain is asserted
for an arbitrary $Q$.  The factor $\eta_*$ is the original displayed good
credit, so marking the component as an owner does not change
$N_{\rm good}$.  For a full $\{4\}$ owner the same remaining majorant still
contains its unique $\eps^{-(d-1)}$ baseline and its $N_4$ charge.

There is also the following common-two-time clause.  Let $\mathcal X$ be an
ordinary support-good type-$\{33{\rm A}\}$ component on the short-gap cell
$|t_1-t_2|<\eta$, and let $R_{\mathcal X}$ contain the $r\in\{1,2\}$ atom
times of $\mathcal X$ which have been selected as distinct birth owners.
Assume in addition $0<\vartheta\le q_0<1$.  After the same pointwise source
majorant has been formed, keep the original indicator and integrate the two
atom times together on their common conditional measure.  Including the two
rate/activity/envelope factors assigned to the two atoms, this contribution
is bounded by
\begin{equation}\label{eq:two-time-owner-good-joint-bound}
 \cL_\eps^C\vartheta^r\eta_*
\end{equation}
times the same remaining source majorant.  The $r$ marked times are consumed
by this joint estimate and are not subsequently used in the projected-time
volume; the component contributes its one pre-existing good identifier.
No pointwise fixed-$(t_1,t_2)$ bound is used on this cell.

This is a standalone one-component conditional lemma.  It neither fixes the
owner time nor invokes the global innovation partition: both the owner cell
and the good cell are measured in the same local conditional measure
$\lambda_{\mathcal C}$.  Simultaneous application to a family of components
occurs only in Proposition~\ref{prop:packet-free-small-window-consumer} and
Lemma~\ref{lem:time-owner}, after the global one-use disintegration has been
established.
\end{lemma}

\begin{proof}
Use the notation fixed in \eqref{eq:owner-good-bartau}.
The source support partition gives a conditional good-cell measure bounded
by $\cL_\eps^C\eta$, while its final ledger retains only $\eta_*$.  Since
$\eta/\eta_*=\eps^{1/(72d)}$, hypothesis
\eqref{eq:owner-good-diagonal-hypothesis}, with $C_{\rm own}$ chosen to
dominate the fixed constant $C$ here, gives
\begin{equation}\label{eq:owner-good-margin}
 C\mathfrak L\frac{\eta}{\eta_*}\le1.
\end{equation}
Consequently
\begin{equation}\label{eq:owner-good-min}
 \min(\tau,C\eta)\le C\bar\tau\eta_*.
\end{equation}
Indeed, if $\tau>C\eta$, divide $C\eta$ by
$\bar\tau$ and use
$C\mathfrak L\eta/\max(1,t_{\rm fin})\le\eta_*$.  If
$\tau\le C\eta$, then
$t_{\rm fin}=\mathfrak L\tau\le C\mathfrak L\eta\le\eta_*<1$;
hence $\bar\tau=\mathfrak L^{-1}$ and
$\tau=\bar\tau t_{\rm fin}\le\bar\tau\eta_*$.  In particular, the
apparently dangerous subcase $t_{\rm fin}\ge1$ and
$\tau\le C\eta$ is empty on the parameter diagonal.

We now work only after the envelope operations just specified have supplied
a pointwise conditional bound by a remaining majorant
$M_{\mathcal C}^{\rm rem}$ which is measurable in the unintegrated variables
and independent of the local coordinates about to be integrated.  Let
$\lambda_{\mathcal C}$ be the original positive conditional
measure of the component, including its exact cell indicator and collision
chart density.  There is no change of the argument of the kernel here.  For
a degree-three singleton, the degree-three chart gives, after forgetting only the good
restriction,
\[
 \lambda_{\mathcal C}(\text{owner cell})
 \le C\cL_\eps^C\tau .
\]
Keeping the same cell and the same measure, the section calculation
\eqref{eq:one-atom-good-section} gives
\[
 \lambda_{\mathcal C}(\text{owner cell}\cap G_{\mathcal C})
 \le C\cL_\eps^C\eta .
\]
Therefore the exact intersection, against the already available pointwise
majorant, costs
\begin{equation}\label{eq:owner-good-common-measure}
 \int_{\text{owner cell}\cap G_{\mathcal C}}
       M_{\mathcal C}^{\rm rem}\,\dd\lambda_{\mathcal C}
 \le C\cL_\eps^C\min(\tau,C\eta)
       M_{\mathcal C}^{\rm rem} .
\end{equation}
This is a conditional measure estimate after pointwise majorization, not a
supremum of the original arbitrary packet kernel.  The velocity-support,
dyadic transported-position and time-support cells are all covered; in the
last case the owner section itself has length $\min(\tau,C\eta)$.

For an ordinary full $\{4\}$ owner, factor its unique $\eps^{-(d-1)}$
normalization into $M_{\mathcal C}^{\rm rem}$ before defining the local
conditional measure.  The remaining normalized one-atom chart has the same
literal time coordinate as the degree-three chart, with its additional
centre variables in a boundary-independent polylogarithmic box.  Therefore
forgetting the good restriction costs at most $C\cL_\eps^C\tau$, while the
degree-four extension of \eqref{eq:one-atom-good-section} costs at most
$C\cL_\eps^C\eta$ on the same normalized measure.  Their exact intersection
again satisfies \eqref{eq:owner-good-common-measure}.  Thus the joint factor
is relative to, and does not cancel or duplicate, the retained full
baseline.

Pair only the unique owner's rate with this common-measure bound.  The
non-owner atom factors stay in the remaining source majorant.  Equations
\eqref{eq:owner-good-bartau}--\eqref{eq:owner-good-min} give
\[
 q_\beta\min(\tau,C\eta)
 \le C\vartheta\eta_* ,
\]
which is \eqref{eq:owner-good-joint-bound}.  Coarsening preserves the
original component identifier, so $N_{\rm good}$ is unchanged and no
measure factor is duplicated.

For the common-two-time clause, use the original two birth charts before the
singular change from a time to a fixed output.  Their atom times are literal
Lebesgue coordinates; after the collision-weight, dyadic and Maxwellian
envelopes have produced the stated pointwise majorant, all complementary
chart densities and cutoff boxes cost at most $\cL_\eps^C$.  On a layer
interval of length $\tau$ the retained indicator has measure
\begin{equation}\label{eq:two-time-short-strip-volume}
 \big|\{(t_1,t_2):|t_1-t_2|<\eta\}\big|
 \le 2\tau\min(\tau,\eta).
\end{equation}
The diagonal hypothesis implies
\begin{equation}\label{eq:two-time-short-strip-diagonal}
 \tau\min(\tau,\eta)\le C\bar\tau^2\eta_*.
\end{equation}
Indeed, if $\tau>\eta$, then
$\eta\le C\bar\tau\eta_*$ because $\bar\tau\ge\mathfrak L^{-1}$.
If $\tau\le\eta$, then
$t_{\rm fin}=\mathfrak L\tau\le\mathfrak L\eta\le C\eta_*<1$,
so $\bar\tau=\mathfrak L^{-1}$ and
$\tau\le C\bar\tau\eta_*$; hence
$\tau^2\le C\bar\tau^2\eta_*^2\le C\bar\tau^2\eta_*$.  Multiplying
\eqref{eq:two-time-short-strip-diagonal} by $q_\beta^2$ gives
\begin{equation}\label{eq:two-time-owner-good-rate-bound}
 q_\beta^2\tau\min(\tau,\eta)
 \le C\vartheta^2\eta_*
 \le C\vartheta^r\eta_*,
\end{equation}
where the last inequality uses $r\le2$ and $\vartheta<1$.  This proves
\eqref{eq:two-time-owner-good-joint-bound} on the same measure on which the
short-gap indicator is defined; no time-strip volume is extracted from an
arbitrary kernel or reused later.
\end{proof}

\begin{proposition}[Packet-free small-window ordinary consumer]
\label{prop:packet-free-small-window-consumer}
Fix a terminal time $T>0$, and put
\begin{align*}
 \tau_T&:=T/\mathfrak L,
 &q_\beta&:=C_\beta\max(1,\alpha)\max(1,A),
 &\vartheta_T&:=q_\beta\frac{\max(1,T)}{\mathfrak L},\\
 \eta&:=\eps^{1/(8d)},
 &\eta_*&:=\eps^{1/(9d)}.
\end{align*}
and assume
\begin{equation}\label{eq:packet-free-diagonal-hypothesis}
 0<\vartheta_T\le q_0<1,
 \qquad C_{\rm own}\mathfrak L\frac{\eta}{\eta_*}\le1,
\end{equation}
where \(C_{\rm own}\) is a fixed constant large enough for the conditional
owner/good estimates below.  Fix a legal ordinary cutting branch and, in its
pre-cut source molecule, fix the actual birth forest and mark its creation
atoms before executing the branch.  Let \(\mathcal M\) be the resulting
current positive Source-B molecule, after any legal ordinary cuts and after
protection of any periodic certificates, but before calling the exceptional
Proposition~\ref{prop:source-exact-replacement}.
Let \(\mathcal N(\mathcal M)\) be its set of nonroot physical particle lines
created by the actual source birth forest.  Thus
\[
 |\mathcal N(\mathcal M)|=\mathsf X.
\]
Then the complete ordinary consumer may be disintegrated so that it retains
one distinct literal creation time for every member of
\(\mathcal N(\mathcal M)\).  With all ordinary and periodic good credits and
all full-component baselines kept in their original ledgers, its owner part
contributes
\begin{equation}\label{eq:packet-free-owner-factor}
 \vartheta_T^{\mathsf X}.
\end{equation}
Consequently the small-window replacement
\eqref{eq:source-weight-projection} gives, with no lower bound on
\(\tau_T\), the complete general estimate
\begin{equation}\label{eq:packet-free-general-ledger}
 |\mathcal I\mathcal N(\mathcal M,H,H')|
 \le \vartheta_T^{\mathsf X}\eps^{\Delta}
 \cL_\eps^{C^*(\rho+N_{\rm good}+|H'\setminus Z|)},
\end{equation}
where \(\Delta\), \(\nu\), \(N_4\), \(N_{\rm ee}\) and
\(N_{\rm del}\) are exactly those of
\eqref{eq:source-Delta}--\eqref{eq:source-nu}.  In particular this
proposition applies to the unchanged high-particle UP route; it invokes no
two-landing packet and no packet recovery map.
\end{proposition}

\begin{proof}
For every \(p\in\mathcal N(\mathcal M)\), let \(b(p)\) be its unique
creation atom in the actual source birth forest and put
\begin{equation}\label{eq:packet-free-owner-set}
 S_0:=\{t_{b(p)}:p\in\mathcal N(\mathcal M)\}.
\end{equation}
Different new lines are different child edges of the birth forest and hence
have different creation atoms.  Therefore
\begin{equation}\label{eq:packet-free-owner-cardinality}
 |S_0|=\mathsf X.
\end{equation}
This selection is made before elementary cutting.  Lemma
\ref{lem:birth-chart-retention}, which is stated for an arbitrary set of
creation atoms in the actual source forest, leaves every member of \(S_0\)
as a genuine Lebesgue output.  It also covers a marked creation atom whose
current representative would otherwise be a standalone degree-two root
substitution.  If two selected atoms lie in one two-atom component, Lemma
\ref{lem:ordinary-fixed-distinguished-times} keeps both times fixed while
estimating the complementary local variables, except that on the type-
$\{33{\rm A}\}$ short-gap cell the two original time coordinates are retained
together until the post-envelope common-measure estimate below.  Thus no
component type can erase an owner selected in
\eqref{eq:packet-free-owner-set}.

For every ordinary support-good one-atom \(\{3\}\) or full-\(\{4\}\)
identifier whose time set meets \(S_0\), put its unique member of \(S_0\)
into \(S_{\rm ord}\), in the fixed reverse forest order.  Thus such an
identifier contributes exactly one owner whenever the intersection is
nonempty, and contributes none otherwise; uniqueness follows because it is
a one-atom identifier.  In particular, no support-good one-atom birth owner
in \(S_0\) is left for a fixed-time good-section estimate.  Let
\(\mathscr G_{\rm sg}\) be the set of ordinary support-good type-
\(\{33{\rm A}\}\) identifiers on the short-gap cell
\(|t_1-t_2|<\eta\) whose atom-time set meets \(S_0\), and put
\[
 S_{\rm sg}:=S_0\cap
 \bigcup_{\mathcal X\in\mathscr G_{\rm sg}}
       \{t_1(\mathcal X),t_2(\mathcal X)\}.
\]
For \(\mathcal X\in\mathscr G_{\rm sg}\) write
\(r_{\mathcal X}:=|S_0\cap
\{t_1(\mathcal X),t_2(\mathcal X)\}|\in\{1,2\}\).
The component identifiers are disjoint, so
\(|S_{\rm sg}|=\sum_{\mathcal X\in\mathscr G_{\rm sg}}r_{\mathcal X}\).
Let \(M_t\) be the indexed set of all atom-time coordinates of
\(\mathcal M\), and put
\begin{equation}\label{eq:packet-free-owner-partition}
 S_{\rm reg}:=S_0\setminus(S_{\rm ord}\cup S_{\rm sg}),
 \qquad \mathcal A_t:=M_t\setminus S_{\rm reg}.
\end{equation}
Every selected birth owner in an ordinary two-atom component belongs to
\(S_{\rm reg}\), except for the members of \(S_{\rm sg}\); in particular the
long-gap type-\(\{33{\rm A}\}\), type-\(\{33{\rm B}\}\) and
type-\(\{44\}\) owners remain in \(S_{\rm reg}\).  Every nongood ordinary
one-atom owner and every birth owner carried by a protected periodic
certificate also belongs to \(S_{\rm reg}\).  A periodic good certificate is
integrated in its conditional velocity coordinate with its time fixed; its
different-lift tube supplies its one good excess, while the time remains for
the common projection.  A periodic or ordinary good full-\(\{4\}\) keeps
its unique \(\eps^{-(d-1)}\) baseline and its \(N_4\) charge throughout.

Keep the source indicators, and first perform the collision-weight, dyadic,
and Maxwellian-envelope bounds.  Only after the complete positive
kernel has been pointwise majorized do we integrate all current ordinary
and periodic good identifiers in one reverse topological order of their
shadow-dependency graph.  For
\(\mathcal X\in\mathscr G_{\rm sg}\) first use the common-two-time clause of
Lemma~\ref{lem:owner-marked-good-component}.  Its factors multiply to
\[
 \cL_\eps^{C|\mathscr G_{\rm sg}|}
 \vartheta_T^{|S_{\rm sg}|}\eta_*^{|\mathscr G_{\rm sg}|},
\]
and the times in \(S_{\rm sg}\) are consumed on that common conditional
measure.  On an owner in \(S_{\rm ord}\), Lemma
\ref{lem:owner-marked-good-component}, with its terminal symbols specialized
as \(t_{\rm fin}=T\), \(\tau=\tau_T\) and
\(\vartheta=\vartheta_T\), gives the common-measure factor
 \(\cL_\eps^C\vartheta_T\eta_*\); it does not multiply two marginal
estimates.  Every other ordinary good identifier uses its usual conditional
good section, with every time in \(S_{\rm reg}\) fixed; any remaining type-
\(\{33{\rm A}\}\) identifier is on the long-gap cell.  A periodic
identifier uses its conditional tube section, again with those times fixed.
Lemma~\ref{lem:global-innovation-partition} puts the pre-coarea velocity
inputs of these blocks, the normal degree-three family and the residual
velocity variables into disjoint output groups of one square full-pair
velocity transformation.  Each local nonvelocity input is adjoined only
afterward in its own equidimensional coarea chart.  Hence the same
integration coordinate is never used twice.  This loop yields
\begin{equation}\label{eq:packet-free-good-loop}
 \cL_\eps^{C N_{\rm good}}
 \vartheta_T^{|S_{\rm ord}|+|S_{\rm sg}|}\eta_*^{N_{\rm good}}
\end{equation}
and leaves every coordinate in \(S_{\rm reg}\) unintegrated.

Let \(\mathcal D_{\mathcal M}^{\rm reg}\) be the exact time domain remaining
after that loop and set
\[
 \mathcal V_{\rm reg}:=
 |\pi_{S_{\rm reg}}\mathcal D_{\mathcal M}^{\rm reg}|.
\]
All preceding operations only retained restrictions, so this projection is
contained in the original existential projection with omitted set
\(\mathcal A_t\).  Lemma~\ref{lem:distinguished-time-weight} therefore gives
\begin{equation}\label{eq:packet-free-projected-volume}
 B_{\mathcal M}\mathcal V_{\rm reg}^{5/6}
 \le \tau_T^{5|S_{\rm reg}|/6}
       \cL_\eps^{C^*|\mathcal A_t|}.
\end{equation}
After taking the \(6/5\) power, multiplication by the unique collision-weight
factor \(B_{\mathcal M}\) leaves
\begin{equation}\label{eq:packet-free-full-volume}
 B_{\mathcal M}\mathcal V_{\rm reg}
 \le \tau_T^{|S_{\rm reg}|}
       \cL_\eps^{C^*|\mathcal A_t|},
\end{equation}
because \(B_{\mathcal M}^{-1/5}\le1\).  The owners are distinct atoms, so
one actual collision-rate/activity/envelope factor, bounded by \(q_\beta\),
may be paired with each displayed time.  Since
\begin{equation}\label{eq:packet-free-tau-to-vartheta}
 q_\beta\tau_T
 \le q_\beta\frac{\max(1,T)}{\mathfrak L}
 =\vartheta_T,
\end{equation}
\eqref{eq:packet-free-good-loop}--\eqref{eq:packet-free-full-volume} and
\eqref{eq:packet-free-owner-cardinality} give
\[
 \vartheta_T^{|S_{\rm ord}|+|S_{\rm sg}|}
 (q_\beta\tau_T)^{|S_{\rm reg}|}
 \le \vartheta_T^{|S_{\rm ord}|+|S_{\rm sg}|+|S_{\rm reg}|}
 =\vartheta_T^{\mathsf X}.
\]

The logarithmic exponent created by the omitted coordinates also closes in
the original ledger.  By construction of the source molecule,
\(|M_t|=|M|=\mathsf Z\), and the actual birth forest has exactly
\(\mathsf X\) creation atoms.  Since
\(S_0=S_{\rm ord}\dot\cup S_{\rm sg}\dot\cup S_{\rm reg}\),
\begin{equation}\label{eq:packet-free-omitted-count}
 |\mathcal A_t|
 =\mathsf Z-|S_{\rm reg}|
 =(\mathsf Z-\mathsf X)+|S_{\rm ord}|+|S_{\rm sg}|
 \le \rho+3N_{\rm good}.
\end{equation}
Here \(\mathsf Z-\mathsf X\) is the source recollision/overlap count and is
at most \(\rho=|H'|+\mathsf Z-\mathsf X\), while every member of
\(S_{\rm ord}\) has a distinct good identifier and
\(|S_{\rm sg}|\le2|\mathscr G_{\rm sg}|\le2N_{\rm good}\).  The fixed factor
three is absorbed in \(C^*\).  Thus the factor
\(\cL_\eps^{C^*|\mathcal A_t|}\) in
\eqref{eq:packet-free-projected-volume} is absorbed by the displayed
\(\rho+N_{\rm good}+|H'\setminus Z|\) exponent after enlarging \(C^*\) and
never becomes an
uncontrolled \(\cL_\eps^{C\mathsf X}\) loss.

The endpoint \(\mathcal A_t=M_t\) is now explicit rather than tacit.  It can
occur only when \(S_{\rm reg}=\varnothing\).  If \(\mathsf X>0\), every birth
owner has then already been consumed either on its one-atom owner/good
conditional measure or in its short-gap common-two-time block by
Lemma~\ref{lem:owner-marked-good-component}, and its
factor \(\vartheta_T\) is present in
\eqref{eq:packet-free-good-loop}; the zero-dimensional projection is not
asked to manufacture a power of \(\tau_T\).  Conversely, a nongood birth
owner necessarily lies in \(S_{\rm reg}\), so when \(\mathsf X>0\) it rules
out \(\mathcal A_t=M_t\).  For \(\mathsf X=0\) the required time factor is
one.  This also covers the one-line stress test and remains uniform as
\(T\downarrow0\).

All other normalizations, good credits and logarithmic losses are exactly
those in the Source-B reverse-cutting calculation.  Combining the owner
factor just proved with its unchanged epsilon ledger yields
\eqref{eq:packet-free-general-ledger}.
\end{proof}

Let $M_t$ be the indexed set of all atom-time coordinates of the complete
molecule and put
\begin{equation}\label{eq:owner-time-partition}
 \begin{aligned}
 S&:=S_0\,\dot\cup\,R_{\rm birth},\\
 S_{\rm ord}&:=\{t_a:a\in\mathfrak F_{\rm own},
        \ a\text{ belongs to an ordinary support-good }\{3\}
        \text{ or full }\{4\}\text{ component}\},\\
 S_{\rm per}&:=\{t_a:a\in\mathfrak F_{\rm own},
        a\text{ carries an }\mathcal L_{\rm par}
        \text{ periodic good certificate}\},\\
 \mathscr G_{\rm sg}&:=\{\mathcal X:
        \mathcal X\text{ is ordinary support-good type-}\{33{\rm A}\},
        \ |t_1(\mathcal X)-t_2(\mathcal X)|<\eta,\\
 &\hspace{42mm}
        S\cap\{t_1(\mathcal X),t_2(\mathcal X)\}\ne\varnothing\},\\
 S_{\rm sg}&:=S\cap
       \bigcup_{\mathcal X\in\mathscr G_{\rm sg}}
       \{t_1(\mathcal X),t_2(\mathcal X)\},\\
 S_{\rm reg}&:=S\setminus(S_{\rm ord}\cup S_{\rm sg}),\\
 \mathcal A_t&:=M_t\setminus S_{\rm reg}.
 \end{aligned}
\end{equation}
The provenance ledgers in \eqref{eq:three-ledgers} make
$S_{\rm ord}$, $S_{\rm sg}$ and $S_{\rm per}$ pairwise disjoint.  A periodic
owner is deliberately
in $S_{\rm reg}$: its time is paid by the ordinary projection, while its
different-lift tube is paid once by the periodic certificate operator.  For
a protected periodic good $\{4\}$ certificate, the creation time is included
in $S_{\rm per}$ whenever that atom is an owner; its full baseline remains in
$N_4$ throughout and is not part of the owner-time projection.  An ordinary
full-$\{4\}$ owner belongs to $S_{\rm ord}$ exactly when it is support-good,
and otherwise to $S_{\rm reg}$; its baseline remains in $N_4$ in either case.
Every selected birth time in an ordinary two-atom $\{33{\rm B}\}$ or
$\{44\}$ component, and in a long-gap type-$\{33{\rm A}\}$ component, lies
in $S_{\rm reg}$, as does every unselected extra periodic birth time.
The short-gap type-$\{33{\rm A}\}$ times in $S_{\rm sg}$ remain literal until
they are consumed together by the common-two-time clause of
Lemma~\ref{lem:owner-marked-good-component} after the pointwise envelope.
For the other two-atom cases, Lemma~\ref{lem:birth-chart-retention} makes
each selected time a genuine output and Lemma
~\ref{lem:ordinary-fixed-distinguished-times} collects the component's single
good gain at its unique vertex while leaving that time fixed.
The
set $\mathcal A_t$ includes all incoming time--normal root pivots and the
ordinary support-good owner times in $S_{\rm ord}\cup S_{\rm sg}$ that will
be consumed by Lemma~\ref{lem:owner-marked-good-component} after the source
envelope is formed.  No time in $S$ is itself a time--normal root pivot.  This partition
prevents both the ordinary time-support gain and the periodic tube gain from
being mistaken for a second copy of an owner-time measure.  Put
\begin{equation}\label{eq:periodic-owner-certificate-family}
 \mathscr C_{\rm per}:=\{C(a):t_a\in S_{\rm per}\},
\end{equation}
where $C(a)$ is the protected periodic certificate carried by $a$.  The
fixed rule selecting at most one atom of $\mathfrak F_{\rm own}$ from each
immutable certificate identifier makes $a\mapsto C(a)$ bijective onto this
family; hence $|\mathscr C_{\rm per}|=|S_{\rm per}|$.  For fixed $S$-times
$s$, put
\begin{equation}\label{eq:time-fibre}
 \mathcal D_H(s)=\{a:(s,a)\in\mathcal D_H\},
 \qquad
 \pi_S\mathcal D_H=\{s:\mathcal D_H(s)\ne\varnothing\}.
\end{equation}

\begin{lemma}[Legal virtual refinement and coarsening of an all-C protected packet]
\label{lem:packet-virtual-refinement}
Fix one positive source term and a protected one-layer packet $K$ consisting
entirely of C-atoms.  Before protection there is a finite sequence of
Source-B legal cuts
\[
 K=K^0\longrightarrow K^1\longrightarrow\cdots\longrightarrow K^N,
 \qquad N\le |K|,
\]
such that $K^N$ is a union of elementary one-atom components.  At every
step the Source-B cutting identity is an equality of positive Radon measures
with the states related by $\mathrm{Sh}$ and $\mathrm{Ba}$.  More precisely,
let $\mu_K$ be the positive packet--complement measure of the original term,
including its source normalization but not the test kernel.  There are a
finite set $\Lambda$, pairwise disjoint half-open refined cells $E_\lambda$,
refined positive Radon measures $\mu^{\rm ref}_\lambda$, and Borel composite
shadow maps
\begin{equation}\label{eq:virtual-refinement-composite-shadow}
 \Sigma_\lambda:(\widetilde z,t,z_{\rm comp})\longmapsto
       (z,t,z_{\rm comp}),
 \qquad
 z_e=\widetilde z_{\operatorname{Sh}_\lambda(e)},
\end{equation}
such that $\Sigma_\lambda$ is the identity on every atom time and every
complement variable and
\begin{equation}\label{eq:virtual-refinement-radon-pushforward}
 \boxed{\quad
 \mu_K=
 \sum_{\lambda\in\Lambda}
  (\Sigma_\lambda)_\#
  \bigl(\mathbf1_{E_\lambda}\mu^{\rm ref}_\lambda\bigr).
 \quad}
\end{equation}
Equivalently, for every nonnegative Borel kernel $Q$ coupling the packet to
its complete complement, with the original time indicator retained,
\begin{align}
 &\int Q(z,t,z_{\rm comp})\mathbf1_{\mathcal D_K}(t)\,
          \dd\mu_K(z,t,z_{\rm comp})
 \notag\\
 &\qquad=
 \sum_{\lambda\in\Lambda}
 \int_{E_\lambda}
 Q\!\left(\Sigma_\lambda(\widetilde z,t,z_{\rm comp})\right)
 \mathbf1_{\mathcal D_K}(t)\,
 \dd\mu^{\rm ref}_\lambda(\widetilde z,t,z_{\rm comp}).
 \label{eq:virtual-refinement-test-kernel-identity}
\end{align}
The refined side is used only to select a Borel shadow/component order.
It is then pushed forward and summed in
\eqref{eq:virtual-refinement-radon-pushforward}; the joint packet coarea is
applied once, and only once, to the coarsened measure $\mu_K$ on the left.
No packet recovery map or packet coarea is applied termwise to
$\mu^{\rm ref}_\lambda$.  No refined component owns a good gain,
degree-four baseline, or independent time credit.
\end{lemma}

\begin{proof}
The exceptional packet $K$ is a finite one-layer C-molecule.  In a virtual
copy of the term, choose the chronologically first atom not yet isolated and
cut that singleton as free.  With no O-atoms inside $K$, this is case (i) of
the Source-B cutting definition: each boundary bond becomes the matched
free/fixed pair and the isolated atom is elementary.  Repeat on the
remaining atoms.  The process ends after at most $|K|$ cuts.

At the $j$th step apply \eqref{eq:source-cutting-identity}, retaining the
cell indicator and the complete complement kernel.  On a compact exhaustion,
one cut and its finite half-open cell partition give, for every nonnegative
Borel $F$,
\[
 \int F\,\dd\mu_{K^{j-1}}
 =\sum_{\lambda_j}
   \int_{E_{\lambda_j}}F\circ\Sigma_{\lambda_j}^{(j)}
       \,\dd\mu_{\lambda_j}^{j}.
\]
Because this holds for all $F$, uniqueness of Radon measures gives the
corresponding one-step push-forward identity.  The one-step map is the
identity on $t$ and $z_{\rm comp}$ and sends every original edge state to
the state of its current shadow, exactly as in
\eqref{eq:source-cutting-identity}.

Iterate the identity.  Functoriality of push-forward composes the one-step
maps into $\Sigma_\lambda$, while positivity and finite Tonelli combine the
cell labels into the finite set $\Lambda$.  This proves
\eqref{eq:virtual-refinement-radon-pushforward} on every compact exhaustion;
monotone convergence proves it for the full Radon measures.  Since every
composite map fixes $t$, multiplication by the original
$\mathbf1_{\mathcal D_K}(t)$ commutes pointwise with the push-forward.
Taking
$F=Q\mathbf1_{\mathcal D_K}$ proves
\eqref{eq:virtual-refinement-test-kernel-identity} for an arbitrary joint
kernel, without a factorization or boundary supremum.

The refinement is performed on a copy solely to define the shadow order
required by the simultaneous-volume proposition; it is not appended to the
actual cleanup sequence and does not violate protection.  After extracting
that order, first form the push-forward sum in
\eqref{eq:virtual-refinement-radon-pushforward}.  It is exactly $\mu_K$, so
the single joint packet coarea used later acts on the original packet
measure, not on any refined branch.  The normalization, gain and time
ledgers are consequently those of $K$, exactly once.
\end{proof}

\begin{lemma}[Fibrewise full-molecule time owner]
\label{lem:time-owner}
Apply \eqref{eq:p38j-packet} to the Source-B positive kernel
\eqref{eq:source-positive-Q}.  Fix the selected Borel conditional
representative of
Definition~\ref{def:appD-selected-conditional-representative}.  For every
fixed value of all $S$-variables, the packet inequality for that
representative is applied to
$Q\mathbf 1_{\mathcal D_H(s)}$ and its output is zero when
$s\notin\pi_S\mathcal D_H$.  Take the union, over the declared countable
positive cell family of $K$, of the fixed Borel null sets selected at the
measure-kernel level in Lemma~\ref{lem:appD-fibre-time-owner}, and denote it
by $N_{S,K}$.  This union is still null, is fixed before $Q$, and both
conditional sides below are zero-extended on it.  The recovery map preserves
the same indicator
$\mathbf 1_{\mathcal D_H}$, the same $S$-variables and the same Maxwellian
energy envelope.  With
$\sigma_K=\cL_\eps^C\eps^{2(d-1)}\epstar^{-2(d-1)}$, the fibre inequality is
\begin{equation}\label{eq:fibrewise-packet-owner}
 \overline J_{K,s}\!\left(Q\mathbf 1_{\mathcal D_H(s)}\right)
 \le \eps^{-(d-1)}\sigma_K\,
 \mathbf 1_{N_{S,K}^c}(s)
 \mathbf 1_{\pi_S\mathcal D_H}(s)
 \int Q_{K,s}^{\sharp}\mathbf 1_{\mathcal D_H}(s,a^\sharp)
 \dd a^\sharp .
\end{equation}
The bar is part of the notation: it denotes the selected, synchronously
zero-extended conditional kernel fixed above.  Removing the bar is permitted
only after integration in \(s\), when one recovers the original complete
operator.
Every point counted on the right carries the displayed witness
$a^\sharp$ for the same source projection.  Moreover the subsequent
full-molecule weight, time-volume, collision-rate, activity and Boltzmann-
envelope calculation retains the complete owner factor
\begin{equation}\label{eq:full-time-owner-factor}
 \vartheta^{|S|}
 \le \vartheta^{\mathsf X+m-\iota},
 \qquad 0<\vartheta\le q_0<1,
\end{equation}
up to the already allowed logarithmic ledger.  Thus the attached branch
retains $\vartheta^{\mathsf X+m}$, while in the disjoint branch the two
independent factors retain respectively $\vartheta^{\mathsf X}$ and
$\vartheta^{m-\iota}$.
\end{lemma}

\begin{proof}
Disintegrate the positive associated measures with respect to the
$S$-variables and use the common fixed zero-extended representatives of
Lemma~\ref{lem:appD-fibre-time-owner}.  This changes no integrated source
operator and is simultaneous for all positive kernels.  For fixed $s$,
include
$\mathbf 1_{\mathcal D_H(s)}$ in the arbitrary nonnegative test function of
Theorem~\ref{thm:two-landing}.  Definition~\ref{def:packet-output} does not
replace this indicator by the product time cell: after every root
substitution it remains
$\mathbf 1_{\mathcal D_H}(s,a)$.  Hence the resulting kernel is supported on
the exact existential projection $\pi_S\mathcal D_H$ in
\eqref{eq:time-fibre}.  The packet coarea replaces $a$ by the recovered live
coordinates $a^\sharp$ but evaluates the original indicator at that recovered
point.  Therefore nonzero output supplies the witness
$(s,a^\sharp)\in\mathcal D_H$, and the fibrewise theorem gives exactly
\eqref{eq:fibrewise-packet-owner}.  No enlargement is made in the
$s$-coordinates or in their existential projection.

At the fibrewise packet stage we stop with
\eqref{eq:fibrewise-packet-owner}; no good-cell volume and no projected-time
volume is pulled through its arbitrary kernel.  Insert that output into the
complete iterated source integral \eqref{eq:source-iterated-integral}, retain
all original indicators, and apply the collision-weight, dyadic and
Maxwellian envelope estimates.  For the normal degree-three family,
\eqref{eq:source-volume} is used here only to construct the common
boundary-independent containing set $Y$; its variables are not yet
integrated and its volume is not yet extracted.  The resulting conditional
integrand is now bounded pointwise by the nonnegative source majorant.  This
is exactly the stage corresponding to Source B's $Q_1$, before its good
components are integrated.  Only at this point do the following
measure-level operations begin.

Let $\mathscr G_{\rm ord}$ be the set of \emph{all} current ordinary
support-good component identifiers in the source decomposition and let
$\mathscr G_{\rm per}$ be the set of \emph{all} current protected periodic
certificate identifiers counted by $N_{\rm good}$.  Their one-atom
owner-marked subset $\mathscr G_{\rm ord}^{\rm own}$ and periodic
owner-marked subset
$\mathscr G_{\rm per}^{\rm own}$ are in bijection with $S_{\rm ord}$ and
$S_{\rm per}$, respectively.  The family $\mathscr G_{\rm sg}$ from
\eqref{eq:owner-time-partition} is a disjoint subset of
$\mathscr G_{\rm ord}$, and each of its identifiers carries one or two
members of $S_{\rm sg}$.  Thus
\begin{equation}\label{eq:all-good-partition}
 \begin{split}
 N_{\rm good}&=|\mathscr G_{\rm ord}|+|\mathscr G_{\rm per}|,\\
 N_{\rm good,own}&=|S_{\rm ord}|+|S_{\rm per}|
                       +|\mathscr G_{\rm sg}|,\qquad
 N_{\rm good,rest}=N_{\rm good}-N_{\rm good,own}.
 \end{split}
\end{equation}
Form one dependency graph on the complete vertex set
$\mathscr G_{\rm ord}\dot\cup\mathscr G_{\rm per}$.  Draw $C\to C'$ when a
free shadow of $C$ is used as a fixed boundary state of $C'$.  The legal
cutting order makes this graph acyclic.  In particular, an unmarked good
vertex is not discarded merely because it lies between two owner-marked
vertices.  Integrate \emph{all} vertices in one global reverse topological
order; neither provenance nor owner status is integrated in a separate
batch.

For $C\in\mathscr G_{\rm sg}$ apply the common-two-time clause of Lemma
~\ref{lem:owner-marked-good-component} after the pointwise envelope, and
consume the two original atom times together.  For
$C\in\mathscr G_{\rm ord}^{\rm own}$ apply Lemma
~\ref{lem:owner-marked-good-component} to the intersection of its owner and
good cells.  For
$C\in\mathscr G_{\rm ord}\setminus
(\mathscr G_{\rm ord}^{\rm own}\cup\mathscr G_{\rm sg})$ use, at
this same post-envelope stage, the ordinary Source-B conditional good-cell
bound \eqref{eq:source-elementary-operator}.  If the component contains a
distinguished birth time in $S_{\rm reg}$, use the fixed-time version in
Lemma~\ref{lem:ordinary-fixed-distinguished-times} and leave every such time
outside the local block for the later projection.  Its
$\cL_\eps^C\eta$ cost is absorbed into the retained
$\cL_\eps^C\eta_*$ factor; any type-$\{33{\rm A}\}$ component remaining in
this case is on the long-gap cell.  For a periodic vertex apply Lemma
~\ref{lem:periodic-certificate-free-coordinate} to its complete local block
$W_C=(w_C,\zeta_C)$: use the tube section in $w_C$ and the
boundary-independent Source-B box or norm-one substitution in the
complementary variables.  Keep every distinguished birth-owner time carried
by that certificate fixed for the later ordinary time projection, whether
or not it is the selected bookkeeping representative.  Thus the whole
non-distinguished protected local
block is consumed once at this vertex.  In every case all boundary states on which the current good
domain depends have already been fixed, and every unintegrated predecessor
remains a parameter.  The owner injection, the at-most-two selected atom
owners per $\{33{\rm A}\}$ identifier and the single use of its component
identifier, and
\eqref{eq:certificate-pivot-disjointness} prevent reuse of a local block.

Collecting the owner-marked ordinary vertices, including the short-gap
common-two-time blocks, gives exactly once
\begin{equation}\label{eq:good-owner-product}
 \cL_\eps^{C(|S_{\rm ord}|+|\mathscr G_{\rm sg}|)}
 \vartheta^{|S_{\rm ord}|+|S_{\rm sg}|}
 \eta_*^{|S_{\rm ord}|+|\mathscr G_{\rm sg}|},
\end{equation}
and their owner times are henceforth only existential witnesses in
$\mathcal A_t$.  Collecting the owner-marked periodic vertices gives,
also exactly once,
\begin{equation}\label{eq:periodic-owner-tube-product}
 \cL_\eps^{C|S_{\rm per}|}
 \eps^{(d-1)|S_{\rm per}|}
 \le \cL_\eps^{C|S_{\rm per}|}\eta_*^{|S_{\rm per}|}.
\end{equation}
For a periodic good $\{4\}$ owner its unique $\eps^{-(d-1)}$ normalization
remains in the $N_4$ part of the remaining source majorant; the displayed
tube is its good excess relative to that baseline.  Only after this single
complete interleaved reverse-order loop is complete have all unmarked
ordinary and periodic vertices also supplied their
$\cL_\eps^C\eta_*$ factors.  Altogether the loop yields
\begin{equation}\label{eq:complete-good-dependency-product}
 \cL_\eps^{C N_{\rm good}}
 \vartheta^{|S_{\rm ord}|+|S_{\rm sg}|}\eta_*^{N_{\rm good}},
\end{equation}
with every good identifier used exactly once.  Only then do we integrate the normal
degree-three velocity family over the containing set $Y$ from
\eqref{eq:source-volume}.  Because $Y$ is uniform in the fixed variables,
every good domain, including every unmarked domain which crosses an owner
dependency, has been used in the legal Source-B dependency order and the
simultaneous normal volume is extracted afterwards.

There is no overlap between these velocity operations.  Invoke
Lemma~\ref{lem:global-innovation-partition} on the present conditional
fibre.  The pre-coarea velocity blocks $U_C$ of \emph{all} periodic
certificates, the velocity inputs of ordinary good blocks, normal velocities
and residual velocity variables are disjoint output groups of one complete
square full-pair transformation.  Each complete periodic output
$W_C=(w_C,\zeta_C)$ is formed only afterward from its disjoint
equidimensional input $(U_C,Z_C)$.  A marked degree-two birth is a transparent
relay: its innovation is spent by whichever one of those groups genuinely
owns its downstream free copy and has no extra group of its own.  This does
not pretend that a transported mixed birth direction is a literal
particle-velocity slot, and it never inverts a one-particle projection of a
C-scattering.  Equation~\eqref{eq:normal-volume-on-innovation-complement} therefore
has exactly the Source-B dimension and bound for the remaining normal
family.  Every complete periodic block is integrated exactly once in
\eqref{eq:complete-good-dependency-product}; the owner-marked subproduct,
and only that subproduct, is displayed separately in
\eqref{eq:periodic-owner-tube-product}.

Let $\mathcal D_H^{\rm reg}$ be the exact remaining time domain and put
\[
 \mathcal V_{\rm reg}
 :=|\pi_{S_{\rm reg}}\mathcal D_H^{\rm reg}|.
\]
For every fixed value of the other remaining variables this set is contained
in the original existential projection
$\widetilde{\mathcal D}_{\mathcal A_t}$: the post-envelope integrations only
retain restrictions and never create a new time vector.  Therefore
\[
 \mathcal V_{\rm reg}
 \le |\widetilde{\mathcal D}_{\mathcal A_t}|.
\]
The scale-correct projection repair, Lemma
~\ref{lem:distinguished-time-weight} with retained set $R=S_{\rm reg}$ and
omitted set $\mathcal A_t$, gives
\begin{equation}\label{eq:time-weight-owner}
 B_{\mathcal M}\mathcal V_{\rm reg}^{5/6}
 \le\tau^{5|S_{\rm reg}|/6}
       \cL_\eps^{C^*|\mathcal A_t|}.
\end{equation}
Taking the $6/5$ power gives
\begin{equation}\label{eq:time-owner-complete-tau}
 \mathcal V_{\rm reg}
 \le B_{\mathcal M}^{-6/5}\tau^{|S_{\rm reg}|}
       \cL_\eps^{C^*|\mathcal A_t|}.
\end{equation}
The collision-weight estimate contributes the unique factor
$B_{\mathcal M}$, so the net scalar is
$B_{\mathcal M}^{-1/5}\le1$.  The remaining one-sixth source integration is
part of this same complete estimate and is not a second projected-volume
factor.

If $s_j\in[a_j,a_j+\tau]$, the exact dilation of the remaining
projection is
\begin{equation}\label{eq:owner-time-dilation}
 \int_{\pi_{S_{\rm reg}}\mathcal D_H^{\rm reg}}G(s)\,\dd s
 =\tau^{|S_{\rm reg}|}
   \int_{\widehat\pi_{S_{\rm reg}}\mathcal D_H^{\rm reg}}
     G(a+\tau u)\,\dd u .
\end{equation}
This identity is used only for the already-majorized conditional integrand;
it is not an operator inequality for the original $Q$.

Every $t_a\in S_{\rm per}$ is among these regularly projected coordinates.
Its local velocity block was integrated in
\eqref{eq:periodic-owner-tube-product} with this time held fixed, and the
time is now paid by the single common projection
\eqref{eq:time-owner-complete-tau}.  Thus the tube and time factors are two
conditional coordinates of one Tonelli integration, not two independent
copies of an owner-time measure.  No ordinary Source-B support gain is
assigned to these certificates.

Every time in $R_{\rm birth}$ remains a literal coordinate because its
marked contact pivots a relative cluster translation, and all incoming
time--normal root pivots belong instead to $\mathcal A_t$.  In particular,
\begin{equation}\label{eq:time-distinguished-owner}
 |R_{\rm birth}|=P_K-1,
 \qquad |S|=|S_0|+P_K-1.
\end{equation}

For the velocity-volume proposition use the legal virtual refinement of
Lemma~\ref{lem:packet-virtual-refinement}.  Its explicit singleton cut list
and measure identities construct the required shadow/component order, after
which all refined pieces are coarsened and joint packet coarea is applied
once to the common positive integrand.  Because
$|K|\le P_{\rm pkt}(d)$, its branch and
logarithmic costs are $O_d(1)$.

The default sharp FCT chart in Theorem~\ref{thm:two-landing} uses
$2(d-1)$ separator-velocity coordinates together with the two lower landing
times.  In the overlapping case the radial speed is only a local exterior
augmentation proving the positive rank-$(d-1)$ Schur factor; it is not an
additional global recovery variable.  The full-$2d$ disjoint chart and the
full-chord overlapping chart remain independent coarse fallbacks.  All of
these charts use complete elastic reflections, free transports, contact roots
and joint coarea.  Elastic energy is
preserved at every atom, the $S$-times, normals and selected free velocities
are not enlarged, and all original velocity cutoffs remain by zero
extension.  The Maxwellian envelope in
\eqref{eq:source-weight-1} and the fixed-boundary independence in
\eqref{eq:source-volume} are therefore unchanged.

It remains to convert the physical time powers into the source-series
parameter without losing an owner.  Set
\[
 q_\beta:=C_\beta\max(1,\alpha)\max(1,A).
\]
The normalized positive source kernel assigns to each collision atom its
collision rate together with the applicable activity and Boltzmann-envelope
factor; the Source-A atom bound is at most $q_\beta$ per atom (the maxima
also cover recollision atoms and normalized roots, where one of these
factors is absent).  The $S_{\rm ord}$ factors were paired with their unique
owner rates in \eqref{eq:good-owner-product}, while each short-gap block
there uses the atom rates entering its common-two-time clause.  The owner maps use distinct
atoms, so one rate may be paired with each of the $|S_{\rm reg}|$ projected
factors of $\tau$; this includes the times in $S_{\rm per}$.  Their tube
factors come separately from \eqref{eq:periodic-owner-tube-product} and do
not consume another time coordinate.  More precisely, after the complete
source operations described above, let $\mathscr R_H$ denote the remaining
nonnegative iterated majorant, with all owner-time variables already
integrated.  If $S_{\rm per}$ contains a protected good-$\{4\}$ owner,
$\mathscr R_H$ still contains that component's single $\eps^{-(d-1)}$ baseline
and the corresponding $N_4$ charge; only its tube excess appears in the
$\eta_*$ factor below.  The owner part of that actual consumer satisfies
\begin{equation}\label{eq:owner-source-majorant}
 \begin{aligned}
 \mathfrak I_{H,\rm owner}
 &\le \cL_\eps^{C(|S_{\rm ord}|+|S_{\rm per}|
                              +|\mathscr G_{\rm sg}|)}
       \vartheta^{|S_{\rm ord}|+|S_{\rm sg}|}
       \eta_*^{|S_{\rm ord}|+|S_{\rm per}|
                              +|\mathscr G_{\rm sg}|}
       B_{\mathcal M}q_\beta^{|S_{\rm reg}|}
       \mathcal V_{\rm reg}\,\mathscr R_H\\
 &\le \cL_\eps^{C^*(|\mathcal A_t|+|S_{\rm per}|
                                +|\mathscr G_{\rm sg}|)}
       \vartheta^{|S_{\rm ord}|+|S_{\rm sg}|}
       \eta_*^{|S_{\rm ord}|+|S_{\rm per}|
                              +|\mathscr G_{\rm sg}|}
       (q_\beta\tau)^{|S_{\rm reg}|}\mathscr R_H,
 \end{aligned}
\end{equation}
because $B_{\mathcal M}B_{\mathcal M}^{-6/5}\le1$.
By
\eqref{eq:effective-layer-parameter},
\begin{equation}\label{eq:tau-to-vartheta-owner}
 q_\beta\tau
 =\frac{C_\beta\max(1,\alpha)\max(1,A)t_{\rm fin}}
        {\mathfrak L}
 \le\vartheta .
\end{equation}
The remaining atoms are handled by the existing source recollision,
component and cutoff ledger and are not counted a second time as owners;
their exceptional costs remain in the displayed
 $\cL_\eps^{C^*(|\mathcal A_t|+|S_{\rm per}|
 +|\mathscr G_{\rm sg}|)}$ and the existing source
 recollision/good ledger.
Consequently the owner contribution contains
\begin{equation}\label{eq:vartheta-owner-count}
 \vartheta^{|S_{\rm ord}|+|S_{\rm sg}|}
 (q_\beta\tau)^{|S_{\rm reg}|}
 \le\vartheta^{|S_{\rm ord}|+|S_{\rm sg}|+|S_{\rm reg}|}
 =\vartheta^{|S|}
 \le\vartheta^{\mathsf X+m-\iota},
\end{equation}
where the last inequality uses
\eqref{eq:complete-owner-cardinality} and $0<\vartheta\le1$.  When $r>0$,
$\iota=0$, giving the attached count $\mathsf X+m$.  In the disjoint branch
the two source measures, owner maps and time projections factor; applying
\eqref{eq:disjoint-owner-cardinality} on the two factors gives
$\vartheta^{\mathsf X}$ below and $\vartheta^{m-\iota}$ above.  The factors
$\eta_*^{|S_{\rm ord}|+|S_{\rm per}|+|\mathscr G_{\rm sg}|}$ in
\eqref{eq:owner-source-majorant} are precisely the owner-marked ordinary and
periodic parts of the final good ledger, with one factor for each short-gap
component identifier.  They join the good factors from
unmarked components and are not multiplied again.  Thus the original
$N_{\rm good}$ exponent is unchanged.  This proves
\eqref{eq:full-time-owner-factor}.  Tonelli
justifies every fibrewise integration because the kernels are nonnegative.
\end{proof}

\subsection{Attached and disjoint Source-B consumers}

Set
\begin{equation}\label{eq:sigma-K}
 \sigma_K=\cL_\eps^C\eps^{2(d-1)}\epstar^{-2(d-1)}.
\end{equation}

\begin{proposition}[Exact top-layer transfer]
\label{prop:v2-consumer}
The packet alternative in Theorem~\ref{thm:p38j} satisfies items
(I1)--(I5) of Proposition~\ref{prop:source-exact-replacement}.  Hence it
supplies the second interpoland
\eqref{eq:source-top-second-interpoland} in both the attached and disjoint
branches.  The first interpoland
\eqref{eq:source-top-first-interpoland} follows from the ordinary uncut
molecule bound together with
Proposition~\ref{prop:packet-free-small-window-consumer} in the attached
branch, from the exact top--lower
factorization in the disjoint nonempty branch, and directly from the packet
bound in the disjoint empty branch.
\end{proposition}

\begin{proof}
Items (I1)--(I2) follow from Definition~\ref{def:packet-output},
Theorem~\ref{thm:two-landing}, Lemma~\ref{lem:time-owner}, and the protected
cleanup with dependency order in
Lemmas~\ref{lem:protected-band-restart}--
\ref{lem:packet-size-partition-ledger}.  We verify the two consumers
explicitly.

In the attached branch, integrate $K$ first and place the complete top
complement and all lower layers inside $Q$.  Lemma~\ref{lem:time-owner}
retains the source time factor.  Since every exceptional-top cleanup
component outside $K$ is nonfull, the top contribution $N_4^1$ in the source
split is zero.  The
lower DOWN algorithm is unchanged and the cut of the attached packet creates
at least one lower fixed end.  Source B's CH-molecule incidence identity
states that the number of connected lower components plus $N_{\rm ee}$ is at
most $|\widetilde H|$.  A DOWN input creates a degree-four component if and
only if that connected input is full.  Since at least one lower component
carries the fixed end just created by the packet cut, $N_4^2$ is strictly
smaller than the total number of lower components.  Combining the two
integer inequalities gives
\begin{equation}\label{eq:attached-baseline-new}
 N_4^2+N_{\rm ee}\le|\widetilde H|-1.
\end{equation}
The packet version of \eqref{eq:source-attached-ledger} is
\begin{multline}\label{eq:attached-consumer}
 \mathcal I\mathcal N(\widetilde{\mathcal M},\widetilde H,H')
 \le\vartheta^{\widetilde{\mathsf X}}
 \eps^{(d-1)(|\widetilde H|-1-N_4^2-N_{\rm ee})}\sigma_K\\
 \times |\log\eps|^{C^*(\widetilde\rho+N_{\rm ee}+N_{\neg N3})}
 \eps^{(9d)^{-1}N_{\rm good}+|H'\setminus Z|/4}.
\end{multline}
The exponent preceding $\sigma_K$ is nonnegative by
\eqref{eq:attached-baseline-new}.  The logarithmic factors involving good
pieces and initial links are absorbed by their displayed positive powers;
more precisely, the unchanged Source-B count gives
\[
 \widetilde\rho+N_{\rm ee}+N_{\neg N3}
 \le C^*(\widetilde\rho+N_4^2+N_{\rm good}+|H'\setminus Z|),
\]
and $N_4^2\le|\widetilde H|-1\le C^*\widetilde\rho$.
The remaining allowed loss therefore depends on $\widetilde\rho$.  In the
attached branch $r>0$, so $\iota=0$ and
$\widetilde{\mathsf X}=\mathsf X+m$.  Choose the source constant $c_0$
small enough (after the fixed combinatorial constant $C$) that
$\vartheta\le q_0<1$.  Then
$\vartheta^m\le\vartheta^{m/10}$, and the layer power in
\eqref{eq:attached-consumer} is stronger than
the one in \eqref{eq:source-top-second-interpoland}.  Moreover the exact
source relation
$\widetilde\rho\le\rho+r+3\Gamma$ and $r>0$ lets us enlarge $C^*$ once and
bound every remaining logarithmic loss by
$|\log\eps|^{C^*(\rho+r)}$.  From \eqref{eq:packet-excess-d},
\[
 \sigma_K\le\eps^{2(d-1)-1/(100d)}
 \le\eps^{d-1+1/(18d)}
\]
for $d\ge4$ and small $\eps$.
Thus \eqref{eq:attached-consumer} implies
\eqref{eq:source-top-second-interpoland} for every fixed $d\ge4$.

In the disjoint branch retain the exact factorization
\eqref{eq:source-disjoint-factorization}.  The lower factor is precisely
\eqref{eq:source-disjoint-lower}, or is one when the lower molecule is empty.
For the top factor, the rooted
normalization \eqref{eq:two-landing-rooted} and the time-owner lemma give
\begin{equation}\label{eq:disjoint-consumer}
 \mathcal I\mathcal N(\widetilde{\mathcal M}_\ell,G,\varnothing)
 \le\vartheta^{m-\iota}
 \eps^{(d-1)(|G|-1)}\sigma_K,
 \qquad |G|=r+\iota.
\end{equation}
If $r=0$, then $\iota=1$, $|G|=1$ and the connected high-recollision graph
has $m\ge2$; hence $\vartheta^{m-1}\le\vartheta^{m/9}$ and the epsilon power is
$2(d-1)-1/(100d)$.  If $r>0$, then $\iota=0$, the time power is stronger than
$\vartheta^{m/9}$ and the epsilon power is
$(d-1)(r+1)-1/(100d)$.  No bound of $r$ by
$\Gamma$ is used here (the source only gives $r+m\le2\Lambda_\ell$).  Choose
$C_{14}^*>4$.  Uniformly for every integer $r\ge1$, the difference from the
target epsilon exponent in \eqref{eq:source-disjoint-top} is
\[
 r\left(d-1-\frac4{C_{14}^*}\right)-\frac{59}{900d}>0.
\]
For $r=0$ the difference is $d-1-59/(900d)>0$.
After $C_{14}^*$ is fixed, the
per-root logarithmic factor is $\eps^{-o(1)r}$ and is absorbed by a fixed
fraction of this linear margin.  Thus both cases are stronger than
\eqref{eq:source-disjoint-top}, uniformly in the full source range of $r$.

We verify the first interpoland by three distinct routes.  In the attached
branch, before choosing the exceptional packet, the ordinary Source-B
theorem, with the packet-free small-window owner input of
Proposition~\ref{prop:packet-free-small-window-consumer}, applies to the
original uncut $\widetilde{\mathcal M}$ and gives
\eqref{eq:source-top-first-interpoland}.

In the disjoint branch with nonempty lower molecule, we do \emph{not} use
that uncut route.  Apply the exact factorization
\eqref{eq:source-disjoint-factorization}; use
\eqref{eq:source-disjoint-lower} on the lower factor and
\eqref{eq:disjoint-consumer} on the top factor.  If $r=0$, then $m\ge2$ and
the top layer power $\vartheta^{m-1}$ dominates $\vartheta^{m/10}$ after the fixed
one-layer choice, while $1/(12d)>1/(18d)$ and
$4(C_{14}^*)^{-1}\rho$ dominates the required
$(C_{14}^*)^{-1}R$; the top factor
$\eps^{2(d-1)-1/(100d)}$ is unused room.  If
$r>0$, the top layer power $\vartheta^m$ dominates $\vartheta^{m/10}$ and
$\eps^{(d-1)(r+1)-1/(100d)}$ absorbs the $r$-part and every bounded
logarithmic loss.
Multiplication gives \eqref{eq:source-top-first-interpoland} for the same
positive integral.

If the disjoint lower molecule is empty, its factor is one and
\eqref{eq:disjoint-consumer} alone has exponent
$2(d-1)-1/(100d)$ for $r=0$ and
$(d-1)(r+1)-1/(100d)$ for $r>0$, together with the stronger time power
just checked.
After $C_{14}^*$ is fixed, the same uniform linear-in-$r$ margin absorbs the
per-root subpower loss and implies the first interpoland directly.  These
mutually exclusive routes exhaust the attached/disjoint cases.
Lemma~\ref{lem:source-two-regime-envelope}, with $C_{15}^*$ chosen after
$C_{14}^*$ as in \eqref{eq:source-constant-separation}, gives
\eqref{eq:source-top-target}; substituting that estimate into
\eqref{eq:source-top-error-reduction} gives
\eqref{eq:source-trunc-output}.  This proves (I3)--(I5) and the proposition.
\end{proof}

\section{Epoch restart and proof of the main theorem}
\label{sec:global}

This section has two tasks.  It closes the quantitative estimate on one epoch,
then restarts and stitches those estimates across the full time range of
Theorem~\ref{thm:main}.  The two-landing geometry is already complete.  The
full-range restart additionally uses
Theorem~\ref{thm:birth-flag-frame-production} and
Corollary~\ref{cor:automatic-fixed-k-birth-production}, proved independently
in Section~\ref{sec:transfer-principle}.  We begin with the quantitative
estimate on one epoch satisfying \eqref{eq:one-epoch-horizon}.  The block
restart then stitches these estimates over the longer range
\eqref{eq:loglog-main}.  On one epoch, we choose the layer hierarchy in the
required order, propagate the three endpoint errors, and aggregate them
without losing the fixed epsilon margin.  The proof
separates the physical layer length $\tau$ from the effective series parameter
$\vartheta$, controls the growing constants, and reruns the fixed-terminal
expansion for every $T\in[0,t_{\rm fin}]$ before summing the cumulants.

Put $L=|\log\eps|$ and choose
\begin{equation}\label{eq:number-layers}
 \mathfrak L=\max\{1,\lfloor\kappa\sqrt{\log L}\rfloor\},
 \qquad \tau=\frac{t_{\rm fin}}{\mathfrak L},
\end{equation}
where $\kappa=\kappa(d,\beta)>0$ is fixed sufficiently small.  After that,
choose the one-epoch constant in \eqref{eq:one-epoch-horizon} sufficiently small
relative to $\kappa$.  For small $\eps$,
$\mathfrak L\ge(\kappa/2)\sqrt{\log L}$, hence
$\tau\le2c_0/\kappa<1$.  There is no lower bound on $\tau$.
The quantitative Source-A estimate shows that the atom-series parameter is
not bare $C\tau$ but
\begin{equation}\label{eq:effective-layer-parameter}
 \vartheta
 :=\frac{C_\beta\max(1,\alpha)\max(1,A)
                    \max(1,t_{\rm fin})}{\mathfrak L}.
\end{equation}
For small $\eps$, the one-epoch hypothesis and the lower bound for
$\mathfrak L$ give $\vartheta\le2C_\beta c_0/\kappa$.  Choose $c_0$ after
$C_\beta$ and $\kappa$ so that
\begin{equation}\label{eq:effective-layer-small}
 0<\vartheta\le q_0<1.
\end{equation}
Physical contact-time volumes continue to use $\tau$; every source
summability, restart and top-consumer power uses $\vartheta$.

\subsection{Cutoffs and exceptional constants}

Use the dimension-dependent hierarchy
\begin{equation}\label{eq:4040-hierarchy}
 A_{\mathfrak L}=L_\eps^\#,\qquad
 \Lambda_\ell=A_\ell^{10d},\qquad
 A_{\ell-1}=\Lambda_\ell^{10d}.
\end{equation}
In particular,
\begin{equation}\label{eq:lambda-gap}
 \Lambda_{\ell-1}=\Lambda_\ell^{100d^2}
 \gg \Lambda_\ell^2A_\ell.
\end{equation}
Appendix~\ref{app:parameters} verifies the fifteen dominance stages, including
the two genuine layer-dependent feedbacks.  It also includes the ordinary
cutoff/lift exponent and the initial-forest exponent, and defines
$C_{\rm all}^*=\max\{C_1^*,\ldots,C_{15}^*,C_{\rm ord}^*,C_{\rm init}^*\}$.
The unified bound is
\begin{equation}\label{eq:Cstar-growth}
 C_{\rm all}^*(\mathfrak L)\le C^{C\Gamma\mathfrak L^2}
 \le L^\delta,
\end{equation}
for any fixed $\delta<1/10$, after choosing
$C\Gamma\kappa^2<\delta$ and reducing $\eps_0$.  This is the quantitative
Source-A Part-7 bound; it explicitly includes
$C_{10}^*\gg(C_7^*)^{\mathfrak L}$ and
$C_6^*\gg(C_5^*)^{\Gamma\mathfrak L}$.  Therefore
\begin{equation}\label{eq:subpower-log}
 \mathscr S_\eps
 =L^{C\mathfrak L C_{\rm all}^*}
 =\exp\!\bigl(O(\mathfrak L L^\delta\log L)\bigr)
 =\eps^{-o(1)}
\end{equation}
because $\mathfrak L L^\delta\log L=o(L)$.  This includes all non-root
occurrences.  Losses repeated per root require the more precise quantity
\begin{equation}\label{eq:qeff-def}
 q_{\rm eff}=\eps^{1/C_{14}^*}L^{C_{\rm all}^*}.
\end{equation}
Equations \eqref{eq:Cstar-growth}--\eqref{eq:qeff-def} give
\begin{equation}\label{eq:qeff-bound}
 \log q_{\rm eff}
 \le-L^{1-\delta}+L^\delta\log L
 \le-\tfrac12L^{1-\delta}
\end{equation}
for small $\eps$.  This is stronger than any subset-counting loss for
$s\le L$.
In the small-$R$ branch of the two-regime envelope the same bound gives
$R\le C L^\delta$ and
$C^*R\log L\le C_{\rm all}^*R\log L
\le C L^{2\delta}\log L=o(L)$; this is the uniform absorption
proved in Lemma~\ref{lem:source-two-regime-envelope}, not a fixed-layer
smallness argument.

\subsection{Iteration of the error}

The cumulative-error recurrence in dimension $d$ is
\begin{equation}\label{eq:err-recurrence}
 \|\Err_\ell^2\|_1
 \le\eps^{-2(d-1)\Lambda_\ell^2A_\ell}
 \|\Err_{\ell-1}\|_1.
\end{equation}
The restart factors are
\begin{equation}\label{eq:restart-factors}
 b_0=\eps^{\Lambda_0/10},
 \qquad b_\ell=\vartheta^{\Lambda_\ell/10}.
\end{equation}

\begin{lemma}[Quantified restart domination]
\label{lem:restart-domination}
Fix $M>0$.  Choose $\kappa$ first, then $c_0$ sufficiently small relative to
$\kappa$, and finally $\eps_0$ sufficiently small relative to
$(M,\beta,\kappa,c_0)$.  For $0<\eps\le\eps_0$,
\begin{equation}\label{eq:restart-domination}
 b_0\prod_{k=1}^{\mathfrak L}
      \eps^{-2(d-1)\Lambda_k^2A_k}
 +\sum_{\ell=1}^{\mathfrak L}
 b_\ell\prod_{k=\ell+1}^{\mathfrak L}
      \eps^{-2(d-1)\Lambda_k^2A_k}
 \le \eps^M.
\end{equation}
The estimate is uniform over every admissible
$(\alpha,A,t_{\rm fin})$ satisfying the one-epoch condition
\eqref{eq:one-epoch-horizon}.
\end{lemma}

\begin{proof}
The choice of $c_0$ in \eqref{eq:effective-layer-small} gives
$\vartheta\le q_0<1$; put $a_0=-\log q_0>0$.  Write
$Y_k=\Lambda_k^2A_k=\Lambda_k^{2+1/(10d)}$.  From
$\Lambda_k=\Lambda_{k+1}^{100d^2}$ and
$\Lambda_{\mathfrak L}=(L_\eps^\#)^{10d}\ge L^{10d}$, successive $Y_k$ decrease so rapidly that,
after reducing $\eps_0$,
\begin{equation}\label{eq:tail-loss-geometric}
 \sum_{k=\ell+1}^{\mathfrak L}Y_k
 \le2Y_{\ell+1}
 =2\Lambda_{\ell+1}^{2+1/(10d)}
 \qquad(0\le\ell<\mathfrak L).
\end{equation}
This is a direct finite geometric comparison: every ratio
$Y_{k+1}/Y_k
=\Lambda_{k+1}^{-(100d^2-1)(2+1/(10d))}$ is at most $1/2$.

For $1\le\ell<\mathfrak L$, the negative logarithm of the restart factor is
at least $a_0\Lambda_\ell/10$, whereas the total later loss has logarithm at
most
\begin{equation}\label{eq:restart-vs-tail}
 2(d-1)L\sum_{k>\ell}Y_k
 \le4(d-1)L\Lambda_{\ell+1}^{2+1/(10d)}
 \le \frac{a_0}{20}\Lambda_{\ell+1}^{100d^2}
 =\frac{a_0}{20}\Lambda_\ell.
\end{equation}
The middle inequality holds uniformly in $\ell$, since the smallest possible
$\Lambda_{\ell+1}$ is $(L_\eps^\#)^{10d}\ge L^{10d}$.  Hence the $\ell$th summand on the left of
\eqref{eq:restart-domination} is at most
$\exp(-a_0\Lambda_\ell/20)$.  This is at most $\eps^{M+1}$ because the
smallest $\Lambda_\ell$ in this range is
$\Lambda_{\mathfrak L-1}=(L_\eps^\#)^{1000d^3}\ge L^{1000d^3}$.

For $\ell=\mathfrak L$ the product is empty and
\[
 b_{\mathfrak L}\le\exp(-a_0L^{10d}/10)\le\eps^{M+1}.
\]
For the initial error, \eqref{eq:tail-loss-geometric} gives
\[
 b_0\prod_{k=1}^{\mathfrak L}\eps^{-2(d-1)Y_k}
 \le \exp\left[-L\left(\frac{\Lambda_0}{10}
             -4(d-1)\Lambda_1^{2+1/(10d)}\right)\right]
 \le \eps^{\Lambda_0/20}\le\eps^{M+1}.
\]
Finally $\mathfrak L+1=O(\sqrt{\log L})\le\eps^{-1/2}$ for small $\eps$;
summing the displayed $\eps^{M+1}$ bounds proves
\eqref{eq:restart-domination}.
\end{proof}

Iterating \eqref{eq:err-recurrence} yields exactly the left side of
\eqref{eq:restart-domination}: the $\ell=0$ term is the initial error and the
$\ell\ge1$ terms are the errors restarted in layer $\ell$.  Taking, for
example, $M=M_0:=1/100$ makes the iterated error negligible relative to the final
target exponent $1/(400d)$.

\begin{lemma}[Uniform terminal-time rerun]
\label{lem:uniform-terminal-rerun}
For every $T\in(0,t_{\rm fin}]$ there is a version-2 expansion constructed
on $[0,T]$, denoted
$(f^{\mathcal A,T},E_H^T,f_s^{{\rm err},T},\Err_{\mathfrak L}^T)$, which
satisfies the exact identities \eqref{eq:source-final-cumulant} and
\eqref{eq:source-iterated-error} with $t_{\rm fin}$ replaced by $T$.
Every source estimate and every estimate proved in
Sections~\ref{sec:ensemble-flow}--\ref{sec:p38j} that is used below holds for
this rerun with the same constants and the same $\eps_0$, uniformly in $T$.
At $T=0$ the corresponding assertion is the initial expansion of Section~3.
\end{lemma}

\begin{proof}
For $T>0$, this is the terminal-time content of
Proposition~\ref{prop:sourceA-wrapper} and Corollary
\ref{cor:diagonal-reduction}.  They apply the fixed structural endpoint
theorem separately to $[0,T]$, rather than interpolating inside the
decomposition of $[0,t_{\rm fin}]$.  Online Resource~1, Section~\ref{app:source-interface}
records the ownership distinction.  In the cited Source-A/Source-B arguments,
the only use of a lower comparison for the layer length is the
$\tau^{-1}\le|\log\eps|$ step behind
\eqref{eq:source-weight-projection-original}; it is replaced here by
\eqref{eq:source-weight-projection} together with
Proposition~\ref{prop:packet-free-small-window-consumer} on every ordinary,
high-particle and attached-first route, and with
Lemma~\ref{lem:time-owner} on the exceptional packet route.  Thus no repaired
projection is asked by itself to supply a missing owner power, and the
estimate applies for arbitrarily small $T/\mathfrak L$, while
$\vartheta_T\le\vartheta$ and all starred constants and thresholds remain
unchanged.  At $T=0$, \eqref{eq:source-initial-cumulant},
\eqref{eq:source-initial-error}, Corollary
\ref{cor:initial-cumulant-L1}, and the object assignment
\eqref{eq:T0-object-assignment} give the endpoint directly: the truncation
error vanishes, and the kinetic approximation error vanishes because
$f^{\mathcal A,0}(0)=n_0=n(0)$, while the unique forest
remainder is $\Err_{\mathfrak L}^0=\Err_0$.
\end{proof}

\subsection{Paid-cumulant restart and epoch stitching}

The preceding rerun is a one-epoch statement.  Its output now becomes boundary
data for a new epoch without exposing the old collision word to the new layer
selector.  This step enlarges the total time range.

The boundary class used below is weighted.  This is forced by the
collision rate: a bare \(L^1\) kernel may be concentrated at arbitrarily
large velocity, whereas the flux at one contact grows linearly in the
relative speed.  The Gaussian coefficient decreases with the global kinetic
clock.  Its radius is independent of the current number of
retained labels, so losing a red port never strengthens the output weight.

Put
\[
 M_*:=\max(1,\alpha)\max(1,A),\qquad
 \lambda_0:=\beta/32,\qquad
 \nu:=C_{\beta,d}M_*,\qquad
 \lambda_t:=\frac{\lambda_0}{1+\nu t},
\]
where \(C_{\beta,d}\) is chosen so that the negative derivative of
\(\lambda_t\) absorbs the high-energy part of the absolute collision rate.
For a nonempty label set \(H\), define
\begin{equation}\label{eq:common-energy-Gaussian-weight}
 E_H(v):=\sum_{i\in H}|v_i|^2,\qquad
 W_{H,t}(v):=
 \exp\!\left\{\lambda_t E_H(v)\right\}.
\end{equation}
If \(H\cap K=\varnothing\), then
\begin{align}
 W_{H\sqcup K,t}(v_H,v_K)
 &= W_{H,t}(v_H)W_{K,t}(v_K),
 \label{eq:common-energy-factorization}\\
 W_{H,t}(v_H)&\le W_{P,s}(v_P)
 \qquad(H\subset P,\ 0\le s\le t),
 \label{eq:common-energy-marginal}
\end{align}
The first identity is additivity of the total energy.  For the second,
\(E_H(v_H)\le E_P(v_P)\) and \(\lambda_t\le\lambda_s\).  It remains valid
after an energy-preserving hard-sphere flow before the variables in
\(P\setminus H\) are integrated.  This marginal monotonicity is the reason
for using one common Gaussian radius rather than a cardinality-dependent
moving threshold.

At a nonzero-time boundary, the exact reorganization below takes one complete
signed kernel as input and never splits it into its coordinates.  All other
bottom leaves retain their one-particle Source-B provenance until the final
integration.

\begin{lemma}[Paired stopped-contact partition identity]
\label{lem:paired-first-contact-partition}
Let \(\pi\) be a partition of a finite label set and let
\(\mathsf S_\pi^\eps(t,s)=\bigotimes_{B\in\pi}\mathsf S_B^\eps(t,s)\)
be the product of the exact periodic hard-sphere groups inside its blocks,
with interblock contacts transparent.  It acts on the product chamber
\(\mathcal D_\pi\), where only intrablock exclusions are imposed.  The
initial datum is supported in the full hard-core chamber
\(\mathcal D\subset\mathcal D_\pi\), and the physical output below is
extended by zero to the common ambient product chamber.  For distinct
blocks \(B,C\), let
\(\mathsf F_{B,C;\pi}^{\rm first}(u,s)\) be the incoming flux density of
product trajectories whose first interblock contact is the \(B,C\) contact
at time \(u\), and let \(\mathsf R_{B,C}^{\rm C}\) and
\(\mathsf R_{B,C}^{\rm O}\) be respectively its reflected and transparent
boundary continuations.  The trace is only weak notation: the
time-integrated compositions below are defined by contact coarea, not by a
fixed-time trace of an \(L^1\) representative.  If
\(\mathsf U_\pi^\eps\) denotes the physical evolution in which an
interblock contact merges its incident blocks, then, first for smooth data
on every finite-energy truncation,
\begin{equation}\label{eq:paired-first-contact-renewal}
 \begin{split}
 \mathsf U_\pi^\eps(t,s)
 ={}&\mathsf S_\pi^\eps(t,s)\\
 &+\sum_{\{B,C\}\subset\pi}\int_s^t
 \left[
 \mathsf U_{\pi/(B\sim C)}^\eps(t,u)\mathsf R_{B,C}^{\rm C}
 -\mathsf S_\pi^\eps(t,u)\mathsf R_{B,C}^{\rm O}
 \right]
 \mathsf F_{B,C;\pi}^{\rm first}(u,s)\,\dd u .
 \end{split}
\end{equation}
Only the C continuation is recursively expanded, so iteration terminates
after at most \(|\pi|-1\) interblock merges.  The paired bracket is one
signed transfer; its O continuation is not placed in the merged hard-core
flow.
\end{lemma}

\begin{proof}
Outside the null set of Proposition~\ref{prop:ae-flow}, two independent
block trajectories have either no interblock contact on \((s,t)\), or a
unique first interblock contact pair and time; simultaneous or grazing
first contacts form a null set.  Denote the stopped no-contact evolution by
\(\mathsf S_\pi^{0,\eps}\).  Partitioning the transparent product path and
the physical path by the same first pair and first time gives, respectively,
\[
 \mathsf S_\pi^\eps
 =\mathsf S_\pi^{0,\eps}
 +\sum_{B<C}\int \mathsf S_\pi^\eps\mathsf R_{B,C}^{\rm O}
                  \mathsf F_{B,C;\pi}^{\rm first},
 \qquad
 \mathsf U_\pi^\eps
 =\mathsf S_\pi^{0,\eps}
 +\sum_{B<C}\int \mathsf U_{\pi/(B\sim C)}^\eps
                  \mathsf R_{B,C}^{\rm C}
                  \mathsf F_{B,C;\pi}^{\rm first}.
\]
Here every suppressed propagator carries the times \(t,u,s\) displayed in
\eqref{eq:paired-first-contact-renewal}.  Subtracting the first equality
from the second proves that identity.  In particular, the transparent O
state is never used as initial data for a merged hard-core group.

For comparison with the layered molecule formula, expand
\(\prod_{B\ne C}\mathbf1_{\{B\not\sim C\}}\), group graphs by their
rooted forest, and then sum the complete fibre of allowed IAL choices.
This is exactly the finite Möbius sum in the Source-B
inclusion--exclusion identity recorded in
Online Resource~1, Section~\ref{app:source-interface}.  Reversing that finite calculation
gives the same paired stopped transfer, including its zero-contact branch.
Thus the identity is not a termwise reinterpretation of a residual special
C/O word.  For smooth truncated data the two displayed path partitions are
ordinary change-of-variables formulae.  The weighted estimate in
Lemma~\ref{lem:sealed-one-hyperedge}, proved on those data, is a bound for
the complete time-integrated signed transfer on the structured class with
one arbitrary weighted-\(L^1\) kernel \(K_P\) and all complementary leaves
in their stated pointwise source class.  It therefore extends that complete
transfer uniquely by density in \(K_P\).  Neither C nor O is separately
asserted to be a fixed-time \(L^1\) boundary operator, and no
\(\eps^{d-1}\) gain is assigned to a \(K_P\)-to-\(K_P\) contact.
\end{proof}

The next estimate is the analytic part of this reorganization.  It also
specifies the stopping rule, so a contact trace is never applied to a bare
\(L^1\) kernel.

\begin{lemma}[Sealed one-hyperedge propagation]
\label{lem:sealed-one-hyperedge}
Fix one genuine layer \([s,t]\), with
\(C_{\beta,d}M_*(t-s)\le q_0\).  At its lower boundary insert one complete
kernel \(K_P\), keep its \(P\) variables as one sealed red hyperedge, and
write every other bottom factor as either \(n(s)\) or a one-particle
remainder with the pointwise Boltzmann envelope of
\eqref{eq:source-fa}.  Reorganize the exact one-layer hierarchy by
\eqref{eq:paired-first-contact-renewal} as follows:
\begin{enumerate}[label=\textup{(\roman*)}]
\item all contacts internal to an already merged influence block are summed
in its exact group \(\mathsf S_B^\eps\);
\item every maximal ordinary acyclic decoration with one active boundary
port is summed by the Fr\'echet derivative of the complete Source-B
one-output polynomial;
\item a paired stopped transfer which first joins a previously independent
output-root block is retained as an interblock merge slot;
\item the stopping slices are free-time slices and every exposed boundary
transfer is the complete C--O bracket in
\eqref{eq:paired-first-contact-renewal}.
\end{enumerate}
Let \(H_R\) be the upper roots in the causal influence closure of the sealed
hyperedge, put \(p=|P|\), \(h_R=|H_R|\), and factor every component outside
that closure by the exact Source-B disconnected-subforest identity.  Then
the corresponding complete red-closure contribution satisfies
\begin{equation}\label{eq:sealed-one-hyperedge-bound}
 \|W_{H_R,t}\mathsf V_{H_R\leftarrow P}^{[s,t]}K_P\|_1
 \le \mathfrak d(t,s)^{p+h_R}
 \chi^{\,(h_R-p)_+}
 \|W_{P,s}K_P\|_1 ,
\end{equation}
where
\begin{equation}\label{eq:sealed-layer-factors}
 \mathfrak d(t,s):=
 \exp\!\left\{C_{\beta,d}
 \bigl((1+\nu t)^3-(1+\nu s)^3\bigr)\right\},
 \qquad
 \chi:=C_{\beta,d}\eps^{d-1}.
\end{equation}
For \(H_R=\varnothing\), the operator is the free-slice marginal
contraction and the root factor is omitted.  The estimate is uniform in
the number of internal periodic collisions.  It is not obtained from the
old \(\Gamma\)-truncation.
\end{lemma}

\begin{proof}
Work first on one finite Source-B molecule prefix and on a finite kinetic
energy set.  Between two skeleton atoms,
Proposition~\ref{prop:ae-flow} and energy conservation give
\[
 \|W_{B,u_2}\mathsf S_B^\eps(u_2,u_1)F_B\|_1
 \le \|W_{B,u_1}F_B\|_1 .
\]
Thus all internal red--red, red--blue and blue--blue periodic repeats cost
one after the blocks carrying them have merged.  In particular, no
quadratic pair rate is taken in absolute value.

We now sum the unmarked source decorations before estimating the marked
kernel.  Let \(\mathsf A_I^{[N]}[f](u)\) denote the \(N\)-th binary Picard
polynomial of the ordinary acyclic one-output Source-B recursion on
\([s,u]\), after the complete twist/C--O fibre has been summed.  Each such
polynomial is finite, and their union contains every finite ordinary rooted
molecule.  This is the finite
identity recorded in Source B by
\[
 \begin{gathered}
 \texttt{eq.fAterm2},\quad \texttt{prop.twistcancel},\\
 \texttt{eq.fB\_3Moleq\_3},\quad \texttt{eq.IMrecursive}.
 \end{gathered}
\]
Its kernels and its full partial-order domain are recorded here in
\eqref{eq:source-associated-integral}--\eqref{eq:source-time-domain}.

We first construct the untruncated ordinary map rather than identifying it
with the finite Source cutoff.  Put \(A_{\rm src}=\max(1,A)\) and
\(\gamma_I(u)=9\beta/80-\mu(u-s)\), where
\(\mu=C_{\beta,d}M_*\).  Choose the fixed genuine-layer constant \(q_0\)
so small that \(C_{\beta,d}M_*(t-s)\le q_0\) implies
\(\gamma_I(u)\ge\beta/16>\lambda_0\).  The binary recursion
\texttt{eq.IMrecursive} is precisely the Picard recursion for the signed
translated collision form
\(\mathcal Q_\eps=\mathcal Q_{\rm C}-\mathcal Q_{\rm O}\).  Its positive
majorant \(\mathsf B_I^{[N]}\) is obtained by using \(|f|\) and replacing
the complete signed C--O bracket by the sum of its two positive kernels.
Thus \(|\mathsf A_I^{[N]}|\le\mathsf B_I^{[N]}\).  With
\(F_N(u)=\sup_{x,v}e^{\gamma_I(u)|v|^2}
 \mathsf B_I^{[N]}(u,x,v)\), it obeys
\begin{equation}\label{eq:scalar-positive-Gaussian-iteration}
 D^+F_N(u)
 \le \sup_{v\in\mathbb R^d}
 \{-\mu|v|^2F_N(u)
   +C_{\beta,d}\alpha F_N(u)^2(1+|v|)\}.
\end{equation}
Indeed, energy conservation transfers the Gaussian to the two incoming
velocities, the unused velocity has an integrable Gaussian, and the torus
translations do not change the spatial supremum.  On the bootstrap region
\(F_N\le2C_{\beta,d}A_{\rm src}\), the identities
\(M_*=\max(1,\alpha)A_{\rm src}\) and \(\mu=C_{\beta,d}M_*\) make the right side at
most \(C_{\beta,d}M_*F_N\).  The bottom source envelope gives
\(F_N(s)\le C_{\beta,d}A_{\rm src}\), because its weakest exponent is
\(9\beta/80\).  A first-exit argument and the layer smallness give,
uniformly in \(N\),
\begin{equation}\label{eq:intermediate-source-envelope}
 \sup_N\sup_{s\le u\le t}\sup_{x,v}
 e^{\gamma_I(u)|v|^2}
 |\mathsf A_I^{[N]}[f](u,x,v)|\le C_{\beta,d}A_{\rm src} .
\end{equation}
Apply the same estimate to the difference of two successive iterates of the
positive fixed-point map defining \(\mathsf B_I\).  In the integrated norm
containing the negative
\(-\mu|v|^2\) term, the high-velocity part is estimated by
\[
 M_*\sup_{r\ge0}r\int_0^{t-s}e^{-\mu\theta r^2}\,\dd\theta
 \le C M_*\sqrt{\frac{t-s}{\mu}}
 \le C\sqrt{M_*(t-s)}.
\]
The bounded-velocity part costs \(CM_*(t-s)\).  Thus the Lipschitz ratio is
at most
\(C_{\beta,d}\{M_*(t-s)+\sqrt{M_*(t-s)}\}<1/2\), after the same fixed
reduction of \(q_0\).
Hence the positive Picard sequence converges, and it absolutely dominates
the signed sequence \(\mathsf A_I^{[N]}\).  The latter therefore converges
in the weaker Gaussian norm, uniformly for \(u\in[s,t]\), to a signed map
\begin{equation}\label{eq:untruncated-ordinary-source-map}
 \mathsf A_I^\infty[f](u):=\lim_{N\to\infty}
 \mathsf A_I^{[N]}[f](u).
\end{equation}
Thus \eqref{eq:intermediate-source-envelope} is an independent
intermediate-time estimate of the same ordinary molecule series; it is not
deduced by interpolating the endpoint assertion \eqref{eq:source-fa} or by
using the weaker molecule-by-molecule \(\tau^{n/2}\) estimate.

The derivative series is absolutely convergent as well.  Write
\(a_I=\mathsf A_I^\infty[f]\) and
\(\mathsf G_I=D\mathsf A_I^\infty[f]\).  At each binary atom its mark
continues through exactly one of
\(\mathcal Q_\eps(g,a_I)\) and
\(\mathcal Q_\eps(a_I,g)\); the number of derivative marks never branches.
For \(r\) simultaneous active ports the generator is the sum of these
one-port generators over the \(r\) coordinates.  Its order-\(n\) positive
Picard majorant has at most \((2r)^n\) slot choices and the ordered marked
times give the simplex factor \(1/n!\).  This is the point at which a raw
two-sided Galton--Watson expansion, which has no such simplex, is avoided.

Here is the exact finite-prefix partition which justifies the
differentiation.  Starting from the future influence closure of the upper
set \(H_R\), retain in the skeleton the lower sealed set \(P\), the upper
set \(H_R\), every atom whose two descendant cones contain distinct active
ports, every selected birth packet, and every additional recollision,
periodic overlap or cycle not belonging to the ordinary paired C/O fibre.
An ordinary O companion is kept with its C mate.  Delete those skeleton
atoms.  Every remaining maximal component has no active port, exactly one
lower and one upper active port, or one lower and no upper active port.  The
first belongs to the pure-source factor, the second is one term of
\(D\mathsf A_I^\infty[f]\), and the third is a free-slice Liouville
marginal contraction.  A zero-lower, one-upper component is a pure-source
birth base and remains a credited skeleton slot.  The assertion follows by
induction from the binary recursion \texttt{eq.IMrecursive}: at the highest
atom either the mark uses the first child or it uses the second, while the
other child is the already summed scalar polynomial.  The complete
twist/C--O fibre is summed before this partition, so the induction is an
identity for the signed polynomial, not a termwise absolute-value
replacement.  A source label which interacts with the red closure later in
the block belongs to the future closure and is an active boundary port; it
is not marginalized.  In particular, if its ordinary source tree crosses
the nominal pure-source cutoff before that contact, the whole ordinary
decoration belongs to \(\mathsf A_I^\infty\), not to a Source defect.  The
Source cutoff is applied only to components which remain future-disconnected
from red throughout the block.  Half-open layer cells and the
first-failure owner make the alternatives disjoint.  Consequently
preclusters, later attachments and alternating raw contact times are all
present, but only the marked spine is ordered by
\eqref{eq:source-time-domain}; unmarked branches have already been summed
into \(a_I\).

We next estimate this exact one-port map.  At a marked collision the
activity \(\alpha\eps^{-(d-1)}\) cancels the contact-sphere factor
\(\eps^{d-1}\).  If the new source velocity is \(w\), energy conservation
gives \(E_P(v'_P)\le E_P(v_P)+|w|^2\), and therefore
\(W_{P,u}(v'_P)\le W_{P,u}(v_P)e^{\lambda_u|w|^2}\), also when the derivative
mark changes collision slot.  The pointwise source envelope
\eqref{eq:intermediate-source-envelope}, with
\(\gamma_I(u)\ge\beta/16>\lambda_u\), absorbs this factor and the
relative speed.  Summing
over the \(p\) marked coordinates gives
\[
 C_{\beta,d}\nu
 \bigl(p+\sqrt p\sqrt{E_P}\bigr).
\]
Since
\[
 \lambda'_u=-\frac{\lambda_0\nu}{(1+\nu u)^2},
\]
completion of the square gives, after increasing \(C_{\beta,d}\),
\begin{equation}\label{eq:decreasing-Gaussian-rate-absorption}
 \lambda'_uE_P+C_{\beta,d}\nu
       \bigl(p+\sqrt p\sqrt{E_P}\bigr)
 \le C_{\beta,d}p\nu(1+\nu u)^2.
\end{equation}
Thus the decreasing coefficient, rather than a cardinality-dependent
threshold, absorbs the retained high velocity.
For the limit, let \(\Pi_R\) cut off both the marked total energy and the
new source velocity, and insert \(\Pi_R\) immediately before and after each
marked collision.  The resulting killed positive generator is bounded on
\(L^1\).  If \(\mathsf G_{I,R}^{\sharp,[N]}\) is its order-\(N\) Picard
sum, then every finite cutoff satisfies
\[
 \frac{\dd}{\dd u}
 \|W_{P,u}\mathsf G_{I,R}^{\sharp,[N]}(u,s)K_P\|_1
 \le C_{\beta,d}p\nu(1+\nu u)^2
 \|W_{P,u}\mathsf G_{I,R}^{\sharp,[N]}(u,s)K_P\|_1 .
\]
Its order-\(N+1\) tail is at most
\([C_RpM_*(u-s)]^{N+1}/(N+1)!\).  First let \(N\to\infty\), then let
\(R\to\infty\).  Positivity makes the latter limit monotone, while the
displayed decreasing-weight estimate is uniform in \(R\).  Weighted-\(L^1\)
density therefore yields
\begin{equation}\label{eq:one-port-derivative-bound}
 \|W_{P,t}\mathsf G_I^\sharp(t,s)K_P\|_1
 \le \mathfrak d(t,s)^p\|W_{P,s}K_P\|_1 .
\end{equation}
Every finite signed derivative prefix is termwise dominated by this full
positive derivative.  Conversely, absolute convergence of
\eqref{eq:untruncated-ordinary-source-map} shows that the union of the
finite prefixes is exactly the signed family of all red-dependent ordinary
one-port histories.  Thus no such history is lost at the pure-source cutoff.
Across genuine layers the chain rule
\[
 D(\mathsf A_m^\infty\circ\cdots\circ\mathsf A_1^\infty)
 =D\mathsf A_m^\infty\circ\cdots\circ D\mathsf A_1^\infty
\]
applies on each maximal one-port decoration and keeps one mark; all
multiport and special atoms remain in the skeleton.  Hence the exponents in
\eqref{eq:one-port-derivative-bound} telescope without mixing the ordinary
map with a Source remainder.  This is a restart of a closed one-particle
map, not a reset of physical grey memory.  Pure-red contacts remain in the
exact Liouville
group, and \(K_P\) enters only on free slices; whenever a skeleton merge is
exposed, the source-owned leaf, time and direction are integrated through
the complete bracket \eqref{eq:paired-first-contact-renewal} before the
variables of \(K_P\).  Thus no fixed-time contact trace is applied to the
arbitrary kernel.  After the transient skeleton half-edges are contracted
as below, only the \(p\) lower and \(h_R\) upper factors remain.  This proves
the factor \(\mathfrak d(t,s)^{p+h_R}\).

It remains to count observed root blocks.  Integrate a complete ordered
merge forest, not an intermediate correlated block.  At its bottom all
nonred leaves still carry the pointwise Maxwellian or one-particle
remainder envelope.  Remove a leaf block of the merge forest.  Its common
torus translation and its first-contact time are integration variables.
The contact coarea formula, with the flux as its normal Jacobian, gives
\(\eps^{d-1}\).  The remaining retained-energy factor is controlled by
the same decreasing Gaussian.  Indeed, with \(y=\sqrt{E_B}\),
\[
 \delta\lambda_{s,t}:=\lambda_s-\lambda_t
 =\frac{\lambda_0\nu(t-s)}{(1+\nu s)(1+\nu t)},
\]
and
\begin{align*}
 &(t-s)(1+y)e^{-\delta\lambda_{s,t}y^2}\\
 &\quad\le C_{\beta,d}\left[1+
 \sqrt{(1+\nu t)^3-(1+\nu s)^3}\right]
 \le C_{\beta,d}\mathfrak d(t,s).
\end{align*}
The first inequality follows by maximizing \(ye^{-\delta\lambda y^2}\)
and using \(\nu\ge1\); it is valid at \(t=s\) by continuity.  More
generally, \eqref{eq:decreasing-Gaussian-rate-absorption} applied between
successive skeleton times controls all ordered contact rates at once.

A finite prefix with \(b\) retained merge/birth credits is
bounded before the final forest sum by
\[
 C^{p+h_R+b}\mathfrak d_0(t,s)^{p+h_R+2b}\chi_0^b
 \|W_{P,s}K_P\|_1,
\]
where \(\mathfrak d_0\) and \(\chi_0\) denote the same two factors with
smaller absolute constants.  Port counting gives
\(b\ge(h_R-p)_+\).  After increasing the constant in
\(\mathfrak d\) and reducing \(\eps_0\),
\(C\mathfrak d_0(t,s)^2\chi_0\le1/2\) on a genuine layer; summing
\(b\ge(h_R-p)_+\) therefore gives exactly
\(\mathfrak d(t,s)^{p+h_R}\chi^{(h_R-p)_+}\).
Consequently the estimate is uniform after integration against the weighted
parent kernel, even when that parent contains the arbitrary joint \(K_P\),
and each retained merge edge supplies its contact-sphere factor.
It also covers a blue block which was internally correlated before meeting
the red block, because all its original source leaves are integrated
simultaneously along the same forest; no pointwise bound is asserted for
the intermediate block.

No \(\chi\)-credit is assigned to a contact for which both incident
variables belong to the arbitrary sealed kernel: such contacts are
internal exact flow.  In the red closure, order the credited slots from the
top down and assign to each its first still-unassigned source-owned leaf.
The forest property makes this assignment injective, so the preceding
coarea integration never traces \(K_P\) on a prescribed contact tube.

To count the red closure, contract the \(p\) sealed input ports.  A grey
event and a contact internal to the exact influence block expose no new
root port, whereas every retained merge/birth credit exposes at most one.
Thus a complete finite prefix ending with \(h_R\) red-influenced roots has
at least \((h_R-p)_+\) credited slots.  This is a port-count invariant, not
an atom-disjoint path matching.  The arbitrary-\(k\) packet input is not
used in this red-closure count; it enters only the pure-reference
exceptional branch below, where its full
\(\eps^{(d-1)(k-1)}\epstar^{-k(d-1)}\) factor is weakened to
\(q_{\rm b}^{k-1}\) in the common ledger.

In particular, the \(p=h_R=2\) pure-red C/O test remains inside the exact
\(P\)-block and has no asserted new gain.  When \(p>h_R\), lost red
variables are integrated on a free-time slice, not by a trace of \(K_P\).
Energy conservation, Liouville measure preservation and
\eqref{eq:common-energy-marginal} give the required weighted contraction
with operator norm one; the output Gaussian radius does not shrink with
\(h_R\).  A blue block which correlates and subsequently
joins the red closure contributes both its internal credited slot and its
attaching slot; if it never joins, it belongs to the factorized pure-source
part and is not included in \(H_R\).  This proves the estimate and, by the
absolute limit \eqref{eq:untruncated-ordinary-source-map} together with the
global first-failure owner for the remaining special skeleton, its exact
infinite-prefix meaning.
\end{proof}

The estimate just proved is stronger than the root factor used by the
Source-B block.  We retain the weaker common factor because it is already
compatible with all ordinary and exceptional source consumers.

\begin{theorem}[Sealed-past block stop line]
\label{thm:sealed-past-block-stop-line}
Let \(I=[a,b]\) contain \(r\le m\) new layers and satisfy the one-epoch
owner condition and \(\vartheta_I\le q_0\).  Let \(C_\bullet\) dominate
the new-block starred constants and the common forest exponent, and put
\begin{equation}\label{eq:block-S-q}
 \Sigma_{\rm b}:=L^{CmC_\bullet},\qquad
 q_{\rm b}:=\eps^{1/C_{14}^*}L^{C_\bullet},\qquad
 \mathfrak R_{\rm b}:=
 \mathfrak d(b,a)\Xi_{\rm kin}.
\end{equation}
Suppose that at time \(a\) a hierarchy is written over the exact
hard-core reference generated by the activity profile \(n(a)\), and that
each nonempty term contains at most one sealed complete kernel \(K_P\);
all remaining factors retain the displayed one-particle source provenance.
For every output \(H\ne\varnothing\), the complete new-block contribution
containing \(K_P\) obeys
\begin{equation}\label{eq:sealed-block-bound}
 \|W_{H,b}\mathsf V_{H\leftarrow P}^{I}K_P\|_1
 \le \Sigma_{\rm b}\mathfrak R_{\rm b}^{|P|+|H|}
 q_{\rm b}^{\,(|H|-|P|)_+}
 \|W_{P,a}K_P\|_1 .
\end{equation}
The pure reference contribution has the weighted estimate
\begin{equation}\label{eq:weighted-new-block-source}
 \|W_{H,b}K_H^{\rm new}\|_1
 \le \Sigma_{\rm b}\eps^{a_*}q_{\rm b}^{(|H|-1)_+},
 \qquad a_*:=\frac7{1000d}.
\end{equation}
The old history does not enter \(C_\bullet\).  The same statements hold
with one localized marked defect in place of \(K_P\).
\end{theorem}

\begin{proof}
Apply Lemma~\ref{lem:sealed-one-hyperedge} successively in the \(r\) new
layers.  Its red-dependent expansion is organized by merge/birth number,
not by the number of internal contacts, so arbitrarily many pure-red
periodic contacts remain inside the exact group.  At every free slice use
the exact disconnected-subforest identity to write
\(H=H_R\sqcup H_0\), where \(H_R\) is the red influence closure and
\(H_0\) is the union of pure-source components.  Lemma
\ref{lem:sealed-one-hyperedge} supplies
\(\chi^{(h_R-p)_+}\) on the first part.  If \(H_0\ne\varnothing\) and its
pure-source component is in the retained main class,
\eqref{eq:weighted-new-block-source} supplies
\(\eps^{a_*}q_{\rm b}^{h_0-1}\), hence \(q_{\rm b}^{h_0}\) after
\(\eps^{a_*}\le q_{\rm b}\).  Thus the total new-root credit is at least
\((h-p)_+\), including the case in which the red hyperedge is completely
marginalized.  If the disconnected \(H_0\) component is instead the
first-cutoff pure-source term, the whole tensor product with the red
closure is routed to the marked defect bucket and is not used in
\eqref{eq:sealed-block-bound}; Lemma
\ref{lem:localized-weighted-defect} below combines that product into one
marked hyperedge without assigning its source credit to \(K_P\).

Intermediate cardinalities are retained throughout this
iteration.  Compress every transient continuation chain into its exact
block group.  Every remaining transient root is injected into its first
incoming credited merge or pure-source base.  Its two adjacent
\(\mathfrak d\)-factors are absorbed using
\[
 \mathfrak d(b,a)^2\chi\le q_{\rm b},\qquad
 \mathfrak d(b,a)^2\eps^{a_*}\le q_{\rm b},
\]
after one common reduction of \(\eps_0\); only the \(p\) bottom and \(h\)
top factors remain in \(\mathfrak R_{\rm b}^{p+h}\).  The half-open owner,
lift and dyadic partitions and the global forest range contribute
\(\Sigma_{\rm b}\).  This finite-prefix radius induction proves
\eqref{eq:sealed-block-bound} without multiplying unrecorded intermediate
sizes.

For the pure reference term run the unchanged Source-B block with initial
profile \(n(a)\), taking its source parameter
\(\beta_{\rm src}=\beta/4\).  The shifted Gaussian and spatial hypotheses
follow from
\eqref{eq:weighted-assumptions} and
Lemma~\ref{lem:quantified-shifted-regularity} below.  In the source weight
estimate reserve \(\lambda_0=\beta/32\) after the collision weights have
been absorbed; this leaves \eqref{eq:weighted-new-block-source}.  Ordinary
and periodic consumers are unchanged.  On an exceptional blue merge the
arbitrary fixed-\(k\) operator estimate of
Corollary~\ref{cor:automatic-fixed-k-birth-production} is applied to the
complete sealed complement, never to a factorization of it.

Here is the equation-level constant and weight check suppressed in that
invocation.  The atom-series majorants use only the pointwise amplitude
envelope for \(n\); consequently their contraction parameter is
\(\vartheta_I=C_{\beta,d}\max(1,\alpha)\max(1,A)
\max(1,|I|)/m\), with no spatial-gradient factor.  In the proof of the
Source-B one-particle comparison, the spatial derivative occurs only after
telescoping one translated collision constraint.  Each resulting summand
has the explicit factor \(\eps\), at most \(\Lambda^2\) choices, and exactly
one bottom factor replaced by \(\nabla_xn(a)\).  Thus, for every source-tree
cutoff \(\Lambda\) used in the block, its additional contribution is bounded
by
\begin{equation}\label{eq:shifted-source-gradient-split}
 C_{\beta,d}\eps\Lambda^2 A\Xi_{\rm kin}\Sigma_{\rm b}
 \le \eps^{1/2}\le\eps^{a_*}.
\end{equation}
Indeed, the global tower gives
\(\log\Lambda+\log\Sigma_{\rm b}
=O(mL^{\delta_0}\log L)\), while
\(\log\Xi_{\rm kin}=O((\log L)^2)\); the logarithm of every factor multiplying
\(\eps\) is therefore \(o(L)\).  This is the finite-\(\eps\) comparison
absorbed in \(f^{\mathcal A}(b)-n(b)\), not an atom-series factor, and it is
used only once per block.  This proves the amplitude/gradient constant split
needed for the restarted source call.  With
\(\beta_{\rm src}=\beta/4\), the Source-B weight recursion gives
\(\beta_\ell\ge9\beta/80\) for its one-particle remainder, while
\eqref{eq:weighted-assumptions} is stronger for \(n\).  Pull \(W_{H,b}\)
backwards along a fixed history.  Free transport leaves velocities
unchanged, each C reflection
and its O companion preserve total kinetic energy and have unit velocity
Jacobian, and \eqref{eq:common-energy-factorization} separates
disjoint source subforests.  At a source atom at time \(u\),
\(\lambda_b\le\lambda_u\le\lambda_0<9\beta/80\), so the unused Maxwellian factor absorbs
every relative-speed polynomial in the contact coarea.  The volume,
owner, lift and good-cell bounds are velocity-independent, while the
exceptional packet estimate is applied after multiplying its complete
nonnegative kernel by this pulled-back weight.  Thus every printed
Source-B \(L^1\) inequality holds with \(W_{H,b}\) on its output and with
the same epsilon power.  The first kinematic factor is absorbed in
\(\Sigma_{\rm b}\), and each further one is absorbed by weakening the
corresponding Source-B per-root epsilon factor to \(q_{\rm b}\); these
inequalities
follow from the same logarithmic comparison used below
\eqref{eq:epoch-q-log}.  This proves
\eqref{eq:weighted-new-block-source}, rather than assuming a weighted
version of the source theorem.

Finally, the Source-B recurrence is algebraic before its boundary
cumulants are inserted.  Substituting one sealed kernel and reversing the
complete partition fibres by
Lemma~\ref{lem:paired-first-contact-partition} proves that the reference
branch and the one-hyperedge branches are disjoint and exhaustive.  A
time-zero cross-block exclusion belongs to the stopped zero-contact branch;
equivalently its initial exclusion edge is included in the same finite
Möbius fibre, so it is not silently replaced by independence.  It also has
the required analytic credit: keeping \(K_P\) intact and integrating the
source endpoint gives
\[
 \int_{\T^d\times\R^d}
 \mathbf1_{\{|x-y|\le\eps\}}\,n(a,y,w)\,\dd y\,\dd w
 \le C_{\beta,d}A\eps^d
 \le L^{C_\bullet}\eps^{d-1}\le q_{\rm b}.
\]
The initial forest owner and the open positive-time first-contact owner are
disjoint, so this edge is not counted twice.  A
first cutoff-saturating reference tree which stays future-disconnected from
red remains in the original global cutoff error.  A red-dependent branch
has no ordinary-tree or
internal-collision cutoff; its one-port decorations are the signed limit
\eqref{eq:untruncated-ordinary-source-map}, while all multiport and special
atoms remain in the exact skeleton of
Lemma~\ref{lem:sealed-one-hyperedge}.  This proves exactness and prevents
the old \(\Gamma\)-defect from being applied to an arbitrary red kernel.
The argument is unchanged when the unique sealed hyperedge is marked as a
defect.
\end{proof}

The root algebra used after \eqref{eq:sealed-block-bound} contains no
hidden spare-root assertion.  If \(p=|P|\), \(h=|H|\),
\(0\le\zeta,q<1\), \(R\ge1\), and \(R\zeta\le1\), then
\begin{equation}\label{eq:coloured-cardinal-algebra}
 R^{p+h}\zeta^{p-1}q^{(h-p)_+}
 \le R^{2h}\bigl(\zeta+q\bigr)^{h-1}.
\end{equation}
If \(p\le h\), the factor
\(\zeta^{p-1}q^{h-p}\) is a degree-\(h-1\) monomial and
\(R^{p+h}\le R^{2h}\).  If \(p>h\), write the left side as
\(R^{2h}\zeta^{h-1}(R\zeta)^{p-h}\).  Thus the two-red-root C/O test with
\(p=h=2\) merely transports its existing \(\zeta\), while lost red ports
are absorbed by the already paid radius rather than assigned a fictitious
new gain.

The shifted profile satisfies the remaining source hypothesis.

\begin{lemma}[Quantified shifted regularity]
\label{lem:quantified-shifted-regularity}
Under the assumptions of Theorem~\ref{thm:main},
\begin{equation}\label{eq:shifted-gradient-bound}
 \sup_{0\le t\le t_{\rm fin}}
 \|e^{\beta|v|^2/2}\nabla_xn(t)\|_\infty
 \le C_{\beta,d}A\,\Xi_{\rm kin},\qquad
 \Xi_{\rm kin}:=
 \exp\!\{C_{\beta,d}(M_*t_{\rm fin})^2+
                   C_{\beta,d}M_*t_{\rm fin}\}.
\end{equation}
In the range \eqref{eq:loglog-main},
\(\log\Xi_{\rm kin}=O((\log L)^2)=o(L)\).
\end{lemma}

\begin{proof}
The case \(t_{\rm fin}=0\) is immediate.  Otherwise apply a spatial
difference quotient \(\delta_h\) to the mild equation, put
\(g_h=|\delta_hn|\), and use the positive absolute collision form
\(\mathcal Q^\sharp\).  In the distributional mild sense,
\[
 (\partial_t+v\cdot\nabla_x)g_h
 \le\alpha\{\mathcal Q^\sharp(g_h,n)
             +\mathcal Q^\sharp(n,g_h)\}.
\]
This is an inequality, not the signed linearized Boltzmann equation.
Set
\[
 \gamma(t):=2\beta-\frac{3\beta t}{2t_{\rm fin}},\qquad
 H_h(t):=\|e^{\gamma(t)|v|^2}g_h(t)\|_\infty .
\]
Energy conservation in the collision change of variables and
\(\|e^{2\beta|v|^2}n(t)\|_\infty\le A\) give, uniformly for
\(\beta/2\le\gamma\le2\beta\),
\[
 e^{\gamma|v|^2}
 \{\mathcal Q^\sharp(g_h,n)+\mathcal Q^\sharp(n,g_h)\}(v)
 \le C_{\beta,d}A(1+|v|)H_h .
\]
Along torus characteristics the derivative of the decreasing Gaussian
contributes
\(-3\beta|v|^2/(2t_{\rm fin})\).  Consequently
\[
 \frac{\dd}{\dd t}H_h(t)
 \le \sup_{r\ge0}
 \left\{C_{\beta,d}\alpha A(1+r)
       -\frac{3\beta r^2}{2t_{\rm fin}}\right\}H_h(t)
 \le C_{\beta,d}(M_*+M_*^2t_{\rm fin})H_h(t).
\]
Since \(H_h(0)\le A\), Gronwall gives the right side of
\eqref{eq:shifted-gradient-bound}, uniformly in \(h\).  At every time
\(\gamma(t)\ge\beta/2\), and Fatou passes to the weak spatial derivative.
This estimate is used only in the finite-\(\eps\) comparison, not in
\(\vartheta\).
\end{proof}

Choose the fixed hierarchy constant so large that
\begin{equation}\label{eq:restart-constant-reserve}
 C_{14}^*\ge \frac{1000d}{7}.
\end{equation}
For a nonempty finite label set \(H\), call a joint kernel \(K_H\)
\emph{\((G,\zeta)\)-paid at time \(t\)} if
\begin{equation}\label{eq:paid-cumulant-norm}
 \|W_{H,t}K_H\|_1
 \le G\eps^{a_*}\zeta^{(|H|-1)_+},
 \qquad a_*:=\frac7{1000d},\qquad 0\le\zeta<1.
\end{equation}
Thus \(\zeta=0\) is allowed: it means that only singleton nonempty kernels
are present.  The first root is paid by the fixed factor
\(\eps^{a_*}\); every additional root is paid by the activity radius.

\begin{lemma}[Weighted complete-operator paid restart]
\label{lem:paid-cumulant-restart}
Let \(I=[a,b]\) and the block data satisfy the hypotheses of
Theorem~\ref{thm:sealed-past-block-stop-line}.  Put
\begin{equation}\label{eq:block-Omega}
 \Omega_{\rm b}:=
 C\Sigma_{\rm b}\max\{1,\mathfrak R_{\rm b}^3\}.
\end{equation}
At time \(a\), take the exact hard-core grand-canonical reference generated
by the activity profile \(n(a)\), and add a signed hierarchy whose nonempty
kernels are \((G,\zeta)\)-paid.  If
\begin{equation}\label{eq:paid-radius-admissibility}
 \Omega_{\rm b}\zeta\le\frac14 ,
\end{equation}
then the exact main hierarchy at \(b\), after the product rebase over
\(n(b)\), has an exact product-coordinate decomposition with
\(K_\varnothing^{\rm prod}=1\) and with signed nonempty kernels that are
\((G^+,\zeta^+)\)-paid, where
\begin{equation}\label{eq:paid-restart-map}
 G^+\le \Omega_{\rm b}(G+1),\qquad
 \zeta^+\le \Omega_{\rm b}(\zeta+q_{\rm b}).
\end{equation}
If \(b\) is an internal epoch boundary preceding another block, the same
hierarchy has an exact decomposition over the fresh hard-core
grand-canonical reference generated by \(n(b)\), with
\(K_\varnothing^{\rm rel}=0\) and with signed nonempty relative kernels
obeying the same paid bounds.  At such an internal boundary the
fresh-reference forest remainder is inserted, with its exact sign, in the
finite-correlation error bucket.
No old collision word and no old internal-collision count enters a starred
constant of the new block.
\end{lemma}

\begin{proof}
Substitute the boundary decomposition in the algebraic one-layer hierarchy.
Every term contains either no old kernel or exactly one complete old kernel;
all complementary factors retain their one-particle source provenance.
Apply \eqref{eq:weighted-new-block-source} to the first class and
\eqref{eq:sealed-block-bound} to the second.  If \(p,h\) are its input and
output sizes, respectively, then
\eqref{eq:coloured-cardinal-algebra} and
\(\Omega_{\rm b}\zeta\le1/4\) give
\[
 \|W_{H,b}K_{H;p}^{+}\|_1
 \le \Omega_{\rm b}(G+1)\eps^{a_*}
       \{\Omega_{\rm b}(\zeta+q_{\rm b})\}^{h-1}.
\]
The half-open label assignments and the sums over a lost red port are
geometric because \(\Omega_{\rm b}\zeta\le1/4\); their remaining fixed
factors are included in the displayed choice of \(\Omega_{\rm b}\).
This proves the map before rebasing.  In particular, a pure-red \(p=h\)
branch transports its old radius and is not assigned a fictitious new
gain.

For the exceptional pure-reference branch use the canonical birth flag
\begin{equation}\label{eq:adaptive-fixed-k}
 k_K:=\max\bigl(\{2\}\cup
 \{k\in\mathbb N:k\ge2,\ 2k-1\le P_K\}\bigr).
\end{equation}
The packet cap makes this a finite dimension-dependent choice.
Corollary~\ref{cor:automatic-fixed-k-birth-production} acts on the complete
sealed complement, and
\[
 \eps^{(d-1)(k_K-1)}\epstar^{-k_K(d-1)}
 \le \eps^{d-1}\epstar^{-2(d-1)}.
\]
Thus it is at least as strong as the two-landing input used in (R5).

Finally write \(f^{\mathcal A}(b)=n(b)+u_b\).  The exact identity
\begin{equation}\label{eq:exact-cumulant-rebase}
 \sum_{H\subset Q}(f^{\mathcal A})^{\tensor(Q\setminus H)}E_H
 =\sum_{K\subset Q}n^{\tensor(Q\setminus K)}
   \sum_{H\subset K}E_Hu_b^{\tensor(K\setminus H)}
\end{equation}
does not create a second arbitrary hyperedge.  The source Gaussian reserve
gives
\(\|W_{\{1\},b}u_b\|_1\le\Sigma_{\rm b}\eps^{a_*}\).
Use \eqref{eq:common-energy-factorization},
\(\eps^{a_*}\le q_{\rm b}\), and the binomial theorem.  Extra singleton
bases pay the extra rebase roots, so this gives the asserted paid bound in
the product \(n(b)^{\tensor}\) coordinates, with
\(K_\varnothing^{\rm prod}=1\).

At an internal boundary preceding another block, it remains to perform,
rather than suppress, the fresh-reference boundary conversion.  Let
\(\mathfrak G_Q[n(b)]\) be the exact hard-core
grand-canonical hierarchy generated by \(n(b)\).  Apply the forest
construction of Section~3 with \(n_0\) replaced by \(n(b)\) and with the
cutoff at this boundary.  It gives the exact identity
\[
 \mathfrak G_Q[n(b)]
 =\sum_{H\subset Q}n(b)^{\tensor(Q\setminus H)}
   C_H^{\rm ref}(b)+\Err_Q^{\rm ref}(b),
 \qquad C_\varnothing^{\rm ref}=1.
\]
The proof of Proposition~\ref{prop:initial-cumulant} is profile-uniform: it
uses only unit mass and the pointwise Maxwellian amplitude.  Repeating its
leaf integrations with \(W_{H,b}\), allowed because
\(\lambda_b\le\beta/32\), yields, for \(h=|H|\ge1\),
\begin{equation}\label{eq:fresh-reference-paid-bound}
 \|W_{H,b}C_H^{\rm ref}(b)\|_1\le\eps^{h/3}
 \le\eps^{a_*}q_{\rm b}^{h-1},
 \qquad
 \|W_{Q,b}\Err_Q^{\rm ref}(b)\|_1
 \le\eps^{\Lambda/10}.
\end{equation}
Here \(\Lambda\) is the boundary forest cutoff, and the second inequality
in the first bound follows from \(C_{14}^*\ge1000d/7\) and the definition
of \(q_{\rm b}\).

If \(K_H^{\rm prod}\) denotes the product-coordinate kernel obtained above,
then the required triangular conversion is the exact equality
\begin{equation}\label{eq:fresh-reference-boundary-conversion}
 \sum_{H\subset Q}n(b)^{\tensor(Q\setminus H)}K_H^{\rm prod}
 =\mathfrak G_Q[n(b)]
  +\sum_{\varnothing\ne H\subset Q}
    n(b)^{\tensor(Q\setminus H)}
    \bigl(K_H^{\rm prod}-C_H^{\rm ref}(b)\bigr)
  -\Err_Q^{\rm ref}(b).
\end{equation}
In particular,
\(K_\varnothing^{\rm rel}:=K_\varnothing^{\rm prod}
-C_\varnothing^{\rm ref}=0\).  Thus the next signed nonempty kernel is
\(K_H^{\rm prod}-C_H^{\rm ref}(b)\), not
\(K_H^{\rm prod}\).  Subtraction creates no second arbitrary hyperedge,
and \eqref{eq:fresh-reference-paid-bound} adds only the already reserved
unit amplitude and \(q_{\rm b}\) radius.  Increasing the absolute constant
in \(\Omega_{\rm b}\) once absorbs these two additions and preserves
\eqref{eq:paid-restart-map}.  The final term in
\eqref{eq:fresh-reference-boundary-conversion} is routed to the stated
finite-correlation error bucket at this internal boundary.  This proves
both the exact input type for the next block and the paid estimate; if there
is no next block, the already paid product-coordinate output is retained.
\end{proof}

The dynamical error must keep the connected subset on which it was created.
Collapsing it to an arbitrary aggregate \(L^1\) kernel would recreate the
high-velocity and contact-trace obstructions.

\begin{lemma}[Localized weighted defect hierarchy]
\label{lem:localized-weighted-defect}
Let
\begin{equation}\label{eq:defect-root-constant}
 C_{\rm def}:=C_{\rm def,0}
 \max\{C_\bullet C_{15}^*,C_{14}^*,1\},\qquad
 q_{\rm def}:=\eps^{1/C_{\rm def}}L^{C_\bullet},
 \qquad
 a_{\rm err}:=\min\{1/(20d),1/100\}.
\end{equation}
After increasing the absolute \(C_{\rm def,0}\), assume also that the
incoming \((G,\zeta)\)-paid radius satisfies
\(\zeta+q_{\rm b}\le q_{\rm def}\).  The dynamical defect created by the
pure-reference factor in one block, including its exact disconnected tensor
product with an incoming paid red closure when present, has an exact
expansion
\begin{equation}\label{eq:localized-defect-expansion}
 \mathcal D_s(t)=
 D_\varnothing(t)n(t)^{\tensor s}
 {}+
 \sum_{\varnothing\ne H\subset[s]}
 n(t)^{\tensor([s]\setminus H)}D_H(t)
\end{equation}
with
\begin{equation}\label{eq:localized-defect-bound}
 |D_\varnothing(t)|
 \le \Sigma_{\rm b}(G+1)\eps^{a_{\rm err}},
 \qquad
 \|W_{H,t}D_H(t)\|_1
 \le \Sigma_{\rm b}(G+1)\eps^{a_{\rm err}}
 q_{\rm def}^{(|H|-1)_+}.
\end{equation}
If a previously marked nonempty family enters with parameters
\((G^{\rm def},\eta)\) at base \(\eps^{a_{\rm err}}\), then the nonempty
part of its exact propagation, including fusion with a new
future-disconnected pure-source cutoff term, is again one marked family with
\begin{equation}\label{eq:marked-defect-fusion-map}
 G^{{\rm def},+}\le\Omega_{\rm b}G^{\rm def},
 \qquad
 \eta^+\le\Omega_{\rm b}(\eta+q_{\rm def}).
\end{equation}
The completely marginalized output is routed to the scalar bucket and has
normalized size at most \(\Omega_{\rm b}G^{\rm def}\).
This assertion is used in the regime
\(\Omega_{\rm b}(\eta+q_{\rm def})<1\).
Under later blocks, each \(D_H\) is kept as the unique sealed marked
hyperedge and is propagated by
Theorem~\ref{thm:sealed-past-block-stop-line}.  The finite-correlation
cutoff remainder remains in the single global recurrence
\eqref{eq:err-recurrence}.  The scalar \(D_\varnothing n^{\tensor s}\) is
propagated on the distinguished pure-reference coordinate and is never
treated as an arbitrary contact-traced hyperedge.
\end{lemma}

\begin{proof}
Retain the fixed-term information in the Source-B truncation proof before
its final sum.  First consider the pure-reference branch.  The bad
last-layer cluster is \(\lambda_{\rm bad}\), of size \(r\), and \(K_{\rm src}\) is
the source-generated connected index.  Exact subforest factorization first
gives
\[
 \mathcal D_s=
 \sum_{H_0\subset[s]}
 (f^{\mathcal A})^{\tensor([s]\setminus H_0)}
 \widetilde D_{H_0},
 \qquad H_0=\lambda_{\rm bad}\cup K_{\rm src};
\]
it does not yet give an \(n\)-based expansion.  The fixed term has the
source factor
\begin{equation}\label{eq:fixed-term-defect-credit}
 \vartheta^{(\mathsf X+m)/20}
 \eps^{d-1+1/(20d)+(\rho+r)/C_{15}^*}
 L^{C_\bullet(\rho+r+1)}
\end{equation}
and the forest count gives
\(|K_{\rm src}|\le C_\bullet\rho\).  Thus
\(|H_0|\le r+C_\bullet\rho\).  After increasing \(C_{\rm def,0}\), one half
of \((\rho+r)/C_{15}^*\) absorbs
\(L^{C_\bullet(\rho+r)}\), the other half gives
\(q_{\rm def}^{(|H_0|-1)_+}\), the fixed \(\eps^{d-1}\) absorbs the final
logarithmic constant, and \(1/(20d)\) gives \(a_{\rm err}\).  The unused
Gaussian reserve gives \(W_{H_0,t}\).

There is no analogous creation estimate with an arbitrary incoming paid
kernel substituted for \(K_{\rm src}\), and none is used.  What can occur
is the following exact disconnected product.  Let \(H_R\) be the output of
the red influence closure carrying an incoming paid kernel, and let
\(H_0\) be the affected set of a future-disconnected pure-source cutoff
term.  Exact subforest
factorization gives, before rebasing and symmetrization,
\begin{equation}\label{eq:red-source-defect-fusion}
 \bigl(\mathsf V_{H_R\leftarrow P}K_P\bigr)
 \tensor \widetilde D_{H_0}.
\end{equation}
We immediately declare \(H_R\cup H_0\) to be the affected set and regard
\eqref{eq:red-source-defect-fusion} as one sealed marked hyperedge; its two
displayed factors are never propagated separately.  By
\eqref{eq:common-energy-factorization},
\eqref{eq:sealed-block-bound}, and
\eqref{eq:fixed-term-defect-credit}, its weighted norm is bounded by the
product of the corresponding two bounds.  If \(H_0\ne\varnothing\), the
root powers already present are
\[
 (p-1)+( |H_R|-p)_+ + (|H_0|-1)_+
 \ge |H_R|+|H_0|-2.
\]
The missing root is paid by the incoming fixed factor
\(\eps^{a_*}\le q_{\rm def}\); all the other root factors are replaced by
\(q_{\rm def}\) using \(\zeta+q_{\rm b}\le q_{\rm def}\).  The pure-source
factor retains \(\eps^{a_{\rm err}}\).  This proves
\[
 \|W_{H_R\cup H_0,t}
   [\mathsf V_{H_R\leftarrow P}K_P\tensor\widetilde D_{H_0}]\|_1
 \le \Sigma_{\rm b}G\eps^{a_{\rm err}}
 q_{\rm def}^{(|H_R|+|H_0|-1)_+},
\]
after absorbing the fixed block constants in \(\Sigma_{\rm b}\).  If
\(H_0=\varnothing\) and \(H_R\ne\varnothing\), the scalar source defect
times the red closure is fused in the same way and enters the nonempty
marked bucket; the red factor already supplies
\(q_{\rm def}^{(|H_R|-1)_+}\).  Only
\(H_R=H_0=\varnothing\) enters the scalar bucket.  If \(H_R=\varnothing\),
the preceding pure-source estimate applies directly.

All cutoff-connected red-dependent histories, by contrast, are propagated
by the exact sealed operator of
Theorem~\ref{thm:sealed-past-block-stop-line}: internal red repeats are in
the exact group, ordinary one-port decorations are summed by
\eqref{eq:one-port-derivative-bound}, and all finite merge forests are
retained.  Thus the only new cutoff
factor in \eqref{eq:red-source-defect-fusion} is
\(\widetilde D_{H_0}\); source credit is never assigned to \(K_P\).
Previously created marked defects are inputs to the exact future sealed
operator, not new Source-B connected indices.

We verify the marked-input variant
\eqref{eq:marked-defect-fusion-map}.  Suppose
\[
 \|W_{P,a}D_P\|_1
 \le G^{\rm def}\eps^{a_{\rm err}}
       \eta^{(|P|-1)_+}.
\]
If no new source cutoff occurs, the sealed-block estimate propagates this
single marked input with root radius \(C(\eta+q_{\rm b})\).  If a new
future-disconnected source cutoff with affected set \(H_0\) occurs, exact
subforest factorization gives
\((\mathsf V_{H_R\leftarrow P}D_P)\tensor\widetilde D_{H_0}\), and we fuse
it immediately on \(H_R\cup H_0\).  For \(H_0\ne\varnothing\) the visible
root powers again total at least
\(|H_R|+|H_0|-2\).  Of the two fixed factors
\(\eps^{a_{\rm err}}\), retain one as the marked base and use the other via
\(\eps^{a_{\rm err}}\le q_{\rm def}\) to pay the missing root.  All
remaining powers are bounded by \(\eta+q_{\rm def}\), since
\(q_{\rm b}\le q_{\rm def}\).  A scalar new source defect times a nonempty
old marked closure is fused into the same nonempty bucket; only the wholly
empty case remains scalar.  In that case
\(H_R=H_0=\varnothing\), and the free-slice marginal contraction in
Lemma~\ref{lem:sealed-one-hyperedge}, followed by the subset and rebase sums,
gives scalar size at most
\(\Omega_{\rm b}G^{\rm def}\eps^{a_{\rm err}}\).  The block, subset and
rebase constants are exactly those included in \(\Omega_{\rm b}\), proving
\eqref{eq:marked-defect-fusion-map}.  Repeating this operation never creates
a second propagated hyperedge: its affected set is replaced by the union at
the first new cutoff slice.

Finally rebase the first display by
\(f^{\mathcal A}=n+u\).  The identity
\eqref{eq:exact-cumulant-rebase},
\(\|W_{\{1\}}u\|_1\le\Sigma_{\rm b}\eps^{a_*}\), and
\(\eps^{a_*}\le q_{\rm def}\) pay every newly added singleton and give
\eqref{eq:localized-defect-expansion}--\eqref{eq:localized-defect-bound}.
The exceptional case \(r=\rho=0\), in which the source formula uses one
unobserved proxy and \(H_0=\varnothing\), is precisely
\(D_\varnothing n^{\tensor s}\); its fixed \(\eps^{d-1+1/(20d)}\) pays the
displayed scalar bound.

In a later block the marked \(D_H\) is precisely the one sealed hyperedge
allowed in Theorem~\ref{thm:sealed-past-block-stop-line}; its future
pure-internal contacts are resummed before absolute values.  Therefore the
same activity recurrence applies and no aggregate contact operator is ever
placed on \(D_H\).  The global forest-cutoff remainder is a different,
two-sided algebraic remainder; its exact coefficient is still
\(2(d-1)\), so \eqref{eq:err-recurrence} is unchanged.
\end{proof}

\begin{lemma}[Structured global error propagation]
\label{lem:epoch-global-error-bucket}
For consecutive epochs using the same global cutoff tower, let
\((G_j^{\rm def},\eta_j)\) be the paid parameters of the sum of all
nonempty localized defects present after epoch \(j\), and let \(\sigma_j\)
be the corresponding normalized scalar-bucket size.  Let \(G_j\) be the
main paid amplitude at the same boundary, and put
\(E_j:=G_j^{\rm def}+\sigma_j\).  With \(\Omega_j\) denoting
\eqref{eq:block-Omega} for that epoch, after increasing its absolute
constant once to absorb the two augmented marked coordinates,
\begin{align}
 E_{j+1}&\le
 \Omega_j(E_j+G_j+1),
 &\eta_{j+1}&\le\Omega_j(\eta_j+q_{\rm def}),
 \label{eq:epoch-weighted-defect-sum}\\
 &&G_0^{\rm def}&=\eta_0=\sigma_0=E_0=0.
 \notag
\end{align}
The finite-correlation remainders obey the single concatenated recurrence
\eqref{eq:err-recurrence}.
\end{lemma}

\begin{proof}
Use Lemma~\ref{lem:localized-weighted-defect} at the creation epoch and
Theorem~\ref{thm:sealed-past-block-stop-line} at every later epoch, always
retaining the marked affected set.  Adjoin the empty kernel as a scalar
coordinate and the one-particle reference as a distinguished coordinate.
On this augmented hierarchy space rebasing and block propagation are
linear (the apparent affine source term is the image of the empty
coordinate).  Hence the exact telescoping identity
\begin{equation}\label{eq:epoch-telescoping}
 \mathcal U_B\cdots\mathcal U_1-\mathcal T_B\cdots\mathcal T_1
 =\sum_{j=1}^B\mathcal U_B\cdots\mathcal U_{j+1}
 (\mathcal U_j-\mathcal T_j)\mathcal T_{j-1}\cdots\mathcal T_1
\end{equation}
shows that there is exactly one first marked creation in every summand.
There are four augmented marked outputs to count.  An old scalar is
propagated as that scalar times the pure-reference block: its empty output
stays scalar, while every nonempty output
\(D_\varnothing K_H^{\rm new}\) is a sealed marked kernel.  The weighted
pure-reference estimate \eqref{eq:weighted-new-block-source},
\(q_{\rm b}\le q_{\rm def}\), and \(\eps^{a_*}\le1\) bound both kinds by
the block constant times \(\sigma_j\) at base
\(\eps^{a_{\rm err}}\).  The other possibility is that this old scalar
meets a later future-disconnected pure-source cutoff.  Exact subforest
factorization then gives
\(D_\varnothing\widetilde D_{H_0}\).  By
\eqref{eq:fixed-term-defect-credit}, retaining the old marked base and using
the new fixed \(\eps^{a_{\rm err}}\le1\), this costs at most the block
constant times \(\sigma_j\); it stays scalar when \(H_0=\varnothing\) and
enters the nonempty marked coordinate with radius \(q_{\rm def}\) otherwise.
Conversely, an old nonempty marked family has a
nonempty output controlled by \eqref{eq:marked-defect-fusion-map}, while its
completely marginalized output is scalar and is controlled by the last
assertion of Lemma~\ref{lem:localized-weighted-defect}; both cost the block
constant times \(G_j^{\rm def}\).  Finally, a new pure-source defect,
including fusion with the current main red closure, contributes at most the
block constant times \(G_j+1\) to the two coordinates by
\eqref{eq:localized-defect-bound}.  Adding the two coordinate bounds and
increasing the absolute constant in \(\Omega_j\) by a fixed factor gives
the first inequality in \eqref{eq:epoch-weighted-defect-sum}.  The radius
inequality is unchanged: old nonempty outputs use the fusion map, and
nonempty descendants of a scalar use \(q_{\rm b}\le q_{\rm def}\).  No
unstructured later \(L^1\) operator is used.  Concatenating the algebraic
forest-cutoff recurrences gives the last assertion because the coefficient
\(2(d-1)\) is independent of epoch boundaries.
\end{proof}

\begin{proposition}[Epoch-stacked endpoint theorem]
\label{prop:epoch-stacked-endpoint}
Assume \eqref{eq:loglog-main}.  Uniformly for every
\(t\in[0,t_{\rm fin}]\) and every integer \(1\le s\le L\), there is an
exact decomposition
\begin{equation}\label{eq:epoch-final-cumulant}
 f_s(t)=\sum_{H\subset[s]}
 n(t)^{\tensor([s]\setminus H)}\mathcal K_H^t
 +\mathcal R_s^t+\mathcal F_s^t,
 \qquad \mathcal K_\varnothing^t=1,
\end{equation}
on the hard-core domain, such that, for \(H\ne\varnothing\),
\begin{align}
 \|\mathcal K_H^t\|_1
 &\le \mathfrak S_{\rm ep}\eps^{a_*}
       q_{\rm ep}^{(|H|-1)_+},
 \label{eq:epoch-K-bound}\\
 \|\mathcal R_s^t\|_1+\|\mathcal F_s^t\|_1
 &\le \mathfrak S_{\rm ep}\eps^{a_{\rm err}},
 \qquad
 a_{\rm err}:=\min\left\{\frac1{20d},\frac1{100}\right\},
 \label{eq:epoch-error-bound}
\end{align}
where
\begin{equation}\label{eq:epoch-subpower-activity}
 \mathfrak S_{\rm ep}=\eps^{-o(1)},
 \qquad q_{\rm ep}=o(L^{-2}).
\end{equation}
All \(o(1)\) terms and the threshold are uniform in the full parameter range
of \eqref{eq:loglog-main}.
\end{proposition}

\begin{proof}
Put \(X=\log L\) and
\[
 M_{\alpha,A}:=\max(1,\alpha)\max(1,A),
 \qquad
 m=\max\{1,\lfloor\kappa_{\rm b}\sqrt X\rfloor\}.
\]
For \(t>0\), choose
\begin{equation}\label{eq:epoch-number}
 B=B(t):=\max\left\{1,
 \left\lceil C_{\rm b}
 \frac{M_{\alpha,A}\max(1,t)}m\right\rceil\right\},
 \qquad \Delta=t/B,
\end{equation}
and divide \([0,t]\) into \(B\) equal epochs, each of which is divided into
exactly \(m\) equal genuine time layers.  Choose \(C_{\rm b}\) large, then
\(\kappa_{\rm b}\) small, and finally \(c_0\) small relative to
\(\kappa_{\rm b}/C_{\rm b}\).  The two inequalities in
\eqref{eq:loglog-main} give
\begin{equation}\label{eq:epoch-count-bounds}
 M_{\alpha,A}\max(1,\Delta)\le
 \max\left\{\frac{2c_0}{\kappa_{\rm b}},\frac2{C_{\rm b}}\right\}m
 =:c_{\rm ep}m,
 \qquad B\le C\sqrt X,
 \qquad N:=Bm\le\kappa_{\rm g}X,
\end{equation}
where \(c_{\rm ep}\) is chosen so that the local effective parameter
\begin{equation}\label{eq:epoch-local-vartheta}
 \vartheta_{\rm b}:=
 \frac{C_\beta M_{\alpha,A}\max(1,\Delta)}{m}\le q_0,
\end{equation}
and with fixed \(\kappa_{\rm g}>0\) as small as required below.
Indeed, if \(B=1\), the first bound follows from the definition of \(B\);
if \(B>1\), then
$M_{\alpha,A}\Delta\le m/C_{\rm b}$ and
$M_{\alpha,A}\le 2c_0m/\kappa_{\rm b}$ for small \(\eps\).
Also \(B\le1+C_{\rm b}c_0X/m\) and
$Bm\le m+C_{\rm b}c_0X$, which prove the remaining two bounds in the
display after the stated choice order.  No zero-length or empty layer is
inserted.

Use one global cutoff tower
\begin{equation}\label{eq:global-epoch-cutoff-tower}
 A_N=L_\eps^\#,
 \qquad \Lambda_j=A_j^{10d},
 \qquad A_{j-1}=\Lambda_j^{10d}\quad(1\le j\le N).
\end{equation}
Thus the output particle range of one epoch is exactly the input range of the
next.  Its largest forest exponent is bounded by
\begin{equation}\label{eq:global-forest-depth}
 C_{\rm init,glob}^*
 \le C\bigl(1+(100d^2)^N\bigr)\le L^{\delta_0},
\end{equation}
after choosing
\(\kappa_{\rm g}\log(100d^2)<\delta_0\).  On the other hand, the exceptional
constants are restarted and see only \(m\) layers.  Lemma
\ref{lem:appE-sourceA-Cj} gives
\begin{equation}\label{eq:block-star-depth}
 \max_j C_{j,{\rm b}}^*\le C^{C\Gamma m^2}\le L^{\delta_0}
\end{equation}
after choosing \(\kappa_{\rm b}\) small.  Enlarge
\(C_\bullet\) to include \eqref{eq:global-forest-depth}; then
\(C_\bullet\le L^{\delta_0}\).

Let \((\Sigma_{\rm b},q_{\rm b},\mathfrak R_{\rm b})\) be
\eqref{eq:block-S-q} and let \(\Omega_{\rm b}\) be
\eqref{eq:block-Omega}, using the worst block value.  Starting with the
exact initial grand-canonical ensemble, apply
Lemma~\ref{lem:paid-cumulant-restart} successively.  After each of the first
\(B-1\) blocks, use the exact fresh-reference conversion
\eqref{eq:fresh-reference-boundary-conversion} before the later block; after
the last block, retain the already paid product-coordinate output.
The recurrence
\[
 G_{j+1}\le\Omega_{\rm b}(G_j+1),\qquad
 \zeta_{j+1}\le\Omega_{\rm b}(\zeta_j+q_{\rm b}),\qquad
 G_0=\zeta_0=0,
\]
gives
\begin{equation}\label{eq:epoch-recursion-solved}
 G_B\le B\Omega_{\rm b}^{B+1},\qquad
 \zeta_B\le B\Omega_{\rm b}^{B+1}q_{\rm b}.
\end{equation}
For the cumulants and localized defects use
\begin{equation}\label{eq:epoch-S-q-definitions}
 \mathfrak S_{\rm ep}:=B^2\Omega_{\rm b}^{B+2},\qquad
 q_{\rm ep}:=B\Omega_{\rm b}^{B+1}q_{\rm b},\qquad
 \eta_{\rm ep}:=B\Omega_{\rm b}^{B+1}q_{\rm def},\qquad
 \sigma_{\rm ep}:=B\Omega_{\rm b}^{B+1}.
\end{equation}
The factors from \(\mathfrak d\) telescope:
\[
 \prod_{j=1}^B\mathfrak d(t_j,t_{j-1})
 =\exp\!\left\{C_{\beta,d}
 \bigl((1+\nu t)^3-1\bigr)\right\}.
\]
For the equal epochs in \eqref{eq:epoch-number}, the mean-value theorem gives
\[
 B\max_j\log\mathfrak d(t_j,t_{j-1})
 \le 3C_{\beta,d}\nu t(1+\nu t)^2
 \le 3C_{\beta,d}\bigl((1+\nu t)^3-1\bigr).
\]
After enlarging the generic constant, using one worst block value in the
radius recurrence therefore loses no more than the same cubic kinetic clock.
By \eqref{eq:loglog-main},
\(1+\nu t\le C(1+\log L)\).
The shifted-profile factor is used at most once per block.  Hence, from
\eqref{eq:epoch-count-bounds}, \eqref{eq:block-star-depth}, and
\(\delta_0=1/20\),
\begin{equation}\label{eq:epoch-total-loss-log}
 \log\mathfrak S_{\rm ep}
 \le C L^{\delta_0}(\log L)^2
      +C(\log L)^3=o(L).
\end{equation}
Moreover
\begin{align}
 \log q_{\rm ep}
 &\le-L^{1-\delta_0}
       +C L^{\delta_0}(\log L)^2
       +C(\log L)^3
 \le-\tfrac12L^{1-\delta_0},
 \label{eq:epoch-q-log}\\
 \log\eta_{\rm ep}
 &\le-L^{1-2\delta_0}
       +C L^{\delta_0}(\log L)^2
       +C(\log L)^3
 \le-\tfrac12L^{1-2\delta_0}.
 \label{eq:epoch-defect-radius-log}
\end{align}
Here \(C_{\rm def}\le CL^{2\delta_0}\) follows from
\eqref{eq:defect-root-constant}.  In particular
\(q_{\rm ep}+\eta_{\rm ep}=o(L^{-2})\), and at every intermediate boundary
\(\Omega_{\rm b}\zeta_j\le1/4\) and
\(\zeta_j+q_{\rm b}\le q_{\rm def}\) for small \(\eps\).  To verify the
latter directly, use \(\zeta_j\le q_{\rm ep}\) and
\[
 \log\frac{q_{\rm ep}}{q_{\rm def}}
 =-L\left(\frac1{C_{14}^*}-\frac1{C_{\rm def}}\right)
   +\log\bigl(B\Omega_{\rm b}^{B+1}\bigr)
 \le-\frac{L}{2C_{14}^*}<0.
\]
Here \(C_{\rm def,0}\) was chosen at least four and
\eqref{eq:block-star-depth}--\eqref{eq:epoch-total-loss-log} make the last
logarithm \(o(L/C_{14}^*)\).  The same comparison gives
\(q_{\rm b}\le q_{\rm def}\), and
\(\Omega_{\rm b}(\eta_j+q_{\rm def})<1\) follows from
\eqref{eq:epoch-defect-radius-log}.  Thus
\eqref{eq:paid-radius-admissibility} is verified before each restart,
without using the conclusion circularly.  This proves
\eqref{eq:epoch-K-bound} and \eqref{eq:epoch-subpower-activity}.
Combining \eqref{eq:epoch-recursion-solved} with the first recurrence in
\eqref{eq:epoch-weighted-defect-sum} gives a sharper reserve which we record
for the singleton terms.  After increasing the absolute constant in
\(\Omega_{\rm b}\), we have \(\Omega_{\rm b}\ge4\) and
\(\Omega_{\rm b}^2\ge8L\).  Induction in the two recurrences gives
\[
 G_j\le2\Omega_{\rm b}^j,\qquad
 E_j\le4j\Omega_{\rm b}^j.
\]
The nonempty marked amplitude satisfies
\(G_B^{\rm def}\le E_B\le4B\Omega_{\rm b}^B\), while
\eqref{eq:epoch-S-q-definitions} and
\eqref{eq:epoch-defect-radius-log} give the stated terminal radius.
Since a singleton defect has no radius factor, the correct subset sum is
\[
 \sum_{\varnothing\ne H\subset[s]}
 \eta_{\rm ep}^{|H|-1}
 =\frac{(1+\eta_{\rm ep})^s-1}{\eta_{\rm ep}}
 \le2s\le2L,
\]
with the value at \(\eta_{\rm ep}=0\) understood by continuity.  Hence the
nonempty defect total is at most
\[
 8LB\Omega_{\rm b}^B\eps^{a_{\rm err}}
 \le B^2\Omega_{\rm b}^{B+2}\eps^{a_{\rm err}}
 =\mathfrak S_{\rm ep}\eps^{a_{\rm err}}.
\]
The scalar buckets satisfy
\(\sigma_B\le E_B\le4B\Omega_{\rm b}^B\le\sigma_{\rm ep}\), and hence
contribute at most
\(\sigma_{\rm ep}\eps^{a_{\rm err}}\le
\mathfrak S_{\rm ep}\eps^{a_{\rm err}}\) by the coupled recurrence.
For the finite-correlation remainder repeat the proof of
Lemma~\ref{lem:restart-domination} for the concatenated tower
\eqref{eq:global-epoch-cutoff-tower}, with \(N\) replacing
\(\mathfrak L\), the local parameter \(\vartheta_{\rm b}\) from
\eqref{eq:epoch-local-vartheta}, and the fixed main-theorem value
\(M=M_0=1/100\).  That proof uses only the adjacent cutoff gap,
\(\vartheta_{\rm b}\le q_0\), and the exact coefficient \(2(d-1)\), so it is
unchanged at epoch boundaries and bounds the global recurrence remainder by
\(\eps^{M_0}\).  In addition, each internal conversion
\(j=1,\ldots,B-1\) in
\eqref{eq:fresh-reference-boundary-conversion} contributes one signed
reference remainder with weighted norm at most
\(\eps^{\Lambda_{jm}/10}\) at that boundary.  The initial-term calculation in
Lemma~\ref{lem:restart-domination}, with its indices shifted to \(jm\),
bounds each such remainder after all later cutoff losses by
\(\eps^{M_0+1}\).
There are at most \(B-1\le C\sqrt X\) of them (and none when \(B=1\)), so
their sum is at most
\(\eps^{M_0}\) after the same reduction of \(\eps_0\).  This proves
\eqref{eq:epoch-error-bound}.  At \(t=0\) use the initial forest expansion.
Finally, the possible values of \(B(t)\) are bounded by \(C\sqrt X\), all
estimates use their worst value, and the time cells are half open.  One
common reduction of \(\eps_0(d,\beta)\) is uniform in \(t\).

\end{proof}

\subsection{The one-epoch output estimates}

Assume throughout this subsection the one-epoch condition
\eqref{eq:one-epoch-horizon}.  The estimates below are retained regression
outputs for one structural run; they are not the full-range ledger of
Proposition~\ref{prop:epoch-stacked-endpoint}.
Fix $t\in[0,t_{\rm fin}]$, apply the rerun of
Lemma~\ref{lem:uniform-terminal-rerun} with $T=t$, and suppress the
superscript $T$ on its expansion objects.  Corollary
\ref{cor:diagonal-reduction}, with the inputs proved in
Sections~\ref{sec:ensemble-flow}--\ref{sec:p38j}, gives the explicit bounds
\begin{align}
 \|E_H(t)\|_1
 &\le \mathscr S_\eps\eps^{1/(15d)}q_{\rm eff}^{|H|},
 \quad H\ne\varnothing,
 \label{eq:global-EH}\\
 \|f_s^{\rm err}(t)\|_1
 &\le\mathscr S_\eps\eps^{1/(20d)},
 \label{eq:global-trunc}\\
 \|f^{\mathcal A}(t)-n(t)\|_1
 &\le\mathscr S_\eps\eps^{7/(1000d)}.
 \label{eq:global-approx}
\end{align}
There is no suppressed layer-dependent factor in these displays.  The
single factor $\mathscr S_\eps$ is defined in \eqref{eq:subpower-log}, and
all per-root losses are in $q_{\rm eff}$.  The exponent in
\eqref{eq:global-approx} is
\[
 \frac12\left(\frac1{200d}+\frac9{1000d}\right)
 =\frac7{1000d}.
\]

\subsection{Final aggregation}
\label{subsec:final-aggregation}

The nonempty cumulants, including the per-root logarithmic loss and the
tensor factor on the complementary labels, sum to
\begin{align}
 \sum_{\varnothing\ne H\subset[s]}
 \|(f^{\mathcal A}(t))^{\tensor([s]\setminus H)}E_H(t)\|_1
 &\le \sum_{h=1}^s \binom{s}{h}
 \|f^{\mathcal A}(t)\|_1^{s-h}
 \mathscr S_\eps\eps^{1/(15d)}q_{\rm eff}^{h}\notag\\
 &\le \max\{1,\|f^{\mathcal A}(t)\|_1\}^{s}
 \mathscr S_\eps\eps^{1/(15d)}
 \big((1+q_{\rm eff})^s-1\big)\notag\\
 &\le 4s\,\mathscr S_\eps\eps^{1/(15d)}q_{\rm eff},
 \label{eq:cumulant-subset-sum}
\end{align}
where \eqref{eq:global-approx}, mass conservation,
\eqref{eq:qeff-bound} and $s\le L$ were used.  Indeed
$\|f^{\mathcal A}(t)\|_1
\le1+\mathscr S_\eps\eps^{7/(1000d)}$, so its $s$th power is
$1+o(1)$ uniformly for $s\le L$; that harmless factor was absorbed on the
second line.

Mass conservation gives $\|n(t)\|_1=1$.  The tensor telescoping identity and
\eqref{eq:global-approx} imply
\begin{equation}\label{eq:tensor-replacement}
 \|(f^{\mathcal A})^{\tensor s}-n^{\tensor s}\|_1
 \le s\,\mathscr S_\eps\eps^{7/(1000d)}
       (1+\mathscr S_\eps\eps^{7/(1000d)})^{s-1}.
\end{equation}
Finally, the union bound over pairs and
$\sup_x\int n(t,x,v)\dd v\le C_\beta A$ give
\begin{equation}\label{eq:excluded-volume}
\|n(t)^{\tensor s}(1-\1_{\cD_s^\circ})\|_1
 \le C_{\beta,d} A s^2\eps^d.
\end{equation}

For the longer range we use Proposition
\ref{prop:epoch-stacked-endpoint}.  Since \(\|n(t)\|_1=1\), its nonempty
cumulants satisfy
\begin{align}
 \sum_{\varnothing\ne H\subset[s]}
 \|n(t)^{\tensor([s]\setminus H)}\mathcal K_H^t\|_1
 &\le \mathfrak S_{\rm ep}\eps^{a_*}
 \sum_{h=1}^s\binom{s}{h}q_{\rm ep}^{(h-1)_+}\notag\\
 &\le 2s\mathfrak S_{\rm ep}\eps^{a_*},
 \label{eq:epoch-cumulant-subset-sum}
\end{align}
because \(s q_{\rm ep}=o(1)\) for \(s\le L\).  The singleton
term is kept explicitly; no nonexistent per-root factor is assigned to the
rebased one-particle error.

\begin{proof}[Proof of Theorem~\ref{thm:main}]
Apply Proposition~\ref{prop:epoch-stacked-endpoint} at the selected time
\(t\).  Extending the microscopic correlation by zero off the hard-core
domain and using the zero-measure contact boundary, its exact identity
\eqref{eq:epoch-final-cumulant} gives
\begin{align}
 &\|f_s(t)-n(t)^{\tensor s}\1_{\cD_s^\circ}\|_1\notag\\
 &\quad\le
 \|\mathcal R_s^t\|_1+\|\mathcal F_s^t\|_1
 +\sum_{\varnothing\ne H\subset[s]}
   \|n(t)^{\tensor([s]\setminus H)}\mathcal K_H^t\|_1
 +\|n(t)^{\tensor s}(1-\1_{\cD_s^\circ})\|_1.
 \label{eq:final-five-term-bound}
\end{align}
The terms are controlled by \eqref{eq:epoch-error-bound},
\eqref{eq:epoch-cumulant-subset-sum}, and \eqref{eq:excluded-volume}.
Proposition~\ref{prop:epoch-stacked-endpoint} makes the bounds uniform for
$0\le t\le t_{\rm fin}$; all constants and $\eps_0$ are independent of
$t$ and $1\le s\le L$.

The exponents \(a_*=7/(1000d)\), \(a_{\rm err}\), and \(d\) have strict
room above \(1/(400d)\).  Equations
\eqref{eq:epoch-total-loss-log}--\eqref{eq:epoch-q-log} give the common
subpower factor and absorb every repeated paid-root loss.  Reducing
$\eps_0$ gives
\[
 \|f_s(t)-n(t)^{\tensor s}\1_{\cD_s^\circ}\|_1
 \le\eps^{1/(400d)}.
\]
\end{proof}
\section{Scope and Further Consequences}
\label{sec:scope}

The proof is organized around six structural choices that determine the scope of
the main theorem and the way its local estimates enter the long-time argument.

\subsection{Exact scope of the proof}

\begin{enumerate}[label=\textup{(\roman*)}]
\item The serial $j=3$ input is the shell estimate
\eqref{eq:j3-shell}.  At the critical radial scale, the higher-dimensional
closure for the long degree-$(3,3)$ component is supplied by the joint two-sublayer packet mechanism.
\item All two-root capacity comes from the explicit physical first-cross
residual \eqref{eq:landing-residual}; the resulting transversality is therefore
tied directly to the collision geometry used in the packet construction.
\item In fixed dimension $d\ge4$, Theorem~\ref{thm:p38j} replaces the third
alternative of \cite[Proposition 3.8]{DHMtorus}, and
Proposition~\ref{prop:v2-consumer} proves the exact operator-level form needed
by the downstream consumer.
\item Theorem~\ref{thm:two-landing} retains the entire connected two-sublayer
owner closure as the analytic object.  This packet-level formulation includes
long triangles and strongly degenerate configurations, while the outer
Source--A dichotomy may dispatch a strongly degenerate branch before the
theorem is invoked.
\item The imported long-time architecture is version~2 of \cite{DHMlong}, the
version cited by \cite{DHMtorus}; its exponent ledger is used consistently
throughout the proof.
\item The epoch interface acts on complete joint kernels in the
time-decreasing common-energy Gaussian paid class
\eqref{eq:paid-cumulant-norm}.  Pure internal contacts are first resummed by
Lemma~\ref{lem:paired-first-contact-partition}; the estimate is then applied
to the paired time-integrated interblock operator.  The packet cap keeps the adaptive integer in
\eqref{eq:adaptive-fixed-k} bounded independently of \(\eps\).  The first
condition in \eqref{eq:loglog-main} retains the one-epoch
\((\alpha,A)\)-range.  The added range is in total time and becomes
\(t_{\rm fin}=O(\log|\log\eps|)\) when \(\alpha,A\) are fixed.
\end{enumerate}

\subsection{Joint-kernel coarea and branch ownership}

On every packet cell, four properties keep the coarea and all ownership
charges inside one joint operator.
\begin{itemize}
\item The raw block \eqref{eq:raw-landing-Jacobian} identifies the two long
scales, but the actual coarea is one joint map.  The default FCT chart uses
$2(d-1)$ separator velocities and the two landing times.  Its radial row is
only a local exterior augmentation; Cauchy--Binet returns ordinary Cartesian
coordinates.  The full $2d$ disjoint chart and the collision-free new-line
position--velocity chart provide complementary coarse bounds.  Every chart
acts on the same joint kernel.
\item The selected maximal minor supplies local transversality, while the complete
collision word is retained as the polynomial system
\eqref{eq:appD-full-polynomial-system}; Theorem
\ref{thm:appD-uniform-algebraic-multiplicity} enumerates all its regular real
solutions by the fixed slots $1\le j\le M(d,P_0)$.  Thus the packet operator is the finite sum over all branch labels.
\item Every boundary state remains a parameter in the same nonnegative integrand,
and the output domain is enlarged once independently of its concrete value.
\item The entire complement is inserted as $Q$ before the packet is integrated, so
each source coordinate is charged within the same joint estimate.
\end{itemize}

Lemma~\ref{lem:upper-tree-landings} chooses the first two lower atoms.  Their
edges are either disjoint, when Proposition~\ref{prop:disjoint-radial-fct}
applies, or overlap on three lines, when
Proposition~\ref{prop:overlap-radial-fct} uses the line first appearing at
$d_2$.  This dichotomy exhausts every packet cell.  The full-chord, short-gap,
three-line and tree-only propositions provide complementary coarse checks.

\subsection{Normalization and the two landing gains}

The full packet has the unique negative baseline
\eqref{eq:unique-baseline-setup}.  The two positive $\eps^{d-1}$ factors in
\eqref{eq:two-landing-excess} belong to two distinct lower collision
contacts.  In the default FCT chart their two fluxes cancel the live landing
time columns exactly once, and the inverse gap cost is
$\eps_*^{-2(d-1)}$.  The disjoint $2d$-velocity coarea and the overlapping
full-chord plus $d$-row landing coarea instead cost $\eps_*^{-2d}$ and are
retained as independent coarse checks.
Upper tree collisions cancel particle-source normalizations.  In the
periodic overlap repair, protection makes the first-failure assignment
injective, and a recut transfers the same certificate from a pair to a
singleton representative.  The final theorem rate is the minimum of the
disjoint error-family exponents.

\subsection{Finite fibres and downstream propagation}

Global recovery combines local rank with uniform finite-fibre control.  For
example, $p\mapsto\sin p$ has nonzero derivative at every zero but several
zeros on a long interval.  Lemma~\ref{lem:appD-maximal-minor-atlas} supplies
the local transversality certificate, and the polynomial fibre bound below
supplies the global branch count.
For each coarse disjoint, coarse overlapping or default radial $\mathrm{FCT}$ cell, the full
hard-sphere word instead supplies a square polynomial system of degree at
most $4P_0+4$ in at most $(2d+1)P_0+2d$ pivot variables.  Every regular real fibre
is therefore finite with the uniform bound
\[
 M(d,P_0)=(4P_0+4)^{(2d+1)P_0+2d}.
\]
The branch label is part of the fixed enlarged output domain.  Proposition
\ref{prop:appD-arbitrary-kernel-coarea} sums every algebraic branch within
the chosen chart, with the chart-specific gap exponent
\(\gamma_\kappa\) defined in Section~\ref{app:coarea}, and Lemma
\ref{lem:appD-fibre-time-owner} retains on each nonzero branch the same
original packet-time witness.  The sharp downstream route chooses the
default radial FCT atlas, for which \(\gamma_\kappa=2(d-1)\); the independent
coarse atlas retains \(\gamma_\kappa=2d\).  Hence the finite-branch multiplier contributes only a constant depending on
the fixed packet cap and carries no additional $\eps$ power.  This
correction is propagated through Propositions
\ref{prop:disjoint-rooted-extension} and
\ref{prop:asymmetric-new-line-peeling}, Theorem
\ref{thm:two-landing}, Theorem~\ref{thm:p38j}, Proposition
\ref{prop:v2-consumer}, and finally Theorem~\ref{thm:main}.

\appendix
\section{Parameter Order and Arithmetic Closure}
\label{app:parameters}

This appendix collects the numerical comparisons used in Sections~7--8.  It
also records the order of choices and each absorption of a growing constant.

\subsection{Order of choices}

Fix an integer $d\ge4$ and the Maxwellian weight parameter \(\beta>0\).
Before any growing
constant is selected, fix once and for all
\[
 \delta_0=\frac1{20},\qquad M_0=\frac1{100},
\]
and the finite list \(\boldsymbol{\eta}_{\rm main}\) of positive subpower
margins actually used in the proof of Theorem~\ref{thm:main}.  This list
includes, in particular, one margin smaller than \(9/(4000d)\) for the final
aggregation.  The generic lemmas below remain valid for other fixed values
of \(M\) or \(\eta\), but their thresholds are not included in the single
main-theorem threshold.  The constants are selected in the following order:
\begin{equation}\label{eq:appE-choice-order}
 \begin{gathered}
 d,\beta\longrightarrow
 \delta_0,M_0,\boldsymbol{\eta}_{\rm main}
 \longrightarrow
 \Gamma\in\mathbb N,\quad \widehat\Gamma:=\Gamma+\tfrac12
 \\
 \longrightarrow
 P_{\rm pkt}(d),(a_j)_{1\le j<15},C_{\rm ord,0},C_{\rm init,0},C_{\rm own}
 \longrightarrow \kappa\longrightarrow C_{\rm b}
 \longrightarrow \kappa_{\rm g}\longrightarrow \kappa_{\rm b}
 \longrightarrow c_0\longrightarrow \varepsilon_0.
 \end{gathered}
\end{equation}
Here the integer $\Gamma$ is the Source-B recollision threshold and is chosen
to satisfy $\Gamma>(60d)^{60d}K_D$ from
\eqref{eq:main-case-2-cutoffs};
$\widehat\Gamma=\Gamma+\tfrac12$ is the strict graph-partition threshold used
only in Theorem~\ref{thm:p38j}, and the bounded packet cap $P_{\rm pkt}(d)$ comes from
the finite cutting algorithms.  Both thresholds are fixed before every
growing parameter.  The $a_j$ are the fixed-power dominance exponents.  The
fixed constant $C_{\rm own}$ dominates the finitely many conditional
owner/good constants in Proposition
\ref{prop:packet-free-small-window-consumer}.  It is chosen before $\kappa$
and is independent of every terminal parameter.  The
constant $\kappa=\kappa(d,\beta)$ is fixed after $\delta_0$, small enough that
the quadratic growth exponent is below $\delta_0$; \(c_0\) is chosen only
after the epoch constants.  More precisely, \(C_{\rm b}\) is made large
enough for the block smallness in \eqref{eq:epoch-count-bounds}; then
\(\kappa_{\rm g}\) is made small enough for
\(\kappa_{\rm g}\log(100d^2)<\delta_0\); then \(\kappa_{\rm b}\) is made
small enough for \eqref{eq:block-star-depth} and relative to
\(\kappa_{\rm g}\).  Finally \(c_0\) is made small relative to
\(\kappa,\kappa_{\rm b}/C_{\rm b}\), and \(\kappa_{\rm g}/C_{\rm b}\), so
that all atom series have ratio below \(q_0\) and the global layer count in
\eqref{eq:epoch-count-bounds} holds.
For a given
$\eps$, only after those choices, the integer $\mathfrak L$ and the actual
hierarchy constants $C_j^*(\mathfrak L)$ are evaluated by
\eqref{eq:appE-Cj-actual-recursion}.  Finally $\varepsilon_0$ is reduced to
validate all asymptotic comparisons simultaneously.  In particular, no
constant is chosen after seeing $t_{\rm fin}$.

\subsection{Growth of the fifteen exceptional constants}

The fixed-time convention in Source B gives an intentionally extravagant
example and cannot itself be placed on the present diagonal.  The
quantitative argument is Source A,
\cite[proof of the main theorem, Part~7]{DHMtorus}.  A line-by-line comparison identifies
exactly two layer-dependent dominance requirements:
\begin{equation}\label{eq:appE-two-layer-dominances}
 C_{10}^*\gg(C_7^*)^{\mathfrak L},
 \qquad
 C_6^*\gg(C_5^*)^{\Gamma\mathfrak L}.
\end{equation}
They arise respectively from the layer selector and from the at most
$\Gamma\mathfrak L$ thin layers.  Every other transition has the form
$C_{j+1}^*\ge (C_j^*)^{a_j}$ with a fixed integer $a_j$, independent of
$\mathfrak L$.  This is precisely the dependence asserted in Source A,
\cite[proof of the main theorem, Part~7]{DHMtorus}, in the displayed
dominance conditions and the ensuing bound
$C_j^*\le C^{C\Gamma\mathfrak L^2}$; it replaces, for the quantitative
theorem, the nonquantitative example in Source B's parameter convention.

Let $\mathscr E_{\rm src}$ be the finite set of fixed powers occurring in
the Source-B v2 dominance conditions for the transitions
$C_j^*\mapsto C_{j+1}^*$, after removing the two layer-dependent powers
displayed in \eqref{eq:appE-two-layer-dominances}.  Equivalently, this is the
finite list denoted by the absolute exponent $C$ in Source A, Part~7,
item~(3).  Fix once and for all
\[
 a_{\rm src}:=\max\bigl(2,\lceil a\rceil:a\in\mathscr E_{\rm src}\bigr)
\]
and take every $a_j\ge a_{\rm src}$.  Thus the symbols $a_j$ below dominate
every fixed-power source requirement and are independent of
$\mathfrak L$ and $\eps$.  Choose a sufficiently large absolute $C_1^*$.
The complete recursion is
\begin{equation}\label{eq:appE-Cj-actual-recursion}
 C_{j+1}^*=\begin{cases}
 \max\{(C_5^*)^{a_5},(C_5^*)^{a_5\Gamma\mathfrak L}\},&j=5,\\
 \max\{(C_9^*)^{a_9},(C_7^*)^{a_9\mathfrak L}\},&j=9,\\
 (C_j^*)^{a_j},&j\notin\{5,9\}.
 \end{cases}
\end{equation}
After increasing one absolute exponent $C_0$, define the single free cutoff
coefficient
\begin{equation}\label{eq:appE-Cstar-definition}
 C^*(\mathfrak L):=C^{C_0\Gamma\mathfrak L^2}.
\end{equation}
It dominates all fifteen $C_j^*(\mathfrak L)$ and is the quantity denoted by
$C^*$ in every logarithmic packet weight below.  Thus no occurrence of
$C^*$ is an unevaluated ``sufficiently large'' constant on the diagonal.

Two further routed exponents have to be included.  Let
\begin{equation}\label{eq:Cord-definition}
 C_{\rm ord}^*(\mathfrak L)
 :=C_{\rm ord,0}\bigl(1+C^*(\mathfrak L)\bigr),
\end{equation}
where the absolute $C_{\rm ord,0}$ dominates the fixed degrees with which
the velocity cutoff, lift count, dyadic cells and fixed enlarged domains
occur in the ordinary operators.  Thus a displayed ordinary loss
$\cL_\eps^C$ is routed globally as at most
$\cL_\eps^{C_{\rm ord}^*}$; this convention does not turn a genuinely finite
chart count into a growing one.  The initial-forest exponent is the explicit
quantity from \eqref{eq:Cinit-definition},
\begin{equation}\label{eq:Cinit-recalled}
 C_{\rm init}^*(\mathfrak L)
 =C_{\rm init,0}\bigl(1+(100d^2)^{\mathfrak L}\bigr).
\end{equation}
Define the single global exponent
\begin{equation}\label{eq:Call-definition}
 C_{\rm all}^*(\mathfrak L)
 :=\max\bigl\{C_1^*,\ldots,C_{15}^*,
               C_{\rm ord}^*,C_{\rm init}^*\bigr\}.
\end{equation}
The absolute prefactors in \eqref{eq:Cord-definition} and
\eqref{eq:Cinit-definition} are chosen once, before $\kappa$, large enough
that $C_{\rm all}^*$ absorbs every fixed multiple of a per-root logarithmic
exponent used in the paper.

\begin{lemma}[Source-A quantitative constant bound]
\label{lem:appE-sourceA-Cj}
For $1\le j\le15$, the construction
\eqref{eq:appE-Cj-actual-recursion} satisfies all source dominance
conditions and
\begin{equation}\label{eq:appE-Cj-quadratic-growth}
 C_j^*\le C^{C\Gamma\mathfrak L^2}.
\end{equation}
Moreover,
\begin{equation}\label{eq:appE-Call-quadratic-growth}
 C_{\rm all}^*(\mathfrak L)
 \le C^{C\Gamma\mathfrak L^2}.
\end{equation}
The bound includes the two actual layer-dependent feedbacks; none is hidden
as a mere cell multiplicity.
\end{lemma}

\begin{proof}
Before stage six, $\log C_j^*=O(1)$.  The first exceptional transition in
\eqref{eq:appE-Cj-actual-recursion} gives
$\log C_6^*=O(\Gamma\mathfrak L)$.  Fixed-power transitions preserve this
order through $C_9^*$.  Since $\log C_7^*=O(\Gamma\mathfrak L)$, the second
exceptional transition gives $\log C_{10}^*=O(\Gamma\mathfrak L^2)$.
The remaining fixed-power transitions preserve the latter order.  Enlarging
the absolute base $C$ proves \eqref{eq:appE-Cj-quadratic-growth} and the
strict inequalities in \eqref{eq:appE-two-layer-dominances}.
Equation \eqref{eq:Cord-definition} has the same quadratic-exponential
upper bound.  Also
$C_{\rm init}^*\le C\exp(\log(100d^2)\,\mathfrak L)$, which is smaller than
$C^{C\Gamma\mathfrak L^2}$ for $\mathfrak L\ge1$ after another absolute
enlargement.  Taking the maximum proves
\eqref{eq:appE-Call-quadratic-growth}.
\end{proof}

\begin{lemma}[Uniform one-epoch quantifier form]
\label{lem:appE-uniform-quantifiers}
After fixing the finite numerical data in \eqref{eq:appE-choice-order}, namely
\[
 (\delta_0,M_0,\boldsymbol{\eta}_{\rm main},\Gamma,P_{\rm pkt},
  (a_j),C_{\rm ord,0},C_{\rm init,0},C_{\rm own},
  \kappa,C_{\rm b},\kappa_{\rm g},\kappa_{\rm b},c_0),
\]
there is one
\(\varepsilon_0=\varepsilon_0(d,\beta)>0\) such that every estimate
\emph{invoked in one application of the fixed structural reduction with
those fixed numerical data} holds simultaneously for every
\(0<\varepsilon\le\varepsilon_0\), every admissible triple
\((\alpha,A,t_{\rm fin})\) satisfying
\[
 \max(1,\alpha)\max(1,A)\max(1,t_{\rm fin})
 \le c_0\sqrt{\log|\log\varepsilon|}.
\]
The simultaneity includes every $T\in[0,t_{\rm fin}]$,
$1\le s\le|\log\varepsilon|$, and every subset $H\subset[s]$ occurring in
the terminal rerun.  It does not assert one threshold uniform over all
possible auxiliary choices $M>0$ or $\eta>0$ in the generic asymptotic
lemmas.
\end{lemma}

\begin{proof}
From
\eqref{eq:appE-Call-quadratic-growth} and
$\mathfrak L^2\le\kappa^2\log L$,
\[
 C_{\rm all}^*\le L^{C_1\Gamma\kappa^2}.
\]
Choose $\kappa$ so that $C_1\Gamma\kappa^2<\delta_0$.  Then all fifteen
source constants, the ordinary routed exponent and the initial-forest
exponent are simultaneously at most $L^{\delta_0}$ after one common reduction
of $\eps_0$.
The same reduction gives
\begin{equation}\label{eq:appE-owner-diagonal}
 C_{\rm own}\mathfrak L\eps^{1/(72d)}\le1,
\end{equation}
because \(\mathfrak L=O(\sqrt{\log L})\).  This is the fixed-dimensional
owner/good margin and is uniform in \((\alpha,A,T)\).
There are finitely many
\emph{forms} of smallness condition,
although their coefficients vary along the diagonal.  The only form in
which the discrete count $R$ grows together with such a coefficient is the
two-regime envelope; Lemma~\ref{lem:source-two-regime-envelope}, especially
\eqref{eq:source-uniform-growing-constants}, proves its threshold uniformly
by the estimate $L^{2\delta_0}\log L=o(L)$.  Repeated one-root and one-layer
losses are uniform by \eqref{eq:appE-growing-subpower} and
\eqref{eq:appE-qeff-log}.  All remaining conditions are monotone smallness
requirements, applications of the subpower ledger with
$\eta\in\boldsymbol{\eta}_{\rm main}$, or the single application of
Lemma~\ref{lem:restart-domination} with $M=M_0$.  This is a finite collection
after the numerical data have been fixed.  Taking the minimum of these
uniform thresholds gives the claim, without taking an infimum over arbitrary
$M$ or over margins $\eta\downarrow0$.
\end{proof}

\subsection{Subpower scales}

Put
\[
 L=|\log\eps|,
 \qquad
 \epstar=\exp(-\sqrt L),
 \qquad
 \mathfrak L=\max\{1,\lfloor\kappa\sqrt{\log L}\rfloor\}.
\]

\begin{lemma}[Finite subpower ledger]
\label{lem:appE-subpower-ledger}
For fixed \(a,b,M>0\) and every \(\eta>0\), after reducing
\(\varepsilon_0\),
\begin{align}
 \epstar^{-a}&\le\eps^{-\eta},
 \label{eq:appE-epstar-subpower}\\
 L^b&\le\eps^{-\eta},
 \label{eq:appE-log-subpower}\\
 \mathfrak L^M L^b\epstar^{-a}&\le\eps^{-\eta}.
 \label{eq:appE-combined-subpower}
\end{align}
If \(C_{\rm all}^*\le L^\delta\), \(0<\delta<1\), and the number of occurrences is
at most \(C\mathfrak L\), then
\begin{equation}\label{eq:appE-growing-subpower}
 \bigl(L^{C_{\rm all}^*}\bigr)^{C\mathfrak L}
 \le\exp\!\left(C\kappa L^\delta(\log L)^{3/2}\right)
 =\eps^{-o(1)}.
\end{equation}
\end{lemma}

\begin{proof}
Take logarithms.  The logarithms of the three factors on the left of
\eqref{eq:appE-combined-subpower} are respectively
\(a\sqrt L\), \(b\log L\), and
\(M\log\mathfrak L\), all \(o(L)=o(|\log\eps|)\).
For \eqref{eq:appE-growing-subpower}, the logarithm of the left side is at
most \(C\mathfrak L L^\delta\log L\), which is
\(O(L^\delta(\log L)^{3/2})=o(L)\).
\end{proof}

The packet scale is therefore
\begin{equation}\label{eq:appE-packet-scale}
 \sigma_K:=\cL_\eps^C\eps^{2(d-1)}\epstar^{-2(d-1)}
 =\eps^{2(d-1)-o(1)}.
\end{equation}
The proof uses a much weaker fixed consequence:
\begin{equation}\label{eq:appE-packet-fixed-credit}
 \sigma_K\le\eps^{2(d-1)-1/(100d)}.
\end{equation}
Indeed, \eqref{eq:appE-packet-fixed-credit} follows immediately from
Lemma~\ref{lem:appE-subpower-ledger}.

\subsection{Time powers in the two consumers}

In the attached branch $r>0$, hence $\iota=0$ and
$\widetilde{\mathsf X}=\mathsf X+m$.  After the physical $\tau$-volumes are
multiplied by the collision rate and Boltzmann-envelope factors, the
quantitative estimate \eqref{eq:effective-layer-parameter} makes the packet
ledger $\vartheta^{\mathsf X+m}$; the second source interpoland asks for
$\vartheta^{(\mathsf X+m)/10}$.

\begin{lemma}[Attached time comparison]
\label{lem:appE-attached-time}
Choose $c_0$ after $C_\beta$ and the fixed combinatorial constants so that
\begin{equation}\label{eq:appE-time-smallness}
 \vartheta\le q_0<1.
\end{equation}
Then, for every integer \(q\ge0\),
\begin{equation}\label{eq:appE-attached-time-comparison}
 \vartheta^q\le\vartheta^{q/10}.
\end{equation}
\end{lemma}

\begin{proof}
This is immediate from $0<\vartheta\le1$.
\end{proof}

This is the exact condition used in Proposition~\ref{prop:v2-consumer}.
Writing only $C\tau$ here would omit the growing $\alpha$ and $A$ factors
identified in Source A, Part~7.

In the disjoint branch the top packet has \(m\ge2\).  If \(r=0\), then
\(\iota=1\), \(|G|=1\), and its time factor is
\(\vartheta^{m-1}\).  Since $0<\vartheta\le1$,
\begin{equation}\label{eq:appE-disjoint-time}
 \vartheta^{m-1}\le\vartheta^{m/9}.
\end{equation}
For \(m\ge2\), this follows because
\(m-1-m/9=(8m-9)/9\ge7/9\).  If \(r>0\), then
\(\iota=0\) and the exponent \(m\) is automatically stronger than \(m/9\).

\subsection{Epsilon powers in the attached branch}

The attached packet ledger is
\begin{multline}\label{eq:appE-attached-ledger}
 \vartheta^{\widetilde{\mathsf X}}\,
 \eps^{(d-1)(|\widetilde H|-1-N_4^2-N_{\rm ee})}\,
 \sigma_K\\
 {}\times
 L^{C^*(\widetilde\rho+N_{\rm ee}+N_{\neg N3})}
 \eps^{N_{\rm good}/(9d)+|H'\setminus Z|/4}.
\end{multline}
The integer inequality
\[
 N_4^2+N_{\rm ee}\le|\widetilde H|-1
\]
makes the first epsilon exponent nonnegative.  The unchanged source
incidence bound gives
\[
 \widetilde\rho+N_{\rm ee}+N_{\neg N3}
 \le C^*(\widetilde\rho+N_4^2+N_{\rm good}
                    +|H'\setminus Z|).
\]
The good and initial-link terms absorb the corresponding logarithms using
their displayed fixed powers.  Since
\(N_4^2\le|\widetilde H|-1\le C^*\widetilde\rho\), the only remaining
logarithmic power is \(L^{C^*\widetilde\rho}\).

The source relation
\begin{equation}\label{eq:appE-rho-window}
 \rho+r\le\widetilde\rho\le\rho+r+3\Gamma
\end{equation}
shows that the bounded \(3\Gamma\) difference is a finite subpower loss,
whereas the \((\rho+r)\)-dependent part is exactly the one permitted in the
second interpoland.  Combining
\eqref{eq:appE-packet-fixed-credit},
\eqref{eq:appE-attached-time-comparison} and
\eqref{eq:appE-rho-window} gives
\[
 \mathcal I\mathcal N
 \le
 \vartheta^{(\mathsf X+m)/10}
 \eps^{d-1+1/(18d)}
 L^{C^*(\rho+r)}.
\]
Indeed,
$2(d-1)-1/(100d)>d-1+1/(18d)$ for $d\ge4$, so this has the epsilon
strength required by
\eqref{eq:source-top-second-interpoland}.  Equality in this arithmetic is
not used: all omitted packet logarithms were absorbed using a strict
fraction of the much larger margin in \eqref{eq:appE-packet-scale}.

\subsection{Epsilon powers in the disjoint branch}

The rooted top packet estimate is
\begin{equation}\label{eq:appE-disjoint-ledger}
 \vartheta^{m-\iota}\eps^{(d-1)(|G|-1)}\sigma_K,
 \qquad |G|=r+\iota.
\end{equation}
If \(r=0\), then \(\iota=1\), the baseline exponent is zero, and
\eqref{eq:appE-disjoint-ledger} has epsilon strength
\(2(d-1)-1/(100d)\).  If
\(r>0\), then \(\iota=0\) and its strength is
\((d-1)(r+1)-1/(100d)\).  The source top target is
\begin{equation}\label{eq:appE-disjoint-target-power}
 d-1+\frac1{18d}+\frac{4r}{C_{14}^*}.
\end{equation}
Choose \(C_{14}^*>4\).  No bound of \(r\) by \(\Gamma\) is available or
needed: the source range is \(r+m\le2\Lambda_\ell\).  For \(r=0\),
the difference is $d-1-59/(900d)>0$, so the packet dominates
\eqref{eq:appE-disjoint-target-power}.  For \(r\ge1\),
\[
 r\left(d-1-\frac4{C_{14}^*}\right)-\frac{59}{900d}>0.
\]
The right-hand side is bounded below by the positive \(r=1\) value and grows
linearly with \(r\).  Once \(C_{14}^*\) is fixed, every per-root logarithmic
factor is \(\eps^{-o(1)r}\) and is absorbed by a fixed fraction of this
margin, uniformly for the full range \(r+m\le2\Lambda_\ell\).  Thus
\eqref{eq:source-disjoint-top} follows uniformly in both cases.

\subsection{Closed form of the cutoff hierarchy}

Set \(A_{\mathfrak L}=L_\eps^\#\) and
\[
 \Lambda_\ell=A_\ell^{10d},
 \qquad A_{\ell-1}=\Lambda_\ell^{10d}.
\]
Writing \(A_\ell=(L_\eps^\#)^{a_\ell}\) gives
\begin{equation}\label{eq:appE-hierarchy-exponents}
 a_{\mathfrak L}=1,
 \qquad
 a_{\ell-1}=100d^2a_\ell.
\end{equation}
Consequently
\begin{equation}\label{eq:appE-hierarchy-closed}
 A_\ell=(L_\eps^\#)^{(100d^2)^{\mathfrak L-\ell}},
 \qquad
 \Lambda_\ell=(L_\eps^\#)^{10d(100d^2)^{\mathfrak L-\ell}}.
\end{equation}
In particular,
\[
 \Lambda_{\ell-1}=\Lambda_\ell^{100d^2}.
\]
The loss exponent in one restart step is
\begin{equation}\label{eq:appE-restart-loss-scale}
 Y_\ell:=\Lambda_\ell^2A_\ell
 =\Lambda_\ell^{2+1/(10d)}.
\end{equation}
Thus
\begin{equation}\label{eq:appE-Y-ratio}
 \frac{Y_{\ell+1}}{Y_\ell}
 =\Lambda_{\ell+1}^{-(100d^2-1)(2+1/(10d))}.
\end{equation}
Since the smallest \(\Lambda_{\ell+1}\) is \((L_\eps^\#)^{10d}\), the ratio in
\eqref{eq:appE-Y-ratio} is below \(1/2\) uniformly in \(\ell\) after
reducing \(\varepsilon_0\).  Hence
\begin{equation}\label{eq:appE-Y-tail}
 \sum_{k>\ell}Y_k\le2Y_{\ell+1}.
\end{equation}

Let \(a_0=-\log q_0>0\), where \(\vartheta\le q_0<1\).  The negative logarithm of
the layer-\(\ell\) restart factor is at least
\(a_0\Lambda_\ell/10\).  The logarithm of all later recurrence losses is,
by \eqref{eq:appE-Y-tail}, at most
\[
 2(d-1)L\sum_{k>\ell}Y_k
 \le4(d-1)L\Lambda_{\ell+1}^{2+1/(10d)}.
\]
Using \(\Lambda_\ell=\Lambda_{\ell+1}^{100d^2}\), their ratio is
\begin{equation}\label{eq:appE-restart-ratio}
 \frac{4(d-1)L\Lambda_{\ell+1}^{2+1/(10d)}}
      {a_0\Lambda_\ell/10}
 =
 \frac{40(d-1)L}{a_0}
 \Lambda_{\ell+1}^{-(100d^2-2-1/(10d))}
 =o(1)
\end{equation}
uniformly in \(\ell\).  This is the quantitative reason one restart gain
dominates every later loss; it is not an informal comparison of names in the
hierarchy.

\subsection{Per-root losses versus finite losses}

Suppose \(C_{\rm all}^*\le L^\delta\) with \(\delta<1/10\).  A factor that occurs
only \(C\mathfrak L\) times is handled by
\eqref{eq:appE-growing-subpower}.  A factor repeated for every cumulant root
must instead be combined with the root credit:
\[
 q_{\rm eff}
 =\eps^{1/C_{14}^*}L^{C_{\rm all}^*}.
\]
Since \(C_{14}^*\le L^\delta\),
\begin{equation}\label{eq:appE-qeff-log}
 \log q_{\rm eff}
 \le-L^{1-\delta}+L^\delta\log L
 \le-\frac12L^{1-\delta}.
\end{equation}
For \(s\le L\), this gives
\begin{equation}\label{eq:appE-qeff-subset}
 (1+q_{\rm eff})^s-1
 \le2s q_{\rm eff}
 \le2L\exp\!\left(-\frac12L^{1-\delta}\right).
\end{equation}
The last quantity is smaller than every fixed negative power of \(L\), and
in particular absorbs all label and subset counts that are repeated per
root.  This calculation is why a per-root loss is never placed into the
undifferentiated symbol \(\eps^{-o(1)}\).

\subsection{Final fixed-exponent margin}

For the full epoch-stacked range, the final decomposition
\eqref{eq:epoch-final-cumulant} has the following three families.  The older
one-epoch quantities in \eqref{eq:global-EH}--\eqref{eq:global-approx}
remain valid regression estimates under the one-epoch hypothesis, but they
are not the final ledger used here.
\begin{center}
\small
\setlength{\tabcolsep}{3pt}
\begin{tabular}{@{}p{0.36\textwidth}p{0.21\textwidth}p{0.29\textwidth}@{}}
\toprule
family & available exponent & distance above \(1/(400d)\)\\
\midrule
rebased nonempty cumulants \(\mathcal K_H^t\)
 & \(7/(1000d)\) & \(9/(2000d)\)\\
structured dynamical and finite-correlation errors
 & \(a_{\rm err}\)
 & \(\ge3/(80d)\)\\
excluded volume
 & \(d\) & \(d-1/(400d)\)\\
\bottomrule
\end{tabular}
\end{center}
Here \(a_{\rm err}=\min\{1/(20d),1/100\}\), as in
\eqref{eq:epoch-error-bound}.
The smallest margin is \(9/(2000d)\), from the rebased nonempty cumulants.
Choose once and for all \(\eta<9/(4000d)\) in
Lemma~\ref{lem:appE-subpower-ledger}.  Every bounded subpower loss then
uses at most half of that smallest margin.  For the full-range main cumulants,
\eqref{eq:epoch-q-log} and \eqref{eq:epoch-cumulant-subset-sum} handle the
per-root factors.  The localized-defect sum retains its singleton and is
controlled by the explicit \(2L\) reserve
in the proof of Proposition~\ref{prop:epoch-stacked-endpoint}.  Polynomial
factors \(s\), \(s^2\), \(A\) and \(\mathfrak L\), as well as the cubic
Gaussian clock \(O((\log L)^3)\), are subpower under the main logarithmic
hypothesis.  Therefore each of the three families is bounded by
\(o(1)\eps^{1/(400d)}\), uniformly in the theorem's full parameter range.
There is no separate approximate-tensor family in this ledger because
\eqref{eq:epoch-final-cumulant} is already rebased over \(n(t)\).

For the hard-core indicator, the pair union bound gives
\[
 \int n^{\otimes s}
 \1_{\{\exists i<j:\ d_{\T^d}(x_i,x_j)<\eps\}}
 \le \binom{s}{2}\|\rho(t)\|_\infty |B_d(1)|\eps^d.
\]
The propagated weighted Boltzmann bound gives
\(\|\rho(t)\|_\infty\le C_\beta A\), so this is the excluded-volume entry in
the table.  These calculations complete the arithmetic closure of the main
theorem.
\clearpage
\section*{Online Resource 1: Supplementary Proofs and Interfaces}
The following sections reproduce Online Resource~1, collecting the detailed
proofs, complementary routes, source concordance, and compatibility tables in
one continuous document.
\renewcommand{\thesection}{\arabic{section}}
\makeatletter\let\sectionname\@empty\makeatother
\renewcommand{\theHsection}{OR.\arabic{section}}
\renewcommand{\theHsubsection}{\theHsection.\arabic{subsection}}
\renewcommand{\theHsubsubsection}{\theHsubsection.\arabic{subsubsection}}
\renewcommand{\theHequation}{\theHsection.\arabic{equation}}
\setcounter{section}{0}
\setcounter{subsection}{0}
\setcounter{subsubsection}{0}
\section{Detailed Coordinate Verifications}
\label{supp:fixed-periodic-proofs}

These complete proofs correspond to the statements in Main Sections~\ref{sec:ensemble-flow}--\ref{sec:parallel-ov} and are presented in source order.

\suppproofheading{Lemma~\ref{lem:fixed-time-33A-j2-j3}}{supp:proof:lem:fixed-time-33A-j2-j3}
\begin{proof}
\emph{Step 1: serial three-output conditional measure.}

We first treat $j=3$.  Write $B=b_m(t_2)$ and
$\Delta=t_1-t_2$.  Disintegrating
\eqref{eq:j3-transformed-density} with respect to $(t_1,t_2)$ gives the
canonical fixed-time measure
\begin{equation}\label{eq:j3-fixed-time-density}
 \dd\nu_{3}^{t_1,t_2}
 =c_d\eps^{d-1}
   \frac{|p\cdot\omega_2|}{|y|^{d-2}|\Delta|^d}
   \dd\omega_2\,
   \dd\mathcal H^{d-1}_{y^\perp}(w),
 \qquad
 y=\frac{B-\eps\omega_2}{\Delta}.
\end{equation}
At the single point (if any) where $B=\eps\omega_2$, set the displayed
density equal to zero and choose the direction of $y$ arbitrarily.  This
does not change the measure.  The resulting formula defines a Borel Radon
kernel for every fixed $(t_1,t_2)$ with $|\Delta|\ge\mu$; by Tonelli and
the area formula it agrees with the ordinary Radon disintegration for almost
every pair.  We use this formula as the selected representative at the
exceptional pairs as well.  On the actual velocity cutoff,
$|p|+|w|\le\cL_\eps^C$.  Thus it remains to put the integrable singularity
 $|B-\eps\omega_2|^{-(d-2)}$ into a boundary-independent angular box.
On every enlarged chart below the pulled-back kernel, including its cutoff
indicator, is extended by zero outside the actual cell.

We give the charts explicitly.  Make the following half-open partition in
the fixed parameter $B$: either
\[
 \big||B|-\eps\big|\ge\eps/4,
 \qquad\hbox{or}\qquad
 \big||B|-\eps\big|<\eps/4,
\]
assigning the equality to the first cell.  In the first cell
$|B-\eps\omega_2|\ge\eps/4$, so the standard finite half-open atlas of
 $\Sph^{d-1}$ gives a fixed box and a density at most $C_d\eps^{-(d-2)}$.

In the second cell put $e=B/|B|$ and divide $\Sph^{d-1}$ into the half-open cap
$\omega_2\cdot e>1/2$ and its complement.  The complement is as above.  On
the cap use
\begin{equation}\label{eq:j3-fixed-cap-chart}
 \omega_2=\sqrt{1-r^2}\,e+r\sigma,
 \qquad
 0\le r<\sqrt3/2,
 \qquad \sigma\in\Sph^{d-2}\subset e^\perp.
\end{equation}
Choose the rotation sending a fixed coordinate vector to $e$ by the first
available half-open coordinate chart; Gram--Schmidt then identifies the
$\sigma$-sphere with one of finitely many fixed $\Sph^{d-2}$ patches.  This is a
Borel finite atlas with unit rotational Jacobian.  In these coordinates
\begin{align}
 \dd\omega_2
 &=\frac{r^{d-2}}{\sqrt{1-r^2}}\dd r\dd\sigma,
 \label{eq:j3-fixed-cap-area}\\
 |B-\eps\omega_2|^2
 &= (|B|-\eps)^2
    +2|B|\eps\bigl(1-\sqrt{1-r^2}\bigr)
 \ge (|B|-\eps)^2+c\eps^2r^2.
 \label{eq:j3-fixed-cap-denominator}
\end{align}
Consequently
\begin{equation}\label{eq:j3-fixed-cap-density-bound}
 \frac{\dd\omega_2}{|B-\eps\omega_2|^{d-2}}
 \le C_d\eps^{-(d-2)}\dd r\dd\sigma.
\end{equation}
This is the fixed-box form of the $d$-dimensional spherical-potential
bound
\begin{equation}\label{eq:j3-fixed-sphere-potential}
 \sup_{B\in\R^d}
 \int_{\Sph^{d-1}}\frac{\dd\omega}{|B-\eps\omega|^{d-2}}
 \le C_d\eps^{-(d-2)}.
\end{equation}

For complete endpoint bookkeeping let
$K_\eps=\lceil|\log_2\eps|\rceil+2$ and refine the cap, without enlarging
it, into the disjoint half-open radial cells
\[
 2^{-k-1}<r\le2^{-k}\quad(0\le k<K_\eps),
 \qquad 0\le r\le2^{-K_\eps},
\]
intersected with $r<\sqrt3/2$.  These are merely indicators in the single
fixed box of \eqref{eq:j3-fixed-cap-chart}; hence they do not create a
factor $K_\eps$ in the operator norm.  The terminal cap contains
$B=\eps\omega_2$, equivalently $y=0$.  At that single radial face the
direction of $y$ may be assigned arbitrarily because the face has zero
$\dd r\dd\sigma$ measure.  Away from it, cover $y/|y|\in\Sph^{d-1}$ by the
usual finite half-open atlas and write
$w=\sum_{i=1}^{d-1}\xi_iE_i(y/|y|)$ with
$|\xi_i|\le\cL_\eps^C$.  This puts the moving hyperplane measure in fixed
$\xi$-boxes with unit Gram density.  The set $p\cdot\omega_2=0$ is kept in
the zero-flux half-open cell; no inverse Jacobian is used there.

Since $|y|^{-(d-2)}=|\Delta|^{d-2}|B-\eps\omega_2|^{-(d-2)}$,
\eqref{eq:j3-fixed-time-density} and
\eqref{eq:j3-fixed-cap-density-bound} now give, chart by chart,
\[
 \dd\nu_{3}^{t_1,t_2}
 \le \cL_\eps^C\eps|\Delta|^{-2}
       \dd r\dd\sigma\dd\xi
 \le \cL_\eps^C\eps\mu^{-2}
       \dd r\dd\sigma\dd\xi.
\]
Composition of $Q$ with the displayed recovery maps proves
\eqref{eq:j3-fixed-time-sphere} for arbitrary $Q$.

\emph{Step 2: four-output latitude estimate.}

We next prove the fixed-time $j=4$ estimate on the full range, without using
the ambient small-$w$ volume in \eqref{eq:j4-direct-ball}.  Put
$q_{\rm sgn}=\eps/\Delta$, $q=|q_{\rm sgn}|$, retain the Carleman variables
$y=\rho e$, $e\in\Sph^{d-1}$, and use $w_*$ from
\eqref{eq:j4-wstar}.  Disintegrating the exact joint measure with respect to
$(t_1,t_2)$ and substituting
$w=w_*+q_{\rm sgn}e$ in the delta form
$|y|^{-(d-3)}\delta(y\cdot w)\dd w\dd y$ gives, up to the bounded residual
chart density, the conditional measure
\begin{multline}\label{eq:j4-fixed-time-latitude-density}
 \dd\nu_4^{t_1,t_2}
 \le C_d\eps^{d-1}|\Delta|^{-d}
 |(v_4-v_7)\cdot\omega_2|\,
 \rho\,\dd\rho\,\dd\omega_2\,\dd\xi\\
 {}\times
 \delta(e\cdot w_*+q_{\rm sgn})\,\dd e .
\end{multline}
Indeed, $\delta(\rho(e\cdot w_*+q_{\rm sgn}))$ supplies $\rho^{-1}$,
which leaves the radial density $\rho\,\dd\rho$.  This formula is first an
identity away from $\rho=0$; assigning an arbitrary direction at $\rho=0$
gives its selected Borel version without changing the measure.

Let $R=|w_*|$.  If $R<q$, the latitude is empty.  If $R\ge q$, spherical
coarea gives, for every nonnegative Borel $F$,
\begin{multline}\label{eq:j4-fixed-time-latitude-coarea}
 \int_{\Sph^{d-1}}\delta(e\cdot w_*+q_{\rm sgn})F(e)\dd e\\
 =\frac1R\left(1-\frac{q^2}{R^2}\right)^{(d-3)/2}
 \int_{\Sph^{d-2}}
 F\left(-\frac{q_{\rm sgn}}R\widehat w_*
 +\sqrt{1-\frac{q^2}{R^2}}\,\sigma\right)\dd\sigma .
\end{multline}
For $d\ge4$ the displayed density is bounded by
$R^{-1}\le q^{-1}=|\Delta|/\eps$, including the limiting latitude $R=q$.
The collision numerator and the $\rho$-range are polylogarithmically
bounded.  A finite half-open atlas for $\omega_2$, $\widehat w_*$ and
$\sigma$, together with the fixed radial interval for $\rho$, gives a
boundary-independent product box.  Compose $Q$ with the latitude recovery
map and extend it by zero outside the original cell.  Equations
\eqref{eq:j4-fixed-time-latitude-density}--\eqref{eq:j4-fixed-time-latitude-coarea}
then yield
\[
 \dd\nu_4^{t_1,t_2}
 \le \cL_\eps^C\eps^{d-2}|\Delta|^{-(d-1)}
       \dd\rho\,\dd\omega_2\,\dd\xi\,\dd\sigma,
\]
which proves \eqref{eq:j4-fixed-time-sphere} on every small- and large-$w$
cell for arbitrary $Q$.

\emph{Step 3: two-output recovery map.}

For $j=2$, use the second-contact coarea before making the direct/transformed
choice used in \eqref{eq:j2-direct-ball}.  The area formula behind
\eqref{eq:safe-output-second-coarea} is applied on each fixed lift and sign
cell before either contact indicator is enlarged.  On such a cell the
recovered $d$-dimensional output $u=v_2-v_1$ is a smooth
fractional-linear function of $(t_2,\omega_2)$ (it is not affine in $t_2$).
More explicitly, after the first-contact coordinates have been fixed and
absorbed into $\xi$, the second-contact equation has the form
\begin{equation}\label{eq:j2-fixed-time-recovery-map}
 u=\Psi_{m,\xi}(t_2,\omega_2)
   =\frac{B_{m,\xi}(t_2)-\eps\omega_2}{\Delta},
 \qquad B'_{m,\xi}(t_2)=v_1-v_7.
\end{equation}
Since $\Delta'= -1$ and $B_{m,\xi}-\eps\omega_2=\Delta u$,
\begin{equation}\label{eq:j2-fixed-time-recovery-derivatives}
 \partial_{t_2}\Psi_{m,\xi}=\frac{v_2-v_7}{\Delta},
 \qquad
 D_{\omega_2}\Psi_{m,\xi}[E_k]
 =-\frac{\eps}{\Delta}E_k\quad(1\le k\le d-1).
\end{equation}
Consequently the absolute $d$-dimensional Jacobian of the time column and
an oriented orthonormal tangent frame $(E_1,\ldots,E_{d-1})$ is
\begin{equation}\label{eq:j2-fixed-time-recovery-jacobian}
 \left|\det D_{(t_2,\omega_2)}\Psi_{m,\xi}\right|
 =\eps^{d-1}\frac{|(v_2-v_7)\cdot\omega_2|}{|\Delta|^d},
\end{equation}
which is exactly \eqref{eq:safe-output-second-coarea}.  The fixed lift and
incoming-sign refinement separates the at most two line--sphere roots, so
the map is one-to-one on each of the resulting half-open cells.

To record all factors rather than hide the first collision weight in the
notation $\dd\xi$, choose the finite first-contact angular and moving-frame
atlas used above and denote its remaining coordinates by $\xi$.  After
zero extension to its fixed cutoff box, its Gram factors and the first
collision flux have a Borel density
\begin{equation}\label{eq:j2-fixed-time-residual-density}
 0\le a_2(t_1,t_2,\omega_2,\xi)\le\cL_\eps^C.
\end{equation}
This density is independent of the test kernel.  Thus the exact joint
measure identity says that, for every nonnegative Borel test function $F$,
\begin{multline}\label{eq:j2-fixed-time-kernel-identity}
 \int F\,\dd\nu_2
 =\int\!\int\eps^{d-1}
   \frac{|(v_2-v_7)\cdot\omega_2|}{|t_1-t_2|^d}
   a_2(t_1,t_2,\omega_2,\xi)
   F(t_1,t_2,\omega_2,\xi)\\
 \dd\omega_2\dd\xi\dd t_2\dd t_1.
\end{multline}
In particular, disintegration with respect to the two atom times has the
selected Borel version
\begin{equation}\label{eq:j2-fixed-time-density}
 \dd\nu_{2}^{t_1,t_2}
 =\eps^{d-1}
   \frac{|(v_2-v_7)\cdot\omega_2|}{|\Delta|^d}
   a_2(t_1,t_2,\omega_2,\xi)\dd\omega_2\dd\xi.
\end{equation}
Equations \eqref{eq:j2-fixed-time-kernel-identity} and
\eqref{eq:j2-fixed-time-density}, rather than a formal division by a root
Jacobian, define the kernel also on rank-drop parameter values.  They agree
with the ordinary conditional measure almost everywhere and supply the
selected zero-extended version at every fixed pair of times.  In particular,
the selected version is fixed before $Q$.

Notice in particular that this argument does not drop the second-contact
indicator on the small-$u$ cell: doing so would control a $d$-dimensional
ambient $u$-volume and would not prove a bound for its fixed-$t_2$
$(d-1)$-dimensional conditional slice.

To make the endpoints explicit, put $r_0=4\eps/\mu$ and partition the
recovered $u$-output into the half-open terminal ball $|u|\le r_0$ and the
dyadic shells
\[
 2^{k-1}r_0<|u|\le2^kr_0
\]
up to the velocity cutoff, assigning the final boundary to the last cell.
All these sets are retained as disjoint indicators in the same
$(\omega_2,\xi)$ box; in particular the small-velocity point $u=0$ is not
discarded and no shell is enlarged separately.  On the cutoff support the
numerator in \eqref{eq:j2-fixed-time-density} is at most $\cL_\eps^C$.
The grazing set where it vanishes is a zero-flux cell, not an inverse pivot.
Thus
\[
 \dd\nu_{2}^{t_1,t_2}
 \le\cL_\eps^C\eps^{d-1}\mu^{-d}
      \dd\omega_2\dd\xi,
\]
which proves \eqref{eq:j2-fixed-time-sphere}.  The lift and sign partitions
are finite and half-open, so summing them costs only $\cL_\eps^C$.
\end{proof}

\section{Complementary Coarse and Conditional Packet Routes}
\label{supp:packet-proofs}

The sharp radial-conjugate proof and every load-bearing closure are contained in the main manuscript.  This section records independent coarse full-chord, root-exchange, and conditional packet routes that are not substituted for the sharp main input.

\suppproofheading{Lemma~\ref{lem:two-sided-bridge}}{supp:proof:lem:two-sided-bridge}
\begin{proof}
There is no collision on the selected physical line between its last upper
atom $u_i$ and its first lower atom $d_i$.  Its endpoint positions are
therefore $X_i-a_iV_i$ and $X_i+b_iV_i$.  After subtracting the frozen
companion positions and fixed lifts, the derivative block is
\[
 \begin{pmatrix}I_d&-a_iI_d\\ I_d&b_iI_d\end{pmatrix}.
\]
Subtracting the first block row from the second gives
\eqref{eq:two-sided-bridge-determinant}.  The two separator states are
distinct because the matching uses two different physical lines, so the
conditional product follows.
\end{proof}

\suppproofheading{Proposition~\ref{prop:aligned-three-line-bridge}}{supp:proof:prop:aligned-three-line-bridge}
\begin{proof}
Translate the separator to time zero and put
\[
 A=s-t_{u_1},\quad a=s-t_{u_2},\quad
 b=t_{d_1}-s,\quad B=t_{d_2}-s.
\]
Thus $A>a\ge h$ and $B>b\ge h$.  Let $n$ be the normal at the later upper
contact $(c,2)$ and let $m$ be the normal at the earlier lower contact
$(0,c)$.  Write $P_n=n\tensor n$ and $P_m=m\tensor m$.  The other two
normals do not enter the position derivative before their own scattering.

Use the contact-residual orientations from the outer line toward the centre
line.  Backward reconstruction through the later upper reflection and
forward reconstruction through the earlier lower reflection give, up to
terms depending only on the retained centre state,
\begin{align*}
 U_1&=x-Au+(A-a)P_nw,&
 U_2&=y-aw,\\
 R_1&=x+bu,&
 R_2&=y+Bw-(B-b)P_m u.
\end{align*}
These formulas are the exact equal-mass reflection rule: a variation of the
outer velocity transfers to the centre velocity by the rank-one projection
onto the contact normal.

Subtract the $R_1$ row from the $U_1$ row and the $U_2$ row from the $R_2$
row, and then pivot the $x,y$ columns.  With
\[
 L_1=A+b,\qquad L_2=B+a,
 \qquad r=\frac{(A-a)(B-b)}{L_1L_2},
\]
the remaining $2d$-dimensional block on $(u,w)$ is, up to harmless row
signs,
\[
 \begin{pmatrix}
  L_1I_d&-(A-a)P_n\\
  -(B-b)P_m&L_2I_d
 \end{pmatrix}.
\]
The matrix determinant lemma and the fact that $P_mP_n$ has sole possible
nonzero eigenvalue $(m\cdot n)^2$ yield the exact identity
\begin{align}
 \left|\det D_{(x,u,y,w)}(U_1,U_2,R_1,R_2)\right|
 &= (L_1L_2)^d\left(1-r(m\cdot n)^2\right)\notag\\
 &\ge (L_1L_2)^d(1-r)\notag\\
 &= (L_1L_2)^{d-1}(A+B)(a+b).
 \label{eq:aligned-three-line-exact}
\end{align}
Here we used
\[
 L_1L_2-(A-a)(B-b)=(A+B)(a+b).
\]
Every factor $L_1,L_2,A+B,a+b$ is at least $2h$, which proves
\eqref{eq:aligned-three-line-lower}.  In particular no angular separation of
$m$ and $n$ is required.

Apply the $d$-dimensional sphere identity at all four contacts and coarea
in $(x,u,y,w)$.  The four sphere factors give $\eps^{4(d-1)}$ and
\eqref{eq:aligned-three-line-lower} costs at most $C_d h^{-2d}$.  Elastic and
transport substitutions have unit Jacobian.  The four collision fluxes are
not denominators in this chart and are bounded by $\cL_\eps^C$ on the
velocity cell.  The three-particle full normalization is
$\eps^{-3(d-1)}$.
Keeping the centre state and all four times and normals inside the same
nonnegative integrand gives \eqref{eq:aligned-three-line-operator} for
arbitrary $Q$.  Zero extension of the two solved outer states gives a fixed
output box, independent of the retained centre values.
\end{proof}

\suppproofheading{Proposition~\ref{prop:all-three-line-bridge}}{supp:proof:prop:all-three-line-bridge}
\begin{proof}
Relabel the particles so that the upper contacts are $(0,1)$ at time $-A$
and $(1,2)$ at time $-a$, where $A>a\ge h$.  Write the lower times as
$b<B$, with $b\ge h$, and put
\[
 p=A-a,\qquad q=B-b,\qquad
 L_1=A+b,\qquad L_2=a+B.
\]
Let $N=n\tensor n$ be the normal projection at the later upper contact and
$M=m\tensor m$ the normal projection at the earlier lower contact.  The
first-lower rule leaves six possible ordered pairs of distinct lower edges.
For the first four, the indicated matching and elimination of its two
separator positions reduce the velocity block to
\begin{equation}\label{eq:three-line-projection-block}
 \begin{pmatrix}
  L_1I_d&-pE_n\\
  -qE_m&L_2I_d
 \end{pmatrix},
\end{equation}
with the following complete table:
\[
\begin{array}{c|c|c|c}
 (d_1,d_2)&\text{selected lines}&E_n&E_m\\ \hline
 ((0,1),(0,2))&(0,2)&N&I_d-M\\
 ((0,1),(1,2))&(0,2)&N&M\\
 ((0,2),(0,1))&(0,1)&I_d-N&I_d-M\\
 ((0,2),(1,2))&(0,1)&I_d-N&M
\end{array}
\]
this table follows directly by applying
\[
 (v_i,v_j)\longmapsto
 ((I_d-P)v_i+Pv_j,\ Pv_i+(I_d-P)v_j)
\]
at the one intervening upper or lower collision.  Thus every $E_n,E_m$ in
the table is an orthogonal projection.

The determinant of \eqref{eq:three-line-projection-block} is
\begin{equation}\label{eq:three-line-projection-determinant}
 (L_1L_2)^d\det(I_d-rE_mE_n),
 \qquad r=\frac{pq}{L_1L_2}.
\end{equation}
The nonzero eigenvalues of a product of two orthogonal projections coincide
with those of the positive contraction $E_nE_mE_n$ on
$\operatorname{Ran}E_n$; hence they belong to $[0,1]$.  Since $0\le r<1$,
\[
 \det(I_d-rE_mE_n)\ge(1-r)^d.
\]
Moreover
\[
 L_1L_2-pq=(A+B)(a+b)\ge4h^2.
\]
Therefore the modulus of every determinant in the table is at least
$(4h^2)^d=2^{2d}h^{2d}$.

The two remaining lower orders begin with $(1,2)$ and end with either
$(0,1)$ or $(0,2)$.  Select line $0$ and, respectively, line $1$ or line
$2$.  Subtracting the upper and lower residuals on the repeated $(1,2)$
edge first makes the block triangular; its determinant is exactly
\[
 (A+B)^d(a+b)^d\ge2^{2d}h^{2d}.
\]
This exhausts the ordered lower paths.  The sphere, normalization,
arbitrary-$Q$ and fixed-output-domain calculation is identical to the last
paragraph of Proposition~\ref{prop:aligned-three-line-bridge}.
\end{proof}

\suppproofheading{Proposition~\ref{prop:three-line-upstream-rooted}}{supp:proof:prop:three-line-upstream-rooted}
\begin{proof}
Write $u_-<u_+$ for the two distinct last-upper atoms selected by the
matching, and let $p$ be the selected physical line whose last upper atom is
$u_-$.  The collision edge at $u_+$ cannot contain $p$, because that would
make $u_+$ a later upper atom on the same line.  With only three particle
lines, the edge at $u_+$ is therefore the pair formed by the other two lines.
The two distinct edges at $u_-,u_+$ form the complete reverse-Kruskal
spanning tree.

Suppose that an unmarked upper contact $a$ occurred in $(u_-,u_+)$.  It
cannot contain $p$, again by the definition of $u_-$.  Hence it has the same
unordered particle pair as $u_+$.  If several such contacts occur, take the
latest one.  No upper contact lies between it and $u_+$ on either of those
two physical lines: a contact with $p$ would contradict the last-upper
property of $p$, and a contact without $p$ would be another contact of the
same pair.  The two atoms are therefore joined by both intervening particle
bonds and form a double bond, contrary to the hypothesis.  There is no
upper contact after $u_+$, so every unmarked upper contact lies strictly
before $u_-$.

Use independent right event states immediately after every contact and the
raw separator state as coordinates.  Order the columns as
\[
 (y_{<u_-},\,h_T,V_s,\,y_{>d_2}),
\]
where $y_{<u_-}$ are the incoming-root variables before $u_-$, $h_T,V_s$
are the two tree translations and separator velocities used by the bridge,
and $y_{>d_2}$ are later root variables.  Before any root is substituted,
the constraint derivative has the causal form
\[
 \begin{pmatrix}
 A_<&0&0&0\\
 *&D_{h_T}\mathcal U&D_{V_s}\mathcal U&0\\
 0&D_{h_T}\mathcal R_D&D_{V_s}\mathcal R_D&0\\
 *&*&*&A_>
 \end{pmatrix}.
\]
The zero in the first column of the selected landing rows is literal in
this chart: the residuals at $u_-,u_+$ and at the two lower landings use
only the event states from $u_-$ through $d_2$.  This is the precise
coordinate/elimination order behind the asserted absence.  Schur-eliminate
$A_<$ first and $A_>$ last.  Their off-diagonal stars cannot change the
middle bridge block.  The diagonal root blocks are
$\eps^{d-1}|g_a\cdot\omega_a|$, while the remaining block is exactly one of the
six three-line bridge matrices in
Proposition~\ref{prop:all-three-line-bridge}.  Thus its determinant is at
least $2^{2d}h^{2d}$.  Every upstream root is neutral against its own sphere and
collision flux by Lemma~\ref{lem:unrooted-neutral-completion}; contacts after
$d_2$ are causally absent from the selected residuals and are neutral in the
same way.  The sphere and normalization ledger is unchanged, and the
arbitrary-$Q$ zero-extension argument in
Proposition~\ref{prop:all-three-line-bridge} proves
\eqref{eq:three-line-upstream-rooted-operator}.
\end{proof}

\suppproofheading{Proposition~\ref{prop:all-tree-only-bridge}}{supp:proof:prop:all-tree-only-bridge}
\begin{proof}
Every physical line must occur at a lower atom, while the two lower atoms
contain only four particle slots.  Hence $P\le4$.  The case $P=3$ is
Proposition~\ref{prop:all-three-line-bridge}; the hypotheses exclude $P<3$.
It remains to treat $P=4$.

For four lines the two lower edges $f_1,f_2$ are disjoint and partition the
particle set.  At the separator, take the $2d$ landing conormals whose
configuration parts are the $d$ coordinate copies of
$b_{f_1}^*$ and $b_{f_2}^*$.  Before the upper constraints are eliminated,
backward free flight from the two landing times gives the configuration and
velocity parts
\[
 Z=(b_{f_1}^*,b_{f_2}^*),\qquad
 W=(\Delta_1b_{f_1}^*,\Delta_2b_{f_2}^*).
\]
Because the two edges are disjoint,
\begin{equation}\label{eq:disjoint-conormal-gram}
 Z^*W=\operatorname{diag}(2\Delta_1I_d,2\Delta_2I_d)>0,
 \qquad
 \sqrt{\det(W^*W)}=2^d(\Delta_1\Delta_2)^d.
\end{equation}
The earlier lower reflection acts on the disjoint pair $f_1$ and therefore
does not alter the $f_2$ frame.

Lift this combined frame chronologically across the upper tree while
annihilating the upper contact-position rows.  In raw tree-edge coordinates
this is precisely Gaussian elimination of those rows: at a free segment the
lift is \eqref{eq:backward-free-conormal}, and at a retained-normal collision
it is the orthogonal part of
\eqref{eq:backward-collision-conormal}.  If the normal is allowed to vary
before the sphere variable is retained, the additional term is the positive
operator $\mathcal H$ in that formula.  Thus the cone condition in
\eqref{eq:disjoint-conormal-gram} is preserved and the $2d$-dimensional
velocity volume cannot decrease.  After the last tree row is eliminated,
the velocity part of the lifted frame is exactly the adjoint of the Schur row
$\mathcal S$ in \eqref{eq:tree-schur-row}.  This is the identification made
explicit in Lemma~\ref{lem:tree-schur-conormal-lift}.  Passage from tree-edge
to orthonormal reduced particle coordinates costs only a constant depending
on $P=4$.  Lemmas~\ref{lem:conormal-frame-identities} and
\ref{lem:tree-schur-conormal-lift} therefore give
\begin{equation}\label{eq:disjoint-tree-row-volume}
 \Jac_{2d}(\mathcal S)
 \ge c(d)|\Delta_1\Delta_2|^d
 \ge c(d)h^{2d}.
\end{equation}

Cauchy--Binet selects $2d$ scalar separator-velocity coordinates with a
minor comparable to \eqref{eq:disjoint-tree-row-volume}.  Apply the sphere
identity at the three upper-tree contacts and the two lower contacts, and
coarea in the $3d$ relative separator-position coordinates and those $2d$
velocity coordinates.  The sphere factor is
$\eps^{5(d-1)}$, the inverse long-gap Jacobian is at most
$C_dh^{-2d}$, and the full four-particle normalization is
$\eps^{-4(d-1)}$.  No
collision time is a pivot in this chart, so all five fluxes remain harmless
numerators bounded by $\cL_\eps^C$.  This gives
$\cL_\eps^C\eps^{d-1}\epstar^{-2d}$ after $h=\epstar/4$ and proves
\eqref{eq:all-tree-only-operator}.  All substitutions occur in the same
nonnegative integrand; zero extension defines the displayed fixed output
$Q_{H,\kappa}^{\sharp,\mathrm{tree}}$.
\end{proof}

\suppproofheading{Proposition~\ref{prop:disjoint-rooted-extension}}{supp:proof:prop:disjoint-rooted-extension}
\begin{proof}
At the separator the combined landing frame is still exactly the disjoint
frame in \eqref{eq:disjoint-conormal-gram}: before $d_2$ the only lower
collision is $d_1$, and it acts on the particle pair disjoint from $f_2$.
Process the upper history backward in its complete chronological order.  At a
marked tree contact use its cluster-translation pivot; at an unmarked upper
contact use its $d$ time--normal pivots.  Lemma~\ref{lem:rooted-tree-schur-lift}
identifies the resulting single Schur complement with the chronological
backward adjoint-pullback free-flight/collision lift.  In particular it is not an
iteration of Lemma~\ref{lem:tree-schur-conormal-lift}, whose hypotheses exclude
unmarked upper contacts.  The positive block Gram matrix in
\eqref{eq:disjoint-conormal-gram} and conormal monotonicity give
\[
 \Jac_{2d}(\mathcal S)\ge c(d,P_0)|\Delta_1\Delta_2|^d.
\]

Contacts strictly after $d_2$ do not occur in either selected residual by
causality.  Lemma~\ref{lem:rooted-tree-schur-lift} also records the neutral
$\eps^{d-1}|g\cdot\omega|$ determinant of every upper root, and the same local
block processes the later lower roots.  Thus every unmarked sphere factor
and collision flux cancels its own inverse root Jacobian.  Only the upper
tree, the two selected lower spheres and the full normalization remain:
\[
 \eps^{-(d-1)P}\eps^{(d-1)(P-1)}\eps^{2(d-1)}=\eps^{d-1}.
\]
Cauchy--Binet chooses one of finitely many $2d$-velocity minors of
$\mathcal S$.  This is only a local submersion statement.  For fixed
retained variables, the complete pivot equations---that minor, all
marked-tree cluster translations, and all unmarked time--normal
variables---form the bounded-degree square hard-sphere system
\eqref{eq:appD-full-polynomial-system}.  Theorem
\ref{thm:appD-uniform-algebraic-multiplicity} enumerates all of its regular
solutions in at most $M(d,P_0)$ boundary-independent branch slots.  Applying
Lemma~\ref{lem:appD-finite-fibre-substitution} branchwise gives the inverse
factor $C(d,P_0)\epstar^{-2d}$.  All substitutions act on the same nonnegative
integrand.  Retaining the minor and algebraic-branch labels and extending
the transformed kernel by zero on the fixed velocity, time, normal and lift
boxes defines $Q_{H,\kappa}^{\sharp,\mathrm{disj}}$ and proves
\eqref{eq:disjoint-rooted-operator}.
\end{proof}

\suppproofheading{Proposition~\ref{prop:asymmetric-new-line-peeling}}{supp:proof:prop:asymmetric-new-line-peeling}
\begin{proof}
Let $u_c$ be the last upper contact on the new line $c$.  Since $c$ does not
occur at $d_1$ and $d_1,d_2$ are consecutive first lower atoms, the line $c$
has no collision in $(u_c,d_2)$.  Choose the reverse-chronological Kruskal
tree $T$ from Lemma~\ref{lem:upper-tree-landings}.  Because $u_c$ is the
last upper contact of the selected crossing line $c$, that lemma proves that
$u_c\in T$ and that every rejected upper edge is rooted through a path of
later retained contacts.  Remove $u_c$ from $T$.  This gives two tree
components.
Use a particle in the component opposite $c$ as the common position
reference, take the separator state $(X_c,V_c)$ as $2d$ pivot coordinates,
and use the remaining particle states as complementary coordinates.  The
reduced incidence matrix of $T\setminus\{u_c\}$ is then invertible: its two
components are pinned respectively at $c$ and at the reference particle.
Its determinant and inverse are bounded in terms of $P_0$.

Use the separator-side outgoing representation at $u_c$ and the incoming
representation at $d_2$.  Lemma~\ref{lem:causal-free-line-endpoints}, first
backward in the upper word and then forward in the lower word, proves in raw
separator coordinates that neither companion state depends on
$(X_c,V_c)$.  Therefore Lemma~\ref{lem:positive-free-chord-reduction}
applies and gives the exact first
pivot
\begin{equation}\label{eq:asymmetric-chord-determinant}
 \left|\det D_{(X_c,V_c)}(U_{u_c},R_{d_2})\right|
 =(t_{d_2}-t_{u_c})^d\ge\epstar^d.
\end{equation}
All other molecule variables, all remaining contact constraints and the
complete nonnegative kernel are retained during this substitution.

It remains to impose the $P-2$ vector contacts of
$T\setminus\{u_c\}$ and the first landing $d_1$.  At the separator the $d$
conormals of $R_{d_1}$ have the one-speed form
\[
 (Z_1,W_1)=(Z_1,bZ_1),
 \qquad b=t_{d_1}-s\ge h=\epstar/4,
\]
so $Z_1^*W_1=bZ_1^*Z_1\ge0$ and
$\Jac_d(W_1)\ge c_d b^d$.  The chord substitution just performed contributes
the nonnegative term \eqref{eq:free-chord-index-term}.  Every incoming root
contributes its nonnegative cylinder-curvature shunt, while fixed-normal
elastic transports are orthogonal.  Thus the remaining forest, including
roots before $u_c$ on line $c$ and an arbitrary dependence of the $d_2$
companion on $d_1$, is exactly a grounded positive network in
Lemma~\ref{lem:grounded-positive-network-lift}.  Its trace-completeness
hypothesis is not inferred from connectivity; it is verified for this
actual word by Lemma~\ref{lem:appC-overlap-trace-certificate}.  Apply
\eqref{eq:grounded-one-speed-volume} with $k=d$.  Since the normalized
incidence map $Z_1$ is an isometry up to a constant depending only on
$P_0$, we obtain
\begin{equation}\label{eq:asymmetric-single-landing-volume}
 \Jac_d(W_{1,\mathrm{Sch}})\ge c(d,P_0)b^d
 \ge c(d,P_0)h^d.
\end{equation}
This is a forest Schur complement, but no new argument is hidden: the
chronological cluster pivots are the edges of
$T\setminus\{u_c\}$, and the two missing common translations are exactly
the two endpoint grounds used in
\eqref{eq:free-chord-index-term}.

Apply the sphere identity at the $P-1$ upper-tree contacts and at
$d_1,d_2$.  The chord pivot
\eqref{eq:asymmetric-chord-determinant} and the forest/single-landing pivot
\eqref{eq:asymmetric-single-landing-volume} cost at most
$C(P_0)\epstar^{-2d}$.  Every unmarked contact is neutral through its own
incoming time--normal determinant and collision flux; later lower contacts
are causally absent from the selected rows.  The complete sphere and
normalization ledger is
\[
 \eps^{-(d-1)P}\,
 \underbrace{\eps^{2(d-1)}}_{u_c,d_2}\,
 \underbrace{\eps^{(d-1)(P-2)}\eps^{d-1}}_{
   T\setminus\{u_c\},\ d_1}
 =\eps^{d-1}.
\]
The chord pivot is globally affine and single-valued for fixed remaining
pivots.  The second $d$-row Schur submersion is not declared globally
injective.  Together with every remaining tree and incoming-root equation it
is the square bounded-degree system
\eqref{eq:appD-full-polynomial-system}; Theorem
\ref{thm:appD-uniform-algebraic-multiplicity} enumerates every regular
solution in at most $M(d,P_0)$ branch slots.  Thus the finite-fibre coarea
formula, not a local implicit-function argument, justifies both elimination
steps with only a $C(d,P_0)$ cost.  Both steps act on one nonnegative
integrand.  The untranslated
reference position, all unpivoted velocities, every retained time and
normal, all root witnesses, all complement variables and the finite branch
label remain outputs.  Zero extension to the fixed half-open position,
velocity, time, normal, lift and branch cells defines
$Q_{H,\kappa}^{\sharp,\mathrm{ch}}$ independently of the boundary values and
proves \eqref{eq:asymmetric-new-line-operator}.
\end{proof}

\suppproofheading{Proposition~\ref{prop:covariant-rooted-bridge}}{supp:proof:prop:covariant-rooted-bridge}
\begin{proof}
Keep the two selected upper-tree rows and their two chronological cluster-
translation pivots until the end, and eliminate all other marked-tree rows
and all unmarked root blocks first.  In the raw reverse-Kruskal cluster
coordinates the eliminated position block is triangular with identity
reduced-incidence diagonal.  The two remaining selected position coordinates
are the kinetic-orthonormal relative translations of the isolated selected
line against its later-edge cluster.  Thus their input and output port maps
are partial isometries; passage back to a fixed orthonormal reduced-particle
basis costs only $c(P_0)$ in the final exterior volume.

Translate the separator to zero.  Relabel the two selected upper contacts so
that their times are $-A<-a<0$, and write the lower times as $0<b<B$.  Hence
\[
 A>a\ge h,\qquad B>b\ge h,
 \qquad p=A-a,\quad q=B-b.
\]
Before the other positions are eliminated, the selected upper and lower
rows have direct coefficients $-AI_d,-aI_d,bI_d,BI_d$.  The hypothesis is
that their simultaneous Schur elimination preserves these diagonal blocks.
For the chronological pairing in which the earlier upper endpoint is matched
to the earlier lower endpoint, the assumed remaining $2d$-velocity block is
\begin{equation}\label{eq:covariant-rooted-bridge-block}
 M=\begin{pmatrix}
       (A+b)I_d&-P\\
       -Q&(a+B)I_d
    \end{pmatrix}.
\end{equation}
Here, by hypothesis, $P$ is the total upper companion compliance compression
between the two matched event ports and $Q$ is the lower companion response
between the consecutive lower ports.  For the crossed chronological pairing,
the criterion assumes that both cross terms lie in the same strict triangular
block and that the diagonal of $M$ is unchanged.

We verify the two norm bounds in the preceding block identity.  In the
covariant event-port
coordinates of Lemma~\ref{lem:rooted-compliance-factorization}, the row
subtraction just made inserts the balanced source--sink current
\[
 j_u=-f,\qquad j_v=O_\gamma f
\]
at the two upper companion ports.  The chronological path $\gamma$ between
them is the fixed-event path from the reverse-Kruskal bridge calculation; it
is contained in $[-A,-a]$, its orthogonal transport is $O_\gamma$, and its
total resistance is $p$.  The endpoint maps from the selected velocity to
the source and from the port displacement to the selected residual are the
partial isometries fixed in the first paragraph.  Operator
Cauchy--Schwarz for the compliance form and
\eqref{eq:rooted-compliance-path-bound} therefore give
\begin{equation}\label{eq:covariant-upper-port-bound}
 \|P\|\le p=A-a.
\end{equation}
This is an identity for the particular row-reduced bridge block, not an
assertion that the complete raw matrix
$-(D_h\mathcal U)^{-1}D_V\mathcal U$ is self-adjoint.  Incoming contacts
inside the path add nonnegative shunts and can only decrease the compliance;
contacts before $-A$ are absent by causality.

The first two lower contacts are consecutive.  Between them the only reset
in the selected history is the retained-normal orthogonal reflection at the
first landing.  Sending the same balanced trial current along $[b,B]$ gives
\begin{equation}\label{eq:covariant-lower-port-bound}
 \|Q\|\le q=B-b.
\end{equation}
There is no earlier lower root, while later lower contacts are causally absent
from both selected residuals.

Apply the block determinant identity to
\eqref{eq:covariant-rooted-bridge-block}.  Equations
\eqref{eq:covariant-upper-port-bound}--
\eqref{eq:covariant-lower-port-bound} imply
\begin{align*}
 |\det M|
 &\ge\big((A+b)(a+B)-(A-a)(B-b)\big)^d\\
 &=\big((A+B)(a+b)\big)^d\ge2^{2d}h^{2d}.
\end{align*}
The triangular pairing has determinant
$(A+B)^d(a+b)^d$ and obeys the same bound.  Restoring the already eliminated
tree-position block changes the determinant only by a factor bounded above
and below in terms of $P_0$.  Cauchy--Binet then gives
\eqref{eq:covariant-rooted-bridge-volume}.
\end{proof}

\suppproofheading{Corollary~\ref{cor:conditional-long-gap-rooted-packet}}{supp:proof:cor:conditional-long-gap-rooted-packet}
\begin{proof}
Use Proposition~\ref{prop:covariant-rooted-bridge} and choose a maximal
$2d$-coordinate separator-velocity minor by Cauchy--Binet, with the usual
lexicographic half-open tie rule.  Jointly solve that minor and all marked-
tree cluster translations.  At every unmarked contact, solve its incoming
time and normal in the same chronological raw chart.  Its determinant is
$\eps^{d-1}|g\cdot\omega|$ and is cancelled by its own collision sphere and
flux.  Elastic and transport maps have unit Jacobian.  The factors not
cancelled locally are therefore exactly
\[
 \eps^{-(d-1)P}\eps^{(d-1)(P-1)}\eps^{2(d-1)}=\eps^{d-1}.
\]
Equation \eqref{eq:covariant-rooted-bridge-volume} costs at most
$C(P_0)h^{-2d}$.  The common untranslated position, the unpivoted separator
velocities, all retained times and normals, every incoming-root witness and
lift, and all complement variables remain in the same nonnegative integrand.
Extend the recovered variables by zero to their fixed fundamental-domain,
velocity-cutoff and half-open cell boxes.  This defines
$Q_{H,\kappa}^{\sharp,\mathrm{port}}$ independently of the boundary values
and proves \eqref{eq:conditional-long-gap-rooted-packet}.
\end{proof}

\suppproofheading{Lemma~\ref{lem:unrooted-neutral-completion}}{supp:proof:lem:unrooted-neutral-completion}
\begin{proof}
Return temporarily to the raw edge-state coordinates of the molecule and
order all upper contacts chronologically.  At an unmarked contact, its $d$
vector rows pivot its own time and $d-1$ normal coordinates.  Pointwise
normalize the tangent columns by
$\delta\widehat\theta=G_\omega^{1/2}\delta\theta$, so that
$D_\theta\omega\,\delta\theta=E_\omega\delta\widehat\theta$; the chart
density $J_\omega=(\det G_\omega)^{1/2}$ is the determinant of the same
column change.  In these normalized columns the normal derivative is
$-\eps I_{d-1}$ on $\omega_a^\perp$, while the time derivative is $g_a$ and has normal
component $g_a\cdot\omega_a$; the local determinant therefore has modulus
$\eps^{d-1}|g_a\cdot\omega_a|$.  At a marked tree contact, the $d$ vector rows
pivot the relative translation of the newly merged chronological tree
cluster.  Its reduced incidence block has determinant one in tree-edge
coordinates and is bounded above and below in orthonormal cluster
coordinates by constants depending only on $P_0$.  Elastic pair maps and
transport translations have determinant one.  Causality makes the resulting
full matrix of the rows $(\mathcal U,C_A)$ against the tree translations and
$y_A$ block lower triangular in these raw coordinates.  Hence its determinant
is the tree-incidence determinant times
$\eps^{(d-1)|A|}\prod_a|g_a\cdot\omega_a|$.

Now eliminate the tree rows by the block determinant identity.  The
tree-incidence determinant cancels, and the remaining Schur block is exactly
$D_{y_A}C_A$ in the reduced coordinates used in the statement.  This proves
\eqref{eq:unrooted-causal-determinant}, including upper co-tree contacts that
occur before the last tree merge.

Apply the block determinant identity to the rows $(C_A,\mathcal R_D)$ and
columns $(y_A,V_I,t_{d_1},t_{d_2})$.  The lower-right Schur complement is
exactly the derivative obtained after solving $C_A=0$ for $y_A$ by the
implicit-function theorem.  Chronological solution first fixes the incoming
time root and then its contact normal, so this derivative is precisely the
total root-substituted matrix $[\mathcal S_I\ \mathcal T]$.  Multiplication
by \eqref{eq:unrooted-causal-determinant} proves
\eqref{eq:unrooted-schur-equivalence}.  Finally, the sphere identity supplies
$\eps^{d-1}$ and the collision kernel supplies $|g_a\cdot\omega_a|$ at every
$a\in A$; these factors cancel their counterparts in the inverse Jacobian.
\end{proof}

\suppproofheading{Lemma~\ref{lem:three-normal-root-exchange}}{supp:proof:lem:three-normal-root-exchange}
\begin{proof}
Order the rows and columns as
\[
 (C_{A\setminus\{a\}},\mathcal U;,C_a,\mathcal R_D),
 \qquad
 (y_{A\setminus\{a\}},h;,t_a,V_I,t_{d_1},t_{d_2}).
\]
In chronological raw edge-state coordinates, the first diagonal block is
lower triangular.  Its marked-tree diagonal is the reduced tree-incidence
matrix.  At every $b\ne a$, after the same pointwise tangent-column
normalization used above, its unmarked-contact diagonal is
\[
 D_{(t_b,\theta_b)}C_b=[g_b\ -\eps E_b],
 \qquad
 \left|\det[g_b\ -\eps E_b]\right|
 =\eps^{d-1}|g_b\cdot\omega_b|.
\]
It follows that the determinant of the first pivot block is the prefactor
in \eqref{eq:root-exchange-block-identity}.  The lower-right Schur
complement is precisely the total derivative of $\mathcal F_a$ displayed
there.  The block determinant identity proves the formula.  Since
$\theta_a$ is not among the pivot columns, none of the $d-1$ columns
$-\eps E_a$ has been used and no selected-normal factor is hidden in the
last determinant.
\end{proof}

\suppproofheading{Lemma~\ref{lem:partial-normal-root-exchange}}{supp:proof:lem:partial-normal-root-exchange}
\begin{proof}
Put
\[
 G_a=(D_{\theta_a}\omega_a)^*D_{\theta_a}\omega_a,\qquad
 L_a=G_a^{-1}(D_{\theta_a}\omega_a)^*.
\]
Before any Schur elimination, replace the $d$ rows of $C_a$ by its normal
component $\omega_a\cdot C_a$ and its dual tangent residual $L_aC_a$.
This pointwise row change has density $J_a^{-1}$, uniformly bounded on the
chosen chart; it is paired with the chart density
$J_a=(\det G_a)^{1/2}$ in $\dd\omega_a$.  On the constraint surface $C_a=0$
the derivative of the coefficient $L_a$ contributes nothing, and
\[
 D_{\theta_a}(L_aC_a)=-\eps L_aD_{\theta_a}\omega_a
 =-\eps I_{d-1}.
\]
Chronological causality therefore makes the columns
$\theta_a^{J^c}$ lower triangular: their first nonzero diagonal block is
\[
 D_{\theta_a^{J^c}}
  \big((L_aC_a)_j\big)_{j\in J^c}
 =-\eps I_{d-1-k}.
\]
Use these tangent rows to clear every later occurrence of the same columns.
This row operation has determinant one and extracts the factor
$\eps^{d-1-k}$.  The selected contact then leaves its normal component and
the $k$ tangent components indexed by $J$, hence $1+k$ rows.  Together with
the $2d$ landing rows this is the $(2d+1+k)$-row map
$\mathcal F_{a,J}$.  Its remaining pivot columns are one root time,
$2(d-1)+k$ separator-velocity coordinates and the two landing times, also
$2d+1+k$ columns.

All other unmarked contacts and the marked tree are still the chronological
triangular pivot block used in
\eqref{eq:root-exchange-block-identity}.  They contribute
\[
 c_T\eps^{(d-1)(|A|-1)}\prod_{b\ne a}|g_b\cdot\omega_b|.
\]
Multiplying this factor by the extracted $\eps^{d-1-k}$ and the last Schur
minor proves \eqref{eq:partial-root-exchange-block-identity}.

When $k=0$, only the normal component of $C_a$ remains.  Its $t_a$ derivative
is $g_a\cdot\omega_a$.  Eliminate that row and column.  The last $2d$ rows
are exactly the total root-substituted landing derivative
$[\mathcal S_I\ \mathcal T]$, which proves
\eqref{eq:partial-root-standard-endpoint}.  When $k=d-1$ no selected normal
column is pivoted, so the construction reduces to
Lemma~\ref{lem:three-normal-root-exchange}.
\end{proof}

\suppproofheading{Proposition~\ref{prop:common-normal-fct}}{supp:proof:prop:common-normal-fct}
\begin{proof}
In each of the $d-1$ directions in $n^\perp$, every fixed-normal reflection
acts as the identity.  Orient the upper tree $T$ and let $\beta_e$ be the
scalar incidence row of an edge $e\in T$.  If the upper contact time is
$t_e\le s$ and a lower edge $f$ lands at time $t_f>s$, elimination of the
tree-position rows gives, in one $n$-tangential component,
\begin{equation}\label{eq:common-normal-path-row}
 S_f^{\tan}=\sum_{e\in P_T(f)}
 \sigma_{f,e}(t_f-t_e)\beta_e,
\end{equation}
where $P_T(f)$ is the oriented tree path joining the endpoints of $f$ and
$\sigma_{f,e}\in\{-1,1\}$.  This follows by writing
$\beta_f=\sum_{e\in P_T(f)}\sigma_{f,e}\beta_e$ and subtracting the unique
linear combination of the upper tree rows from the lower row.

The paths for the two lower edges are different, because a tree has a
unique path between two vertices and the selected unordered lower pairs are
different.  If neither path contains the other, choose one private tree
edge from each path.  The corresponding two-by-two minor of
$(S_{f_1}^{\tan},S_{f_2}^{\tan})$ is diagonal up to signs and has modulus at
least $\Delta_1\Delta_2$.  If one path is properly contained in the other,
choose one edge in the smaller path and one private edge of the larger
path; the minor is triangular and has the same lower bound.  Passing between
tree-edge coordinates and orthonormal reduced particle coordinates costs
only a constant depending on $P_0$.

The $d-1$ tangential spatial directions are orthogonal copies of this
minor, and hence supply $2(d-1)$ velocity columns with determinant at least
$c(d,P_0)(\Delta_1\Delta_2)^{d-1}$.  Subtracting their contributions from the two
landing-time columns does not change the resulting $2d$-by-$2d$
determinant.  In the remaining normal direction those columns have the
lower-triangular diagonal $g_i\cdot\omega_i$.  Their determinant is
$\phi_1\phi_2$, proving \eqref{eq:tree-exterior-coercivity}.
Here the $2(d-1)$ selected velocity columns are tangential to $n$ and therefore
have zero entries in the two remaining normal landing rows; this is the
block-zero statement that makes the preceding column subtraction preserve
the normal $2$-by-$2$ determinant.
\end{proof}

\section{Source Concordance and Imported Interfaces}
\label{supp:dependency-source-ledger}

\subsection{Dependency interfaces}
\label{app:ledger}

The seven-part reduction is summarized below in dependency order.

\subsubsection{Logical flow}

\begin{enumerate}[label=\textbf{Part \arabic*.}]
\item The grand-canonical correlation identity uses
\eqref{eq:correlation-def}, Proposition~\ref{prop:ae-flow}, and the exact
forest expansion \eqref{eq:initial-expansion}.
\item Ordinary molecule formation uses the degree-two/three/four operators
and Proposition~\ref{prop:ordinary-subset}.
\item Weight and velocity-volume bounds are dimension-general; the
$d$-dimensional error recurrence changes the exponent to the factor $2(d-1)$ in
\eqref{eq:err-recurrence}.
\item The first two bearing periodic overlap consumers are replaced by
Lemma~\ref{lem:first-failure}.
\item The lower cyclic, 2CONNUP and final DOWN consumers use
\eqref{eq:lower-cyclic-count}, \eqref{eq:2conn-count}, and the protected
restart rule.
\item The exceptional top layer uses Theorem~\ref{thm:p38j},
Lemma~\ref{lem:time-owner}, and Proposition~\ref{prop:v2-consumer}.
\item The full-range chain is Lemmas~\ref{lem:paired-first-contact-partition}
and~\ref{lem:sealed-one-hyperedge},
Theorem~\ref{thm:sealed-past-block-stop-line},
Corollary~\ref{cor:automatic-fixed-k-birth-production}, proved later from
Theorem~\ref{thm:birth-flag-frame-production} independently of the hierarchy,
Lemmas~\ref{lem:paid-cumulant-restart}--\ref{lem:epoch-global-error-bucket},
and Proposition~\ref{prop:epoch-stacked-endpoint};
\eqref{eq:number-layers}--\eqref{eq:qeff-bound} provide the one-epoch inputs.
\end{enumerate}

\subsubsection{Certificate priority rule}
\label{appA:certificate-priority}

The deterministic priority convention behind
Lemma~\ref{lem:certificate-ledger} can be written as the following ordered
rule.  For each legal cut in the fixed atom/slot order, the first applicable
case is used.
\begin{enumerate}[label=\textup{(\arabic*)}]
\item If the cut is a Source-B step-1 or degree-four step-2 exception, place
its original identifier in $\mathcal L_{\rm old}$.
\item Otherwise, if the ordinary adjacent-pair rule succeeds, create one
six-port record in $\mathcal L_{33}$ and protect its active representatives.
\item Otherwise apply Lemma~\ref{lem:first-failure}, create exactly one
identifier in $\mathcal L_{\rm par}$, pad its port list to six, and protect
all active representatives.
\item At a later permanent recut, retain the existing identifier, owner and
port list, and apply \eqref{eq:simultaneous-ledger-transfer} using the
pre-recut value of $g_C$.  A temporary owner or conditional-measure
refinement changes no counter.
\end{enumerate}
This order is independent of the drawings and of later cleanup conventions.

\subsubsection{Epsilon-power bookkeeping}

\begin{center}
\small
\begin{tabular}{@{}lll@{}}
\toprule
Input & Bound & Role\\
\midrule
Initial nonroot link
 & $C_\beta\alpha A\eps$
 & Mayer forest\\
Ordinary serial $j=3$
 & $\cL_\eps^C\eps^2\mu^{-2}$
 & Parts 2--5\\
Two landing contacts
 & $\cL_\eps^C\eps^{2(d-1)}\epstar^{-2(d-1)}$
 & joint relative excess\\
Full packet baseline
 & $\eps^{-(d-1)}$
 & exactly once\\
Rebased nonempty cumulants
 & $\eps^{7/(1000d)}$
 & epoch-stacked main sum\\
Structured dynamical and finite-correlation errors
 & $\eps^{a_{\rm err}}$
 & $\mathcal R_s^t+\mathcal F_s^t$\\
Excluded volume
 & $C_{\beta,d} A s^2\eps^d$
 & hard-core indicator\\
Final rate
 & $\eps^{1/(400d)}$
 & retained minimum\\
\bottomrule
\end{tabular}
\end{center}
Here $a_{\rm err}:=\min\{1/(20d),1/100\}$.

\subsubsection{Root and time identities}

The two normalization identities that must not be conflated are
\[
 I_{H,R}=\eps^{(d-1)|R|}J_H,
 \qquad
 \eps^{-(d-1)P}\eps^{(d-1)(P-1)}=\eps^{-(d-1)}.
\]
The first converts Source-A and Source-B associated integrals; the second is
the unique full-component baseline.  The time gain comes separately from
the projected full time-order domain in \eqref{eq:time-weight-owner}.

\subsubsection{Growth versus fixed powers}

The complete one-epoch growing exponent is
\[
 C_{\rm all}^*
 =\max\{C_1^*,\ldots,C_{15}^*,C_{\rm ord}^*,C_{\rm init}^*\}.
\]
Lemma~\ref{lem:appE-sourceA-Cj} and the diagonal choice give
$C_{\rm all}^*\le L^\delta$.  Here $C_{\rm ord}^*$ owns the cutoff, lift,
dyadic and ordered-cell losses, while $C_{\rm init}^*$ owns the initial
forest hierarchy.  A bounded number of their logarithmic losses is
$\eps^{-o(1)}$.  A one-epoch per-root loss is included in $q_{\rm eff}$ and
controlled by \eqref{eq:qeff-bound}.  The full range instead uses $q_{\rm ep}$
for cumulants and $\eta_{\rm ep}$ with the explicit $2L$ reserve for defects;
see \eqref{eq:epoch-q-log}, \eqref{eq:epoch-defect-radius-log}, and
\eqref{eq:epoch-cumulant-subset-sum}.
\subsection{Imported Source-A / Source-B interfaces}
\label{app:source-interface}

This section fixes the imported part of the argument at equation level.  The
reference sources are arXiv:2408.07818v2, dated 14 November 2024
(\emph{Source B}), and arXiv:2503.01800v1, dated 3 March 2025
(\emph{Source A}).  A statement below attributed to a source means the
statement with precisely that version's hypotheses and notation.  No later
version, and no unlabelled summary in the present paper, is used to enlarge
the imported result.

Source B also has arXiv:2408.07818v3, dated 18 July 2025.  The proof continues
to use version~2 because Source A v1 predates version~3 and its
exceptional-layer argument cites the version-2 proposition numbers, component
classes and exponent ledger.  Version~3 reorganizes that argument
(in particular, v2 Proposition~13.1 becomes the materially different v3
Proposition~15.1).  The complete v3 source is retained for a separate
comparison, but every bearing implication in this paper is checked
against v2.  Both version identifiers and dates are recorded in
\cite{DHMlong,DHMtorus}.  To avoid presenting two arXiv versions of the
same work as different articles, both Source-B versions are recorded under
the single bibliography entry \cite{DHMlong}; the present section carries
the version-specific dates and proposition correspondence.

\subsubsection{Version 2 and the current public version}
\label{appB:version-crosswalk}

The Annals of Mathematics article page records the underlying Source-B
paper as accepted on 31 October 2025 and ``to appear.''  As of the date
of this manuscript, that page does not expose a separate accepted-manuscript
source or final pagination; the latest public arXiv text available for
line-by-line comparison is version~3 of 18 July 2025.  We therefore do not
label v3 as the accepted manuscript.  The following table separates the
v2--v3 semantic correspondence from publication-status metadata.  The
official status page is
\href{https://annals.math.princeton.edu/articles/22284}{the Annals article
record}.

\begin{center}
\scriptsize
\begin{tabular}{@{}p{0.16\textwidth}p{0.20\textwidth}p{0.24\textwidth}
 p{0.30\textwidth}@{}}
\toprule
role used here & reference v2 location & public v3 location & semantic
relation and use here\\
\midrule
initial cumulant package
 & Proposition 5.1 and equations (5.1)--(5.2)
 & Proposition 3.25 plus the v3 initial-link expansion
 & v3 incorporates the initial package into the unified expansion and has
   zero initial truncation error; this paper retains the v2 interface because
   Source A invokes it and proves it independently in Section~3\\
\addlinespace
Boltzmann approximation
 & Proposition 5.2
 & Proposition 6.1
 & same analytic role; v3 is reorganized through Proposition 3.25\\
\addlinespace
nonempty cumulant gain
 & Proposition 5.3, gain $1/(15d)$
 & Proposition 6.2, gain $3^{-d-2}$
 & exponent contract changed; no v3 exponent is inserted into the v2 ledger\\
\addlinespace
truncation error
 & Proposition 5.4, gain $1/(20d)$ and the two-interpoland CH-molecule
   reduction in Section 13
 & Proposition 6.3, gain $3^{-d-2}$ and a unified good/bad reduction in
   Section 15
 & consumer interface materially changed; Proposition~\ref{prop:v2-consumer}
   verifies exactly the v2 route used by Source A\\
\addlinespace
exceptional top layer
 & Proposition 13.1, label \texttt{prop.comb\_est\_extra}
 & Proposition 15.1, same label
 & v2 returns two alternatives in $\{33B\},\{44\},\{4\}$ counts; v3 returns
   one good/bad inequality.  They are not formal renamings\\
\addlinespace
cycle hypothesis
 & tree plus $\gamma$ extra edges, $\gamma>\Gamma$; the application earlier
   proves $\gamma\le2\Gamma$
 & $\Gamma<\rho(\mathcal M)\le2\Gamma$
 & on the connected application graph $\rho=\gamma$; realizability remains
   an inherited support condition in both applications\\
\addlinespace
good-object scale
 & $\eps^{1/(8d)}$, credit $1/(10d)$
 & $\upsilon=3^{-d-1}$, credit $\upsilon/2$
 & quantitatively different; the present fixed-$d$ argument uses the v2 thresholds\\
\addlinespace
Source-A torus argument
 & Source A explicitly cites Source-B Propositions 5.1--5.4 and the
   Section-13/Proposition-13.1 route
 & no public v3-shaped torus argument is asserted
 & the logical arrow used here is Source A $\to$ Source B v2; v3 is a
   comparison source, not a replacement\\
\bottomrule
\end{tabular}
\end{center}

The three reference/comparison files have the following SHA256 values:
\begin{center}
\small
\texttt{4D8CB6042B3C1E312570B4FEE099C090935EFE902F76DC6F2218B8D65CBD0838}
\quad (Source-B v2 TeX),
\\[2pt]
\texttt{8CF03F6BF754A96CC14203B3B5F291C3824D3F34B3C3BD1C1CE79D83369A409F}
\quad (Source-B v3 TeX),
\\[2pt]
\texttt{83C06AA48806934499A995D65C0A942A153D06D042DD922402825F3C0B8A444C}
\quad (Source-A v1 TeX).
\end{center}
The proposition statements, the definitions of good and bad objects,
proof dependencies, and the uses of Source B, Proposition~5.4, were compared against
these files.  If a separately paginated accepted manuscript becomes public,
the bibliographic comparison may be updated without altering the version-2
implication proved here.

\subsubsection{Source-file location notes}
\label{appB:source-location-notes}

The precise source-file locations underlying three comparisons
in the body are as follows.
\begin{enumerate}[label=\textup{(\roman*)}]
\item In Source-B v2, the close-root statement is at
\texttt{main.tex}, lines 2081--2082, and the choice $Z=S$ occurs at lines
2180--2184, immediately after the initial-cumulant formula.  This is the
comparison discussed in Remark~\ref{rem:source-close-repair}.
\item The four Source-A calls corresponding to
$\mathsf U_1,\mathsf U_2,\mathsf{CY},\mathsf{2C}$ occur in
arXiv:2503.01800v1, \texttt{main.tex}, lines 1431, 1437, 1438, and 1439.
The TeX hash is the Source-A value displayed above; the containing archive has
SHA256
\[
\texttt{194D7D050D3DCC8F475B7BAA95DC0F4343D5CA2FD637D2ADC84BFD07D3578138}.
\]
\item In Source-B v2, cutting cases (0), (2), (3), and (4) occur in
\texttt{main.tex}, lines 2879--2883.  The fixed/free mate statement used in
Lemma~\ref{lem:protected-band-restart} is labelled
\texttt{prop.match\_ends} and occurs at lines 2915--2941.
\end{enumerate}

The purpose of recording the interfaces is twofold.  First, it specifies the
positive operator into which the $d$-dimensional packet has to be inserted.
Second, it separates identities already proved in the sources from the new
geometric and combinatorial assertions that must be proved here.  In
particular, the final implication in Proposition~\ref{prop:reduction} may be
used only after all obligations in
Proposition~\ref{prop:source-exact-replacement} have been verified.

\subsubsection{Layer parameters and the four terminal estimates}

Let $L$ be the number of time layers and put
\begin{equation}\label{eq:source-tau}
 \tau=t_{\rm fin}/L,
 \qquad [0,t_{\rm fin}]=\bigcup_{\ell=1}^{L}
 [ (\ell-1)\tau,\ell\tau].
\end{equation}
For exact cardinal stopping, use the Source-B recursion with integer base
\begin{equation}\label{eq:source-hierarchy}
 A_L=L_\eps^\#,
 \qquad \Lambda_\ell=A_\ell^{10d},
 \qquad A_{\ell-1}=\Lambda_\ell^{10d},
\end{equation}
and
\begin{equation}\label{eq:source-theta}
 \theta_0=(100d)^{-1},\qquad
 \theta_\ell=\left(\frac9{10}+\frac{L-\ell}{10L}\right)\theta_0,
 \qquad \theta=\theta_0/2.
\end{equation}
The constants $C$ are independent of $L$ and the constants $C^*$ may depend
on $L$; the increasing family $C_j^*$ is selected before $\eps$ is made
small.  This order is part of the interface.

For $|H|\le A_0$, Source B, Proposition 5.1, gives the initial expansion
with $f^{\mathcal A}(0)=n_0$ and
\begin{equation}\label{eq:source-initial-cumulant}
 |E_H(0,z_H)|
 \le |\log\eps|^{C^*|H|}
 \sum_{Z\subset H}\eps^{|H\setminus Z|/2}
 \prod_{p\in H}n_0(z_p)\,
 \mathbf 1_{\rm close}^{(H,Z)}(z_H),
\end{equation}
where $\mathbf 1_{\rm close}^{(H,Z)}=1$ precisely when every $p\in Z$ has
a distinct partner $p'\in Z$ satisfying
$|x_p-x_{p'}|\le\Lambda_0\eps$.  Its forest remainder obeys
\begin{equation}\label{eq:source-initial-error}
 \|\Err_0\|_{L^1}\le\eps^{\Lambda_0/10}.
\end{equation}
There is a small but consequential gap between the printed statement and
the selector used in the Source-B proof: the latter initially takes $Z$ to
be the singleton-tree roots and only proves closeness to a root in $H$.
Lemma~\ref{lem:corrected-close-selector} repairs this by adjoining the root
of every tree hit by a singleton root.  It proves the printed internal-$Z$
condition and preserves $|F|-|H|\ge|H\setminus Z|$; see
Remark~\ref{rem:source-close-repair}.  Thus
\eqref{eq:source-initial-cumulant}, rather than the insufficient
intermediate selector, is the interface used below.
Source B, Proposition 5.2, defines $f^{\mathcal A}$ through its one-layer
recursion and proves
\begin{equation}\label{eq:source-fa}
 f^{\mathcal A}(\ell\tau)=n(\ell\tau)
 +(f^{\mathcal A})_{\rm rem}(\ell\tau),
 \qquad
 \|(f^{\mathcal A})_{\rm rem}(\ell\tau)\|_{\mathrm{Bol}^{\beta_\ell}}
 \le\eps^{\theta_\ell}.
\end{equation}
The printed outputs of Source B, Propositions 5.3 and 5.4, are
\begin{align}
 \|E_H(\ell\tau)\|_{L^1}
 &\le \eps^{1/(15d)+(C_{14}^*)^{-1}|H|},
 &&H\ne\varnothing,\label{eq:source-cumulant-output}\\
 \|\Err^1_\ell\|_{L^1}&\le\tau^{\Lambda_\ell/10},
 \label{eq:source-layer-error}\\
 \|f_s^{\rm err}(t_{\rm fin})\|_{L^1}
 &\le\eps^{1/(20d)}.\label{eq:source-trunc-output}
\end{align}
The quantitative Source-A estimate used here has
$0<\tau\le\vartheta\le q_0<1$ (after the common choice of constants).
Consequently the printed restart estimate implies the weaker but uniform
form
\begin{equation}\label{eq:source-layer-error-vartheta}
 \|\Err^1_\ell\|_{L^1}
 \le\vartheta^{\Lambda_\ell/10}.
\end{equation}
This replacement is a deduction in the present paper; it is not attributed
to the printed v2 display.

Define the truncated correlation by
$\widetilde f_s:=f_s-f_s^{\rm err}$; this is the exact decomposition used
by the truncation algorithm.  The final aggregation is not a separate black
box.  Source B gives the exact cumulant identity
\begin{equation}\label{eq:source-final-cumulant}
 \widetilde f_s(t_{\rm fin},z_s)
 =\sum_{H\subset[s]}
 (f^{\mathcal A}(t_{\rm fin}))^{\otimes([s]\setminus H)}
 E_H(t_{\rm fin},z_H)+\Err_{L}(z_s),
\end{equation}
The recurrence used below is obtained by combining three facts from the
specified source versions: the initial term
\eqref{eq:source-initial-error} from Source B, Proposition 5.1; the exact
two-sided propagation coefficient $2(d-1)$ in Source B's multilayer cluster
multiplication formula, equation
\texttt{eq.cluster\allowbreak\_mult\allowbreak\_err2}; and the uniformized restart bound
\eqref{eq:source-layer-error-vartheta}.  It is
\begin{equation}\label{eq:source-iterated-error}
 \|\Err_L\|_{L^1}
 \le \eps^{\Lambda_0/10}
       \prod_{k=1}^{L}\eps^{-2(d-1)\Lambda_k^2A_k}
 +\sum_{\ell=1}^{L}\vartheta^{\Lambda_\ell/10}
       \prod_{k=\ell+1}^{L}\eps^{-2(d-1)\Lambda_k^2A_k}.
\end{equation}
Here the first summand is the initial forest remainder, and the remaining
summands are the propagated one-layer restart errors.  The two-sided
coefficient is $2(d-1)$, exactly as in
\eqref{eq:err-recurrence}.  Notice the distinction: the source's local
restart bound is first printed with $\tau$, while the present diagonal
recurrence legitimately weakens it to $\vartheta$ because
$\tau\le\vartheta$; the initial term remains an epsilon power.
Only after \eqref{eq:source-cumulant-output}--\eqref{eq:source-iterated-error}
and \eqref{eq:source-fa} are available does the subset sum reduce to the
Boltzmann tensor product.  Inserting the hard-core indicator costs
\begin{equation}\label{eq:source-excluded-volume}
 \left\|n(t_{\rm fin})^{\otimes s}
 (1-\mathbf 1_{\cD_s^\circ})\right\|_{L^1}
 \le Cs^2\eps^d.
\end{equation}

Equations \eqref{eq:source-fa}--\eqref{eq:source-excluded-volume} have, at
this point, exactly the time ownership present in Source B: the first two
families are asserted at layer endpoints and the truncation and final
cumulant identities are asserted at the chosen terminal time
$t_{\rm fin}$.  They do not by themselves give a statement at every
interior time of one previously fixed decomposition.  For a prescribed
$T\in(0,t_{\rm fin}]$, the present paper instead reruns the same finite
structural expansion on $[0,T]$, with the same number of layers and with
$\tau_T=T/L$.  We denote its objects by
$f^{\mathcal A,T},E_H^T,f_s^{{\rm err},T},\Err_L^T$ when the distinction is
needed.  The small-window estimate
\eqref{eq:source-weight-projection}, rather than the source-only
fixed-window formula \eqref{eq:source-weight-projection-original}, is what
makes this rerun uniform as $T\downarrow0$.  The uniform quantifiers and the
$T=0$ endpoint are proved in Lemma~\ref{lem:uniform-terminal-rerun}; hence no
all-time conclusion is being attributed to Source B.

\subsubsection{The molecule integral and the positive kernel}

For a layered interaction diagram let $\mathsf X$ be the number of new
particle lines, $\mathsf U$ the total number of incoming particle lines and
$\mathsf Z$ the number of collisions plus overlaps.  In Source B notation,
\begin{equation}\label{eq:source-XUZ}
 \mathsf U=\sum_{\ell'}s_{\ell'},\qquad
 \mathsf X=\sum_{\ell'}(s_{\ell'}^+-s_{\ell'}),\qquad
 \mathsf Z=|E(\mathcal G)|.
\end{equation}
The total recollision count is $\mathsf Z-\mathsf X$, and the parameter that
enters the molecule bounds is
\begin{equation}\label{eq:source-rho}
 \rho=|H'|+\mathsf Z-\mathsf X.
\end{equation}

Let $\mathcal M=(M,E,P)$ be the molecule, $E_*\subset E$ the bonds and free
ends, and $E^-_{\rm end}$ its bottom ends.  At an atom $a$, write
$(e_1,e_2)$ for its bottom pair and $(e_1',e_2')$ for the serial top pair.
With
\[
 r_a=x_{e_1}-x_{e_2}+t_a(v_{e_1}-v_{e_2}),
 \qquad \omega_a=r_a/\eps,
\]
the C-atom kernel is the product of the two transport deltas, the contact
delta and flux, and the two elastic-scattering deltas:
\begin{align}
 \Delta_a^{\rm C}
 &=\delta^{(d)}(x_{e_1'}-x_{e_1}+t_a(v_{e_1'}-v_{e_1}))
   \delta^{(d)}(x_{e_2'}-x_{e_2}+t_a(v_{e_2'}-v_{e_2}))\notag\\
 &\quad\times\delta(|r_a|-\eps)
 [(v_{e_1}-v_{e_2})\cdot\omega_a]_-\notag\\
 &\quad\times\delta^{(d)}\!\left(v_{e_1'}-v_{e_1}
 +[(v_{e_1}-v_{e_2})\cdot\omega_a]\omega_a\right)\notag\\
 &\quad\times\delta^{(d)}\!\left(v_{e_2'}-v_{e_2}
 -[(v_{e_1}-v_{e_2})\cdot\omega_a]\omega_a\right).
 \label{eq:source-C-kernel}
\end{align}
For an O-atom the last two factors in \eqref{eq:source-C-kernel} are replaced
by $\delta^{(d)}(v_{e_1'}-v_{e_1})
\delta^{(d)}(v_{e_2'}-v_{e_2})$.  On a torus lift cell, $r_a$ is replaced
by $r_a-m_a$ for one fixed $m_a\in\Z^d$.

The time domain is the full partial-order domain
\begin{equation}\label{eq:source-time-domain}
 \mathcal D_M=\left\{t_M:
 (\ell[a]-1)\tau<t_a<\ell[a]\tau,\quad
 t_a<t_{a^+}\text{ whenever }a^+\text{ is a parent of }a\right\}.
\end{equation}
The Source-B associated integral is
\begin{multline}\label{eq:source-associated-integral}
 \mathcal I\mathcal N(\mathcal M,H,H')
 =\eps^{-(d-1)(|E|-2|M|-|H|)}
 \int \mathbf 1_{\mathcal D_M}(t_M)
 \prod_{a\in M}\Delta_a\\
 \times\prod_{e\in E^-_{\rm end}\setminus H'}
 \left|f^{\mathcal A}
 ((\ell_1[e]-1)\tau,x_e-(\ell_1[e]-1)\tau v_e,v_e)\right|
 \left|E_{H'}(0,(x_e,v_e)_{e\in H'})\right|
 \dd z_E\dd t_M.
\end{multline}
This formula fixes the normalization, the live variables and the original
time owner.  The identity $|E|-2|M|=\mathsf X+|H|$ converts it into
$\eps^{(d-1)|H|}$ times the full-molecule integral used in Source B,
Proposition 8.10.

After using \eqref{eq:source-initial-cumulant} and
\eqref{eq:source-fa}, Source B decomposes the integrand into positive terms.
For fixed $Z\subset H'$, sets $\mathcal W,\mathcal V$ of bottom ends and a
fixed collection of generalized links, the positive factor has the form
\begin{align}
 Q_{Z,\mathcal W,\mathcal V}
 &=\eps^{|H'\setminus Z|/2}\mathbf 1_{\mathcal D_M}(t_M)
 \prod_{e\in\mathcal W}|n(\ell_e\tau,\Phi_ez_e)|
 \prod_{e\in\mathcal V}|(f^{\mathcal A})_{\rm rem}
 (\ell_e\tau,\Phi_ez_e)|\notag\\
 &\quad\times\prod_{e\in H'}n_0(z_e)
 \prod_{(e,e')\in\mathcal L}
 \mathbf 1_{|x_e-x_{e'}|\le\eps^{1-1/(8d)}}.
 \label{eq:source-positive-Q}
\end{align}
Here $\Phi_e$ is the prescribed free transport to the bottom of the initial
layer.  Every $e\in Z$ belongs to at least one and at most two generalized
links.  The number of positive terms is bounded by
$C^{\mathsf X}|\log\eps|^{C^*\rho}$.  Consequently, an estimate valid only
for a factorized or boundary-supremized test function is insufficient: the
new local estimate must hold for every nonnegative function of all variables
in \eqref{eq:source-positive-Q}, with the same
$\mathbf 1_{\mathcal D_M}$.

\subsubsection{Cutting, elementary operators, weights and volume}

If $\mathcal M'$ is obtained from $\mathcal M$ by a Source-B cutting
sequence and $\operatorname{Sh}:E\to E_*'$ is the shadow map, Source B,
Proposition 8.10, proves the exact positive-kernel identity
\begin{multline}\label{eq:source-cutting-identity}
 \eps^{-(d-1)(|E|-2|M|)}\int\prod_{a\in M}\Delta_a\,Q(z_E,t_M)
 \dd z_E\dd t_M\\
 =\eps^{-(d-1)(|E_*'|-2|M|)}
 \int\prod_{a\in M}\Delta_a(\widetilde z)
 Q((\widetilde z_{\operatorname{Sh}(e)})_{e\in E},t_M)
 \dd\widetilde z_{E_*'}\dd t_M.
\end{multline}
It is essential that $Q$ is not replaced by a product in this identity.
Deleting an O-atom preserves $|E|-2|M|$ and gives the corresponding
inequality after dropping its contact indicator.

Source B, Proposition 9.1, supplies the following one-atom operations.
A degree-two component is an exact collision-root substitution of norm at
most one.  With one fixed end, a degree-three component is represented by
\begin{equation}\label{eq:source-degree-three}
 \int [(v_1-v_2)\cdot\omega]_-Q\,
 \dd t\dd\omega\dd v_2,
 \qquad \omega\in\mathbb S^{d-1}.
\end{equation}
A full degree-four component has the same integral with the additional
variables $(x_1,v_1)$ and the single normalization $\eps^{-(d-1)}$.
The proposition also states the small-time, small-relative-velocity and
small-position support gains, together with the exclusions forced by serial
O-atom slots.  Source B, Proposition 9.2, gives the two-atom
$\{33{\rm A}\}$, $\{33{\rm B}\}$ and $\{44\}$ operators.  Each is a finite
sum
\begin{equation}\label{eq:source-elementary-operator}
 \mathcal J(\mathcal X;Q)
 \le\kappa_{\mathcal X}\int_{\Omega_{\mathcal X}}Q^\flat\dd w
 \le\kappa_{\mathcal X}\int_{\Omega'_{\mathcal X}}Q^\flat\dd w,
\end{equation}
where $\Omega'$ is independent of the fixed boundary values,
$|\Omega'|\le|\log\eps|^{C^*}$.  For each of the proposition's permitted
support exponents $0<\upsilon<1/2$, the corresponding good support gives
$|\kappa_{\mathcal X}||\Omega_{\mathcal X}|
\le\eps^\upsilon|\log\eps|^{C^*}$.  The cutting algorithm later chooses
$\upsilon=1/(8d)$, and after the fixed logarithmic absorption the final
ledger deliberately retains only $\eps^{1/(9d)}$ per good component.  The
full list of slot restrictions is part of Proposition 9.2; it may not be
replaced by the component name alone.  Importantly, the product
$|\kappa_{\mathcal X}||\Omega_{\mathcal X}|$ is a conditional measure
bound, not an $L^1$ operator norm for an arbitrary kernel.  In particular,
\eqref{eq:source-elementary-operator} does \emph{not} imply
\[
 \mathcal J(\mathcal X;Q)
 \le \eps^\upsilon|\log\eps|^{C^*}
       \int_{\Omega'_{\mathcal X}}Q^\flat
\]
before a pointwise source majorant has been established.  In every consumer
below the complete positive kernel and its indicators are therefore retained
through cutting, packet coarea, weight absorption and the Maxwellian bound.
Only then is the conditional measure of the good cell used.  This order is
the measure-level content of the source proof and prevents a small support
volume from being pulled through an arbitrary $Q$.

Periodic certificates from $\mathcal L_{\rm par}$ are not silently placed
among these ordinary Source-B support alternatives.  If such a certificate
carries a new-line owner, that time remains in the distinguished projected
set.  After the same pointwise source majorant has been formed, the
different-lift restriction is integrated in its conditional velocity
variable using \eqref{eq:multiple-contact-tube}.  Thus its ordinary time
factor and its one periodic good factor come from two explicitly different
conditional coordinates and each is used once.  A periodic good $\{4\}$
singleton is nevertheless a full degree-four Source-B component: it
increments $N_4$ and retains the unique $\eps^{-(d-1)}$ baseline.  Its tube
is the good excess relative to that baseline; good status never deletes the
baseline.

Source B, Proposition 9.3, absorbs all collision weights.  For every atom
set $A\subset M$ and $|\log\eps|^{-1}\le\gamma\le1$, there is a scalar
$B_M\ge1$ such that its weight estimate and its printed fixed-window time
estimates are
\begin{align}
 \prod_{a\in M}(1+|v_{e_1(a)}-v_{e_2(a)}|)
 &\le(C\gamma^{-1/2})^{|M|}|\log\eps|^{C^*\mathsf R}
 (C^*)^{\mathsf U}B_M
 \exp\!\left(\gamma\sum_{e\in E^-_{\rm end}}|v_e|^2\right),
 \label{eq:source-weight-1}\\
 B_M\left(\int_{\widetilde{\mathcal D}_A}1\dd t_{M\setminus A}
 \right)^{5/6}
 &\le\tau^{5|M|/6}|\log\eps|^{C^*|A|},
 \label{eq:source-weight-projection-original}\\
 B_M\left(\int_{\mathcal D_M}1\dd t_M\right)^{5/6}
 &\le\tau^{5|M|/6}.
 \label{eq:source-weight-2}
\end{align}
The domain $\widetilde{\mathcal D}_A$ is the projection of the full domain
\eqref{eq:source-time-domain}: the omitted $A$-times must exist.  Source B
obtains \eqref{eq:source-weight-projection-original} in its fixed-window
regime by using the extra comparison $\tau^{-1}\le |\log\eps|$.  It is not a
scale-uniform statement when the terminal window is allowed to shrink.

The present paper replaces only that projected-time line by the following
small-window estimate:
\begin{equation}\label{eq:source-weight-projection}
 B_M\left(\int_{\widetilde{\mathcal D}_A}1\dd t_{M\setminus A}
 \right)^{5/6}
 \le\tau^{5(|M|-|A|)/6}|\log\eps|^{C^*|A|}.
\end{equation}
This is not attributed to Source B.  It is proved in this paper by the tree
simplex and projection argument of Lemmas~\ref{lem:tree-time-owner}--
\ref{lem:alternating-tree-decomposition}; its precise consumer form is
Lemma~\ref{lem:distinguished-time-weight}.  The existential qualification is
the time-owner constraint.  Every consumer below designates, before coarea,
the birth/forest times which own the required line credits and retains them
as output coordinates.  A time used as a contact-root pivot is never restored
by a bare finite-fibre argument.  Thus
\eqref{eq:source-weight-projection-original} records the source formula,
whereas \eqref{eq:source-weight-projection} records the new uniform repair.

The repaired projection is not used with an arbitrary omitted set and then
asked to reproduce the old time ledger.  Before elementary integration,
Proposition~\ref{prop:packet-free-small-window-consumer} marks the unique
creation time of every new physical line.  Each mark is either retained for
the common existential projection or is consumed, after pointwise source
majorization, by an owner/good estimate on the same conditional measure.
The two classes are disjoint and exhaust all \(\mathsf X\) new lines.  This
admissibility statement is essential in the extreme case \(A=M\): the
zero-dimensional projection supplies no \(\tau\) power, so every new-line
credit must already have been supplied by a common-measure owner.  A finite
coarea fibre is never used to manufacture a missing time power.

Source B, Proposition 9.4, controls simultaneously the free velocity in
every normal degree-three component.  If $\mathcal F$ is their selected set
of free ends and $X_e$ are the dyadic bounds of original free ends, then
there is a set $Y\subset(\R^d)^{|\mathcal F|}$, independent of the concrete
fixed velocities, such that
\begin{equation}\label{eq:source-volume}
 |Y|\le \prod_{e\in E_{\rm free}(\mathcal M)}X_e^d
 C^{|M|}|\log\eps|^{C^*K},
\end{equation}
where $K$ counts components outside the selected normal degree-three family.

For a packet-free ordinary branch,
Proposition~\ref{prop:packet-free-small-window-consumer} is inserted into the
same reverse cutting order.  It converts the retained physical time volumes
and their distinct atom rates into \(\vartheta_T^{\mathsf X}\), while leaving
the epsilon, good-component and full-baseline ledgers unchanged.  After the
components are integrated in reverse order, the resulting uniform estimate
is the following; at the original terminal time \(T=t_{\rm fin}\) we write
\(\vartheta_T=\vartheta\):
\begin{equation}\label{eq:source-iterated-integral}
 |\widetilde{\mathcal I}(\mathcal M')|
 \le\eps^{-(d-1)(|E_*'|-2|M|)}
 \int\kappa_{\rm tot}Q_1\,
 \dd w_{r+q}\cdots\dd w_1.
\end{equation}
The domain of a good component at position $j$ may depend only on variables
$w_k$ with $k<j$.  Reordering the normal degree-three variables is legal only
with this dependency property.  The final ledger is
\begin{align}
 |\mathcal I\mathcal N(\mathcal M,H,H')|
 &\le\vartheta^{\mathsf X}\eps^{\Delta}
 |\log\eps|^{C^*(\rho+N_{\rm good}+|H'\setminus Z|)},
 \label{eq:source-general-ledger}\\
 \Delta&=\frac14(N_{\rm ee}+|H'\setminus Z|)
 +\frac1{9d}N_{\rm good}
 -\max(d\nu,(d-1)\nu)-N_{\rm del},
 \label{eq:source-Delta}\\
 \nu&=N_4+N_{\rm ee}-|H|.\label{eq:source-nu}
\end{align}
The ordinary cutting theorem supplies the lower bound on $\Delta$.  Thus a
new exceptional component must be inserted before this ledger is evaluated;
it cannot be justified merely by comparing its final epsilon exponent with
the target.

\subsubsection{The exceptional top-layer consumer}

The following is the exact point at which Source A modifies Source B.  The
Source-B first-threshold selector has two exhaustive exceptional routes:
the unchanged high-particle route $r+m>\Lambda_\ell$, with
$r+m\le2\Lambda_\ell$, and the high-recollision route
$\Gamma<\gamma\le2\Gamma$.  In the high-particle route Source B first proves
$\mathsf X+m>\Lambda_\ell$ and uses the ordinary UP estimate.  This route is
combinatorially unchanged, but its small-window analytic input is now
Proposition~\ref{prop:packet-free-small-window-consumer}.  For the joined
ordinary molecule set \(N:=\mathsf X+m\).  Its nonroot birth count is
\(\widetilde{\mathsf X}=N-\iota\), and the ordinary UP molecule estimate,
after the standard weakening used in the first source interpoland, contains
\(\vartheta_T^{(N-\iota)/9}\).  The part needed for the target and the spare
part separate exactly because
\[
 \frac{N-\iota}{9}-\frac N{20}-\frac N{100}
 =\frac{46N-100\iota}{900}\ge0.
\]
Indeed this is immediate when \(\iota=0\); when \(\iota=1\), the integer
inequality \(N>\Lambda_\ell\ge2\) gives \(N\ge3\).  Hence
\[
 \vartheta_T^{(N-\iota)/9}
 \le \vartheta_T^{N/20}q_0^{N/100}
 \le \vartheta_T^{N/20}q_0^{\Lambda_\ell/100}.
\]
The hierarchy has
\(\min_\ell\Lambda_\ell=(L_\eps^\#)^{10d}\ge L_\eps^{10d}\), so one common
main-theorem threshold makes
\[
 q_0^{\Lambda_\ell/100}\le\eps^{2d}
\]
for every layer.  The dimension-uniform credit $2d$ more than pays the
displayed $\eps^{d-1}$ target and the finite subpower losses, while the
ordinary UP epsilon gain pays the
remaining \(1/(20d)+(C_{15}^*)^{-1}(\rho+r)\) part exactly as in the cited
source calculation.  Thus this route remains in (R3),
and already implies
\eqref{eq:source-top-target}; it never calls the two-landing packet.  The
replacement below is restricted to the second, high-recollision route.

In that route the exceptional layer $\ell$ contains $r+m$ particle lines
and has a collision graph equal to a tree plus $\gamma$ edges, with
$\Gamma<\gamma\le2\Gamma$.  Joining this layer to
the lower diagram gives a molecule $\widetilde{\mathcal M}$.  Let $m$ be the
number of newly inserted lines, $r$ the size of the distinguished root
cluster, $\iota=0$ for $r>0$ and $\iota=1$ for $r=0$, and retain the lower
parameters $(\mathsf X,\rho)$.  Source B reduces the truncation error to
\begin{equation}\label{eq:source-top-error-reduction}
 \|f_s^{\rm err}(t_{\rm fin})\|_{L^1}
 \le C^{\mathsf X+m}|\log\eps|^{C^*(r+\rho+1)}
 \eps^{-\iota(d-1)}
 \sup_{\widetilde{\mathcal M}}
 \mathcal I\mathcal N(\widetilde{\mathcal M},\widetilde H,H').
\end{equation}
The required molecule estimate is
\begin{equation}\label{eq:source-top-target}
 \mathcal I\mathcal N(\widetilde{\mathcal M},\widetilde H,H')
 \le\vartheta^{(\mathsf X+m)/20}\eps^{d-1}
 \eps^{1/(20d)+(C_{15}^*)^{-1}(\rho+r)}.
\end{equation}
It is obtained from the following two estimates:
\begin{align}
 \mathcal I\mathcal N(\widetilde{\mathcal M},\widetilde H,H')
 &\le\vartheta^{(\mathsf X+m)/10}
 \eps^{1/(18d)+(C_{14}^*)^{-1}(\rho+r)},
 \label{eq:source-top-first-interpoland}\\
 \mathcal I\mathcal N(\widetilde{\mathcal M},\widetilde H,H')
 &\le\vartheta^{(\mathsf X+m)/10}\eps^{d-1+1/(18d)}
 |\log\eps|^{C^*(\rho+r)}.
 \label{eq:source-top-second-interpoland}
\end{align}

The word ``interpolation'' in the source is a two-regime minimum argument,
not a fixed-weight geometric mean.  We record the elementary step because a
literal geometric mean would lose half of the factor $\eps^{d-1}$.

\begin{lemma}[Two-regime source envelope]
\label{lem:source-two-regime-envelope}
Put $R=\rho+r$, $a=1/(18d)$, $b=1/(20d)$ and
$\delta=a-b=1/(180d)$.  Choose $C_{15}^*$ after $C_{14}^*$ so that, with
$c=(C_{14}^*)^{-1}$ and $c'=(C_{15}^*)^{-1}$,
\begin{equation}\label{eq:source-constant-separation}
 c'<c,
 \qquad \frac{(c-c')\delta}{2c'}\ge d-1.
\end{equation}
Let $L_\eps=|\log\eps|$.  The implication is uniform along the final
diagonal choice provided, for some fixed $0<\eta<1/10$,
\begin{equation}\label{eq:source-uniform-growing-constants}
 C^*\le L_\eps^\eta,
 \qquad C_{15}^*=(c')^{-1}\le L_\eps^\eta,
\end{equation}
where $C^*$ is the coefficient in the logarithmic exponent of the second
interpoland.  More precisely, there is $\eps_0(d,\eta)>0$, independent of
$R$ and of the allowed constants in
\eqref{eq:source-uniform-growing-constants}, such that for
$0<\eps<\eps_0$ the two bounds
\eqref{eq:source-top-first-interpoland}--
\eqref{eq:source-top-second-interpoland} imply
\eqref{eq:source-top-target}.
\end{lemma}

\begin{proof}
If $R\le\delta/(2c')$, then
$a-b-c'R\ge\delta/2$.  Since $R$ is bounded by a constant independent of
$\eps$ only in the fixed-layer formulation, we instead use the uniform
diagonal bounds.  Equations
\eqref{eq:source-uniform-growing-constants} give
\[
 R\le\frac{\delta}{2}C_{15}^*\le C_dL_\eps^\eta,
 \qquad
 \log\bigl(L_\eps^{C^*R}\bigr)
 \le C_dL_\eps^{2\eta}\log L_\eps=o(L_\eps).
\]
Because $2\eta<1$, for all sufficiently small $\eps$, uniformly in every
allowed $R,C^*,C_{15}^*$, the last expression is at most
$(\delta/3)L_\eps$.  Hence the polylogarithmic factor is at most
$\eps^{-\delta/3}$.  The remaining gap
$a-b-c'R\ge\delta/2$ proves that the second bound gives the epsilon power in
\eqref{eq:source-top-target} with room $\delta/6$.

If $R>\delta/(2c')$, then \eqref{eq:source-constant-separation} gives
\[
 a+cR-(d-1+b+c'R)
 =\delta-(d-1)+(c-c')R
 >\delta-(d-1)+\frac{(c-c')\delta}{2c'}
 \ge\delta>0,
\]
so the first bound gives the required epsilon power.  In either regime,
  $0<\vartheta\le1$ makes $\vartheta^{(\mathsf X+m)/10}$ no larger than
$\vartheta^{(\mathsf X+m)/20}$.  This proves the claim.
\end{proof}

In the attached case, the first interpoland is the ordinary molecule bound
\begin{equation}\label{eq:source-attached-first}
 \mathcal I\mathcal N(\widetilde{\mathcal M},\widetilde H,H')
 \le\vartheta^{\widetilde{\mathsf X}/9}
 \eps^{1/(12d)+4(C_{14}^*)^{-1}\widetilde\rho},
 \quad
 \widetilde{\mathsf X}=\mathsf X+m-\iota,
\end{equation}
where $\rho+r\le\widetilde\rho\le\rho+r+3\Gamma$.  Its small-window owner input is
Proposition~\ref{prop:packet-free-small-window-consumer}; the remaining
ordinary source argument and its standard interpolation weaken the full
owner power to the displayed \(\widetilde{\mathsf X}/9\) exponent.  Thus
\eqref{eq:source-attached-first} does not call the exceptional packet-time
lemma.  The second interpoland
comes from the exact ledger
\begin{multline}\label{eq:source-attached-ledger}
 \mathcal I\mathcal N(\widetilde{\mathcal M},\widetilde H,H')
 \le\vartheta^{\widetilde{\mathsf X}}
 \eps^{(d-1)(|\widetilde H|-N_4-N_{\rm ee})}
 |\log\eps|^{C^*(\widetilde\rho+N_{\rm ee}+N_{\neg N3})}\\
 \times\eps^{(9d)^{-1}N_{\rm good}+|H'\setminus Z|/4}.
\end{multline}
Here $N_4=N_4^1+N_4^2$ separates the top-layer and lower-layer degree-four
components.  The top combinatorics controls $N_4^1$ by good components; the
lower algorithm and the fact that cutting the attached top creates at least
one lower fixed end give
\begin{equation}\label{eq:source-attached-baseline}
 N_4^2+N_{\rm ee}\le|\widetilde H|-1.
\end{equation}
The strict $-1$ is the unique remaining full-component baseline.

In the disjoint case the integral factorizes as
\begin{equation}\label{eq:source-disjoint-factorization}
 \mathcal I\mathcal N(\widetilde{\mathcal M},\widetilde H,H')
 \le \mathcal I\mathcal N(\widetilde{\mathcal M}_\ell,G,\varnothing)
 \mathcal I\mathcal N(\widetilde{\mathcal M}_{<\ell},K,H'),
\end{equation}
where $|G|=r+\iota$ and $K$ is the lower root set.  The lower factor obeys
\begin{equation}\label{eq:source-disjoint-lower}
 \mathcal I\mathcal N(\widetilde{\mathcal M}_{<\ell},K,H')
 \le\vartheta^{\mathsf X/9}
 \eps^{1/(12d)+4(C_{14}^*)^{-1}\rho},
\end{equation}
or equals one if the lower molecule is empty.  In the nonempty case the
ordinary small-window owner input is
Proposition~\ref{prop:packet-free-small-window-consumer}, applied only to the
lower factor; the top packet variables do not enter this projection.  The
required top estimate is
\begin{equation}\label{eq:source-disjoint-top}
 \mathcal I\mathcal N(\widetilde{\mathcal M}_\ell,G,\varnothing)
 \le\vartheta^{m/9}\eps^{d-1+1/(18d)+4(C_{14}^*)^{-1}r}.
\end{equation}

The top combinatorial interface in Source B is modeled on Source A,
Proposition~3.8.  Its printed input is a full, one-layer C-molecule whose bonds form a tree of at
most $|\log\eps|^{C^*}$ atoms plus $\gamma$ bonds,
$\Gamma<\gamma<2\Gamma$.  The Source-B error selector also permits the
endpoint $\gamma=2\Gamma$.  Theorem~\ref{thm:p38j} supplies the required
fixed-$d\ge4$ result by rerunning the graph partition and replacing every
analytic leaf by an internal all-dimensional estimate.  It fixes an integer
Source-B threshold \(\Gamma\), sets
\(\widehat\Gamma=\Gamma+\tfrac12\), and uses the elementary implication
\[
 \Gamma<\gamma\le2\Gamma
 \Longrightarrow
 \widehat\Gamma<\gamma<2\widehat\Gamma
\]
for the integer \(\gamma\).  After at most
$C^{|M|}|\log\eps|^{C^*}$ positive indicator pieces, its three alternatives
are: a good $\{44\}$ component and no $\{33{\rm B}\}$ or degree-four
component; the quantitative good-versus-degree-four inequality; or exactly
one degree-four component and at most ten excess-bearing components with
\begin{equation}\label{eq:source-A-third-alternative}
 \prod_j\sigma_j\le\eps^{d-1+1/(15d)}.
\end{equation}
Every replacement of the third alternative must retain the input class,
indicator count, unique degree-four baseline, nonfull cleanup and the
arbitrary positive kernel in \eqref{eq:source-associated-integral}.

\clearpage
\subsubsection{Cross-reference map for the imported argument}
\label{appB:complete-call-matrix}

The next table makes the source correspondence in
Proposition~\ref{prop:reduction} explicit without
requiring the reader to infer which source call is being replaced.  Line
locations refer to the TeX sources distributed with
arXiv:2408.07818v2 and arXiv:2503.01800v1; the arXiv proposition/equation
labels are the cited-source identifiers.  The variables column lists the
variables live at the call, not merely the variables displayed in the final
estimate.

\begingroup
\fontsize{7.5}{8.4}\selectfont
\setlength{\tabcolsep}{3pt}
\begin{longtable}{@{}p{0.14\textwidth}p{0.16\textwidth}p{0.20\textwidth}
 p{0.21\textwidth}p{0.20\textwidth}@{}}
\toprule
source call & present replacement & variables & normalization/time owner
& exact output\\
\midrule
\endfirsthead
\toprule
source call & present replacement & variables & normalization/time owner
& exact output\\
\midrule
\endhead
\midrule
\multicolumn{5}{r}{\textit{continued on the next page}}\\
\endfoot
\bottomrule
\endlastfoot
Source B (2.1)--(2.2), ensemble definition
& \eqref{eq:partition}--\eqref{eq:grand-density} and
  \eqref{eq:correlation-def}
& all \(z_N\), sector \(N\), activity \(a_\eps\)
& \(a_\eps=\alpha\eps^{-(d-1)}\);
  correlation scaling \((\alpha^{-1}\eps^{d-1})^s\)
& normalized periodic correlation and hard-core sector measure\\
\addlinespace
Source B, Proposition~5.1
& Proposition~\ref{prop:initial-cumulant} and
  Lemma~\ref{lem:corrected-close-selector}
& initial roots, active-forest nonroots, Mayer tree links and cross-pair
  selectors
& no time owner; each nonroot insertion carries one activity factor
& \eqref{eq:source-initial-cumulant} and
  \eqref{eq:source-initial-error}\\
\addlinespace
Source B, Propositions~5.2--5.3
& Sections~3--5
& complete layered molecule, fixed ends, all live atom times
& \(\eps^{-(d-1)(|E|-2|M|-|H|)}\);
  original full partial-order domain
& \eqref{eq:source-fa}, \eqref{eq:source-cumulant-output},
  \eqref{eq:source-layer-error}\\
\addlinespace
Source B, Proposition~5.4 and Section~13.1--13.2
& packet-free owner Proposition~\ref{prop:packet-free-small-window-consumer}
  and unchanged UP/high-particle combinatorics; Theorem~\ref{thm:p38j},
  Proposition~\ref{prop:v2-consumer}, and Section~\ref{app:consumer-verification} for high recollision
& joined top/lower molecule, \(G,K,H'\), complement variables
& unique attached baseline or rooted factor
  \(\eps^{(d-1)(|G|-1)}\); exact projected time domain
& \eqref{eq:source-top-target}, hence
  \eqref{eq:source-trunc-output}\\
\addlinespace
Source B, Proposition~8.10
& \eqref{eq:source-cutting-identity} and
  Lemma~\ref{lem:appF-shadow-compatibility}
& every original edge state and its current shadow
& preserves \(|E|-2|M|\); no time integration
& identical nonnegative joint kernel after cutting\\
\addlinespace
Source B, Propositions~9.1--9.2
& Proposition~\ref{prop:ordinary-subset}
& fixed ends, component atom times/normals, free component velocities
& degree-four baseline \(\eps^{-(d-1)}\) once;
  original component time interval
& all degree \(2,3,4\) operators;\newline
  \(\{33A\},\{33B\},\{44\}\) operators\\
\addlinespace
Source B, Proposition~9.3
& Lemmas~\ref{lem:tree-time-owner},
  \ref{lem:distinguished-time-weight}, and
  Proposition~\ref{prop:packet-free-small-window-consumer};
  Lemma~\ref{lem:time-owner} only for the exceptional packet
& ordinary kept birth times \(t_S\), eliminated times \(t_A\), all weights
& existential projection
  \(\{t_S:\exists t_A,\ (t_A,t_S)\in\mathcal D_M\}\)
& source formula \eqref{eq:source-weight-projection-original}; present
  small-window repair \eqref{eq:source-weight-projection};
  $\vartheta^{\mathsf X}$ consumer with physical projections retaining
  $\tau$\\
\addlinespace
Source B, Proposition~9.4
& Proposition~\ref{prop:ordinary-subset} and
  Lemma~\ref{lem:global-innovation-partition} for the ordinary route;
  Proposition~\ref{prop:appF-I1} for the exceptional packet
& all selected normal degree-three velocities and fixed boundary velocities
& branch labels are finite counting variables, not velocity volume
& boundary-independent simultaneous containing set\\
\addlinespace
Source B, equations (9.50), (9.53)--(9.54)
& \eqref{eq:source-iterated-integral}--\eqref{eq:source-nu}
& reverse cutting variables \(w_{r+q},\ldots,w_1\)
& dependency order retained; one full baseline only
& general epsilon ledger \(\Delta\)\\
\addlinespace
Source B, Propositions~10.2(1)--(2)
& Proposition~\ref{prop:periodic-up-local}
& periodic first-failure ports and current analytic representatives
& no new baseline under a protected recut
& unchanged UP local alternatives\\
\addlinespace
Source B, Proposition~10.9
& Proposition~\ref{prop:cyclic-certificate-injection}
& cyclic lower components, bearing record and one of six ports
& combinatorial only
& \(N_{33A}+N_{\rm sing}\ge|\operatorname{Cyc}|/6\)\\
\addlinespace
Source B, Propositions~11.2 and 2CONNUP
& Proposition~\ref{prop:periodic-2connup}
& \(\mathcal L_0,\mathcal L_1,\mathcal P_1,\mathcal P_2\), certificate owner
& old exceptions and periodic certificates are disjoint
& exact coefficient \(50d^2\) in \eqref{eq:2conn-count}\\
\addlinespace
Source B, Proposition~12.1; Source A, Part~5
& Proposition~\ref{prop:periodic-down-tail}
& original recollisional blocks and child-assigned connectors
& choose one protected charge per recollisional block; the chosen charges
  have distinct block owners
& one distinct good object per recollisional block\\
\addlinespace
Source A, Proposition~3.8
& Theorem~\ref{thm:p38j}
& complete special-class molecule, two sublayers, arbitrary joint \(Q\)
& packet factor
  \(\eps^{-(d-1)}\eps^{2(d-1)}\epstar^{-2(d-1)}\);
  fibrewise time owner retained
& operator-level replacement, not a componentwise third alternative\\
\addlinespace
Source B, final cumulant identity
& Section~\ref{sec:global} and Appendix~\ref{app:parameters}
& all subsets \(H\subset[s]\), \(s\le|\log\eps|\), all layers
& v2 exponents only; \(q_{\rm eff}\) absorbs repeated losses
& Theorem~\ref{thm:main} with rate \(\eps^{1/(400d)}\)\\
\end{longtable}
\endgroup

Every bearing imported call appears in the table.  The rows involving the
new packet are expanded further in Section~\ref{app:consumer-verification}; the remaining rows are
proved in Sections~3--5 or are literal algebraic identities displayed in
this section.  In particular, there is no unstated fifth
dimension-sensitive input hidden behind Proposition~\ref{prop:reduction}.

\section{Compatibility of the Exceptional-Top Insertion}
\label{app:consumer-verification}

Items (I1)--(I5) of Proposition~\ref{prop:source-exact-replacement} are
verified here in a fixed order.  The proofs give the expanded equation-level
form of
Proposition~\ref{prop:v2-consumer}, while keeping the operator replacement
distinct from a comparison of final epsilon powers.

\subsection{Dictionary for the exceptional top call}

The notation in this section has the following fixed meaning.
\begin{center}
\small
\begin{tabular}{@{}p{0.17\textwidth}p{0.73\textwidth}@{}}
\toprule
symbol & meaning in the Source-B version-2 consumer\\
\midrule
\(\widetilde{\mathcal M}\)
 & the complete molecule after adjoining the exceptional top layer to all
 lower layers;\\
\(\widetilde H\)
 & the root set of the complete molecule;\\
\(H'\)
 & the initial cumulant root set at the bottom;\\
\(\mathsf X\)
 & the number of new particle lines in the lower molecule;\\
\(m\)
 & the number of new lines inserted in the exceptional top layer;\\
\(r\)
 & the size of the distinguished root cluster already present at the top
 interface;\\
\(\iota\)
 & \(0\) when \(r>0\) and \(1\) when \(r=0\);\\
\(\widetilde{\mathsf X}\)
 & \(\mathsf X+m-\iota\), the total number of new lines in the joined
 molecule;\\
\(\rho\)
 & the lower generalized-link and recollision parameter
 \( |H'|+\mathsf Z-\mathsf X\);\\
\(\widetilde\rho\)
 & the same parameter for the complete joined molecule;\\
\(N_4^1,N_4^2\)
 & the numbers of degree-four full components produced respectively in the
 top and lower cutting calls, including periodic good $\{4\}$ singletons;
 good status never deletes their $\eps^{-(d-1)}$ baseline;\\
\(N_{\rm ee}\)
 & the number of empty-end charges in the source ledger;\\
\(N_{\neg N3}\)
 & the number of elementary components outside the selected normal
 degree-three family;\\
\(N_{\rm good}\)
 & the number of protected ordinary or periodic good objects;\\
\(G\)
 & the top root set in the disjoint factorization, with
 \( |G|=r+\iota\);\\
\(K\)
 & the lower root set in the disjoint factorization; it is unrelated to the
 Jacobi response operator \(K_{\mathcal N}\) of the main manuscript,
 Section~\ref{app:positive-network}.\\
\bottomrule
\end{tabular}
\end{center}

All counters are evaluated after every \emph{permanent} protected
pair-to-singleton recut has performed the simultaneous transfer
\eqref{eq:simultaneous-ledger-transfer}, including the possible first
$N_{\rm good}$ credit of a formerly nongood ordinary pair.  A temporary
$\{33{\rm A}\}\leftrightarrow\{3\}+\{2\}$ refinement used only inside an
owner or conditional-measure calculation changes no counter.  Thus no atom
is present in two analytic counters.  Provenance remains stored by the
immutable certificate identifier but does not create an additional epsilon
factor.

\subsection{The insertion operator}

Let \(K_{\rm pkt}\subset\widetilde{\mathcal M}\) be the unique protected
two-sublayer packet returned by Theorem~\ref{thm:p38j}.  We reserve the
subscript \({\rm pkt}\) to avoid confusion with the lower root set \(K\).
On one positive cell \(\kappa\), write
$N_{\partial,K_{\rm pkt},\kappa}$ for the fixed source-boundary null set
selected in Definition~\ref{def:appD-selected-conditional-representative},
and define
\begin{equation}\label{eq:appF-packet-insertion}
 \begin{aligned}
 \mathfrak P_\kappa Q
 &:=
 \1_{N_{\partial,K_{\rm pkt},\kappa}^c}(z_{\rm comp})
 \sum_{j=1}^{M(d,P_0)}
 \int_{\mathfrak B_{K_{\rm pkt},\kappa}}
 Q(\Phi_{K_{\rm pkt},\kappa,j}(b),z_{\rm comp})\\
 &\qquad{}\times
 \1_{\mathcal D_{\widetilde{\mathcal M}}}
  (\Phi_{K_{\rm pkt},\kappa,j}(b),z_{\rm comp})
 W_{\kappa,j}^{\rm rem}(b)\dd b .
 \end{aligned}
\end{equation}
The complement variables \(z_{\rm comp}\) are parameters in this
definition; here this symbol denotes the full tuple of complement and
fixed-end boundary parameters at the insertion.  At a fixed
$z_{\rm comp}$, both the packet operator and the
output operator below mean the selected Borel conditional representatives
of Definition~\ref{def:appD-selected-conditional-representative}.  After
integration in $z_{\rm comp}$ they are the original complete Source-B
operators.  With this convention, choose the default radial FCT atlas used
by every sharp downstream bound in this section.  Then
Proposition~\ref{prop:appD-arbitrary-kernel-coarea}, with
\(\gamma_\kappa=2(d-1)\), gives
\begin{equation}\label{eq:appF-packet-insertion-bound}
 \overline J_{K_{\rm pkt},\kappa}
  (Q\1_{\mathcal D_{\widetilde{\mathcal M}}})
 \le
 \cL_\eps^C\eps^{-(d-1)}\eps^{2(d-1)}\epstar^{-2(d-1)}
 \mathfrak P_\kappa Q .
\end{equation}
Here the bar denotes the selected conditional packet kernel at the fixed
\(z_{\rm comp}\); after integration in \(z_{\rm comp}\) it is the original
unbarred complete Source-B operator.
This is a positive linear operator inequality.  No absolute value, supremum
over a boundary state, or factorization of \(Q\) is introduced.  The finite
sum is not a uniqueness convention: it contains every regular solution in
the retained fibre by
Theorem~\ref{thm:appD-uniform-algebraic-multiplicity}.

\begin{lemma}[Compatibility with cutting shadows]
\label{lem:appF-shadow-compatibility}
Every recovery branch \(\Phi_{K_{\rm pkt},\kappa,j}\) commutes with the
Source-B shadow substitution in \eqref{eq:source-cutting-identity}.  The
nontrivial assertion concerns every original bond or free end \(e\) whose
state occurs in the packet kernel and every fixed end \(f\) created at the
packet--complement interface.  For every existing branch \(j\),
\[
 z_e\circ\Phi_{K_{\rm pkt},\kappa,j}
 =
 z_{\operatorname{Sh}(e)}
 \circ\Phi_{K_{\rm pkt},\kappa,j},
 \qquad
 z_f\circ\Phi_{K_{\rm pkt},\kappa,j}
 =
 z_{\operatorname{Sh}(f)}
 \circ\Phi_{K_{\rm pkt},\kappa,j}.
\]
An original edge wholly in the complement satisfies the same commutation
identity trivially because \(\Phi_{K_{\rm pkt},\kappa,j}\) is the identity
on its state.
No assertion is made that a free end belonging to a Source-B simple pair
lies in the range of \(\operatorname{Sh}\).
Consequently \(\mathfrak P_\kappa\) may be inserted at the position of the
packet in the reverse cutting order without changing the arguments of the
positive source kernel.
\end{lemma}

\begin{proof}
We verify one legal cut.  In the no-overlap case
\(\mathtt{cutting\_0}\), a broken original bond is represented by the new
free end, and the paired fixed end has that same free end as its shadow.
This is exactly the substitution used in
\eqref{eq:source-cutting-identity}.

For an overlap subsegment, use the notation of Source B,
\(\mathtt{def.cutting}\).  In case \(\mathtt{cutting\_2}\), all original
bonds \(e_0,\ldots,e_q\) of the subsegment have the new bypass bond \(e_+\)
as their next shadow.  The two ends created at each intermediate O-atom
form a simple pair.  The free member of that pair need not itself be in the
range of \(\operatorname{Sh}\); Source B,
\(\mathtt{prop.match\_ends}\), instead gives
\[
 \operatorname{Sh}(\operatorname{Ba}(n,e))
 =\operatorname{Sh}(f)=e_+,
\]
where \(f\) is its fixed mate.  Hence every occurrence of an old subsegment
state, including an occurrence read through a simple-pair incidence, is
read from \(z_{e_+}\).

In case \(\mathtt{cutting\_3}\), the broken chain
\(e_0,\ldots,e_q\) has free representative \(e_+\) at the atom in the
free side and fixed end \(e_-\) at the terminal C-atom.  Directly from the
definition,
\[
 \operatorname{Sh}(e_i)=e_+,\qquad
 \operatorname{Sh}(e_-)=e_+,
 \]
up to subsequent shadow updates.  Intermediate simple pairs are handled by
the fixed-mate formula just displayed.  In case
\(\mathtt{cutting\_4}\), the same statement holds for
\(e_0,\ldots,e_{q-1}\), while the old terminal end \(e_q\) is also updated
to \(e_+\); every new fixed end again has shadow \(e_+\).  These are all
possibilities in \(\mathtt{cutting\_2}\)--\(\mathtt{cutting\_4}\).

Now compose the one-cut statements in cutting order.  A later cut merely
replaces the current representative by its next shadow, so induction gives
the final \(\operatorname{Sh}\).  On each side the branch
\(\Phi_{K_{\rm pkt},\kappa,j}\) performs the same transport translations
and elastic involutions on the packet states after this naming
substitution.  Complement and fixed-end states are not pivots and remain
parameters.  Therefore composition with \(\Phi_{K_{\rm pkt},\kappa,j}\)
preserves every displayed equality, including all three overlap cases.
\end{proof}

\subsection{Proof of item (I1)}

\begin{proposition}[Complete-kernel and time-owner property]
\label{prop:appF-I1}
The packet satisfies item (I1) of
Proposition~\ref{prop:source-exact-replacement}.
\end{proposition}

\begin{proof}
Take for \(Q\) one entire positive term
\(Q_{Z,\mathcal W,\mathcal V}\) from
\eqref{eq:source-positive-Q}, including its generalized-link indicators,
Maxwellian factors, remainders, fixed-end variables and the complete
complement.  Inequality \eqref{eq:appF-packet-insertion-bound} applies
because it assumes only nonnegativity.  Definition~\ref{def:appD-output-kernel}
uses zero extension into a product box and the fixed branch set
\(\{1,\ldots,M(d,P_0)\}\), both chosen before the fixed-boundary values are
supplied.  Thus Source-B sees a finite direct sum of admissible positive
kernels.  Its constants absorb \(M(d,P_0)\), while its velocity boxes and
boundary-independence requirement are unchanged.

Put \(S=S_0\cup R_{\rm birth}\) as in
Lemma~\ref{lem:new-line-time-owner-matching}, and at the packet-fibre stage
put
\[
 \mathcal A_{\rm pkt}:=M_t\setminus S.
\]
For fixed \(t_S\), using the selected conditional representative fixed
before the positive kernel is chosen,
Lemma~\ref{lem:appD-fibre-time-owner} retains the indicator of the exact
existential projection
\[
 \{t_S:\exists\,t_{\mathcal A_{\rm pkt}}\text{ with }
   (t_{\mathcal A_{\rm pkt}},t_S)\in
   \mathcal D_{\widetilde{\mathcal M}}\},
\]
where, at the fibrewise packet stage, all owner times are literal fixed
coordinates and $\mathcal A_{\rm pkt}\cap S=\varnothing$.  In the later
measure-level estimate the one-atom ordinary support-good owner times
$S_{\rm ord}$ and the short-gap two-atom owner times $S_{\rm sg}$ have been
consumed by the post-envelope common-measure estimates.  The retained
set is then $S_{\rm reg}=S\setminus(S_{\rm ord}\cup S_{\rm sg})$ and the omitted set is
$\mathcal A_t=M_t\setminus S_{\rm reg}$; the latter consists of the incoming
time--normal root pivots, the other nonowner times, and the already consumed
$S_{\rm ord}\cup S_{\rm sg}$ times.  Thus the two projections occur at different stages and
never use two overlapping sets as simultaneous coordinate blocks.  The
matching count
\eqref{eq:complete-owner-cardinality} and the weight--time estimate
Lemma~\ref{lem:time-owner} then give the complete
\(\vartheta^{\mathsf X+m-\iota}\) owner factor; in particular the original
lower factor \(\vartheta^{\mathsf X}\) is present.  All retained packet
birth times are literal output coordinates, while new lines outside the
packet tree are charged to the ordinary owners \(S_0\), so no unproved
birth-time injection is being used.  This proves every part of (I1).  The
analytic order inside that lemma is the Source-B order: collision weights and
the Maxwellian envelope first produce the pointwise majorant while the normal
degree-three velocities remain live; ordinary one-atom owner-good
intersections, short-gap common-two-time blocks and periodic free-coordinate
blocks are then integrated in one interleaved
reverse topological order of the complete cutting dependency graph, using
Lemma~\ref{lem:periodic-certificate-free-coordinate} at the periodic vertices.
Only then are the simultaneous normal velocity volume and the remaining
projected owner times paid.  Thus no conditional good domain loses a variable
on which it depends.  The normal family is formed after the certificate
ledger.  Lemma~\ref{lem:global-innovation-partition} splits one complete
square full-pair velocity transformation into exactly the four groups
$\mathcal H_{\rm ord}\oplus\mathcal H_{\rm per}\oplus
\mathcal H_{\rm nor}\oplus\mathcal H_{\rm rem}$.  A marked degree-two birth
is a transparent relay: it owns and consumes no birth innovation.  Its
coefficient $h_p$ remains live through that relay and is assigned, by the
complete $2d$-dimensional pair joint group, to exactly one of the four
displayed groups; it is consumed once by that group's genuine downstream
owner and nowhere else.  This includes the temporary degree-three member of
a refined $\{33{\rm A}\}$ certificate without identifying its transported
mixed direction with a particle slot.  Hence the birth relay,
good-component factors and Source-B simultaneous volume never integrate the
same coordinate twice.
\end{proof}

\subsection{Proof of item (I2)}

Order all protected calls by the moment at which their owner is first cut.
Inside one top call use the full order constructed in
Theorem~\ref{thm:p38j}; Lemma~\ref{lem:packet-special-induction} supplies the
initial order inside the selected connected component.  After protecting
\(K_{\rm pkt}\), the cleanup consists of:
\begin{enumerate}[label=\textup{(\roman*)}]
\item the remaining components inside the selected component, in the
DOWN/UP restart order of Lemma~\ref{lem:packet-special-induction};
\item the transversal set, processed by DOWN;
\item the remaining upper components, processed by UP;
\item the lower components, processed by DOWN;
\item any good recut, performed at the identifier's original position and
then removed before restart.
\end{enumerate}

\begin{lemma}[Acyclic dependency order]
\label{lem:appF-dependency-order}
Let \(\mathcal C_1,\ldots,\mathcal C_J\) be the elementary components
produced by the packet and cleanup, in cutting order.  Direct a dependency
edge \(\mathcal C_i\to\mathcal C_j\) when the integration domain of
\(\mathcal C_j\) uses a fixed end whose free shadow lies in
\(\mathcal C_i\).  Then \(i<j\), so the dependency graph is acyclic.  In
reverse integration order, the domain of a good component depends only on
variables that have not yet been integrated.
\end{lemma}

\begin{proof}
A fixed end is created only when an already selected component is cut free.
Its free shadow therefore belongs to a component selected strictly earlier.
Protection removes every active representative of the protected atoms
before a restart, so a later component cannot create a backward dependency
through the same certificate.  The TOP/BOTTOM orientation in the cleanup
ensures that a boundary bond is always on the causal side required by UP or
DOWN.  Hence every dependency edge follows the cutting order.  Reversing
that order gives exactly the condition preceding
\eqref{eq:source-iterated-integral}.
\end{proof}

\begin{proposition}[Admissible protected cleanup]
\label{prop:appF-I2}
The packet satisfies item (I2) of
Proposition~\ref{prop:source-exact-replacement}.
\end{proposition}

\begin{proof}
Theorem~\ref{thm:p38j} constructs the packet by legal Source-B free cuts.
Lemma~\ref{lem:appF-shadow-compatibility} identifies the corresponding
shadow variables, and Lemma~\ref{lem:appF-dependency-order} gives the
required integration order.  Every exceptional-top cleanup component is
nonfull by Theorem~\ref{thm:p38j}(ii), whose proof uses
Lemma~\ref{lem:packet-special-induction} on the selected connected component;
hence no second
degree-four full baseline is created by the top cleanup.  Lower DOWN
components are accounted for separately by $N_4^2$ below.  The indicator decomposition count is
\(C^{|M|}\cL_\eps^C\), as stated in Theorem~\ref{thm:p38j}.  These are
precisely the assertions in (I2).
\end{proof}

\subsection{Proof of item (I3): attached branch}

In the attached branch, cutting the top packet creates at least one fixed
end in the lower molecule.  Let its connected lower components be
\(\mathcal L_1,\ldots,\mathcal L_q\).  The Source-B CH incidence identity is
\begin{equation}\label{eq:appF-lower-incidence}
 q+N_{\rm ee}\le|\widetilde H|.
\end{equation}
A connected DOWN input produces a degree-four full component if and only if
it has no fixed end.  At least one \(\mathcal L_j\) contains the fixed end
created by the packet cut; therefore
\begin{equation}\label{eq:appF-lower-full-strict}
 N_4^2\le q-1.
\end{equation}
Adding \eqref{eq:appF-lower-incidence} and
\eqref{eq:appF-lower-full-strict} gives
\begin{equation}\label{eq:appF-attached-strict-baseline}
 N_4^2+N_{\rm ee}\le|\widetilde H|-1.
\end{equation}

\begin{lemma}[Exactly one packet baseline]
\label{lem:appF-one-baseline}
The normalized packet contribution in the attached branch contains exactly
one factor \(\eps^{-(d-1)}\).  No exceptional-top cleanup component contributes
a second top negative baseline; lower full components remain in $N_4^2$.
\end{lemma}

\begin{proof}
If the packet has \(P\) physical lines, its source normalization is
\(\eps^{-(d-1)P}\).  The \(P-1\) upper-tree contacts give
\(\eps^{(d-1)(P-1)}\), leaving \(\eps^{-(d-1)}\).  The two selected lower contacts
are the two positive factors in the packet excess and are not upper-tree
normalizations.  By Proposition~\ref{prop:appF-I2}, every exceptional-top
cleanup component is nonfull.  The lower full components are not discarded:
\eqref{eq:appF-attached-strict-baseline} includes them through $N_4^2$ and
makes the remaining source baseline exponent nonnegative.
\end{proof}

\begin{proposition}[Attached consumer]
\label{prop:appF-I3}
The packet satisfies item (I3) of
Proposition~\ref{prop:source-exact-replacement}.
\end{proposition}

\begin{proof}
Insert \eqref{eq:appF-packet-insertion-bound} into the complete source
ledger.  Lemma~\ref{lem:appF-one-baseline} and
\eqref{eq:appF-attached-strict-baseline} give
the following bound.  Here the displayed good exponent is assembled only
once: write
\[
 N_{\rm good}=N_{\rm good,own}+N_{\rm good,rest},
 \qquad
 N_{\rm good,own}=|S_{\rm ord}|+|S_{\rm per}|
                         +|\mathscr G_{\rm sg}|.
\]
The $S_{\rm ord}$ factors come from the common conditional-measure
intersection in Lemma~\ref{lem:owner-marked-good-component}; the
$\mathscr G_{\rm sg}$ factors come from the common-two-time clause of that
lemma, one per short-gap component identifier; the
$S_{\rm per}$ factors come once from the different-lift tube
\eqref{eq:periodic-owner-tube-product}, while their time coordinates remain
in the regular projection.  Together these are exactly the
$\eta_*^{|S_{\rm ord}|+|S_{\rm per}|+|\mathscr G_{\rm sg}|}$ factors in
\eqref{eq:owner-source-majorant}.  The ordinary Source-B and periodic
certificate integrations supply only the remaining $N_{\rm good,rest}$
factors in the same complete reverse-topological loop
\eqref{eq:complete-good-dependency-product}; they are not postponed past the
normal-family volume extraction.  Their product is denoted below by
$\eps^{N_{\rm good}/(9d)}$; it is not a second application of a good estimate
to an owner-marked component.
An ordinary support-good full $\{4\}$ owner belongs to $S_{\rm ord}$: the
same common-measure lemma is applied to its baseline-normalized local measure,
while its unique $\eps^{-(d-1)}$ normalization and $N_4$ charge remain in the
source exponent.  An ordinary full $\{4\}$ owner which is not good belongs to
$S_{\rm reg}$ and receives only its projected time factor.  Thus all ordinary
full components are represented in the owner count without manufacturing a
good credit.
If an owner-marked periodic certificate is a full good $\{4\}$ singleton,
its tube supplies the displayed good excess but its unique $\eps^{-(d-1)}$
normalization remains in the $N_4$ exponent.  In the attached consumer the
top cleanup is nonfull, so such a singleton, if present, is a lower component
and is counted in $N_4^2$.  In the disjoint consumer it belongs to the lower
factor; the rooted top cleanup again contains no full singleton.  Hence the
owner split never hides a degree-four baseline.
\begin{multline*}
 \mathcal I\mathcal N
 \le
 \vartheta^{\widetilde{\mathsf X}}
 \eps^{(d-1)(|\widetilde H|-1-N_4^2-N_{\rm ee})}
 \sigma_K\\
 {}\times
 L^{C^*(\widetilde\rho+N_{\rm ee}+N_{\neg N3})}
 \eps^{N_{\rm good}/(9d)+|H'\setminus Z|/4}.
\end{multline*}
This is \eqref{eq:attached-consumer}.  The first epsilon exponent is
nonnegative.  The source incidence count absorbs the good and initial-link
logarithms, leaving only \(L^{C^*(\rho+r)}\) and a bounded
\(3\Gamma\)-loss.  Lemma~\ref{lem:appE-attached-time} handles the time
power, and \eqref{eq:appE-packet-fixed-credit} handles the packet scale.
The result is \eqref{eq:source-top-second-interpoland} for every fixed
integer \(d\ge4\).
\end{proof}

\subsection{Proof of item (I4): disjoint branch}

If the exceptional top molecule has no bond to the lower molecule, the
source variables split into disjoint sets.  The time indicator is the
product of the two layer indicators because all top times already lie in a
strictly later layer.  The collision kernels, bottom data and source
normalizations also split.  Tonelli therefore gives the exact factorization
\begin{equation}\label{eq:appF-disjoint-factorization}
 \mathcal I\mathcal N(\widetilde{\mathcal M},\widetilde H,H')
 \le
 \mathcal I\mathcal N(\widetilde{\mathcal M}_\ell,G,\varnothing)
 \mathcal I\mathcal N(\widetilde{\mathcal M}_{<\ell},K,H').
\end{equation}
The inequality, rather than equality, records only the positive
decompositions already made before this step.

\begin{lemma}[Top rooted normalization]
\label{lem:appF-top-rooted-normalization}
For \(G\) with \(|G|=r+\iota\), the top factor in
\eqref{eq:appF-disjoint-factorization} satisfies
\begin{equation}\label{eq:appF-top-rooted-bound}
 \mathcal I\mathcal N(\widetilde{\mathcal M}_\ell,G,\varnothing)
 \le
 \vartheta^{m-\iota}
 \eps^{(d-1)(|G|-1)}\sigma_K.
\end{equation}
\end{lemma}

\begin{proof}
Apply the rooted packet estimate
\eqref{eq:two-landing-rooted}.  Each root other than the unique translation
root contributes \(\eps^{d-1}\), giving \(\eps^{(d-1)(|G|-1)}\).
The top time power is not charged solely to packet birth contacts.
Lemma~\ref{lem:new-line-time-owner-matching} splits the genuine top
new-line owners between the retained packet set \(R_{\rm birth}\) and the
ordinary top set \(S_{0,\rm top}\), and proves the exact disjoint count
\eqref{eq:disjoint-owner-cardinality}:
\[
 |S_{0,\rm top}|+|R_{\rm birth}|\ge m-\iota .
\]
Lemma~\ref{lem:time-owner}, applied to the top factor of the disjoint
source measure, keeps the exact projected time indicator and assigns one
literal Lebesgue owner coordinate to each member of these two disjoint
sets.  Its ledger \eqref{eq:vartheta-owner-count} therefore supplies
\(\vartheta^{m-\iota}\).  Multiplication with the rooted spatial
normalization gives \eqref{eq:appF-top-rooted-bound}.  Thus no injection of
all top new lines into the possibly smaller packet tree, and no
unsupported product-domain replacement or enlargement of the exact
projected time set is used; the scale-correct projected-volume estimate in
Lemma~\ref{lem:time-owner} is used exactly once.
\end{proof}

\begin{proposition}[Disjoint consumer]
\label{prop:appF-I4}
The packet satisfies item (I4) of
Proposition~\ref{prop:source-exact-replacement}.
\end{proposition}

\begin{proof}
The lower factor in \eqref{eq:appF-disjoint-factorization} is untouched, so
it obeys \eqref{eq:source-disjoint-lower}, or equals one when the lower
molecule is empty.  For the top factor use
\eqref{eq:appF-top-rooted-bound}.  If \(r=0\), then
\(\iota=1\), \(m\ge2\), and its strength is
\(\vartheta^{m-1}\eps^{2(d-1)-1/(100d)}\).  If \(r>0\), then
\(\iota=0\) and its strength is
\(\vartheta^m\eps^{(d-1)(r+1)-1/(100d)}\).  The exact comparisons in
\eqref{eq:appE-disjoint-time} and
\eqref{eq:appE-disjoint-target-power} prove
\eqref{eq:source-disjoint-top} in both cases.
\end{proof}

\subsection{Proof of item (I5)}

The source first interpoland has three routes.  Only the attached route uses
the ordinary theorem on the original uncut molecule.  In the disjoint
nonempty route the two factors have disjoint variable sets and are multiplied
only through the exact Source-B factorization.

\begin{lemma}[Three first-interpoland routes]
\label{lem:appF-first-interpoland}
The first interpoland \eqref{eq:source-top-first-interpoland} is obtained as
follows.
\begin{enumerate}[label=\textup{(\roman*)}]
\item Attached: apply the ordinary molecule theorem and its packet-free
small-window owner input,
Proposition~\ref{prop:packet-free-small-window-consumer}, to the original
uncut molecule.
\item Disjoint, lower nonempty: combine the exact factorization
\eqref{eq:appF-disjoint-factorization}, the lower estimate
\eqref{eq:source-disjoint-lower}, and the top estimate
\eqref{eq:appF-top-rooted-bound}.
\item Disjoint, lower empty: apply
\eqref{eq:appF-top-rooted-bound} alone.
\end{enumerate}
\end{lemma}

\begin{proof}
For (i), return to \eqref{eq:source-associated-integral} before choosing the
exceptional packet.  The ordinary theorem has the source root set, fixed
ends, and parameters of the attached first interpoland.  Proposition
\ref{prop:packet-free-small-window-consumer} supplies its complete owner
exhaustion before the ordinary interpolation weakens the time exponent to
the one in \eqref{eq:source-attached-first}.

For (ii), the variables and indicators split exactly as proved before
\eqref{eq:appF-disjoint-factorization}.  If $r=0$, then $R=\rho$, $m\ge2$,
and the lower factor supplies
$\vartheta^{\mathsf X/9}\eps^{1/(12d)+4(C_{14}^*)^{-1}\rho}$, while the top
factor supplies
$\vartheta^{m-1}\eps^{2(d-1)-1/(100d)}$.  These powers dominate
$\vartheta^{(\mathsf X+m)/10}
 \eps^{1/(18d)+(C_{14}^*)^{-1}R}$.
If $r>0$, the top factor instead supplies
$\vartheta^m\eps^{(d-1)(r+1)-1/(100d)}$; its positive margin is linear in $r$ and
absorbs both the required $r$-coefficient and the per-root logarithmic
losses, uniformly for $r+m\le2\Lambda_\ell$.  Thus the product gives the first
interpoland.  This multiplication is legitimate because it is the exact
factorization, not a product of overlapping marginal bounds.

For (iii), the complete integral is the rooted top packet.
Lemma~\ref{lem:appF-top-rooted-normalization} gives
$\eps^{2(d-1)-1/(100d)}$ when $r=0$ and
$\eps^{(d-1)(r+1)-1/(100d)}$ when $r>0$, with the
stronger time power.  These exponents dominate the first interpoland after
$C_{14}^*$ is selected.
\end{proof}

\begin{proposition}[Completion of the version-2 replacement]
\label{prop:appF-I5-completion}
The packet satisfies item (I5).  Consequently all hypotheses of
Proposition~\ref{prop:source-exact-replacement} hold and the replacement
proves \eqref{eq:source-trunc-output}.
\end{proposition}

\begin{proof}
Propositions~\ref{prop:appF-I3} and~\ref{prop:appF-I4} give the second
interpoland in the attached and disjoint branches.  Lemma~\ref{lem:appF-first-interpoland}
gives the first interpoland independently.  All finite packet chart counts
and bounded exceptional constants are covered by
Lemma~\ref{lem:appE-subpower-ledger}; per-root losses are handled by
\eqref{eq:appE-qeff-log}.  Lemma~\ref{lem:source-two-regime-envelope}
therefore yields \eqref{eq:source-top-target}.  Substitution in
\eqref{eq:source-top-error-reduction} gives
\(\|f_s^{\rm err}\|_1\le\eps^{1/(20d)}\), which is
\eqref{eq:source-trunc-output} for every fixed integer \(d\ge4\).
\end{proof}

\subsection{Sharp packet scale and complementary coarse bound}

Every default bound in this section uses the sharp packet scale
\(\sigma_K=\cL_\eps^C\eps^{2(d-1)}\epstar^{-2(d-1)}\), supplied by
Theorem~\ref{thm:packet-fct}.  The comparisons in Appendix~\ref{app:parameters} were already
valid for the larger coarse scale with \(\epstar^{-2d}\), so the downstream
hierarchy is a monotone consequence and no parameter has been retuned.  The
full collision-free-line chord argument remains an independent proof of that
coarse scale; it is not used to justify the radial-conjugate FCT theorem.

\section*{Statements and declarations}

\noindent\textbf{Funding.} None.

\medskip
\noindent\textbf{Competing interests.} The author declares no competing interests.

\medskip
\noindent\textbf{Data availability.} No datasets were generated or analysed during the current study.

\medskip
\noindent\textbf{Materials availability.} The mathematical source materials and exact version locks used by the proof are identified in the manuscript and Online Resource~1.

\medskip
\noindent\textbf{Code availability.} No research software or numerical code was used to obtain the mathematical results.

\medskip
\noindent\textbf{Use of large-language-model tools.} Large-language-model tools were used to assist with language editing, structural and document reorganization, literature and submission-guideline lookup, cross-reference checking, mechanical invariance comparisons, adversarial consistency review, and preparation of submission and archival materials. The author reviewed the outputs and takes full responsibility for all mathematical statements, proofs, citations, and conclusions.

\end{document}